\documentclass{amsart}

\usepackage{amsmath, amssymb, amsfonts, mathtools} 
\usepackage{bm} 
\usepackage{mathbbol} 
\usepackage{mathrsfs} 
\usepackage{stmaryrd} 
\usepackage{verbatim} 
\usepackage{tabularx} 

\usepackage[all]{xy} 
\usepackage{tikz} 
\usetikzlibrary{decorations.pathmorphing} 
\tikzset{snake it/.style={decorate, decoration=snake}}

\usepackage{tikz-cd} 
\usepackage{standalone} 

\usepackage{caption}
\usepackage{subcaption}

\usepackage{babel}

\usepackage[
	bookmarks=true,
	linktocpage=true,
	bookmarksnumbered=true,
	breaklinks=true,
	pdfstartview=FitH,
	hyperfigures=false,
	plainpages=false,
	naturalnames=true,
	colorlinks=true,
	pagebackref=true,
	pdfpagelabels
	]{hyperref}

\hypersetup{
	linkcolor=blue,
	citecolor=blue,
	urlcolor=blue,
	colorlinks=true
}

\usepackage[capitalise,noabbrev,nosort]{cleveref}

\makeindex

\let\oldtocsection=\tocsection
\let\oldtocsubsection=\tocsubsection
\let\oldtocsubsubsection=\tocsubsubsection
\renewcommand{\tocsection}[2]{\hspace{0em}\oldtocsection{#1}{#2}}
\renewcommand{\tocsubsection}[2]{\hspace{1em}\oldtocsubsection{#1}{#2}}
\renewcommand{\tocsubsubsection}[2]{\hspace{2em}\oldtocsubsubsection{#1}{#2}}

\usepackage{todo} 
\usepackage{ulem} 

\newcommand{\greg}[1]{\todo[G.\!]{#1}}

\DeclareSymbolFontAlphabet{\mathbba}{bbold}

\newtheorem{theorem}{Theorem}[section]
\newtheorem*{theorem*}{Theorem}
\newtheorem{lemma}[theorem]{Lemma}
\newtheorem{corollary}[theorem]{Corollary}
\newtheorem{proposition}[theorem]{Proposition}

\theoremstyle{definition}
\newtheorem{definition}[theorem]{Definition}
\newtheorem{remark}[theorem]{Remark}
\newtheorem{example}[theorem]{Example}
\newtheorem{convention}[theorem]{Convention}
\newtheorem{notation}[theorem]{Notation}

\newcommand{\defeq}{\stackrel{\mathrm{def}}{=}}

\newcommand{\mc}[1]{\mathcal{#1}}
\newcommand{\ms}[1]{\mathscr{#1}}
\newcommand{\R}{\mathbb{R}}
\newcommand{\C}{\mathbb{C}}

\newcommand{\Z}{\mathbb{Z}}

\newcommand{\id}{\mathrm{id}}
\newcommand{\e}{\mathbf{e}}
\newcommand{\f}{\mathbf{f}}
\newcommand{\bb}{\mathbf{b}}
\newcommand{\interval}{\mathbb{I}}

\newcommand{\Or}{{\rm Det}}

\newcommand{\im}{\text{im}}
\newcommand{\into}{\hookrightarrow}

\newcommand{\sms}{\smallsmile}
\newcommand{\pf}{\pitchfork}
\newcommand{\bd}{\partial}
\newcommand{\td}{\tilde}
\newcommand{\Hom}{\textup{Hom}}
\newcommand\red[1]{\textcolor{red}{#1}}
\newcommand\blue[1]{\textcolor{blue}{#1}}

\newcommand{\uW}{\underline{W}}
\newcommand{\udW}{\underline{\partial W}}
\newcommand{\uV}{\underline{V}}
\newcommand{\uM}{\underline{M}}
\newcommand{\uN}{\underline{N}}
\newcommand{\xr}{\xrightarrow}
\newcommand{\xl}{\xleftarrow}
\newcommand{\uX}{\underline{X}}
\newcommand{\cok}{{\rm cok}}
\newcommand{\Cre}{\mathrm{Cre}}
\newcommand{\aug}{\mathbf{a}}
\newcommand{\vertices}{{\rm Vert}} 
\newcommand{\jinterval}{\mathbb{J}}
\newcommand{\onto}{\twoheadrightarrow}
\newcommand{\diag}{\mathbf{d}}
\newcommand{\Ext}{\text{Ext}}
\newcommand{\mf}{\mathfrak}
\newcommand{\mr}[1]{\mathring{#1}}
\newcommand{\pipe}{|} 

\title[Foundations of geometric cohomology]{Foundations of geometric cohomology:\\ From co-orientations to product structures}

\author[G. Friedman]{Greg Friedman}
\address{Department of Mathematics, Texas Christian University, USA}
\email{\href{mailto:g.friedman@tcu.edu}{g.friedman@tcu.edu}}
\thanks{This work was partially supported by a grant from the Simons Foundation (\#839707 to Greg Friedman)}

\author[A. Medina-Mardones]{Anibal M. Medina-Mardones}
\address{Mathematics Department, Western University, Canada}
\email{\href{mailto:anibal.medina.mardones@uwo.ca}{anibal.medina.mardones@uwo.ca}}

\author[D. Sinha]{Dev Sinha}
\address{Mathematics Department, University of Oregon, USA}
\email{\href{mailto:dps@uoregon.edu}{dps@uoregon.edu}}

\begin{document}
	\begin{center}
		||PREPRINT||
		\vskip 20pt
	\end{center}


\begin{abstract}
	This manuscript develops a geometric approach to ordinary cohomology of smooth manifolds, constructing a cochain complex model based on co-oriented smooth maps from manifolds with corners.
	Special attention is given to the pull-back product of such smooth maps, which provides our geometric cochains with a partially defined product structure inducing the cup product in cohomology.
	A parallel treatment of homology is also given allowing for a geometric unification of the contravariant and covariant theories.
\end{abstract}
	\maketitle
	\tableofcontents
	\index{Table of Contents}

\section{Introduction}\label{S: introduction}

We fully develop a geometric approach to ordinary homology and cohomology on smooth manifolds, with classes represented by smooth maps from manifolds to our target manifold of interest.
Such a development is in line with thinking about homology dating back to Poincar\'e and Lefschetz, but is also a time-honored approach to cohomology through intersection theory and Thom classes.
Thom's seminal work on bordism theory showed that not all homology classes can be represented by pushing forward fundamental classes of manifolds, and this necessitates a broader notion of representing manifold.
The one we utilize here is that of manifolds with corners, the smallest category containing manifolds with boundary and which is closed under transverse pullback.
The former property is needed to define homologies, and we use the latter property to define multiplicative structures.
The study of these multiplicative structures -- a partially defined ring structure on cochains and a module over it on chains -- is also in line with the classical perspective, and their rigorous development is a central goal of this project.
While these structures give rise to the usual cup and cap products on cohomology and homology, they are markedly different at the chain and cochain level when compared to simplicial, cubical, or singular approaches, providing greater geometric access and enjoying strict graded commutativity when defined.

Our approach using manifolds with corners was initiated by Lipyanskiy in the unpublished manuscript \cite{Lipy14}, and we extend that work to a full treatment of cochains and multiplicative structures, along the way filling in details needed.
Chains of degree $i$ are represented by maps from compact oriented $i$-dimensional manifolds with corners, while cochains are represented geometrically by proper and co-oriented smooth maps from manifolds with corners, with associated degree of the cochain given by the codimension of the map.
While it is expected that ordinary homology can be captured in this way -- after all simplices and cubes are manifolds with corners -- the technicalities in this setting are surprising.
For example, the boundary of a boundary of a manifold with corners is not empty or even naively zero as a chain or cochain.
Following Lipyanskiy, we simply quotient by the sorts of ``trivial'' chains that arise in the image of the boundary squared, but then our chains and cochains are themselves equivalence classes.
Another quotient by ``degenerate'' chains and cochains is needed to ensure homology and cohomology theories that satisfy the dimension axiom.
Once one is working with such equivalence classes and needs, for example, transversality to define products, truly substantial difficulties arise.
Indeed, at one point we doubted the existence of a well-defined multiplication.

Approaches to ordinary homology and cohomology were an active area of development eighty years ago, and indeed we highlight that our work is in some sense parallel to de Rham theory, allowing one to make calculations and invoke geometry for manifolds, rather than rely on transcendental approaches to cochains as formal linear duals.
Compared with de Rham theory, our present work has the key advantage of being defined over the integers.
Slightly more recent developments with similar goals of describing homology, and especially cohomology, more geometrically include Goresky's work on geometric homology and cohomology of Whitney stratified objects \cite{goresky1981stratified} and the book by Buoncristiano, Rourke, and Sanderson \cite{buoncristiano1976homology}.
But this project is more contemporary than one might assume.
Symplectic geometers have been revisiting these ideas as a parallel to work on Floer theory, with both Lipyanskiy \cite{Lipy14} and Joyce \cite{Joyc15} offering versions.
In a similar vein, Kreck's ``differential algebraic topology'' \cite{Krec10} provides homology and cohomology on smooth manifolds using maps from \textit{stratifolds}, a certain kind of singular space.
Even the foundations of a theory of transversality for manifolds with corners suitable for our work has only been worked out in recent decades by Margalef-Roig and Outerelo Dominguez \cite{MaDo92} and Joyce \cite{Joy12}.

While the idea of homology as represented by fundamental classes of submanifolds or, more generally, manifold mappings is quite familiar, historically it was only in some corners of 20th century geometric topology, such as those noted above, that \textit{co}homology classes were also represented by appropriate maps from manifolds.
Such cochains are geometric objects in their own right, which partially evaluate on chains through intersection.
A great benefit to such thinking is that the classical operations of algebraic topology, such as cup and cap products, can be described \textit{at the level of chains and cochains} by simple geometric operations based on intersection, without recourse to chain approximations to the diagonals, Alexander--Whitney maps, or other such combinatorics.
This again is reminiscent of the original thinking about such products in terms of intersection, and parallels modern work such as intersection theory in the PL category as in McClure's \cite{McC06}.

The trade-off for such a pleasant description is that these intersections are not always defined; they require transversality.
This limitation is also classically anticipated by the famous commutative cochain problem.
Loosely speaking, no integral cochain construction computing ordinary cohomology can be made canonically into a (graded) commutative ring.
Since the process of forming intersections is commutative, the ring structure it induces in our theory cannot be fully defined.
We find the trade-off worthwhile, and in work building on these foundations \cite{FMS-flows} we have already ``married'' multiplicatively the theory we develop here, which is commutative and partially defined, and the one defined by cubical cochains with the Serre product, which is not commutative but everywhere defined.
For concrete applications of ideas in line with our viewpoint, we mention Cochran's work relating Seifert surface intersections with Massey products in the context of Milnor invariants for links \cite{cochran1990milnor}, as well as the work of our third-named author and his collaborator on group cohomology through configuration spaces, where key calculations are made through intersections \cite{giusti2012symmetric, giusti2021alternating}.

Lipyanskiy's manuscript \cite{Lipy14}, on which we build, gives a fairly thorough account of geometric homology, but a much more lightly sketched account of geometric cohomology, which leaves several major theorems unproven.
Some other expected results are not stated at all, including an isomorphism between geometric and ordinary cohomology, either as graded abelian groups or as rings.
So one of our main goals is to give a thorough account, with detailed proofs, of geometric homology and cohomology, with our primary focus on geometric cohomology, both because Lipyanskiy's account of this requires more filling in and also because cohomology with its algebra structure is of more interest to us.
In addition to research applications, we find these ideas helpful in teaching graduate students as, for example, cohomology of projective spaces follows from linear algebra, and pushforward or umkehr maps are defined just by taking images.
At the research level, we found Lipyanskiy's work thanks to Mike Miller, while working on \cite{FMS-flows} and looking for a rigorous foundation to geometrically model the cup product.
Ultimately our main aim is to obtain a full -- but partially defined -- $E_\infty$-algebra structure on geometric cochains, with the $E_\infty$-structure ``resolving'' partial definedness (though itself being partially defined!) rather than non-commutativity.
We also see plenty of room for development to make geometric cochains more broadly applicable, for example conjecturally to CW complexes with smooth attaching maps.
Our careful treatment here is offered to facilitate such work.

\section*{Acknowledgments}\index{acknowledgments}

The authors thank Mike Miller for pointing us to \cite{Lipy14} and Dominic Joyce for answering questions about his work.

\section*{Conventions}\index{conventions|textbf}

Throughout we will denote the dimension of a manifold represented by an upper case character by the corresponding lower case character, for example, $\dim(M) = m$, $\dim(V) = v$, etc. \index{conventions!dimension labels}

All maps are assumed to be smooth, in the sense to be defined in \cref{D: smooth}, unless stated otherwise.
For a map $f \colon W \to M$ and $x \in W$, when we want to specify the derivative map restricted to $T_xW$ and it is not otherwise clear from context, we write $D_xf$.
As noted below in \cref{N: implicit notation}, once we have established a specific map $f \colon W \to M$, we will often abuse notation by writing simply $W$ to refer to $W$ together with its map to $f$; in this sense, $W$ will be treated as an object over $M$.

For evaluation of a tensor products of cochains on a tensor product of chains, we follow the convention $(\alpha \otimes \beta)(x\otimes y) = \alpha(x)\beta(y)$.\index{conventions!evaluation of tensors}
This convention is used to define the cup product in Munkres \cite[Section 60]{Mun84}, Hatcher \cite[Section 3.2]{Hatc02}, and Spanier \cite[Section 5.6]{Span81}, though it disagrees with the conventions in some other sources, such as Dold \cite[Section VII.7]{Dol72}.

The symbol $\sqcup$ denotes disjoint union.\index{$Z$@$\sqcup$} 

\section{Manifolds with corners}\label{S: manifolds with corners}

In this section, we provide an overview of manifolds with corners, which are the main geometric objects in the definitions of geometric chains and cochains.
Our main references for this material are Joyce \cite{Joy12} and Margalef--Roig and Outerelo Dominguez \cite{MaDo92}. In \cite{MaDo92}, manifolds with corners (which they simply call ``manifolds'') may be infinite-dimensional and may fail to be Hausdorff, so there is some work involved in seeing that their definition agrees with Joyce's when one assumes, as we shall, that all spaces are finite-dimensional, Hausdorff, and second countable.

By \textbf{smooth} we always mean differentiable of all orders.
Throughout the paper, all manifolds and maps will be in the smooth category unless noted otherwise.

\begin{definition}\label{D: smooth}
	If $A \subset \R^n$ and $B \subset \R^m$ are any subsets, we say that $f \colon A \to B$ is \textbf{smooth}\index{smooth!map of subset of Euclidean space|textbf} if it extends to a smooth map from an open neighborhood of $A$ to $\R^m$.
	We say that $f$ is a \textbf{diffeomorphism}\index{diffeomorphism|textbf} if $f$ is a homeomorphism and both $f$ and $f^{-1}$ are smooth\footnote{Joyce requires $n = m$ in his definition for $f$ to be a diffeomorphism.
	But suppose $f$ is a diffeomorphism as defined here and, without loss of generality, $n>m$.
	Let $\R^m \subset \R^n$ in the usual way.
	Then as $f$ extends to a smooth map from a neighborhood of $A$ to $\R^m$, we can certainly consider $f$ as extending to a smooth map from a neighborhood of $A$ to $\R^n$.
	Similarly, as $f^{-1}$ extends to a smooth map from a neighborhood of $B$ in $\R^m$ to $\R^n$, we can extend $f^{-1}$ to a neighborhood $B$ in $\R^n$ by precomposing with the projection $\R^n \to \R^m$.
	So diffeomorphisms in our sense can be made into diffeomorphisms in Joyce's sense, and obviously vice versa.}.
\end{definition}

With this definition, the notions of smooth charts and atlases for manifolds in the standard setting can be extended to define (smooth) manifolds with corners; see \cite[Section 2]{Joy12}: An $n$-dimensional chart $(U,\phi)$ of the space $W$ has domain $U$ an open subset of $\R^n_k = [0,\infty)^k \times \R^{n-k} \subset \R^n$,\index{$R$@$\R^n_k$} and $\phi$ is a homeomorphism from $U$ to $\phi(U) \subset W$.
Two $n$-dimensional charts $(U,\phi)$ and $(V,\psi)$ are compatible if $\psi^{-1}\phi \colon \phi^{-1}(\phi(U) \cap \psi(V)) \to \psi^{-1}(\phi(U) \cap \psi(V))$ is a diffeomorphism of subsets of $\R^n$.
An $n$-dimensional atlas for $W$ is a family of pairwise compatible $n$-dimensional charts that cover $W$, and an $n$-dimensional manifold with corners is a second countable Hausdorff space with a maximal $n$-dimensional atlas.
Our condition that $W$ be second countable is slightly more restrictive than Joyce's definition, which requires only that $W$ be paracompact.

The reason for our additional restriction is that we choose to work entirely with subspaces of $\R^\infty$ in order to have a set of such objects, so rather than directly employ Joyce's definition, we define manifolds with corners for the purposes of this text as follows.

\begin{definition}\label{D: MWC}
	An $n$-dimensional \textbf{manifold with corners}\index{manifold with corners|textbf} $W$ is a subspace of some $\R^N \subset \R^\infty$ such that each point of $W$ possesses a neighborhood diffeomorphic to an open subset of $\R^n_k$ for some $k$.
\end{definition}

We note that these local diffeomorphisms provide charts $\phi \colon U \subset \R^n_k \to W$, and these charts are compatible, with $\psi^{-1}\phi \colon \phi^{-1}(\phi(U) \cap \psi(V)) \to \psi^{-1}(\phi(U) \cap \psi(V))$ being a diffeomorphism, even using Joyce's more restrictive notion of diffeomorphism.
These charts collectively give an atlas, and every atlas extends to a unique maximal one.
So our manifolds with corners are also manifolds with corners by Joyce's definition.

On the other hand, Joyce notes in \cite[Remark 2.11 (see also Remark 6.3)]{Joy12} that his version of manifolds with corners agrees with that of Margalef--Roig and Outerelo Dominguez \cite{MaDo92} (except that they consider the much more general case of possibly infinite-dimensional manifolds and possibly of differentiable class $p<\infty$).
They show that when their manifolds with corners (which they simply call manifolds) are second countable, Hausdorff, and of the same finite dimension at all points, they have smooth closed embeddings into finite-dimensional Euclidean spaces \cite[Corollary 11.3.10]{MaDo92} (in particular, the embeddings are proper\footnote{Recall that a map between locally compact Hausdorff spaces is \textbf{proper}\index{proper map|textbf} if the inverse image of any compact set is compact.} maps \cite[Proposition 3.3.4]{MaDo92}).
Joyce already assumes his manifolds with corners to be paracompact Hausdorff in \cite[Definition 2.1]{Joy12}, so if we add the second countable condition to his definition,\footnote{\label{F: countable}For us this will not be a terrible restriction, as in this case being second countable is the same as being paracompact with a countable number of components.
Indeed, every regular second countable space is metrizable \cite[Theorem 34.1]{Mu00}, and hence paracompact, and clearly a locally connected second countable space must have a countable number of connected components.
Conversely, each connected component of a paracompact manifold with corners will be closed, paracompact \cite[Theorem 41.2]{Mu00}, and locally compact (via our charts), so it will be $\sigma$-compact \cite[Appendix A.1]{Sp79}.
A $\sigma$-compact space is Lindelof, and this implies that our manifold with corners can be covered by a countable number of charts, each of which is second countable, which implies second countability.
And then a countable union of second countable spaces is second countable.
So the only effect of our second countability condition is to rule out manifolds with corners in Joyce's sense that have uncountably many connected components.}
any one of Joyce's manifolds with corners is diffeomorphic to a manifold with corners in our sense.

As noted in \cite[Remark 2.11]{Joy12}, Joyce's definition of manifolds with corners agrees with those of Cerf \cite{Ce61}, Douady \cite{Doua61}, and Margalef--Roig and Outerelo Dominguez \cite{MaDo92}, and they also correspond to Melrose's t-manifolds in \cite{Melrose}.
As an alternative to citing \cite{MaDo92} in the preceding paragraph, we can also observe that Melrose shows in \cite[Proposition 1.14.1]{Melrose} that t-manifolds can be embedded into manifolds without boundary; so, when second countable, these can then be embedded into Euclidean space by the Whitney Embedding Theorem.

By modeling on $\R^n_k$, our category includes manifolds ($k = 0$ in all charts) and manifolds with boundary ($k \leq 1$ in all charts), as well as cubes and simplices, but not the octahedron, for example, as the cone on $[0,1] \times [0,1]$ is not modeled by any $\R^n_k$.

\begin{definition}
	A map $f \colon W \to M$ between manifolds with corners is {(weakly) \bf smooth}\index{smooth!map of manifold with corners|textbf} if whenever $(U,\phi)$, $(V,\psi)$ are charts for $W$ and $M$ respectively then
	$$\psi^{-1}f \phi \colon (f\phi)^{-1}\psi(V) \to V$$
	is smooth.

	The \textbf{tangent bundle}\index{tangent bundle|textbf} of a manifold with corners is the space of derivations of the ring of smooth real-valued functions.
	Analogously to smooth manifolds with boundary, if $(U,\phi)$ is a chart of $W$ then $D\phi$ takes the tangent space to $\R^n$ at $x \in U$ isomorphically to the tangent spaces of $W$ at $\phi(x)$.
	In particular, if $W$ is $n$-dimensional and $x \in W$, then the tangent space at $x$, denoted $T_xW$, is isomorphic to $\R^n$.
\end{definition}

\begin{remark}\label{R: weakly smooth}
	In \cite{Joy12}, Joyce reserves the word \textit{smooth} for weakly smooth maps that also satisfy an additional condition concerning how they interact with the boundaries of their codomains.
	When the codomain is a manifold without boundary, which will be our primary situation of interest, the weaker and stronger notions coincide.
	In the few other cases in which we need to consider maps whose codomains have boundary, in particular boundary immersions and projections of pullbacks, the maps will also satisfy Joyce's stronger condition, though we will never need to utilize this explicitly.
	Thus we will feel justified in simply using the word ``smooth'' throughout, referring the reader to \cite[Definition 3.1]{Joy12} for the full definition.
\end{remark}

While the notion of embeddings of manifolds with corners into manifolds with corners is not straightforward (and it is not clear that these notions in \cite{Joy12} and \cite{MaDo92} agree), the following observation concerning embeddings of manifolds with corners into manifolds without boundary is reassuring.
Suppose a manifold with corners $V$ is embedded in a manifold without boundary $M$, meaning, as in the classical sense, that there is a map $f \colon V \to M$ that is a homeomorphism onto its image and such that $Df$ is injective at each point of $V$; see \cite[Definition 3.3.1 and Theorem 3.2.6]{MaDo92}.
Let $x$ be a boundary point of $V$, and let $(U,\phi)$ be a chart for $V$ at $x$ such that $\phi(0)=x$.
Then, by the definition of smooth maps, $f \phi$ extends smoothly to a neighborhood of the origin in $\R^v$.
As $f$ is an embedding, $D(\phi f)_0$ must be injective, so by classical differential topology there is a neighborhood of $0$ in $\R^v$ on which $f\phi$ extends to an embedding with $x$ in the interior \cite[Section 1.3]{GuPo74}.
In other words, the embedding of $V$ can be extended (locally in some neighborhood of $x$) to an embedding of an ordinary manifold without boundary, and the neighborhood of $x$ in $V$ sits inside a Euclidean neighborhood of this embedding as the diffeomorphic image of $\R^v_k$ for some $k$.

\begin{notation}\index{conventions!dimension labels}
	Our default notation for manifolds with corners will be capital letters with the corresponding lower case letter denoting the dimension.
	In other words, $\dim(V) = v$, $\dim(W) = w$, $\dim(M) = m$, etc.
	We generally reserve $M$ for a manifold without boundary.
\end{notation}

\subsection{Boundaries}\label{S: boundaries}

We next need to describe boundaries of manifolds with corners.
Again see \cite[Section 2]{Joy12} for further details.

\begin{definition}
	A point $x$ in an $n$-dimensional manifold with corners $W$ has \textbf{depth}\index{depth|textbf} $k$ if there is a chart from an open subset of $\R^n_k$ which sends the origin to $x$.
	The set $S^k(W) \subseteq W$ of elements having depth~$k$ is called the \textbf{stratum of depth $k$}.\index{stratum|textbf}\index{depth|textbf}
	By \cite[Proposition 2.4.]{Joy12}, $S^k(W)$\index{$S^k(W)$} is an $n-k$ manifold without boundary.
\end{definition}

\begin{example}
	If $W$ is a smooth manifold with boundary in the classical sense, then $S^0(W)$ is its interior, $S^1(W)$ is its boundary, and $S^k(W) = \emptyset$ for $k>1$.
	If $S^k(W) = \emptyset$ for all $k>0$, then $W$ is a manifold without boundary.
\end{example}

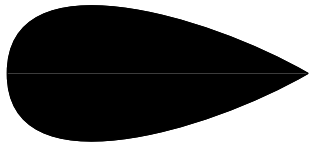
\begin{figure}[h]
		\begin{tikzpicture}[scale=4, rotate=90]
	\coordinate (a) at (0,1.2);
	\coordinate (b) at (0,0.2);

	\draw[out=120, in=180, fill, opacity=.2] (b) to (a);
	\draw[out=120, in=180] (b) to (a);
	\draw[out=60, in=0, fill, opacity=.2] (b) to (a);
	\draw[out=60, in=0] (b) to (a);
	\end{tikzpicture}
	\caption{The 2-dimensional ``teardrop'' is a manifold with corners whose boundary inclusion is not injective.}
	\label{F: teardrop}
\end{figure}

When $W$ is a general manifold with corners, the boundary is more naturally a space equipped with a map to $W$, rather than a subspace of $W$.
The reason can be seen, for example, in the teardrop space of \cref{F: teardrop}.
The intuitive boundary here would be homeomorphic to the circle but have a corner at the point of the teardrop.
However, this would not be a smooth manifold with corners.
Instead, we consider the boundary to be homeomorphic to the closed interval together with a map taking both endpoints to the vertex of the teardrop.
To explain in more detail, we have the following definitions.

\begin{definition}
	A \textbf{local boundary component}\index{boundary!local boundary component (of manifold with corners)} $\beta$ of $W$ at $x \in W$ is a consistent choice of connected component $\mathbf{b}_V$ of $S^1(W) \cap V$ for any neighborhood $V$ of $x$, with consistent meaning that $\mathbf{b}_{V \cap V'} \subset \mathbf{b}_{V} \cap \mathbf{b}_{V'}$.
\end{definition}

Since this notion is local, the number of such components is determined by the depth.
Considering the origin in $\R^n_k$, for any $k \geq 0$, points having depth $k$ have exactly $k$ local boundary components.
As another example, letting $\interval = [0,1]$ and so $\interval^3$ the 3-cube, the set $S^1(\interval^3)$ consists of the interiors of two-dimensional faces, $S^3(\interval^3)$ is the set of eight corners, and any sufficiently small neighborhood of a corner intersects exactly three of the two-dimensional faces.

\begin{definition}\label{D: MWC boundary}
	Let $W$ be a manifold with corners.
	Define its \textbf{boundary}\index{boundary!of manifold with corners|textbf} $\bd W$ to be the space of pairs $(x, \bb)$ with $x \in W$ and $\bb$ a local boundary component\footnote{Note that if $x \in S^0(W)$ then $x$ does not have any local boundary components and so does not appear in such a pair.} of $W$ at $x$.
	Define $i_{\bd W} \colon \bd W \to W$ by sending $(x,\bb)$ to $x$.
\end{definition}

In the teardrop example, $i_{\bd W} \colon \bd W \to W$ takes $\interval$ to the topological manifold boundary with both endpoints going to the unique point of $S^2(W)$.
For $\interval^3$, the boundary consists of six closed two-dimensional squares each mapping homeomorphically to a face of the cube.
In general, $|i_{\bd W}^{-1}(x)|$ is equal to the depth of $x$.

As in Douady \cite{Doua61}, $\bd W$ is itself a manifold with corners, and the boundary map $i_{\bd W}$ is a smooth immersion \cite[Theorem 3.4]{Joy12}. Note that $S^1(W)$ will always be diffeomorphic to the interior of $\bd W$, i.e.\ $S^0(\bd W) = S^1(W)$.
Inductively, we let $\bd^k W$\index{$bd$@$\bd^k$|textbf} denote $\bd (\bd^{k-1} W)$ with $\bd^0 W = W$, and we let $i_{\bd^k W}$ denote the composite of the boundary maps sending $\bd^k W$ to $W$. The map $i_{\bd^k W}$ takes $S^a(\bd^k W)$ onto $S^{a+k}(W)$, but not, in general, injectively.

\begin{remark}\label{R: bd diff}
	If $f \colon V \to W$ is a diffeomorphism of manifolds with corners then it must, in particular, take components of $S^i(V)$ diffeomorphically onto corresponding components of $S^i(W)$, and, consequently, given a local boundary component $\mathbf{b}$ at a point $x \in V$, the map $f$ picks out a corresponding local boundary component, say $f_*(\mathbf{b})$ in $W$. We obtain a diffeomorphism $f_\bd \colon \bd V \to \bd W$ by $(x,\mathbf{b})\mapsto (f(x),f_*(\mathbf{b}))$ with $fi_{\bd V} = i_{\bd W} f_\bd$.
\end{remark}

If $W$ is a classical smooth manifold with boundary, then we have $\bd^2 W=\emptyset$, but this will not generally be the case for a manifolds with corners.
However, given a map $W \to M$ with either $W$ oriented or the map co-oriented, then $\bd^2 W$ does come equipped with a natural orientation- or co-orientation-reversing $\Z_2$ action.
This will be explained below and is a key component in showing that $\bd$ is suitable for defining boundary maps of geometric chain and cochain complexes.

N.B. When no confusion is likely to occur, we will sometimes abuse notation and use $\bd W$ also to refer to the image $i_{\bd W}(\bd W) \subset W$, which is the boundary of $W$ as a topological manifold in the usual sense.

In the following example, which is a mild generalization of \cite[Example 7.3]{Joy12}, we see $\bd^2 W$ immerse into $W$ in an interesting way.

\begin{example}\label{boundary}
	Consider the quotient $Q$ of $S^n \times \R^2_2$ by the diagonal $\Z_2$ action, where $\Z_2$acts antipodally on $S^n$ and acts by permuting the two coordinates
	of $\R^2_2$.
	The projection of $Q$ onto $S^n / \Z_2 = \R P^n$ defines a fiber bundle with fiber $\R^2_2$.
	Local coordinates can then be used to endow $Q$ with the structure of a manifolds with corners.

	The subspace $S^1(Q)$ is diffeomorphic to $S^n \times (0,\infty)$, and $S^2(Q)$ is diffeomorphic to $\R P^n$.
	The boundary $\bd Q$ is the quotient of $S^n \times \bd \R^2_2$ by $\Z_2$, which is diffeomorphic to
	$S^n \times [0,\infty)$.
	Thus $\bd^2 Q$ is $S^n$, which maps to $S^2(Q)$ by the standard quotient by antipodal action.
\end{example}

In general, Proposition~2.9 of \cite{Joy12} identifies $\bd^k W$ with the set of points $(x,\bb_1,\ldots,\bb_k)$ with $x \in W$ and the $\bb_i$ providing an ordered $k$-tuple of distinct local boundary components of $W$.

Special cases of manifolds with corners, including
manifolds with faces or manifolds with embedded corners \cite{Joy12}, steer clear of the interesting boundary phenomena of \cref{boundary}.
Although, a more restrictive notion should suffice for our application, Lipyanskiy develops geometric cohomology in the current generality, and we appreciate Joyce's careful treatment of transversality in this category, so we use their definitions.

We conclude this section by observing that the product $V \times W$ of two manifolds with corners is naturally a manifold with corners with, by \cite[Proposition 2.12]{Joy12},
$$\bd (V \times W) = (\bd V \times W) \sqcup (V \times \bd W).$$

\subsection{Transversality}

Transversality of smooth maps will play a key role, as this is the condition that assures that intersections or, more generally, pullbacks of manifolds (with corners) are also manifolds (with corners).
Recall that in the classical setting if $f \colon V \to M$ and $g \colon W \to M$ are smooth maps of manifolds without boundary then we say that $f$ and $g$ are \textbf{transverse}\index{transversality!classical} if whenever $f(x) = g(y) = z$ for some $x \in V$, $y \in W$, and $z \in M$ we have $Df(T_xV)+Dg(T_yW) = T_z M$.
We here discuss the extension of transversality to manifolds with corners, though only in the case where $M$ is without boundary; as noted in \cite[Remark 6.3]{Joy12}, in this special case Joyce's transversality of manifolds with corners is equivalent to the formulation of transversality in \cite[Section 7.2]{MaDo92}.
More general versions of transversality can be found in \cite[Section 6]{Joy12}.

\begin{definition}{\cite[Special case of Definition 6.1]{Joy12}}
	Let $f \colon V \to M$ and $g \colon W \to M$ be smooth maps of manifolds with corners to a manifold without boundary.
	We say $f$ and $g$ are \textbf{transverse}\index{transversality!manifolds with corners|textbf} if whenever $x \in V$ and $y\in W$ with $f(x) = g(y) = z$, if $x\in S^j(V)$ and $y \in S^k(W)$ then $Df|_{S^j(V)}(T_xS^j(V))+Dg|_{S^k(W)}(T_yS^k(W)) = T_zM$.
	In this case, we often abuse notation by saying that $V$ and $W$ are transverse, leaving the maps tacit.
	Note that the empty map $\emptyset \to M$ is transverse to all other maps.
\end{definition}

Note that in our version of the definition, we use maps of the form $f|_{S^j(V)}$, which are classical smooth maps, as the strata $S^j(V)$ are smooth manifolds. Joyce takes a slightly different approach, first noting that if $f:V\to M$ is a smooth map of manifolds with corners then there is an induced map $Df$ ``in the usual way'' \cite[Definition 3.2]{Joy12} and then restricting $Df$ to $T_xS^j(V)\subset T_xV$ (see also \cite[Proposition 2.4]{Joy12}). However, some thinking through the definitions shows that the resulting spaces $D(f|_{S^j(V)})(T_xS^j(V))$ and $Df(T_xS^j(V))$ are identical, both being the derivations that act on functions $h:M\to \R$ by taking the derivatives of $hf$ in the directions tangent to $S^k(V)$. We somewhat prefer our formulation as it expresses transversality entirely in terms of properties of maps of smooth manifolds (without corners).

The following special case of the definition will also sometimes be useful:

\begin{definition}\label{D: regular value}
If $M = \R$ and $g$ is the inclusion of the point $p$ into $\R$, we will say that $p$ is a \textbf{regular value}\index{regular value|textbf} for $f \colon V \to M$ when $f$ and $g$ are transverse.
\end{definition}

While this definition of transversality is given in terms of the behavior of $f$ and $g$ on the strata $S^j(V)$ and $S^k(W)$, it is sometimes useful to have a reformulation in terms of the boundaries $\bd^j V$ and $\bd^kW$.
This is the content of \cref{L: simple trans} below, for which we need a further definition.

\begin{definition}\label{D: naive transversality}
Let $f \colon V \to M$ and $g \colon W \to M$ be maps from manifolds with corners to a manifold without boundary.
We say that $f$ and $g$ are \textbf{naively transverse}\index{transversality!naive|textbf} if whenever $x \in V$ and $y\in W$ with $f(x) = g(y) = z$ we have $Df(T_xV)+Dg(T_yW) = T_z M$. Note that $x$ and $y$ may be boundary points of $V$ and $W$, the maps $Df$ and $Dg$ being those of \cite[Definition 3.2]{Joy12}.
\end{definition}

\begin{lemma}\label{L: simple trans}
	Let $f \colon V \to M$ and $g \colon W \to M$ be maps from manifolds with corners to a manifold without boundary.
	Then $f$ and $g$ are transverse if and only if $fi_{\bd^j} \colon \bd^jV \to M$ and $gi_{\bd^k} \colon \bd^kW \to M$ are naively transverse for all $j,k$.
\end{lemma}

Note that, a priori, the latter is a stronger condition as it imposes conditions not just on the interior of strata but on their closures.

\begin{proof}
	First suppose $fi_{\bd^j}$ and $gi_{\bd^k}$ are naively transverse for all $j,k$.
	Suppose $x \in S^j(V)$ and $y \in S^k(W)$ for some fixed $j,k$ and $f(x) = g(y)$.
	The preimage of $x$ under $i_{\bd^j}$ consists of $j!$ points in $S^0(\bd^jV)$, and $i_{\bd^j}$ maps a neighborhood of each such preimage point diffeomorphically to a neighborhood of $x$ in $S^j(V)$, and similarly for $y$.
	Let $\psi_x$ be the inverse diffeomorphism from a neighborhood of $x$ in $S^j(V)$ to a neighborhood of one of the preimages.
	Then $f|_{S^j(V)} = fi_{\bd^j}\psi_x$ in a neighborhood of $x$, and similarly for $y$.
	Since $fi_{\bd^j}$ and $gi_{\bd^k}$ are naively transverse, and $\psi_x$ and $\psi_y$ are diffeomorphisms, it follows that $f|_{S^j(V)}$ and $g|_{S^k(W)}$ are transverse at $f(x) = g(y)$.

	Conversely, suppose $f|_{S^j(V)}$ and $g|_{S^k(W)}$ are transverse for all $j,k$, and suppose $x \in \bd^jV$ and $y \in \bd^kW$ for some fixed $j,k$ with $fi_{\bd^j}(x) = gi_{\bd^k}(y) = z$.
	Furthermore, suppose $i_{\bd^j}(x) \in S^a(V)$ and $i_{\bd^k}(y) \in S^b(W)$, which implies $x \in S^{a-j}(\bd^j V)$.
	By focusing on local charts, there is a neighborhood of $x$ in $\bd^j(V)$ whose intersection with $S^{a-j}(\bd^j V)$ maps diffeomorphically via $i_{\bd^j}$ onto a neighborhood of $i_{\bd^j}(x)$ in $S^a(V)$, and analogously for $y$.
	Thus, $f|_{S^a(V)}i_{\bd^j}|_{S^{a-j}(\bd^j V)}$ and the analogous $g|_{S^b(W)}i_{\bd^k}|_{S^{b-k}(\bd^k W)}$ are transverse at $f(x) = g(y)$ as they precompose transverse maps with local diffeomorphisms.
	But the image of $D_x(fi_{\bd^j})$ contains the image of $D_x(f|_{S^j(V)}i_{\bd^j}|_{S^{a-j}(\bd^j V)})$ and similarly for $y$, and thus the images of $D_x(fi_{\bd^j})$ and $D_y(gi_{\bd^k})$ must also span $T_{z}M$.
	Therefore, $fi_{\bd^j}$ and $gi_{\bd^k}$ are naively transverse at $f(x)$.
\end{proof}

\subsubsection{Achieving transversality}

Throughout this text, we will need a series of increasingly more general results guaranteeing that we can make certain maps transverse to each other.
We begin here with two relatively simple cases that will be used in \cref{S: basic properties} to show that geometric cohomology is contravariantly functorial with respect to continuous maps of manifolds without boundary.
We put these first transversality results here, as some of the fundamental ideas will be used in a number of places throughout the text.

Our transversality theorems will be built primarily using some basic tools that can be found in \cite[Section 2.3]{GuPo74}.
In particular, we record the following results, referring the reader to \cite[Section 2.3]{GuPo74} for the proofs\footnote{We rephrase the statements of these theorems slightly to better fit our context and notation.}:

\begin{theorem}[Transversality Theorem]\label{T: GP transversality}\index{transversality!Transversality Theorem!Guillemin-Pollack}
	Suppose $F \colon X \times S \to Y$ is a smooth map of manifolds, where only $X$ has boundary, and let $Z$ be any boundaryless submanifold of $Y$.
	If both $F$ and $F|_{\bd X \times S}$ are transverse to $Z$, then for almost every $s \in S$, both $F(-,s) \colon X \to Y$ and $F(-,s)|_{\bd X} \colon \bd X \to Y$ are transverse to $Z$.
\end{theorem}

\begin{theorem}[$\epsilon$-Neighborhood Theorem]\label{T: epsilon neighborhood}\index{epsilon@$\epsilon$-neighborhood theorem}
	For a compact boundaryless manifold $Y$ in $\R^M$ and a positive number $\epsilon$, let $Y_\epsilon$ be the open set of points in $\R^M$ with distance less than $\epsilon$ from $Y$.
	If $\epsilon$ is sufficiently small, then each point $w \in Y_\epsilon$ possesses a unique closest point in $Y$, denoted $\pi(w)$.
	Moreover, the map $\pi \colon Y_\epsilon \to Y$ is a submersion.
	When $Y$ is not compact but properly embedded in $\R^m$, these statements remain true after replacing the constant $\epsilon$ with a smooth function $\epsilon \colon Y \to (0, \infty)$ and letting $Y_\epsilon = \{w \in \R^m \mid |w-y|<\epsilon(y) \text{\ for some\ } y \in Y\}$.
\end{theorem}

\begin{remark}\label{R: epsilon neighborhood}
	The closest point characterization of the retraction map $\pi$ in \cref{T: epsilon neighborhood} is not discussed in \cite{GuPo74}.
	However, it is observed in \cite[Problem 6-5]{Lee13} (see \cite{KK24} for the argument) that there exists a tubular neighborhood $U$ of $Y$ in $\R^M$ with retraction $r \colon U \to Y$ such that for each $z \in U$, $r(z)$ is the unique point of $Y$ closest to $z$.
	The construction of $\epsilon$-neighborhoods in \cite{GuPo74} essentially finds them within any tubular neighborhood using the same retraction map for both the tubular neighborhood and the $\epsilon$-neighborhood, so by using tubular neighborhoods as described above, we may assume the closest-point characterization holds for our $\epsilon$-neighborhoods even when $Y$ is not compact.
\end{remark}

These theorems are used in \cite{GuPo74} to prove the following basic transversality result:

\begin{theorem}[Transversality Homotopy Theorem]\index{transversality!Transversality Homotopy Theorem}
	Let $X$ be a smooth manifold with boundary.
	For any smooth map $f \colon X \to Y$ and any boundaryless submanifold $Z$ of the boundaryless manifold $Y$, there exists a smooth map $g \colon X \to Y$ homotopic to $f$ such that $g$ is transverse to $Z$ and $g|_{\bd X}$ is transverse to $Z$.
\end{theorem}

Among other generalizations as we progress, we will extend these results to maps that go from manifolds with corners to manifolds without boundary.
We will also often require that the homotopies take a special form, i.e.\ that they are \textit{universal} homotopies as defined in \cref{D: universal homotopy}.
See, for example, \cref{P: perturb transverse to map}.

When working with transversality, it is often much easier to work with the case where one of the maps is an embedding.
The following technical results will facilitate that and help us to extend the results of \cite{GuPo74} from transversality with respect to submanifolds to transversality with respect to maps.

\begin{lemma}\label{L: proper product}
	Let $f \colon X \to Y$ be an arbitrary map and $g \colon X \to Z$ a proper map.
	Then the map $(f,g) \colon X \to Y \times Z$ given by $(f,g)(x)=(f(x),g(x))$ is proper.
\end{lemma}

\begin{proof}
	Let $K \subset Y \times Z$ be compact.
	Let $\pi_Y,\pi_Z$ be the projections of $Y \times Z$ to $Y$ and $Z$. Then
	\begin{align*}
		(f,g)^{-1}(K) &= \{x \in X \mid (f(x),g(x)) \in K \}\\
		&\subset \{ x \in X \mid g(x) \in \pi_Y(K)\}\\
		&= g^{-1}(\pi_Y(K)).
	\end{align*}
	As $K$ is compact and $g$ is proper, $g^{-1}(\pi_Y(K))$ is compact.
	So $(f,g)^{-1}(K)$ is a closed subset of a compact set, and so compact.
\end{proof}

\begin{corollary}\label{C: embed V}
Let $g \colon W \to M$ be a map from a manifold with corners to a manifold without boundary.
Then for some $k$ there exists a closed embedding $e \colon W \to M \times \R^k$ such that, if $\pi \colon M \times \R^k \to M$ is the projection, then $\pi e = g$.
\end{corollary}

\begin{proof}
	By \cite[Corollary 11.3.10]{MaDo92}, there is a closed embedding $\alpha \colon W \to \R^k$ for some $k$.
	The map $\alpha$ is proper, as if $K$ is a compact subset of $\R^k$, the inverse image $\alpha^{-1}(K)$ is homeomorphic to the intersection of $K$ with the closed set $\alpha(W)$.
	So, by \cref{L: proper product}, the map $(g,\alpha) \colon W \to M \times \R^k$ is proper.
	Further, $\alpha$ is an immersion, so is $(g,\alpha)$.
	In fact, by \cite[Theorem 3.2.6]{MaDo92} and the fact that $M$ and $\R^k$ do not have boundary, a map is an immersion in this context if and only if it is injective on tangent spaces.
	If $\pi_{\R^k} \colon M \times \R^k \to \R^k$ is the projection, the map $\pi_{\R^k}(g,\alpha)$ equals the embedding $\alpha$, so $D(g,\alpha)$ must be injective at each point.
	Similarly, the map $(g,\alpha)$ is itself injective because $\alpha$ is injective.
	So $(g,\alpha)$ is an injective proper immersion, and just as for ordinary manifolds, an injective proper immersion is a closed embedding \cite[Proposition 3.3.4]{MaDo92}.
	Finally, it is clear from the construction that $\pi (g,\alpha) = g$.
\end{proof}

\begin{lemma}\label{L: all transversality is wrt embeddings}
	Let $f \colon V \to M$ and $g \colon W \to M$ be smooth maps from manifolds with corners to a manifold without boundary.
	Let $e \colon W \to M \times \R^n$ be an embedding such that $\pi e = g$, where $\pi$ is the projection $M \times \R^n \to M$.
	Then $f$ and $g$ are transverse if and only if $e$ is transverse to $f \times \id_{\R^n} \colon V \times \R^n \to M \times \R^n$.
\end{lemma}

\begin{proof}
	It suffices to assume that $V$ and $W$ are without boundary.
	Otherwise we can apply the following argument to each pair of strata of $V$ and $W$.

	Suppose that $f$ and $g$ are transverse, i.e.\ that if $f(v) = g(w)$ then $Df(T_vV)+Dg(T_wW) = T_{f(v)}M$.
	For each $w \in W$, we can write $e(w) = (g(w),e_\R(w)) \in M \times \R^n$.
	Now suppose $w \in W$ and $(v,z) \in V \times \R^n$ such that $e(w) = (f \times \id_{\R^n})(v,z)$.
	Then we have $(g(w),e_\R(w)) = (f(v),z)$.
	The image of the derivative of $f \times \id_{\R^n}$ at such a point will span $Df(T_vV) \oplus T_z(\R^{n}) = Df(T_vV) \oplus \R^{n}$, while the derivative of $e$ will take $a \in T_w(W)$ to $Dg(a)+ De_{\R}(a)$.
	But the image of $D(f \times \id_{\R^n})$ already includes $0 \oplus \R^{n}$, so subtracting off the second summand, $D(f \times \id_{\R^{n}})(T_{(v,z)}(V \times \R^n))+De(T_wW)$ contains $Dg(a)$.
	It follows that $D(f \times \id_{\R^{n}})(T_{(v,z)}(V \times \R^n))+De(T_wW)$ contains $Df(T_vV) \oplus 0$, $Dg(T_wW) \oplus 0$, and $0 \oplus \R^n$.
	Since $f$ and $g$ are transverse and $D(f \times \id_{\R^{n}})(T_{(v,z)}(V \times \R^n))+De(T_wW)$ is a vector space, it therefore contains all of $T_{f(v)}M \oplus \R^n = T_{e(w)}(M \oplus \R^n)$.
	So $f \times \id_{\R^n}$ and $e$ are transverse.

	Next suppose $f \times \id_{\R^n}$ and $e$ are transverse and that $f(v) = g(w) \in M$.
	Suppose $e(w) = (g(w),z)$.
	Then $e(w) = (f \times \id_{\R^n})(v,z)$.
	So, by definition and assumption,
	\begin{equation}\label{E: Quillen transverse}
		D(f \times \id_{\R^{n}})(T_{(v,z)}(V \times \R^n))+De(T_wW) = T_{e(w)}(M \times \R^n) = T_{f(v)}M \oplus \R^n.
	\end{equation}
	As $\pi$ is a submersion, the image of this tangent space under $D\pi$ is all of $T_{f(v)}M$.
	But $(D\pi)(De) = D(\pi e) = Dg$, so $(D\pi \circ De)(T_wW) = Dg(T_wW)$.
	Furthermore, letting $\pi_V \colon V \times \R^n \to V$ be the projection, we have $(D\pi)(D(f \times \id_{\R^{n}})) = D(\pi(f \times \id_{\R^{n}})) = D(f\pi_V) = (Df)(D\pi_V)$, so, as $D\pi_V \colon T_{(v,z)}(V \times \R^n) \to T_vV$ is surjective, we have $(D\pi)(D(f \times \id_{\R^{n}}))(T_{(v,z)}(V \times \R^n)) = Df(T_vV)$.
	So applying $D\pi$ to equation \eqref{E: Quillen transverse}, we get $Df(T_vV)+Dg(T_wW) = T_{f(v)}M$, and $f$ is transverse to $g$.
\end{proof}

We now need one more particular technical result that will allow us to achieve transversality while keeping proper maps proper.

\begin{lemma}\label{L: nearby proper homotopy}\index{proper map!homotopy criterion}
	Let $X$ be a topological space and $M$ a manifold without boundary possessing a metric $d$.
	Suppose $f \colon X \to M$ is a proper map.
	Then there exists a continuous function $\varepsilon \colon X \to (0,\infty)$ such that if $g \colon X \to M$ satisfies $d(f(x),g(x)) < \varepsilon(x)$ for all $x \in X$, then $g$ is proper and homotopic to $f$ by a proper homotopy.
	Furthermore, if $X$ is a manifold with corners and $f$ and $g$ are smooth and homotopic by a proper homotopy, then they are homotopic by a smooth proper homotopy.
\end{lemma}

\begin{proof}
	Let $C(X,M)$ be the set of continuous maps from $X$ to $M$.
	By \cite[Proposition 9.2.28]{MaDo92}, there is an open neighborhood $V^f$ of $f$ in the topology $T_s$ on $C(X,M)$ (defined in \cite[Proposition 9.2.1]{MaDo92}) such that every $g \in V^f$ is proper and homotopic to $f$ by a proper homotopy.
	Furthermore, by \cite[Proposition 9.3.9]{MaDo92} there is a local basis at $f$ in the $T_s$ topology consisting of the sets
	$$B_\varepsilon^d = \{g \in C(X,M) \mid d(f(x),g(x)) < \varepsilon (x) \text{ for all } x \in X\}$$
	as $\varepsilon$ ranges over all continuous functions $X \to (0,\infty)$.
	So, in particular, there is an $\varepsilon$ such that $B_\varepsilon^d \subset V^f$, which proves the first part of the lemma.
	When $X$ is a manifold with corners and $f$ and $g$ are smooth, we can take the proper homotopy between them to be proper and smooth by \cite[Proposition 9.2.35]{MaDo92}; note that $X$ satisfies the hypotheses of having partitions of unity by \cite[Corollary 1.5.14]{MaDo92}.
\end{proof}

We now apply the above results to demonstrate two transversality theorems that we will need in \cref{S: cohomology pullback} to obtain contravariant functoriality of geometric homology and cohomology.
We will see various fancier versions of these theorems below.
See, for example, \cref{P: ball stability,P: perturb transverse to map}.

\begin{theorem}\label{T: basic trans}\index{transversality!Transversality Theorem}
	Let $f \colon M \to N$ be a continuous map of manifolds without boundary, and let $g \colon V \to N$ be a smooth map from a manifold with corners.
	Then there is a homotopy $h \colon M \times I \to N$ such that $h(-,0) = f$ and $h(-,1)$ is smooth and transverse to $g$.
	Furthermore, if $f$ is proper, then there is such an $h$ so that both $h$ and $h(-,1)$ are proper.
\end{theorem}

\begin{proof}
	By \cref{C: embed V}, there is a closed embedding $e \colon V \to N \times \R^k$ that satisfies the hypotheses of \cref{L: all transversality is wrt embeddings}, so by that lemma it suffices to find a homotopy such that $h(-,0) = f$ and $h(-,1) \times \id_{\R^k} \colon M \times \R^k \to N \times \R^k$ is smooth and transverse to $e(V)$.

	We first find a homotopy from $f$ to a smooth map.
	This can be done by the smooth approximation theorem; see \cite[Corollary 9.2.31]{MaDo92}.
	Furthermore, if $f$ is proper, then we can take the smooth map and the homotopy between it and $f$ to be proper by \cite[Corollary 9.2.36]{MaDo92}.
	So we assume for the rest of the proof that this first homotopy has been completed and that $f$ is a smooth map.

	Next, we can think of $N$ as properly embedded in some $\R^K$ by the Whitney Embedding Theorem \cite[Section 1.8]{GuPo74}, and we give $N$ the corresponding metric.
	We let $N_\epsilon$ be an $\epsilon$-neighborhood of $N$ in $\R^K$ with submersion $\pi \colon N_\epsilon \to N$ (see \cref{T: epsilon neighborhood}).
	Let $D$ be the open unit ball in $\R^K$.
	We define a smooth composite map $H \colon M \times D \to \R^K \to N$ by
	$$H(x,s) = \pi(f(x)+ \eta(x)s)$$
	where $\eta \colon M \to (0,\infty)$ is defined as follows.

	First, so that $H(x,s)$ will be well defined, we want $f(x)+ \eta(x)s$ to land in $N_\epsilon$ for all $x$, and for this it is sufficient to have $0< \eta(x) < \epsilon(f(x))$.
	If $f$ is not proper or we do not need the last clause of the theorem to hold, any such $\eta$ will be sufficient.
	If $f$ is proper and we want the last clause of the theorem to hold, then we next choose a function $\varepsilon \colon M \to (0,\infty)$ to meet the requirements of \cref{L: nearby proper homotopy} with $d$ being the distance on $N$ chosen above.
	Letting $|\cdot|$ be the norm on $\R^K$, we then want to choose $\eta$ so that $|f(x) - H(x,s)| < \varepsilon(x)$.
	Clearly, $|f(x) - (f(x)+ \eta(x)s)| = |\eta(x)s| < \eta(x)$, and by definition $\pi$ takes points of $N_\epsilon$ to the nearest point in $N$. So, assuming $f(x)+ \eta(x)s \in N_\epsilon$, then $|f(x)+ \eta(x)s - \pi(f(x)+ \eta(x)s)| < \eta(x)$, as $f(x) \in N$ and so the distance from $f(x)+ \eta(x)s$ to $N$ must be $\leq |\eta(x)s| < \eta(x)$.
	So by the triangle inequality, $|f(x) - H(x,s)| < 2\eta(x)$.
	Thus if we take $\eta(x) < \min(\epsilon(f(x)), \varepsilon(x)/2)$, we will have $H$ well defined and $|f(x) - H(x,s)| < \varepsilon(x)$ for any $x \in M$ and $s \in D$ by \cref{L: nearby proper homotopy}.

	Now, we have $H(x,0) = \pi(f(x)) = f(x)$, and, since $\eta(x)>0$ for all $x$, the map $(x,s) \to f(x)+ \eta(x)s$ is a submersion onto its image in $\R^K$; in fact, for $x$ fixed, $H(x,s)$ maps $D$ linear onto a ball neighborhood of $f(x)$.
	As $\pi$ is also a submersion, so is $H$, and then $H \times \id_{\R^k} \colon M \times D \times \R^k \to N \times \R^k$ is a submersion.
	In particular, $H \times \id_{\R^k}$ is transverse to $e(V)$ and to $e(S^i(V))$ for each stratum $S^i(V)$ of $V$.
	We can now apply the Transversality Theorem, (\cref{T: GP transversality}, though with the parameter space $D$ in a slightly nontraditional location in the ordering of factors) to obtain that $H(-,s) \times \id_{\R^k}$ is transverse to $e(V) \times \R^k$ for almost every $s \in D$ and similarly for each $S^i(V)$.
	Since $V$ has a finite number of strata, and since the intersection of a finite number of sets whose complements have measure zero is another set whose complement has measure zero, it holds for almost every $s \in D$ that $H(-,s) \times \id_{\R^k}$ is transverse to every stratum of $V$.
	Now let $s_0$ be one such $s$, and let $h \colon M \times I \to N$ be given by $h(x,t) = H(x,ts_0)$.
	Then $h(x,0) = H(x,0) = \pi(f(x)) = f(x)$, and $h(-,1) = H(-,s_0)$ is smooth and transverse to every stratum of $V$ by \cref{L: all transversality is wrt embeddings}.
	Furthermore, if $f$ is proper, then by \cref{L: nearby proper homotopy} and our choice of $\eta$, the map $h(-,1)$ is proper and if necessary we can replace $h$ with a proper homotopy from $f$ to $h(-,1)$. \qedhere

\end{proof}

\begin{remark}\label{R: countable trans}
	In the statement of \cref{T: basic trans} we make $f$ transverse to $g \colon V \to N$, but the construction depends on $g$ only in determining the set of almost all possible $s\in D$ that can be chosen as $s_0$ to construct the homotopy.
	So as the countable intersection of sets of full measure will still have full measure, we can generalize the theorem to see that $f$ is homotopic to a map that is simultaneously transverse to any given countable collection of maps $V_i \to N$.
\end{remark}

Our second theorem will be applied below when $M$ is of the form $Y \times I$ for $Y$ a manifold without boundary in order to improve a homotopy between maps that are both transverse to another given map $g$ to a homotopy that is transverse to $g$.

\begin{theorem}\label{T: homotopy trans}\index{transversality!Transversality Theorem!relative}
	Let $f \colon M \to N$ be a smooth map from a manifold with boundary to a manifold without boundary, and let $g \colon V \to N$ be a smooth map from a manifold with corners.
	Suppose $\bd f \colon \bd M \to N$ is transverse to $g$.
	Then there is a smooth homotopy $h \colon M \times I \to N$ such that $h(-,0) = f$, $h(-,1)$ is smooth and transverse to $g$, and $h$ fixes $\bd M$, i.e.\ $h(x,t)=f(x)$ for all $x \in \bd M$ and all $t \in I$.
	Furthermore, if $f$ is proper, then there is such an $h$ so that both $h$ and $h(-,1)$ are proper.
\end{theorem}

\begin{proof}
	We continue with the same notation as in the proof of \cref{T: basic trans}.
	In particular, we again assume a closed embedding $e \colon V \to N \times \R^k$ that satisfies the hypotheses of \cref{L: all transversality is wrt embeddings}, a proper embedding of $N$ in $\R^K$ with an $\epsilon$ neighborhood $N_\epsilon$, a submersion $\pi \colon N_\epsilon \to N$, and, in the proper case, a function $\varepsilon \colon M \to (0,\infty)$ to meet the requirements of \cref{L: nearby proper homotopy}.

	In this case $f$ is already assumed to be smooth, but we alter our definition of the map $H \colon M \times D \to \R^K \to N$ as follows.
	Let $\rho \colon M \to [0,1]$ be a smooth function that is $0$ on $\bd M$ and $>0$ on $M-\bd M$ \cite[Proposition 5.43]{Lee13}.
	Then we define
	$$H(x,s) = \pi(f(x)+ \rho(x)\eta(x)s),$$
	where $\eta$ is defined as above.
	Then when $x \in \bd M$, we have $H(x,s) = \pi(f(x)) = f(x)$, so $H$ does not depend on $s$ when $x \in \bd M$.
	The remainder of the argument is now identical to the rest of the proof of \cref{T: basic trans}.
	In particular, for almost all $s$, the restriction of $H(-,s)$ to $M - \bd M$ will be transverse to $g$, while the restriction of $H(-,s)$ to $\bd M$ is just $f(x)$ and so transverse to $g$ by assumption.
	We then define $h$ from $H$ exactly as in the proof of \cref{T: basic trans}, and the conditions in the proper case are satisfied as $\rho(x)\leq 1$.
\end{proof}

\subsection{Pullbacks and fiber products}

When two embedded submanifolds of a manifold meet transversely, their intersection is again a submanifold.
More generally, if two smooth maps of manifolds are transverse, we can form their pullback, also called their fiber product, which is again a manifold.
This construction extends to manifolds with corners mapping into a manifold without boundary.

\begin{definition}\label{D: top pullback}
	Let $f \colon V \to M$ and $g \colon W \to M$ be transverse smooth maps from manifolds with corners to a manifold without boundary.
	Define the \textbf{pullback}\index{pullback|textbf} or \textbf{fiber product}\index{fiber product|textbf} $V \times_M W$ by
	$$V \times_M W = \{(x, y) \in V \times W \mid f(x) = g(y)\}.$$

	There are canonical maps from $V \times_M W$ to $V$, $W$, and $M$ that respectively take $(x,y)$ to $x$, $y$, and $f(x) = g(y)$.
	\[
	\begin{tikzcd}
		V \times_M W \arrow[r, "g^*"] \arrow[d, "f^*"] & V \arrow[d, "f"] \\
		W \arrow[r, "g"] & M.
	\end{tikzcd}
	\]
	We typically suppress the maps from the notation, though we sometimes label them as in the diagram and sometimes write $f \times_M g \colon V \times_M W \to M$ for the composition around the square.
	We also sometimes write $f^*$ as $\pi_W$ and $g^*$ as $\pi_V$, as these maps are induced by restricting to $V \times_MW$ the projections from $V \times W$ to $V$ and $W$.
	We also often abbreviate $V \times_M W$ as $P$.

	We will generally use the term \textit{pullback} when we want to emphasize $V \times_M W$ with its map to $V$ or $W$, while the \textit{fiber product} is to be considered as mapping to $M$.
	When treating $V \times_M W$ as a pullback, we also sometimes use the notation $g^*V \to W$ or $f^*W \to V$; this notation is consistent with the analogous notation for pullbacks of fiber bundles, which is a special case.

	Note that if one of $V$ or $W$ is empty, then $V \times_M W = \emptyset$.
\end{definition}

The reason for the transversality requirement in \cref{D: top pullback} is that this ensures that $V \times_M W$ is itself a manifold with corners, as follows from \cite[Theorem 7.2.7]{MaDo92} or \cite[Theorem 6.4]{Joy12}. In fact, we have the following.

\begin{theorem}\label{pullback}\index{pullback!manifold with corners}\index{pullback!stratification}
	Let $f \colon V \to M$ and $g \colon W \to M$ be transverse smooth maps from manifolds with corners to a manifold without boundary.
	Then $V \times_M W$ is a manifold with corners of dimension $v+w-m$ with
	\begin{equation*}
		S^i(V \times_M W) = \bigsqcup_{k + \ell = i} S^k(V) \times_M S^\ell(W).
	\end{equation*}
	Moreover, the maps from the fiber product to $V$, $W$, and $M$ are smooth in the strong sense of \cite[Definition 3.1]{Joy12} (see \cref{R: weakly smooth}).
\end{theorem}

\begin{proof}
	By \cite[Theorem 6.4]{Joy12} or \cite[Theorem 7.2.7]{MaDo92}, the pullback is a manifold with corners; see \cite[Remark 6.3]{Joy12} for the equivalence of the definitions of transversality in the two papers given our assumptions.
	Beyond Joyce's definition of manifolds with corners, we have also assumed that our manifolds with corners are second countable.
	However, products and subspaces of second countable spaces are also second countable \cite[Theorem 30.2]{Mu00}, so if $V$ and $W$ are second countable then so is $V \times_M W$.

	As we have realized $V \times_M W$ as a subset of $V \times W$, we meet our condition of \cref{D: MWC} that our manifolds with corners be contained in some Euclidean space, as by hypothesis this is the case for $V$ and $W$.

	The statements that the projections are smooth also follows from \cite[Theorem 6.4]{Joy12}, which states that $V\times_MW$ is a fiber product in the category of manifolds with corners (with smooth maps), and this includes the statement that the projection maps are maps in the category and hence smooth; see the introduction to \cite[Section 6]{Joy12}.

	For the decomposition, we refer to \cite[Theorem 7.2.7.a]{MaDo92}, which says that $V\times_MW$ is a \textit{totally neat} submanifold of $V\times W$.
	By \cite[Definition 3.1.10]{MaDo92}, the definition of ``totally neat submanifold'' includes the condition that the \textit{index} of a point in the submanifold is equal to its index in the ambient manifold.
	But the notion of index in \cite{MaDo92} is precisely our notion of depth. So $V\times_MW$ being a totally neat submanifold implies that $S^i(V\times_MW)$ must consist of those points of $V\times_MW$ in $S^i(V\times W)$.
	But clearly from the definitions,
	$$S^i(V\times W)=\bigsqcup_{k + \ell = i} S^k(V) \times S^\ell(W).$$
	So a point of $S^i(V\times_M W)$ consists of a pair $(x,y)\in V\times W$ with $f(x)=g(y)$ and such that $x\in S^k(V)$ and $y\in S^\ell(W)$ for some $k,\ell$ with $k+\ell=i$.
	But this means precisely that $(x,y)\in S^k(V) \times_M S^\ell(W)$ for some $k,\ell$ with $k+\ell=i$.
	Conversely, if $(x,y)\in S^k(V) \times_M S^\ell(W)$ and $k+\ell=i$, then $(x,y)\in V\times_MW$ and, since $(x,y)$ has index $i$ in $V\times W$, it must in fact be in $S^i(V\times_M W)$.
\end{proof}

To generalize this theorem when $M$ is also a manifold with corners requires substantial additional hypotheses in the definition of transverse smooth maps.
Such a generalization is a central result in \cite{Joy12}.

There is a Leibniz rule for taking boundaries of fiber products of transverse maps \cite[Proposition~6.7]{Joy12}:

\begin{proposition}\label{P: product boundary}\index{pullback!boundary formula}
	If $f \colon V \to M$ and $g \colon W \to M$ are transverse maps from manifolds with corners to a manifold without boundary, then
	\[
	\bd(V \times_M W) = (\bd V \times_M W) \sqcup (V \times_M \bd W),
	\]
	recalling that if $g \colon W \to M$ then we interpret $\bd W$ as equipped with the map $gi_{\bd W} \colon \bd W \to M$ and similarly for $V$.
\end{proposition}

We will see versions of this formula below that take into account orientations and co-orientations.

\subsubsection{Some further properties of transversality and fiber products}

In this section we mostly collect some well-known, though not always easy to find in the literature, results about transversality and fiber products.
We state these results mainly in the classical setting of manifolds without boundary, though they generally extend to the case of transverse maps of manifolds with corners mapping to a manifold without boundary, either by applying them to the pairwise transverse strata or by thinking in terms of the ``naive transversality'' of \cref{D: naive transversality}.

\begin{lemma}[Tangent spaces of pullbacks]\label{L: tangent of pullbacks}\index{pullback!tangent space}
	Let $f \colon V \to M$ and $g \colon W \to M$ be transverse smooth maps of manifolds with corners to a manifold without boundary.
	Suppose $x \in V$ and $y \in W$ with $f(x) = g(y) = z$ so that $(x,y) \in V \times_MW$.
	Then the tangent space of $V \times_MW$ at $(x,y)$ as a subspace of $T_{(x,y)}(V \times W) = T_xV \oplus T_yW$ consists of those vectors $(\mathbf v,\mathbf w)$ such that $Df(\mathbf v) = Dg(\mathbf w) \in T_zM$.
	In other words,
	$$T_{(x,y)}(V \times_M W) = T_xV \times_{T_zM} T_yW,$$
	or ``the tangent space of the fiber product is the fiber product of the tangent spaces.''
\end{lemma}

\begin{proof}
For manifolds without boundary, a proof can be found in \cite[Theorem 5.47]{Wed16}.
Wedhorn proves this theorem for ``premanifolds,'' which are essentially manifolds minus the Hausdorff and second countability conditions; these conditions will be automatic in our setting, so Wedhorn's proof applies.
For manifolds with corners, see \cite[Proposition 7.2.7.c]{MaDo92}.
\end{proof}

In many circumstances, this lemma allows us to reduce arguments about fiber products of maps of manifolds to arguments about fiber products of linear maps.

\begin{corollary}\label{C: tangent map of pullbacks}
With the notation of \cref{L: tangent of pullbacks}, $D_{(x,y)}(f \times_M g)=D_xf \times_{T_zM} D_yg$ as maps from $T_{(x,y)}(V \times_M W) = T_xV \times_{T_zM} T_yW$ to $T_z(M)$.
\end{corollary}

\begin{proof}
	By definition, the map $f \times_M g$ is either of the compositions $f \pi_V = g \pi_W$. So $D_{(x,y)} (f \times_M g) = D_yg \circ D_{(x,y)} \pi_W$. On the other hand, considering the fiber product $D_xf \times_{T_z M} D_y g$ of the maps $D_x f \colon T_xV \to T_zM$ and $D_y g \colon T_yW \to T_zM$, it is the composition of the projection $T_xV \times_{T_zM} T_yW \to T_yW$ and $D_y g$.
	But this projection is precisely $D_{(x,y)} \pi_W$, so the maps agree.
\end{proof}

\begin{lemma}\label{L: normal pullback}\index{pullback!normal bundle}
	Let $W$ be a manifold with boundary and $V$ and $M$ manifolds without boundary. Suppose $f \colon V \to M$ and $g \colon W \to M$ are transverse smooth maps and that $f$ is an embedding so that $V$ is a submanifold of $M$ with normal bundle $\nu V$.
	Then the pullback map $g^*V = V \times_MW \to W$ is an embedding onto the neat submanifold $g^*V \cong g^{-1}(V)$ and the normal bundle of $g^*V$ in $W$ is isomorphic to the pullback of $\nu V$.
	In other words,
	$$\nu(g^*V) \cong \left(g|_{g^*V}\right)^*(\nu V),$$
	or ``the pullback of the normal bundle is the normal bundle of the pullback.''
\end{lemma}

Here ``neat'' means that $\bd (g^*V) \subset \bd W$ and at each such point there is a chart $(U,\phi)$ of $W$ such $\phi$ is a diffeomorphism between $(U, U\cap g^*V)$ and $(\R^w_1, \R^{v+w-m}_1)$; see \cite[Definition II.2.2]{Kos93}.

\begin{proof}
	See \cite[Proposition IV.1.4]{Kos93}.
\end{proof}

As we like to think geometrically where possibly, we will generally perform the common abuse of thinking of normal bundles as subbundles of the tangent bundle of the ambient space over the embedded submanifold via a (usually tacit) choice of splitting of the quotient map to the true normal bundle.

\begin{lemma}\label{L: fiber product of embeddings}\index{fiber product!of submanifolds}
	Let $f \colon V \to M$ and $g \colon W \to M$ be transverse smooth maps of manifolds without boundary, and suppose $f$ and $g$ are both embeddings.
	Then, identifying $V$ and $W$ as submanifolds of $M$, the fiber product $V \times_MW$ is simply the intersection $V \cap W$, which is a smooth submanifold of $M$.
	In this case, we can identify the restriction of $\nu V$ to $V \cap W$ with a sub-bundle of the restriction of $TW$ to $V \cap W$.
\end{lemma}

\begin{proof}
	By \cref{L: normal pullback}, the pullback is a smooth submanifold of $W$, but $g$ is itself an embedding, so, identifying the manifolds with their embedded images, we obtain smooth submanifolds $V \times_MW \subset W \subset M$.
	That $V \times_MW = V \cap W$ in this situation follows immediately from the definition of the fiber product.

	The second claim is now a special case of \cref{L: normal pullback}, using the splitting of surjective vector bundle maps \cite[Theorem 3.9.6]{Hus75} to identify the normal bundle of $V \cap W$ in $W$ with a subbundle of $TW$.
\end{proof}

The preceding statements are not so evident in the case of manifolds with corners, as even how to extend the notions of immersions and embeddings when the target space is a manifold with corners is not so obvious. Joyce enforces his stronger version of smoothness in \cite{Joy12}, while the definition of immersion in \cite[Section 3.2]{MaDo92} has conditions beyond the classical injectivity of the tangent map. However, we have the following useful facts.

\begin{lemma}\label{L: immersions}\index{immersion}
	Let $f \colon V \to M$ be a map from a manifold with corners to a manifold without boundary.
	Then:
	\begin{enumerate}
		\item $f$ is an immersion at $x \in V$ if and only if $D_xf$ is injective.
		\item If $f$ is an immersion at all points and a homeomorphism onto its image, then it is an embedding.
		\item If $f$ is an immersion at $x$, then it restricts to an embedding on an open neighborhood of $x$.
	\end{enumerate}
\end{lemma}

\begin{proof}
	See \cite[Theorem 3.2.6, Proposition 3.2.13, and Definition 3.3.1]{MaDo92}.
	Note that when working with finite dimensional manifolds, as we always do, the ``topological supplement'' condition that appears often in \cite{MaDo92} is automatic \cite[page 7]{MaDo92}.
\end{proof}

The third statement of \cref{L: immersions} holds even when $M$ is a manifold with corners but requires the more sophisticated definition of immersion in \cite{MaDo92}.

\begin{lemma}\label{L: immersion pullback}\index{pullback!of immersion}
	Let $f \colon V \to M$ be an immersion of a manifold with corners into a manifold without boundary, let $W$ be a manifold without boundary, and let $g \colon W \to M$ be transverse to $V$.
	Then the pullback $\pi_W: V \times_M W \to W$ is an immersion and so, locally in a neighborhood of each point, an embedding.
\end{lemma}

\begin{proof}
	By \cref{L: tangent of pullbacks}, the tangent space of the pullback is the pullback of the tangents space, so if we write $P = V \times_M W$ and suppose $(v,w) \in P$, then $T_{(v,w)}P = T_vV \times_{T_{f(v)}M} T_wW$.
	Suppose $(a,b) \in T_{(v,w)}P$ with $a \in T_vV$ and $b \in T_wW$.
	Let $\pi_V \colon P \to V$ and $\pi_W \colon P \to W$ be the restrictions to $P$ of the projections from $V \times W$, and recall that $\pi_W = f^*$, the pullback of $f$.
	Then, suppressing the basepoints for the remainder of the argument, $D(f \times_M g) = D(g\pi_W) = D(f \pi_V)$.
	Suppose $D\pi_W(a,b)= b = 0$.
	Then $$0 = Dg(D\pi_W(a,b))= D(g\pi_W)(a,b) = D(f \pi_V)(a,b) = Df(a).$$
	But $f$ is an immersion, so $Df_v$ is injective by \cref{L: immersions}, and hence $a=0$.
	So $D\pi_W$ is injective, and therefore $\pi_W$ is an immersion, and so a local embedding, again by \cref{L: immersions}.
\end{proof}

\begin{lemma}\index{fiber product!of immersions}
	The fiber product of two transverse immersions $f \colon V \to M$ and $g \colon W \to M$ from manifolds with corners to a manifold without boundary is an immersion and so, locally in a neighborhood of each point, an embedding onto its image.
\end{lemma}

\begin{proof}
	Even though $W$ may have corners, the argument in the proof of \cref{L: immersion pullback} still holds to show that $D \pi_W$ is an injection at each point of $V \times_M W$.
	Then, as $g$ is an immersion, the composition $D(g \pi_W)$ is also injective at each point.
	So the lemma follows by \cref{L: immersions}.
\end{proof}

The next lemma, concerning transversality with respect to pullbacks, will be useful later in proving certain functoriality formulas.
By \cref{pullback} and the fact that transversality of manifolds with corners reduces to transversality of manifolds without corners (looking at all the pairs of strata), it suffices to state the lemma just for manifolds without boundary.

\begin{lemma}\label{L: transverse to pullback}\index{transversality!with respect to pullbacks}\index{pullbacks!transversality of}
	Let $f \colon V \to M$, $g \colon W \to M$, and $h \colon X \to W$ be smooth maps of manifolds without boundary, and suppose $f$ is transverse to $g$.
	Then $gh \colon X \to M$ is transverse to $f$ if and only if $h$ is transverse to the pullback $\pi_W \colon V \times_M W \to W$.
\end{lemma}

\begin{proof}
	For simplicity of notation, let $P = V \times_MW$, and let $\pi_V \colon P \to V$ and $\pi_W \colon P \to W$ be the maps induced by the projections from $V \times W$ to $V$ and $W$.
	So we have the diagram
	\[
	\begin{tikzcd}
		& P = V \times_M W \arrow[r, "\pi_V"] \arrow[d, "\pi_W"] & V \arrow[d, "f"] \\
		X \arrow[r, "h"] & W \arrow[r, "g"]& M.
	\end{tikzcd}
	\]
	First, suppose $h$ is transverse to the pullback $P \to W$.
	Note that the existence of $P = V \times_M W$ uses the assumption that $f$ and $g$ are transverse.
	Suppose $x \in X$ and $y \in V$ with $gh(x) = f(y)$.
	Then $(y,h(x)) \in P$, and by assumption we have both $Df(T_yV)+Dg(T_{h(x)}W) = T_{f(y)}M$ and
	$Dh(T_xX)+D\pi_W(T_{(y,h(x))}P) = T_{h(x)}W$.
	Applying $Dg$ to the second formula, we have
	$$Dg(T_{h(x)}W) = D(gh)(T_xX)+D(g\pi_W)(T_{(y,h(x))}P).$$
	So from the first formula,
	$$T_{f(y)}M = Df(T_yV)+D(gh)(T_xX)+D(g\pi_W)(T_{(y,h(x))}P).$$
	But $g\pi_W = f\pi_V$, so $D(g\pi_W)(T_{(y,h(x))}P) = D(f\pi_V)(T_{(y,h(x))}P) \subset Df(T_yV)$.
	It follows that
	$$T_{f(y)}M = Df(T_yV)+D(gh)(T_xX),$$
	i.e.\ $f$ is transverse to $gh$.

	Next suppose $f$ is transverse to $gh$ and that $x \in X, p \in P$ with $h(x) = \pi_W(p)$.
	As $gh(x) = g\pi_W(p) = f\pi_V(p)$, we have $gh(x)$ in the image of $f$, so in particular $gh(x)$ is in the intersection of the images of $V$ and $W$ in $M$.
	Now suppose $\mathbf w \in T_{h(x)}W$.
	By assumption $T_{hg(x)}M = Df(T_{\pi_V(p)}V)+D(gh)(T_xX)$, so we can write $Dg(\mathbf w) = Df(\mathbf a)+Dgh(\mathbf b)$ for some $\mathbf a \in T_{\pi_V(p)}V$ and $\mathbf b \in T_x X$.
	Now consider $\mathbf z = \mathbf w-Dh(\mathbf b) \in T_{h(x)}W$.
	Applying $Dg$ to both sides, we have $Dg(\mathbf z) = Dg(\mathbf w)-Dgh(\mathbf b) = Df(\mathbf a)$.
	So $(\mathbf a,\mathbf z)$ is in the pullback $T_{\pi_V(p)}V \times_{T_{gh(x)}M} T_{h(x)}W$.
	But by \cref{L: tangent of pullbacks}, this is precisely the tangent space of $V \times_MW$ at $(h(x),\pi_V(p))$.
	Furthermore, we have $D\pi_W(\mathbf a,\mathbf z) = \mathbf z$.
	Thus $\mathbf w = \mathbf z+Dh(\mathbf b) = D\pi_W(\mathbf a,\mathbf z)+Dh(\mathbf b)$.
	As $\mathbf w$ was arbitrary, $W = D\pi_W(T_{(h(x),\pi_V(p))}P)+Dh(T_xX)$, as desired.
\end{proof}

\subsubsection{Fiber products with more than two inputs}\index{fiber product!more than two inputs}

Finally for this section, we briefly consider transversality of more than two maps. \cref{R: multiproducts} highlights some of the issues and difficulties involved, while \cref{P: 3 out of 4 trans} shows that there are still some simplifications that can be observed.

\begin{remark}\label{R: multiproducts}
	In this long remark we briefly discuss transversality and fiber products that involve more than two maps.
	This is relevant, for example, when considering associativity of fiber products or pullbacks of fiber products.
	As transversality of maps of manifolds is defined in terms of the behavior of the maps of tangent spaces, it is useful to first recall some notions about transversality in the setting of linear maps of vector spaces.

	If $f \colon V \to M$ and $g \colon W \to M$ are linear maps of vector spaces, then transversality of $f$ and $g$ can be expressed in a number of equivalent ways \cite[Section 4.7]{RamBas09}:
	\begin{itemize}
		\item $f(V)$ and $g(W)$ span $M$,
		\item the map $\Delta \colon V \times W \to M$ given by $\Delta(v,w) = f(v)-g(w)$ is surjective,
		\item $\dim(V \times_MW) = v+w-m$.
	\end{itemize}
	The last two formulations easily generalize to $n$-tuples of maps $f_i \colon V_i \to M$.
	Such an $n$-tuple is considered transverse (as an $n$-tuple) when either of the following equivalent conditions hold:
	\begin{itemize}
		\item the map $\Delta \colon \prod V_i \to M^{n-1}$ given by $\Delta(v_1,\ldots,v_n) = (f_2(v_2)-f_1(v_1),\ldots, f_n(v_n)-f_{n-1}(v_{n-1}))$ is surjective,
		\item the fiber product given by $\{(v_1,\ldots,v_n) \in V_1 \times \cdots \times V_n \mid f_1(v_1) = \cdots = f_n(v_n)\}$ has dimension $\sum_{i = 1}^n\dim(V_i) -n\dim(M)$.
	\end{itemize}
	This version of transversality behaves very well in that this ``$n$-transversality'' is equivalent to the iterated transversality conditions that are required when taking fiber products two at a time.
	In other words, an $n$-tuple is transverse if and only if for any $1 \leq i<j \leq n$,
	\begin{itemize}
		\item the $j-i+1$-tuple $\{f_i,\ldots,f_j\}$ is transverse, and
		\item letting $P$ denote the fiber product of the maps $\{f_i,\ldots,f_j\}$, the $n-(j-i)$-tuple consisting of the fiber product map $P \to M$ and the $f_k$ with $k\notin\{i,\ldots,j\}$ is transverse.
	\end{itemize}
	Iterating this fact, we can see that this is equivalent to having $V_1$ transverse to $V_2$, then $V_3$ transverse to $V_1 \times_MV_2$ and so on.
	In particular, we can form the $n$-tuple fiber product if and only if we can form the \textit{transverse} iterated fiber products such as $(((V_1 \times_MV_2) \times_M V_3) \times_M\ldots) \times_M V_n$ as well as in any other order of association.
	See\footnote{Proposition 8-1 of \cite{RamBas09} actually concerns a more general situation of ``mixed associativity'' in which the data consists of zig-zags
	$$V_1 \xr{f_1}S_1\xl{g_2}V_2 \xr{f_2}\cdots \xr{f_{n-1}}S_{n-1}\xl{g_n}V_n.$$
	However, this reduces to our setting by taking $S_i = M$ and $g_i = f_i$ for all $i$.}
	\cite[Propositions~4-9 and 8-1]{RamBas09}.

	Unfortunately, in the setting of maps of manifolds, the situation is less well behaved.
	Suppose $f_i \colon V_i \to M$ is an $n$-tuple of smooth maps of manifolds.
	There are notions of transversality for such an $n$-tuple that imply the fiber product $P = \{(v_1,\ldots, v_n) \in V_1 \times \cdots \times V_n \mid f_1(v_1) = \cdots = f_n(v_n)\}$ is a smooth manifold of dimension $\sum\dim(V_i)-(n-1)\dim(M)$; e.g.\ \cite[Theorem 7.2.7]{MaDo92}.
	However, in this setting the equivalence between such $n$-ary transversality and ``iterated transversality'' cannot hold in general, because the $n$-ary fiber product does not involve potentional problem points that would arise in intermediate stages of an iterated fiber product.

	For example, consider maps $f \colon V \to M$, $g \colon W \to M$, and $h \colon Z \to M$.
	Suppose that $f$ and $g$ are not transverse; in particular, suppose that the fiber product $V \times_MW$ is not a manifold.
	Further, suppose that $h(Z)$ is disjoint from $f(V)$ and $g(W)$.
	In this case, the triple of maps is transverse (vacuously) and its 3-ary fiber product is well defined as $\emptyset$.
	However, the iterated fiber product $(V \times_MW) \times_M Z$ is not well defined in the category of smooth manifolds and maps.
	Given that we define fiber products only of transverse maps of smooth manifolds, in this case $V \times_M W$ is not properly defined in our category, and it is further impossible then to define $(V \times_MW) \times_M Z$ in this context.

	The upshot of all this is that when considering situations involving fiber products of more than two maps, we shall have to be careful about the transversality assumptions.
\end{remark}

Despite the preceding remark, the transversality conditions involved in associativity of fiber products are not completely independent, as the following proposition shows. In later sections, we will usually simply assume all needed transversality exists, but the following proposition can be useful in practice.

\begin{proposition}\label{P: 3 out of 4 trans}\index{fiber product!more than two inputs!redundancy of transversality conditions}
	Let $f \colon V \to M$, $g \colon W \to M$, and $h \colon Z \to M$ be maps from manifolds with corners to a manifold without boundary. Suppose that $W$ is transverse to $Z$ and that $V$ is transverse to $W$ and to $W \times_MZ$. Then $V \times_MW$ is transverse to $Z$. In particular, if $V \times_M(W \times_MZ)$ and $V \times_MW$ are well defined, then so is $(V \times_MW) \times_M Z$.
\end{proposition}

\begin{proof}
	As above, it suffices to assume that $V$, $W$, and $M$ are manifolds without boundary. We must show that $V \times_MW$ is transverse to $Z$, so we consider points $(v,w) \in V \times_M W$ and $z \in Z$ such that $h(z)$ is equal to $(f \times_Mg)(v,w)$, which by definition is equal to $f(v) = g(w)$. In other words, we consider $(v,w,z) \in V \times W \times Z$ such that $f(v) = g(w) = h(z)$.

	So suppose $(v,w,z)$ is such a triple, and denote the common image by $m \in M$. By the transversality assumptions, we know that the images of $D_wg \colon T_wW \to T_mM$ and $D_z h \colon T_zZ \to T_mM$ span $T_mM$, i.e.\ that $D_wg$ and $D_zh$ are transverse as linear maps, and similarly that $D_{(w,z)}(g \times_M h) \colon T_{(w,z)}(W \times_M Z) \to T_mM$ is transverse to $D_vf \colon T_vV \to T_mM$. Furthermore, by \cref{L: tangent of pullbacks,C: tangent map of pullbacks}, the tangent space of a fiber product is the fiber product of the tangent spaces, so $T_{(w,z)}(W \times_M Z) = T_wW \times_{T_mM}T_zZ$ and $D_{(w,z)}(g \times_M h) = D_wg \times_{T_mM}D_zh$.

	Now by \cite[Propositions~4-9]{RamBas09}, the triple of linear maps $(D_vf,D_wg,D_zh)$ is transverse as a triple of maps if and only if both $D_wg$ is transverse to $D_zh$ and $D_wg \times_{T_mM}D_zh$ is transverse to $D_vf$. As such statements are independent of how we order the terms, the transversality established in the preceding paragraph also implies that $D_vf$ and $D_wg$ are transverse (which already follows from the hypotheses of the proposition), and $D_vf \times_{T_mM}D_wg$ is transverse to $D_zh$. But this implies, again using \cref{L: tangent of pullbacks,C: tangent map of pullbacks}, that $h$ is transverse to $f \times_Mg$, as desired.
\end{proof}

\begin{remark}
	The end of the preceding proof at first seems to imply that if $g$ and $h$ are transverse and $g \times_M h$ is transverse to $f$, then $f$ is transverse to $g$ and $f \times_Mg$ is transverse to $h$. Indeed, \cite[Propositions~4-9]{RamBas09} says this is the case for the linear maps of the tangent spaces. Unfortunately, however, as \cite[Propositions~4-9]{RamBas09} applies only to linear maps, we can apply it only at those points $(v,w,z) \in V \times W \times Z$ where we know that $f(v) = g(w) = h(z)$ so that all three tangent space maps are well defined. So such a result would hold if the only intersections among the maps were such triple intersections. However, as noted in \cref{R: multiproducts}, there could be pairs $(v,w) \in V \times W$ with $f(v) = g(w)$, but with this common image in $M$ not in the image of $h$. At such points, \cite[Propositions~4-9]{RamBas09} cannot tell us anything about the transversality of $V$ and $W$, and so $V \times_MW$ might not be well defined as a manifold with corners due to failure of transversality, even if $V \times_M (W \times_M Z)$ is.
\end{remark} 

\section{Orientations and co-orientations}\label{S: orientations and co-orientations}

Manifolds with corners are, in particular, topological manifolds (with boundary), and so they carry the standard notions of orientability and orientation.
As in singular or simplicial homology, orientations carry sign information in geometric versions of homology theory.
For geometric cohomology, however, it turns out that the natural structures to carry sign information are co-orientations, sometimes called orientations of maps.
Unlike orientations, co-orientation can be ``pulled back.''

Co-orientations are less familiar than orientations, so it is helpful to keep the following central example in mind: If $W \to M$ is an immersion of manifolds, a co-orientation is equivalent to an orientation of the normal bundle of the image; see \cref{normal co-or} below.
Notice that this condition does not require the orientability of either $W$ or $M$.
In fact, an important case is when neither $W$ nor $M$ is orientable but the monodromies of their orientation bundles around loops in $W$ are either both orientation-preserving or both orientation-reversing; it is in this sense that we have a \textit{co}-orientation.

The special case of (local) immersions is important both for intuition and in practice; for example, the only geometric cochains which evaluate nontrivially on fixed collections of chains through the intersection homomorphisms of \cite{FMS-flows} are local immersions.
However, for complete definitions that cover all cases of interest, it is critical to co-orient more general maps and to do so in a way that provides key properties when forming pullbacks, such as a Leibniz rule when taking boundaries, graded commutativity of cochains, and a simple expression when the maps are immersions.
While co-orientations can be found in many places in the literature, we could not find a careful treatment that guaranteed these key properties.
Therefore, we develop co-orientations in depth in this (long) section.
Some readers may prefer to skip ahead, either considering geometric cohomology only with mod~2 coefficients in a first reading, or perhaps take orientation of the normal bundle as a temporary partial definition, coming back later to understand the general setting.

\subsection{Orientations}\label{S: orientations}\index{orientation|textbf}

If $W$ is a manifold with corners then it is a topological manifold with boundary and so in the interior of $W$, i.e.\ on $S^0(W)$, we can consider orientability and orientations of $W$ in the usual sense, from either the topological or smooth manifold points of view, which are equivalent \cite[Theorem VI.7.15]{Bred97}.
Following standard conventions, we typically refer to an oriented manifold with corners $W$ with the orientation tacit.
If an orientation on $W$ is understood, then $-W$ refers to $W$ with the opposite orientation.

When $W$ is orientable, so is its (topological) boundary \cite[Lemma 6.9.1]{Bred97}, and since $S^0(W) \cup S^1(W)$ is a smooth manifold with boundary, we can allow an orientation of $W$ to determine an orientation of $S^1(W)$ using standard smooth manifold conventions.
This gives an orientation for $\bd W$, as we can identify $S^1(W)$ with the interior of $\bd W$.
In particular, we choose the following convention, which agrees with that of Joyce \cite[Convention 7.2.a]{Joy12}:

\begin{convention}\label{Con: oriented boundary}\index{orientation!boundary}
	For a smooth oriented manifold with boundary $N$, we orient $\bd N$ by stipulating that an outward normal vector followed by an oriented basis of $\bd N$ yields an oriented basis for $N$.

	When $W$ is an oriented manifold with corners, we identify $S^1(W)$ with the interior of $\bd W$, and this convention determines an orientation of $\bd W$.
\end{convention}

When we wish to work with orientations symbolically, the following interpretation will be extremely useful.

\begin{definition}\label{D: det bundle}\index{determinant line bundle}
	Let $E \to B$ be a rank $d$ real vector bundle.
	Define the \textbf{determinant line bundle}
	$\Or(E)$ to be $\bigwedge^d E$.
	If $d = 0$ this is interpreted to be the trivial rank one ``bundle of coefficients.''
	We call the principal $O(1) \cong C_2$ bundle associated to $\Or(E)$ the \textbf{orientation cover}\index{orientation cover} of $E$.

	An \textbf{orientation}\index{determinant line bundle!orientation of} of $E$, if it exists, is then a section of its orientation cover or, equivalently, an equivalence class of non-zero sections of $\Or(E)$ such that two sections are equivalent if they differ by multiplication by an everywhere positive scalar function.
\end{definition}

In particular, we thus think of orientations of a manifold $M^m$ as (equivalence classes of) non-zero sections of $\Or(TM)$.
We typically use the notation $\beta_M$ to stand for such a section, and, as we are only ever interested in such sections to represent orientations, we systematically abuse notation by not distinguishing between a section and its equivalence class.
Thus an expression such as $\beta_M = \beta_V \wedge \beta_W$ should be interpreted as an equality of equivalence classes.
Such formulas are also often meant to be interpreted locally over some particular point or subspace that will be understood from context, and in this context we refer to expressions such as $\beta_V$ as \textbf{local orientations}\index{orientation!local}.
This notation turn out to be extremely useful in working with orientations, and we will use it frequently.
Whenever we form such wedge products, if one of the terms is an element of $\bigwedge^0 TV$ for some $V$ we treat that term as a scalar function and interpret $\wedge$ as the fiberwise scalar product.

Of course, over a point $x \in M$, we identify the oriented basis $(e_1,\ldots, e_m)$ of $T_xM$ with the section $e_1 \wedge \cdots \wedge e_m$ of $\Or(TM)$.
We also use the standard identification $\bigwedge (A\oplus B)=(\bigwedge A )\otimes (\bigwedge B)$, letting $\bigwedge$ always denote the top-dimensional exterior product.
So, for example, if $V$ and $W$ are two submanifolds of a manifold $M$ intersecting transversely at a point $x$, we would have $T_xM = T_xV\oplus T_xW$, and the local orientation formula $\beta_M = \beta_V \wedge \beta_W$ would indicate that concatenating an oriented basis for $T_xV$ with an oriented basis for $T_xW$, in that order, gives an oriented basis for $T_xM$.
More generally, if $E$ is an oriented bundle over $M$, we can write $\beta_E$ for the corresponding local orientation at a point or in the neighborhood of a point.
Any abuses of notation involved in this calculus are very much justified by its extreme usefulness in working with orientations, as we shall see.

\begin{example}
	We have said that for a smooth oriented manifold with boundary $N$, we orient $\bd N$ so that an outward normal vector followed by an oriented basis of $\bd N$ yields an oriented basis for $N$.
	In our notation, we can let $\nu$ denote the normal bundle to the boundary and $\beta_{\nu}$ the section of the normal bundle tangent bundle corresponding to the outward-pointing orientation.
	Then letting $\beta_N$ and $\beta_{\bd N}$ denote the local orientations of $N$ and $\bd N$ at a boundary point, we can write our boundary orientation convention by saying that at each boundary point $\beta_N=\beta_\nu\wedge \beta_{\bd N}$.

	Of course, the normal bundle is technically a quotient bundle, but we can use the splitting of exact sequence of vector bundles \cite[Theorem 3.9.6]{Hus75} to identify it with a subbundle of the tangent bundle of $N$.
	Such choices are not unique, but they are unique up to elements of $T(\bd N)$, so the (equivalence class of the) exterior product $\beta_\nu\wedge \beta_{\bd N}$ does not depend on such choices.
	Such identifications will be used regularly and tacitly when working with normal bundles.
\end{example}

\subsubsection{Orientations of fiber products}\label{S: orientation of fiber products}\index{orientation!of fiber product|(}

If $V$ and $W$ are oriented manifolds, we orient $V \times W$ in the standard way by concatenating oriented bases of tangent spaces of $V$ with those of $W$.
More generally, if $f \colon V \to M$ and $g \colon W \to M$ are transverse maps with $V$, $W$, and $M$ all oriented and $M$ without boundary, Joyce defines an orientation on the pullback $V \times_M W$ as follows \cite[Convention 7.2b]{Joy12}.
Consider the short exact sequence of vector bundles over $P = V \times_M W$ given by
\begin{equation}\label{E: fiber orientation}
	0 \to TP \xr{D\pi_V \oplus D\pi_W} \pi_V^*(TV) \oplus \pi_W^*(TW) \xr{\pi_V^*(Df)-\pi_W^*(Dg)} (f\pi_V)^*TM \to 0.
\end{equation}
Here $\pi_V$ and $\pi_W$ are the projections of $V \times_M W$ to $V$ and $W$, and $D\pi_V$ is being treated as a map $TP \to \pi_V^*(TV)$ and similarly for $D\pi_W$.
Analogously, $\pi_V^*(Df)$ is the pullback of the map $Df \colon TV \to TM$ obtained first by treating it as a map $TV \to f^*(TM)$ and then pulling back functorially by $\pi_V^*$, and similarly for $\pi_W^*(Dg)$; note that $f\pi_V = g\pi_W$ at points of $P$.
By choosing a splitting, this sequence determines an isomorphism
\begin{equation*}
	TP \oplus (f\pi_V)^*TM\cong\pi_V^*(TV) \oplus \pi_W^*(TW).
\end{equation*}
The choices of orientations on $V$, $W$, and $M$ determine orientations on all summands in this expression except $TP$.
The orientation on $TP$ is then chosen so that the two direct sums differ in orientation by a factor of $(-1)^{wm}$, recalling that $w=\dim(W)$, etc.

Much more about the orientations of fiber products can be found in the technical report of Ramshaw and Basch \cite{RamBas09}.
While the focus there is on manifolds without boundary, and sometimes just fiber products of linear maps of vector spaces, the results about orientations extend to manifolds with corners by employing them on the top-dimensional stratum and utilizing their stability property, by which orientation properties of fiber products of linear maps extend to properties of fiber products of transverse manifolds (see \cite[Sections 6.3, 9.1.2, and 9.3]{RamBas09}).
Their orientation of fiber products agrees with Joyce's.
This can be checked directly from the definitions\footnote{Their multiplicative ``fudge factor'' in \cite[Theorem 9.14]{RamBas09} at first appears to be different from Joyce's, but this is only because their conventions utilize what in our notation would be the map $\pi_W^*(Dg)-\pi_V^*(Df)$ rather than $\pi_V^*(Df)-\pi_W^*(Dg)$.} or, as Joyce notes in \cite[Remark 7.6.iii]{Joy12}, axiomatically, as Ramshaw and Basch show that theirs is the unique choice of orientation convention satisfying certain basic expected properties.
It is these properties that determine the sign in the definition.
We state these properties in the following two propositions.

\begin{proposition}\label{P: oriented fiber product basic properties}
	Let $f \colon V \to M$ and $g \colon W \to M$ be transverse maps from oriented manifolds with corners to an oriented manifold without boundary.
	\begin{enumerate}
		\item\index{orientation!of fiber product!over a point} When $M$ is a point, the oriented fiber product $V \times_M W$ is simply $V \times W$, and in this case the fiber product orientation is consistent with the basic concatenation rule for products.
		\item{orientation!of fiber product!one map is the $\id$} When one of the maps is the identity $\id_M \colon M \to M$, the projection map to the other factor is an orientation preserving diffeomorphisms
		\begin{equation*}
			M \times_M V = V\quad\text{and}\quad V \times_M M = V.
		\end{equation*}
	\end{enumerate}
\end{proposition}

\begin{proposition}\label{P: oriented fiber mixed associativity}\index{orientation!of fiber product!mixed associativity}
	Let $V$, $W$, and $Z$ be oriented manifolds with corners, and let $M$ and $N$ be oriented manifolds without boundary.
	Then the ``mixed associativity'' formula for oriented fiber products
	\begin{equation}\label{E: mixed associativity fiber orientation}
		(V \times_M W) \times_N Z = V \times_M (W \times_N Z)
	\end{equation}
	holds when given maps
	$$V \xr{f} M\xleftarrow{g} W \xr{h} N \xl{k} Z$$
	and assuming sufficient transversality for all the fiber products in \eqref{E: mixed associativity fiber orientation} to be well defined (see \cref{R: multiproducts}).
	In this case, the map $V \times_M W \to N$ is given by composing the projection from $V \times_M W$ to $W$ with $h$, and similarly for the map $W \times_N Z \to M$.
\end{proposition}

These propositions are evident at the level of spaces.
When taking orientations into account, the first property in \cref{P: oriented fiber product basic properties} is proven in \cite[Sections 9.3.9]{RamBas09} as the ``concatenation axiom,'' and the second is proven in \cite[Sections 9.3.5 and 9.3.6]{RamBas09} as the ``left and right identity axioms.''
The mixed associativity property is proven in \cite[Sections 9.3.7]{RamBas09}.
An important special case of this associativity that we will need below occurs when $M = N$ and $g = h$, so that our initial data is three maps all to $M$.
In this case we have the ordinary associativity\index{orientation!of fiber product!associativity}
\begin{equation}\label{E: oriented fiber associativity}
	(V \times_M W) \times_M Z = V \times_M (W \times_M Z).
\end{equation}
That these properties determine the orientation rule for fiber products is the content of \cite[Theorem 9-10]{RamBas09}.
Technically, Ramshaw and Basch require for uniqueness two other properties: an Isomorphism Axiom, which says that the construction is consistent across oriented homeomorphisms, and a Stability Axiom, which implies that the orientation can be determined pointwise in a globally consistent manner.
These properties are both implicit in Joyce's global definition of the fiber product orientation.

There is also a commutativity rule proven in \cite[Sections 9.3.8]{RamBas09} that follows from the other properties:

\begin{proposition}\label{P: commute oriented fiber}
	Let $f \colon V \to M$ and $g \colon W \to M$ be transverse maps from oriented manifolds with corners to an oriented manifold without boundary.
	Then, as oriented manifolds,\index{orientation!of fiber product!commutativity}
	\begin{equation*}
		V \times_M W = (-1)^{(m-v)(m-w)}W \times_M V.
	\end{equation*}
\end{proposition}

Recalling that we write $\dim(M)=m$, etc., this means that the canonical diffeomorphism taking $(x,y) \in V \times_M W \subset V \times W$ to $(y,x) \in W \times_M V \subset W \times V$ takes a positively-oriented basis of the tangent space of $V \times_M W$ to a $(-1)^{(m-v)(m-w)}$-oriented basis of the tangent space of $W \times_M V$.
We note that these signs, while note quite in line with the Koszul conventions, agree with those for the intersection product of homology classes in Dold \cite[Section VIII.13]{Dol72}.

Furthermore, with our convention for oriented boundaries, one obtains the following useful identity; see \cite[Propositions 7.4 and 7.5]{Joy12}

\begin{proposition}\label{P: oriented fiber boundary}\index{orientation!of fiber product!boundary formula}
	Let $f \colon V \to M$ and $g \colon W \to M$ be transverse maps from oriented manifolds with corners to an oriented manifold without boundary.
	Then, as oriented manifolds,
	\begin{equation*}
		\bd (V \times_M W) = (\bd V \times_M W) \sqcup (-1)^{m-v}(V \times_M \bd W).
	\end{equation*}
\end{proposition}

\subsubsection{Fiber products of immersions}\index{orientation!of immersion}

The special case of fiber products with $f \colon V \to M$ and $g \colon W \to M$ embeddings or, a bit more generally, immersions, is of particular interest, especially for developing intuition.
Once again, since we are concerned primarily with orientations in this section, it is sufficient to consider $V$ and $W$ to be manifolds without boundary, and then the results we obtain extend directly to manifolds with corners.
In the case of embeddings, $V \times_M W$ is simply the intersection of $V$ and $W$ as submanifolds of $M$ by \cref{L: fiber product of embeddings}, and in the immersed case this is true locally, i.e.\ restricting attention to submanifolds of $V$ and $W$ on which $f$ and $g$ are embeddings.
Let us try to understand the orientation of $V \times_M W$ determined by orientations of $V$, $W$, and $M$ in this setting.

For convenience of notation, let us assume we have transverse embeddings so that $V$ and $W$ are submanifolds of $M$.
Then, by \cref{L: fiber product of embeddings}, we know that $P = V \times_M W$ is just the intersection $V \cap W$.
Consider a point $x \in P$.
As $V$, $W$, and $P$, are all submanifolds, $T_xV$, $T_xW$, and $T_xP$ are all subspaces of $T_xM$. Furthermore, as orientations are defined via the tangent bundles, it suffices to consider the relation among the orientations of $T_x V$, $T_xW$, $T_xM$, and $T_xP$, and, as the tangent space of a fiber product is the fiber product of the tangent spaces by \cref{L: tangent of pullbacks}, we have $T_xP = T_xV \times_{T_xM} T_xW = T_xV \cap T_xW$.
For simplicity of notation in the following, we drop the ``$T_x$'' from the notation and treat the maps $f$ and $g$ as linear embeddings of vector spaces.
We can then work with orientations of vector spaces expressed in our exterior power notation.

So we have a vector space $M$ with subspaces $V$, $W$, and $P = V \times_M W = V \cap W$.
By the transversality assumption, $V$ and $W$ span $M$.
Let $\nu W \subset V$ be a complementary subspace to $P$ in $V$ so that $V = P \oplus \nu W$; the notation is meant to suggest that $\nu W$ is a choice of normal subspace to $W$ in $M$.
Similarly, let $\nu V \subset W$ be a complementary subspace to $P$ in $W$ so that $W = P \oplus \nu V$.
Then we have $M = \nu W \oplus P \oplus \nu V$. \greg{I think this would be a good place for a picture. Maybe have $V$ and $W$ be 2-dimensional subspaces of $\R^3$? In fact, see \cref{Ex: intersection orientation}.}

In our current context, the exact sequence \eqref{E: fiber orientation} becomes a sequence of vector spaces
\begin{equation*}
	0 \to P \to V \oplus W \xr{f-g} M \to 0
\end{equation*}
with the first non-trivial map being the direct sum of inclusions and with $f$ and $g$ the emeddings of $V$ and $W$ into $M$.
Using our direct sum decompositions of the preceding paragraph, we can choose a splitting $M \to V \oplus W$ as
\[
\nu W \oplus P \oplus \nu V \to P \oplus \nu W \oplus P \oplus \nu V,
\]
given by $(x,p,y) \mapsto (0, x, -p, -y)$.
The signs are necessary due to the sign in $f-g$.
With this splitting, our resulting isomorphism $P \oplus M \to V \oplus W$ can be written in block matrix form as
\begin{equation}\label{E: orientation matrix}
	\begin{pmatrix*}[r]
		I&0&0&0\\
		0&I&0&0\\
		I&0&-I&0\\
		0&0&0&-I
	\end{pmatrix*},
\end{equation}
which has determinant $(-1)^{w}$.

By definition, the orientation of $P$ determined by the orientations of $V$, $W$, and $M$ is the orientation so that this matrix takes the concatenation of the orientation of $P$ with an orientation of $M$ (i.e.\ an ordered basis representing this concatenation orientation) to $(-1)^{wm}$ times the concatenation of the orientations of $V$ and $W$ (i.e.\ an ordered basis of $V \oplus W$ that differs from the concatenation orientation by a permutation of sign $(-1)^{wm}$).

This call can be expressed via our calculus of local orientations as follows.

\begin{proposition}\label{P: orient intersection}\index{orientation!of fiber product!of immersions}
	Let $V$ and $W$ be transverse oriented submanifolds of the oriented manifold $M$.
	Let $\beta_V$, $\beta_W$, and $\beta_M$ be local orientations of $V$, $W$, and $M$ at a point of $P = V \times_M W = V \cap W$.
	Then, using our notation established just above, the fiber product orientation of $P$ is the unique orientation $\beta_P$ such that if we choose orientations $\beta_{\nu W}$ and $\beta_{\nu V}$ for $\nu W$ and $\nu V$ such that $\beta_P \wedge \beta_{\nu W} = \beta_V$ and $\beta_P \wedge \beta_{\nu V} = \beta_W$ then $\beta_{\nu W} \wedge \beta_P \wedge \beta_{\nu V} = (-1)^{w(m+1)}\beta_M$ or, alternatively, $$\beta_P \wedge \beta_{\nu V} \wedge \beta_{\nu W} = \beta_M.$$
\end{proposition}

\begin{proof}
	We may work with $V$, $W$, and $P$ as linear subspaces of a vector space $M$ as above. Note that if we replace $\beta_P$ with its opposite orientation, then this must also reverse the orientations $\beta_{\nu W}$ and $\beta_{\nu V}$ and hence altogether we get the opposite orientation for $\beta_P \wedge \beta_{\nu V} \wedge \beta_{\nu W}$.
	Thus there is a unique choice of orientation $\beta_P$ as described in the lemma, and we must show that this is the orientation as defined in \cref{S: orientation of fiber products}.

	It will be more convenient to prove the lemma in the first form, but the second form follows by observing that $$(-1)^{w(m+1)} = (-1)^{wm+w} = (-1)^{wm-w^2} = (-1)^{w(m-w)},$$
	and then
	\[
	\beta_{\nu W} \wedge \beta_P \wedge \beta_{\nu V} =
	(-1)^{w(m-w)} \beta_P \wedge \beta_{\nu V} \wedge \beta_{\nu W}
	\]
	as $\dim(P \oplus \nu V) = \dim(W) = w$ and $\dim(\nu W) = m-w$.

	To prove the first statement, let $(p_1,\cdots,p_a)$ be an ordered basis for $P$ consistent with the orientation described in the lemma; so we can write $\beta_P = p_1 \wedge\cdots\wedge p_a$.
	When we consider each $p_i$ as a vector in $V$, $W$, or $M$, we write $p_i^V$, $p_i^W$, or $p_i^M$.
	We employ a similar convention with the other bases we will consider.
	Let $(x_1,\cdots,x_b)$ and $(y_1,\cdots,y_c)$ be corresponding ordered bases for $\nu W$ and $\nu V$ as described in the lemma, and we can write $\beta_{\nu W}$ and $\beta_{\nu V}$ analogously. Recall that, by assumption, we have $\beta_P \wedge \beta_{\nu W} = \beta_V$, $\beta_P \wedge \beta_{\nu V} = \beta_W$, and $\beta_M = (-1)^{w(m+1)} \beta_{\nu W} \wedge \beta_P \wedge \beta_{\nu V}$, and we must show that this is consistent with the definition of the fiber product orientation for $P$.

	So our orientation of $P \oplus M$ obtained by concatenation is
	$$(-1)^{w(m+1)} p_1 \wedge\cdots\wedge p_a \wedge x^M_1 \wedge\cdots\wedge x^M_b \wedge p^M_1 \wedge\cdots\wedge p^M_a \wedge y^M_1 \wedge\cdots\wedge y^M_c.$$
	When we apply the matrix \eqref{E: orientation matrix}, we obtain the form in $V \oplus W$ given by
	$$(-1)^{w(m+1)} (p^V_1+p^W_1) \wedge\cdots\wedge (p^V_a+p^W_a) \wedge x^V_1 \wedge\cdots\wedge x^V_b \wedge (- p^W_1) \wedge\cdots\wedge (-p^W_a) \wedge (-y^W_1) \wedge\cdots\wedge(- y^W_c).$$
	As the number of terms with a negative sign is $w$, this expression simplifies to
	$$(-1)^{wm} p^V_1 \wedge\cdots\wedge p^V_a \wedge x^V_1 \wedge\cdots\wedge x^V_b \wedge p^W_1 \wedge\cdots\wedge p^W_a \wedge y^W_1 \wedge\cdots\wedge y^W_c.$$
	But this is now precisely $(-1)^{wm}$ times the concatenation orientation of $V \oplus W$ as desired for the definition of the fiber product orientation of $P$.
\end{proof}

\begin{example}\label{Ex: intersection orientation}
	As an example, let $M = \R^3$ oriented by the standard ordered basis $(e_x,e_y,e_z)$.
	Let $V$ be the $z = 0$ plane oriented by the ordered basis vectors $(e_x,e_y)$, and let $W$ be the $x = 0$ plane oriented by the ordered basis vectors $(e_y,e_z)$.
	The intersection $P$ is the $y$ axis.
	We claim that $P$ has the fiber product orientation by $-e_y$.
	Indeed, assuming so we have $\beta_V = e_x \wedge e_y = -e_y \wedge e_x = (-e_y) \wedge e_x$, so $\beta_{\nu W} = e_x$, and $\beta_W = e_y \wedge e_z = (-e_y) \wedge (-e_z)$, so $\beta_{\nu V} = -e_z$.
	And then $$\beta_P \wedge \beta_{\nu V} \wedge \beta_{\nu W} = -e_y \wedge (-e_z) \wedge e_x = e_x \wedge e_y \wedge e_z = \beta_M,$$
	as required.
\end{example}

\begin{corollary}\label{C: orient complementary intersection}
	Suppose $V$ and $W$ have complementary dimensions so that they intersect in a point.
	Then the fiber product orientation of the point is positive if and only if $\beta_{W} \wedge \beta_{V} = \beta_M$.
\end{corollary}

\begin{proof}
	In this case, $\nu W = V$, $\nu V = W$, and $\beta_P = \pm 1 \in \R$.
	If $\beta_P = 1$, then the formula from \cref{P: orient intersection} becomes exactly the formula of the corollary.
\end{proof}

\begin{remark}
	The corollary shows that the fiber product orientation is \textit{not} necessarily the expected concatenation orientation in the case of transverse complementary embeddings.
\end{remark}

\begin{example}
	Let $M = \R^2$ with the standard orientation that we can write $e_x \wedge e_y$.
	Let $V$ be the $x$-axis with orientation $e_x$ and $W$ be the $y$-axis oriented by $e_y$.
	Then $\beta_W \wedge \beta_V = e_y \wedge e_x$, while $\beta_M = e_x \wedge e_y = -e_y \wedge e_x$.
	So the fiber product orientation of the intersection point is the negative one.
	This runs against the standard convention for transverse intersections of manifolds of complementary dimension, but we nonetheless favor this overall convention for orienting fiber products due to the properties and uniqueness result of \cite{RamBas09}.
\end{example}

\index{orientation!of fiber product|)}

\subsection{Co-orientations}\label{S: co-orientations}

To define co-orientations, we recall our definition of an orientation of a bundle from \cref{D: det bundle} as an equivalence class, up to positive scalar multiplication, of an everywhere non-zero section of the top exterior power of the bundle.
This motivates the following.

\begin{definition}\label{D: co-orientations}\index{co-orientation|textbf}
	A \textbf{co-orientation} $\omega_g$ of a \textit{continuous}\footnote{We will most often be interested in the case of $g$ smooth, but continuous co-orientable maps do come up, for example in \cref{S: basic properties} where we consider covariant functoriality of geometric cohomology with respect to continuous maps.} map $g \colon W \to M$ of manifolds with corners is an equivalence class, up to positive scalar multiplication, of a nowhere zero section of the line bundle $\Hom(\Or(TW), \Or(g^*TM)) \cong \Hom(\Or(TW), g^*\Or(TM))$.
	Equivalently, a co-orientation is a choice of isomorphism between the associated orientation cover $\Or(TW)$ and the pullback of the associated orientation cover $\Or(TM)$.

	In particular, a co-orientation exists if and only if $\Hom(\Or(TW), \Or(g^*TM))$ is a trivial line bundle, in which case we say that $g$ is \textbf{co-orientable}.\index{co-orientability|textbf}
\end{definition}

Thus, if $W$ is connected and $g \colon W \to M$ is co-orientable, there are exactly two co-orientations, which are \textbf{opposite}\index{co-orientation!opposite co-orientation} to one another; we write the opposite of $\omega_g$ as $-\omega_g$.
In particular, for connected $W$ a choice of co-orientation at a single point determines a co-orientation globally when $g$ is co-orientable (analogously to orientations).
Also, just as most manifolds do not possess a preferred orientation, most maps $g \colon W \to M$ do not carry a natural choice of co-orientation, with the following notable exception.

\begin{definition}\label{D: tautological co-orientation}\index{co-orientation!tautological}
	Suppose $g$ is a diffeomorphism, or more generally a codimension-0 immersion. In this case, the top exterior power of $Dg$,
	$$\textstyle{\bigwedge^w} Dg \colon \textstyle{\bigwedge^w} TW \to \textstyle{\bigwedge^w} g^*(TM),$$
	provides a \textbf{tautological co-orientation}.
\end{definition}

The local triviality of the determinant line bundle of a manifold means being able to choose a consistent basis vector over sufficiently small neighborhoods.
Again, we call such a choice of basis vector around a point in $W$ a \textbf{local orientation}\index{orientation!local}, and, as for global orientations, often denote a local orientation by $\beta_W$.\index{$B$@$\beta_W$}
Again, abusing notation, we also often allow $\beta_W$ to refer to its equivalence class up to multiplication by a positive scalar.
We identify $\beta_W$ with a local choice of (equivalence class of) non-zero section of\footnote{As usual, if $\dim(W) = 0$ we identify $\bigwedge^0 TW$ with the trivial $\R$ bundle, and, when forming exterior products, multiplication by a section is treated as scalar multiplication.} $\bigwedge^w TW$ in a neighborhood of a point $x$ in $W$ or, equivalently, a local smoothly varying ordered basis for the fibers of $TW$.
We typically do not specify the point $x$, though when necessary we write $\beta_{W,x}$.
We then use ordered-pair notation for co-orientation homomorphisms, with $\omega_g = (\beta_W, \beta_M)$ being the \textbf{local co-orientation}\index{co-orientation!local}\index{co-orientation!local notation} that sends the local orientation $\beta_W$ at $x \in W$ to the local orientation $\beta_M$ for $g^*(TM)$.
We will often further abuse notation by neglecting the pullback and treating $\beta_M$ as a local orientation at $g(x)$ in $M$.
We write the opposite co-orientation $(\beta_W,-\beta_M) = (-\beta_W,\beta_M)$ as $-(\beta_W,\beta_M)$.
As a co-orientation at a point completely determines the co-orientation of a co-orientable map for connected $W$, it is useful to cheat further and write $\omega_g = (\beta_W,\beta_M)$ for appropriate $\beta_W$ and $\beta_M$ with the chosen points $x$ implicit.

A manifold is orientable if and only if the orientation cover is trivial.
So if $M$ is orientable, $\Or(g^*(TM))$ is trivial, and co-orientability of $g \colon W \to M$ implies that $W$ is orientable.
Moreover, an orientation on $M$ along with a co-orientation of $g$ gives rise to an \textbf{induced orientation}\index{orientation!induced by a co-orientation} of $W$.
Explicitly, if $\beta_M$ denotes the global orientation of $M$, then we orient $W$ at each point by the $\beta_W$ such that $\omega_g = (\beta_W,\beta_M)$.
Conversely, if $M$ and $W$ are both oriented, say by $\beta_M$ and $\beta_W$ respectively, we have the \textbf{induced co-orientation}\index{co-orientation!induced by map of oriented manifolds} given by $\omega_g = (\beta_W,\beta_M)$ at each point of $W$.
On the other hand, it is not true that if we have an orientation of $W$ and a co-orientation of $g \colon W \to M$ then we obtain an orientation of $M$.
For example, if $W$ is orientable, any constant map to $M$ is co-orientable, regardless of whether or not $M$ is orientable.

\begin{notation}\label{N: hat check}
	Suppose $M$ is an oriented manifold with corners and $g \colon W \to M$ is a map from another manifold with corners.
	By convention, below we will often simply write $W$ to refer to $W$ together with its given map to $M$.
	If $W$ comes equipped with an orientation, we may write $\hat W$\index{$W$@$\hat W$} to refer to $W$ with the induced co-orientation on $g$.
	If $W$ (technically the map $g \colon W \to M$) comes equipped with a co-orientation, we may write $\check W$\index{$W$@$\check W$} to refer to $W$ with its induced orientation.
	For spaces with longer expressions, such as $V \times_M W$, we will sometimes write instead $(-)\,\hat{\vrule height1.3ex width0pt}$ or $(-)\,\check{\vrule height1.3ex width0pt}$.
	These are inverse operations in the sense that if $W$ is oriented then $(\hat W)\,\check{\vrule height1.3ex width0pt} = W$ and if $g$ is co-oriented then  $(\check W)\,\hat{\vrule height1.3ex width0pt} = W$.

	The notation is meant to evoke raising and lowering indices, as maps from oriented objects will eventually play the role of chains and co-oriented maps will eventually play the role of cochains.
\end{notation}

Here is another point of view on co-orientations, immediate from the definition.
Recall that the fundamental group of a manifold acts on classes of local orientations as the deck transformations of the orientation cover.
A map is co-orientable if it holds that a loop in $W$ acts nontrivially on a local orientation of $W$ if and only if its image in $M$ acts nontrivially on a local orientation of $M$.
Explicitly, if $g \colon W \to M$ is co-orientable, following a loop in $W$ may change the local orientation pair
$(\beta_W, \beta_M)$ to $(-\beta_W, -\beta_M)$, but these pairs define equivalent co-orientations.

Similarly, to compare local constructions at different points, it is useful to use paths.
Suppose $\gamma \colon I \to W$ is a path with $\gamma(0) = x$ and $\gamma(1) = y$.
We can choose a lift $\td \gamma$ of $\gamma$ to the complement of the $0$-section of $\Or(TW)$ such that $\td \gamma(0)$ is in the equivalence class of $\beta_{W,x}$.
We then define $\gamma_*\beta_{W}$ to be the equivalence class of $\td \gamma(1)$.
Likewise, we define $(g\gamma)_*\beta_M$ by a lift of $g\gamma$ to the complement of the zero section of $\Or(TM)$.
Then $\gamma_*\beta_{W}$ and $(g\gamma)_*\beta_M$ depend on $\gamma$, but if $g \colon W \to M$ is co-orientable the pair $(\gamma_*\beta_{W}, (g\gamma)_*\beta_M)$ is independent of $\gamma$ as this data also determines a non-vanishing lift of $\gamma$ in $\Hom(\Or(TW),g^*\Or(TM))$, which is trivial if $g$ is co-orientable.
In particular, if $g$ is co-oriented and $(\beta_{W}, \beta_M)$ represents the choice of co-orientation locally at $x$, then $\gamma_*(\beta_W,\beta_M) \defeq (\gamma_*\beta_{W}, (g\gamma)_*\beta_M)$ will represent the same co-orientation locally at $y$.
We refer to this as \textbf{transporting}\index{co-orientation!transportation} the co-orientation from $x$ to $y$.

\begin{example}
	Let $g$ be any map $g \colon S^1 \to S^2$.
	As $S^1$ and $S^2$ are orientable, $g$ is co-orientable.
	If we choose a local orientation vector $e_{\theta}$ at any point $x \in S^1$ and latitude/longitude coordinates $\phi,\psi$ at $g(x)$ so that $e_\phi \wedge e_\psi$ is a local orientation in a neighborhood of $g(x)$, then the two possible co-orientations for $g$ can be written $(e_\theta, e_\phi \wedge e_\psi)$ and $-(e_\theta, e_\phi \wedge e_\psi) = (-e_\theta, e_\phi \wedge e_\psi) = (e_\theta,- e_\phi \wedge e_\psi)$.
	While the notation explicitly references a local orientation at a point, this is sufficient to determine the co-orientation globally.
	In what follows we will often demonstrate properties of co-orientations by showing that they hold locally at an arbitrary point but do not depend on the choice of point.

	As another example, consider the standard embedding $g \colon \R P^2 \into \R P^4$.
	Choosing local orientations $e_1 \wedge e_2$ at some $x \in \R P^2$ and $f_1 \wedge f_2 \wedge f_3 \wedge f_4$ at $g(x)$, the two co-orientations are $(e_1 \wedge e_2, f_1 \wedge f_2 \wedge f_3 \wedge f_4)$ and its opposite.
	If $\gamma$ is a loop that reverses the orientation of $\R P^2$ then it also reverses the orientation of $\R P^4$, so $\gamma_*(e_1 \wedge e_2, f_1 \wedge f_2 \wedge f_3 \wedge f_4) = (-e_1 \wedge e_2,- f_1 \wedge f_2 \wedge f_3 \wedge f_4) = (e_1 \wedge e_2, f_1 \wedge f_2 \wedge f_3 \wedge f_4)$, reflecting that $g$ is co-orientable.

	By contrast, no embedding of the M\"obius strip in $\R^3$ is co-orientable.
\end{example}

\begin{remark}\label{R: cooriented composition}
	Co-oriented maps compose\index{co-orientation!compositions of} in an immediate way, forming a category.
	Namely, given $V \xr{f} W \xr{g} M$ and co-orientations $\Or(TV) \to \Or(f^*TW)$ and $\Or(TW) \to \Or(g^*TM)$, we simply compose the former with the pullback of the latter via $f^*$, recalling that $f^*(\Or(E)) = \Or(f^*E)$ in a natural way.
	We will refer to this simply as composing co-orientations and write the composition in symbols as $\omega_f*\omega_g$.
	Warning: note that we write the terms in the order $\omega_f*\omega_g$ for the map $g \circ f$.
	This is more convenient when writing out co-orientations using the local orientations as we obtain expressions such as $(\beta_V, \beta_W)*(\beta_W,\beta_M) = (\beta_V,\beta_M)$.
\end{remark}

\begin{notation}\label{N: implicit notation}
	It will be useful in notation to sometimes leave the maps, codomains, and co-orientations all implicit once they have already been established and just write $V$ to represent the co-oriented map $f \colon V \to M$.
	In this, case we write $-V$\index{$V$@$-V$}\index{co-orientation!notation for opposite} to refer to the same map with the opposite co-orientation.
\end{notation}

Just as is the case for orientations, when working with co-orientations it suffices to consider the interior of the manifold, in this case of the domain.

\begin{proposition}\label{P: interior co-orientation}\index{co-orientability!only depends on interior}
	Let $g \colon W \to M$ be a map of manifolds with corners. Then $g$ is co-orientable if and only if its restriction to $S^0(W)$, the interior of $W$, is co-orientable.
\end{proposition}
\begin{proof}
	It is clear from the definition that if $g$ is co-orientable then so is its restriction to any open set of $W$.

	Conversely, suppose the restriction of $g$ to $S^0(W)$ is co-orientable. As $W$ is a topological manifold with boundary, which we will here denote $bd(W)$, the collaring theorem tells us that $bd(W)$ possesses a collar \cite[Theorem 2]{Bro62}. Let $C$ be a closed collar of $bd(W)$. Then $C$ is homeomorphic to $bd(W)\times I$ with $bd (W)$ identified under the homeomorphism with $bd(W) \times 0$. Let $B$ be the image of $bd(W) \times 1$ in $W$, and let $int(C)$ be the image of $bd(W) \times [0,1)$. Since $g|_{S^0(W)}$ is co-orientable, any choice of co-orientation restricts to a bundle isomorphism $\bigwedge TW|_B \to \bigwedge g^*\Or(TM)|_B$. Due to the product structure of the collar, the restriction of $\bigwedge TW$ to $C$ is isomorphic to $\bigwedge TW|_B \times I$ and similarly for $\bigwedge g^*\Or(TM)$ \cite[Theorem 3.4.4]{Hus94}, and so the bundle isomorphism over $B$ can be extended to a bundle isomorphism over $C$. As we chose the bundle isomorphism over $B$ to be the restriction of an isomorphism over $g|_{S^0(W)}$, we can glue together the bundle isomorphisms over $C$ and over $W-int(C)$ to obtain an isomorphism over all of $W$.
\end{proof}

\subsection{Normal co-orientations of immersions and co-orientations of boundaries}\label{S: normal orientation}

Once again, a key example is when $g$ is an immersion, which is co-orientable if and only if its normal bundle is orientable\footnote{Recall that technically all bundles are over $W$, though our convention is to elide that in the notation.
	Hence we can consider $W$ to have a normal bundle even if $g$ is merely an immersion and not actually an embedding.
	The normal bundle can be identified with $g^*(TM)/TW$ after identifying $TW$ with a sub-bundle of $g^*(TM)$ using the differential.
	In any case, locally in the neighborhood of any point of $W$ one has the usual identification of the normal bundle with a tubular neighborhood of the image, which suffices for our purposes here.}.
Specifically, if $g \colon W \to M$ is an immersion, letting $\nu W$ denote the normal bundle, we have $TW \oplus \nu W \cong g^*TM$.
So, taking $w = \dim(W)$ and $m = \dim(M)$, a co-orientation is a nowhere-zero map from $\bigwedge^w TW$ to $\bigwedge^m g^*TM = \bigwedge^m (TW \oplus \nu W) \cong \bigwedge^w TW \otimes \bigwedge^{m-w}\nu W$.
Such a nowhere-zero map exists if and only $\bigwedge^{m-w}\nu W$ is a trivial line bundle, i.e.\ if $\nu W$ is orientable.

Given a specific orientation of $\nu W$, we specify a standard choice of \textbf{normal co-orientation}\index{co-orientation!normal} for the immersion $g \colon W \to M$ by the following local construction:

\begin{definition}\label{normal co-or}
	Let $g \colon W \to M$ be an immersion with normal bundle $\nu$ locally oriented (at some point of $W$) by $\beta_\nu$.
	Define the \textbf{normal co-orientation} associated to $\beta_\nu$ locally by the pair $\omega_{\nu} = (\beta_W, \beta_W \wedge \beta_\nu)$, where $\beta_W$ is any choice of a local orientation of $W$.
\end{definition}

This construction is independent of the choice of $\beta_W$, as reversing the orientation of $\beta_W$ gives
$$(-\beta_W, -\beta_W \wedge \beta_\nu) = (\beta_W, \beta_W \wedge \beta_\nu).$$
If the normal bundle to $W$ is oriented globally on $W$ then the construction is also independent of the point at which it is carried out since if $\gamma$ is a path from $x$ to $y$ with $\beta_W$ and $\beta_\nu$ constructed at $x$ then
\begin{equation*}
	(\gamma_*\beta_{W,x}\, , (g\gamma)_* (\beta_{W,x} \wedge \beta_{\nu,x})) =
	(\gamma_*\beta_{W,x}\, , (\gamma_* \beta_{W,x}) \wedge (g\gamma)_*\beta_{\nu,x}) =
	(\gamma_*\beta_{W,x}\, , (\gamma_* \beta_{W,x}) \wedge \beta_{\nu,y}),
\end{equation*}
using that $\nu$ is assumed oriented and that, via the immersion, we can treat a path in $W$ as a path in $M$ and tangent vectors of $W$ as tangent vectors of $M$.
Now taking $\beta_{W,y}$ to be $\gamma_* \beta_{W,x}$, the expression on the right again has the prescribed form.
So if the normal bundle to $W$ is oriented, these local choices determine a global co-orientation of $W \to M$.
If the normal bundle is orientable, one can conversely orient the normal bundle via this formula if one is given a co-orientation: choose $\beta_\nu$ so that $(\beta_W, \beta_W \wedge \beta_\nu)$ is the co-orientation of the immersion.

As we shall see, signs in the theory of co-orientations are highly dependent on choices.
One such choice in this definition is whether to append the local normal orientation before or after the local tangent orientation.

\subsubsection{Quillen co-orientations}\label{S: Quillen}

There is a useful alternative, though equivalent, definition of co-orientations due to Quillen \cite{Quil71} that only involves orientations of normal bundles as in the preceding section.\footnote{Quillen's context was slightly different.
	He assumed the normal bundles to have complex structures and so called these \textit{complex orientations}.}

\begin{lemma}\label{L: Quillen}\index{map, decomposition as an embedding then a projection}
	A map $g \colon W \to M$ is co-orientable if and only if for some $N \in \Z_{\geq 0}$ it factors as the composition $W \into M \times \R^N \to M$ of an embedding and a projection such that the image of $W$ in $M \times \R^N$ has an orientable normal bundle.
\end{lemma}

\begin{proof}
	We first note that such a smooth factorization always exists by \cref{C: embed V}.

	Next we observe that $T(M \times \R^N) \cong \pi^*(TM) \oplus \underline{\R}^N$ where $\underline{\R}^N$ is the trivial $\R^N$ bundle over $M \times \R^N$.
	As $\Or( \underline\R^N)$ is a trivial line bundle,
	\begin{equation*}
		\Or(T(M \times \R^N)) \cong \Or(\pi^*(TM)) \otimes \Or(\underline\R^N) \cong \Or(\pi^*TM).
	\end{equation*}
	Thus $\pi$ is always co-orientable.
	Furthermore, we see that
	\begin{equation*}
		\Or(g^*TM) \cong \Or(e^*\pi^*TM) \cong e^*\Or(\pi^*TM) \cong e^*\Or(T(M \times \R^N)) \cong \Or(e^*T(M \times \R^N)).
	\end{equation*}
	So $e$ is co-orientable if and only if $g$ is co-orientable.
	But $e$ is an immersion, and so it is co-orientable if and only if the normal bundle of the image is orientable by the discussion preceding \cref{normal co-or}.
\end{proof}

\begin{definition}\label{D: Quillen normal or}
	If $(e_1, \ldots, e_n)$ is the standard ordered basis of $\R^N$, we denote the corresponding orientation from by $\beta_E = e_1 \wedge \cdots \wedge e_N$.
	Then the projection $\pi \colon M \times \R^N \to M$ has a canonical co-orientation $(\beta_M \wedge \beta_E, \beta_M)$ that is well defined as $\R^N$ is contractible.
	If $g \colon W \to M$ is co-oriented and $e \colon W \into M \times \R^N$ is an immersion with normal bundle $\nu$, we then define the \textbf{compatible normal orientation} or \textbf{Quillen normal orientation}\index{co-orientation!normal|textbf}\index{co-orientation!Quillen|textbf} of $\nu$ so that the composition of co-orientations $(\beta_W,\beta_W \wedge \beta_\nu)$ with $(\beta_M \wedge \beta_E,\beta_M)$ is the given co-orientation of $g$.
	In other words, if $(\beta_W,\beta_M)$ is the given co-orientation of $g$, then the compatible normal orientation of $\nu$ is such that $\beta_W \wedge \beta_\nu = \beta_M \wedge \beta_E$.

	We sometimes speak of the entire structure $W \into M \times \R^N \to M$ with an orientation of $\nu$ as a \textbf{Quillen co-orientation}\index{co-orientation!normal|textbf}\index{co-orientation!Quillen|textbf} for $W \to M$ or as a ``compatible Quillen co-orientation'' if we have already specified a co-orientation for $W \to M$ and we wish to choose the Quillen co-orientation that agrees with it.
\end{definition}

\begin{remark}\label{R: immersion}
	In particular, if $g \colon W \to M$ is an immersion, then, by taking $N = 0$, a co-orientation of $g$ is equivalent to a Quillen normal orientation of the normal bundle $\nu$ of $W$ in $M$.
	In particular, the co-orientation is given locally by $(\beta_W, \beta_M)$ if and only if $\nu$ is oriented so that $\beta_W \wedge \beta_\nu = \beta_M$.
	If $g$ is a codimension-$0$ immersion, then $\nu$ will be $0$-dimensional, and if the co-orientation is the tautological one then $\beta_\nu$ will be the positive orientation at each point.
\end{remark}

Lipyanskiy's definition of co-orientation in \cite{Lipy14} factors a proper map through a map which is surjective onto $TM$, rather than injective from $TW$ as in the Quillen approach.
We discuss Lipyanskiy's co-orientations in an appendix at the end of this section; see \cref{S: Lipyanskiy co-orientations}.

\subsubsection{Co-orientations of boundaries}

Given a co-oriented map $g \colon W \to M$ where $W$ is a manifold with corners, we can use a normal co-orientation of $\bd W$ in $W$ together with the composition of co-orientations noted in \cref{R: cooriented composition} to define ``boundary co-orientations'':

\begin{definition}\label{D: boundary co-orientation}
	Let $\nu$ be the normal bundle of $\bd W$ as an immersed manifold with corners of $W$ (recall that $i_{\bd W} \colon \bd W \to W$ is a smooth immersion by \cite[Theorem 3.4]{Joy12}).
	The \textbf{standard co-orientation of a boundary immersion}\index{co-orientation!standard boundary co-orientation} $i_{\bd W} \colon \bd W \into W$ is the normal co-orientation (see \cref{normal co-or}) associated to the \textit{inward}-pointing\footnote{The outward-pointing normal would also work to provide a co-orientation convention for which the Leibniz formula of \cref{S: co-orient pullbacks} holds.
		However, using an outward normal is not consistent with the intersection map $\mc I$ of \cref{S: intersection map} being a chain map with our other conventions, while using the inward normal does make $\mc I$ a chain map.} orientation of $\nu$.

	If $g \colon W \to M$ is co-oriented, the \textbf{induced co-orientation}\index{co-orientation!induced} or \textbf{boundary co-orientation}\index{co-orientation!boundary} of the composition $gi_{\bd W}$ is the composition of the standard co-orientation of $i_{\bd W}$ with the given co-orientation of $g \colon W \to M$.
	We write $\bd g \colon \bd W \to M$ to denote $gi_{\bd W}$ with its induced co-orientation.
\end{definition}

In \cref{S: geometric cochains}, we will use induced co-orientations on boundaries to define the differential in the geometric cochain complex.

\begin{example}\label{E: splitting example 1}\index{co-orientation!codimension 0 embedding}
	Suppose $g \colon W \into M$ is the embedding of a codimension-$0$ submanifold of $M$.
	In this case $TW$ is the pullback of $TM$, and we have the tautological co-orientation of \cref{D: tautological co-orientation} that, slightly abusing notation, we can write as $(\beta_M,\beta_M)$.
	A particularly important pair of examples is given by the inclusions $g^- \colon (-\infty, 0] \into \R$ and $g^+ \colon [0,\infty) \into \R$, each tautologically co-oriented by $(e_1,e_1)$, where $e_1$ is the standard unit vector in $\R$.

	Next\index{co-orientation!codimension 1 embedding} consider the submanifold consisting of the point $0 \in \R$.
	Its tangent space has trivial determinant line bundle, and we can choose the basis element to be the element $1 \in \bigwedge^0 T0 \cong \R$.
	If we give the normal bundle to $0$ in $\R$ the standard orientation in the positive direction, denoted by the standard basis vector $e_1$, then the corresponding normal orientation of the inclusion of $0$ into $\R$ will be $(1, 1 \wedge e_1) = (1, e_1)$.
	By \cref{normal co-or,D: boundary co-orientation}, the standard co-orientation of the boundary inclusion $\{0\} \into (-\infty, 0]$ is $(1, 1 \wedge -e_1) = (1, -e_1)$.
	The boundary co-orientation of the inclusion $\{0\} \to \R$ induced by the tautological co-orientation of the inclusion $g^- \colon (-\infty, 0] \to \R$ is then the composition of $(1, -e_1)$ with $(e_1, e_1)=(-e_1,-e_1)$, which is again $(1,-e_1) = -(1,e_1)$ (all bases interpreted in the appropriate spaces).
	As the inward normal to $[0,\infty)$ at $0$ is $e_1$, the inclusion $g^+ \colon [0,\infty) \to \R$ induces the opposite co-orientation $(1, e_1)$ on the inclusion $\{0\} \to \R$.
	So the inclusion $\{0\} \into \R$ has opposite co-orientations as the boundary of $(-\infty, 0] \into \R$ versus as the boundary of $[0,\infty) \into \R$.
\end{example}

The following example will be very important later for breaking geometric chains and cochains into pieces.

\begin{example}\label{E: manifold decomposition}\index{co-orientation!codimension 1 embedding}
	Suppose we have a smooth map $\varphi \colon W \to \R$ and that $0$ is a regular value\index{regular value} of $\varphi$, meaning by \cref{D: regular value} that $\varphi$ is transverse to the inclusion of $0$ into $\R$.
	Consider the spaces\index{$W^0$} $W^0 = (\varphi)^{-1}(0)$, $W^- = (\varphi)^{-1}((-\infty,0])$, and $W^+ = (\varphi)^{-1}([0,\infty))$.
	We can identify $W^0$ with the fiber product $0 \times_{\R} W$ of $\varphi$ with the embedding of $0$ into $\R$, as our canonical realization of this fiber product is
	$$\{(0,x) \in 0 \times W \mid \varphi (x)= 0\}.$$
	Similarly, $W^-$\index{$W^-$} and $W^+$\index{$W^+$} are $(-\infty,0] \times_{\R} W$ and $[0,\infty) \times_{\R} W$,
	as $$[0,\infty) \times_{\R} W = \{(t,x) \in [0,\infty) \times W \mid \varphi(x) = t\},$$
	which is simply the graph of $\varphi$ over $W^+$;
	so there is a diffeomorphism between the fiber product and $W^+$ determined by the inverse maps $(t,x) \mapsto x$ and $x \mapsto (\varphi(x), x)$.
	The case of $W^-$ is analogous.
	By \cref{pullback}, the spaces $W^0$, $W^+$, and $W^-$ are manifolds with corners.
	In fact, $W^0$ is a submanifold (with corners) of $W$ by \cite[Proposition 4.2.9]{MaDo92}.
	By the boundary formula of \cref{P: product boundary} (which is purely topological, we do not yet take into account co-orientations),
	$$\bd (W^-)= \left(0 \times_{\R} W\right) \sqcup \left((-\infty,0] \times_{\R} \bd W\right) = W^0 \sqcup (\bd W)^-,$$
	so $W^0$ is a boundary component of $W^-$ and similarly for $W^+$.
	We also observe that as the inclusions of $0$, $(-\infty,0]$, and $[0,\infty)$ into $\R$ are all proper, it will follow from \cref{L: co-orientable pullback}, below, that the inclusions of $W^0$, $W^-$, and $W^+$ into $W$ are all proper.

	Now, at least over its interior, the normal bundle of $W^0$ in $W$ is the pullback via $\varphi$ of the normal bundle of $0$ in $\R$ by \cref{L: normal pullback}. Giving the normal bundle of $0$ in $\R$ the standard orientation in the positive direction, denoted by the basis vector $e_1$, this orientation pulls back to a natural orientation of the normal bundle to the interior of $W^0$ in $W$. Applying \cref{normal co-or,P: interior co-orientation}, we obtain a normal co-orientation for $W^0$ in $W$, which we call the \textbf{co-orientation of $W^0$\index{co-orientation!of $W^0$|(}\index{$W^0$!induced co-orientation|(} in $W$ induced by $\varphi$}. At a point $x$ in the interior of $W^0$, this is the co-orientation $(\beta_{W^0}, \beta_{W^0} \wedge v)$, where $v$ is any vector in $T_x W$ such that $D_x \varphi (v) \in T_0 \R$ points in the positive direction.

	We also have the co-orientations of the inclusion of $W^0$ into $W$ arising from taking the boundary of the inclusions of $W^{\pm}$ into $W$.
	As these boundary co-orientations are the normal co-orientations associated with inward pointing normal vectors, we have completely analogously to \cref{E: splitting example 1} that the boundary co-orientation for the composite $W^0 \to W^- \into W$ is the opposite of the co-orientation of $W^0\into W$ induced by $\varphi$, while the boundary co-orientation from the composite $W^0 \to W^+ \into W$ agrees with the co-orientation of $W^0\into W$ induced by $\varphi$.

	Next we suppose that $\varphi$ is the composition $\varphi = \phi g$ of two smooth maps $W \xr{g} M \xr{\phi}\R$, where $M$ is a manifold without boundary.\footnote{This may at first seem a bit artificial, but such a setting will be important below when $W \to M$ represents a geometric cochain; see in particular the splitting and creasing constructions of \cref{S: splitting and creasing} that will be used for breaking up geometric chains and cochains into smaller pieces.}
	We further suppose $0$ is a regular value for $\phi$ and for $\phi g$; see \cref{D: regular value}.
	Note that when $0$ is a regular value for $\phi$, then \cref{L: transverse to pullback} implies that $0$ being a regular value for $\varphi = \phi g$ is equivalent to $g$ being transverse to $M^0 = \phi^{-1}(0)$, which is a codimension-one submanifold of $M$ by the preceding discussion.
	By the applying the preceding discussion to $\phi$, we obtain spaces $M^0$ and $M^\pm$, and applying the discussion to $\varphi = \phi g$, we obtain spaces $W^0$ and $W^\pm$.
	Furthermore, by considering $M^0$ and $M^\pm$ as pullbacks over $M$ as above, we observe that $W^0$ and $W^\pm$ are diffeomorphic to the respective fiber products $M^0 \times_M W$ and $M^{\pm} \times_M W$ by \cref{P: pullback functoriality}.\index{co-orientation!of $W^0$|)}\index{$W^0$!induced co-orientation|)}

	Now suppose that $g$ is co-oriented with local representatives $(\beta_W, \beta_M)$.
	If we precompose this co-orientation of $g$ with the tautological co-orientations of the inclusions $W^\pm \into W$, we can co-orient $g|_{W^\pm}$ also by $(\beta_W, \beta_M)$.
	This is simply the restriction to $W^\pm$ of the co-orientation of $g$.
	Similarly, by composing the co-orientation of $g$ with the $\varphi$-induced co-orientation of $W^0 \into W$, we obtain the \textbf{co-orientation of $g|_{W^0} \colon W^0 \to M$ induced by $\phi$}, which in symbols is just $(\beta_{W^0}, \beta_{W^0} \wedge v)*(\beta_W,\beta_M)$, where $D\varphi(v)$ points in the positive direction.
	Furthermore, as in \cref{E: splitting example 1}, the co-orientation of $g|_{W^0} \colon W^0 \to M$ induced by $\phi$ disagrees with its boundary co-orientation as a component of $\bd(g|_{W^-})$, while it agrees with its co-orientation as a component of $\bd(g|_{W^+})$.

	When $g$ is co-oriented, our co-orientation computations here for the restrictions of $g$ to $W^0$ and $W^\pm$ will later be seen to be consistent with the co-orientations we define for fiber products of co-oriented maps; see \cref{E: codim 0 and 1 co-or as fiber products,E: codim 1 pullbacks}.
	In \cref{C: co-orient W0}, we will use that technology to see that $\bd (W^0)$ and $(\bd W)^0$ agree as spaces, but their maps to $M$ will have opposite co-orientations under our conventions, i.e.\ ``$\bd (W^0) = -(\bd W)^0$,'' eliding the maps.
\end{example}

\subsubsection{Co-orientations of boundaries of boundaries}\index{co-orientation!of $\bd^2$}

In order to form a chain complex of geometric cochains in \cref{S: geometric cochains}, we will need a result about co-orientations of $\bd^2 W$.
Recall from \cref{S: boundaries} that Proposition 2.9 of \cite{Joy12} identifies $\bd^2 W$ with the set of points $(x,\bb_1,\bb_2)$ with $x \in W$ and the $\bb_i$ encoding distinct local boundary components.
The map $i_{\bd^2 W} \colon \bd^2 W \to W$ takes $(x,\bb_1,\bb_2)$ to $x$.
The manifold with corners $\bd^2 W$ is equipped with a canonical diffeomorphism $\rho$ defined by $(x,\bb_1,\bb_2) \mapsto (x,\bb_2,\bb_1)$.

\begin{lemma}\label{L: boundary2}
	Suppose $i_{\bd^2 W} \colon \bd^2 W \to W$ is co-oriented via the composition of boundary co-orientations $\bd^2 W \to \bd W \to W$, and suppose $\rho \colon \bd^2 W \to \bd^2 W$ is given its tautological co-orientation (see \cref{D: tautological co-orientation}).
	Then $i_{\bd^2 W}$ and $i_{\bd^2 W}\rho$ have opposite co-orientations.
\end{lemma}

\begin{proof}
	It suffices to consider points $(x,\bb_1,\bb_2) \in \bd^2 W$ with $x \in S^2(W)$, as such points fill out the interior of $\bd^2 W$.
	In $W$, such $x$ have neighborhoods of the form $[0,\infty)^2 \times \R^{w-2}$ with $x$ at the origin.
	We identify $[0,\infty)^2$ with the first quadrant of $\R^2$, letting $X$ and $Y$ denote the non-negative $x$ and $y$ axes.
	We let $\bb_X$ and $\bb_Y$ be the corresponding local boundary components.
	Then the preimage in $\bd^2 W$ of a small neighborhood $U$ of $x$ in $S^2(W)$ consists of two copies of $U$ that we can write $(U,\bb_X,\bb_Y)$ and $(U,\bb_Y,\bb_X)$.
	The notation indicates that we think of the first copy of $U$ as embedding into $X \times \R^{w-2} \subset \bd W$ and the second as embedding into $Y \times \R^{w-2} \subset \bd W$.
	The map $i_{\bd W} \colon \bd W \to W$ then identifies the two copies.
	The map $\rho$ simply interchanges them.\greg{Picture here?}

	Let $\beta_X$ and $\beta_Y$ correspond to the positively-directed tangent vectors in $X$ and $Y$, and let $\beta_U$ be an arbitrary local orientation of $U$.
	Abusing notation, we also write $\beta_U$ for the corresponding local orientations of $(U,\bb_X,\bb_Y)$ and $(U,\bb_Y,\bb_X)$.
	The induced co-orientation on $\rho$ can then be written $(\beta_U,\beta_U)$.
	Up to identifying neighborhoods in $W$ with their local models, the boundary co-orientation of $i_{\bd^2 W}$ on $(U,\bb_X,\bb_Y)$ comes from first mapping it into $X \times \R^{w-2}$ and then into $X \times Y \times \R^{w-2}$.
	So from the definition of boundary co-orientations, this co-orientation is $(\beta_U, \beta_U \wedge \beta_X \wedge \beta_Y)$.
	Analogously, the boundary co-orientation of $i_{\bd^2 W}$ on $(U,\bb_Y,\bb_X)$ is $(\beta_U, \beta_U \wedge \beta_Y \wedge \beta_X)$.
	By composition, the co-orientations of $i_{\bd^2 W}\rho$ on $(U,\bb_X,\bb_Y)$ and $(U,\bb_Y,\bb_X)$ are respectively $(\beta_U, \beta_U \wedge \beta_Y \wedge \beta_X)$ and $(\beta_U, \beta_U \wedge \beta_X \wedge \beta_Y)$ as first we interchange then embed.
	But $\beta_U \wedge \beta_X \wedge \beta_Y = -\beta_U \wedge \beta_Y \wedge \beta_X$, which establishes the lemma.
\end{proof}

\begin{remark}\label{R: bd2 oriented}\index{orientation!of $\bd^2$}
	A similar, though more familiar, argument using \cref{Con: oriented boundary} shows that if $W$ is oriented then $\bd^2 W$ possess an orientation reversing diffeomorphism.
	In this case, we observe that of our two copies of $U$, one is oriented by \cref{Con: oriented boundary} so that $\beta_X \wedge \beta_Y \wedge \beta_U$ is the local orientation of $W$ and the other is oriented so that $\beta_Y \wedge \beta_X \wedge \beta_U$ is the orientation of $W$.
	Thus the two copies of $U$ have opposite orientations, and again the diffeomorphism simply interchanges the two copies.
\end{remark}

\subsection{Co-orientation of homotopies}\label{S: co-oriented homotopy}\index{co-orientation!of homotopy|(}

In this section we develop co-orientations related to homotopies.
As the product of two spaces is the same as their fiber product over a point, we have by \cref{P: product boundary} and rearranging the order of components the topological formula:
\begin{equation*}
	\bd(W \times I) =
	(\bd W \times I) \sqcup (W \times \bd I) =
	(W \times 1) \sqcup (W \times 0) \sqcup (\bd W \times I).
\end{equation*}

Now recall that in general if we have a map $f \colon V \to M$ then we write $\bd f$ for the composition
$$\bd V \xr{i_{\bd V}} V \xr{f} M.$$
We make the following definition:

\begin{definition}\label{D: co-oriented homotopy}
	If $G \colon W \times I \to M$ is a co-oriented map, we say that $G$ is a \textbf{co-oriented homotopy}\index{co-orientation!co-oriented homotopy|textbf} (or simply a \textbf{homotopy} when working with co-orientations is understood) from $g_0 \colon W \to M$ to $g_1 \colon W \to M$ if $\bd G = g_1 \amalg -g_0\amalg H$, where $g_1$, $-g_0$, and $H$ correspond respectively to the compositions of $G$ with the inclusions into $W \times I$ of $W \times 1$, $W \times 0$, and $\bd W \times I$, taking each with its boundary co-orientation as in \cref{D: boundary co-orientation}.
	In particular, then, $g_0$ is the composition
	$$W=W \times 0 \into W \times I \xr{G} M$$
	with the opposite of the boundary co-orientation coming from $G$.
\end{definition}

Note that, by analogy with homotopies involving oriented manifolds, a homotopy from $g_0$ to $g_1$ involves the oppositely co-oriented $-g_0$ in the boundary formula.
In the oriented case, this arises because if we orient $W$ by, say, $\beta_W$ then to orient $W \times I$ we consider $\beta_W \wedge \beta_I$.
Then at one end of the cylinder $\beta_W$ agrees with the boundary orientation of $\bd (W \times I)$ while at the other end it disagrees.
The situation for co-orientations is analogous.

Although we will not need it, we note that by employing appropriate smoothing near the boundaries in order to accomplish transitivity, co-oriented homotopy can be shown to be an equivalence relation on maps $W \to M$.

In our most common use of homotopies, we begin with a co-oriented map $g \colon W \to M$ and want to construct a homotopic co-oriented map.
For this the following lemma is useful.

\begin{lemma}\label{L: co-orientable homotopies}
	If $g \colon W \to M$ is co-orientable and $G \colon W \times I \to M$ is a homotopy with $g = G(-,t_0)$ for some $t_0 \in I$, then $G$ is co-orientable.
	Conversely, if $G \colon W \times I \to M$ is co-orientable, then so is $g = G(-,t_0) \colon W \to M$ for any $t_0 \in I$.
\end{lemma}

\begin{proof}
	Over $W \times t_0$, we have $\Or(T(W \times I)) \cong \Or(TW \oplus TI) \cong \Or(TW) \otimes \Or (TI) \cong \Or(TW)$, while the restriction of $G^*\Or(TM)$ over $W \times t_0$ is just $g^*\Or(TM)$.
	If $G$ is co-orientable then there is a nowhere-vanishing map of line bundles $\Or(T(W \times I)) \to G^*\Or(TM))$, so restricting to $W \times t_0$ and using the above identifications we obtain a nowhere-vanishing map of line bundles $\Or(TW) \to g^*\Or(TM))$ over $W \times t_0$, hence $g$ is co-orientable.
	Conversely,
	if $g$ is co-orientable, there is a nowhere-vanishing map of line bundles $\Or(TW) \to g^*\Or(TM)$ over $W \times t_0$.
	By general bundle theory \cite[Theorem 3.4.4]{Hus94}, any vector bundle $E$ over $W \times I$ is isomorphic to $E_{t_0} \times I$, where $E_{t_0}$ is the restriction of $E$ to $W \times \{t_0\}$.
	So our nowhere-vanishing map of line bundles over $W \times t_0$ extends to a nowhere vanishing map of line bundles $\Or(T(W \times I)) \to G^*\Or(TM))$ over $W \times I$.
	This implies the co-orientability of $G$.
\end{proof}

\begin{definition}\label{D: homotopy co-orientation}
	Suppose $g_0 \colon W \to M$ is co-oriented and $G \colon W \times I \to M$ is a smooth homotopy with $G(-,0) = g_0$.
	Then, by the above lemma, $G$ is co-orientable and clearly there is exactly one choice of co-orientation for $G$ for which the $W \times 0$ component of $\bd G$ is $-g_0$.
	We call this co-orientation the \textbf{co-orientation on $G$ induced by $g_0$}.\index{co-orientation!induced on a homotopy}
	The map $G$ then determines a co-oriented homotopy from $g_0$ to a co-oriented map $g_1 \colon W \to M$ (in the sense of \cref{D: co-oriented homotopy}).
	We call this co-orientation on $g_1 = G(-,1)$ the \textbf{induced co-orientation on $g_1$.}\index{co-orientation!induced by a homotopy}
\end{definition}

\begin{remark}\label{R: stationary homotopy}\index{co-orientation!of homotopy!local description}
	In the above scenario, if $g_0$ is co-oriented locally at $x \in W$ by $(\beta_W,\beta_M)$, then at $(x,0) \in W \times I$, the corresponding local co-orientation of $G$ that yields $-g_0$ as a boundary component of $G$ is $(\beta_W \wedge -\beta_I, \beta_M)$, where $\beta_I$ corresponds to the standard orientation of $I$.
	This follows from $(\beta_W,\beta_W \wedge \beta_I)$ being the boundary co-orientation of $W \times 0 \into W \times I$ as $\beta_I$ corresponds to the inward pointing normal at $0 \in I$.
	As we can take $\beta_W \wedge -\beta_I$ to be a consistent orientation along the path given by $\gamma(t) = (x,t)$, we have, recalling the notation from \cref{S: co-orientations}, that $\gamma_*(\beta_W \wedge -\beta_I, \beta_M) = (\beta_W \wedge -\beta_I,\gamma_*\beta_M)$, and at this end of the homotopy the induced local co-orientation of $g_1$ at $x$ is $(\beta_W,\gamma_*\beta_M)$.
	If $G$ is stationary along $x \times I$, then the co-orientation for $g_1$ at $x$ is again $(\beta_W,\beta_M)$ so that the co-orientations of $g_0$ and $g_1$ agree at $x$.
	This observation will be useful below in showing that pullback co-orientations are well defined.
\end{remark}

\begin{lemma}\label{L: co-oriented homotopy}\index{co-orientation!co-oriented homotopy!boundary of}
	Suppose $G \colon W \times I \to M$ is a co-oriented homotopy from $g_0$ to $g_1$ so that $\bd G = g_1 \amalg -g_0\amalg H$ as in \cref{D: co-oriented homotopy}.
	Then $H$ is a homotopy from $-\bd g_0$ to $-\bd g_1$.
\end{lemma}

\begin{proof}
	By definition, $H$ is co-oriented as a boundary component of the co-oriented map $G$, so it remains to check that the induced co-orientations of the ends of $H$ have the expected signs.
	This could be done directly, but rather we use \cref{L: boundary2}, noting that each copy of $\bd W$ (at the top and bottom of the cylinder) can be considered to be a piece of $\bd^2(W \times I)$.
	In particular, applying \cref{L: boundary2} we see the maps $\bd W \to M$ take opposite co-orientations depending on whether we think of them as first mapping $\bd W$ into $W$ and then identifying $W$ as one end of the cylinder versus first including $\bd W$ into $(\bd W) \times I$ and then mapping this to $W \times I$.
	In both cases, we follow with the map $G$.
	As we think of $g_1$ as defined on $W \subset W \times I$ and let $\bd g_1$ denote its boundary, we see that the corresponding map from the top of the cylinder $(\bd W) \times I$ must be $-\bd g_1$.
	Similarly, the map at the bottom of the cylinder is $\bd g_0$.
	The lemma now follows from \cref{D: co-oriented homotopy}.
\end{proof}

\index{co-orientation!of homotopy|)}

\subsection{Co-orientations of pullbacks and fiber products}\label{S: co-orient pullbacks}

In this section, we define a convention for co-orientations of pullbacks and fiber products.
More specifically, if $f \colon V \to M$ and $g \colon W \to M$ are transverse smooth maps from manifolds with corners to a manifold without boundary and $f$ is co-oriented, we define a co-orientation of the pullback $f^* \colon V \times_M W \to W$.
This does not require $g$ to be co-oriented, but if it is, we can compose with $g$ to also get a co-orientation of the fiber product $V \times_M W \to M$.
Ultimately, this will allow us to define cup products of geometric cochains, and the various properties we will demonstrate for co-orientations of fiber products will be reflected in the standard properties for cup products.

Recall that our canonical realization of the topological pullback $P = V \times_M W$ is defined to be $P = \{(x,y) \in V \times W \mid f(x) = g(y)\}$.
By Joyce \cite[Section 6]{Joy12}, the projections $P \to V$ and $P \to W$ are smooth, and hence so is $f \times_M g \colon P \to M$ given by $(x,y) \to f(x) = g(y)$.
It is not obvious how to define the co-orientations of pullbacks and fiber products, and any such definition will depend on choices of convention.
Our goal in this section is to provide a definition such that co-orientations of fiber products of co-oriented maps possess the following desirable properties:

\begin{enumerate}
	\item Embedding property: If $f \colon V \to M$ and $g \colon W \to M$ are transverse co-oriented embeddings of manifolds without boundary into a manifold without boundary, then their fiber product is locally just the (embedding of the) intersection of the images of $V$ and $W$ in $M$ \cite[Proposition 7.3.5]{MaDo92}.
	If $f$ and $g$ are normally co-oriented (see \cref{normal co-or}), then the intersection should be normally co-oriented with the orientation of the normal bundle of the intersection given by concatenating the orientation for the normal bundle of $V$ followed by the orientation for the normal bundle of $W$.

	\item Associativity: If $V$, $W$, and $X$ are all manifolds with corners mapping to the manifold without boundary $M$ and if all the required transversality assumptions hold for the following statement to involve only well-defined fiber products, then, using Notation \ref{N: implicit notation} to allow the domain to represent also the co-oriented map, we should have
	$$(V \times_M W) \times_M X = V \times_M (W \times_M X).$$

	\item Graded commutativity: We should have
	$$V \times_M W = (-1)^{(m-v)(m-w)}W \times_M V$$
	as fiber products, using Notation \ref{N: implicit notation}.
	See \cref{R: precise commutativity}, below, for further remarks on how to interpret this formula.

	\item Leibniz rule: We should have
	$$\bd (V \times_M W) = (\bd V \times_M W)\amalg (-1)^{m-v}(V \times_M \bd W),$$
	again using Notation \ref{N: implicit notation}.
	This formula will hold for pullbacks as well as fiber products.
\end{enumerate}
Of course the latter two properties will later correspond to the graded commutativity and boundary formulas for geometric cup products.

Before getting into the specifics of the construction, we conclude this introductory section with some important general observations.
In the following subsections, we first show that pullbacks of co-orientable maps are co-orientable; the proof will then become the roadmap for defining a specific pullback co-orientation convention.
We then demonstrate that our convention yields pullback and fiber product co-orientations satisfying a number of desireable properties, including the above Leibniz rule.
Further properties of co-orientations of fiber products, including associativity and graded commutativity, will be proven in \cref{S: exterior products}, using a co-orientation we will define there for direct (exterior) products of co-oriented maps.

\begin{remark}\label{R: pullback representative}
	While our canonical pullback\index{pullback!canonical vs. categorical} $P$ has been defined as $V \times_M W = \{(x,y) \in V \times W \mid f(x) = g(y)\}$, categorically the pullback $P$ is technically only well defined up to canonical diffeomorphisms.
	In particular, if $P$ and $P'$ are two specific representatives of the pullback categorically, we have commutative diagrams
	\begin{equation*}
		\begin{tikzcd}[column sep=small]
			P \arrow[rr, "\cong"] \arrow[dr] & & P' \arrow[dl] \\
			& W. &
		\end{tikzcd}
	\end{equation*}
	But, as we have observed in \cref{D: tautological co-orientation}, diffeomorphisms come equipped with tautological co-orientations, and so a co-orientation of $P \to M$ determines\index{co-orientation!of pullback!equivalence across realizations} a unique co-orientation of $P' \to M$ by composition and vice versa.
	Thus, when working with co-orientations of pullbacks, we typically think of selecting a fixed representative $P \to W$ to work with for computations, though not always the canonical one.
	This observation shows that we are free to do so, and typically we will do so tacitly.
	This foreshadows the notion of isomorphic representatives of geometric chains and cochains; see \cref{D: equiv triv and small}.
\end{remark}

\begin{remark}\label{R: precise commutativity}\index{co-orientation!of fiber product commutativity}
	This is also a good place to point out exactly what we mean by writing $V \times_M W = (-1)^{(m-v)(m-w)}W \times_M V$ in our graded commutativity statement, as, using the canonical pullbacks, $V \times_M W \subset V \times W$ and $W \times_M V \subset W \times V$ are different spaces, though canonically identified via the map $\tau \colon V \times W \to W \times V$ that switches the coordinates.
	This map fits into a commutative diagram of fiber products
	\[
	\begin{tikzcd}[column sep=tiny]
		V \times_M W \arrow[rr, "\tau"] \arrow[dr] & & W \times_M V \arrow[dl] \\
		& M. &
	\end{tikzcd}
	\]
	Again, $\tau$ has a tautological co-orientation from being a diffeomorphism, and so the statement means that co-orientation of the fiber product $V \times_M W \to M$ and the composite co-orientation of $\tau$ and then the fiber product co-orientation of $W \times_M V \to M$ should differ by the sign $(-1)^{(m-v)(m-w)}$.
\end{remark}

\subsubsection{Co-orientability of pullbacks and fiber products}

Before defining pullback and fiber product co-orientations, we first want to ensure that pullbacks and fiber products of co-orientable maps are themselves co-orientable.
The following argument about co-orientability will also provide a roadmap to defining such co-orientations.
We also take the opportunity to observe that pullbacks of proper maps are proper.

\begin{lemma}\label{L: co-orientable pullback}
	Suppose $f \colon V \to M$ and $g \colon W \to M$ are transverse maps of manifolds with corners to a manifold without boundary.
	Then:
	\begin{enumerate}
		\item\index{co-orientability!of pullback} If $f$ is co-orientable, the pullback $f^* \colon P = V \times_M W \to W$ is co-orientable.
		\item\index{proper map!pullback of} If $f$ is proper, the pullback $f^* \colon P = V \times_M W \to W$ is proper.
	\end{enumerate}
\end{lemma}

Note that $g$ need not be co-orientable or proper for this lemma to apply.

\begin{proof}
	We first show that the pullback of a proper map is proper.
	Recall that if we use the canonical version of $V \times_M W$ as $P = \{(x,y) \mid f(x) = g(y )\}$ then the pullback map $f^*:V\times_MW\to W$ can be identified with the restriction to $V \times_M W$ of the projection $\pi_W: V\times W\to W$. Let us write $f^*=\pi_W$ so that we label our maps
	\[
	\begin{tikzcd}
		P \arrow[r, "\pi_V"] \arrow[d, "\pi_W"] & V \arrow[d, "f"] \\
		W \arrow[r, "g"] & M.
	\end{tikzcd}
	\]
	Suppose $K \subset W$ is compact.
	We have
	\begin{align*}
		\pi_W^{-1}(K)& = \{z \in P \mid \pi_W(z) \in K\}\\
		& \subset \{z \in P \mid g\pi_W(z) \in g(K)\} \\
		& = \{z \in P \mid f\pi_V(z) \in g(K)\} \\
		& = \{z \in P \mid \pi_V(z) \in f^{-1}(g(K))\}.
	\end{align*}
	So $\pi_W^{-1}(K) \subset K \times f^{-1}(g(K)) \subset V \times W$.
	But this is a product of compact sets as $f$ is proper.
	So $\pi_W$ is proper.

	For co-orientability, it suffices by \cref{P: interior co-orientation} to consider the restriction to the interior of $V \times_M W$. Equivalently, we assume that $V$ and $W$ are manifolds without boundary for the remainder of the argument.
	By \cref{L: Quillen}, it suffices to utilize Quillen's approach to co-orientability as discussed in \cref{S: Quillen}.
	We factor $f$ as $V \xhookrightarrow{e} M \times \R^N \to M$, and then we have the pullback diagram
	\begin{equation}\label{D: pullback}
		\begin{tikzcd}
			P \arrow[r, "\pi_V"] \arrow[d] & V \arrow[hookrightarrow, d, "e"] \\
			W \times \R^N \arrow[r, "g \times \id"] \arrow[d] & M \times \R^N \arrow[d, "\pi_M"] \\
			W \arrow[r, "g"] & M.
		\end{tikzcd}
	\end{equation}

	The bottom square is evidently a pullback.
	Thus by elementary topology the top square is a pullback diagram if and only if the composite rectangle is a pullback diagram.
	So by choosing $P$ so that the top square is a pullback diagram, we obtain also the pullback of $W \xr{g} M\xleftarrow{f} V$.

	Since $f$ is transverse to $g$, we have $g \times \id$ transverse to $e$ by \cref{L: all transversality is wrt embeddings}.
	As $e$ is an embedding, it follows that $P = (g \times \id)^{-1}(e(V))$ is a submanifold of $W \times \R^N$ by \cite[Proposition IV.1.4]{Kos93}.
	Furthermore, by \cref{L: Quillen}, $e(V)$ has an orientable normal bundle in $M \times \R^N$, and since the pullback of the normal bundle is the normal bundle of the pullback, again by \cite[Proposition IV.1.4]{Kos93}, it follows that the normal bundle of $P$ in $W \times \R^N$ is also orientable.
	Applying \cref{L: Quillen} again, the map $f^* = \pi_W \colon P \to W$ is co-orientable.
\end{proof}

\begin{remark}\label{R: pullback representative 2}
	As foreshadowed in \cref{R: pullback representative}, we here use a different realization of $P = V \times_M W$.
	Thinking of the top square of the diagram as a pullback square, this $P$ is concretely the subset $\{(v,(w,z)) \in V \times (W \times \R^N) \mid e(v) = (g(w),z)\}$.
	If $\pi_1 \colon M \times \R^N \to M$ and $\pi_2 \colon M \times \R^N \to \R^N$ are the projections, we know by definition that $\pi_1(e(v)) = f(v)$, and this is also $g(w)$, so the points in this realization of $P$ also satisfy $f(v) = g(w)$.
	In fact, there is a canonical diffeomorphism between this realization of $P$ and our standard realization $\{(v,w) \in V \times W \mid f(v) = g(w)\}$ given by $(v,(w,z)) \mapsto (v,w)$ with inverse given by $(v,w) \mapsto (v,(w,\pi_2(e(v))))$.

	We also already observed in the above proof that $P$ can be identified with $(g \times \id)^{-1}(e(V)) \subset W \times \R^N$, though here we need to be a bit more careful as it is not completely clear that $(g \times \id)^{-1}(e(V))$ has the structure of an embedded submanifold with corners.
	In fact, the notions of immersion and embedding are complex when working with manifolds with corners; for example, see \cite[Chapter 3]{MaDo92}.
	However, we know by \cref{L: normal pullback} that $(g|_{S^0(W)} \times \id)^{-1}(e(S^0(V)))$ will be an embedded submanifolds in $S^0(W)\times \R^N$, and due to \cref{pullback,P: interior co-orientation} this will generally suffice for working with co-orientations.
	So, unless noted otherwise, we will below assume such restrictions, or simply that we are working with manifolds without boundary in the first place, whenever we identify $P$ with $(g \times \id)^{-1}(e(V))$.
	Otherwise, we assume that $P$ has its structure as a manifold with corners as a pullback as given in \cite{Joy12}.
	Of course when we do identify $P$ with $(g \times \id)^{-1}(e(V))$, it is via the map $\{(v,(w,z)) \in V \times (W \times \R^N) \mid e(v) = (g(w),z)\} \to W \times \R^N$ given by $(v,(w,z)) \mapsto (w,z)$.
\end{remark}

\begin{corollary}\index{co-orientability!of fiber product}
	If $f \colon V \to M$ and $g \colon W \to M$ are transverse and co-orientable, their fiber product $V \times_M W \to M$ is also co-orientable.
\end{corollary}

\begin{proof}
	By the preceding lemma, the pullback $V \times_M W \to W$ is co-orientable, and the map $g \colon W \to M$ is co-orientable by assumption.
	Now choose co-orientations and compose to get a co-orientation of $P \to M$.
\end{proof}

\subsubsection{Co-orientations of pullbacks and fiber products}\label{S: co-orientation of pullbacks}

The construction in the proof of \cref{L: co-orientable pullback} provides a roadmap to define the co-orientations of pullbacks and fiber products.
For the following definition, refer again to Diagram \eqref{D: pullback}. Again thanks to \cref{P: interior co-orientation}, we can assume that our spaces are all smooth manifolds.

\begin{definition}\label{D: pullback coorient}
	Suppose $f \colon V \to M$ and $g \colon W \to M$ are transverse maps from manifolds without boundary to a manifold without boundary such that $f$ is co-oriented and the normal bundle $\nu V$ of $e(V) \subset M \times \R^N$ is given its Quillen normal orientation as defined in \cref{D: Quillen normal or}.
	Then the pullback $P = V \times_M W = (g \times \id_{\R^N})^{-1}(e(V)) \subset W \times \R^N$ has an oriented normal bundle that is the pullback of $\nu V$, which, by abuse of notation, we also label $\nu V$.
	Let $\beta_P$ and $\beta_W$ be local orientations of $P$ and $W$, and let $\beta_E$ be the standard orientation of $\R^N$.
	Define the \textbf{pullback co-orientation}\index{co-orientation!of pullback|(textbf} on $P \to W$ to be the composition of the normal co-orientation $(\beta_P,\beta_P \wedge \beta_{\nu V})$ with the canonical co-orientation $(\beta_W \wedge \beta_E,\beta_W)$.
	In other words, the pullback co-orientation is $(\beta_P,\beta_W)$ if $\beta_P \wedge \beta_{\nu V} = \beta_W \wedge \beta_E$ and $-(\beta_P,\beta_W)$ otherwise.
	Equivalently, the co-orientation of the pullback is the one for which $\beta_{\nu V}$ is the compatible Quillen normal orientation.

	If $V$ and $W$ are manifolds with corners, then we define the \textbf{pullback co-orientation} on $P \to W$ by applying the above construction to the restrictions of $f$ and $g$ to $S^0(V)$ and $S^0(W)$; this is sufficient to determine the co-orientation by \cref{pullback,P: interior co-orientation}.\index{co-orientation!of pullback|)textbf}

	Following \cref{D: top pullback}, we sometimes write $f^* \colon P \to W$.
	We also sometimes write $P = g^*V$ to emphasize that $P$ is the pullback of $V$ by $g$ to a manifold over $W$.

	If $g$ is co-oriented, define the \textbf{fiber product co-orientation}\index{co-orientation!of fiber product} on $P \to M$ as the composition of the pullback co-orientation with the co-orientation of $g \colon W \to M$.
\end{definition}

In the definition, note that the Quillen orientation of $\nu V$ is determined by the co-orientation of $f$, and the orientation $\beta_E$ of $\R^N$ is taken to be canonically fixed across all instances.
The other orientations appearing in the definition are $\beta_P$ and $\beta_W$, but the co-orientation of the pullback $f^* \colon P \to W$ does not depend on the particular choices.
For example, suppose we choose $\beta_P$ and $\beta_W$ so that $\beta_P \wedge \beta_{\nu V} = \beta_W \wedge \beta_E$ and hence the pullback co-orientation is $(\beta_P,\beta_W)$.
If we replace $\beta_P$ with $\beta_P' = -\beta_P$, then
$\beta_P' \wedge \beta_{\nu V} = -\beta_P \wedge \beta_{\nu V} = -\beta_W \wedge \beta_E$, so the co-orientation is $-(\beta_P',\beta_W) = -(-\beta_P,\beta_W) = (\beta_P,\beta_W)$.
So the co-orientation is unchanged.
Similarly, the definition is independent of our choice of $\beta_W$.
We will show just below that the definition is independent of $N$ and $e$, as well.

\begin{remark}\label{R: co-or restriction or switch}
	It follows from the definition that reversing the co-orientation of $f \colon V \to M$ reverses the co-orientation of $f^* \colon V \times_M W \to W$.
	Furthermore, if $g \colon W \to M$ is co-oriented, then reversing the co-orientation of either $f$ or $g$ reverses the co-orientation of the fiber product $f \times_M g \colon V \times_M W \to M$.

	It is also clear that the definition is consistent under restrictions to open sets.
	In other words if $x \in V$, $y \in W$ with $f(x) = g(y)$, then replacing $V$, $W$, and $M$ with neighborhoods of $x$, $y$, and $f(x) = g(y)$ yields a co-orientation of the restriction of $f^*$ to a neighborhood of $(x,y) \in V \times_M W$ that is consistent with the co-orientation of all of $f^*$, at least so long as we use the same $N$ and a restriction of $e$, though we will now show independence of these choices as well.
\end{remark}

\begin{lemma}\label{L: pullback co well defined}
	The pullback and fiber product co-orientations do not depend on the choices of $N$, $e$, or the choices of local orientations of $P$, $V$, $W$, or $M$ used in the definition (note that local orientations of $V$ and $M$ are implicit in choosing a Quillen normal orientation for $\nu V$).
\end{lemma}

\begin{remark}
	The pullback co-orientations \textit{do} depend on such choices as the choice to use the standard orientation for $\R^N$ and the choice for the co-orientation of the projection $M \times \R^N \to M$ to be $(\beta_M \wedge \beta_E, \beta_M)$, but these are universal choices that we make once and for all.
	The point is that the pullback and fiber product co-orientations only depend on $f$, $g$, and their co-orientations, after fixing such universal choices that do not depend on $f$ or $g$.
\end{remark}

\begin{proof}[Proof of \cref{L: pullback co well defined}]
	It suffices to assume that $P$, $V$, $W$, and $M$ are all manifolds without boundary.
	As the local orientations of $P$, $V$, $W$, and $M$ used in the construction all come in pairs (e.g.\ $\beta_P$ in $(\beta_P,\beta_P \wedge \beta_{\nu V}))$, the construction is independent of those choices.

	Next, suppose we are given an embedding $e \colon V \into M \times \R^N$ and extend it to $e' = (e,0) \colon V \into M \times \R^N \times \R^n$.
	In the construction involving $e$, if we choose $\beta_V$, $\beta_M$ so that $(\beta_V,\beta_M)$ is the co-orientation for $f$, then by the definition of the Quillen orientation, $\nu V$ will be such that $\beta_V \wedge \beta_{\nu V} = \beta_M \wedge \beta_E$.
	If we now increase the dimension of the Euclidean factor to $\R^{N+n}$ and write its canonical local orientation as $\beta_{E^N} \wedge \beta_{E^n}$ while extending $e$ to $e'$, we see that $\nu V$ becomes $\nu V \oplus \underline{\R}^n$ so that $\beta_{\nu V}$ becomes $\beta_{\nu V} \wedge \beta_{E^n}$.
	Pulling back over $W$ we obtain the pullback co-orientation $(\beta_P,\beta_P \wedge \beta_{\nu V} \wedge \beta_{E^n})*(\beta_W \wedge \beta_{E^N} \wedge \beta_{E^n},\beta_W)$.
	This is $(\beta_P,\beta_W)$ if and only if $\beta_P \wedge \beta_{\nu V} \wedge \beta_{E^n} = \beta_W \wedge \beta_{E^N} \wedge \beta_{E^n}$, but this condition is equivalent to having $\beta_P \wedge \beta_{\nu V} = \beta_W \wedge \beta_{E^N}$.
	So the pullback co-orientation is unchanged.

	Now suppose that $e_0 \colon V \to M \times R^{N_0}$ and $e_1 \colon V \to M \times R^{N_1}$ are any two embeddings over $f$.
	By the preceding paragraph, by adding Euclidean factors we can assume $N_0 = N_1 = N$ for some sufficiently large $N$ without changing the pullback co-orientations associated to $e_0$ and $e_1$.
	Let $\pi \colon M \times \R^N \to M$ be the projection to $M$.
	As $\pi e_0 = \pi e_1$ by assumption, the maps $e_0$ and $e_1$ are homotopic over $f$, say by linear homotopies in the Euclidean fibers.
	Let $H \colon V \times I \to M \times \R^N$ be the chosen homotopy.
	Next, by the same argument by which embeddings such as $e_0$ and $e_1$ exist (see the proof of \cref{L: Quillen}), there is an embedding $\td H \colon V \times I \into M \times \R^N \times \R^Q$ for some $Q \geq 0$ so that if $\td \pi \colon M \times \R^N \times \R^Q \to M \times \R^N$ is the projection then $\td \pi \td H = H$.
	If we let $\td e_0, \td e_1$ denote respectively $\td H(-,0), \td H(-,1) \colon V \to M \times \R^N \times \R^Q$, then $\td \pi \td e_0 = e_0$, $\td \pi \td e_1 = e_1$, and $\td H$ is a homotopy from $\td e_0$ to $\td e_1$.
	If then $(e_0,0) \colon V \to M \times \R^N \times \R^Q$ denotes the map $x \mapsto (e_0(x),0)$, then there is also a homotopy from $\td e_0$ to $(e_0,0)$; in fact as $e_0$ is an embedding and $\td \pi \td e_0 = e_0$, we can let these homotopies be linear in the $\R^Q$ factor and constant in the other factors and this homotopy will be an embedding of $V \times I$ into $M\times \R^N \times \R^Q$.
	We can define $\td e_1$, $(e_1,0)$, and an embedded homotopy between them similarly.

	So we have a sequence of three embedded homotopies, say $F_0$, $\td H$, and $F_1$, respectively from $(e_0,0)$ to $\td e_0$, from $\td e_0$ to $\td e_1$, and from $\td e_1$ to $(e_1,0)$.
	Additionally, as for $\td H$, we have $\pi \td \pi F_j(x,t) = f(x)$, so each of the three homotopies is constant in $I$ when projected down to $M$, and in particular each of $(e_0,0)$, $(e_1,0)$, $\td e_0$, and $\td e_1$ is an embedding $V \into M \times \R^{N+Q}$ over $f \colon V \to M$.
	We know from above that the pullback co-orientation obtained from using $(e_0,0)$ and $(e_1,0)$ agree with those obtained from $e_0$ and $e_1$.
	So it suffices to use these homotopies to show successively that $(e_0,0)$, $\td e_0$, $\td e_1$, and $(e_1,1)$ all provide the same pullback co-orientation of $f^*$. This is the content of the following lemma.
\end{proof}

\begin{lemma}\label{L: homotopy pullback independence}
	Suppose that $f \colon V \to M$ and $g \colon W \to M$ are transverse with $f$ co-oriented.
	Let $\pi \colon M \times \R^N \to M$ be the projection, and let $h \colon V \times I \to M \times \R^N$ be an embedding such that $\pi h(x,t) = f(x)$.
	Say that $h$ is a homotopy between the embedding $h_0$ and $h_1$.
	Then the co-orientations of $V \times_M W$ obtained from the above pullback co-orientation construction, using $h_0$ and $h_1$ as the embeddings of $V$, are the same.
\end{lemma}

\begin{proof}
	The composition $\pi h \colon V \times I \to M$ is the constant homotopy with $\pi h(x,t) = f(x)$, so $\pi h$ is also transverse to $g$, and
	\begin{align*}
		(V \times I) \times_M W
		&= \{ (x,t,y)\in V \times I \times M \mid \pi h(x,t) = g(y) \} \\
		&= \{ (x,t,y)\in V \times I \times M \mid f(x) = g(y) \} \\
		&= (V \times_M W) \times I.
	\end{align*}
	So, within the setup of the lemma, we have $V \times I$ embedded in $M \times \R^N$ and $(V \times_M W) \times I$ embedded in $W \times \R^N$ with $g \times \id_{\R^N}$ mapping the latter to the former.
	Let us fix through the proof a particular pair $(x,y)$ with $f(x) = g(y)$.
	Then identifying the spaces with their embeddings, the map $g \times \id$ takes $(x,y,t) \subset (V \times_M W) \times I$ to $(x,t) \subset V$; in fact, as $\pi h$ is a constant homotopy, these correspond to the same path in the fiber $\R^N$.
	To simplify notation, we write $x_0 = (x,0) \in V_0= V \times 0$ and $x_1 = (x,1) \in V_1 = V \times 1$ and identify these points and spaces also with their images under $h$.

	Suppose $\beta_V$ is a local orientation at $x$ and $\beta_M$ is a local orientation at $f(x)$ such that $(\beta_V,\beta_M)$ is the given co-orientation of $f$.
	We can also think of these as co-orientations of the restrictions of $\pi h$ to each $V \times t \subset V \times I$.
	We co-orient the composite $\pi h \colon V \times I \to M$ with the co-orientation $(\beta_V \wedge \beta_I, \beta_M)$ at each $(x,t)$, with $\beta_I$ corresponding to the standard positive orientation of the interval.
	So then if we give each inclusion $V \to V \times t \subset V \times I$ the normal co-orientation $(\beta_V, \beta_V \wedge \beta_I)$ at $x$, the composition of this normal co-orientation with the co-orientation of $\pi h$ is the co-orientation of $f$.

	Now, considering $V \times I$ as embedded via $h$ in $M \times I$, we have the Quillen orientation $\beta_{\nu(V \times I)}$ of the normal bundle $\nu(V \times I)$.
	If we let $\beta_{V \times I} = \beta_V \wedge \beta_I$ at some $(x,t)$, then by definition the Quillen orientation of the normal bundle is the one such that
	\[
	(\beta_V \wedge \beta_I, \beta_V \wedge \beta_I \wedge \beta_{\nu(V \times I)}) * (\beta_M\wedge \beta_E,\beta_M) =
	(\beta_V\wedge \beta_I,\beta_M),
	\]
	i.e.\ such that $\beta_V \wedge \beta_I \wedge \beta_{\nu(V \times I)} = \beta_M \wedge \beta_E$.
	Similarly, if $\nu V_0$ is the normal bundle to $V_0$ at $x_0$, then its Quillen orientation is such that
	\[
	(\beta_V, \beta_V \wedge \beta_{\nu V_0})*(\beta_M \wedge \beta_E,\beta_M) =
	(\beta_V, \beta_M),
	\]
	i.e.\ $\beta_V \wedge \beta_{\nu V_0} = \beta_M \wedge \beta_E$.
	But on $V_0$, the normal bundle $\nu V_0$ is just the sum of $\nu(V \times I)$ and a line bundle tangent to $V \times I$ in the $I$ direction, which we write as $\nu V_0 = \nu(V \times I) \oplus TI$.
	And since the above computations imply $\beta_V \wedge \beta_I \wedge \beta_{\nu(V \times I)} = \beta_V \wedge \beta_{\nu V_0}$ at $x_0$, the relation among the orientations is that $\beta_I \wedge \beta_{\nu(V \times I)} = \beta_{\nu V_0}$.
	By the same argument, we have the equivalent relation at $x_1 \in V_1$.
	The point of all this is that we have now related the orientations of the Quillen normal bundles at $x_0 \in V_0$ and $x_1 \in V_1$ to the orientations there of the single oriented bundle $\nu(V \times I)$.

	Now, following the recipe for the pullback co-orientation of $(V \times I) \times_M W \cong (V \times_M W) \times I$, we pull $V \times I \subset M \times \R^N$ back by $g \times \id$, and the oriented normal bundles of $(V \times_M W) \times I$, $V_0 \times_M W$, and $V_1 \times_M W$ will be the pullbacks of the normal bundles of $V \times I$, $V_0$, and $V_1$.
	Again, we abuse notation and use the same notations for the pulled back bundles and their orientations.
	As we have noted that the path $(x,y) \times I \in (V \times_M W) \times I$ maps to the path $x \times I \in V \times I$, preserving the orientation of $I$, the identification $\nu V_0 = \nu(V \times I) \oplus TI$ continues to hold under the pullback, and similarly at the other end of the homotopy.

	At $(x_0,y) \in (V \times_M W) \times I$, let us now choose a local orientation of the form $\beta_P \wedge \beta_I$ with $\beta_P$ a local orientation of $V \times_M W$.
	By definition, the pullback co-orientation of the projection $(V \times I) \times_M W \to W$ at this point will be $(\beta_P \wedge \beta_I, \beta_P \wedge \beta_I \wedge \beta_{\nu(V \times I)})* (\beta_W \wedge \beta_E,\beta_W)$.
	For convenience, let use suppose we choose $\beta_W$ so that $\beta_P \wedge \beta_I \wedge \beta_{\nu(V \times I)} = \beta_W \wedge \beta_E$, in which case the co-orientation of $(V \times I) \times_M W \to W$ is $(\beta_P \wedge \beta_I, \beta_W)$.
	Similarly, continuing to use the same $\beta_P$ for $V \times_M W$, the co-orientation at $(x_0,y)$ of $V_0 \times_M W \to W$ will be $(\beta_P,\beta_P \wedge \beta_{\nu V_0})* (\beta_W \wedge \beta_E,\beta_W)$.
	As the relationship $\beta_I \wedge \beta_{\nu(V \times I)} = \beta_{\nu V_0}$ is maintained under the pullback, this latter co-orientation is
	$(\beta_P,\beta_P \wedge \beta_I \wedge \beta_{\nu(V \times I)}) * (\beta_W \wedge \beta_E, \beta_W)$, which is then $(\beta_P,\beta_W)$ by our preceding assumption that $\beta_P \wedge \beta_I \wedge \beta_{\nu(V \times I)} = \beta_W \wedge \beta_E$.
	So the pullback co-orientation of $V_0 \times_M W \to W$ is just $(\beta_P, \beta_W)$.
	But now the computations concerning $V_1$ are equivalent.
	Furthermore, as the pullback map $(V \times_M W) \times I \to W$ is also a constant homotopy, the choice of $\beta_W$ so that $\beta_P \wedge \beta_I \wedge \beta_{\nu(V \times I)} = \beta_W \wedge \beta_E$ will be the same: as we transport the orientations $\beta_P$, $\beta_I$, and $\beta_{\nu(V \times I)}$ along $(x,y) \times I$ in $(V \times_M W) \times I$, the corresponding $\beta_W$ will stay constant.
	In particular, the entire equality $\beta_P \wedge \beta_I \wedge \beta_{\nu(V \times I)} = \beta_W \wedge \beta_E$ transports along the image of $(x,y) \times I$ in $(V \times_M W) \times I$, so the co-orientation of $V_1 \times_M W$ at $(x_1,y)$ is also $(\beta_P, \beta_W)$.
\end{proof}

\begin{remark}\label{R: local pullback co-orientations}
	The pullback co-orientation is determined locally in the sense that if $U$ is an open subset of $M$ then the pullback co-orientation of $f^{-1}(U) \times_U g^{-1}(U) \to g^{-1}(U)$ will just be the restriction of the pullback co-orientation of $V \times_M W \to W$.
	This is clear from the construction if we co-orient the local pullback using the Quillen co-orientation of $f^{-1}(U) \to U$ obtained from $f^{-1}(U) \xhookrightarrow{e|_{f^{-1}(U)}} U \times \R^N \to U$, the restriction of the Quillen co-orientation for $f$ obtained using the map $V\xhookrightarrow{e}M \times \R^N \times M$.
	But \cref{L: pullback co well defined} says that we are free to make such a choice.
\end{remark}

\begin{remark}\label{R: what products exist}
	We have just shown that, after choosing conventions, the fiber product of two transverse co-oriented maps is co-oriented, and this will eventually lead us to the cup product of geometric cochains.
	Analogously, if $f \colon V \to M$ is co-oriented and $W$ is oriented, then the pullback co-orientation $f^* \colon P = V \times_M W \to W$ provides a way to orient $P$, namely if $\beta_W$ is the given globally-defined orientation of $W$ then we can choose $\beta_P$ so that $(\beta_P, \beta_W)$ is the co-orientation of $f^*$ (this is just the induced orientation discussed in \cref{S: co-orientations}).
	This observation will be utilized below in our construction of the cap product.
	However, somewhat surprisingly, the fiber product of two maps with oriented domains cannot necessarily be oriented, and so there is in general no product of geometric chains and hence, in general, no homology product.
	Such oriented fiber products can be formed if the codomain $M$ is oriented, as in this case there is an equivalence between orientations of domains and co-orientations of maps.
	But this is not always possible when $M$ is not orientable.
	For example, we recall that the intersection of two orientable $\R P^3$s in the non-orientable $\R P^4$ can be a non-orientable $\R P^2$.
\end{remark}

\subsubsection{Functoriality of pullbacks}

The co-oriented pullback construction is functorial in the following sense.

\begin{proposition}\label{P: pullback functoriality}\index{co-orientation!of pullback!functoriality}
	Suppose $f \colon V \to M$ is a co-oriented map from a manifold with corners to a manifold without boundary.
	Then the pullback of $f$ by the identity $\id_M \colon M \to M$ is (diffeomorphic to) $f \colon V \to M$ with the same co-orientation.

	Suppose further that $X$ is a manifold with corners, that $W$ is a manifold without boundary, that $g \colon W \to M$ is transverse to $f$ and that $h \colon X \to W$ is transverse to $V \times_M W \to W$ (or, equivalently by \cref{L: transverse to pullback}, that $gh$ is transverse to $f$).
	Then $(gh)^*V \cong h^*g^*V$ as co-oriented manifolds over $X$.
\end{proposition}

\begin{proof}
	By \cref{pullback,P: interior co-orientation} it suffices to assume that all manifolds are without boundary.

	We first note that there is a diffeomorphism between $V$ and $V \times_M M = \{(v,x) \in V \times M \mid f(v) = x\}$ given by $v \mapsto (v,f(v))$ with inverse $(v,x) \mapsto v$.
	Then, given a compatible Quillen co-orientation of $f$, we can form the pullback diagram as
	\[
	\begin{tikzcd}[column sep=large]
		V \arrow[r, "\id_V"] \arrow[d,"e"] & V \arrow[d,"e"] \\
		M \times \R^N \arrow[r, "\id_{M \times \R^N}"] \arrow[d,"\pi_M"] & M \times \R^N \arrow[d,"\pi_M"] \\
		M \arrow[r, "\id_M"] & M,
	\end{tikzcd}
	\]
	and the conclusion is evident.

	For the second claim, there is a diffeomorphism between $V \times_M X = \{(v,x) \in V \times X \mid f(v) = g(h(x))\}$ and $(V \times_M W) \times_W X = \{((v,w),x) \in (V \times_M W) \times X \mid w = h(x)\}$ given by $(v,x) \mapsto (v,h(x),x)$ and $((v,w),x) \mapsto (v,x)$.
	To see that the last map is well defined notice that $f(v) = g(h(x))$ as $h(x) = w$, and $f(v) = g(w)$ from the assumption $(v,w) \in V \times_M W$.
	Alternatively, these two pullbacks must be diffeomorphic by general category theory, as our pullbacks are pullbacks in the category of manifolds with corners by \cite[Section 6]{Joy12}.

	Compatibility of the co-orientations now follows by considering the following diagram.
	We may assume $V$ and $W$ are manifolds without boundary by \cref{P: interior co-orientation}, and then the map labeled $e^*$ is an embedding, as recalled in the proof of \cref{L: co-orientable pullback}.
	We then note that it is equivalent to pull back the normal bundle $\nu V$ to $X \times \R^N$ either in two steps or all at once.
	\[
	\begin{tikzcd}[column sep=large]
		(V \times_M W) \times_W X \arrow[r, "\pi_{V \times_M W}"] \arrow[d] & V \times_M W \arrow[r, "\pi_V"] \arrow[d,"e^*"] & V \arrow[d, "e"] \\
		X \times \R^N \arrow[r, "h \times \id"] \arrow[d] & W \times \R^N \arrow[r, "g \times \id"] \arrow[d] & M \times \R^N \arrow[d, "\pi_M"] \\
		X \arrow[r, "h"] & W \arrow[r, "g"] & M.
	\end{tikzcd}
	\]
\end{proof}

\subsubsection{Fiber products of immersions}\label{S: co-or product immersion}

Pullbacks have particularly nice descriptions when one or both of the maps are embeddings or immersions.
In addition, these special cases are good for building intuition about the more general situation.

\begin{example}\label{E: V embedded}
	When $f \colon V \to M$ is a co-oriented embedding, the pullback co-orientation is particularly easy to describe. Again, for co-orientation purposes, we can restrict to considering the case where $V$ and $W$ are manifolds without boundary.

	We know from \cref{S: normal orientation} that in this case a co-orientation is equivalent to an orientation $\beta_{\nu V}$ of the normal bundle to $V$ in $M$.
	Then, as $f$ is already an embedding, we can take $N = 0$ in \cref{D: pullback coorient}.
	So the pullback $V \times_M W$ is just the submanifold $g^{-1}(V) \subset W$, co-oriented by $(\beta_P,\beta_W)$, where $\beta_P \wedge \beta_{\nu V} = \beta_W$, the $\nu V$ here being the pullback of the normal bundle to $g^{-1}(V)$ in $W$.
	In other words, the co-orientation of the pullback is just the normal co-orientation corresponding to the pulled back orientation of $\nu V$.

	The case where $g$ is an embedding instead also has a nice description but requires some more technology.
	We will discuss that case below in \cref{E: embedded}.
\end{example}

In the key example when both $f \colon V \to M$ and $g \colon W \to M$ are immersions, we know by \cref{L: fiber product of embeddings,S: normal orientation} that the co-orientations correspond locally to orientations of the normal bundles $\nu V$ and $\nu W$ and the fiber product $V \times_M W \to M$ corresponds locally to the intersection of the images of $V$ and $W$.
In this case our fiber product co-orientation of $V \times_M W \subset M$ is easily determined in terms of the orientations of $\nu V$ and $\nu W$.
Again, \cref{pullback,P: interior co-orientation} allow us to focus on the case where $V$ and $W$ lack boundaries.

\begin{proposition}\label{P: normal pullback}\index{co-orientation!of fiber product!of immersions}
	Let $f \colon V \to M$ and $g \colon W \to M$ be transverse co-oriented immersions from manifolds without boundary to a manifold without boundary.
	Let $\nu V$ and $\nu W$ denote the respective normal bundles.
	Choose local Quillen orientations $\beta_{\nu V}$ and $\beta_{\nu W}$ so that the normal co-orientations $(\beta_V, \beta_V \wedge \beta_{\nu V})$ and $(\beta_W, \beta_W \wedge \beta_{\nu W})$ agree with the given co-orientations of $f$ and $g$.
	Then, decomposing the normal bundle of the fiber product $P = V \times_M W \to M$ at any point as $\nu V \oplus \nu W$ and giving it the orientation $\beta_{\nu V} \wedge \beta_{\nu W}$, the fiber product co-orientation agrees with the normal co-orientation, i.e.\
	$$\omega_{f \times_M g} = (\beta_P,\beta_P \wedge \beta_{\nu V} \wedge \beta_{\nu W}).$$
\end{proposition}

That is, if one orients the normal bundle of the intersection by following an oriented basis of the normal bundle of $V$ with an oriented basis of the normal bundle of $W$, the associated normal co-orientation is the fiber product co-orientation.

\begin{proof}
	It suffices to demonstrate this property in the neighborhood of any intersection point, so we may assume that $f$ and $g$ are embeddings of manifolds without corners and consider $x \in V$, $y \in W$ with $f(x) = g(y) = z \in M$.
	Locally, for our Quillen co-orientation of $f$ we can apply the definition of the pullback co-orientation with $N = 0$ and the embedding $e \colon V \into M \times \R^N$ being simply $f$ itself.
	As $N = 0$, in this case $\nu V$ is itself the oriented normal bundle of $e(V) = f(V)$ in $M \times \R^N = M$.
	Pulling back via $g$ to $W$, we obtain the oriented pullback of $\nu V$ (which we also call $\nu V$) as the normal bundle of $P = g^{-1}(V)$ in $W$.
	By definition, the co-orientation of $P \to W$ is then the composition of $(\beta_P,\beta_P \wedge \beta_{\nu V})$ with the standard co-orientation of the projection $W \times \R^N$ to $W$, which in this case is the identity.
	The co-orientation of the fiber product is thus the composition of $(\beta_P,\beta_P \wedge \beta_{\nu V})$ with the co-orientation $(\beta_W, \beta_W \wedge \beta_{\nu W})$ of $g$.
	But this last co-orientation is independent of the choice of $\beta_W$, so we can take $\beta_W = \beta_P \wedge \beta_{\nu V}$.
	Thus we see that the fiber product co-orientation of $P \to M$ is $(\beta_P, \beta_P \wedge \beta_{\nu V} \wedge \beta_{\nu W})$, as desired.
\end{proof}

\subsubsection{The Leibniz rule}

We now verify the Leibniz rule.
We recall from \cref{D: boundary co-orientation} that if $g \colon W \to M$ is co-oriented, then the boundary co-orientation of the composite $\bd W \xr{i_{\bd W}} W \xr{g} M$ is the composite of the boundary co-orientation $(\beta_{\bd W}, \beta_{\bd W} \wedge \beta_{\nu \bd W})$, with $\beta_{\nu \bd W}$ corresponding to the inward pointing normal vector, and the co-orientation of $g$.
Recall also that we often abuse notation by letting $W$ stand for the co-oriented map $g$, and in this case we write $\bd W$ to stand for the co-oriented composite.
We also write $-W$ for $g$ with the opposite co-orientation.
This notation makes the statement of the Leibniz rule, stated just below, comprehensible.
Establishing the rule directly for immersions, for which we can use the normal co-orientations, is a quick exercise; the general case requires more care.

\begin{proposition}[Leibniz rule]\label{leibniz}\index{co-orientation!of pullback!boundary formula (Leibniz rule)}\index{co-orientation!of fiber product!boundary formula (Leibniz rule)}
	Let $f \colon V \to M$ and $g \colon W \to M$ be transverse maps from manifolds with corners to a manifold without boundary, and suppose $f$ co-oriented.
	Let $V \times_M W \to W$ be the co-oriented pullback.
	Then
	$$\bd (V \times_M W) = \left((\bd V) \times_M W\right) \bigsqcup (-1)^{m-v} \left(V \times_M (\bd W)\right),$$
	interpreting each of these pullback spaces as representing its co-oriented map to $W$; see Notation \ref{N: implicit notation}. Here we interpret $V \times_M \bd W \to W$ as the composition of the co-oriented pullback $V \times_M \bd W \to \bd W$ with the boundary immersion $\bd W \to W$ with its boundary co-orientation.

	If $g$ is also co-oriented then this formula also holds as fiber products mapping to $M$.
\end{proposition}

We first need a lemma.

\begin{lemma}\label{L: pullback boundary normal}
	Let $V$ be a smooth manifold with boundary embedded in a manifold without boundary $M$, and suppose $g \colon W \to M$ is a map from a manifold without boundary and transverse to $V$.
	In this case we know by \cref{L: immersion pullback} that $V \times_M W \to W$ is a local embedding onto $g^{-1}(V)$ and the pullback of the normal bundle $\nu V$ of $V$ in $M$ is the normal bundle of $g^{-1}(V)$ in $W$, and similarly replacing $\bd V$ with $V$.
	Let $w\in \bd (V \times_M W) = (\bd V) \times_M W$, which we identify as a subset of $W$.
	Then there is a vector $b \in T_w(V \times_M W)$, not contained in $T_w((\bd V) \times_M W)$, such that $Dg(b) \in T_{g(w)}V$ but $Dg(b) \notin T_{g(w)}(\bd V)$. Furthermore, both $b$ and $Dg(b)$ can be taken to be inward pointing, toward $V \times_M W$ and $V$, respectively.
	So, roughly speaking, there is a correspondence between the normal direction to $(\bd V) \times_M W$ in $V \times_M W$ and the normal direction to $\bd V$ in $V$ (up to the usual ambiguities in the choices of splittings for normal bundles).
\end{lemma}
\begin{proof}
	Let $v = g(w)$, and let $a \in T_v V \subset T_v M$ such that $a \notin T_v (\bd V)$.
	As $g$ is transverse to $\bd V$, there are vectors $b \in T_w W$ and $c \in T_v(\bd V)$ such that $a= Dg(b) + c$.
	Now rewriting as $Dg(b) = a - c$, the righthand side is contained in $T_v(V)$, but not in $T_v(\bd V)$ or else $a$ would be in $T_v(\bd V)$.
	As the tangent space of the pullback at $w$ is the pullback of the tangent spaces of $T_v V$ and $T_w W$ mapping to $T_v M$ by \cref{L: tangent of pullbacks}, we see that the pair $(a-c, b) \in T_v V \times_{T_v M} T_w W$ is in the tangent space of the pullback, and as the derivative of the pullback map $\pi_W \colon V \times_M W \to W$ is just projection in the $W$ co-ordinate, we see that $b$ must be a tangent vector to at $w$ mapping to $Dg(b)$ via $Dg$ as desired.
	Note that $b$ cannot be in $T_w((\bd V) \times_M W)$ as such a vector would map to $T_v (\bd V)$.
	It is clear that $Dg$ must take inward pointing vectors to inward pointing vectors.
\end{proof}

The lemma will let us identify a normal direction to $\bd V$ in $V$ with a normal direction to $(\bd V) \times_M W$ in $V \times_M W$ in the following argument.

\begin{proof}[Proof of \cref{leibniz}]
	The statement at the level of underlying manifolds with corners is \cite[Proposition 6.7]{Joy12}, so we focus on co-orientations.
	The second statement follows from the first by composing each map with the co-oriented map $g \colon W \to M$ and taking the composite co-orientations.
	We will write $\bd P$ when considering the boundary of the pullback map $P = V \times_M W \to W$ with its boundary co-orientation, and we write $(\bd V) \times_M W$ or $V \times_M (\bd W)$ when considering the maps from these boundary components with their pullback co-orientations, and in general we speak respectively of the ``boundary'' and ``pullback'' co-orientations.
	However, in both cases we write $\beta_{\bd P}$ when speaking of local orientations to simplify the notation.
	In the following arguments, it suffices to consider points in the interiors of $\bd V$ or $\bd W$, as knowing a co-orientation at one such point of each component determines it globally; in other words, we can avoid corners.
	So working locally and separating the two cases, we can assume for each case that only one of $V$ or $W$ has a boundary.

	We first consider the boundary co-orientation.
	By \cref{D: pullback coorient}, at a point of $P = V \times_M W$, the co-orientation $\omega_{f^*}$ of $P \to W$ is $(\beta_P,\beta_W)$ if and only if $\beta_P$ and $\beta_W$ are chosen so that $\beta_P \wedge \beta_{\nu V} = \beta_W \wedge \beta_E$, where $\nu V$ is the pullback to $P$ of the Quillen-oriented normal bundle of $e(V)$ in $M \times \R^N$ for some appropriate embedding $e$ that we fix throughout the following.
	The normal bundle $\nu V$ is given its Quillen orientation corresponding to the co-orientation of $f$.
	By \cref{D: boundary co-orientation}, if $\nu (\bd P)$ is an inward pointing normal of $\bd P$ in $P$ then the boundary co-orientation $\bd P \to W$ is obtained by composing the normal co-orientation of the boundary immersion $(\beta_{\bd P},\beta_{\bd P} \wedge \beta_{\nu (\bd P})$ with $\omega_{f^*}$.
	If we choose $\beta_W$ and $\beta_P$ to satisfy the condition above that $\beta_P \wedge \beta_{\nu V} = \beta_W \wedge \beta_E$ and then choose $\beta_{\bd P}$ so that $\beta_P = \beta_{\bd P} \wedge \beta_{\nu(\bd P)}$, then we have that $\omega_{\bd P \to W}$ is $(\beta_{\bd P},\beta_W)$.
	We assume in what follows that we have made such choices.

	Now, consider a point in the interior of $(\bd V) \times_M W \subset V \times_M W = P$.
	In the construction of the pullback co-orientation for $(\bd V) \times_M W$, we can take $e \colon \bd V \to M \times \R^N$ to be the restriction to $\bd V$ of our fixed embedding $e$.
	As there are two objects we could reasonable notate $\nu (\bd V)$, throughout the proof we will use $\nu (\bd V)$ for the normal bundle of $\bd V$ in $V$, and we will write $\nu^s(\bd V)$ ($s$ for stable) for the normal bundle of $e(\bd V)$ in $M \times \R^N$.
	Analogously, $\nu(\bd W)$ and $\nu(\bd P)$ will be the normal bundles as submanifolds of $W$ and $P$, respectively.
	Note that $\nu^s(\bd V) \cong \nu(\bd V) \oplus \nu V$.
	Furthermore, by employing \cref{L: pullback boundary normal}, when we pull back $\nu^s(\bd V)$ to the normal bundle of $(\bd V) \times_M W$ in $W \times \R^N$, we maintain this decomposition, identifying the pullback of $\nu V$ with the normal bundle to $V \times_M W$, as usual, and the pullback of $\nu(\bd V)$ with $\nu (\bd P)$, the normal to $\bd P$ in $P$ at our point.
	Employing our standard abuses of notation, we thus write $\beta_{\nu(\bd P)} = \beta_{\nu(\bd V)}$, identifying the inward pointing normal orientations.

	Next, recall the Quillen orientation $\beta_{\nu V}$ was chosen so that if $(\beta_V, \beta_M)$ is the co-orientation of $V$ (at an appropriate point) then $\beta_V \wedge \beta_{\nu V} = \beta_{M} \wedge \beta_E$.
	Let us fix such $\beta_V$ and $\beta_M$.
	Then if we choose $\beta_{\bd V}$ so that $\beta_V = \beta_{\bd V} \wedge \beta_{\nu(\bd V)}$ then the boundary co-orientation of $\bd V \to M$ will be $(\beta_{\bd V}, \beta_{\bd V} \wedge \beta_{\nu(\bd V)}) * (\beta_V, \beta_M) = (\beta_{\bd V},\beta_M)$, so we can then perform the pullback co-orientation construction of \cref{D: pullback coorient} using this co-orientation of $\bd V$.
	Note that as $\beta_{\nu V}$ is chosen so that $\beta_V \wedge \beta_{\nu V} = \beta_M \wedge \beta_E$, we will have also $\beta_M \wedge \beta_E = \beta_{\bd V} \wedge \beta_{\nu(\bd V)} \wedge \beta_{\nu V}$ and so the corresponding Quillen orientation of $\nu^s(\bd V)$ is $\beta_{\nu^s(\bd V)} = \beta_{\nu(\bd V)} \wedge \beta_{\nu V}$.

	So now applying \cref{D: pullback coorient}, the pullback co-orientation of $(\bd V) \times_M W \to W$ is $(\beta_{\bd P},\beta_W)$ (for our previously chosen $\beta_{\bd P}, \beta_W$) if and only if $\beta_{\bd P} \wedge \beta_{\nu^s(\bd V)} = \beta_W \wedge \beta_E$.
	But $\beta_{\bd P} \wedge \beta_{\nu^s(\bd V)} = \beta_{\bd P} \wedge \beta_{\nu(\bd V)} \wedge \beta_{\nu V}$.
	So the pullback co-orientation is $(\beta_{\bd P},\beta_W)$ if and only if $\beta_{\bd P} \wedge \beta_{\nu(\bd V)} \wedge \beta_{\nu V} = \beta_W \wedge \beta_E$.
	But we have previously identified $\beta_{\nu(\bd V)}$ with $\beta_{\nu(\bd P)}$, and $\beta_{\bd P} \wedge \beta_{\nu(\bd P)} = \beta_P$ by assumption.
	So this condition reduces to $\beta_P \wedge \beta_{\nu V} = \beta_W \wedge \beta_E$, which also holds by previous assumption.
	So the boundar and pullback co-orientations agree at points of $(\bd V) \times_M W$.

	Next, consider a point $x$ in $V \times_M \bd W$.
	As our pullback $V \times_M W$ is embedded neatly in $W \times \R^N$ \cite[Proposition IV.1.4]{Kos93}, it is immediate that, at $x$, the normal bundle of $V \times_M \bd W = (g \times \id)|_{(\bd W) \times \R^N}^{-1}(V)$ in $W \times \R^N$ can be decomposed into the direct sum of the pullback of $\nu V$, which can be identified with a subbundle of $T_x((\bd W) \times \R^N)$, and a 1-dimensional summand that is normal to $(\bd W) \times \R^N$. By projection, we can identify this summand with a normal direction to $\bd W$ in $W$ and again write $\beta_{\nu(\bd P)}=\beta_{\nu(\bd W)}$ for the orientation determined by these inward pointing normals.

	Now, we continue to assume that $(\beta_P, \beta_W)$ is the co-orientation of $V \times_M W \to W$ (and so $\beta_P \wedge \beta_{\nu V} = \beta_W \wedge \beta_E$).
	We choose $\beta_M$ so that $(\beta_W,\beta_M)$ is the co-orientation of $g$ and $\beta_{\bd W}$ so that $(\beta_{\bd W},\beta_W)$ co-orients $i_{\bd W}$, which implies $\beta_W = \beta_{\bd W} \wedge \beta_{\nu(\bd W)}$.
	We also continue to choose $\beta_{\bd P}$ so that $\beta_P = \beta_{\bd P} \wedge \beta_{\nu(\bd P)} = \beta_{\bd P} \wedge \beta_{\nu(\bd W)}$ and hence the boundary co-orientation of $\bd P \to W$ is $(\beta_{\bd P}, \beta_W)$.

	Using these local orientations and applying \cref{D: pullback coorient} to $\bd g \colon \bd W \to M$, the co-orientation of the pullback $V \times_M \bd W \to \bd W$ is $(\beta_{\bd P},\beta_{\bd W})$ (and so the composite co-orientation to $W$ is $(\beta_{\bd P}, \beta_W)$) if and only if $\beta_{\bd P} \wedge \beta_{\nu V} = \beta_{\bd W} \wedge \beta_E$ as local orientations at the image of $x$ in $\bd W \times \R^N$.
	Considering $\bd W \times \R^N \subset W \times \R^N$, this condition holds if and only if $\beta_{\bd P} \wedge \beta_{\nu V} \wedge \beta_{\nu(\bd W)} = \beta_{\bd W} \wedge \beta_E \wedge \beta_{\nu(\bd W)}$ in $W \times \R^N$.
	But as the dimensions as bundles are $\dim(\nu(\bd W)) = 1$ and $\dim(\nu V) = m+N-v$,
	$$\beta_{\bd P} \wedge \beta_{\nu V} \wedge \beta_{\nu(\bd W)} = (-1)^{m+N-v}\beta_{\bd P} \wedge \beta_{\nu(\bd W)} \wedge \beta_{\nu V} = (-1)^{m+N-v}\beta_{\bd P} \wedge \beta_{\nu(\bd P)} \wedge \beta_{\nu V} = (-1)^{m+N-v}\beta_P \wedge \beta_{\nu V} ,$$
	and
	$$\beta_{\bd W} \wedge \beta_E \wedge \beta_{\nu(\bd W)} = (-1)^N\beta_{\bd W} \wedge \beta_{\nu(\bd W)} \wedge \beta_E = (-1)^N\beta_{W} \wedge \beta_E.$$
	So the pullback co-orientation of $V \times_M \bd W \to W$ is $(\beta_{\bd P},\beta_W)$ if and only if $\beta_P \wedge \beta_{\nu V} = (-1)^{m-v}\beta_{W} \wedge \beta_E$.
	But in defining the boundary co-orientation we assumed that $\beta_P \wedge \beta_{\nu V} = \beta_{W} \wedge \beta_E$, so two co-orientations agree or disagree according to the sign $(-1)^{m-v}$ as claimed.
\end{proof}

\subsubsection{Codimension $0$ and $1$ pullbacks}\label{S: codim 0 and 1 co-or}

The results in this section should be compared with, and justify, the discussion and choices in \cref{E: splitting example 1,E: manifold decomposition}.
They will be useful when working with the splitting and creasing constructions for geometric cochains in \cref{S: splitting and creasing}.

\begin{proposition}\label{P: codim 0 pullback}\index{co-orientation!of pullback!of codimension $0$ embedding}\index{co-orientation!of fiber product!of codimension $0$ embedding}\index{pullback!of codimension $0$ embedding}
	Let $V$ be an embedded codimension $0$ submanifold with corners in the manifold without boundary $M$, and let $f \colon V \to M$ be the embedding, co-oriented by the tautological co-orientation.
	Let $W$ be a manifold with corners and suppose $g \colon W \to M$ is transverse to $f$.
	Then the co-oriented pullback $V \times_M W \to W$ is the inclusion of the codimension $0$ manifold with corners $g^{-1}(V) \into W$, co-oriented by the tautological co-orientation.
	Consequently, if $g$ is co-oriented, the co-oriented fiber product $V \times_M W \to M$ is just the restriction of the co-oriented map $g$ to $g^{-1}(V)$.
\end{proposition}

\begin{proof}
	It is clear topologically that the pullback is $g^{-1}(V)$, and it must be a codimension $0$ manifold with corners in $W$ by Joyce \cite[Theorem 6.4]{Joy12}.
	So we consider co-orientations, for which we can assume $V$ and $W$ are without boundary by \cref{pullback,P: interior co-orientation}.
	As $f$ is an embedding we may choose $N = 0$ and $e = f$ in \cref{D: pullback coorient}.
	As $V \to M$ is tautologically co-oriented, we can identify $\beta_V$ with $\beta_M$ at any point via the embedding.
	The Quillen normal bundle of $V$ in $M$ is then the positively oriented $\R^0$-bundle.
	So then the definition says that the pullback co-orientation at any point is $(\beta_P, \beta_W)$ when $\beta_P = \beta_W$, i.e.\ the pullback co-orientation is $(\beta_P, \beta_P)$, the tautological co-orientation.
\end{proof}

\begin{corollary}\label{C: cup with identity}\index{fiber product!with identity map}
	Let $f \colon V \to M$ be a co-oriented map from a manifold with corners to a manifold without boundary, and let $\id_M \colon M \to M$ be the identity with its tautological co-orientation.
	Then both co-oriented fiber products $V \times_M M \to M$ and $M \times_M V \to M$ are again $f \colon V \to M$ with the given co-orientation.
\end{corollary}

\begin{proof}
	The case of $M \times_M V \to M$ follows from the preceding lemma, and the other is the first statement of \cref{P: pullback functoriality}.
\end{proof}

\begin{example}\label{E: codim 0 and 1 co-or as fiber products}
	Let $\phi \colon M \to \R$ be a smooth function from a manifold without boundary to $\R$ that is transverse to $0$.
	By \cref{E: manifold decomposition}, we know (ignoring co-orientations for the moment) that $[0,\infty) \times_\R M = \phi^{-1}([0,\infty))$, that $(-\infty,0] \times_\R M = \phi^{-1}((-\infty,0])$, and that $0 \times_\R M = \phi^{-1}(0)$.
	As in \cref{E: manifold decomposition}, we denote these respectively as $M^+$, $M^-$, and $M^0$\index{$M^0$}\index{$M^+$}\index{$M^-$}.
	By \cite[Proposition 4.2.9]{MaDo92}, the inclusion of $M^0$ into $M$ is a closed embedding and by \cite[Proposition 6.7]{Joy12}, we have $\bd M^\pm = M^0$ as spaces.

	Suppose the inclusions $M^\pm \into M$ are given the tautological co-orientations of \cref{D: tautological co-orientation}, and let $g \colon W \to M$ be transverse to $M^\pm$, which in this case is equivalent to being transverse to $M^0$, which is also equivalent to $\phi g$ being transverse to $0$ in $\R$ by \cref{L: transverse to pullback}.
	Then let $W^\pm = M^\pm \times_M W$.\index{$W^+$}\index{$W^-$}
	By \cref{P: codim 0 pullback}, the pullback map $W^\pm \to W$ is just the inclusion $g^{-1}(M^\pm) \into W$ with its tautological co-orientation.
	If $g$ is co-oriented, the compositions $W^\pm \into W \to M$ are then the fiber products $M^\pm \times_M W \to M$, and by \cref{P: codim 0 pullback} their co-orientations are just the restrictions of the co-orientation of $g$ to $W^\pm$.
	This agrees with the co-orientations discussed in \cref{E: manifold decomposition}.
\end{example}

\begin{proposition}\label{P: codim 1 co-orient}
	Suppose $V \subset M$ is a closed codimension $1$ submanifold without boundary in the manifold without boundary $M$.
	Further suppose $V$ has oriented normal bundle $\nu$.
	Let the embedding $f \colon V \into M$ be co-oriented by the normal co-orientation $(\beta_V, \beta_V \wedge \beta_\nu)$.
	Let $W$ be a manifold with corners and suppose $g \colon W \to M$ is transverse to $f$.
	Then $W^0 \defeq g^{-1}(V)$\index{pullback!of codimension $1$ submanifold}\index{$W^0$} is a codimension $1$ submanifold with corners of $W$ with oriented pullback normal vector bundle $\nu_0$, and the co-oriented pullback $V \times_M W \to W$ is the embedding $W^0 \into W$, co-oriented by $(\beta_P, \beta_P \wedge \beta_{\nu_0})$.
\end{proposition}

\begin{proof}
	Since the normal bundle of $V$ is $1$-dimensional and oriented, it is the trivial line bundle.
	Embedding the normal bundle as a tubular neighborhood, we can then construct a map $\phi: M \to \R$ so that $V = \phi^{-1}(0)$.
	Then $g$ being transverse to $V$ is equivalent to $\phi g$ being transverse to $0$, so by \cite[Proposition 4.2.9]{MaDo92} the pullback $g^{-1}(V) = W^0$ is a codimension $1$ submanifold with corners of $W$.
	For the co-orientation, if we take $N = 0$ and $e = f$ in \cref{D: pullback coorient}, then $\nu$ is just our normal bundle $\nu V$ and $\nu_0$ is simply the pullback.
	Then, by definition, the pullback co-orientation is $(\beta_P, \beta_W)$ (at interior points) if and only if $\beta_P \wedge \beta_{\nu_0} = \beta_W$, as claimed.
\end{proof}

\begin{example}\label{E: codim 1 pullbacks}
	We continue with the assumptions and notation of \cref{E: codim 0 and 1 co-or as fiber products}, but now let $V = M^0 = \phi^{-1}(0)$ with normal bundle oriented by the pullback of the standard (positive-direction) orientation of the normal bundle of $0 \in \R$.
	This determines a normal co-orientation of the embedding $M^0 \to M$.
	Then, by \cref{P: codim 1 co-orient}, the pullback co-orientation of $W^0 = M^0 \times_M W \into W$ agrees with the $\phi$-induced co-orientation of $W^0$ defined in \cref{E: manifold decomposition}.
	We can also confirm now, using the Leibniz rule and that the codimension of $M^-$ in $M$ is $0$, that as spaces with co-oriented maps to $W$ we have
	\begin{multline*}
		\bd(W^-) = \bd(M^- \times_M W) = \left( (\bd (M^-)) \times_M W \right) \bigsqcup \left( M^- \times_M \bd W \right)\\
		= \left( -(M^0) \times_M W \right) \bigsqcup \left( M^- \times_M \bd W \right)= -(W^0) \bigsqcup (\bd W)^-.
	\end{multline*}
	Here we also use that the orientation of the normal bundle to $M^0$ is outward pointing from $M^-$ and so disagrees with the inward-pointing normal used to co-orient the boundary inclusion; hence $\bd(M^-) = -(M^0)$ as spaces with co-oriented maps to $M$.
	We also note that, by \cref{E: codim 0 and 1 co-or as fiber products}, $(\bd W)^- \to W$ is co-oriented by the tautological co-orientation $(\bd W)^- \to \bd W$ followed by the boundary co-orientation of $i_{\bd W} \colon \bd W \to W$.

	Analogously,
	$$\bd (W^+) = W^0 \bigsqcup (\bd W)^+,$$
	using that the orientation of the normal bundle to $M^0$ is inward pointing for $M^+$.
\end{example}

We can now prove the claim from the end of \cref{E: manifold decomposition}.
We express the following corollary using Notation \ref{N: implicit notation}.

\begin{corollary}\label{C: co-orient W0}
	Suppose the hypotheses and notation of \cref{P: codim 1 co-orient} and suppose $V$ is without boundary.
	Then $(\bd W)^0 = -\bd (W^0)$\index{$bd$@$\bd (W^0)$ vs. $(\bd W)^0$} as co-oriented maps to $W$, with $W^0$ and $(\bd W)^0 = (gi_{\bd W})^{-1}(V)$ co-oriented as in \cref{P: codim 1 co-orient} as the pullbacks $V \times_M W \to W$ and $V \times_M \bd W \to W$.
\end{corollary}

\begin{proof}
	By \cref{leibniz} and the preceding examples, we have that
	$$\bd (W^0) = \bd (V \times_M W) = (-1)^{m-v} V \times_M \bd W = -V \times_M \bd W = -(\bd W)^0$$
	as spaces mapping to $W$.
\end{proof}

\subsection{Exterior products and their relations with fiber products}\label{S: exterior products}

In this section, we consider products of maps that will eventually become the exterior products in geometric homology and cohomology, as well as their relations to fiber products.
While fiber products, which will eventually be used to define cup and intersection products, require special transversality conditions in order to be defined, exterior products are always fully defined.
Of course products of oriented manifolds are familiar objects, so we treat them only briefly in the next section.
Then we consider products of co-oriented maps of manifolds.
In \cref{S: product relations}, we show that the co-oriented fiber product is the pullback by the diagonal map of the co-oriented exterior product, foreshadowing the classical cohomology relation between exterior products and cup products.
This will also allow us to prove associativity of fiber products.

\subsubsection{Exterior products of oriented manifolds}\label{S: oriented product}

Recall that we defined the oriented fiber product of oriented manifolds with corners in \cref{S: orientation of fiber products}.\index{exterior product!orientation}\index{orientation!of external product}
In particular, if $V$ and $W$ are oriented manifolds, we saw in \cref{P: oriented fiber product basic properties} that the oriented fiber product of the maps from $V$ and $W$ to a point is just the standard product $V \times W$ oriented with the usual concatenation convention.
In other words, if $V$ and $W$ are oriented by $\beta_V$ and $\beta_W$, then $V \times W$ is oriented at any point by $\beta_V \wedge \beta_W$.

We observe the following interplay between fibered and exterior products of maps of oriented manifolds more generally.
In our notational shorthand (see \cref{N: implicit notation}), if $f \colon V \to M$ and $h \colon X \to N$, then we let $V \times X$ represent the product map $f \times h \colon M \times N$.

\begin{proposition}\label{P: oriented interchange}\index{orientation!of fiber product!of external product}\index{orientation!of external product!of fiber product}\index{criss cross}
	Suppose $f \colon V \to M$ and $g \colon W \to M$ are transverse maps of oriented manifolds with corners to an oriented manifold without boundary and similarly for $h \colon X \to N$ and $k \colon Y \to N$.
	Then
	$$(V \times X)\times_{M \times N} (W \times Y) = (-1)^{(m-w)(n-x)}(V \times_M W) \times (X \times_N Y)$$
	as oriented manifolds.
\end{proposition}

\begin{proof}
	We first note that the transversality assumptions ensure also that $f \times h$ will be transverse to $g \times k$.
	It is straightforward to verify that these are diffeomorphic spaces, so we focus on the orientations.
	For simplicity, let us write
	$P = (V \times X)\times_{M \times N} (W \times Y)$ and $P' = (V \times_M W) \times (X \times_N Y)$ as oriented manifolds.
	We then write local orientations symbolically as $\beta_P$, etc.
	By the construction of fiber product orientations in \cref{S: orientation of fiber products}, and omitting the pullbacks from the notation, $P$ is oriented so that
	$$\beta_P \wedge \beta_{M \times N} = (-1)^{(w+y)(m+n)}\beta_{V \times X} \wedge \beta_{W \times Y},$$
	or, as (non-fiber) products are oriented by concatenation, we have
	$$\beta_P \wedge \beta_M\wedge\beta_N = (-1)^{(w+y)(m+n)}\beta_V \wedge \beta_X \wedge \beta_W \wedge \beta_Y.$$
	Recall that here we identify $T(M \times N)$ as a summand of $T(V \times X \times W \times Y)$ over $P$ by splitting the derivative $D(f \times h) - D(g \times k)$.
	Similarly, for $V \times_M W$ and $X \times_N Y$ we have
	\begin{align*}
		\beta_{V \times_M W} \wedge \beta_M & = (-1)^{wm}\beta_V \wedge \beta_W \\
		\beta_{X \times_N Y} \wedge \beta_N & = (-1)^{yn}\beta_X \wedge \beta_Y,
	\end{align*}
	using the splittings of $Df-Dg$ and $Dh-Dk$.
	We note that the signs of $Df$, $Dg$, $Dh$, and $Dk$ in all the splitting formulas are consistent in computing the orientations for $P$ and $P'$.

	As $\beta_{P'} = \beta_{V \times_M W} \wedge \beta_{X \times_N Y}$, we have
	\begin{align*}
		\beta_{P'} \wedge \beta_M \wedge \beta_N
		&= \beta_{V \times_M W} \wedge \beta_{X \times_N Y} \wedge \beta_M \wedge \beta_N\\
		& = (-1)^{m(x+y-n)}\beta_{V \times_M W} \wedge \beta_M \wedge \beta_{X \times_N Y} \wedge \beta_N\\
		& = (-1)^{m(x+y-n)+wm+ny}\beta_V \wedge \beta_W \wedge \beta_X \wedge \beta_Y \\
		& = (-1)^{m(x+y-n)+wm+ny+xw}\beta_V \wedge \beta_X \wedge \beta_W \wedge \beta_Y.
	\end{align*}
	So $\beta_P$ differs from $\beta_{P'}$ by $-1$ to the power
	$$m(x+y-n) + wm + ny + xw - (w+y)(m+n).$$
	An elementary computation now shows that this is $(m-w)(n-x)$ as desired.
\end{proof}

\subsubsection{Exterior products of co-oriented maps}

Next we define and study a co-oriented exterior product for co-oriented maps.
In the next subsection, we will see that such products are intimately related to fiber products, and this will allow us to easily prove some properties about fiber products that we have delayed.

\begin{lemma}\index{co-orientability!of external product}
	If $f \colon V \to M$ and $g \colon W \to N$ are co-orientable maps of manifolds with corners then the product map $f \times g \colon V \times W \to M \times N$ is co-orientable.
\end{lemma}

\begin{proof}
	We recall that, by definition, a co-orientation of $f$ is equivalent to a choice of isomorphism between the orientation cover $\Or(TV)$ and the pullback $f^*\Or(TM)$ of the orientation cover $\Or(TM)$, and similarly for $g$.

	If we let $\pi_V, \pi_W$ denote the projections of $V \times W$ to $V$ and $W$, then $T(V \times W) \cong \pi_V^*(TV) \oplus \pi_W^*(TW)$, and so $$\Or(T(V \times W)) \cong \Or(\pi_V^*(TV))\otimes\Or(\pi_W^*(TW)) \cong \pi_V^*\Or(TV)\otimes\pi_W^*\Or(TW).$$ Similarly
	\begin{multline*}(f \times g)^*T(M \times N) \cong (f \times g)^*(\pi_M^*(TM) \oplus \pi_N^*(TN))\\
		\cong (f \times g)^*\pi_M^*(TM) \oplus (f \times g)^*\pi_N^*(TN)) \cong \pi_V^*f^*(TM) \oplus \pi_W^*g^*(TN),
	\end{multline*}
	using that $\pi_M(f \times g) = f\pi_V \colon V \times W \to M$ and $\pi_N(f \times g) = g\pi_W \colon V \times W \to N$.
	So
	\begin{multline*}
		(f \times g)^*\Or(T(M \times N)) \cong \Or((f \times g)^*T(M \times N))\\ \cong \Or(\pi_V^*f^*TM) \otimes \Or(\pi_W^*g^*TN) \cong \pi_V^*f^*\Or(TM) \otimes \pi_W^*g^*\Or(TN).
	\end{multline*}
	Thus if $\Or(TV) \cong f^*\Or(TM)$ and $\Or(TW) \cong g^*\Or(TN)$, we can construct an isomorphism $\Or(T(V \times W)) \cong (f \times g)^*\Or(T(M \times N))$.
\end{proof}

\begin{definition}\label{D: co-oriented exterior}
	If $f \colon V \to M$ and $g \colon W \to N$ are co-oriented maps of manifolds with corners with co-orientations given by isomorphisms $\phi \colon \Or(TV) \to f^*\Or(TM)$ and $\psi \colon \Or(TW) \to g^*\Or(TN)$, we define the \textbf{product co-orientation}\index{co-orientation!of external product|textbf} of $f \times g \colon V \times W \to M \times N$ by the isomorphism $(-1)^{(m-v)w}\pi_V^*\phi \otimes \pi_W^*\psi$.
	In particular, if at $x \in V$ the co-orientation of $f$ is given locally by $(\beta_V,\beta_M)$ and at $y \in W$ the co-orientation of $g$ is given locally by $(\beta_W,\beta_N)$, then the product co-orientation is locally represented at $(x,y)$ by $$(-1)^{(m-v)w}(\beta_V \wedge \beta_W,\beta_M \wedge \beta_N).$$

	Following our standard convention from \cref{N: implicit notation}, we often write simply $V \times W$ to represent the co-oriented product.
\end{definition}

\begin{remark}
	The sign in the definition is not at first obvious, though it will be justified in the following lemmas.
	One way to think of it is as follows: If we we take $V$ and $W$ as immersed submanifolds co-oriented by orienting their normal bundles as in \cref{normal co-or}, then $V \times W$ is also immersed, and at an image point we have $T(M \times N) \cong TV \oplus \nu V \oplus TW \oplus \nu W$, letting $\nu V$ and $\nu W$ stand for the normal bundles of $V$ and $W$ in $M$ and $N$, respectively.
	The sign $(-1)^{(m-v)w}$ is the sign needed in the local orientation to permute this to $TV \oplus TW \oplus \nu V \oplus \nu W$ so that we can properly utilize the normal co-orientation for $\nu(V \times W) \cong \nu V \oplus \nu W$.
	While this argument is essentially heuristic, it is borne out in the computations below.
\end{remark}

The following example will be useful in the proof of \cref{T: intersection is cup product}.

\begin{example}\label{E: sphere product}
	Let $S^p$ and $S^q$ be oriented spheres with $p,q>0$.
	Let $V = W = pt$, and let $f \colon V \to S^p$ and $g \colon W \to S^q$ be embeddings to points $x \in S^p$, $y \in S^q$.
	Let $f$ be co-oriented by $(1,\beta_{S^p})$; in other words $V$ is normally co-oriented by the orientation of its normal bundle that agrees with the orientation of $S^p$.
	Let $g$ be co-oriented similarly.
	Then $V \times W$ is represented by the embedding of the point to $(x,y) \in S^p \times S^q$ with normal bundle oriented consistently with the product orientation of $S^p \times S^q$.
	There is no extra sign in this case as $\dim(W) = 0$.
\end{example}

\begin{proposition}\label{P: co-oriented exterior unit}\index{product with a point}
	Let $f \colon V \to M$ be a co-oriented map of manifolds with corners, and let $g:pt \to pt$ be the unique map with the tautological co-orientation.
	Then $f \times g \colon V \times pt \to M \times pt$ and $g \times f:pt \times V \to pt \times M$ are each isomorphic as co-oriented maps of manifolds with corners to $f \colon V \to M$.
\end{proposition}

\begin{proof}
	This is obvious ignoring co-orientations.
	Considering co-orientations, if $f$ is co-oriented at a point by $(\beta_V,\beta_M)$, then the co-orientation of $f \times g$ is simply $(\beta_V \wedge 1,\beta_M \wedge 1) = (\beta_V,\beta_M)$, noting that
	the sign $(-1)^{(m-v)\cdot 0} = 1$ in this case.
	The case $g \times f$ is similar, though due to the transposition the sign is now $(-1)^{(0-0)v}$, which is again $1$.
\end{proof}

\begin{proposition}\label{P: boundary of exterior product}\index{co-orientation!of external product!boundary (Leibniz formula)}
	Let $f \colon V \to M$ and $g \colon W \to N$ be co-oriented maps of manifolds with corners and suppose $f \times g \colon V \times W \to M \times N$ is given the product co-orientation.
	Then the boundary co-orientation of $V \times W$ as co-oriented maps to $M \times N$ is $$\bd(V \times W) = (\bd V) \times W \bigsqcup (-1)^{m-v}V \times \bd W.$$
\end{proposition}

\begin{proof}
	We know that this expression is an identity ignoring co-orientations, so we must establish the agreement of the co-orientations for each component.
	As usual, it suffices to consider points in the top dimensional strata of $\bd(V \times W)$.
	In what follows, we fix $\beta_V$, $\beta_W$, $\beta_M$, and $\beta_N$, so that $(\beta_V,\beta_M)$, $(\beta_W,\beta_N)$, and $(-1)^{(m-v)w}(\beta_V \wedge \beta_W,\beta_M \wedge \beta_N)$ denote the co-orientations of $V$, $W$, and $V \times W$ at the point under consideration.

	Let $\nu$ denote an inward pointing normal to $V \times W$ at such a point.
	Then the inclusion $\bd(V \times W) \to V \times W$ is co-oriented at that point by $(\beta_{\bd(V \times W)},\beta_{\bd(V \times W)} \wedge \beta_\nu)$ for any $\beta_{\bd(V \times W)}$.
	If we choose $\beta_{\bd(V \times W)}$ so that $(\beta_{\bd(V \times W)}\wedge\beta_\nu,\beta_M \wedge \beta_N)$ represents the co-orientation of $V \times W \to M \times N$, then from the definition of the boundary co-orientation, the boundary $\bd(V \times W) \to M \times N$ is co-oriented by $(\beta_{\bd(V \times W)},\beta_M \wedge \beta_N)$.
	We fix such a choice in what follows.

	Now suppose our point is more specifically in the top-dimensional stratum of $(\bd V) \times W$.
	If we choose $\beta_{\bd V}$ so that $\beta_{\bd V} \wedge \beta_\nu = \beta_V$, then $\bd V \to M$ is co-oriented by $(\beta_{\bd V},\beta_M)$ and so $(\bd V) \times W \to M \times N$ is co-oriented by $(-1)^{(m-v+1)w}(\beta_{\bd V} \wedge \beta_W,\beta_M \wedge \beta_N)$.
	On the other hand, the co-orientation of $V \times W$ can then be written $(-1)^{(m-v)w}(\beta_{\bd V} \wedge \beta_\nu \wedge \beta_W,\beta_M \wedge \beta_N) = (-1)^{(m-v)w+w}(\beta_{\bd V} \wedge \beta_W \wedge \beta_\nu,\beta_M \wedge \beta_N)$, so the boundary co-orientation of $\bd(V \times W)$ is $$(\beta_{\bd V} \wedge \beta_W,\beta_{\bd V} \wedge \beta_W \wedge \beta_\nu)*(-1)^{(m-v)w+w}(\beta_{\bd V} \wedge \beta_W \wedge \beta_\nu,\beta_M \wedge \beta_N) = (-1)^{(m-v)w+w}(\beta_{\bd V} \wedge \beta_W,\beta_M \wedge \beta_N),$$
	which agrees with our co-orientation for $(\bd V) \times W$.

	Next consider a point in the top-dimensional stratum of $V \times \bd W$.
	If we choose $\beta_{\bd W}$ so that $\beta_{\bd W} \wedge \beta_\nu = \beta_W$ then we have $\bd W$ co-oriented by $(\beta_{\bd W},\beta_N)$ and so $V \times \bd W$ is co-oriented by $(-1)^{(m-v)(w-1)}(\beta_{V} \wedge \beta_{\bd W},\beta_M \wedge \beta_N)$.
	On the other hand, the co-orientation of $V \times W$ can now be written $(-1)^{(m-v)w}(\beta_{V} \wedge \beta_{\bd W} \wedge \beta_\nu,\beta_M \wedge \beta_N)$, so the boundary co-orientation of $\bd(V \times W)$ is $(-1)^{(m-v)w}(\beta_{V} \wedge \beta_{\bd W},\beta_{M} \wedge \beta_N)$, which differs from that of $V \times \bd W$ by a factor of $(-1)^{m-v}$.
\end{proof}

\begin{proposition}\label{P: exterior associativity}\index{co-orientation!of external product!associativity}
	Let $f \colon V \to M$, $g \colon W \to N$, and $h \colon X \to Q$ be co-oriented maps of manifolds with corners.
	Then the co-orientations of $(V \times W) \times X \to M \times N \times Q$ and $V \times (W \times X) \to M \times N \times Q$ agree.
	In other words, forming co-oriented products is associative.
\end{proposition}

\begin{proof}
	If $f,g,h$ are co-oriented by $(\beta_V,\beta_M)$, etc., then both products are co-oriented up to sign by $(\beta_V \wedge \beta_W \wedge \beta_X,\beta_M \wedge \beta_N \wedge \beta_P)$.
	In forming $(V \times W) \times X$ we first have the sign $(-1)^{(m-v)w}$ from $V \times W$, then taking the product with $X$ on the right multiplies by $(-1)^{(m+n-v-w)x}$.
	So the total sign is $(-1)^{(m-v)w+(m+n-v-w)x}$.
	Alternatively, forming $W \times X$ has the sign $(-1)^{(n-w)x}$ and then taking the product with $V$ on the left contributes $(-1)^{(m-v)(w+x)}$.
	So the total sign is $(-1)^{(n-w)x+(m-v)(w+x)}$.
	One readily verifies that these signs agree.
\end{proof}

The following lemma provides a nice description of the Quillen co-orientation of a product of co-oriented maps.
Among other things, it will help us to next demonstrate a commutativity property for exterior products of co-oriented maps.
We assume for convenience that our Euclidean factors are even dimensional, which simplifies the computations and will be sufficient for what follows; of course if $V$ embeds in $M \times \R^n$, it also embeds in $M \times \R^{n+1}$, and we have shown that our fiber product co-orientations do not depend on such choices.

\begin{lemma}\label{L: Quillen product co-orientation}\index{co-orientation!of external product!Quillen co-orientation}
	Let $f \colon V \to M$ and $g \colon W \to N$ be co-oriented maps from manifolds with corners to manifolds without boundary.
	Consider Quillen co-orientations representing $f$ and $g$ via embeddings $e_V \colon V \into M \times \R^a$ and $e_W \colon W \into N \times \R^b$ with $a$ and $b$ even.
	Denote the normal bundles of $V$ and $W$ in $M \times \R^a$ and $N \times \R^b$ by $\nu V$ and $\nu W$.
	Let $T \colon M \times \R^a \times N \times \R^b \to M \times N \times \R^{a+b}$ be the diffeomorphism that interchanges the middle two factors.
	Then
	$T(e_V \times e_W)$ gives an embedding $V \times W \to M \times N \times \R^{a+b}$ with normal bundle isomorphic to the sum of the pullbacks of $\nu V$ and $\nu W$ by the projections of $V \times W \to M \times N \times \R^a \times \R^b$ to either the first and third factor or the second and fourth factors.
	For simplicity, we simply write $\nu V \oplus \nu W$.

	Then the Quillen normal orientation of the normal bundle of $f \times g \colon V \times W \to M \times N$ is given by $$\beta_{\nu V \oplus \nu W} = \beta_{\nu V} \wedge \beta_{\nu W},$$
	suitably interpreting the relevant coordinates in $M \times N \times \R^a \times \R^b$.
\end{lemma}

The last line simply means that if an element of $\nu V$ has coordinates $(x,y)$ in $M \times \R^a$, then this corresponds to a normal vector with coordinates $(x,0,y,0)$ in $M \times N \times \R^a \times \R^b$ for which we do not create a new notation, and similarly for $\nu_W$.

\begin{proof}
	Let $\beta_a$ and $\beta_b$ denote the standard orientations for $\R^a$ and $\R^b$.
	By definition, $\nu V$ and $\nu W$ are oriented so that $\beta_V \wedge \beta_{\nu V} = \beta_M \wedge \beta_a$ and $\beta_W \wedge \beta_{\nu W} = \beta_N \wedge \beta_b$.

	By definition, the Quillen orientation of $\nu V \oplus \nu W$ corresponding to the product co-orientation of $V \times W$ is the local orientation $\beta_{\nu V \oplus \nu W}$ such that
	$$(\beta_{V \times W}, \beta_{V \times W} \wedge \beta_{\nu V \oplus \nu W})*(\beta_{M \times N} \wedge \beta_{a+b},\beta_{M \times N}) = (-1)^{(m-v)w}(\beta_V \wedge \beta_W,\beta_M \wedge \beta_N).$$
	Taking $\beta_{M \times N} = \beta_M \wedge \beta_N$ and $\beta_{V \times W} = \beta_V \wedge \beta_W$ and noting $\beta_{a+b} = \beta_a \wedge \beta_b$, this formula becomes
	$$(\beta_V \wedge \beta_W, \beta_V \wedge \beta_W \wedge \beta_{\nu V \oplus \nu W})*(\beta_M \wedge \beta_N \wedge \beta_a \wedge \beta_b,\beta_M \wedge \beta_N) = (-1)^{(m-v)w}(\beta_V \wedge \beta_W,\beta_M \wedge \beta_N).$$
	We also have $$\beta_M \wedge \beta_N \wedge \beta_a \wedge \beta_b = \beta_M \wedge \beta_a \wedge \beta_N \wedge \beta_b,$$ as $a$ is even, so using $\beta_V \wedge \beta_{\nu V} = \beta_M \wedge \beta_a$ and $\beta_W \wedge \beta_{\nu W} = \beta_N \wedge \beta_b$, we require
	$$(\beta_V \wedge \beta_W, \beta_V \wedge \beta_W \wedge \beta_{\nu V \oplus \nu W})*(\beta_V \wedge \beta_{\nu V}\wedge\beta_W \wedge \beta_{\nu W} ,\beta_M \wedge \beta_N) = (-1)^{(m-v)w}(\beta_V \wedge \beta_W,\beta_M \wedge \beta_N).$$
	But now $$\beta_V \wedge \beta_{\nu V}\wedge\beta_W \wedge \beta_{\nu W} = (-1)^{(m-v)w} \beta_V \wedge \beta_{W}\wedge\beta_{\nu V} \wedge \beta_{\nu W},$$
	so, after all that, we see that the Quillen orientation of the normal bundle to $V \times W$ is simply $$\beta_{\nu V \oplus \nu W} = \beta_{\nu V} \wedge \beta_{\nu W}.$$
\end{proof}

\begin{proposition}\label{P: exterior commutativity}\index{co-orientation!of external product!commutativity}
	Let $f \colon V \to M$ and $g \colon W \to N$ be co-oriented maps from manifolds with corners to manifolds without boundary and suppose $f \times g \colon V \times W \to M \times N$ is given the product co-orientation.
	Let $\tau \colon N \times M \to M \times N$ be the diffeomorphism that interchanges coordinates.
	Denote the pullback of $f \times g$ by $\tau$ as $\tau^*(V \times W) \to N \times M$. Then this pullback is diffeomorphic, as co-oriented maps, to $(-1)^{(m-v)(n-w)}W \times V \to N \times M$. In other words,
	$$\tau^*(V \times W) = (-1)^{(m-v)(n-w)}W \times V.$$
\end{proposition}
We assume $M$ and $N$ to be without corners so that we can properly use the pullback construction, which requires transversality, in the proposition and its proof.
However, the pullback is by a diffeomorphism, so this result should extend without problem to more general settings.
\begin{proof}
	This is clear at the level of spaces, so we focus on co-orientations.
	Let $(\beta_V,\beta_M)$ and $(\beta_W,\beta_N)$ be local representations of the co-orientations at some points.
	The product co-orientation of $V \times W \to M \times N$ is $(-1)^{(m-v)w}(\beta_V \wedge \beta_W,\beta_M \wedge \beta_N)$.

	As in \cref{L: Quillen product co-orientation}, we consider Quillen co-orientations representing $f$ and $g$ via embeddings $e_V \colon V \into M \times \R^a$ and $e_W \colon W \into N \times \R^b \to N$ with $a$ and $b$ even.
	This is sufficient as we know that the pullback construction is independent of $a$ and $b$ for sufficiently large dimensions by \cref{L: pullback co well defined}.
	Assuming the other notation from \cref{L: Quillen product co-orientation}, that lemma tells us that the normal co-orientation of $V \times W$ in $M \times N \times \R^a \times \R^b$ is $$\beta_{\nu V \oplus \nu W} = \beta_{\nu V} \wedge \beta_{\nu W}.$$

	Now using this Quillen co-orientation for $f \times g$, we pull back by the diffeomorphism $\tau \colon N \times M \to M \times N$, obtaining the composition we can write $W \times V \into N \times M \times \R^a \times \R^b \to N \times M$.
	The pulled back normal bundle is still oriented in each fiber as $\beta_{\nu V} \wedge \beta_{\nu W}$ (though of course the order of actual local coordinates have now been jumbled around).
	By definition, the pullback co-orientation is $(\beta_W \wedge \beta_V,\beta_N \wedge \beta_M)$ if and only if $$\beta_W \wedge \beta_V \wedge \beta_{\nu V} \wedge \beta_{\nu W} = \beta_N \wedge \beta_M \wedge \beta_{a+b},$$
	and as $a$ and $b$ are even this last expression is equal to
	$\beta_N \wedge \beta_b \wedge \beta_M \wedge \beta_{a}.$ But by the previous choices, $\beta_V \wedge \beta_{\nu V} = \beta_M \wedge \beta_a$ and $\beta_W \wedge \beta_{\nu W} = \beta_N \wedge \beta_b$.
	So
	\begin{align*}
		\beta_N \wedge \beta_b \wedge \beta_M \wedge \beta_{a}
		& = \beta_W \wedge \beta_{\nu W} \wedge \beta_V \wedge \beta_{\nu V} \\
		& = (-1)^{v(n+b-w)}\beta_W \wedge \beta_V \wedge \beta_{\nu W} \wedge \beta_{\nu V} \\
		& = (-1)^{v(n+b-w)+(m+a-v)(n+b-w)}\beta_W \wedge \beta_V \wedge \beta_{\nu V} \wedge \beta_{\nu W} \\
		& = (-1)^{v(n-w)+(m-v)(n-w)}\beta_W \wedge \beta_V \wedge \beta_{\nu V} \wedge \beta_{\nu W}, \\
	\end{align*}
	where again we use that $a$ and $b$ are even.

	Therefore, the pullback co-orientation is $(-1)^{v(n-w)+(m-v)(n-w)}(\beta_W \wedge \beta_V,\beta_N \wedge \beta_M)$, which is $(-1)^{(m-v)(n-w)}$ times the product co-orientation of $W \times V$, as claimed.
\end{proof}

The next result concerns co-oriented products in which one map is the identity.
We show that such products are simply pullbacks by projections.

\begin{proposition}\label{P: projection pullbacks}\index{co-orientation!of pullback!of a projection is product co-orientation}
	Let $f \colon V \to M$ be a co-oriented map from a manifolds with corners to a manifold without boundary, and let $\id_N \colon N \to N$ be the identity map of a manifold with corners with the tautological co-orientation.
	\begin{enumerate}
		\item The co-oriented pullback of $V$ by the projection $\pi_1 \colon M \times N \to M$ is $f \times \id_N \colon V \times N \to M \times N$ with its product co-orientation, i.e.\ $\pi_1^*V = V \times N$.
		\item The co-oriented pullback of $V$ by the projection $\pi_2 \colon N \times M \to M$ is $\id_N \times f \colon N \times V \to N \times M$ with its product co-orientation, i.e.\ $\pi_2^*V = N \times V$.
	\end{enumerate}
\end{proposition}

\begin{proof}
	As the projections are submersions, the required transversality conditions to ensure the existence of the pullbacks are met.
	These claims are then clear concerning maps of topological spaces, so we need only verify the co-orientations.

	As in the preceding argument, we start again with an embedding $e \colon V \into M \times \R^a$ to establish a Quillen co-orientation for $f$.
	We again may assume $a$ to be even for simplicity.
	We write the co-orientation of $f$ as $(\beta_V,\beta_M)$, and we let $\nu V$ denote the normal to $e(V)$ and orient $\nu V$ so that $\beta_V \wedge \beta_{\nu V} = \beta_M \wedge \beta_a$, writing $\beta_a$ for the standard orientation of $\R^a$.

	For the second statement, the product co-orientation of $\id_N \times f \colon N \times V \to N \times M$ is $(\beta_N \wedge \beta_V,\beta_N \wedge \beta_M)$, as the domain and codomain of $\id_N$ have the same dimension.
	The pullback by the projection $N \times M \to M$ gives us the embedding/projection sequence $N \times V\xhookrightarrow{\id_N \times e} N \times M \times \R^a \to N \times M$, and the orientations of the pullback of the normal bundle $\nu V$ by $\pi_2 \times \id_{\R^a}$ is again $\beta_{\nu V}$ at each point of $N \times e(V)$.
	So now from \cref{D: pullback coorient}, the pullback has the product co-orientation if and only if $\beta_N \wedge \beta_V \wedge \beta_{\nu V} = \beta_N \wedge \beta_M \wedge \beta_a$.
	But $ \beta_V \wedge \beta_{\nu V} = \beta_M \wedge \beta_a$ by assumption, so this holds.

	For the first statement, the product co-orientation of $f \times \id_N \colon V \times N \to M \times N$ is $(-1)^{(m-v)n}(\beta_V \wedge \beta_N,\beta_M \wedge \beta_N)$.
	The pullback by the projection $M \times N \to M$ gives us an embedding/projection sequence $V \times N \into M \times N \times \R^a \to M \times N$ (where the first arrow is the composition of $e \times \id_N$ with a permutation of coordinates), and the orientation of the pullback of the normal bundle $\nu V$ by $\pi_1 \times \id_{\R^a}$ is again $\beta_{\nu V}$.
	So now from the definition, the pullback has the product co-orientation if and only if $\beta_V \wedge \beta_N \wedge \beta_{\nu V} = (-1)^{(m-v)n}\beta_M \wedge \beta_N \wedge \beta_a$.
	But $ \beta_V \wedge \beta_{\nu V} = \beta_M \wedge \beta_a$ by assumption, so
	\begin{align*}
		\beta_V \wedge \beta_N \wedge \beta_{\nu V}& = (-1)^{(m+a-v)n}\beta_V \wedge \beta_{\nu V} \wedge \beta_N\\
		& = (-1)^{(m+a-v)n}\beta_M \wedge \beta_a \wedge \beta_N\\
		& = (-1)^{(m-v)n}\beta_M \wedge \beta_N \wedge \beta_a.\qedhere
	\end{align*}
\end{proof}

The next proposition shows that the exterior product construction is natural.

\begin{proposition}\label{P: natural exterior}\index{pullback!of external product}\index{co-orientation!of pullback!of external product map}
	Let $f \colon V \to M$ and $g \colon W \to N$ be co-oriented maps of manifolds with corners with $M$ and $N$ having no boundaries.
	Let $h \colon X \to M$ and $k \colon Y \to N$ be maps of manifolds with corners that are transverse to $f$ and $g$ respectively.
	Then $(h \times k)^*(V \times W) = h^*V \times k^*W$ as spaces with co-oriented maps to $X \times Y$.
\end{proposition}

\begin{proof}
	It is easy to show that $h \times k$ is transverse to $f \times g$, so we focus on co-orientation.
	As in the preceding proofs, we write the Quillen orientation of the normal bunle to the embedded image of $V \times W \into M \times N \times \R^a \times \R^b$ as $\beta_{\nu V \oplus \nu W} = \beta_{\nu V} \wedge \beta_{\nu W}$.
	Then the pullback $P = (h \times k)^*(V \times W)$ is co-oriented by $(\beta_P,\beta_{X \times Y})$ if and only if we choose $\beta_P$ and $\beta_{X \times Y}$ so that
	$$\beta_{P} \wedge \beta_{\nu V} \wedge \beta_{\nu W} = \beta_{X \times Y} \wedge \beta_{a+b}.$$

	On the other hand, we know $h^*V$ is co-oriented by $(\beta_{h^*V},\beta_X)$ if an only if $\beta_{h^*V} \wedge \beta_{\nu V} = \beta_X \wedge \beta_a$, and $k^*W$ is co-oriented by $(\beta_{k^*W},\beta_Y)$ if an only if $\beta_{k^*W} \wedge \beta_{\nu W} = \beta_Y \wedge \beta_b$.
	Assuming these hold, then $h^*V \times k^*W$ is co-oriented by $$(-1)^{(x-(v+x-m))(w+y-n)}(\beta_{h^*V} \wedge \beta_{k^*W},\beta_X \wedge \beta_Y) = (-1)^{(m-v)(w+y-n)}(\beta_{h^*V} \wedge \beta_{k^*W},\beta_X \wedge \beta_Y).$$

	Now, continuing to assume the equalities of the last paragraph and taking $a$ and $b$ even as usual, we have
	\begin{align*}
		\beta_{h^*V} \wedge \beta_{k^*W} \wedge \beta_{\nu V} \wedge \beta_{\nu W}
		& = (-1)^{(w+y-n)(m-v)}\beta_{h^*V} \wedge \beta_{\nu V} \wedge \beta_{k^*W} \wedge \beta_{\nu W}\\
		& = (-1)^{(w+y-n)(m-v)}\beta_X \wedge \beta_a \wedge \beta_Y \wedge \beta_b \\
		& = (-1)^{(w+y-n)(m-v)}\beta_X \wedge \beta_Y \wedge \beta_a \wedge \beta_b.
	\end{align*}
	So if we take $\beta_P = (-1)^{(w+y-n)(m-v)}\beta_{h^*V} \wedge \beta_{k^*W}$ and $\beta_{X \times Y} = \beta_X \times \beta_Y$, then this also gives us $\beta_{P} \wedge \beta_{\nu V} \wedge \beta_{\nu W} = \beta_{X \times Y} \wedge \beta_{a+b}$.
	Therefore,
	$(h \times k)^*(V \times W)$ is also co-oriented by
	$$(\beta_P,\beta_X \wedge \beta_Y) = ((-1)^{(w+y-n)(m-v)}\beta_{h^*V} \wedge \beta_{k^*W},\beta_X \wedge \beta_Y).$$ We conclude $(h \times k)^*(V \times W) = h^*V \times k^*W$.
\end{proof}

\subsubsection{Applications of co-oriented exterior products to co-oriented fiber products}\label{S: product relations}

Having established some elementary properties for our exterior product, we can now relate co-oriented exterior products to co-oriented fiber products.
These relationships correspond to those in singular cohomology between the exterior and cup products of cochains, though one very nice feature is that in our context these relationships all hold ``on the nose'' at the cochain level and can be proven without any need for Alexander-Whitney maps or any other approximations to the diagonal.
This will also be useful for proving some properties of the co-oriented fiber product that we have deferred so far, including associativity and graded commutativity.

We start with the following version of a well-known fact.

\begin{lemma}\label{L: alternative transversality}\index{transversality!relation to the diagonal}\index{diagonal map!relation to transversality}\index{diagonal map}\index{$D$@$\diag$}
	Suppose $f \colon V \to M$ and $g \colon W \to M$ are maps from manifolds with corners to a manifold without boundary.
	Let $\diag \colon M \to M \times M$ be the diagonal map $\diag(x) = (x,x)$.
	Then $f$ and $g$ are transverse if and only if $f \times g$ is transverse to $\diag$.
\end{lemma}

\begin{proof}
	As transversality is determined stratum by stratum, it suffices to suppose $V$ and $W$ are manifolds without boundary.
	We briefly recall the argument in this case.

	First suppose $f$ and $g$ are transverse and that $f(x) = g(y) = z \in M$.
	Then $Df(T_xV)+Dg(T_yW) = T_zM$.
	Now suppose $(a,b) \in T_{(z,z)}(M \times M) \cong T_zM \oplus T_zM$.
	Write $a = v_a+w_a$ with $v_a \in Df(T_xV)$ and $w_a \in Dg(T_yW)$.
	Similarly, write $v = v_b+w_b$.
	Then
	\begin{align*}
		(a,b)& = (v_a+w_a,v_b+w_b)\\
		& = (v_a-v_b+v_b+w_a, v_b+w_a-w_a+w_b)\\
		& = (v_a-v_b,-w_a+w_b)+(v_b+w_a, v_b+w_a),
	\end{align*}
	which is in $D(f \times g)(T_xV \times T_yW)+D\diag(T_zM)$.
	So $f \times g$ is transverse to $\diag$.

	Conversely, suppose $f \times g$ is transverse to $\diag$, and continue to assume that $f(x) = g(y) = z \in M$ so that $(f \times g)(x,y) = \diag(z)$.
	Let $v \in T_zM$.
	Then there are $a\in T_x V$, $b \in T_y W$, and $c \in T_zM$ such that $$(v,0) = D(f \times g)(a,b) + D\diag(c) = (Df(a), Dg(b)) + (c,c).$$
	Then $Df(a) + c = v$ and $Dg(b) + c = 0$, so $Df(a) - Dg(b) = Df(a) + Dg(-b) = v$.
	So $f$ and $g$ are transverse.
\end{proof}

\begin{proposition}\label{P: cross to cup}\index{pullback!by the diagonal}\item{fiber product!as diagonal pullback}
	Suppose $f \colon V \to M$ and $g \colon W \to M$ are transverse co-oriented maps from manifolds with corners to a manifold without boundary.
	Then the pullback of $V \times W \to M \times M$ by $\diag$ is the co-oriented fiber product $V \times_M W \to M$, i.e.\ $$V \times_M W = \diag^*(V \times W).$$
\end{proposition}

\begin{proof}
	By \cref{L: alternative transversality}, we know $f \times g$ is transverse to $\diag$, so the pullback is defined, and it is a manifold with corners by \cref{pullback}.
	This pullback by $\diag$ is $$P' = \{(v,w,z) \in V \times W \times M \mid (f(v),g(w)) = \diag(z) = (z,z)\},$$ which is diffeomorphic to $$P = V \times_M W = \{(v,w) \in V \times W \mid f(v) = g(w)\}$$ via the projection $(v,w,z) \mapsto (v,w)$ with inverse $(v,w) \mapsto (v,w,f(v))$.
	So we consider the co-orientations, for which we may assume as usual that all manifolds are without boundary by \cref{pullback,P: interior co-orientation}.

	By \cref{L: Quillen product co-orientation}, now with $M = N$, if $f$ and $g$ are co-oriented (at appropriate points) by $(\beta_V,\beta_M)$ and $(\beta_W,\beta_M)$ and if we we take Quillen co-orientations coming from $V \xhookrightarrow{e_V} M \times \R^a \to M$ and $W \xhookrightarrow{e_W} M \times \R^b \to M$ (with $a$ and $b$ assumed even) by orienting $\nu V$ and $\nu W$ in $M \times \R^a$ and $M \times \R^b$, then the co-orientation of the product $f \times g$ has Quillen co-orientation with the normal bundle to $V \times W$ in $M \times M \times \R^{a+b}$ oriented by $\beta_{\nu V} \wedge \beta_{\nu W}$.
	Recall that we interpret this expression so that $\nu V$ is our standard normal bundle to $V$ in $M \times \R^a$ embedded into the first and third factors of $M \times M \times \R^a \times \R^b$ and analogous for $\nu W$.

	Pulling back by the diagonal, we thus obtain from \cref{D: pullback coorient} that the pullback co-orientation is $(\beta_{P'},\beta_M)$ if and only if $\beta_{P'} \wedge \beta_{\nu V} \wedge \beta_{\nu W} = \beta_M \wedge \beta_{a+b}$.
	On the other hand, using the same Quillen co-orientation for $V$, the co-orientation of the fiber product $V \times_M W \to W \to M$ is $(\beta_P,\beta_M)$ if and only if $\beta_P \wedge \beta_{\nu V} = \beta_W \wedge \beta_a$.

	Consider now the following computation, where we use that $\beta_W \wedge \beta_{\nu W} = \beta_{M\times \R^b} = \beta_M \wedge \beta_b$ from the definition of the Quillen orientation of the normal bundle:
	\begin{align*}
		\beta_P \wedge \beta_{\nu V} \wedge \beta_{\nu W}& = \beta_W \wedge \beta_a \wedge \beta_{\nu W}\\
		& = \beta_W \wedge \beta_{\nu W} \wedge \beta_a\\
		& = \beta_M \wedge \beta_b \wedge \beta_a \\
		& = \beta_M \wedge \beta_{a+b}\\
		& = \beta_{P'} \wedge \beta_{\nu V} \wedge \beta_{\nu W}.
	\end{align*}
	This appears to demonstrate that, if we fix $\beta_M$, we get corresponding $\beta_P$ and $\beta_P'$.
	The trouble with this argument, however, is that we've been very free about what spaces exactly all these various local orientations live over.
	Even with our given identification of $P$ and $P'$, the symbol $\nu V$ here represents a normal bundle in at least two or three different spaces.
	To solidify the argument, we need a bit more care.

	First, we recall that our model of $P = V \times_M W$ when computing pullback co-orientations is technically $V \times_{M \times \R^a} (W \times_M \R^a)$, identified with $V \times_M W$ in the obvious way, and similarly for $P'$; see \cref{R: pullback representative 2}.
	It will be convenient here to use these representations, so, relabelling for convenience and letting $e_V(v) = (f(v),e_a(v))$ and $e_W(w) = (g(w),e_b(w))$, we let
	\begin{align*}
		P &= \{(v,w,s)\in V\times W\times \R^a\mid (g(w),s)=e_V(v)\}\\
		P' &= \{(v,w,z,s,t)\in V\times W\times M \times \R^a \times \R^b \mid (z,z,s,t)=(f(v),g(w),e_a(v),e_b(w))\}.
	\end{align*}
	Then $P$ is diffeomorphic to our standard $V \times_M W$ via $(v,w) \leftrightarrow (v,w,e_a(v))$, and $P'$ is diffeomorphic to $(V \times W) \times_{M \times M} M$ via $(v,w,z) \leftrightarrow (v,w,z,e_a(v),e_b(w))$. Furthermore, $P \cong P'$ via $(v,w,s)\leftrightarrow (v,w,g(w),s,e_b(w))$.

	We now consider the following diagram, in which the righthand horizontal maps are those occurring in the pullback co-orientation construction:
	\[
	\begin{tikzcd}
		P \arrow[r, hook, "\pi_{W \times \R^a}"] \arrow[d, "\cong" swap] &[1cm]
		W \times \R^a \arrow[r, "g \times \id_{\R^a}"] \arrow[d, hook, "h"] &[1.3cm]
		M \times \R^a \\
		P' \arrow[r, hook, "\pi_{M\times \R^a \times \R^b}"] &
		M \times \R^a \times \R^b \arrow[r, hook, "\diag \times \id_{\R^a} \times \id_{\R^b}"] &
		M \times M \times \R^a \times \R^b. \arrow[u, "\Pi"]
	\end{tikzcd}
	\]
	The lefthand vertical map is our diffeomorphism $P \cong P'$. We let $h$ be the map $h(w,s) = (g(w),s,e_b(w))$, and we let $\Pi$ be the projection onto the first and third factors.
	One readily checks that this diagram commutes and that the maps with hooks are embeddings.

	Let $p = (v,w,s)\in P$, let $p' = (v,w,z,s,t)$ be the corresponding point in $P'$, so in particular $z = f(v) = g(w)$.
	Then taking derivatives at these points and their images, we obtain the diagram
	\[
	\begin{tikzcd}
		T_pP \arrow[r, hook, "D\pi_{W \times \R^a}"] \arrow[d, "\cong" swap] &[1cm]
		T_wW \oplus T_s\R^a \arrow[r, "D(g \times \id_{\R^a})"] \arrow[d, hook, "Dh"] &[1.6cm]
		T_zM \oplus T_s\R^a \\
		T_{p'}P' \arrow[r, hook, "D\pi_{M \times \R^a \times \R^b}"] &
		T_zM \oplus T_s\R^a \oplus T_t\R^b \arrow[r, hook, "D(\diag \times \id_{\R^a} \times \id_{\R^b})"] &
		T_zM \oplus T_zM \oplus T_s\R^a \oplus T_t\R^b. \arrow[u, "\Pi"]
	\end{tikzcd}
	\]
	Now the normal bundle $\nu V$ to $e_V(V) \subset M \times \R^a$ is our standard normal bundle for the construction of a co-orientation of $V \times_M W$, and we know it pulls back to the normal bundle of $P$, identified with $(g \times \id_{\R^a})^{-1}(e_V(V)) \subset W \times \R^a$.

	By construction, we typically think of choosing a splitting so that we can identify a normal space $\nu V$ of $P$ in $W \times \R^a$ at $p$ with a subspace $N$ of $T_{w,s}(W \times \R^a)=T_wW \oplus T_s\R^a$ such that the restriction to $N$ of the composition of $D(g \times \id_{\R_a})$ with the quotient map to $T_{(g(w),s)}(M \times \R^a)/ T_{(g(w),s)}V$ is an isomorphism.
	With such a choice, it follows that the composition the other way around the righthand square, followed by the quotient map, is also an isomorphism onto $T_{(g(w),s)}(M \times \R^a)/ T_{(g(w),s)}V$, and so we can identify the image of $N$ under $D(\diag \times \id_{\R^a}\times \id_{\R^b}) Dh$ with the fiber of the pullback of $\nu V$ by $\Pi$. But now we observe that if $(\xi,\eta) \in T_wW \oplus T_s\R^a$, then
	$$D(\diag \times \id_{\R^a}\times \id_{\R^b}) Dh(\xi,\eta) = (Dg(\xi),Dg(\xi), \eta, De_b(\xi)).$$
	But $(Dg(\xi),De_b(\xi)) = De_W(\xi)$, so $(0, Dg(\xi),0, De_b(\xi))$ is in the tangent space of our chosen embedding of $V \times W$ into $M \times M \times \R^a \times \R^b$.
	So, the image $D(\diag \times \id_{\R^a}\times \id_{\R^b})Dh(N)$ represents our standard choice of a $\nu V$ summand in the normal bundle to $V \times W$, it simply corresponds to a different choice of splitting.
	In particular, this image and the $\nu V$ summand from \cref{L: Quillen product co-orientation} project to isomorphic images in $T(M \times M \times \R^a \times \R^b)/T(V \times W)$.

	The upshot is that we can thus identify our versions of $\nu V$ in $W \times \R^a$ and $M \times \R^a \times \R^b$ via $Dh$, and similarly identify $T_pP$ and $T_{p'}P'$ via the diagram. In fact, the lefthand square consist entirely of embeddings, we can now properly interpret the above computation as a computation entirely with elements of $T_zM \oplus T_s\R^a \oplus T_t\R^b$ with $\nu W$ being the normal from the embedding of $W$ into the first and third coordinates of $M \times \R^a \times \R^b$.
\end{proof}

\begin{corollary}\label{C: fiber natural pullback}\index{pullback!of fiber product}
	Suppose $f \colon V \to M$ and $g \colon W \to M$ are transverse co-oriented maps of manifolds with corners to a manifold without boundary, that $N$ is a manifold with corners, and that $h \colon N \to M$ is transverse to $f$, $g$, and $f \times_M g \colon V \times_M W \to M$.
	Then
	$$h^*(V \times_M W) = h^*V \times_N h^*W$$
	as manifolds with co-oriented maps to $N$.
\end{corollary}

\begin{proof}
	Let $\diag_M \colon M \to M \times M$ and $\diag_N \colon N \to N \times N$ denote the diagonal maps.
	Using \cref{P: cross to cup,P: pullback functoriality,P: natural exterior}, and that $\diag_M h = (h \times h) \circ \diag_N$,
	we compute
	\begin{align*}
		h^*(V \times_M W)& = h^*\diag_M^*(V \times W)\\
		& = (\diag_M h)^*(V \times W)\\
		& = ((h \times h) \circ \diag_N)^*(V \times W)\\
		& = \diag_N^*(h \times h)^*(V \times W)\\
		& = \diag_N^*(h^*V \times h^*W)\\
		& = h^*V \times_N h^*W.
	\end{align*}
	To apply \cref{P: pullback functoriality} in the second line, we note that $h$ is transverse to $\diag_M^*(V \times W) = V \times_M W$ by assumption.
	For the fourth line we observe that $h \times h$ is transverse to $f \times g$ because $h$ is transverse to $f$ and $g$, and the composite $(h \times h) \circ \diag_N = \diag_M h$ is transverse to $f \times g$ by our assumptions and \cref{L: transverse to pullback}.
\end{proof}

\begin{corollary}[Associativity of co-oriented fiber products]\label{C: fiber assoc}\index{fiber product!associativity}\index{co-orientation!of fiber product!associativity}
	Suppose $f \colon V \to M$, $g \colon W \to M$, and $h \colon X \to M$ are co-oriented maps from manifolds with corners to a manifold without boundary such that the following pairs are transverse (see \cref{R: multiproducts}): $(V,W)$, $(W,X)$, $(V \times_M W,X)$, and $(V,W \times_M X)$.
	Then $$(V \times_M W) \times_M X = V \times_M (W \times_M X)$$ as co-oriented fiber products mapping to $M$.
\end{corollary}

\begin{proof}
	We compute using \cref{P: cross to cup,P: pullback functoriality,P: exterior associativity,P: natural exterior} and that $(\id_M \times \diag)\diag = (\diag \times \id_M)\diag$:
	\begin{align*}
		(V \times_M W) \times_M X& = \diag^*(\diag^*(V \times W) \times X)\\
		& = \diag^*(\diag \times \id_M)^*((V \times W) \times X)\\
		& = ((\diag \times \id_M)\diag)^*((V \times W) \times X)\\
		& = ((\id_M \times \diag)\diag)^*(V \times (W \times X))\\
		& = \diag^*(\id_M \times \diag)^*(V \times (W \times X))\\
		& = \diag^*(V \times \diag^*(W \times X))\\
		& = V \times_M (W \times_M X).\qedhere
	\end{align*}
\end{proof}

\begin{proposition}[Graded commutativity of co-oriented fiber products]\label{P: graded comm}\index{fiber product!commutativity}\index{co-orientation!of fiber product!commutativity}
	Suppose $f \colon V \to M$ and $g \colon W \to M$ are transverse co-oriented maps from manifolds with corners to a manifold without boundary.
	Then as co-oriented fiber products over $M$ we have $\omega_{g \times_M f} = (-1)^{(m-v)(m-w)} \omega_{f \times_M g}$, or, using Notation \ref{R: precise commutativity},
	$$V \times_M W = (-1)^{(m-v)(m-w)} W \times_M V.$$
\end{proposition}

\begin{proof}
	By \cref{P: cross to cup}, the fiber products $V \times_M W \to M$ and $W \times_M V \to M$ are respectively the pullbacks of $f \times g \colon V \times W \to M \times M$ and $g \times f \colon W \times V \to M \times M$ by the diagonal map $\diag \colon M \to M \times M$.
	Meanwhile, by \cref{P: exterior commutativity}, we have
	$$\tau^*(V \times W) = (-1)^{(m-v)(m-w)}W \times V,$$
	where $\tau \colon M \times M \to M \times M$ is the map that interchanges the factors.
	We observe that $\tau \diag = \diag$, so using the functoriality of pullbacks from \cref{P: pullback functoriality} we can compute
	\begin{align*}
		W \times_M V &= \diag^*(W \times V)\\
		&= (-1)^{(m-v)(m-w)} \diag^* \tau^*(V \times W)\\
		&= (-1)^{(m-v)(m-w)}(\tau \diag)^*(V \times W)\\
		&= (-1)^{(m-v)(m-w)}(\diag)^*(V \times W)\\
		&= (-1)^{(m-v)(m-w)} V \times_M W.\qedhere
	\end{align*}
\end{proof}

\begin{example}\label{E: embedded}\index{pullback!of embedding}\index{co-orientation!of pullback!of embedding}
	In \cref{E: V embedded}, we considered pullback co-orientations $V \times_M W \to W$ when $V \into M$ was embedded.
	In this example, we discuss the case where $W$ is embedded.

	Let $f \colon V \to M$ and $g \colon W \to M$ be transverse maps from manifolds with corners to a manifold without boundary.
	Suppose $f$ is co-oriented and $g \colon W \to M$ is an embedding, which we use throughout the example to identify $W$ with its image.
	Let $(x,y) \in V \times_M W$.
	For the remainder of the argument, we fix a local orientation $\beta_M$ at $f(x) \in M$, and we choose $\beta_V$ at $x$ so that $(\beta_V,\beta_M)$ is the co-orientation of $V$ at $x$.
	Furthermore, even though $g$ might not be co-oriented, let us choose a Euclidean neighborhood $U$ of $y$ in $W$ and an arbitrary co-orientation $(\beta_U,\beta_M)$ on the restriction of $g$ to $U$.
	This can be done as $U$ is contractible and so $g|_U$ is co-orientable.

	Although we are interested in $V \times_M U$, the definition of the pullback co-orientation makes it easier to work with $U \times_M V$ when $U$ is embedded; see \cref{E: V embedded}.
	As we have chosen a co-orientation for $U \into M$ and as $V \to M$ comes with a co-orientation, we can consider the fiber product $P_U = U \times_M V \to M$.
	As $U$ is embedded in $M$, we have $P_U = f^{-1}(U) \subset V$.
	If we choose an orientation $\beta_{\nu U}$ of the normal bundle to $U$ in $M$ at $y$ so that $\beta_U \wedge \beta_{\nu U} = \beta_M$, then by \cref{E: V embedded} and the definition of fiber product co-orientation, the map $P_U \xr{f} M$ is co-oriented by $(\beta_P,\beta_M)$, where $\beta_P$ is chosen so that $\beta_P \wedge \beta_{\nu U} = \beta_V$, as usual letting $\nu U$ here also stand for its pullback as a normal bundle of $P_U$ in $V$.
	We note that if we had chosen the opposite co-orientation for $U \into M$ then, as we have fixed $\beta_V$ and $\beta_M$, the result would be to replace our current $\beta_U$ with $-\beta_U$, which would result in also reversing the signs of $\beta_{\nu U}$ and $\beta_P$.
	In particular, the fiber product would have the opposite co-orientation.

	Next, we will apply graded commutativity.
	By \cref{P: graded comm}, we have $U \times_M V = (-1)^{(m-v)(m-w)}V \times_M U$ as co-oriented fiber products over $M$.
	By \cref{R: precise commutativity}, this means that the fiber product $V \times_M U \to M$ with its fiber product co-orientation corresponds to the fiber product $U \times_M V = f^{-1}(U) \xr{f} M$ with co-orientation at $(x,y)$ given by $(-1)^{(m-v)(m-w)}(\beta_P, \beta_M)$.
	But decomposing the fiber product as the pullback $V \times_M U \to U$ and the inclusion $U \into M$, we write $(-1)^{(m-v)(m-w)}(\beta_P, \beta_M)$ as the composite co-orientation
	$$(-1)^{(m-v)(m-w)}(\beta_P, \beta_M) = (-1)^{(m-v)(m-w)}(\beta_P, \beta_U)*(\beta_U, \beta_M).$$
	As $(\beta_U, \beta_M)$ is our chosen co-orientation for $g|_U$, the pullback co-orientation of $V \times_M U \to U$ must be $(-1)^{(m-v)(m-w)}(\beta_P, \beta_U)$.
	But we have already observed that if we had chosen the opposite co-orientation for $g|_U$, that would reverse the signs of both $\beta_U$ and $\beta_P$, so in either case we obtain the same pullback co-orientation for $V \times_M U \to U$.
	In other words, this description of the co-oriented pullback $V \times_M U \to U$ is independent of our choice of co-orientation for $g|_U$, and so it extends globally to give our pullback co-orientation of $V \times_M W \to W$.

	Summarizing then, as a space we have $P = V \times_M W = f^{-1}(W)$, and the pullback and fiber product are the first map and composite of $f^{-1}(W) \xr{f} W \into M$ (analogously to the case where $f$ was an embedding in \cref{E: V embedded}).
	Furthermore, fixing the co-orientation of $f$ as $(\beta_V, \beta_M)$, the co-orientation of the pullback $V \times_M W \to W$ is by $(-1)^{(m-v)(m-w)}(\beta_P, \beta_W)$, where, if we choose any local orientation $\beta_{\nu W}$ for the normal bundle of $W$ in $M$, then $\beta_W \wedge \beta_{\nu W} = \beta_M$ in $M$ and $\beta_P \wedge \beta_{\nu W} = \beta_V$ after pulling back to $V$.

	It is a nice exercise to confirm that this agrees with the computation of \cref{P: normal pullback} when $V$ and $W$ are both embedded.
\end{example}

\begin{corollary}\label{C: criss cross}\index{co-orientation!of fiber product!of external products}\index{co-orientation!of external product!of fiber products}\index{criss cross}
	Suppose $f \colon V \to M$ and $g \colon W \to M$ are transverse co-oriented maps from manifolds with corners to a manifold without boundary and similarly for $h \colon X \to N$ and $k \colon Y \to N$.
	Then$$(V \times X)\times_{M \times N} (W \times Y) = (-1)^{(m-w)(n-x)} (V \times_M W) \times (X \times_N Y) $$
	as co-oriented fiber products over $M \times N$.
\end{corollary}

\begin{proof}
	With our given transversality assumptions, $f \times h$ is transverse to $g \times k$, so the expression on the left is well defined.
	We then compute using \cref{P: cross to cup,P: pullback functoriality,P: exterior associativity,P: exterior commutativity,P: natural exterior} and letting $\tau$ here be the interchange of the interior $N$ and $M$ in the quadruple product:
	\begin{align*}
		(V \times X)\times_{M \times N} (W \times Y)& = \diag_{M \times N}^*(V \times X \times W \times Y)\\
		& = \diag_{M \times N}^*(\id \times \tau^* \times \id)^*(-1)^{(n-x)(m-w)}(V \times W \times X \times Y)\\
		& = (\diag_M \times \diag_N)^*(-1)^{(n-x)(m-w)}(V \times W \times X \times Y)\\
		& = (-1)^{(n-x)(m-w)}\diag_M^*(V \times W) \times \diag_N^*( X \times Y)\\
		& = (-1)^{(m-w)(n-x)} (V \times_M W) \times (X \times_N Y).
	\end{align*}
	For the third equality, we use that $\diag_M \times \diag_N = (\id_M \times \tau \times \id_N)\diag_{M \times N}$.
\end{proof}

\begin{corollary}\label{C: cross is cup}\index{exterior product!as fiber product of pullbacks of projections}\index{fiber product!of pullback by projection is external product}\index{co-orientation!of external product!as fiber product of pullbacks of projections}\index{co-orientation!of fiber product!of pullbacks by projections is external product}
	Let $V \to M$ and $W \to N$ be maps from manifolds with corners to manifolds without boundary.
	Let $\pi_M \colon M \times N \to M$ and $\pi_N \colon M \times N \to N$ be the projections.
	Then $$V \times W = \pi_M^*(V)\times_{M \times N}\pi_N^*(W)$$ as co-oriented manifolds mapping to $M \times N$.
\end{corollary}

\begin{proof}
	By \cref{P: projection pullbacks}, $\pi_M^*(V) = V \times N$ and $\pi_N^*(W) = M \times W$, so these are transverse as spaces mapping to $M \times N$.
	Then by \cref{C: criss cross,C: cup with identity}, we have
	\begin{align*}
		\pi_M^*(V)\times_{M \times N}\pi_N^*(W)& = (V \times N)\times_{M \times N} (M \times W)\\
		& = (V \times_M M) \times (N \times_N W)\\
		& = V \times W.\qedhere
	\end{align*}
\end{proof}

\subsection{Properties mixing orientations and co-orientations}\label{S: mixing}

In this section we study properties that involve both orientations and co-orientations.
In particular, we are mostly interested in the pullback of a co-oriented map $V \to M$ by a map $W \to M$ with $W$ oriented, in which case the co-orientation of the pullback $V \times_M W \to W$ together with the orientation of $W$ produces an induced orientation on $V \times_M W$ as described in \cref{S: co-orientations}.
As $V \times_M W \to M$ with this orientation will eventually correspond to the cap product when we get to geometric homology and cohomology, we will here refer to this orientation as the \textbf{cap orientation}\index{cap orientation|textbf}\index{orientation!cap}.

The next results concern cap orientations on $V \times_M W$.
We note that, by construction, this oriented manifold comes equipped with a map to $W$ and, by composing with the given map $W \to M$, a map to $M$.

We start with the next result, which is yet another Leibniz formula.
This will lead us eventually to a boundary formula for the cap product of geometric chains and cochains that agrees with that for the boundary of cap products of singular chains and cochains in Spanier \cite[Section 5.6.15]{Span81} and Munkres \cite[Section 66]{Mun84}.
It will also be used in \cref{S: intersection map} to demonstrate that our intersection map $\mc I$, relating geometric cochains to cubical cochains of a cubulation, is a chain map.
This map is a critical component in relating geometric cohomology to other cohomology theories and is also central to the main result about cup products in \cite{FMS-flows}.

\begin{proposition}\label{P: Leibniz cap}\index{cap orientation!boundary formual (Leibniz rule)}
	Let $f \colon V \to M$ and $g \colon W \to M$ be transverse maps of manifolds with corners to a manifold without boundary.
	Suppose $f$ is co-oriented and $W$ is oriented.
	Then $$\bd(V \times_M W) = \left[(-1)^{v+w-m} (\bd V) \times_M W\right] \bigsqcup (V \times_M \bd W)$$
	as oriented manifolds, giving $V \times_M W$, $(\bd V) \times_M W$, and $V \times_M \bd W$ each their cap orientations.
	In the last case, this is the cap orientation of the pullback of $V$ by $\bd W \xr{i_{\bd}} W \xr{g} M$, with $\bd W$ having its boundary orientation as the boundary of $W$.
\end{proposition}

\begin{proof}
	We compute and compare these orientations by first considering the pullback co-orientations as defined in \cref{D: pullback coorient}.
	We proceed by analogy to the proof of the Leibniz rule for the pullback of co-oriented maps in \cref{leibniz}, utilizing the computations already performed there.
	Recall that, as needed, we can assume by working locally that $V$ and $W$ are manifolds without boundary or manifolds without strata of greater depth.

	Recall, in brief, from \cref{D: pullback coorient} that to co-orient the pullback $P = V \times_M W \to W$ we first construct a composition $V \xhookrightarrow{e} M \times \R^N \to M$ and find a Quillen orientation for the normal bundle $\nu V$ of $e(V) \subset M \times \R^N$ as determined by the co-orientation of $V \to M$.
	Then we pull back via $W \times \R^N \to M \times \R^N$ to obtain a normal bundle, also labeled $\nu V$, of $P \subset W \times \R^N$.
	Then we co-orient $P \to W$ locally by $(\beta_P,\beta_W)$ so that $\beta_P \wedge \beta_{\nu V} = \beta_W \wedge \beta_E$, where $\beta_E$ represents the standard orientation of $\R^N$.
	In the case at hand, we can assume $\beta_W$ to represent the global orientation of $W$, and then $\beta_P$ becomes the global cap orientation for $P = V \times_M W$.

	Let $\nu\bd P$ denote an outward pointing normal vector in the tangent bundle to $P$ at a boundary point of $P$, and let $\beta_{\nu\bd P}$ denote the corresponding orientation.
	Then, by definition, $\bd P$ is oriented at that point by $\beta_{\bd P}$ so that $\beta_{\nu\bd P} \wedge \beta_{\bd P} = \beta_P$.
	In other words, with $\beta_P$, $\beta_W$, $\beta_{\nu V}$, and $\beta_E$ given as above, $\bd P$ is oriented by $\beta_{\bd P}$ such that $\beta_{\nu\bd P} \wedge \beta_{\bd P} \wedge \beta_{\nu V} = \beta_W \wedge \beta_E$.

	Now, recall that $\bd(V \times_M W) = (\bd V) \times_M W \bigsqcup V \times_M \bd W$ as spaces and consider a point in $(\bd V) \times_M W$.
	By \cref{leibniz}, at such a point the pullback co-orientation of $(\bd V) \times_M W \to W$ (as the pullback of $\bd V \to M$ by $g$) agrees with boundary co-orientation of the pullback $P = V \times_M W \to W$.
	So continuing to let $(\beta_P,\beta_W)$ denote the pullback co-orientation of $P \to W$ and recalling that the boundary co-orientation utilizes the \textit{inward} normal, the boundary co-orientation of $(\bd V) \times_M W \to W$ is the composite $(\beta_{\bd P}, \beta_{\bd P} \wedge -\beta_{\nu\bd P})*(\beta_P,\beta_W)$ for any $\beta_{\bd P}$.
	But if we choose $\beta_{\bd P}$ to represent the orientation of $\bd P$ found above by orienting $P$ and then taking its boundary orientation, we have $\beta_P = \beta_{\nu\bd P} \wedge \beta_{\bd P} = (-1)^{\dim(\bd P)}\beta_{\bd P} \wedge \beta_{\nu\bd P}$.
	So the boundary co-orientation of $(\bd V) \times_M W \to W$ is the composite
	$$(\beta_{\bd P}, \beta_{\bd P} \wedge -\beta_{\nu\bd P})*((-1)^{\dim(\bd P)}\beta_{\bd P} \wedge \beta_{\nu\bd P},\beta_W) = (-1)^{\dim(\bd P)+1}(\beta_{\bd P},\beta_W).$$
	Thus the resulting cap orientation of $(\bd V) \times_M W$ is $(-1)^{\dim(P)}\beta_{\bd P}$, which is $(-1)^{\dim(P)} = (-1)^{v+w-m}$ times the orientation of $\bd P$ obtained by taking the oriented boundary of $V \times_M W$.

	Next we consider a point in $V \times_M \bd W$.
	Note that by \cref{L: normal pullback} we can consider our outward pointing normal vector $\nu\bd P$ at a point of $\bd P$ to also be an outward pointing normal of $\bd W$; below we write $\nu \bd P = \nu \bd W$.
	Again from the definition of pullback co-orientations, the co-orientation of the pullback $V \times_M \bd W \to \bd W$ is $(\beta_{\bd P},\beta_{\bd W})$ when $\beta_{\bd P} \wedge \beta_{\nu V} = \beta_{\bd W} \wedge \beta_E$.
	If we here give $\bd W$ its boundary orientation, this determines the cap orientation as $\beta_{\bd P}$.
	Wedging with $\beta_{\nu\bd P}$ and using the definition of boundary orientation of $W$, we thus have
	$$\beta_{\nu\bd P} \wedge \beta_{\bd P} \wedge \beta_{\nu V} = \beta_{\nu\bd P} \wedge \beta_{\bd W} \wedge \beta_E = \beta_{\nu\bd W} \wedge \beta_{\bd W} \wedge \beta_E = \beta_W \wedge \beta_E.$$
	But this was exactly our condition above for $\beta_{\bd P}$ to be the boundary orientation of $P$ when $P$ is given its cap orientation.
\end{proof}

Next we describe the cap orientation when $V \to M$ and $W \to M$ are embeddings.
As we've observed in the cases where either both maps are oriented or both maps are co-oriented, this is often an instructive and important example.

\begin{proposition}\label{P: cap of immersions}\index{cap orientation!of embeddings}
	Let $f \colon V \to M$ and $g \colon W \to M$ be transverse embeddings from manifolds with corners to a manifold without boundary.
	Suppose $f$ is co-oriented and $W$ is oriented.
	Then $P = V \times_M W$ is just the intersection of $V$ and $W$ in $M$.
	If $\beta_W$ is the orientation of $W$ and $\beta_{\nu V}$ is the Quillen orientation of the normal bundle to $V$ in $M$, which at points of $P$ we can identify\footnote{See \cref{L: normal pullback}.} with the normal bundle to $P$ in $W$, then the cap orientation $\beta_P$ of $P$ satisfies $\beta_P \wedge \beta_{\nu V} = \beta_W$.
	If $f$ and $g$ are immersions, then this description holds locally.
\end{proposition}

\begin{proof}
	As $f$ is an embedding, we can take $N = 0$ in the definition of the pullback co-orientation, \cref{D: pullback coorient}.
	Then the pullback is just the inclusion of $P = g^{-1}(V) = V \cap W$ into $W$, and by definition the pullback co-orientation has the form $(\beta_P,\beta_W)$, where $\beta_P \wedge \beta_{\nu V} = \beta_W$ and $\nu V$ here is the pullback of the normal bundle of $V$ in $M$ to be the normal bundle of $V \cap W$ in $W$.
	Furthermore, if we take $\beta_W$ to be the given orientation of $W$, then $\beta_P$ is the cap orientation on the intersection by definition.
	The last statement about immersions follows as we can compute the co-orientations locally.
\end{proof}

The following corollary is particularly important and follows immediately from \cref{P: cap of immersions}.

\begin{corollary}\label{C: complementary cap}\index{cap orientation!of embeddings!of complementary dimension}
	Let $f \colon V \to M$ and $g \colon W \to M$ be transverse embeddings from manifolds with corners to a manifold without boundary.
	Suppose $f$ is co-oriented, $W$ is oriented, and $\dim(V) + \dim(W) = \dim(M)$.
	Then $V \times_M W$ is the union of intersection points of $V$ and $W$.
	Such a point $x \in V \cap W$ is positively oriented if and only if the Quillen orientation of the normal bundle $\nu V$ of $V$ at $x$ agrees with the orientation of $W$ at $x$, identifying the fiber of $\nu V$ at $x$ with $T_xW$.
	If $f$ and $g$ are immersions, then this description holds locally.
\end{corollary}

The next two propositions will eventually correspond to the unital identities for the cap product for geometric chains and cochains.
The analogues for singular chains and cochains are the cap product with the cochain $1$ and the cap product with a chain representing the fundamental class, though in this case our underlying spaces do not need to be compact.

\begin{proposition}\label{P: cap with 1}\index{cap orientation!with the identity map}\index{fiber product!with identity map}
	Let $g \colon W \to M$ be a map from an oriented manifold with corners to a manifold without boundary, and consider $M \to M$ as the identity with the tautological co-orientation.
	Then $M \times_M W = W$ as oriented manifolds.
\end{proposition}

\begin{proof}
	By definition, there is a Quillen co-orientation for $M$ consisting of the sequence of identity maps $M \into M \to M$ with the normal bundle to $M$ in itself being the $0$-dimensional vector bundle, which we consider to have positive orientation at each point.
	It follows from the definition of the pullback that the corresponding Quillen co-orientation for $M \times_M W$ comes from the sequence $g^{-1}(M) = W \into W \to W$ with $W$ also having a $0$-dimensional positively-oriented normal bundle in itself.
	Consequently, the pullback co-orientation for $W \to W$ is the tautological one, and so the induced orientation on $W$ is the given one.
\end{proof}

\begin{proposition}\label{P: cap with identity M}\index{cap orientation!with the identity map}\index{fiber product!with identity map}\index{orientation!induced!is cap orientation of fiber product with identity}
	Let $f \colon V \to M$ be a co-oriented map from a manifold with corners to an oriented manifold without boundary, and consider $M$ equipped with its identity map $M \to M$.
	Then $V \times_M M$ with the cap orientation is simply $V$ with its induced orientation.
	In other words, following \cref{N: hat check}, $$V \times_M M = \check V.$$
\end{proposition}

\begin{proof}
	By definition, the cap orientation of $V \times_M M$ is the induced orientation arising from the co-oriented pullback of $V \to M$ by the identity map $M \to M$.
	But by \cref{P: pullback functoriality} this co-oriented pullback is simply $f \colon V \to M$ again.
\end{proof}

The next property relates products of pullbacks with pullbacks of products. Again, our signs agree with the cap product formulas in Spanier, in this case \cite[Section 5.6.21]{Span81}.

\begin{proposition}\label{P: cap cross}\index{cap orientation!of external product!}\index{exterior product!of fiber products with cap orientations}
	Let $f \colon V \to M$ and $g:X \to N$ be co-oriented maps from manifolds with corners to manifolds without boundary, and let $h \colon W \to M$ and $k \colon Y \to N$ be maps with $W$ and $Y$ oriented manifolds with corners.
	Suppose that $V$ is transverse to $W$ and that $X$ is transverse to $Y$.
	Then,
	$$(V \times X)\times_{M \times N} (W \times Y) = (-1)^{(x+y-n)(m-v)} (V \times_M W) \times (X \times_N Y),$$
	as oriented manifolds with the pullbacks given their cap orientations.
\end{proposition}

\begin{proof}
	Let $\beta_W$ and $\beta_Y$ denote the orientations of $W$ and $Y$.
	Then $W \times Y$ is oriented by $\beta_W \wedge \beta_Y$.

	Now let $P = V \times_M W$ and $P' = X \times_N Y$.
	By definition, $P$ and $P'$ are oriented by the orientations $\beta_P$ and $\beta_{P'}$ such that $(\beta_P,\beta_W)$ and $(\beta_{P'},\beta_Y)$ are the pullback co-orientations for $P \to W$ and $P' \to Y$.
	Then $P \times P'$ is oriented by $\beta_P \wedge \beta_{P'}$.

	Furthermore, using our construction of pullback co-orientations, $\beta_P$ and $\beta_{P'}$ are such that $\beta_P \wedge \beta_{\nu V} = \beta_W \wedge \beta_a$ and $\beta_{P'} \wedge \beta_{\nu X} = \beta_Y \wedge \beta_b$, where $\beta_a$ and $\beta_b$ are the standard orientations of the Euclidean spaces $\R^a$ and $\R^b$ and we are free to take $a$ and $b$ to be even integers.

	By \cref{L: Quillen product co-orientation}, we have that the Quillen co-orientation of $V \times X \to M \times N$ is represented by an embedding $V \times X \into M \times N \times \R^a \times \R^b$ with normal bundle $\nu V \oplus \nu X$ suitably interpreted.
	So letting $Q = (V \times X)\times_{M \times N} (W \times Y)$, the orientation $\beta_Q$ is the one such that $(\beta_Q,\beta_W \wedge \beta_Y)$ is the pullback co-orientation, i.e.\ the one such that $\beta_Q \wedge \beta_{\nu V} \wedge \beta_{\nu X} = \beta_W \wedge \beta_Y \wedge \beta_a \wedge \beta_b$.
	But then we compute, using $a$ and $b$ even,
	\begin{align*}
		\beta_W \wedge \beta_Y \wedge \beta_a \wedge \beta_b& = \beta_W \wedge \beta_a \wedge \beta_Y \wedge \beta_b\\
		& = \beta_P \wedge \beta_{\nu V} \wedge \beta_{P'} \wedge \beta_{\nu X}\\
		& = (-1)^{|P'||\nu V|}\beta_P \wedge \beta_{P'} \wedge \beta_{\nu V} \wedge \beta_{\nu X}\\
		& = (-1)^{(x+y-n)(m-v)}\beta_P \wedge \beta_{P'} \wedge \beta_{\nu V} \wedge \beta_{\nu X}.
	\end{align*}
	So $\beta_Q = (-1)^{(x+y-n)(m-v)}\beta_P \wedge \beta_{P'} = (-1)^{(x+y-n)(m-v)}\beta_{P \times P'}$.
\end{proof}

The following technical lemma will be used to prove \cref{P: OC mixed associativity}, which will eventually become the associativity relation among cup and cap products, i.e.\ $(a \smile b) \frown x = a \frown (b \frown x)$.

\begin{lemma}\label{L: same induced}\index{cap orientation!associativity|(}
	Let $f \colon V \to M$ and $g \colon W \to M$ be transverse maps from manifolds with corners to a manifold without boundary.
	Suppose that $f$ is co-oriented and that $W$ and $M$ are oriented, with respective (global) orientations $\beta_W$ and $\beta_M$.
	Suppose we co-orient $g$ by $(\beta_W,\beta_M)$.
	Then the cap orientation of $V \times_M W$ (i.e.\ that induced from the pullback co-orientation of $V \times_M W \to W$ and the orientation of $W$) is the same as the orientation induced on $V \times_M W$ by the fiber product co-orientation of $V \times_M W \to M$ and the orientation of $M$.

	In other words, following \cref{N: hat check}, $$V \times_M W = (V \times_M \hat W)\,\check{\vrule height1.3ex width0pt}.$$

	In particular, the orientation of $(V \times_M \hat W)\,\check{\vrule height1.3ex width0pt}$ does not depend on the orientation of $M$.
\end{lemma}

\begin{proof}
	By definition, the orientation of $V \times_M W$ induced from the orientation of $W$ and the pullback co-orientation of $V \times_M W \to W$ is the orientation $\beta_P$ such that $(\beta_P,\beta_W)$ is the pullback co-orientation.
	But then the fiber product co-orientation is the composite $(\beta_P,\beta_W)*(\beta_W,\beta_M) = (\beta_P,\beta_M)$.
	So the orientation induced by the orientation of $M$ and the composite co-orientation is again $\beta_P$.

	The last statement holds because $V \times_M W$ does not depend on the orientation of $M$.
\end{proof}

\begin{proposition}[Mixed associativity 1]\label{P: OC mixed associativity}\index{fiber product!mixed associativity}
	Let $f \colon V \to M$ and $g \colon W \to M$ be co-oriented maps from manifolds with corners to a manifold without boundary.
	Let $h \colon Z \to M$ be a map with $Z$ an oriented manifold with corners.
	Then, assuming sufficient transversality for all terms to be defined (see \cref{R: multiproducts}),
	$$(V \times_M W) \times_M Z = V \times_M (W \times_M Z),$$
	as oriented manifolds.
	Here, on the left, $V \times_M W$ has its fiber product co-orientation, so we can form the cap orientation of the pullback over $Z$.
	On the right we first give $W \times_M Z$ its cap orientation as a pullback over $Z$ and then use that to form the cap orientation of $V \times_M (W \times_M Z)$.\index{cap orientation!associativity|)}
\end{proposition}

\begin{proof}
	First suppose $M$ is orientable and that we have given it an arbitrary, but fixed, orientation.
	Then applying \cref{L: same induced}, the cap orientation of $(V \times_M W) \times_M Z$ is the same as the orientation induced by the orientation of $M$ and the fiber product co-orientation $(V \times_M W) \times_M Z \to M$, after co-orienting $Z \to M$ with the co-orientation induced by the orientations of $Z$ and $M$, i.e.\ $$(V \times_M W) \times_M Z = ((V \times_M W) \times_M \hat Z)\,\check{\vrule height1.3ex width0pt}.$$
	Similarly, applying \cref{L: same induced} twice, the cap orientation of $V \times_M (W \times_M Z)$ is the same as that induced by the orientation of $M$ and the iterated fiber product co-orientation $V \times_M (W \times_M Z) \to M$ again coming from the canonical co-orientation of $Z \to M$:
	\begin{align*}
		V \times_M (W \times_M Z) &= V \times_M (W \times_M \hat Z)\,\check{\vrule height1.3ex width0pt}\\
		&=(V \times_M (W \times_M \hat Z))\,\check{\vrule height1.3ex width0pt}.
	\end{align*}
	But now the co-oriented fiber products $(V \times_M W) \times_M \hat Z$ and $V \times_M (W \times_M \hat Z)$ are the same by \cref{C: fiber assoc}.

	Next, suppose $M$ is not necessarily orientable.
	We know from their constructions that $(V \times_M W) \times_M Z$ and $V \times_M (W \times_M Z)$ are oriented manifolds, and it is not difficult to see that they are diffeomorphic, both being diffeomorphic to $\{(v,w,z) \in V \times W \times Z \mid f(v) = g(w) = h(z)\}$.
	So it suffices to consider these as identical spaces and to show that their induced orientations agree at any arbitrary point.
	If $(v,w,z)$ is such a point, consider its image $a = f(v) = g(w) = h(z) \in M$.
	Let $U$ be a Euclidean neighborhood of $a$, and consider the restrictions of $f$, $g$, and $h$ to $f^{-1}(U)$, $g^{-1}(U)$, and $h^{-1}(U)$.
	The resulting products over $U$ give us the pieces of $(V \times_M W) \times_M Z$ and $V \times_M (W \times_M Z)$ over $U$, and the resulting orientations will be compatible with those of the full manifolds $(V \times_M W) \times_M Z$ and $V \times_M (W \times_M Z)$, as orientations and co-orientations of fiber products are determined locally (see \cref{R: local pullback co-orientations} and the construction of fiber product orientations).
	But as $U$ is orientable, the preceding argument shows that these orientations must agree with each other.
\end{proof}

The following property will eventually manifest itself in geometric (co)homology as the familiar naturality formula for cap products $f_*(f^*(\alpha)\frown x)) = \alpha\frown f_*(x)$.

\begin{proposition}\label{P: natural cap}\index{cap orientation!naturality}
	Let $f \colon V \to M$ and $h \colon N \to M$ be transverse maps with $f$ co-oriented, $V$ a manifold with corners and $M$ and $N$ manifolds without boundary.
	Furthermore, let $g \colon W \to N$ be a map from an oriented manifold with corners that is transverse to the co-oriented pullback $V \times_M N \to N$.
	Then the cap orientation induced on $(V \times_M N) \times_N W$ by pulling back the co-oriented map $V \times_M N \to N$ over $W \to N$
	is the same as the cap orientation obtained by pulling back $V \to M$ by the composite $hg \colon W \to M$.
	In other words, $(V \times_M N) \times_N W = V \times_M W$ as oriented manifolds.
\end{proposition}

\begin{proof}
	Note that $V$ is transverse to $hg \colon W \to M$ by \cref{L: transverse to pullback}, so both expressions are defined.
	It follows directly from \cref{P: pullback functoriality} that the two pullback co-orientations we have described for $(V \times_M N) \times_N W \to W$ agree.
	Therefore, the induced cap orientations agree.
\end{proof}

\subsubsection{Comparing the oriented and co-oriented products}\label{S: c vs o}

Suppose $f \colon V \to M$ and $g \colon W \to M$ are two transverse co-oriented maps from manifolds with corners to a manifold without boundary.
Further, suppose $M$ oriented.
Then we know from the discussion in \cref{S: co-orientations} that there is a bijection between co-orientations of $f$ and orientations of $V$; an orientation of $V$ induces a co-orientation of $f$ and vice versa.
Of course the same is true of $W$ and $g$.
In this scenario, we have two different ways to orient $V \times_M W$, depending on whether we start by thinking of $V$ and $W$ as oriented or by thinking of $f$ and $g$ as co-oriented.
If we think of $V$ and $W $as oriented, we have the fiber product orientation of $V \times_M W$
discussed in \cref{S: orientations}.
Alternatively, if we think of $f$ and $g$ as co-oriented, we can form the fiber product co-orientation of $V \times_M W \to M$ as in \cref{S: co-orientation of pullbacks} and then consider the induced orientation $(V \times_M W)\,\check{\vrule height1.3ex width0pt}$ given the orientation of $M$ (recall \cref{N: hat check}).

Our goal in this section is to compare these two orientations on $V \times_M W$.
To attempt to avoid confusion, let us consider $V$ and $W$ given as oriented, and we will write $V \times_M^o W$\index{$Z$@$\times_M^o$} for the fiber product orientation of \cref{S: orientations} and $V \times_M^c W$\index{$Z$@$\times_M^c$} for $(\hat V \times_M \hat W)\,\check{\vrule height1.3ex width0pt}$.

The reader might have noticed that a third way to orient $V \times_M W$ is to consider $f$ to be co-oriented and $W$ to be oriented and then form the cap orientation that we studied in detail in the preceding section.
However, we already know this to be identical to $V \times_M ^cW$ by \cref{L: same induced}.
By contrast, the following proposition shows this is not always the same as $V \times_M^o W$.

When we move on to geometric homology, $V \times_M ^oW$ will correspond to the classical intersection product of homology classes, as described for example in \cite[Section VI.11]{Bred97}, while $V \times_M^c W$ will correspond to (the Poincar\'e dual of) the cup product of cohomology classes.
When $M$ is closed and oriented, switching between thinking of $V$ and $W$ as oriented vs. co-oriented will be precisely the Poincar\'e duality isomorphism.
So this proposition will ultimately demonstrate that the intersection product is Poincar\'e dual to the cup product, up to a sign; see \cref{S: PD}.

\begin{proposition}\label{P: compare cup and intersection orientations}\index{fiber product!co-oriented vs. oriented}\index{co-orientation!of fiber product!vs. orientation}\index{orientation!of fiber product!vs. co-orientation}
	Let $f \colon V \to M$ and $g \colon W \to M$ be transverse maps from manifolds with corners to an oriented manifold without boundary with $V$, $W$, and $M$ all oriented.
	Then $$V \times_M^o W = (-1)^{(m-v)(m-w)} V \times_M^c W$$ as oriented manifolds with corners.
\end{proposition}

The proof will take a bit of work.
Our strategy will be as follows.
We recall by \cref{pullback,P: interior co-orientation} that co-orientations are determined entirely by what happens on the pullbacks of the interiors, and this is also a standard fact for orientations of topological manifolds with boundary.
So it suffices in what follows to consider only manifolds without boundary.

First, we will prove in \cref{L: compare cup and intersection for immersions} that the result holds when $f$ and $g$ are immersions.
Then we will show, first for co-orientations and then for orientations, that we can replace the fiber product of
$$V \xr{f} M \xleftarrow{g} W$$
with the fiber product
$$V \times \R^b \xhookrightarrow{e \times \id_{\R^b }} M \times \R^a \times \R^b \hookleftarrow W \times \R^a ,$$
where $a$ and $b$ are even, the map $e$ is the embedding $V \into M \times \R^a $ of a Quillen co-orientation of $f$, and the leftward arrow is the identity between the $\R^a $ factors and takes the $W$ factor into $M \times \R^b $ by the embedding map of a Quillen co-orientation of $g$.
By ``replace,'' we mean in the oriented case that the two fiber products are canonically oriented diffeomorphic.
In the co-oriented case we mean that we have a canonical oriented diffeomorphism
between the domains of the two co-oriented fiber products, oriented with their induced orientations coming respectively from the orientation of $M$ and from the concatenation orientation of $M \times \R^a \times \R^b $ using the standard orientations of the Euclidean terms.
The maps of this second fiber product are embeddings for which the proposition holds by \cref{L: compare cup and intersection for immersions}, and so the general case will follow.
The even dimensions of the Euclidean factors are chosen to avoid some extraneous signs in the arguments below.

Before proceeding, let us explain in more detail what we mean by ``canonical'' here and below.
Recall that, as a space, $P = V \times_M W$ can be identified with $\{(v,w) \in V \times W \mid f(v) = g(w)\}$.
Below we will see various fancier embeddings of $P$ in spaces of the form $V \times X \times W \times Y$, with $X$ and $Y$ Euclidean or $I$.
For each such embedding, the projection to $V \times W$ will take the embedding of $P$ back to $P$.
In this way, all versions of $P$ can be canonically identified, and it is these identifications that will yield our orientation preserving diffeomorphisms.
Such identifications have already been discussed in \cref{R: pullback representative,R: pullback representative 2}.
In the latter remark, we provide exactly such a canonical identification between our standard realization of $V \times_M W$ as a subset of $V \times W$ and the version used for co-orienting pullbacks and fiber products.

\medskip

We begin with the case where $V$ and $W$ are immersions and then work toward the general case.

\begin{lemma}\label{L: compare cup and intersection for immersions}\index{fiber product!co-oriented vs. oriented!for immersions}
	If $f$ and $g$ are transverse immersions of oriented manifolds without boundary into an oriented manifold without boundary, then $V \times_M ^oW = (-1)^{(m-v)(m-w)} V \times_M ^cW$.
\end{lemma}

\begin{proof}
	It suffices to consider small neighborhoods of interior points on which $f$ and $g$ are embeddings.
	Let $\beta_V$, $\beta_W$, and $\beta_M$ denote the local orientations of $V$, $W$, and $M$, respectively, at such a point.
	Assuming $f$ and $g$ to be co-oriented with the induced co-orientations $(\beta_V,\beta_M)$ and $(\beta_W,\beta_M)$, we have the resulting Quillen orientations $\beta_{\nu V}$ and $\beta_{\nu W}$ of the normal bundles of $V$ and $W$.
	Recall that these are defined so that $(\beta_V,\beta_V \wedge \beta_{\nu V})$ and $(\beta_W,\beta_W \wedge \beta_{\nu W})$ are the co-orientations of $f$ and $g$, respectively.
	In this scenario, with $\beta_V$ and $\beta_W$ fixed as the orientations of $V$ and $W$, this is equivalent to requiring $\beta_V \wedge \beta_{\nu V} = \beta_M$ and $\beta_W \wedge \beta_{\nu W} = \beta_M$.
	Again, to keep the contexts clear, for the remainder of the argument we will write $\beta^c_{\nu V}$ and $\beta^c_{\nu W}$ for the Quillen orientations of $\nu V$ and $\nu W$.

	By \cref{P: normal pullback}, the co-orientation of the fiber product $\hat V \times_M \hat W$ is $(\beta_P,\beta_P \wedge \beta^c_{\nu V} \wedge \beta^c_{\nu W})$.
	So, taking the induced orientation, $V\times_M^c W$ has the orientation $\beta_P^c$ such that $\beta_P^c \wedge \beta^c_{\nu V} \wedge \beta^c_{\nu W} = \beta_M$.

	On the other hand, in \cref{P: orient intersection}, $\beta^o_P$ is such that if $\beta^o_P \wedge \beta^o_{\nu W} = \beta_V$ and $\beta^o_P \wedge \beta^o_{\nu V} = \beta_W$ then $\beta^o_P \wedge \beta^o_{\nu V} \wedge \beta^o_{\nu W} = \beta_M.$ A priori these may be different orientations of $\nu V$ and $\nu W$ than those of the preceding paragraph, so we use these alternate labels.
	In fact, let us suppose $\beta^o_P$, $\beta^o_{\nu V}$, and $\beta^o_{\nu W}$ chosen so that these expressions all hold.
	Then we have
	$$\beta_M = \beta^o_P \wedge \beta^o_{\nu V} \wedge \beta^o_{\nu W} = \beta_{W} \wedge \beta^o_{\nu W},$$
	so
	$$\beta_{\nu W}^o = \beta_{\nu W}^c.$$
	Similarly,
	$$\beta_M = \beta^o_P \wedge \beta^o_{\nu V} \wedge \beta^o_{\nu W} = (-1)^{(m-v)(m-w)}\beta^o_P \wedge \beta^o_{\nu W} \wedge \beta^o_{\nu V} = (-1)^{(m-v)(m-w)}\beta_V \wedge \beta^o_{\nu V},$$
	so
	$$\beta^o_{\nu V} = (-1)^{(m-v)(m-w)}\beta^c_{\nu V}.$$

	Thus
	$$\beta_M = \beta^o_P \wedge \beta^o_{\nu V} \wedge \beta^o_{\nu W} = (-1)^{(m-v)(m-w)}\beta^o_P \wedge \beta^c_{\nu V} \wedge \beta^c_{\nu W}.$$

	We conclude that $\beta^c_P = (-1)^{(m-v)(m-w)}\beta^o_P$.
\end{proof}

Now we show how to replace a general co-oriented fiber product with a co-oriented fiber product whose maps are embeddings.
In our remaining constructions in this section we take all introduced Euclidean spaces to be even-dimensional to simplify the signs in our computations.

\begin{lemma}
	Let $f \colon V \to M$ and $g \colon W \to M$ be transverse co-oriented maps from manifolds without boundary to an oriented manifold without boundary.
	Let $V\xhookrightarrow{e}M \times \R^a \to M$ be a Quillen co-orientation for $f$ with $a$ even.

	Then $(V \times_M W)\,\check{\vrule height1.3ex width0pt}$ (with the orientation induced by the fiber product co-orientation and the orientation of $M$) is canonically oriented diffeomorphic to $(V\times_{M \times \R^a }(W \times \R^a ))\,\check{\vrule height1.3ex width0pt}$ (with the orientation induced by the fiber product co-orientation and the orientation of $M \times \R^a $).
	Here the maps for the second fiber product are $e \colon V \to M \times \R^a $ and $g \times \id_{\R^a } \colon W \times \R^a \to M \times \R^a $.
	As in the construction of the Quillen co-orientation (\cref{D: Quillen normal or}), we assume $e \colon V \into M \times \R^a $ to be co-oriented so that its composition with the canonical co-orientation $(\beta_{M}\wedge\beta_a,\beta_M)$ is the co-orientation of $f$.
	We also take $g \times \id_{\R^a }$ to be co-oriented by the product co-orientation $(\beta_W \wedge \beta_a,\beta_M \wedge \beta_a)$ if $(\beta_W,\beta_M)$ is the co-orientation of $g$.
	Finally, $M \times \R^a $ is given the product orientation with $\R^a $ having the standard orientation.
\end{lemma}

\begin{proof}
	By the definitions of the induced orientation and the fiber product co-orientation, the co-orientation of the fiber product is obtained by identifying $V \times_M W$ with $(g\times \id_{\R^a})^{-1}(e(V)) \subset W \times \R^a$, and then induced orientation is $\beta_P$, where $\beta_P \wedge \beta_{\nu V} = \beta_W \wedge \beta_a$.
	We recall that the $\nu V$ in this formula is actually the pullback of the normal bundle of $e(V) \subset M \times \R^a $ via the map $g \times \id_{\R^a }: W \times \R^a \to M \times \R^a $.
	But this is also exactly the description of the induced fiber product orientation from the co-oriented fiber product $V\times_{M \times \R^a }(W \times \R^a )$, treating $e \colon V \into M \times \R^a $ as its own Quillen co-orientation with $N = 0$ in \cref{D: pullback coorient}.
	In fact, the co-orientation of \cref{D: pullback coorient} is obtained precisely by identifying these two forms of the pullback.
\end{proof}

\begin{corollary}\label{C: co-oriented full transition to embedded}
	Let $f \colon V \to M$ and $g \colon W \to M$ be transverse co-oriented maps from manifolds without boundary to an oriented manifold without boundary.
	Let $V \xhookrightarrow{r} M \times \R^a \to M$ and $W\xhookrightarrow{s} M \times \R^b \to M$ be Quillen co-orientations compatible with $f$ and $g$.
	Then the fiber product $(V \times_M W)\,\check{\vrule height1.3ex width0pt}$ is canonically oriented diffeomorphic to the fiber product $((V \times \R^b ) \times_{M \times \R^a \times \R^b }(W \times \R^a))\,\check{\vrule height1.3ex width0pt}$, in which the first map is $r \times \id_{\R^b }$ and the second map takes $(w,z) \in W \times \R^a $ to $s(w)$ in the first and third coordinates and $z$ in the second coordinate.
\end{corollary}

\begin{proof}
	By the preceding lemma we have $(V \times_M W)\,\check{\vrule height1.3ex width0pt}$ canonically oriented diffeomorphic to $(V\times_{M \times \R^a}(W \times \R^a))\,\check{\vrule height1.3ex width0pt}$ with $V \to M \times \R^a $ an embedding.
	By \cref{P: graded comm}, the transposition map from this space to $((W \times \R^a ) \times_{M \times \R^a }V)\,\check{\vrule height1.3ex width0pt}$ is $(-1)^{(m-v)(m-w)}$-orientation preserving, using that all our Euclidean spaces are taken even-dimensional.
	We next observe the map described for $W \times \R^a \to M \times \R^a \times \R^b $ is an embedding whose composition with the projection to $M \times \R^a $ gives a Quillen co-orientation for $g \times \id_{\R^a }$.
	Now we can apply the lemma again and then transpose again with the same sign, so that the signs cancel out.
\end{proof}

The preceding lemma and corollary concerned orientations obtained from co-orientations.
Next we consider oriented manifolds and maps that are not necessarily co-oriented and show how to replace the oriented fiber product with an oriented fiber product whose maps are embeddings.
Again, we take all introduced Euclidean spaces to be even-dimensional to simplify the signs.

\begin{lemma}
	Let $f \colon V \to M$ and $g \colon W \to M$ be transverse maps from oriented manifolds without boundary to an oriented manifold without boundary.
	Let $V \xhookrightarrow{e}M \times \R^a \to M$ be a factorization of $f$ with $e$ an embedding.
	Then $V \times_M W$ is canonically oriented diffeomorphic to the oriented fiber product $V\times_{M \times \R^a } (W \times \R^a )$ of $e \colon V \to M \times \R^a $ and $g \times \id_{\R^a } \colon W \times \R^a \to M \times \R^a $.
\end{lemma}

\begin{proof}
	We will first show that $V \times_M W$ is canonically oriented diffeomorphic to the fiber product $V\times_{M \times \R^a } (W \times \R^a )$ with the map $V \to M \times \R^a $ being the composition of $f$ with the inclusion $M = M \times \{0\} \into M \times \R^a $.
	We will write this composite as $f_0$.
	The maps $f_0$ and $g \times \id_{\R^a }$ are transverse as $f$ and $g$ are transverse in $M$ and the $\id_{\R^a }$ factor takes care of the $\R^a $ factor of the tangent spaces.
	We also observe that the two fiber products are canonically the same as spaces, as $V \times_M W = \{(v,w) \in V \times W \mid f(v) = g(w)\}$, while the other fiber product is $\{(v,w,0) \in V \times W \times \R^a \mid f(v) = g(w)\}$.

	Now we consider the orientations.
	Let us choose a point $p = (v,w) \in P = V \times_M W$ and let $z=f(v)=g(w)\in M$.
	As the tangent space of the pullback is the pullback of the tangent spaces by \cref{L: tangent of pullbacks}, we have $T_pP = T_vV \times_{T_z(M)} T_wW \subset T_vV \oplus T_wW$.
	By the definition of the fiber product orientation in \cref{S: orientation of fiber products}, we consider the map
	\begin{equation*}
		\Phi \colon T_pP \oplus T_zM \to T_vV \oplus T_wW
	\end{equation*}
	given by $\Phi((a,b),c) =(a,b)+s(c)$, where $s$ is a splitting of the map $T_vV \oplus T_wW \to T_zM$ that takes $(a,b)$ to $Df(a)-Dg(b)$.
	We recall from \cref{S: orientation of fiber products} that the orientation of $P$ is chosen so that $\Phi$ is an orientation preserving isomorphism up to a sign of $(-1)^{wm}$.

	In the case of $V\times_{M \times \R^a } (W \times \R^a )$, if we write $P'$ for the fiber product,
	we have $T_{(p,0)}P' = T_vV \times_{T_z(M)} T_{(w,0)}(W \oplus \R^a) \subset T_vV \oplus T_{(w,0)}(W \oplus \R^a)$.
	As $Df_0 = (Df,0)$, we have $T_{(p,0)}P' = \{(a,b,0) \in T_vV \oplus T_{(w,0)}(W \oplus \R^a) \mid Df(a)=Dg(b)\}$, so $T_{(p,0)}P'$ is canonically isomorphic to $T_pP = \{(a,b) \in T_vV \oplus T_wW \mid Df(a)=Dg(b)\}$.
	In this case, we consider the map
	\begin{equation}\label{E: fiber plus euclidean}
		\Psi: T_{(p,0)}P' \oplus T_zM \oplus T_0\R^a \to T_vV \oplus T_wW \oplus T_0\R^a ,
	\end{equation}
	with the restriction of $\Psi$ to $T_{(p,0)}P'$ again being the projection maps, while the restriction to $T_zM \oplus T_0\R^a$ must be a splitting of the map $\Upsilon: T_vV \oplus T_{(w,0)}(W \oplus \R^a) \to T_{(z,0)}(M \oplus \R^a)$ that is $Df_0$ on the first factor and $-D(g \oplus \id_{\R^a })$ on the last two factors.
	We claim that we can take $\Psi((a,b,0),c,u) = (a,b,0)+(s(c),0)-(0,0,u)$, with the splitting map $s$ as above.
	This is certainly correct on the $P$ factor.
	For the $M \oplus \R^a $ factor, we must show $\Upsilon \Psi(0,c,u) = (c,u)$.
	We have $\Psi(0,c,u) = (s(c),0)-(0,0,u)$, noting that $s(c) \in T_vV \oplus T_wW$.
	If we write $s(c) = (s_V(c),s_W(c))$, then by definition $Df(s_V(c))-Dg(s_W(c)) = c$.
	So we have
	\begin{align*}
		\Upsilon \Psi(0,c,u) & =
		\Upsilon((s(c),0)-(0,0,u))\\& =
		\Upsilon(s_V(c),s_W(c),-u)\\& =
		Df_0(s_V(c))-(Dg(s_W(c)),-u)\\& =
		(Df(s_V(c)),0)-(Dg(s_W(c)),-u)\\& =
		(Df(s_V(c))-Dg(s_W(c)),u)\\& =
		(c,u).
	\end{align*}
	So our definition of $\Psi$ suffices for determining the orientation of the fiber product.

	As $\Phi$ and $\Psi$ agree in the first two factors (identifying $T_pP$ and $T_{(p,0)}P'$ in the obvious way) and the dimension of $\R^a $ is even, we see that $\Psi$ is orientation-preserving if and only if $\Phi$ is.
	Furthermore, we have $(-1)^{w(m+a)} = (-1)^{wm}$, so the two fiber product orientations agree in this case.

	Next we must generalize from $f_0 \colon V \to M \times \R^a $ to the general case of $e \colon V \to M \times \R^a $.
	By assumption, $f \colon V \to M$ is the composition of $e$ with the projection $M \times \R^a \to M$, so we may write $e(v) = (f(v), e_{\R}(v))$, and there is a fiberwise homotopy $H \colon V \times I \to M \times \R^a $ from $e$ to $f_0$ given by $H(v,t) = (f(v), te_{\R}(v))$.
	We note that $H$ and $e$ are each transverse to $g \times \id_{\R^a }$.
	Indeed, if $e(v) = (g \times \id_{\R^a })(w,u)$, then $(f(v),e_{\R}(v)) = (g(w),u)$.
	The image of the derivative of $g \times \id_{\R^a }$ at such a point will span $Dg(T_wW) \oplus T(\R^a )$, while the derivative of $e$ will have the form $(Df,De_{\R})$.
	But the image of $D(g \times \id_{\R^a })$ already spans $0\oplus \R^a $, so $De(T_vV)+D(g \times \id_{\R^a })T_{(w,u)} = (Df,0)T_vV+(Dg,0)T_wW+(0,\id_{\R^a })T_u\R^a = T_{(z,u)}(M \times \R^a )$.
	The same argument holds for $H(-,t)$ for any fixed $t$, replacing $De_{\R}$ with $tDe_{\R}$.
	But if each $H(-,t)$ is transverse to $g \times \id_{\R^a }$ then so is $H$.

	It follows that we can form the oriented fiber product of $H$ and $g \times \id_{\R^a }$ over $M \times \R^a $.
	In fact, this fiber product is diffeomorphic to $P \times I$: Noting that we have $H(v,t) = (f(v),te_{\R}(v)) = (g(w),u)$ if and only if $f(v) = g(w)$ and $te_{\R}(v) = u$, we obtain a diffeomorphism $P \times I \to (V \times I)\times_{M \times \R^a }(W \times \R^a )$ given by $((v,w),t) \mapsto ((v,t),(w, te_{\R}(v))$ with inverse $((v,t),(w,u)) \mapsto ((v,w),t)$.
	So this space is a cylinder and the two ends correspond to our two versions of $V\times_{M \times \R^a } (W \times \R^a )$, one mapping $V$ by $f_0$ and the other by $e$.
	Since we have a cylinder, these two end spaces are oriented diffeomorphic, and canonically so by our construction.

	Putting this oriented diffeomorphism together with the one constructed above gives the desired oriented canonical diffeomorphism with the original $V \times_M W$.
\end{proof}

\begin{corollary}\label{C: oriented full transition to embedded}
	Let $f \colon V \to M$ and $g \colon W \to M$ be transverse maps from oriented manifolds without boundary to an oriented manifold without boundary.
	Let $V\xhookrightarrow{r} M \times \R^a \to M$ and $W\xhookrightarrow{s} M \times \R^b \to M$ be factorizations of $f$ and $g$ with $r$ and $s$ embeddings.
	Then the oriented fiber product $V \times_M W$ is canonically oriented diffeomorphic to the fiber product $(V \times \R^b )\times_{M \times \R^a \times \R^b }(W \times \R^a )$, in which the first map is $r \times \id_{\R^b }$ and the second map takes $(w,z) \in W \times \R^a $ to $s(w)$ in the first and third coordinates and $z$ in the second coordinate.
\end{corollary}

\begin{proof}
	As in the proof of \cref{C: co-oriented full transition to embedded},
	we apply the preceding lemma to get $V \times_M W$ canonically oriented diffeomorphic to $V\times_{M \times \R^a }(W \times \R^a )$.
	Then we use the graded commutativity rule for oriented fiber product, as given by \cref{P: commute oriented fiber},
	by which
	the transposition map to $(W \times \R^a )\times_{M \times \R^a }V$ is $(-1)^{(m-v)(m-w)}$-orientation preserving, using that all our Euclidean spaces are taken even-dimensional.
	Then we observe the map described for $W \times \R^a \to M \times \R^a \times \R^b $ is an embedding whose composition with the projection to $M \times \R^a $ is $g \times \id_{\R^b }$.
	Now we apply the lemma again and then transpose again.
\end{proof}

\begin{proof}[Proof of \cref{P: compare cup and intersection orientations}]
	By \cref{C: co-oriented full transition to embedded,C: oriented full transition to embedded}, the proposition reduces to \cref{L: compare cup and intersection for immersions}.
\end{proof}

As a corollary to \cref{P: compare cup and intersection orientations} and some of our previous formulas, we have the following, which will eventually provide associativity between the cap and intersection products.
We continue to use the notation of \cref{P: compare cup and intersection orientations}.

\begin{corollary}[Mixed associativity 2]\label{C: cap/intersect}\index{fiber product!mixed associativity}
	Let $M$ be an oriented manifold without boundary, let $f \colon V \to M$ be a co-oriented map from a manifold with corners, and let $g \colon W \to M$ and $h \colon Z \to M$ be maps from oriented manifolds with corners.
	Then, assuming sufficient transversality for all terms to be defined (see \cref{R: multiproducts}), $$V \times_M (W \times_M^o Z) = (-1)^{(m-v)(m-z)} (V \times_M W) \times_M^o Z,$$
	where in each case the $\times_M$ symbol on the left indicates that the output has the cap orientation.
\end{corollary}

\begin{proof}
	Let $\hat W$ and $\hat Z$ denote $W$ and $Z$ when considered with the induced co-orientations, and let $M$ denote the identity map $M \to M$, thinking of $M$ as oriented.
	We apply several of our previous results in the following computation:
	\[
	\begin{aligned}
		V &\times_M (W \times_M^o Z)
		= (-1)^{(m-w)(m-z)} V \times_M (W \times_M^c Z)
		&& \text{Prop.~\ref{P: compare cup and intersection orientations}} \\
		&= (-1)^{(m-w)(m-z)} V \times_M ((\hat W \times_M \hat Z)\times_M M)
		&& \text{def.\ of $\times_M^c$, Prop.~\ref{P: cap with identity M}} \\
		&= (-1)^{(m-w)(m-z)} (V \times_M (\hat W \times_M \hat Z)) \times_M M
		&& \text{Prop.~\ref{P: OC mixed associativity}} \\
		&= (-1)^{(m-w)(m-z)} ((V \times_M \hat W) \times_M \hat Z)\times_M M
		&& \text{Prop.~\ref{C: fiber assoc}} \\
		&= (-1)^{(m-w)(m-z)+(m-(v+w-m))(m-z)} ((V \times_M \hat W)\times_M M) \times_M^o (\hat Z \times_M M)
		&& \text{Prop.~\ref{P: OC mixed associativity},~\ref{P: cap with identity M}} \\
		&= (-1)^{(m-v)(m-z)} (V \times_M (\hat W\times_M M)) \times_M^o (\hat Z \times_M M)
		&& \text{Prop.~\ref{P: OC mixed associativity}} \\
		&= (-1)^{(m-v)(m-z)} (V \times_M W) \times_M^o Z
		&& \text{Prop.~\ref{P: cap with identity M}.}\qedhere
	\end{aligned}
	\]
\end{proof}

We can similarly consider the exterior products.

Suppose $f \colon V \to M$ and $g \colon W \to N$ with $V$ and $W$ are oriented.
We write $V \times^o W$ for $V \times W$ with the product orientation.
On the other hand, let $V \times^c W = (\hat V \times \hat W)\,\check{\vrule height1.3ex width0pt}$, where the product $\hat V \times \hat W$ on the right is the product co-orientation of \cref{D: co-oriented exterior}.

\begin{proposition}\label{P: compare exterior orientations}\index{exterior product!orientation vs. co-orientation}
	Let $f \colon V \to M$ and $g \colon W \to N$ be maps from oriented manifolds with corners to oriented manifolds with corners.
	Then $$V \times^o W = (-1)^{(m-v)w} V \times^c W$$ as oriented manifolds with corners.
\end{proposition}
\begin{proof}
	Let $\beta_M$ and $\beta_N$ be the orientations of $M$ and $N$, and let $\beta_V$ and $\beta_W$ be the orientations on $V$ and $W$.
	Then the orientation of $V \times^o W$ is $\beta_V \wedge \beta_W$, and the orientation for $M \times^o N$ is $\beta_M \wedge \beta_N$.
	The induced co-orientations of $f$ and $g$ are then, respectively, $(\beta_V, \beta_M)$ and $(\beta_W, \beta_N)$.
	By \cref{D: co-oriented exterior}, the product co-orientation is $(-1)^{(m-v)w}(\beta_V \wedge \beta_W,\beta_M \wedge \beta_N),$ and so the resulting orientation on $V \times^c W$ is $(-1)^{(m-v)w}(\beta_V \wedge \beta_W)$.
\end{proof}

\subsection{Appendix: Lipyanskiy's co-orientations}\label{S: Lipyanskiy co-orientations}\index{co-orientation!Lipyanskiy definition|(}

In \cite{Lipy14}, Lipyanskiy uses a different notion of co-orientation from the one we have used to define geometric cochains.
We here discuss Lipyanskiy's co-orientations, which he initially refers to as \textit{orientations of maps}, and show that for a smooth map $f \colon M \to N$, his definition is equivalent to our definition, up to possible sign conventions.
In other words, we show that a smooth map is co-orientable in our sense if and only if it is co-orientable in Lipyanskiy's sense.
We will not explore the precise differences between the specific co-orientation conventions.
We will also see that this alternative framework is in some sense dual to Quillen's approach to co-orientations that we presented in \cref{S: Quillen}: while Quillen's formulation involves replacing arbitrary maps with embeddings, the formulation here involves replacing arbitrary bundle maps with surjective bundle maps.

To define co-orientations, Lipyanskiy utilizes the determinant line bundles of Donaldson and Kronheimer in \cite[Section 5.2.1]{DoKr90}.
A key point throughout our discussion will be the following lemma, which is \cite[Lemma 5.2.2]{DoKr90}.
Donaldson and Kronheimer actually state (without proof) a stronger version of the lemma---that the isomorphism is canonical---but we will not need that for our purposes.
For the statement of the lemma, recall our definition of $\Or(V)$ in \cref{D: det bundle}.

\begin{lemma}\label{L: det sequence}
	Given an exact sequence of vector bundles
	\[
	0 \to V_1 \xr{d_1} \cdots \xr{d_{m-1}} V_m \to 0,
	\]
	there is an isomorphism
	\[
	\bigotimes_{i\ \text{odd}} \Or(V_i) \ \cong
	\bigotimes_{i\ \text{even}} \Or(V_i).
	\]
\end{lemma}

\begin{proof}
	It is implicit in the hypothesis of exactness that the kernels and image of the maps in the sequence are well-defined vector bundles so that it makes sense to say $\im(d_{i-1}) = \ker(d_{i})$ as objects in the category of vector bundles.
	We thus have for each $i$ a short exact sequence
	\[
	0 \to \ker(d_i) \to V_i \to \im(d_i) \to 0,
	\]
	and since short exact sequences of vector bundles split \cite[Theorem 3.9.6]{Hus75}, we have $V_i \cong \ker(d_i) \oplus \im(d_i) = \im(d_{i-1}) \oplus \im(d_i)$.
	Consequently,
	\[
	\bigotimes_{i\ \text{odd}} \Or(V_i)\, \cong
	\bigotimes_{i\ \text{odd}} \Or(\im(d_{i-1}) \oplus \im(d_i))\, \cong
	\bigotimes_{i\ \text{odd}} \Or(\im(d_{i-1})) \otimes \Or(\im(d_i))\, \cong
	\bigotimes_{\text{all}\ i} \Or(\im(d_{i})),
	\]
	and similarly for the other tensor product.
\end{proof}

We can now define the Donaldson-Kronheimer determinant line bundles as in \cite[Section 5.2.1]{DoKr90}.
Donaldson and Kronheimer work in a more general setting, but we will confine ourselves to considering a map of vector bundles $F \colon E \to E'$ over $M$.
At first, we also assume that $\ker(F)$ and $\cok(F)$ are well-defined vector bundles.
Then the Donaldson-Kronheimer determinant line bundle is defined to be
$$\Or(\ker(F)) \otimes \Or(\cok(F))^*,$$
where the $*$ over $\Or(\cok(F))$ denotes the dual line bundle.
Below we will consider that $\ker(F)$ and $\cok(F)$ are not always vector bundles, but for now we see that the determinant bundle is morally related to the index of an operator.
We refer to \cite[Section 5.2.1]{DoKr90} for a more precise statement of the relationship.

To relate the Donaldson-Kronheimer determinant line bundle to our notion of co-orientation, consider the exact sequence of vector bundles
\begin{equation*}
	0 \to \ker(F) \to E \to E' \to \cok(F) \to 0.
\end{equation*}
Applying \cref{L: det sequence}, we have $\Or(\ker(F)) \otimes \Or(E') \cong \Or(E) \otimes \Or(\cok(F))$.
Next we use that for a line bundle $L$ we have $L \otimes L^* \cong \underline{\R}$, the trivial line bundle.
So multiplying both sides by $\Or(\cok(F))^*$ and $\Or(E')^*$, we get
$$\Or(\ker(F)) \otimes \Or(\cok(F))^* \cong \Or(E) \otimes \Or(E')^*.$$
The latter is isomorphic to $\Hom(\Or(E'), \Or(E))$, which is dual to $\Hom(\Or(E), \Or(E'))$.
In particular, $\Hom(\Or(E), \Or(E'))$ is trivial, and so admits a non-zero section, if and only if the Donaldson-Kronheimer determinant bundle $\Or(\ker(F)) \otimes \Or(\cok(F))^*$ is trivial.

In the setting of a smooth map $f \colon M \to N$, we can think of the derivative $Df$ as a map $Df \colon TM \to f^*(TN)$, and then the above demonstrates that $\Hom(\Or(TM), \Or(f^*(TN)))$ is trivial if and only if the determinant bundle $\Or(\ker(Df)) \otimes \Or(\cok(Df))^*$ is trivial.
We recall that the triviality of
\[
\Hom(\Or(TM), \Or(f^*(TN)))
\]
is the condition for co-orientability of $f$ in the sense of \cref{D: co-orientations}.
Our co-orientations in this setting are equivalence classes of non-zero sections of $\Hom(\Or(TM), \Or(f^*(TN)))$ up to positive scalars or, equivalently, orientations of this line bundle.
Lipyanskiy's co-orientations are orientations of $\Or(\ker(Df)) \otimes \Or(\cok(Df))^*$.
As orientations of line bundles exist if and only if the line bundle is trivial, the two notions of co-orientability coincide.
We leave it to the reader to define the isomorphisms in sufficient detail to carry a particular co-orientation as defined in \cref{S: co-orientations} to one of Lipyanskiy's co-orientations.

The problem with the preceding analysis is that in general $\ker(Df)$ and $\cok(Df)$ do not necessarily have the same dimensions from fiber to fiber, and so $\ker(Df)$ and $\cok(Df)$ are not necessarily well defined as vector bundles.
The solution is to reframe the definition of the determinant line bundle as in \cite{DoKr90} so that it is always well defined and such that it is isomorphic to $\Or(\ker(F)) \otimes \Or(\cok(F))^*$ when it is also well defined.

For this, let $\underline{\R}^n$ be the trivial $\R^n$ bundle over $M$, and suppose we have a map $\psi \colon \underline{\R}^n \to E'$ such that $F \oplus \psi \colon E \oplus \underline{\R}^n \to E'$ is surjective\footnote{Donaldson and Kronheimer work with complex vector bundles, so \cite{DoKr90} features $\underline{\C}^n$ rather than $\underline{\R}^n$.}; here we write $F \oplus \psi$ for the map $(x,y) \to F(x) + \psi(y)$.
Such a map will always exist in our setting, as we defined manifolds with corners to be embedded in finite dimensional Euclidean space.
Hence tangent bundles are subbundles of trivial bundles and so also the images of projections of trivial bundles (or, up isomorphism, quotients of the trivial bundle by their orthogonal complements after endowing the trivial bundle with a Riemannian structure).
The gain is that $F \oplus \psi$ then has trivial cokernel and a kernel that is a vector bundle, as now the fibers of the kernel have a fixed dimension.
We then define the determinant line bundle to be
$$\mathscr L = \Or(\ker(F \oplus \psi)) \otimes \Or(\underline{\R}^n)^* \cong \Or(\ker(F \oplus \psi)).$$

In the case where $\ker(F)$ and $\cok(F)$ were already vector bundles, $\mathscr L$ is isomorphic to the earlier Donaldson-Kronheimer determinant line bundle using \cref{L: det sequence} and the following lemma:

\begin{lemma}
	If $F \colon E \to E'$ and $\psi \colon \underline{\R}^n \to E'$ are bundle maps with $F \oplus \psi \colon E \oplus \underline{\R}^n \to E'$ surjective and $\ker(F)$ and $\cok(F)$ well-defined vector bundles, then the following sequence is exact\footnote{This exact sequence appears incorrectly in \cite{DoKr90} with the $\psi$ in place of $F$ in the first and last terms.}:
	\[
	\begin{tikzcd}[column sep=small]
		0 \arrow[r] & \ker(F) \arrow[r] & \ker(F \oplus \psi) \arrow[r] & \underline{\R}^n \arrow[r] & \cok(F) \arrow[r] & 0.
	\end{tikzcd}
	\]
\end{lemma}

\begin{proof}
	This exact sequence is simply the snake lemma exact sequence obtained from the commutative diagram of exact sequences
	\[
	\begin{tikzcd}
		0 \arrow[r] & E \arrow[r] \arrow[d, "F"] & E \oplus \underline{\R}^n \arrow[r] \arrow[d, "F \otimes \psi"] & \underline{\R}^n \arrow[r] \arrow[d] & 0 \\
		0 \arrow[r] & E' \arrow[r, "="] & E' \arrow[r] & 0 \arrow[r] & 0.
	\end{tikzcd}
	\]
	The category of vector bundles over a space is not technically an abelian category, but one can check by hand for this diagram that, with our assumptions, all the maps of the exact sequence are well defined and the exactness then holds fiberwise by the classical snake lemma.
	In particular, the map $\underline{\R}^n \to \cok(F)$ is the composition of the splitting map $\underline{\R}^n \to E \oplus \underline{\R}^n$, the map $F \oplus \psi$, and the projection $E'$ to $\cok(F)$.
\end{proof}

Combining this lemma with \cref{L: det sequence} gives us an isomorphism
$$\Or(\ker(F)) \otimes \Or(\underline{\R}^n) \cong \Or(\ker(F \oplus \psi)) \otimes \Or(\cok(F)).$$
Multiplying both sides by $\Or(\underline{\R}^n)^* \otimes \Or(\cok(F))^*$ and using again that for a line bundle $L$ we have $L \otimes L^* \cong \underline{\R}$, we obtain
$$\Or(\ker(F)) \otimes \Or(\cok(F))^* \cong \Or(\ker(F \oplus \psi)) \otimes \Or(\underline{\R}^n)^*.$$
So, as promised, the two definitions of the Donaldson-Kronheimer determinant line bundle agree (up to canonical isomorphisms) when $\ker(F)$ and $\cok(F)$ are defined.

Finally, we should observe that the construction of $\mathscr L$ is independent, at least up to isomorphism, of the choice of $\psi$ and $n$.
Clearly $\Or(\underline{\R}^n) \cong \Or(\underline{\R}^n)^* \cong \underline{\R}$ for all $n$, so we must only show that if $\psi_1 \colon \underline{\R}^n \to E'$ and $\psi_2 \colon \underline{\R}^m \to E'$ are two maps satisfying the requirement of the definition then $\Or(\ker(F \oplus \psi_1)) \cong \Or(\ker(F \oplus \psi_2))$.
Adapting an argument in \cite[Section 5.1.3]{DoKr90}, we note that the bundle maps
\begin{align*}
	F \oplus \psi_1 \oplus 0 \colon &E \oplus \underline{\R}^{n+m} \to E'\\
	F \oplus 0 \oplus \psi_2 \colon &E \oplus \underline{\R}^{n+m} \to E'
\end{align*}
are each homotopic through fiberwise linear homotopies to $F \oplus \psi_1 \oplus \psi_2$, and so they are homotopic to each other.
Furthermore, these are homotopies through surjective bundle maps, so we can write the homotopy between the two maps as a surjective bundle map $(E \oplus \underline{\R}^{n+m}) \times I \to E' \times I$ over $M \times I$.
Thus we have a well-defined kernel bundle over $M \times I$, which implies that the kernel bundles over $M \times \{0\}$ and $M \times \{1\}$ are isomorphic.
But, up to reordering the summands, these are, respectively, $\ker(F \oplus \psi_1) \oplus \underline{\R}^m$ and $\ker(F \oplus \psi_2) \oplus \underline{\R}^n$.
Therefore,
$$\Or(\ker(F \oplus \psi_1) \oplus \underline{\R}^m) \cong \Or(\ker(F \oplus \psi_2) \oplus \underline{\R}^n).$$
But
$$\Or(\ker(F \oplus \psi_1) \oplus \underline{\R}^m) \cong \Or(\ker(F \oplus \psi_1)) \otimes \Or(\underline{\R}^m) \cong \Or(\ker(F \oplus \psi_1)),$$
and similarly for the other bundle.
So $\Or(\ker(F \oplus \psi_1)) \cong \Or(\ker(F \oplus \psi_2))$.
\index{co-orientation!Lipyanskiy definition|)} 

\section{Geometric chains and cochains}\label{S: geometric cochains}

Geometric homology and cohomology are homology/cohomology theories for smooth manifolds defined through submanifolds or, more generally, maps from manifolds with corners.
They agree with singular homology and cohomology, but having different representatives
at the chain/cochain level, they provide geometric approaches to both theory and calculations.
They are thus akin to de Rham theory in the sense that chains and cochains are defined through the smooth structure rather than continuous maps or open sets, but they are defined over the integers and not just the real numbers.

Our definitions of geometric chains and cochains will be a modification of those given by Lipyanskiy in \cite{Lipy14}, though his results will continue to hold with our modified definitions.
As Lipyanskiy primarily focuses on geometric chains and geometric homology, our focus, where the theories diverge, will primarily be on geometric cochains and cohomology, for which Lipyanskiy's account is less complete despite the cohomological setting having its own subtleties.
We also take the opportunity to fill in some of the details missing from Lipyanskiy's account more generally, especially utilizing the more thorough foundations on manifolds with corners provided by \cite{Joy12,MaDo92}.

\subsection{Preliminary definitions}

We first identify certain types of ``manifolds over $M$.'' In future sections $M$ will most typically be a manifold without boundary, as that is the case where we obtain agreement between geometric (co)homology and singular (co)homology, but in this section we allow it to be a manifold with corners anywhere that the extra generality may be useful but so long as it does not create more technical work.

\begin{definition}\label{V: maps are co-oriented}
	Let $M$ be a smooth manifold, possibly with corners.
	A \textbf{manifold over $\mathbf{M}$}\index{manifold over a manifold|textbf} is a manifold with corners $W$ with a smooth map $r_W \colon W \to M$,
	called the \textbf{reference map}\index{reference map|textbf}\index{$R$@$r_W$}.
	We freely and almost always abuse notation
	by using the domain $W$ to refer to the manifold over $M$, not $r_W$ or some other symbol, letting
	context determine whether we are referring to the entire data or the domain; see \cref{N: implicit notation}.

	We say a manifold over $M$, is:
	\begin{itemize}
		\item \textbf{compact}\index{manifold over a manifold!compact|textbf} if its domain is compact,
		\item \textbf{proper}\index{manifold over a manifold!proper|textbf} if its reference map is proper (i.e. the preimage of compact subsets is compact),
		\item \textbf{oriented}\index{manifold over a manifold!oriented|textbf} if its domain is oriented, and
		\item \textbf{co-oriented}\index{manifold over a manifold!co-oriented|textbf} if its reference map is co-oriented.
	\end{itemize}
	If $W$ is a manifold over $M$ then so is $\bd W$ using the reference map $r_W \circ i_{\bd W} \colon \bd W \to M$.
	If $W$ is oriented or co-oriented, then $\bd W$ inherits an orientation or co-orientation as in \cref{Con: oriented boundary} or \cref{D: boundary co-orientation}.
	By \cite[Lemma 2.8]{Joy12}, $i_{\bd W}$ is proper, so, as the composition of proper maps is proper, if $W$ is a proper manifold over $M$ then so is $\bd W$.

	If $W$ is oriented or co-oriented, we write $-W$ to refer to the manifold over $M$ with the opposite orientation or co-orientation; it should always be clear from context which structure we are replacing with its opposite.

	If $V$ and $W$ are manifolds over $M$ that are disjoint\footnote{Recall that all of our manifolds with corners are subsets of $\R^\infty$.} as submanifolds of $\R^\infty$, we let the disjoint union $V \sqcup W$ be a manifold over $M$ with reference map $r_{V \sqcup W}$ such that $r_{V \sqcup W}(x) = r_V(x)$ when $x \in V$ and $r_{V \sqcup W}(x) = r_W(x)$ when $x \in W$.
	We call this construction the \textbf{disjoint union}\index{disjoint union|textbf}\index{$Z$@$\sqcup$} of $V$ and $W$.
	If $V$ and $W$ are both oriented, co-oriented, compact, or proper, then so is $V \sqcup W$.
	In \cref{E: copies,D: prechain sum} we consider disjoint union for the case where $V$ and $W$ are not necessarily disjoint in $\R^\infty$ using the notion of isomorphism classes, which we introduce next.
\end{definition}


Geometric chains and cochains will be equivalence classes of manifolds over $M$ under an equivalence relation we define using the following concepts, which are taken from or modify
the definitions of \cite{Lipy14}.

\begin{definition}\label{D: equiv triv and small}
	Let $V, W$ be manifolds over $M$ with reference maps $r_V$ and $r_W$.
	\begin{itemize}
		\item Suppose $V$ and $W$ are co-oriented and $\phi \colon W \to V$ is a diffeomorphism such that $r_V \circ \phi = r_W$.
		We say that $\phi$ is \textbf{co-orientation preserving (or co-orientation reversing)}\index{co-orientation preserving|textbf}\index{co-orientation reversing|textbf} if the composition of the tautological co-orientation of $\phi$ (\cref{D: tautological co-orientation}) with the co-orientation of $r_V$ agrees with (or disagrees with) the co-orientation of $r_W$.
		\item If $V$ and $W$ are oriented, we say they are \textbf{(oriented) isomorphic}\index{manifold over a manifold!isomorphic|textbf} if there is an orientation-preserving diffeomorphism $\phi \colon W \to V$ such that $r_V \circ \phi = r_W$.
		\item If $V$ and $W$ are co-oriented, we say they are \textbf{(co-oriented) isomorphic}\index{manifold over a manifold!isomorphic|textbf} if there is a co-orientation-preserving diffeomorphism $\phi \colon W \to V$.
		\item If $W$ is oriented then $W$ is \textbf{(oriented) trivial}\index{manifold over a manifold!trivial|textbf}\index{trivial|textbf} if there is an orientation-reversing diffeomorphism $\rho \colon W \to W$ such that $r_W \circ \rho = r_W$.
		\item If $W$ is co-oriented then it is \textbf{(co-oriented) trivial}\index{manifold over a manifold!trivial|textbf}\index{trivial|textbf} if there is a co-orientation-reversing diffeomorphism $\rho \colon W \to W$.
		\item $W$ has \textbf{small rank}\index{manifold over a manifold!small rank|textbf}\index{small rank|textbf} if the differential $D r_W$ is everywhere less than full rank, in other words if the rank of $D_w r_W$ is less than $\dim(W)$ for all $w\in W$.
		\item If $W$ is (co\nobreakdash-)oriented, then it is \textbf{((co\nobreakdash-)oriented) degenerate}\index{manifold over a manifold!degenerate|textbf}\index{degenerate|textbf} if it has small rank and ${\bd W}$ is the disjoint union of a trivial (co\nobreakdash-)oriented manifold over $M$ and one with small rank.
	\end{itemize}
	We declare the empty map $\emptyset \to M$ to be both trivial and degenerate.
\end{definition}

Note that while we can speak of a diffeomorphism being orientation-preserving or -reversing without it necessarily commuting with the reference maps, the notion of being co-orientation-preserving or -reversing is only defined when $r_V \circ \phi = r_W$.

We observe that isomorphism is an equivalence relation and that it preserves triviality and degeneracy.

Rather than small rank, Lipyanskiy uses the condition of \textbf{small image}\index{small image}.
In our notation, a map $r_W \colon W \to M$ has small image if there is another map $r_T \colon T \to M$ such that $r_T(T)\supset r_W(W)$ but $\dim(T)<\dim(W)$.
However, the small rank condition turns out to be more manageable for our purposes while still providing geometric homology and cohomology theories that are equivalent to singular homology and cohomology.
Roughly speaking, the geometric chains and cochains of a manifold $M$ will consist of isomorphism classes of oriented or proper co-oriented manifolds with corners over $M$ modulo the trivial and degenerate chains or cochains.

The most obvious use of the notion of triviality will be to ensure that $\bd^2 = 0$ so that we have a chain complex.
The degeneracy comes into play in ensuring that geometric homology and cohomology satisfy the dimension axiom; see \cref{E: dimension} and \cref{R: degen1}.

\begin{example}
	Let $S^1$ be the unit circle in the plane with the standard counterclockwise orientation.
	Let $\pi \colon S^1 \to \R$ be the projection of the circle onto the $x$-axis.
	This map is trivial.
	Indeed, if $\rho$ is the reflection of the circle across the $x$-axis, then $\rho$ is orientation-reversing and $\pi \rho = \pi$.

	Alternatively, let $\pi \colon S^1 \to \R$ be co-oriented by $(e_\theta,e_x)$, where $e_x$ is the standard positively-directed unit vector in $\R$ and $e_\theta$ is the counterclockwise tangent vector in $S^1$.
	This map is trivial as a co-oriented map, again via the reflection $\rho$ across the $x$-axis.
	Indeed, we still have $\pi \rho = \pi$, and the tautological co-orientation of $\rho$ is $(e_\theta,-e_\theta)$ so that the composite co-orientation of $\pi\rho$ is $-(e_\theta,e_x)$.
\end{example}

\begin{example}\label{E: copies}
	If $r_V \colon V \to M$ is any oriented or co-oriented map and $-r_V \colon V \to M$ is the same map with the opposite orientation or co-orientation, then $r_V \sqcup -r_V \colon V \sqcup V \to M$ is trivial, taking $\rho$ to be the map that switches the two copies of $V$.

	To be technically accurate, the space $V \sqcup V$ is not clearly defined as a manifold with corners, as we have defined our manifolds with corners in \cref{D: MWC} to be submanifolds of some $\R^N$.
	We solve this problem below by working with isomorphism classes of maps.
	In this case we could take $V \sqcup V$ up to isomorphism to be $V \sqcup W$ where $V$ and $W$ are isomorphic as manifolds over $M$.
	See \cref{D: prechain sum} below for precise details.
\end{example}

\begin{example}
	Any co-oriented map of the interval to a point has small rank, but its boundary does not have small rank; it is nonetheless (co-oriented) degenerate because its boundary is trivial.
\end{example}

\begin{example}\label{E: projected triangle}
	Let $V$ be the $2$-simplex in $\R^2$ with vertices at $(1,0)$, $(-1,0)$, and $(0,1)$, and let $\pi \colon V \to \R$ be the projection to the $x$-axis.
	This map has small rank, but the boundary does not have small rank and is not trivial.
\end{example}

\begin{example}
	Let $V = W = M = \R^1$.
	Let $r_W \colon W \to M$ be the identity map of $\R^1$ with its tautological co-orientation, which we can write $(e_1,e_1)$, letting $e_1$ be a positively-oriented tangent vector to $\R^1$.
	Let $r_V \colon V \to M$ be given by $r_V(t) = -t$ with co-orientation $(-e_1,e_1)$.
	Let $\phi \colon W \to V$ be given by $\phi(t) = -t$.
	Then the tautological co-orientation of $\phi$ is $(e_1,-e_1)$.
	So $r_V\phi = r_W$ as co-oriented maps, and $\phi$ provides a co-oriented isomorphism between $r_W$ and $r_V$ even though they are very different maps.
\end{example}

The notions of orientation-preserving and -reversing diffeomorphisms are fairly standard, and oriented isomorphisms and triviality of oriented manifolds over $M$ are simply such diffeomorphisms that commute with the reference maps.
By contrast, the corresponding notions for co-orientation are less familiar, but the following statement shows how the situation for co-orientations is locally similar to that for orientations.

\begin{lemma}\label{L: co-or preserving/reversing}\index{manifold over a manifold!isomorphism!local (co-)orientation criteria}
	Let $r_V \colon V \to M$ and $r_W \colon W \to M$ be co-oriented manifolds over $M$, and let $\phi \colon W \to V$ be a diffeomorphism such that $r_W = r_V \phi$.
	For $x \in W$, let $(\beta_{W,x},\beta_{M,r_W(x)})$ be a local representation of the co-orientation of $r_W$ and let $(\beta_{V, \phi(x)},\beta_{M,r_W(x)})$ be a local representation of the co-orientation of $r_V$, noting $r_V\phi(x)=r_W(x)$.
	Then $\phi$ provides a co-oriented isomorphism if and only if for all $x\in W$ the derivative $D_x\phi$ takes the orientation $\beta_{W,x}$ of $T_xW$ to the orientation $\beta_{V, \phi(x)}$ of $T_{\phi(x)}V$.
	Similarly, the map $\phi$ is a co-orientation-reversing diffeomorphism if and only if for all $x\in W$ the derivative $D_x\phi$ takes the orientation $\beta_{W,x}$ of $T_xW$ to the orientation $-\beta_{V, \phi(x)}$ of $T_{\phi(x)}V$.
\end{lemma}

In the statement of the lemma, by $D_x\phi$ taking an orientation of $T_xW$ to an orientation of $T_{\phi(x)}V$, we mean that $D_x\phi$ takes an ordered basis consistent with the given orientation of $T_xW$ to an ordered basis consistent with the given orientation of $T_{\phi(x)}V$.
We could equivalently phrase this in terms of the map $\textstyle{\bigwedge^w}D_x\phi$ as in \cref{D: tautological co-orientation}.

\begin{proof}[Proof of \cref{L: co-or preserving/reversing}]
	By definition, $\phi$ is an isomorphism if and only if $\omega_{r_W} = \omega_{\phi} * \omega_{r_V}$, where $\omega_\phi$ is the tautological co-orientation of $\phi$.
	This will be the case if and only if for each $x$ we have $\omega_\phi = (\beta_{W,x}, \beta_{V, \phi(x)})$.
	But by the definition of the tautological co-orientation, if $e_1 \wedge \cdots \wedge e_w$ represents an oriented basis for $T_xW$, then $\omega_\phi$ is represented at $x$ by the pair $(e_1 \wedge \cdots \wedge e_w, D_x\phi(e_1) \wedge \cdots \wedge D_x\phi(e_w))$.
	So $\phi$ is an isomorphism if and only if $D_x\phi$ takes an ordered basis representing the local orientation $\beta_{W,x}$ to an ordered basis representing the local orientation $\beta_{V, \phi(x)}$.

	The orientation-reversing case is equivalent, considering instead the requirement $\omega_{r_W} = - \omega_{\phi} * \omega_{r_V}$.
\end{proof}

The following corollaries are now immediate from \cref{L: co-or preserving/reversing} and the definitions.

\begin{corollary}\label{C: co-or preserving is or preserving}\index{co-orientation!reversing/preserving vs. orientation reversing/preserving}\index{orientation!reversing/preserving vs. co-orientation reversing/preserving}
	Let $r_V \colon V \to M$ and $r_W \colon W \to M$ be co-oriented manifolds over the oriented manifold $M$, and let $\phi \colon W \to V$ be a diffeomorphism such that $r_W = r_V \phi$.
	Then $\phi$ is co-orientation preserving (reversing) if and only if it is orientation preserving (reversing) with respect to the orientations of $V$ and $W$ induced by the co-orientations and the orientation of $M$ (\cref{S: co-orientations}).
\end{corollary}

\begin{corollary}\label{C: induced triviality}\index{triviality!oriented vs. co-oriented}
	Let $r_V \colon V \to M$ be a map from a manifold with corners to an oriented manifold $M$, and let $\check V$ denote $V$ with the induced orientation; equivalently, we could start with an oriented $\check V$ over $M$ and write $V$ when considering the induced co-orientation.
	Then $V$ is co-oriented trivial if and only if $\check V$ is oriented trivial.
\end{corollary}

\subsection{Geometric homology and cohomology}

We can now move toward defining geometric homology and cohomology, starting with ``prechains'' and ``precochains.''

\begin{definition}\label{D: PC}
	Let $M$ be a smooth manifold with corners.
	Denote by $PC^\Gamma_*(M)$\index{$PC^\Gamma_*(M)$} the set of oriented isomorphism classes of compact oriented manifolds over $M$, $r_W \colon W \to M$,
	graded by dimension $\dim(W)$.
	Denote by $PC_\Gamma^*(M)$\index{$PC_\Gamma^*(M)$} the set of co-oriented isomorphism classes of proper co-oriented manifolds over $M$, $r_W \colon W \to M$,
	graded by \textbf{codimension}\index{codimension} $\dim(M) - \dim(W)$.
	We declare the empty manifold of each dimension to be orientable and the empty maps from the empty manifolds to $M$ to be co-orientable.

	As per \cref{V: maps are co-oriented}, we will often write $W \in PC^\Gamma_*(M)$ or $W \in PC_\Gamma^*(M)$, letting $W$ represent both its reference map and its isomorphism class as a manifold over $M$.
	In these respective cases we write $-W$ for $W$ with the opposite orientation or co-orientation.

	We refer to elements of $PC^\Gamma_*(M)$ or $PC_\Gamma^*(M)$ as \textbf{prechains}\index{prechain|textbf} or \textbf{precochains}\index{precochain|textbf}, respectively.
\end{definition}

When necessary, as in the following definition, we may write prechains or precochains as $[W]$ to emphasize that they are isomorphism classes.

\begin{definition}\label{D: prechain sum}\index{prechain!disjoint union}\index{precochain!disjoint union}
	If $V,W$ are compact oriented manifolds over $M$ (respectively proper co-oriented manifolds over $M$), then define $[V] \sqcup [W] = [V' \sqcup W']$, where on the right $\sqcup$ denotes disjoint union and $V'$ and $W'$ are compact oriented manifolds over $M$ (respectively proper co-oriented manifolds over $M$) such $[V]=[V']$, $[W]=[W']$, and $V'$ and $W'$ are disjoint in $\R^\infty$ (cf.
	\cref{D: MWC}); note that the images of $V'$ and $W'$ are not necessarily disjoint in $M$.
\end{definition}

The last clause in the definition is due to our requirement that all manifolds with corners be subsets of $\R^\infty$.
As given, the definition allows constructions like $[W]\sqcup[W]$ to be well defined.
We observe that this definition is itself well defined, since if $V''$ and $W''$ are two other disjoint manifolds over $M$ isomorphic to $V$ and $W$ (in the appropriate sense), then we can compose the diffeomorphisms $V'\xleftarrow{\phi'}V \xr{\phi''} V''$ and $W'\xleftarrow{\psi'}W \xr{\psi''} W''$ to obtain an isomorphism $\phi''(\phi')^{-1} \sqcup \psi''(\psi')^{-1} \colon V' \sqcup W' \to V'' \sqcup W''$.

With the operation $\sqcup$, the sets $PC^\Gamma_*(M)$ and $PC_\Gamma^*(M)$ become commutative monoids in each degree with the empty maps $r_\emptyset \colon \emptyset \to M$ as the identities.

We now return to denoting isomorphism classes by their representatives, noting again that triviality and small rank are properties of the isomorphism classes.
In particular, when we write $V \sqcup W$, we assume that $V$ and $W$ are disjoint in $\R^\infty$.

\begin{definition}\label{D: Q}
	Let $Q_*(M) \subset PC^\Gamma_*(M)$\index{$Q_*(M)$} denote the set of (isomorphism classes of) compact oriented manifolds over $M$ of the form $V \sqcup W$ with $V$ trivial and $W$ degenerate.
	Let $Q^*(M) \subset PC_\Gamma^*(M)$\index{$Q^*(M)$} denote the set of (isomorphism classes of) proper co-oriented manifolds over $M$ of the form $V \sqcup W$ with $V$ trivial and $W$ degenerate.
	In either case $V$ or $W$ may be empty.

	We will sometimes write $W \in Q(M)$\index{$Q(M)$} to mean $W \in Q_*(M)$ or $W \in Q^*(M)$ for arguments that are analogous in the two cases.
	When we do so, we assume a consistent choice of $Q_*(M)$ or $Q^*(M)$ throughout the discussion; see for instance \cref{L: bd defined} and its proof.
\end{definition}

Our first lemma is immediate:

\begin{lemma}\label{L: sum of trivial/degenerate}\index{disjoint union!of trivial is trivial}\index{disjoint union!of degenerate is degenerate}\index{trivial!disjoint union of trivial is trivial}\index{degenerate!disjoint union of degenerate is degenerate}
	\begin{enumerate}
		\item If $V$ and $W$ are trivial, then $V \sqcup W$ is trivial.
		\item If $V$ and $W$ are degenerate, then $V \sqcup W$ is degenerate.
		\item If $V,W\in Q_*(M)$ (or $Q^*(M)$), then $V \sqcup W \in Q_*(M)$ (or $Q^*(M)$).
	\end{enumerate}
\end{lemma}

We also have the following directly from \cref{C: induced triviality} and the fact that the small rank property does not depend on (co-)orientation information.

\begin{lemma}\label{L: Q switch}
	Suppose $V \in PC^*_\Gamma(M)$ with $V$ compact and $M$ oriented, and let $\check V \in PC_*^\Gamma(M)$ denote $V$ with the induced orientation.
	Then $V \in Q^*(M)$ if and only if $\check V \in Q_*(M)$.
\end{lemma}

The following useful basic properties are proven in \cite{Lipy14}; for completeness we provide versions of the arguments here, occasionally augmenting those of \cite{Lipy14}.
In each case, ``isomorphic,'' ``trivial,'' or ``degenerate'' should be read consistently to refer to the compact oriented case or the proper co-oriented case.
Lipyanskiy's proofs assume small image rather than small rank, but the proofs are equivalent.

Our first lemma is contained within the proof of \cite[Lemma 10]{Lipy14}, though it will be useful to have an independent statement.

\begin{lemma}\label{L: trivial structure}
	Suppose $W$ is trivial.
	Then
	$$W = \left(\bigsqcup_a A_a \right) \sqcup \left(\bigsqcup_b B_b \sqcup -B_b \right),$$
	where each $A_a$ and $B_b$ is connected and each $A_i$ is trivial.
\end{lemma}
\begin{proof}
	We can write $W$ as the disjoint union of a (possibly infinite) number of isomorphism classes of connected components and then group the isomorphism classes together up to (co\nobreakdash-)orientation as $W = W_1 \sqcup W_2 \sqcup \cdots$ so that all connected components of each $W_i$ are isomorphic (ignoring (co\nobreakdash-)orientations) for each $i$ and so that no connected component of $W_i$ is isomorphic to a connected component of $W_j$ for $i\neq j$.
	As either $W$ is compact or $r_W$ is proper, each $W_i$ is the union of a finite number of connected components; if not, then there would be a $W_i$ with an infinite number of connected components and, since the components of $W_i$ are isomorphic to each other, a point $x\in M$ such that $r_W^{-1}$ has infinite components, violating that $r_W$ is proper.

	Any diffeomorphism of $W$ preserves the decomposition into $W_i$.
	So because $W$ is trivial, either $W_i$ consists of connected components each isomorphic to some $A$  that possesses a (co\nobreakdash-)orientation reversing diffeomorphism, and so $A$ is trivial, or $W_i$ has zero components when counting with (co\nobreakdash-)orientation, in which case the components of $W_i$ can be paired up into $B \sqcup -B$ pairs.
\end{proof}

\begin{lemma}[Lipyanskiy Lemma 10]\label{L: Lip L10}\index{cancellation property of triviality}\index{trivial!cancellation property}
	If $V$ is trivial and $V \sqcup W$ is trivial, then $W$ is trivial.
\end{lemma}
\begin{proof}
	As in the proof of \cref{L: trivial structure}, as $V$ is trivial, it can be decomposed as $V = V_1 \sqcup V_2 \sqcup \cdots$ where each $V_i$ is a union of isomorphism classes of connected components up to (co\nobreakdash-)orienta\-tion and such that either $V_i$ has zero components when counting with sign or all components of $V_i$ have (co\nobreakdash-)orientation reversing automorphisms.
	And similarly for $V \sqcup W$.
	But then forming $V \sqcup W$ from $W$ contributes zero components when counting with sign or components with (co\nobreakdash-)orientation reversing automorphisms.
	So if $V \sqcup W$ is trivial, $W$ must have already been trivial.
\end{proof}

\begin{lemma}[Lipyanskiy Lemma 11]\label{L: bd defined}\index{$Q(M)$!preserves boundary}
	If $W$ is in $Q(M)$ then so is $\bd W$.
\end{lemma}

\begin{proof}
	We first check that if $W$ is trivial then so is its boundary $r_{\bd W} = r_Wi_{\bd W} \colon \bd W \to M$.
	If $\rho \colon W \to W$ is (co\nobreakdash-)orientation reversing, then we can consider $\rho_\bd$ as defined in \cref{R: bd diff}.
	As $r_W\rho = r_W$ we also have $r_Wi_{\bd W} = r_W\rho i_{\bd W} = r_W i_{\bd W} \rho_\bd$.
	Thus we only need see that $\rho_\bd$ is (co\nobreakdash-)orientation reversing.
	It is sufficient to consider what happens at points on the interior of $\bd W$ so working locally we may identify such points of $\bd W$ with points of $W$ itself and similarly identify $\rho_\bd$ with $\rho$.
	In the oriented case, the orientation of $W$ determines orientations of $T_xW$ and $T_{\rho(x)}W$, and by assumption $D\rho: T_xW \to T_{\rho(x)}W$ takes the orientation of $T_xW$ to the opposite of the orientation of $T_{\rho(x)}W$.
	But also $D\rho$ must preserve inward/outward pointing vectors.
	Thus $D\rho$ must restrict to a map $T_x (\bd W) \to T_{\rho(x)}(\bd W)$ that also reverses the orientation.
	The co-oriented situation is analogous using local orientation pairs $\left(\beta_{W,x}, \beta_{M,r_W(x)}\right)$ and $\left(\beta_{W,\rho(x)}, \beta_{M,r_W\rho(x)}\right)$ and noting $\beta_{M,r_W\rho(x)} = \beta_{M,r_W(x)}$ as $\rho$ is a diffeomorphism over $M$.

	Now suppose $W$ degenerate.
	By definition $\bd W = A \sqcup B$ with $A$ trivial and $B$ of small rank.
	Then $\bd^2 W = \bd A \sqcup \bd B$.
	As noted, $\bd A$ is trivial, and $\bd^2 W$ is trivial for all $W$ by \cref{L: boundary2}.
	It follows from \cref{L: Lip L10} that $\bd B$ is trivial, and so $B$ is degenerate.
	Thus $\bd W$ is degenerate.
\end{proof}

\begin{lemma}[Lipyanskiy Lemma 12]\label{L: Lipy12}\index{$Q(M)$!cancellation property}\index{cancellation property!in $Q(M)$}
	If both $V$ and $V \sqcup W$ are in $Q(M)$ then so is $W$.
\end{lemma}

\begin{proof}
	As in the proof of \cref{L: trivial structure}, decompose $W$ as $W = W_1 \sqcup W_2 \sqcup \cdots$ and $V$ as $V = V_1 \sqcup V_2 \sqcup \cdots$.
	As $V$ and $V \sqcup W$ are in $Q(M)$, each $V_i$ has small rank or, also as in the proof of \cref{L: trivial structure}, $V_i$ is trivial, and similarly for $V \sqcup W$, from which it follows again by counting with signs in the trivial components that each $W_i$ is either trivial or has small rank.
	By grouping terms of the decomposition we can write $W = A \sqcup B$ with $A$ trivial and $B$ of small rank.

	We then have $\bd V$ and $\bd (V \sqcup W) = \bd V \sqcup \bd A \sqcup \bd B$ in $Q(M)$ by \cref{L: bd defined}, and also $\bd A$ is trivial as the boundary of a trivial manifold over $M$ by the proof of \cref{L: bd defined}.
	Now by the same argument as in the previous paragraph, replacing $V$ with $\bd V \sqcup \bd A$ and $W$ with $\bd B$, we have that $\bd B$ can be decomposed into the disjoint union of a trivial manifold over $M$ and one with small rank.
	But this shows $B$ is degenerate, so $W \in Q(M)$.
\end{proof}

\begin{lemma}[Lipyanskiy Lemma 13]\label{L: cancel Q}\index{$Q(M)$!determines equivalence relation}
	The relation given by $V\sim W$ if $V \sqcup -W$ is in $Q_*(M)$ (respectively $Q^*(M)$) is an equivalence relation on $PC^\Gamma_*(M)$ (respectively $PC_\Gamma^*(M)$).
\end{lemma}

\begin{proof}
	Reflexivity: For any $W$, we have $W \sqcup -W$ trivial via the map that interchanges the two copies of $W$.

	Symmetry: If $V \sqcup -W$ is the union of trivial and degenerate manifolds over $M$ then certainly so is $-(V \sqcup -W) = W \sqcup -V$.

	Transitivity: If $V \sqcup -W$ and $W \sqcup -U$ are in $Q_*(M)$ (or $Q^*(M)$), then so is $V \sqcup -W \sqcup W \sqcup -U \cong V \sqcup -U \sqcup W \sqcup -W$.
	We know $W \sqcup -W$ is trivial, and so $V \sqcup -U$ is in $Q_*(M)$ (or $Q^*(M)$) by \cref{L: Lipy12}.
\end{proof}

These lemmas allow us to follow Lipyanskiy in defining geometric chains and cochains.
We will show that the claims of the following definition hold in \cref{L: co/chains well defined} just below.

\begin{definition}\label{D: chains and cochains}
	Let $M$ be a smooth manifold with corners.
	The \textbf{geometric chains}\index{geometric chain|textbf} of $M$, denoted $C^\Gamma_*(M)$\index{$C^\Gamma_*(M)$}, are the $\sim$ equivalence classes in $PC^\Gamma_*(M)$.
	The \textbf{geometric cochains}\index{geometric cochain|textbf} of $M$, denoted $C_\Gamma^*(M)$,\index{$C_\Gamma^*(M)$} are the $\sim$ equivalence classes in $PC_\Gamma^*(M)$.
	In either case, we denote the equivalence class of $W$ by $\uW$ and say that $W$ represents (or is a representative of) $\uW$.

	These are chain complexes with group operation $\uV + \uW = \underline{V \sqcup W}$ and boundary map $\bd \uW = \underline{\bd W}$.
	Furthermore, $\uW = 0$ in $C^\Gamma_*(M)$ (respectively $C_\Gamma^*(M)$) if and only if $W$ is in $Q_*(M)$ (respectively $Q^*(M)$), and in either case, the inverse of $\uW$ is $-\uW = \underline{-W}$.
	In particular, the empty manifold $\emptyset$ represents $0$ in $C^\Gamma_*(M)$ and the empty map $\emptyset \to M$ represents $0$ in $C^*_\Gamma(M)$.

	We define the \textbf{geometric homology}\index{geometric homology|textbf} of $M$, written $H_*^\Gamma(M)$, to be $H_*(C^\Gamma_*(M))$,\index{$H_*^\Gamma(M)$} and
	we define the \textbf{geometric cohomology}\index{geometric cohomology|textbf} of $M$, written $H^*_\Gamma(M)$ to be $H^*(C_\Gamma^*(M))$.\index{$H^*_\Gamma(M)$}
\end{definition}

We note that the formula $\uV + \uW = \underline{V \sqcup W}$ implies that every geometric chain or geometric cochain can be represented by a single map $r_Z \colon Z \to M$ for appropriate (possibly disconnected) $Z$.
In particular, $k \uV$ for $k>0$ can be represented by a single map consisting of the disjoint union of $k$ (co-)oriented isomorphic pre(co)chains, while $k \uV$ for $k<0$ is represented by the same map with the opposite (co-)orientations; see \cref{E: copies,D: prechain sum}.
We will show below in \cref{E: uncountable} that the $C^i_\Gamma(M)$ are not necessarily free abelian groups.
We do not know if the $C^\Gamma_i(M)$ are free abelian groups; if so this is not evident from their construction.
However, we will show that these groups are all torsion free below in \cref{L: flat}.

\begin{lemma}\label{L: co/chains well defined}\index{geometric chain!negative of}\index{geometric chain!representing $0$}\index{geometric cochain!negative of}\index{geometric cochain!representing $0$}
	The chain complexes $C_\Gamma^*(M)$ and $C_\Gamma^*(M)$ are well defined with the properties claimed in the definition.
	In particular, $\uW = 0$ in $C^\Gamma_*(M)$ (respectively $C_\Gamma^*(M)$) if and only if $W$ is in $Q_*(M)$ (respectively $Q^*(M)$), and the additive inverse of $\uW$ is $\underline{-W}$.
\end{lemma}

\begin{proof}
	We will apply the preceding lemmas.
	For simplicity we work with $Q_*(M)$ and $C^\Gamma_*(M)$, but the identical arguments hold with $Q^*(M)$ and $C_\Gamma^*(M)$.

	Suppose $\uV,\uW \in C^\Gamma_*(M)$ with $V$ and $V'$ in the class $\uV$ and $W$ and $W'$ in the class $\uW$.
	Then the sum $\uV+\uW$ is well defined as the class $\underline{V \sqcup W}$ because $(V \sqcup W) \sqcup -(V' \sqcup W') = (V \sqcup -V') \sqcup (W \sqcup -W')$, and $V \sqcup -V'$ and $W \sqcup -W'$ are both in $Q_*(M)$ by assumption.

	The identity in each degree is represented by $\emptyset$ with the unique empty map to $M$ (with either orientation or co-orientation).
	In fact, every element of $Q_*(M)$ represents $0$ in $C^\Gamma_*(M)$ as elements of $Q_*(M)$ are all equivalent to $\emptyset$.
	Conversely, if $\uW = 0$ then $W \in Q_*(M)$, as if $\uV+\uW = \uV$ then $V \sqcup W \sqcup -V \in Q_*(M)$.
	But $V \sqcup -V \in Q_*(M)$, so by \cref{L: Lipy12}, $W \in Q_*(M)$.
	We also see that the additive inverse of $\uW$ is $\underline{-W}$, as $W \sqcup -W$ is trivial.

	That the boundary map is well defined with $\bd \uW = \udW$ is due to \cref{L: bd defined} and \cite[Lemma 2.8]{Joy12}, which implies that $\bd W$ is proper.
	That $\bd^2 = 0$ in the case of cochains follows from \cref{L: boundary2}, which shows that $\bd^2 W$ is always trivial.
	Similarly, to obtain $\bd^2 = 0$ for chains see \cref{R: bd2 oriented}.
\end{proof}

\begin{remark}\label{R: cycles and boundaries}\index{geometric chain!cycle representative}\index{geometric chain!boundary representative}\index{geometric cochain!cocycle representative}\index{geometric cochain!coboundary representative}
	Suppose $W \in PC_*^{\Gamma}(M)$.
	It follows directly from \cref{D: chains and cochains,L: co/chains well defined} that $W$ represents a cycle in $C_*^{\Gamma}(M)$, i.e.\ $\bd \uW=0$, if and only if $\bd W \in Q_*(M)$.
	We also observe that $W$ represents a boundary precisely when there is some $Z \in PC_*^\Gamma(M)$ such that $\bd \underline{Z} = \underline{\bd Z} = \uW$, which translates to $(\bd Z) \sqcup -W \in Q_*(M)$.
	The corresponding statements for cohomology are analogous.
\end{remark}

In Section 6 of \cite{Lipy14}, Lipyanskiy shows that the homology theory based on geometric chains satisfies some of the Eilenberg-Steenrod axioms.
This is enough to state in Section 10 of \cite{Lipy14} that geometric homology is isomorphic to singular homology on the fixed manifold $M$, though we provide our own proofs below in \cref{T: geometric is singular,T: hom iso map}.
Lipyanskiy does not provided a detailed treatment of geometric cohomology, which is different from geometric homology in several respects, though we will also show that it is isomorphic to singular cohomology in \cref{T: geometric is singular,T: intersection qi}.

We next present some immediate examples and observations.
The first is that geometric homology and cohomology satisfy a very strong form of Poincar\'e duality.

\begin{theorem}[Poincar\'e Duality]\label{T: PD}\index{Poincar\'e Duality}
	If $M$ is closed and oriented then tautologically $C_*^\Gamma(M) = C_\Gamma^{m-*}(M)$.
	Consequently, also $H_*^\Gamma(M) = H_\Gamma^{m-*}(M)$.
\end{theorem}
\begin{proof}
	This following directly from the definitions using that the domain of a proper map to a compact space must be compact and that there is an induced orientation on the domain of a co-oriented map with oriented codomain; see the discussion following \cref{D: tautological co-orientation}.
\end{proof}

\begin{example}\label{E: first examples}\index{geometric homology!in top and bottom degrees|(}\index{geometric cohomology!in top and bottom degrees|(}
	We will see in the next chapter that geometric homology and cohomology are isomorphic to ordinary singular homology and cohomology with $\Z$-coefficients.
	So the following calculations are not surprising, though we hope they may be illuminating and help build some intuition.
	We begin by considering homology and cohomology classes represented by maps from $0$-manifolds.

	For any manifold without boundary $M$, an element of $PC_0^\Gamma(M)$ is a map to $M$ from a finite disjoint union of signed points.
	As disjoint union in $PC_0^\Gamma(M)$ translates to addition in $C_0^\Gamma(M)$ the map of any oriented point to $M$ represents a generator of $C_0^\Gamma(M)$.
	Analogously to the computation for singular homology, any two points with the same orientation mapping to the same component of $M$ represent the same element of $H_0^\Gamma(M)$ as can be seen by joining them with a smooth path and applying our boundary conventions.
	On the other hand, the boundary of such an interval mapping to $M$ consists of two points with opposite orientations, and the only compact $1$-dimensional manifolds with corners consist of closed intervals and circles, which have no boundaries.
	So it follows that a positively-oriented point mapping to $M$ generates an infinite cyclic subgroup of $H_0^\Gamma(M)$.
	Altogether, $H_0^\Gamma(M) \cong \oplus \Z$, where the sum is taken over the connected components of $M$.

	Next, suppose $M$ is a closed, connected manifold.
	The considerations for $H^m_\Gamma(M)$ are similar to those for $H_0^\Gamma(M)$, noting that elements of $PC^m_\Gamma(M)$ are represented by proper co-oriented maps of $0$-manifolds to $M$.
	However, in this case if we are given a map $r \colon [0,1] \to M$ with $r(0) = r(1) = z \in M$, the signs of the boundary components will depend on whether or not the loop determined by $r$ preserves or reverse the orientation of $M$.
	More specifically, let $r_0$ denote the map that takes the point to $z \in M$ with the co-orientation $(1, \beta_M)$.
	Treating $r$ as a homotopy, by \cref{R: stationary homotopy} the boundary of $r$ will be the disjoint union of $-r_0$ and $r_1$, where $r_1$ takes the point to $z$ with co-orientation $(1, r_*\beta_M)$.
	Recall $r_*\beta_M$ is the orientation of $M$ at $z$ obtained by starting with the orientation $\beta_M$ at $z$ and traveling around the loop $r$; see \cref{S: co-orientations}.
	If $r_*\beta_M = \beta_M$, then $r_1 = r_0$ and we have $\bd r = r_1 - r_0 = 0$, which is analogous to the homology computation.
	But if $r_*\beta_M = -\beta_M$, then $r_1 = -r_0$ and we have $\bd r = r_1 - r_0 = -2 r_0$, so twice $r_0$ is a boundary.
	Furthermore, for any two maps $r_0$ and $s_0$ taking the point to the same component of $M$ (with either co-orientation), by joining them with a smooth path, we see that the maps will be cohomologous up to sign.
	So we conclude that if $M$ is connected then
	\[H^m_\Gamma(M) =
	\begin{cases}
		\Z, & M \text{ is orientable},\\
		\Z_2, & M \text{ is not orientable}.
	\end{cases}
	\]

	On the other hand, if $M$ is connected, without boundary, and not compact, then any co-oriented map $r_0$ from the point to $M$ represents $0$ in $H_\Gamma^m(M)$ because any proper smooth path $(-\infty,0] \to M$ with the restriction to $0$ being $r_0$ gives a null-cohomology of $r_0$ with an appropriate choice of co-orientation.

	Altogether, if $M$ is a manifold without boundary, not necessarily connected, then $H^m_\Gamma(M)$ is the direct product of $\Z$ and $\Z_2$ factors contributed respectively by the closed orientable and closed non-orientable components of $M$.

	At the opposite end of the spectrum, the groups $H_m^\Gamma(M)$ and $H^0_\Gamma(M)$ are represented by maps of $m$-dimensional manifolds with corners.
	If $M$ is connected, then $H^0_\Gamma(M) \cong \Z$ with a generator represented by the identity map $M \to M$ with its tautological co-orientation.
	This will follow most readily from \cref{E: coho 0 generator,T: intersection qi}, below.
	Consequently, if $M$ is closed and oriented, then \cref{T: PD} implies $H_m^\Gamma(M) \cong \Z$ generated by the identity map.\index{geometric homology!in top and bottom degrees|)}\index{geometric cohomology!in top and bottom degrees|)}
\end{example}

\begin{example}\label{E: dimension range}\index{geometric homology!vanishing theorem}\index{geometric cohomology!vanishing theorem}
	If $\dim(M) = m$, then $H_i^\Gamma(M) = H^i_\Gamma(M) = 0$ if $i < 0$ or $i > m$.
	In fact, $C_i^\Gamma(M) = 0$ when $i<0$ and $C^i_\Gamma(M) = 0$ when $i > m$, as in these cases the only elements of $PC_i^\Gamma(M)$ or $PC^i_\Gamma(M)$ are respectively the empty manifold $\emptyset$ and the empty map $\emptyset \to M$.
	For homology when $i>m+1$, every element of $PC_i^\Gamma(M)$ must have small rank and its boundary must also have small rank, so it is in $Q_i(M)$ and hence $C_i^\Gamma(M) = 0$.
	When $i = m+1$, $C_{m+1}^\Gamma(M)$ may be non-zero, but every element will have small rank.
	If $r_W \colon W \to M$ represents a degree $m+1$ cycle, then its boundary must be in $Q_m(M)$ by \cref{R: cycles and boundaries}.
	So $W$ is degenerate and represents $0$ in $C_{m+1}^\Gamma(M)$ and hence in $H_{m+1}^\Gamma(M)$.
	The argument for $H^i_\Gamma(M)$ for $i<0$ is the same.
\end{example}

\begin{example}[Dimension axiom]\label{E: dimension}\index{geometric homology!dimension axiom}\index{geometric cohomology!dimension axiom}
	By \cref{E: dimension range}, $H_i^\Gamma(pt) = H^i_\Gamma(pt) = 0$ unless $i = 0$.
	When $i = 0$, we have $H^0_\Gamma(M) \cong H_0^\Gamma(M) \cong \Z$ by the computations of \cref{E: first examples}.
	In fact, this is an equality by \cref{T: PD}.
\end{example}

\begin{remark}\label{R: degen1}
	It is in these examples that we most obviously see the need to include degenerate chains and cochains in $Q(M)$.
	On the other hand, the formulation of degeneracy as given will create some difficulty for us in \cref{S: products} when it comes to consider notions of transversality for geometric chains and cochains.
	The reason is that degeneracy causes much of the problematic ambiguity in choosing representatives for chains and cochains.
	For example, consider a connected prechain $V \in PC_*(M)$ with small rank but a boundary that is not in $Q_*(M)$.
	If $V'$ is any other such prechain with small rank and $\bd V = \bd V'$, then $V \sqcup -V' \in Q_*(M)$ since it will have small rank and trivial boundary.
	So $V$ and $V'$ represent the same chain but could behave wildly differently aside from their boundaries.
	By contrast, we will see in \cref{S: products} that trivial chains are less of an issue (they can generally be ignored) and that non-trivial components of prechains without small rank are ``essential'' in a sense we will make precise in \cref{D: essential}.
	In particular, essential components appear in any representative of the same geometric chain or cochain.

	Given the headaches thus caused by the degenerate chains and cochains, it is tempting to ask for a simpler definition of degeneracy.
	One variant that comes to mind would be defining degeneracy so that each individual component must have a boundary consisting of trivial and small rank (co)chains.
	This would seem to be sufficient for the dimension axiom and would eliminate the difficulty described above.
	Unfortunately, with such an alternative definition of degeneracy, it will not generally be true that if $V \in Q(M)$ then fiber products $V \times_M W$ are also in $Q(M)$.
	This is an important property that will arise in the next section and then be needed both to construct cup and cap products and to prove the existence of Mayer--Vietoris sequences.
	See \cref{R: degen2} for further discussion of this point.
\end{remark}

The following algebraic property will be useful below as we consider homological algebra with geometric cochains:

\begin{lemma}\label{L: flat}\index{geometric chain complex is torsion-free}\index{geometric cochain complex is torsion-free}
	Each $C_i^\Gamma(M)$ or $C_\Gamma^i(M)$ is torsion-free and hence flat as a $\Z$-module.
\end{lemma}

\begin{proof}
	The second statement follows from the first as $\Z$ is a Dedekind domain.
	The first statement is proven for geometric chains in \cite[Lemma 34]{Lipy14}.
	The proof for geometric cochains, even accounting for our different definition of degeneracy, is the same.
	The argument is as follows.

	Suppose $\uW$ satisfies $n \uW = 0$ for some $n \in \Z$ with $n \neq 0$.
	Representing $\uW$ by $W$, we can decompose $W$ as in the proofs of \cref{L: trivial structure,L: Lipy12} so that $W = W_1 \sqcup W_2 \sqcup \cdots$ with each $W_i$ being the union of isomorphic-up-to-sign connected manifolds over $M$ with the isomorphism classes of components in $W_i$ and $W_j$ distinct if $i\neq j$.
	Even in the cochain case, each $W_i$ contains only finitely many connected components due to the properness condition.
	Then $n \uW$ is represented by taking $n$ copies of each $W_i$ if $n > 0$ or, if $n < 0$, taking $|n|$ copies of each $-W_i$.
	By \cref{L: co/chains well defined}, $n W \in Q(M)$, so $nW$ is the union of a trivial pre(co)chain and a degenerate pre(co)chain, which we can write as $nW = T \sqcup D$.
	A connected component of $W_i$ has small rank if and only if they all do, so if $W_i$ consists of components that are not of small rank, then all of $nW_i$ must be part of $T$.
	Furthermore, as the components of $nW_i$ and $nW_j$ are distinct for $i \neq j$, each such $nW_i$ must itself be trivial.
	It follows that each $n W_i$ is trivial or has small rank.

	Suppose $n W_i$ is trivial for some $i$.
	Then each connected component of $n W_i$ is trivial or the total number of components of $n W_i$ is zero when counted with sign.
	But this implies that each component of $W_i$ is trivial or that $W_i$ must also have the total number of components be zero when counted with sign.
	So $W_i$ is trivial.
	It is also clear that if $W_i$ is trivial, then so is $n W_i$, so $W_i$ is trivial if and only if $n W_i$ is trivial.

	Now let $W'$ be the union of those $W_i$ that are not trivial.
	Then $nW$ is the union of a trivial pre(co)chain and $n W'$, so by \cref{L: Lipy12}, we have $n W' \in Q(M)$.
	Furthermore, $n W'$ must have small rank, which implies that $W'$ must have small rank.
	Further, since $n W' \in Q(M)$, we have $\bd (n W') = n (\bd W') \in Q(M)$ by \cref{L: bd defined}.
	So by the same argument as above, $\bd W'$ must be a union of pre(co)chains that are trivial or of small rank.
	So $W'$ is degenerate.

	Altogether, we have now shown that $W$ is the union of a trivial pre(co)chain and a degenerate pre(co)chain, so $W \in Q(M)$.
\end{proof}

\subsection{Products of manifolds over \texorpdfstring{$M$}{M}}

In this section we define various products of elements of $PC_*^\Gamma(M)$ and $PC^*_\Gamma(M)$, all coming from the external products or fiber products defined above.
When we need transversality, $M$ will be without boundary.
These products will ultimately become our cup, cap, intersection, and exterior products, but we introduce them here as products on $PC_*^\Gamma(M)$ and $PC^*_\Gamma(M)$ and derive some further properties that we will need in the following sections.

We begin with the products coming from fiber products.
These are only partially defined, as fiber products require transversality.

\begin{definition}\label{D: PC products}\index{manifold over a manifold!fiber product|(}\index{fiber product|(}\index{geometric precochain!cup product|(}\index{geometric prechain!cap product|(}\index{geometric precochain!cap product|(}\index{geometric prechain!intersection product|(}\index{cup product|(}\index{cap product|(}\index{intersection product|(}
	Given a manifold without boundary $M$, the fiber products of transverse manifolds over $M$ determine partially-defined products of the following forms:
	\begin{align*}
		PC^*_\Gamma(M) \times PC^*_\Gamma(M)& \to PC^*_\Gamma(M)\\
		PC^*_\Gamma(M) \times PC_*^\Gamma(M)& \to PC_*^\Gamma(M).
	\end{align*}
	If, furthermore, $M$ is oriented, then there is also a partially-defined product
	$$PC_*^\Gamma(M) \times PC_*^\Gamma(M) \to PC_*^\Gamma(M).$$
	In each case, the product is defined when the reference maps $r_V \colon V \to M$ and $r_W \colon W \to M$ are transverse.

	\begin{itemize}
		\item If $V,W\in PC^*_\Gamma(M)$ are transverse, we define their product to be the fiber product $V \times_M W$ with its fiber product co-orientation of \cref{D: pullback coorient}.
		Note that the fiber product is proper by \cref{L: co-orientable pullback} as the composition of proper maps is proper.

		\item If $V \in PC^*_\Gamma(M)$ and $W \in PC_*^\Gamma(M)$ are transverse, we define their product to be $V \times_M W$ with the orientation induced from the orientation of $W$ and the co-orientation of the pullback $V \times_M W \to W$; see \cref{D: pullback coorient} and the discussion following \cref{D: co-orientations}.
		This is precisely the cap orientation of \cref{S: mixing}.
		Note that $V \times_M W$ is compact, as $W$ is compact and the pullback map is proper by \cref{L: co-orientable pullback} and the properness of $r_V$.

		\item If $M$ is oriented and $V,W \in PC_*^\Gamma(M)$, we define their product to be $V \times_M W$ with the fiber product orientation of \cref{S: orientation of fiber products}.
		Note that $V \times_M W$ is compact as a closed subset of the compact set $V \times W$.
	\end{itemize}

	We show below that these products are all well defined when the transversality conditions are met, as transversality is preserved under isomorphisms of manifolds over $M$, and isomorphic manifolds over $M$ have isomorphic fiber products over $M$.

	In all cases we continue to denote the product by $V \times_M W$, allowing context to determine which (co\nobreakdash-)orientations apply.
\end{definition}

In the last case, we really need $M$ to be oriented in general, as we have observed in \cref{R: what products exist} that if $M$ is not orientable the fiber product of orientable manifolds over $M$ may not be orientable.

\begin{lemma}\label{L: product preserves iso}
	When the transversality conditions are met, the products of \cref{D: PC products} are well defined.
	In particular, if $\phi_V \colon V \to V'$ and $\phi_W \colon W \to W'$ are orientation- or co-orientation-preserving diffeomorphisms of transverse manifolds over $M$, then $\phi_V \times \phi_W$ restricts, as appropriate to the products of \cref{D: PC products}, to a corresponding orientation- or co-orientation-preserving diffeomorphism $V \times_M W \to V' \times_M W'$.

	Similarly, if exactly one of $\phi_V \colon V \to V'$ or $\phi_W \colon W \to W'$ is an orientation- or co-orientation-\textit{reversing} diffeomorphism, then $\phi_V \times \phi_W$ restricts, as appropriate to the products of \cref{D: PC products}, to a corresponding orientation- or co-orientation-\textit{reversing} diffeomorphism $V \times_M W \to V' \times_M W'$.
\end{lemma}
\begin{proof}
	It is clear that the transversality conditions are preserved if representatives $r_V \colon V \to M$ or $r_W \colon W \to M$ are replaced with (oriented or co-oriented) isomorphic manifolds over $M$.
	It is also clear from the universal property that isomorphic manifolds over $M$ lead to diffeomorphic fiber products.
	So our main challenge is to verify the correct behavior of orientations and co-orientations.

	We first consider $V$, $W$, and $M$ all oriented.
	Suppose we have $f \colon V \to M$ isomorphic to $f' \colon V' \to M$ via $\phi_V \colon V \to V'$ and similarly for $W$; see \cref{D: equiv triv and small}.
	The map $\phi_V \times \phi_W$ restricts to give our diffeomorphism from $P = V \times_M W$ to $P' = V' \times_M W'$, so, considering the definition of the fiber product orientation in \cref{S: orientation of fiber products}, we obtain bundle map isomorphisms covering $\phi_V \times \phi_W \colon P \to P'$ of the following exact sequences (or we can think of the bottom sequence as pulled back to $P$ by $\phi_V \times \phi_W$).
	\begin{equation}
		\begin{tikzcd}
			0 \arrow[r] & TP \arrow[r,"D\pi_V \oplus D\pi_W"] \arrow[d] &[1cm] \pi_V^*(TV) \oplus \pi_W^*(TW) \arrow[r,"\pi_V^*(Df)-\pi_W^*(Dg)"] \arrow[d] &[1.7cm] (f\pi_V)^*TM \arrow[r] \arrow[d] & 0 \\
			0 \arrow[r] & TP' \arrow[r,"D\pi_{V'} \oplus D\pi_{W'}"] & \pi_{V'}^*(TV') \oplus \pi_{W'}^*(TW') \arrow[r,"\pi_{V'}^*(Df')-\pi_{W'}^*(Dg')"] & (f'\pi_{V'})^*TM \arrow[r] & 0
		\end{tikzcd}
	\end{equation}
	As the maps $\phi_V$ and $\phi_W$ are orientation preserving, it follows that the two righthand vertical maps are oriented bundle isomorphisms.
	Consequently, if we orient $TP$ and $TP'$ as in \cref{S: orientation of fiber products} and use the isomorphisms to make the splittings of the sequences compatible, then the lefthand map must also be an isomorphism of oriented bundles.
	(Looking at the level of individual tangent space fibers, the reader can also compare with \cite[Sections 9.1.1 and 9.3.1]{RamBas09}.)

	For the other cases, we first observe that being an orientation or co-orientation preserving isomorphism is a local property.
	In particular, a diffeomorphism $V \to V'$ is orientation preserving if and only its derivative is orientation preserving at each $x \in V$, and by \cref{L: co-or preserving/reversing} a similar statement holds for co-orientations.
	Consequently, to show that a diffeomorphism over $M$ is (co-)orientation preserving, it suffices to show that this property holds over each set of an open cover of $M$.

	So let $U$ be any Euclidean subset of $M$, and suppose we give $U$ an arbitrary orientation.
	Considering the case of co-oriented precochains, by \cref{C: co-or preserving is or preserving} the restrictions of $\phi_V$ and $\phi_W$ to $f^{-1}(U)$ and $g^{-1}(U)$ are co-orientation preserving if and only if they are orientation preserving with respect to the induced orientations.
	Note that this statement is independent of the choice of fixed orientation for $U$.
	Since we assume that $\phi_V$ and $\phi_W$ are co-orientation preserving, they are thus orientation preserving over $U$ with respect to the induced orientations.
	It follows from the preceding argument that $\phi_V \times \phi_W$ restricts to an orientation-preserving diffeomorphism of the fiber products over $U$.
	By \cref{P: compare cup and intersection orientations}, there is a fixed sign relation depending only on the dimensions between the fiber product orientations obtained by the above construction and the orientations of the fiber products induced from the orientation of $U$ and the co-oriented fiber product (this is the distinction considered there between $V \times^o_M W$ and $V \times^c_M W$).
	Since $\phi_V \times \phi_W$ provides an orientation-preserving diffeomorphism with respect to the fiber product oriented the first way, it also provides an orientation-preserving diffeomorphism with respect to the fiber product oriented the other way.
	Hence by \cref{C: co-or preserving is or preserving} again, $\phi_V \times \phi_W$ restricts to a co-orientation-preserving diffeomorphism $f^{-1}(U) \times_U g^{-1}(U) \to (f')^{-1}(U) \times_U (g')^{-1}(U)$.
	Since $U$ was arbitrary, $\phi_V \times \phi_W$ gives a co-orientation-preserving diffeomorphism on all of $V \times_M W$.

	Lastly, for the fiber product $PC^*_\Gamma(M) \times PC_*^\Gamma(M) \to PC_*^\Gamma(M)$, we again begin by considering what happens for the preimage of a Euclidean set $U$ in $M$ with an arbitrary orientation $\beta_U$.
	In this case the restriction $g \colon g^{-1}(U) \to U$ is a map between oriented manifolds and so if $W$ has orientation $\beta_W$, this map has an induced co-orientation $(\beta_W, \beta_U)$.
	Now, by definition, $\beta_P$ is the desired orientation of $f^{-1}(U) \times_U g^{-1}(U)$ for this kind of product if and only $(\beta_P,\beta_W)$ is the pullback co-orientation of $f^{-1}(U) \times_U g^{-1}(U) \to g^{-1}(U)$.
	And this is the case if and only if $(\beta_P,\beta_W)*(\beta_W,\beta_U) = (\beta_P, \beta_U)$ is the co-orientation of the co-oriented fiber product of $f$ and $g$ restricted to $f^{-1}(U)$ and $g^{-1}(U)$, obtained from $g$ being given the induced co-orientation as above.
	Furthermore, as $\phi_W$ is orientation-preserving, its restriction to a diffeomorphism $g^{-1}(U)$ to $(g')^{-1}(U)$ is also co-orientation preserving with respect to the induced co-orientations by \cref{C: co-or preserving is or preserving}.
	But now we have seen that the product $\phi_V \times \phi_W$ of co-orientation-preserving diffeomorphisms gives a co-orientation preserving diffeomorphism on $f^{-1}(U) \times_U g^{-1}(U)$.
	So now again by \cref{C: co-or preserving is or preserving}, $\phi_V \times \phi_W$ is an orientation-preserving diffeomorphism with respect to the orientations induced by the orientation of $U$ and the fiber product co-orientations.
	But we have just recalled that these are exactly the orientations that we want for these products.
\end{proof}

Given the preceding lemma, from here on we will usually not explicitly distinguish between (co\nobreakdash-)oriented maps $V \to M$ and their isomorphism classes when discussing products.

The next lemma is critical and will be used repeatedly.
Among its applications, it will be used in \cref{S: products} toward showing that our fiber products products extend to well-defined, though only partially-defined, products of geometric chains and cochains.
It will also be needed much sooner to show that the creasing construction is well defined; this construction is used, in turn, to demonstrate the existence of Mayer--Vietoris sequences.

For the statement, recall that we use $Q(M)$ to stand for $Q_*(M)$ or $Q^*(M)$ as appropriate.

\begin{lemma}\label{L: pullback with Q}\index{fiber product!trivial}\index{fiber product!degenerate}\index{fiber product!small rank}\index{trivial!fiber product}\index{degenerate!fiber product}\index{$Q(M)$!in fiber product}
	For any of the products of \cref{D: PC products}, if either $V$ or $W$ is in $Q(M)$ then so is $V \times_M W$.
	In fact, if $V$ or $W$ is trivial then $V \times_M W$ is trivial, if $V$ or $W$ has small rank then $V \times_M W$ has small rank, and if $V$ or $W$ is degenerate then $V \times_M W$ is degenerate.
\end{lemma}

\begin{proof}
	We provide the proof if $V \in Q(M)$; the other case is similar.
	By assumption $V$ is the disjoint union of trivial and degenerate chains or cochains, so it suffices to consider independently the possibilities that $V$ is trivial or degenerate.

	If $\rho$ is a (co\nobreakdash-)orientation reversing diffeomorphism of $V$ over $M$, then $\rho \times_M \id_W$ is a (co\nobreakdash-)orientation reversing diffeomorphism of $V \times_M W$ by \cref{L: product preserves iso}.

	Next assume that $V$ is degenerate, so in particular it has small rank.
	Recall that the tangent bundle of a fiber product is the fiber product of the tangent bundles by \cref{L: tangent of pullbacks}, and so the derivative is the fiber product of derivatives.
	Note that the fiber product of two linear maps, one with a non-trivial kernel, must also have a non-trivial kernel: If $A,B$ are linear maps with a common codomain and $v \in \ker(A)$, then $(v,0)$ is in the kernel of the fiber product of $A$ and $B$.
	So if the differential of $r_V$ has non-trivial kernel everywhere so will the derivative of any fiber product with $r_V$.
	Thus $V \times_M W$ has small rank.

	Now we recall that $\bd(V \times_M W)$ is, up to (co\nobreakdash-)orientations, the disjoint union of $(\bd V) \times_M W$ and $V \times_M (\bd W)$.
	We have just shown that $V \times_M (\bd W)$ must have small rank.
	As $V$ is degenerate, $\bd V$ is a disjoint of trivial and small rank manifolds over $M$, and so by the preceding arguments $(\bd V) \times_M W$ will be a union of trivial and small rank manifolds over $M$.
	Altogether, $V \times_M W$ is degenerate.
\end{proof}

\begin{remark}\label{R: degen2}
	As noted in \cref{R: degen1}, it is this lemma that fails if we attempt to simplify the definition of degeneracy by requiring each connected component of a degenerate prechain or precochain to have small rank and a boundary that is a union of trivial and small rank pre(co)chains.
	In fact, it is possible to construct a $V$ and $W$ such that $V$ is a non-trivial prechain that is degenerate in this stronger sense but such that $V \times_M W$ has multiple components that are each non-trivial and of small rank but such that the boundary of each component is non-trivial and not of small rank.
	So $V \times_M W$ would not be in a version of $Q(M)$ defined using this stronger, but simpler, notion of degeneracy.
	Of course it is in $Q(M)$ with our actual definitions by the preceding lemma.
\end{remark}

\index{manifold over a manifold!fiber product|)}\index{fiber product|)}\index{geometric precochain!cup product|)}\index{geometric prechain!cap product|)}\index{geometric precochain!cap product|)}\index{geometric prechain!intersection product|)}\index{cup product|)}\index{cap product|)}\index{intersection product|)}

\index{manifold over a manifold!exterior product|(}\index{exterior product|(}\index{geometric precochain!cross product|(}\index{geometric prechain!cross product|(}
We have similar results for the exterior products studied in \cref{S: exterior products}, which are always fully defined:

\begin{lemma}\label{L: ext product preserves iso}
	Suppose $f \colon V \to M$ and $g \colon W \to N$ are maps of manifolds with corners with either $V$ and $W$ both oriented or $f$ and $g$ both co-oriented.
	Then the corresponding exterior product $PC_*^\Gamma(M) \times PC_*^\Gamma(N) \to PC_*^\Gamma(M \times N)$ or $PC^*_\Gamma(M) \times PC^*_\Gamma(N) \to PC^*_\Gamma(M \times N)$ is well defined.

	In particular, if $\phi_V \colon V \to V'$ and $\phi_W \colon W \to W'$ are orientation- or co-orientation-preserving diffeomorphisms over $M$ and $N$, respectively, then $\phi_V \times \phi_W$ is a corresponding orientation- or co-orientation-preserving diffeomorphism $V \times W \to V' \times W'$.

	Similarly, if exactly one of $\phi_V \colon V \to V'$ or $\phi_W \colon W \to W'$ is an orientation- or co-orientation-\textit{reversing} diffeomorphism, then $\phi_V \times \phi_W$ is a corresponding orientation- or co-orientation-reversing diffeomorphism $V \times W \to V' \times W'$.
\end{lemma}
\begin{proof}
	The oriented case is standard.
	For the co-oriented case, we can proceed analogously to the proof of \cref{L: product preserves iso}.
	In particular, we can choose Euclidean neighborhoods $A \subset M$ and $B \subset N$ given arbitrary orientations.
	Then, as in that preceding proof, isomorphisms of product co-orientations over $A \times B$ follow from isomorphisms of product orientations over $A \times B$, now using \cref{P: compare exterior orientations} rather than \cref{P: compare cup and intersection orientations}.
\end{proof}

\begin{lemma}\label{L: exterior Q}\index{exterior product!trivial}\index{exterior product!degenerate}\index{exterior product!small rank}\index{trivial!external product}\index{degenerate!external product}\index{$Q(M)$!in external product}
	Consider two pre(co)chains $f \colon V \to M$ and $g \colon W \to N$.
	If $V \in Q(M)$ or $W \in Q(N)$ then $f \times g \colon V \times W \to M \times N$ is in $Q(M \times N)$.
\end{lemma}

\begin{proof}
	We provide the proof if both maps are in $PC^*$ and $V \in Q^*(M)$; the other cases are similar.
	By assumption $V$ is the disjoint union of trivial and degenerate chains or cochains, so it suffices to consider independently the possibilities that $V$ is trivial or degenerate.

	If $\rho$ is a co-orientation-reversing diffeomorphism of $V$ over $M$, then $\rho \times \id_W$ is a co-orientation-reversing diffeomorphism of $V \times_M W$ by \cref{L: ext product preserves iso}.
	So if $V$ is trivial so is $V \times W$.

	Next assume that $V$ is degenerate, so in particular it has small rank.
	At any point, the derivative of $f \times g$ is a matrix with $Df$ and $Dg$ on the block diagonals, so $f \times g$ has small rank.
	Now we recall that $\bd(V \times W)$ is, up to co-orientations, the disjoint union of $(\bd V) \times_M W$ and $V \times_M (\bd W)$.
	We have just shown that $V \times_M (\bd W)$ must have small rank.
	As $V$ is degenerate, $\bd V$ is a disjoint union of trivial and small rank manifolds over $M$, and so by the preceding arguments $(\bd V) \times_M W$ will be a union of trivial and small rank manifolds over $M$.
	Altogether, $V \times_M W$ is degenerate.
\end{proof}

\index{manifold over a manifold!exterior product|)}\index{exterior product|)}\index{geometric precochain!cross product|)}\index{geometric prechain!cross product|)}

\subsection{Essential decompositions}\label{S: essential decomp}\index{essential decomposition|(}

The material in this subsection will not be needed until \cref{S: products}, but we present it here both because it involves fundamental properties of chains and cochains and because it may be illuminating for the reader concerning the ambiguities involved when considering prechain and precochain representatives of chains and cochains.

It will be useful to adopt the following conventions.
Let $W$ be a connected manifold with corners over $M$.
We say that ``$V$ is isomorphic to $\pm W$'' if $V$ and $W$ are either both elements of $PC_*^\Gamma(M)$ or both elements of $PC^*_\Gamma(M)$ and $V$ is oriented or co-oriented isomorphic\footnote{Recall \cref{D: equiv triv and small}.} either to $W$ or to $W$ with its opposite orientation or co-orientation, as appropriate.
We will say that ``$V$ has $n$ components isomorphic to $W$, counting without sign'' if $V$ has exactly $n$ components, each of which is isomorphic to $\pm W$.
We will say that ``$V$ has $n$ components isomorphic to $W$, counting with sign'' if $n$ is the number of components of $V$ isomorphic to $W$ minus the number of components of $V$ isomorphic to $-W$.
We will only count with sign in contexts in which $W$ is not isomorphic to $-W$, so there should be no ambiguity.
We also note that $W$ being isomorphic to $-W$ is the same as $W$ being trivial.

Also, analogously to our ambiguous notation $Q(M)$ established in \cref{D: Q}, we will sometimes write $W \in PC(M)$ to mean $W \in PC^\Gamma_*(M)$ or $W \in PM_\Gamma^*(M)$ for definitions or arguments that are analogous in the two cases.

\begin{definition}\label{D: essential}
	Suppose $V \in PC(M)$ and $V_1$ is a connected component of $V$.
	\begin{itemize}
		\item We call $V_1$ \textbf{trivially inessential}\index{trivially inessential|textbf}\index{geometric prechain!trivially inessential}\index{geometric precochain!trivially inessential} if either
		$V_1$ is trivial or if the number of components of
		$V$ isomorphic to $\pm V_1$ is $0$, counting with signs.

		\item We call ${V_1}$ \textbf{non-trivially inessential}\index{non-trivially inessential|textbf}\index{geometric prechain!non-trivially inessential}\index{geometric precochain!non-trivially inessential} if it is not trivially inessential but it is of small rank.

		\item We call ${V_1}$ \textbf{essential}\index{essential|textbf}\index{geometric prechain!essential}\index{geometric precochain!essential} if it is not (trivially or non-trivially) inessential.
	\end{itemize}

	Each connected component of $V$ falls into exactly one of these categories by definition.

	The \textbf{essential decomposition} of $V$ is the unique decomposition of $V$ into essential, trivially inessential, and non-trivially inessential components,
	written $V = V_E \sqcup V_{TI} \sqcup V_{NI}$.\index{$V_E$}\index{$V_{TI}$}\index{$V_{NI}$}
\end{definition}

In an essential decomposition, $V_{TI} \in Q(M)$, and in fact $V_{TI}$ is trivial as it decomposes as a union of connected trivial components and pairs of connected components of the form $V_1 \sqcup -V_1$.
On the other hand, $V_{NI}$ may or may not be in $Q(M)$ depending on whether or not $\bd V_{NI}$
contains any essential components.
Nonetheless, a component $V_1$ of $V_{NI}$ is ``inessential to $V$'' in the sense that if $V_1'$ is another small rank element of $PC(M)$ with $\bd V_1 = \bd V_1'$ then $\underline{V_1} = \underline{V_1'}$ in $C_\Gamma^*(M)$ or $C_*^\Gamma(M)$ as $V_1 \sqcup -V_1' \in Q(M)$.
Thus $V_1$ is not required to appear as a component in a representative of $\uV$; for example, in any representative of $\uV$ we could replace each $V_1$ with $V_1'$ and still have the same (co)chain.
By contrast, \cref{L: essential}, which we will prove momentarily, shows that essential components really are essential in that they appear in any representation of $\uV$.

\begin{example}
	Let $V$ be any connected oriented manifold with corners that does not possess an orientation reversing diffeomorphism, and let $r_V \colon V \to M$ be any map that is not of small rank.
	Then $V = V_E$ is essential.
\end{example}

\begin{example}\label{E: bad transversality}
	Consider \cref{E: projected triangle}, which consisted of the projection of the 2-simplex $V \subset \R^2$ with vertices at $(1,0)$, $(-1,0)$, and $(0,1)$ to the $x$-axis.
	This map has small rank but is not trivial.
	In this case $V = V_{NI}$.

	We can make this example a bit more interesting as follows.
	Instead of the codomain being $\R^1$, we let the codomain be $\R^2$.
	We continue to let most of the map $r_V \colon V \to \R^2$ be the projection to the $x$-axis, but let us choose a Euclidean disk in the interior of the simplex and draw its image out into a 1-dimensional ``thread'' in the plane, analogously to how one changes basepoints for an element of some $\pi_2(M)$.
	With some care, this can be done smoothly.
	Any two such maps represent the same chain or cochain (depending on whether we assign orientations or co-orientations).
\end{example}

We now proceed with the following lemma, which will have several useful corollaries illuminating the structure of geometric chains and cochains.
These culminate in \cref{D: minimal essential,T: minimal rep}, which shows that chains and cochains have certain prechain and precochain representatives that are in some sense minimal, although they generally fail to be unique.
The theorem further identifies that non-uniqueness as being mainly vested in the $V_{NI}$.
We also highlight \cref{C: essential trans}, which concerns transversality and will be useful in \cref{S: products}.

\begin{lemma}\label{L: Q essential}
	Suppose $\bigsqcup_i W_i \in Q(M)$, and let $W_i = W_{i,E} \sqcup W_{i,TI} \sqcup W_{i,NI}$ be the essential decomposition for each $W_i$.
	Then $\bigsqcup_i W_{i,E}$ and $\bigsqcup_i W_{i,TI}$ are trivial and $\bigsqcup_i W_{i,NI} \in Q(M)$.
\end{lemma}

\begin{proof}
	We write the proof for cochains, but the argument for chains is equivalent.
	We know each $W_{i,TI}$ is trivial and thus in $Q(M)$, so by \cref{L: Lipy12}, $\bigsqcup_i (W_{i,E} \sqcup W_{i,NI}) \in Q(M)$.
	By definition, each $W_{i,NI}$ has small rank, while each component of $W_{i,E}$ does not have small rank.
	So in any decomposition of $\bigsqcup_i (W_{i,E} \sqcup W_{i,NI})$ into a trivial precochain and a degenerate precochain, which is possible as we know this is in $Q(M)$, the components of the $W_{i,E}$ must be part of the trivial precochain (the components of $W_{i, NI}$ may be part of either the trivial precochain or the degenerate precochain).
	By \cref{L: trivial structure}, each connected component, say $\mc W$, appearing in one of the $W_{i,E}$ either has a co-orientation-reversing self-diffeomorphism or appears zero times in all of $\bigsqcup_i (W_{i,E} \sqcup W_{i,NI})$ when counting with sign.
	If $\mc W$ has a co-orientation-reversing self-diffeomorphism, then $\mc W$ is trivial.
	Otherwise, for each occurrence of $\mc W$ in some $W_{i,E}$, there is an occurrence of $-\mc W$ in some $W_{j,E}$, and $\mc W \sqcup -\mc W$ is trivial.
	So altogether, $\bigsqcup_i W_{i,E}$ decomposes into a union of trivial precochains.
	It follows that $\bigsqcup_i W_{i,E} \in Q(M)$ is trivial, and so again by \cref{L: Lipy12}, $\bigsqcup_i W_{i,NI} \in Q(M)$.
\end{proof}

\begin{corollary}\label{C: Q essential}\index{minimial essential prechain or precochain|(}\index{geometric prechain!minimal essential|(}\index{geometric precochain!minimal essential|(}
	Let $\uW$ be an element of $C^*_\Gamma(M)$ (or $C_*^\Gamma(M)$).
	Then there is some prechain (or precochain) $Z$ so that for any prechain (or precochain) $W$ representing $\uW$, we have $W_E = Z \sqcup T$ for some trivial $T$.
	Furthermore, for every connected component $\mc W$ of $T$, either $\mc W$ or $-\mc W$ is a connected component of $Z$.
\end{corollary}

\begin{proof}
	Suppose $W_1, W_2$ both represent $\uW$.
	So $W_1 \sqcup -W_2 \in Q(M)$, and by \cref{L: Q essential}, $W_{1,E} \sqcup -W_{2,E}$ is trivial.
	Let $\mc W$ be a connected component of $W_{1,E} \sqcup -W_{2,E}$.
	Then $\mc W$ cannot itself be trivial since then it would be in $W_{1,TI}$ or $-W_{2,TI}$ instead of $W_{1,E}$ or $-W_{2,E}$.
	So by \cref{L: trivial structure}, $W_{1,E} \sqcup -W_{2,E}$ must consist of pairs $\mc W \sqcup -\mc W$.
	If there is a pair $\mc W \sqcup -\mc W$ in $W_{1,E}$, then we assign the pair to $T_1$ and similarly for $-W_{2,E}$.
	This leaves only pairs $\mc W \sqcup -\mc W$ so that $\mc W \in W_{1,E}$ and $-\mc W \in -W_{2,E}$ (reversing the sign of $\mc W$ if necessary).
	We let $Z$ be the disjoint union of the $\mc W$ in this last category, and then $W_{1,E} = Z \sqcup T_1$ and $W_{2,E} = Z \sqcup T_2$.
	We note that for any specific $\mc W$ that appears in $T_1$ or $T_2$, it must be the case that $\mc W$ also appears in $Z$, otherwise by definition either $W_{1,E}$ or $W_{2,E}$ would only contain $\mc W$ in cancelling pairs; but in that case these components would live in $W_{1,TI}$ or $W_{2,TI}$ by definition, not $W_{1,E}$ or $W_{2,E}$.

	We must show that this construction of $Z$ did not depend on our choices of $W_1$ and $W_2$.

	Next, consider $V = Z \sqcup W_{1,NI}$.
	From the preceding paragraph, we know that $W_{1,E} = Z \sqcup T_1$, so
	\begin{align*}
			W_1 \sqcup -V &= Z \sqcup T_1 \sqcup W_{1,NI} \sqcup W_{1,TI} \sqcup - (Z \sqcup W_{1,NI})\\
			&=Z \sqcup -Z \sqcup T_1 \sqcup W_{1,NI} \sqcup -W_{1,NI} \sqcup W_{1,TI},
	\end{align*}
	which is in $Q(M)$.
	So $V$ also represents $\uW$, and $V_E = Z$, as $Z$ contains no trivial or small rank components or cancelling pairs.
	But now let us compare $V$ with any other $W_3$, also representing $\uW$.
	Proceeding as in the first paragraph of the proof, there are $Z'$, $T_1'$, and $T_2'$ such that $V_E = Z = Z' \sqcup T'_1$ and $W_{2,E} = Z' \sqcup T'_2$.
	But by the construction of $Z$ above, $Z$ cannot contain any pair $\mc W \sqcup -\mc W$.
	So it follows from that construction that $T'_1 = \emptyset$, and so $Z' = Z$ and $W_{3,E} = Z \sqcup T'_2$.

	Thus any prechain or precochain representing $\uW$ has the form $Z \sqcup T$ for the same $Z$ as initially constructed and for some trivial $T$.
\end{proof}

\begin{definition}\label{D: minimal essential}\index{minimial essential prechain or precochain}
We refer to the prechain or precochain $Z$ of \cref{C: Q essential} as the \textbf{minimal essential prechain or precochain of $\uW$}.
It is constructed from the representative $W$ of $\uW$ by starting with $W_E$ and removing any trivial pairs of components $\mc W \sqcup - \mc W$.
\end{definition}

The preceding work now implies the following:

\begin{theorem}\label{T: minimal rep}
	Suppose $\uW \in C_*^\Gamma(M)$ or $\uW \in C^*_\Gamma(M)$ can be represented by $W = W_E \sqcup W_{TI} \sqcup W_{NI}$.
	Then $\uW$ can also be represented by $Z \sqcup W_{NI}$, where $Z$ is the minimal essential prechain or precochain of $\uW$.
\end{theorem}
\begin{proof}
	By \cref{C: Q essential}, we know $W_E = Z \sqcup T$ for some trivial $T$.
	So
	\begin{align*}
		Z \sqcup W_{NI} \sqcup -(W) &= 	Z \sqcup W_{NI} \sqcup - (W_E \sqcup W_{TI} \sqcup W_{NI})\\
		&= Z \sqcup W_{NI} \sqcup - (Z \sqcup T \sqcup W_{TI} \sqcup W_{NI})\\
		&= Z \sqcup -Z \sqcup W_{NI} \sqcup - W_{NI} \sqcup T \sqcup W_{TI}.
	\end{align*}
	As this consists of trivial objects and pairs, it is in $Q(M)$, which demonstrates the claim.
\end{proof}
\index{minimial essential prechain or precochain|)}\index{geometric prechain!minimal essential|)}\index{geometric precochain!minimal essential|)}

The following further corollaries will be useful below.
The first is an immediate consequence of \cref{C: Q essential}.

\begin{corollary}\label{L: essential}
	Let $\uV = \uW \in C^*_\Gamma(M)$ (or $\uV = \uW \in C_*^\Gamma(M)$) be represented by $V,W \in PC^*_\Gamma(M)$ (or $PC^\Gamma_*(M)$).
	Then $V_E \sqcup -W_E$ is trivial and if $\mc W$ is a connected component of $V_E$ then $\mc W$ or $-\mc W$ is a connected component $W_E$.
\end{corollary}

\begin{corollary}\label{C: essential trans}\index{essential!and transversality}
	Suppose $\uV = \uV' \in C^*_\Gamma(M)$ (or $C_*^\Gamma(M)$) and that $V_E$ is transverse to some $r_W \colon W \to M$.
	Then $V'_E$ is also transverse to $W$.
\end{corollary}

\begin{proof}
	By \cref{L: essential}, any component of $V_E$ is also a component of $V'_E$ (up to sign) and vice versa.
\end{proof}

\begin{corollary}\label{L: same NI}
	If $\uV = \uW \in C^*_\Gamma(M)$ (or $C_*^\Gamma(M)$) then
	$\underline{ V_{NI}} = \underline{ W_{NI}}$ and
	$\underline{\bd V_{NI}} = \underline{\bd W_{NI}}$.
\end{corollary}

\begin{proof}
	By assumption $V\sqcup-W = (V_E \sqcup V_{TI} \sqcup V_{NI}) \sqcup -(W_E \sqcup W_{TI} \sqcup W_{NI})$ is in $Q(M)$.
	We know $V_{TI} \sqcup -W_{TI} \in Q(M)$ as each of $V_{TI}$ and $W_{TI}$ is in $Q(M)$, and $V_E \sqcup -W_E \in Q(M)$ by \cref{L: essential}.
	Thus $V_{NI} \sqcup -W_{NI} \in Q(M)$ by \cref{L: Lipy12}, and the second equality follows by \cref{L: bd defined}.
\end{proof}
\index{essential decomposition|)}

\subsection{Splitting and creasing}\label{S: splitting and creasing}

In this section we discuss the closely related notions of splitting and creasing, both of which are concerned with breaking chains and cochains into smaller pieces.
The idea is somewhat analogous to the role subdivisions play in classical singular homology and cohomology theory, and both will be essential in our discussion of Mayer--Vietoris sequences.

The idea of splitting is to take a manifold with corners $W$ over a manifold without boundary $M$ and split it along a codimension one submanifold with corners into two manifolds with corners, $W^+$ and $W^-$.
This has various uses.
In particular, when $W$ represents a cycle or cocycle, the creasing construction then shows that $\uW$ is homologous (or cohomologous) to $\underline{W^+} + \underline{W^-}$.

\subsubsection{Splitting}\label{S: splitting}\index{splitting|(}
In this section, we begin by consolidating and extending some results previously encountered in \cref{E: manifold decomposition,S: codim 0 and 1 co-or}, namely the splitting of a manifold with corners over a manifold without boundary $r_W \colon W \to M$ into pieces $W^+$\index{$W^+$} and $W^-$\index{$W^-$} with common boundary $W^0$\index{$W^0$} using a function $\phi \colon M \to \R$\index{splitting function}.
In this subsection, we use our previous results to establish the facts we will need in the context of prechains and precochains, and in the next subsection we will utilize these spaces to perform creasing.

Our standing assumptions throughout this section will be that $M$ is a manifold without boundary and $\phi \colon M \to \R$ is a smooth map having $0$ as a regular value in the classical sense, i.e.\ for all $x \in \phi^{-1}(0)$ the differential $D_x\phi$ is nonzero.
We also assume an element of $PC_*^\Gamma(M)$ or $PC^*_\Gamma(M)$ represented by $r_W \colon W \to M$ such that $0$ is also a regular value for $\phi r_W$.
By \cref{D: regular value} this means that $\phi r_W$ is transverse to $0$, which is equivalent to assuming $0$ is a regular value in the classical sense for the restriction of $\phi r_W$ to each stratum of $W$.
By \cref{L: transverse to pullback}, $0$ is a regular value for $\phi r_W$ precisely when $r_W$ is transverse to $0 \times_\R M$, which is simply $\phi^{-1}(0)$, as observed in \cref{E: manifold decomposition}.

We let
\begin{align*}
	M^0 = \phi^{-1}(0) && M^- = \phi^{-1}((-\infty,0]) &&M^+ = \phi^{-1}([0,\infty)).
\end{align*}
Analogously,
\begin{align*}
	W^0 = (\phi r_W)^{-1}(0) && W^- = (\phi r_W)^{-1}((-\infty,0]) && W^+ = (\phi r_W)^{-1}([0,\infty)).
\end{align*}
We sometimes write $M^\pm$ for statements that could involve either $M^+$ or $M^-$, and similarly for $W$.

\begin{remark}
	For simplicity of notation, we primarily use $0$ as our value in $\R$ at which to perform splittings, though it should be clear that we could split $W$ using any regular value of $\phi$ and $\phi r_W$ to obtain analogous results.
\end{remark}

\begin{lemma}\label{L: pm0 as fiber products}\index{splitting!as fiber product}
	With our standing assumptions, there are diffeomorphisms
	\begin{align*}
		M^0 & \cong 0 \times_{\R} M &W^0 & \cong 0 \times_{\R} W & W^0& \cong M^0 \times_M W \\
		M^- & \cong (-\infty,0] \times_{\R} M &W^- & \cong (-\infty,0] \times_{\R} W& W^- &\cong M^- \times_M W\\
		M^+ & \cong [0,\infty) \times_{\R} M & W^+ & \cong [0,\infty) \times_{\R} W& W^+ &\cong M^+ \times_M W.
	\end{align*}
	In particular, these spaces are all manifolds with corners.
	Note that it is possible for some of these spaces to be empty.
\end{lemma}
\begin{proof}
	These diffeomorphisms are discussed in \cref{E: manifold decomposition}.
	The first two columns hold by direct observation, while the rightmost column is a consequence of \cref{P: pullback functoriality}.
\end{proof}

The rightmost column above demonstrates $W^0$, $W^-$, and $W^+$ as fiber products over $M$.
We can use these descriptions to realize these spaces as prechains or precochains.
Note: depending on $\phi$ and $r_W$, some of these spaces may be empty, in which case appropriate versions of the following statements hold vacuously.

\begin{lemma}\label{L: W0 cochain}\index{splitting!of geometric precochain}\index{splitting!boundary formula}
	Suppose $W \in PC^*_\Gamma(M)$, and $0$ is a regular value of $\phi \colon M \to \R$ and $\phi r_W \colon W \to \R$.
	Let the inclusions $M^\pm \into M$ have their tautological co-orientations, and give the inclusion $M^0 \into M$ the co-orientation determined by its normal vector field oriented by pulling back the positively-oriented normal vector field over $0$ in $\R$.
	Then $W^0$, $W^-$, and $W^+$, defined as the fiber products of $M^0$, $M^-$, and $M^+$ with $W$ over $M$, are elements of $PC^*_\Gamma(M)$.
	The co-orientations of $W^\pm \to M$ are the restrictions to $W^\pm$ of the co-orientation of $W$.
	Furthermore,
	\begin{align*}
		\bd(W^-) &= -(W^0) \bigsqcup (\bd W)^- \\
		\bd (W^+) &= W^0 \bigsqcup (\bd W)^+\\
		(\bd W)^0 &= -\bd (W^0).
	\end{align*}
\end{lemma}
\begin{proof}
	As $M^\pm$ and $M^0$ are all closed subsets of $M$, their inclusions are all proper maps, so with the co-orientations assigned above, they all represent elements of $PC^*_\Gamma(M)$.
	By \cref{L: transverse to pullback}, the map $r_W \colon W \to M$ is transverse to $M^0$ and $M^\pm$.
	It follows that $W^0$, $W^-$, and $W^+$, defined as the fiber products of $M^0$, $M^-$, and $M^+$ with $W$, are elements of $PC^*_\Gamma(M)$ by \cref{D: PC products,L: product preserves iso}.
	The co-orientations have been discussed previously in \cref{S: codim 0 and 1 co-or}, and in particular we have the boundary computations from \cref{E: codim 1 pullbacks,C: co-orient W0}.
\end{proof}

\begin{lemma}\label{L: W0 chain}\index{splitting!of geometric prechain}\index{splitting!boundary formulas}
	Suppose $W \in PC_*^\Gamma(M)$, and $0$ is a regular value of $\phi \colon M \to \R$ and $\phi r_W \colon W \to \R$.
	Let $\R$, $(-\infty, 0]$, and $[0,\infty)$ have their standard orientations, and give $0$ its positive orientation.
	Now orient $W^0$, $W^-$, and $W^+$ by realizing them as oriented fiber products of $0$, $(-\infty, 0]$, and $[0,\infty)$ with $\phi r_W \colon W \to \R$ over $\R$.
	Then the restrictions of $r_W \colon W \to M$ to $W^0$, $W^-$, and $W^+$ realize elements of $PC_*^\Gamma(M)$.
	When $W^\pm$ are nonempty, their orientations agree with the orientation of $W$.
	Furthermore, as elements of $PC_*^\Gamma(M)$, we have\footnote{Compare the signs with those in \cref{L: W0 cochain}.}
	\begin{align*}
		\bd(W^-) &= W^0 \bigsqcup (\bd W)^- \\
		\bd(W^+) &= (-W^0) \bigsqcup (\bd W)^+\\
		(\bd W)^0 &= -\bd (W^0).
	\end{align*}
\end{lemma}

\begin{proof}
	Since $W$ is compact, so will be $W^0$, $W^-$, and $W^+$ as closed subsets of $W$.
	So giving each space the orientation from the statement of the lemma, the restrictions of $r_W$ to $W^0$, $W^-$, and $W^+$ are elements of $PC_*^\Gamma(M)$.

	We consider a point $x$ in the interior of $W^+$ with an open neighborhood $N$ also in the interior of $W^+$ and given the orientation restricted from $W$.
	As the orientations of the fiber products are determined locally, the orientation of $W^+$ at $x$ will be consistent with its orientation in the restricted fiber product $[0,\infty) \times_\R N$, which will be the same as its orientation in $(0,\infty) \times_{(0,\infty)} N$.
	But this is the same as the initial orientation of $N$ as a subset of $W$ by \cref{P: oriented fiber product basic properties}.
	The same argument holds for $W^-$.

	We then have by \cref{P: oriented fiber boundary}, our conventions for oriented boundaries, and the standard computations for the boundaries of $(-\infty,0]$ and $[0,\infty)$ that
	\begin{equation*}
		\bd W^- = \bd ((-\infty,0] \times_\R W) = (0 \times_\R W) \sqcup ((-\infty,0] \times_\R \bd W) = W^0 \sqcup (\bd W)^-,
	\end{equation*}
	and
	\begin{equation*}
		\bd W^+ = \bd ([0,\infty) \times_\R W) = (-0 \times_\R W) \sqcup ([0,\infty) \times_\R \bd W) = (-W^0) \sqcup (\bd W)^+.
	\end{equation*}
	We also compute
	\begin{equation*}
		\bd W^0 = \bd (0 \times_\R W) = - 0 \times_\R \bd W = -(\bd W)^0.\qedhere
	\end{equation*}
\end{proof}

Lastly, we will occasionally need to consider the preimage of an interval $[p,q] \subset \R$.
If $p < q$ are regular values for $\phi$, then $\phi^{-1}([p,q])$ will be an embedded manifold with boundary in $M$ that we denote $M^{[p,q]}$.
If further $\phi r_W$ is transverse to $p$ and $q$ then we can form $W^{[p,q]} = M^{[p,q]} \times_M W$\index{$W^{[p,q]}$}.
We also let $W^p = \phi^{-1}(p) \times_M W$\index{$W^p$}, and similarly for $q$.

\begin{lemma}\label{L: Wpq cochain}\index{$W^{[p,q]}$!boundary formula}
	If $p < q$ are regular values for $\phi \colon M \to \R$, if $W \in PC^*_\Gamma(M)$, and if the inclusion $M^{[p,q]} \into M$ is given its tautological co-orientation, then
	$W^{[p,q]} \in PC^*_\Gamma(M)$ and
	$$\bd W^{[p,q]} = W^p \sqcup -W^q \sqcup (\bd W)^{[p,q]}.$$
	Similarly, if $W \in PC_*^\Gamma(M)$, then so is $W^{[p,q]}$, its orientation agrees with that of $W$, and
	$$\bd W^{[p,q]} = -W^p \sqcup W^q \sqcup (\bd W)^{[p,q]}.$$
\end{lemma}
\begin{proof}
	As $M^{[p,q]}$ is a closed subset of $M$, its inclusion into $M$ is proper, so with the tautological co-orientation, $M^{[p,q]} \into M$ is an element of $PC^*_\Gamma(M)$.
	We note that we only need to check transversality at $p$ and $q$, as the inclusion of the interior of $M^{[p,q]}$ to $M$ is transverse to any map.

	In the precochain case, as $M^{[p,q]}$ is a closed subset of $M$, its inclusion into $M$ is proper, so with the tautological co-orientation, $M^{[p,q]} \into M$ is an element of $PC^*_\Gamma(M)$.
	Thus the fiber product $W^{[p,q]} = M^{[p,q]} \times_M W$ is an element of $PC^*_\Gamma(M)$ by \cref{D: PC products,L: product preserves iso}, and the boundary formula follows from the Leibniz rule analogously to our computations above for $W^\pm$.

	In the prechain case, we observe that $W^{[p,q]} \cong (\phi r_W)^{-1}([p,q])$ by similar arguments as for $W^\pm$ in \cref{E: manifold decomposition}, so $W^{[p,q]}$ is compact.
	The agreement of orientation and the boundary formula are similar to the arguments of \cref{L: W0 chain}.
\end{proof}

\index{splitting|)}

\subsubsection{Creasing}\label{S: creasing}\index{creasing|(}

We now discuss the creasing construction of \cite[Section 2.4]{Lipy14}, though we use different orientation conventions and also consider versions involving co-orientations.

Suppose we are given a function $\varphi \colon W \to (-1,1)$ with $0$ a regular value, i.e.\ that $\varphi$ is transverse to the inclusion of $0$ into $(-1,1)$.
Below $\varphi$ will typically be a composition $W \xr{r_W} M \xr{\phi} (-1,1)$.
As in the preceding section, we define $W^+$\index{$W^+$} to be $\varphi^{-1}( [0,1))$, we define $W^-$ to be $\varphi^{-1} ((-1, 0])$\index{$W^-$}, and we define $W^0 = \varphi^{-1}(0)$\index{$W^0$}.
Note that here our codomain is $(-1,1)$ for reasons that will become apparent momentarily, though this is certainly still consistent with our maps being to $\R$ as above, using $(-1,1) \subset \R$.

By \cite[Lemma~9]{Lipy14}, the idea of creasing is that the topological space $W \times [0,1]$ can be given the structure of a manifold with corners, which we will write $\Cre(W)$\index{$C$@$\Cre(W)$}, satisfying $\bd \Cre(W) = W \sqcup -W^+ \sqcup -W^- \sqcup - \Cre(\bd W)$.
We call the manifold with corners $\Cre(W)$ the \textbf{creasing homotopy}\index{creasing homotopy|textbf} of $W$.
The creasing homotopy depends on $\varphi$, but we usually leave it tacit in the notation.
We will define $\Cre(W)$ by using pullbacks, which will provide an arguably simpler description of $\Cre(W)$ than found in \cite{Lipy14}.

To define $\Cre(W)$ via pullbacks, we need another map in addition to $\varphi$.
We let $D$ be the semi-open pentagonal region of the plane given by
$$D = \{(x,y) \in \R^2 \mid -1<x<1, 0 \leq y \leq 2-|x|\}.$$
Then $D$ is a manifold with corners with a smooth proper projection map $\pi \colon D \to (-1,1)$ given by $\pi(x,y) = x$.
We see that $\bd D$ has three pieces, say $D_x$, $D_+$, and $D_-$, corresponding respectively to the intersection of $D$ with the $x$-axis, the graph of $y = 2-x$ over $[0,1)$, and the graph of $y = 2+x$ over $(-1, 0]$.
We orient all three pieces by rightward pointing tangent vectors in the plane and give $D$ itself the standard planar orientation so that $\pi$ restricts to oriented diffeomorphisms from $D_x$, $D_+$, and $D_-$ onto their images in $(-1,1) = D_x$.
Then as oriented manifolds with corners, $$\bd D = D_x \sqcup -D_- \sqcup -D_+.$$

To obtain analogous boundary formulas for creasings of cochains, we let the projections of $D_x$, $D_+$, and $D_-$ to $(-1,1)$ be co-oriented by taking the rightward orientations to the rightward orientations, and $\pi \colon D \to (-1,1)$ to be co-oriented by taking the standard planar orientation to the rightward orientation.
In this case, as co-oriented manifolds with corners over $(-1,1)$, we have
$$\bd D = D_x \sqcup -D_- \sqcup -D_+.$$

The projection $\pi \colon D \to (-1,1)$ restricts to submersions from $D_x$, $D_+$, and $D_-$ onto their images, and so a map $\varphi \colon W \to (-1,1)$ from a manifold with corners is transverse to $\pi$ if and only if it is transverse to the map from the point at the tip of the pentagon to $(-1,1)$.
This is equivalent to the requirement that the restriction of $\varphi$ to every stratum of $W$ have $0$ as a regular value, i.e.\ that $0$ is a regular value for $\varphi$ by \cref{D: regular value}.

We can now officially define creasing.

\begin{definition}
	Suppose $W$ is a manifold with corners and $\varphi \colon W \to (-1,1)$ is a smooth map with a regular value at $0$.
	We define the \textbf{creasing homotopy}\index{creasing homotopy|textbf} to be the pullback
	\[
	\Cre(W) = D\times_{(-1,1)} W \to W.
	\]
	We note that $\Cre(W)$ does depend on the map $\varphi$, though we typically omit it from the notation as the choice of $\varphi$ will usually be clear from context.
	When necessary for clarity, we may write $\Cre_\varphi(W)$.

	Typically in practice, $W$ will arise with a map $r_W \colon W \to M$ representing an element of $PC_*^\Gamma(M)$ or $PC^*_\Gamma(M)$, and in this case our map $\varphi \colon W \to (-1,1)$ will generally be given as a composition $W \xr{r_W}M \xr{\phi} (-1,1)$ for some smooth $\phi \colon M \to (-1,1)$ with $0$ as a regular value.
	In this case, the condition that $\varphi$ have $0$ as a regular value is equivalent to the condition that $r_W$ be transverse to the submanifold $M^0 = \phi^{-1}(0)$ by \cref{L: transverse to pullback}, and as the creasing is governed by $\phi$, we sometimes write $\Cre_\phi(W)$ when we want to emphasize the dependence on $\phi$.
	When $r_W \colon W \to M$ represents an element of $PC_*^\Gamma(M)$, we treat $D\times_{(-1,1)} W$ as oriented by the convention for fiber products of oriented manifolds from \cref{S: orientation of fiber products}, using our fixed orientation of $D$, and we will show in a moment that if $W$ is compact so is $\Cre(W)$.
	Thus composing the pullback map $\Cre(W) \to W$ with $r_W$, we obtain another element of $PC_*^\Gamma(M)$ given by $D\times_{(-1,1)} W \to W \xr{r_W}M$.
	Similarly, when $r_W \colon W \to M$ represents an element of $PC^*_\Gamma(M)$, we treat $\Cre(W) \to W$ as co-oriented by the pullback conventions of \cref{D: pullback coorient}, using our given co-orientation of $D \to (-1,1)$, and this is a proper map by \cref{L: co-orientable pullback}.
	Thus composing the pullback map with the co-oriented map $r_W$, we obtain another element of $PC^*_\Gamma(M)$ given by $\Cre(W) \to W \xr{r_W} M$.
	In either case, we denote the composition $\Cre(W) \to W \xr{r_W}M$ by $r_{\Cre(W)}$.

	We will regularly abuse notation by allowing $\Cre(W)$ to refer to the space $D\times_{(-1,1)} W$, the pullback $D\times_{(-1,1)} W \to W$, or the element of $PC_*^\Gamma(M)$ or $PC^*_\Gamma(M)$ formed by the preceding constructions.
	It should usually be clear from context which is meant at any point.
\end{definition}

To justify the notion of creasing as a type of homotopy, at least topologically, we have the following lemma and corollary.
Note that these results are purely in the topological category, not the smooth category.

\begin{lemma}\index{fiber product!with projection}
	Suppose given a projection $\pi \colon X \times Y \to X$ and a map $g \colon W \to X$.
	Then the fiber product $(X \times Y) \times_X W$ is homeomorphic to $Y \times W$.
\end{lemma}

\begin{proof}
	Note that $g$ is transverse to $\pi$ as $\pi$ is a projection.
	We have $(X \times Y) \times_X W = \{(x,y,w) \in X \times Y \times W \mid x = g(w)\}$.
	A homeomorphism $(X \times Y) \times_X W \to Y \times W$ is then given by $(x,y,w) \mapsto (y,w)$ with inverse given by $(y,w) \mapsto (g(w),y,w)$.
\end{proof}

\begin{corollary}
	$\Cre(W)$ is homeomorphic to $I \times W$.
	In particular, if $W$ is compact then so is $\Cre(W)$.
\end{corollary}

\begin{proof}
	We observe that $D$ is homeomorphic to $(-1,1) \times I$ via a homeomorphism that preserves $\pi$ by taking $(x,y)$ to $\left(x, \frac{y}{2-|x|}\right)$.
	As a homeomorphism over a space $X$ induces a homeomorphism of fiber products, the corollary now follows from the preceding lemma.
\end{proof}

\begin{figure}[h!]
		\begin{tikzpicture}[scale=4]
	\coordinate (a) at (0,1.2);
	\coordinate (b) at (0,0.2);
	\coordinate (c) at (1.2,1);
	\coordinate (d1) at (1.1,.5);
	\coordinate (d2) at (1.4,.5);
	\coordinate (e) at (1.2,0);

	\draw[out=120, in=180] (b) to (a);
	\draw[out=60, in=-10, dotted] (b) to (a);

	\draw (a) -- (c);
	\draw (b) -- (e);

	\draw[out=180, in=130] (c) to (d1);
	\draw[out=-10, in=120] (c) to (d2);
	\draw[out=-100, in=120] (d1) to (e);
	\draw[out=-100, in=40] (d2) to (e);

	\draw (d1) to (d2);

	\node at (0, .7) {$W$};
	\node at (1.2, .75) {$W^+$};
	\node[scale=.7] at (1.23, .55) {$\phi^{-1}(0)$};
	\node at (1.22, .33) {$W^-$};
	\end{tikzpicture}
	\caption{Creasing homotopy}
	\label{F: creasing}
\end{figure}
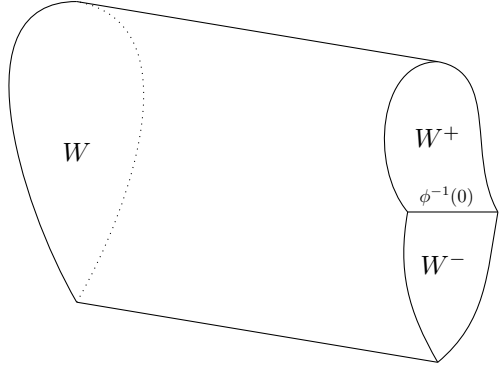

See \cref{F: creasing} for a sketch of a creasing homotopy of the teardrop manifold.

\begin{lemma}\label{E: bd crease}\index{$C$@$\Cre(W)$!boundary formula}
	Let $M$ be a manifold without boundary.
	Suppose given $\phi \colon M \to (-1,1)$ with $0$ a regular value and $r_W \colon W \to M$ transverse to $M^0 = \phi^{-1}(0)$, representing an element of $PC^\Gamma_*(M)$ or $W \in PC_\Gamma^*(M)$.
	Then
	\begin{equation*}
		\bd(\Cre_\phi(W)) = (D_x \sqcup -D_- \sqcup -D_+)\times_{(-1,1)} W \bigsqcup -D\times_{(-1,1)}\bd W = W \sqcup -W^- \sqcup -W^+ \sqcup -\Cre_\phi(\bd W),
	\end{equation*}
	again interpreting these formulas in $PC^\Gamma_*(M)$ or $PC_\Gamma^*(M)$, respectively, by composing the pullback maps to $W$ with $r_W \colon W \to M$ and the pullback maps to $\bd W$ with $\bd r_W \colon \bd W \to M$.
\end{lemma}

\begin{proof}
	The first equality comes from our boundary formulas for $D$ and our Leibniz rules for pullbacks from \cref{P: oriented fiber boundary,leibniz}.
	For the second equality we use, for example, that $\pi$ is an orientation preserving diffeomorphism over $(-1,1)$ from $D_+ \cong [0,1)$ onto its image in $(-1,1)$, and so we have an orientation-preserving diffeomorphism over $W$ given by $D_+\times_{(-1,1)} W \cong [0,1)\times_{(-1,1)} W = W^+$, and similarly in the co-oriented case.
	The argument is analogous for the $W^-$ and $W$ terms.
\end{proof}

To next promote the construction $\Cre(-)$ to an operator on $C_*^{\Gamma}(M)$ or $C^*_\Gamma(M)$, we will need the following corollary of \cref{L: pullback with Q}.

\begin{corollary}\label{C: creasing Q}\index{creasing!in $Q(M)$}
	Let $M$ be a manifold without boundary.
	Suppose given $\phi \colon M \to (-1,1)$ with $0$ a regular value and $r_T \colon T \to M$ transverse to $M^0 = \phi^{-1}(0)$.
	If $T \in Q(M)$, then so are $T^+$, $T^-$, $T^0$, and $\Cre(T)$.
\end{corollary}

\begin{proof}
	As observed in \cref{E: manifold decomposition}, with our assumptions about $\phi$ and $r_T$, the spaces $T^\pm$ and $T^0$ are fiber products over $M$ of $T$ with $M^\pm$ and $M^0$.
	So in this case the claim follows from \cref{L: pullback with Q}.

	For $\Cre(T)$, the pullback projection $\pi \colon \Cre(T) \to T$ is of small rank, as $\dim(\Cre(T))>\dim(T)$, from which it follows that $r_{\Cre(T)} = r_T\pi$ is of small rank.
	We also have $\bd \Cre(T) = T \sqcup - T^+ \sqcup -T^- \sqcup -\Cre(\bd T)$ by \cref{E: bd crease}.
	By assumption and the preceding paragraph, $T, T^\pm \in Q(M)$, and by the preceding sentence $\Cre(\bd T)$ is of small rank.
	Thus all components of $\bd \Cre(T)$ are trivial or of small rank, and so $\Cre(T)$ is degenerate.
\end{proof}

\begin{proposition}\index{creasing!of geometric (co)chains}\index{geometric chain!creasing}\index{geometric cochain!creasing}
	Let $M$ be a manifold without boundary.
	Suppose given $\phi \colon M \to (-1,1)$ with $0$ a regular value and that all creasing is done with respect to $\phi$.
	If $V$ and $W$ are any two representatives of $\uW \in C_*^\Gamma(M)$ whose reference maps are transverse to $M^0$, then $\Cre_\phi(V)$ and $\Cre_\phi(W)$ represent the same element of $C_*^\Gamma(M)$.
	Thus if the equivalence class $\uW$ contains any representative that is transverse to $M^0$, there is a well-defined element $\underline{\Cre(W)} \in C_*^\Gamma(M)$.
	Similarly for $C^*_\Gamma(M)$.
\end{proposition}

\begin{proof}
	The proofs for chains and cochains are the same, so we provide that with chains.

	If $r_V \colon V \to M$ and $r_W \colon W \to M$ represent the same class in $C_*^\Gamma(M)$ then $V \sqcup -W \in Q_*(M)$, and if $V$ and $W$ are both transverse to $M^0$ then so is $V \sqcup -W$.
	So $\Cre(V\sqcup-W)$ is defined, and by \cref{C: creasing Q}, $\Cre(V \sqcup -W) \in Q_*(M)$.
	Using our assumption that all creasing maps are compositions of the fixed $\phi \colon M \to (-1,1)$ with the reference maps and the properties of fiber product orientations, we have $\Cre(V \sqcup -W) = \Cre(V) \sqcup \Cre(-W) = \Cre(V) \sqcup -\Cre(W)$.
	Thus $\Cre(V)$ and $\Cre(W)$ represent the same element of $C_*^\Gamma(M)$.
\end{proof}

Finally, we come to the punchline of creasing:

\begin{theorem}\label{T: cohomology creasing}\index{creasing!in geometric (co)homology}
	Let $M$ be a manifold without boundary.
	Suppose $\uW \in H_*^\Gamma(M)$ and that $\uW$ has a representative $r_W \colon W \to M$ that is transverse to $M^0$, defined with respect to some $\phi \colon M \to (-1,1)$ with $0$ a regular value.
	Then $\uW = \underline{W^+} + \underline{W^-} \in H_*^\Gamma(M)$.
	Similarly for $H^*_\Gamma(M)$.
\end{theorem}

\begin{proof}
	Again the proofs for homology and cohomology are the same so we focus on homology.

	We have seen that, with our assumptions, $\uW$ yields a well-defined element $\underline{\Cre(W)}$ represented by $\Cre(W)$.
	Computing in $C_*^\Gamma(M)$ we have
	\begin{align*}
		\bd \underline{\Cre(W)}& = \underline{\bd \Cre(W) }\\
		& = \underline{W \sqcup -W^+ \sqcup -W^- \sqcup -\Cre(\bd W)}\\
		& = \uW -\underline{W^+}-\underline{W^-}.
	\end{align*}
	In the last line we have used that $\underline{\bd W} = 0$ so that $\bd W \in Q_*(M)$; hence $\Cre(\bd W) \in Q_*(M)$ by \cref{C: creasing Q} and so $\underline{\Cre(\bd W)} = 0 \in C_*^\Gamma(M)$.
	The theorem follows.
\end{proof}

\index{creasing|)}

\section{Basic properties of geometric (co)homology and equivalence with singular (co)homology}\label{S: basic properties}

In this section, we establish the basic properties of geometric homology and cohomology, including functoriality with respect to \textit{continuous} maps in \cref{S: functoriality} and Mayer-Vietoris sequences in \cref{S: MV}.
In \cref{S: homology is homology}, this allows us to show that geometric homology and cohomology are equivalent to (absolute) singular homology and cohomology on manifolds, using a result of Kreck and Singhof.
In the case of homology, we will also construct in \cref{S: homology direct} an explicit isomorphism from smooth singular homology to geometric homology.
It is more challenging to find an explicit isomorphism in the case of cohomology, and, in fact, the bulk of \cref{S: transversality} will be dedicated to such a construction.

\subsection{Functoriality and homotopy properties}\label{S: functoriality}\index{geometric homology!functoriality|(}\index{geometric cohomology!functoriality|(}

Given a \textit{continuous} map of smooth manifolds $f \colon M \to N$, we define in this section the induced maps $f_* \colon H_*^\Gamma(M) \to H_*^\Gamma(N)$ and $f^* \colon H^*_\Gamma(N) \to H^*_\Gamma(M)$ and show that they are independent of $f$ up to homotopy.
In subsequent sections, we will sometimes write the homology map simply as $f$ when the context is clear.

We also treat ``wrong way'' maps, though they require extra hypotheses.
In particular, if $f$ is proper and co-oriented it induces a covariant map on cohomology $f_* \colon H^*_\Gamma(M) \to H^{*+n-m}_\Gamma(N)$ independent of $f$ up to proper homotopy,
while if $f$ is proper and $M$ is oriented, we have a contravariant homology map $f^* \colon H_*^\Gamma(N) \to H_{*+m-n}^\Gamma(M)$.

In the covariant cases, $M$ and $N$ may both be manifolds with corners if we assume all maps and homotopies to be smooth (primarily to avoid smoothing arguments for maps of manifolds with corners).
In the contravariant case, we need to use both transversality and smoothing, and so for this case we consider only $M$ and $N$ without boundary.

The construction of $f_*$ for homology is relatively straightforward and can essentially be found in \cite[Section 6]{Lipy14}, though we provide additional details.
The covariant case of cohomology is similar to that for homology.
For the contravariant cases, there is slightly more work, and our argument here for cohomology parallels Kreck's in \cite{Krec10}, which is slightly different than the sketch in Lipyanskiy \cite[Section 6]{Lipy14} in that we choose to perturb $f$ rather than the reference map for the cohomology class in order to obtain transversality for the purpose of performing pullbacks.

As part of our constructions, we will see that if $F \colon M \times I \to N$ is a homotopy (proper or co-oriented as need be) and $r_W \colon W \to M$ represents a geometric cycle or cocycle, then the composition $F \circ (r_W \times \id_I)$ provides a homology or cohomology from $F(-,0)r_W \colon W \to N$ to $F(-,1)r_W \colon W \to N$.
Such homotopies, which we dub ``universal homotopies,'' will be critical in later sections for making geometric chains and cochains transverse to each other.

\subsubsection{Covariant functoriality of geometric homology and cohomology}\label{S: covariant functoriality}

In this section we consider the covariant behavior of geometric chains and cochains under maps and homotopies.
Both $M$ and $N$ may have corners in this section unless noted otherwise.
We begin with smooth maps but will generalize later to continuous maps via smooth approximation, at which time we will assume $M$ and $N$ are without boundary.

If $r_W \colon W \to M$ is in $PC_*^\Gamma(M)$ and $f \colon M \to N$ is a smooth map, then the composition $fr_W \colon W \to N$ is in $PC^\Gamma_*(N)$.
Consistent with our notation of writing $W \in PC_*^\Gamma(M)$, we write the image in $PC_*^\Gamma(N)$ as $f(W)$.
Similarly, if $f$ is proper and co-oriented we obtain a map $PC^*_\Gamma(M) \to PC^{*+n-m}_\Gamma(N)$ that we also write as $W \to f(W)$.
In this case the change of degree is because $W$ represents an element of degree $m-w$ in $PC^*_\Gamma(M)$ and of degree $n-w$ in $PC^{*}_\Gamma(N)$.
Functoriality is clear, and both $\bd(f(W))$ and $f(\bd W)$ are represented by $fr_Wi_{\bd W}$.
To obtain chain maps\footnote{We will not require here the signs that sometimes accompany chain maps of non-zero degree.} $C_*^\Gamma(M) \to C_*^\Gamma(N)$ and $C^*_\Gamma(M) \to C^{*+n-m}_\Gamma(N)$
(that we also write as $f$), it suffices to show that if $W \in Q(M)$ then $f(W) \in Q(N)$.
This is the content of the following lemma.

\begin{lemma}\label{L: Q preservation}
	If $r_W \colon W \to M$ represents an element of $Q_*(M)$ (or $Q^*(M)$) and $f \colon M \to N$ is any (co-oriented proper) smooth map, then $fr_W \colon W \to N$ is in $Q_*(N)$ (or $Q^{*+n-m}(N)$, co-orienting $fr_W$ with the composition co-orientation).
\end{lemma}

\begin{proof}
	First consider $W \in Q_*(M)$.
	By assumption, $W$ is the disjoint union of a trivial manifold over $M$ given by $r_T \colon T \to M$ and a degenerate manifold over $M$ given by $r_D:D \to M$.
	If $\rho \colon T \to T$ is an orientation-reversing diffeomorphism of $W$ such that $r_T\rho = r_T$ then also $fr_T\rho = fr_T$, so $fr_T \colon T \to N$ is also trivial.
	Furthermore, if $r_D$ has small rank then certainly so does $fr_D$.
	Similarly, $\bd(fr_D)$ is a union of trivial and small rank manifolds over $N$, and so $fr_D:D \to M$ is degenerate.

	The proof for $W \in Q^*(M)$ is the same using that the composition of proper maps is proper and the composition of co-oriented maps is co-oriented.
\end{proof}

\begin{corollary}[Lipyanskiy's Theorem 4]\label{C: homology chain map}
	Given a smooth map $f \colon M \to N$ of manifolds (possibly with corners), there is an induced chain map $C_*^\Gamma(M) \to C_*^\Gamma(N)$ given by $\uW \to \underline{f(W)}$, and this construction gives a covariant functor from the category of smooth manifolds and smooth maps to chain complexes over $\Z$.
\end{corollary}

\begin{corollary}\label{C: proper cofunctoriality}
	Given a smooth proper co-oriented map $f \colon M \to N$ of manifolds (possibly with corners), there is an induced chain map $C^*_\Gamma(M) \to C^{*+n-m}_\Gamma(N)$ given by $\uW \to \underline{f(W)}$, and this construction gives a covariant functor from the category of smooth manifolds and smooth maps to cochain complexes over $\Z$ with chain maps of degree $n-m$.
\end{corollary}
\index{geometric homology!functoriality|)}\index{geometric cohomology!functoriality|)}

Next we consider behavior with respect to homotopies.
For the remainder of the section, and for simplicity of notation, we leave degree shifts tacit in the notation; for example, for a chain complex $A_*$, we write $x \in A_*$ to mean an element of any degree and simply write $f \colon A_* \to B_*$ even when $f$ is a chain map of non-zero degree.

\begin{convention}[Naive prehomotopy co-orientation convention]\label{homotopy product co-orientation convention}\index{co-orientation!naive prehomotopy convention}
	To determine the interaction of homotopies $M \times I \to N$ with geometric chains and cochains, we will often need to consider compositions of the form $W \times I \xr{r_W \times \id_I} M \times I \to N$.
	The maps $W \times I \xr{r_W \times \id_I} M \times I$ are technically external products, but they function more as a sort of ``prehomotopy'' on the way to forming the composite homotopy $W \times I \to N$.
	As such, in this section (and when using its results later on) we follow a naive co-orientation convention for prehomotopies rather than the external product convention of \cref{D: co-oriented exterior} in order to avoid some unpleasant signs that are not consequential to the ultimate results.

	Specifically: If $r_W \colon W \to M$ is co-oriented, then we co-orient $r_W \times \id_I \colon W \times I \to M \times I$ so that if $(\beta_W,\beta_M)$ is the co-orientation at $x \in W$ then $(\beta_W \wedge \beta_I,\beta_M \wedge \beta_I)$ is the co-orientation at $(x,t)$ for any $t \in I$, with $\beta_I$ being the standard orientation of the interval.
\end{convention}

We begin with the following observations for the co-oriented cases.

\begin{lemma}
	Let $r_W \colon W \to M$ be co-oriented, and let $r_W \times \id \colon W \times I \to M \times I$ be co-oriented by \cref{homotopy product co-orientation convention}.
	For $j = 0,1$, let $i_j \colon W = W \times j \into W \times I$ and $k_j \colon M = M \times j \into M \times I$ be the inclusions, co-oriented with the usual boundary co-orientation (\cref{D: boundary co-orientation}).
	Then for $j = 0,1$, the following diagram of co-oriented maps commutes:
	\[
	\begin{tikzcd}[column sep=large]
		W \arrow[r, "r_W"] \arrow[hookrightarrow,d, "i_j"'] & M \arrow[hookrightarrow,d, "k_j"] \\
		W \times I \arrow[r, "r_W \times \id"] & M \times I.
	\end{tikzcd}
	\]
\end{lemma}

\begin{proof}
	Commutativity as maps is clear, so we focus on the co-orientations.

	Let $w \in W$, and suppose the co-orientation of $r_W$ at $w$ is represented by $(\beta_W, \beta_M)$.
	Let $\nu_W$ and $\nu_M$ denote inward pointing normal vector fields to $W \times j$ and $M \times j$ in $W \times I$ and $M \times I$, respectively. If $j=0$, then the co-orientation of $i_0$ at $w$ is $(\beta_W, \beta_W \wedge \beta_{\nu_W}) = (\beta_W, \beta_W \wedge \beta_I)$, where $\beta_I$ is the standard orientation of $I$, while if $j=1$, the co-orientation of $i_1$ is $(\beta_W, \beta_W \wedge \beta_{\nu_W}) = (\beta_W, \beta_W \wedge (-\beta_I))$. The co-orientations of $k_0$ and $k_1$ are analogous.

	So for $j=0$, the co-orientation of the composition down then right is $$(\beta_W, \beta_W \wedge \beta_I) * (\beta_W \wedge \beta_I, \beta_M \wedge \beta_I) = (\beta_W, \beta_M \wedge \beta_I).$$
	The composition right then down is $$(\beta_W, \beta_M) * (\beta_M, \beta_M \wedge \beta_I) = (\beta_W , \beta_M \wedge \beta_I).$$
	So we have commutativity.
	When $j=1$, the computation is the same except that $\beta_I$ is replaced with $-\beta_I$ in both formulas.
\end{proof}

\begin{corollary}\label{C: universal homotopy boundary co-orientation}
	Let $F \colon M \times I \to N$ be a co-oriented homotopy from $f$ to $g$ (see \cref{D: co-oriented homotopy}), and let $r_W \colon W \to M$ be a co-oriented map.
	If $r_W \times \id_I$ is co-oriented according to \cref{homotopy product co-orientation convention}, then $F \circ (r_W \times \id_I)$ is a co-oriented homotopy from $fr_W$ to $gr_W$.
\end{corollary}

\begin{proof}
	Again, the thing to check is the co-orientations. We continue to utilize the notation of the preceding lemma.

	By \cref{D: co-oriented homotopy}, the composition of the co-oriented boundary inclusion $k_0 \colon M = M \times 0 \into M \times I$ with $F$ is the co-oriented map $-f$, while the composition of the co-oriented boundary inclusion $k_1 \colon M = M \times 1 \into M \times I$ with $F$ is the co-oriented map
	$g$. On the other hand, the ``top and bottom'' components of the boundary of $F \circ (r_W \times \id_I)$ are the co-oriented compositions $F (r_W\times \id) i_0$ and $F (r_W\times \id) i_1$.
	But by the preceding lemma, for $j=0,1$, we have $F (r_W\times \id) i_j = F k_j r_W$.
	So when $j=0$, we have
	$$ F (r_W\times \id) i_0 = F k_0 r_W = -f r_W,$$
	and when $j=1$, we have
	$$ F (r_W\times \id) i_1 = F k_1 r_W = g r_W.$$

	So $F (r_W\times \id)$ is a co-oriented homotopy from $fr_W$ to $gr_W$ as desired.
\end{proof}

The next lemma shows that a homotopy $M \times I \to N$ of the above form cannot ``promote'' a manifold over $M$ that is in $Q(M)$ to one that is not in $Q(N)$.

\begin{lemma}\label{L: dessicated homotopy}
	Let $W \in PC_*^\Gamma(M)$ (or $W \in PC^*_\Gamma(M)$).
	Let $F \colon M \times I \to N$ be a smooth homotopy; if $W \in PC^*_\Gamma(M)$ suppose further that $F$ is proper and co-oriented.
	Then $F \circ (r_W \times \id_I) \colon W \times I \to N$ is in $PC_*^\Gamma(N)$ (or $PC^*_\Gamma(N)$, using \cref{homotopy product co-orientation convention}).
	Furthermore,
	if $r_W \colon W \to M$ is trivial, of small rank, or degenerate then so are
	$F \circ (r_W \times id_I) \colon W \times I \to N$ and $F(-,j) \circ r_W \colon W \to N$ for $j=0,1$ (for either co-orientation of $F(-,j)$, which is co-orientable by \cref{L: co-orientable homotopies}).
\end{lemma}

\begin{proof}
	The last statement concerning $F(-,j) \circ r_W \colon W \to N$ holds by \cref{L: Q preservation} and its proof.

	For	$F \circ (r_W \times id_I)$, we first suppose $W \in PC_*^\Gamma(M)$.
	In this case it is clear that $F \circ (r_W \times \id_I) \in PC_*^\Gamma(N)$ using the product orientation on $W \times I$.

	If $\rho \colon W \to W$ is an orientation-reversing diffeomorphism such that $r_W \circ \rho = r_W$, then
	$\rho \times \id_I \colon W \times I \to W \times I$ is an orientation-reversing self-diffeomorphism of $W \times I$ such that $F \circ (r_W \times id_I) \circ (\rho \times \id_I) = F \circ (r_W \times id_I)$.
	So $W \times I \to N$ is trivial.

	Next, if the derivative $Dr_W$ of the reference map $r_W$ has non-trivial kernel at a point $x \in W$ then so does the derivative $D(r_W \times \id_I)$ at each $(x,t) \in W \times I$, and thus so
	will $D(F \circ (r_W \times id_I)) = DF \circ D(r_W \times id_I)$.
	So $W \times I \to N$ has small rank if $r_W$ does.

	By definition, if $W \to M$ is degenerate then it has small rank and $\bd W = T \sqcup S$ with $T$ trivial and $S$ small rank.
	We have shown that $W \times I \xr{F\circ(r_W\times\id_I)} N$ then has small rank, so its suffices to consider its boundary,
	which is the union (up to signs) of $\bd W \times I$ and $W \times \bd I$.
	But $\bd W \times I = (T \times I) \sqcup (S \times I)$, which by our previous
	arguments are trivial and of small rank respectively.
	And the components of $W \times \bd I$ can be written for $j = 0, 1$ as $W \times j = W \xr{r_W} M \xr{F(-,j)} N$, and thus are of small rank since $r_W$ is.
	This completes the argument for $PC_*^\Gamma(M)$.

	If $W \in PC^*_\Gamma(M)$, then $F$ is proper and co-oriented by assumption, and \cref{homotopy product co-orientation convention} gives our co-orientation convention for $r_W \times \id_I$.
	To see the latter is proper, if $K$ is compact in $M \times I$, then $K \subset \pi_M(K) \times I$, which is also compact.
	Then $(r_W \times \id_I)^{-1}(K) \subset (r_W \times \id_I)^{-1}(\pi_M(K) \times I) = r_W^{-1}(\pi_M(K)) \times I$, which is compact as $r_W$ is proper.
	As compositions of proper co-oriented maps are proper and co-oriented, $F \circ (r_W \times \id_I) \colon W \times I \to N$ is well defined in $PC^*_\Gamma(N)$.
	The remaining arguments are analogous to the arguments above for $PC_*^\Gamma(M)$.
\end{proof}

\begin{corollary}\label{C: chain homotopy}
	Let $V, W \in PC_i^\Gamma(M)$ or $V, W \in PC^i_\Gamma(M)$.
	Let $F \colon M \times I \to N$ be a smooth homotopy; if $V,W \in PC^i_\Gamma(M)$ suppose further that $F$ is proper and co-oriented.
	If $V,W$ represent the same element of $C_i^\Gamma(M)$ and $j=0,1$, then $F(-,j) \circ r_V: V \to N$ and	$F(-,j) \circ r_W \colon W \to N$ represent the same element of $C_i^\Gamma(N)$ (for either chosen co-orientation of $F(-,j)$) and $F \circ (r_V \times id_I) \colon V \times I \to N$ and $F \circ (r_W \times id_I) \colon W \times I \to N$ represent the same element of $C_{i+1}^\Gamma(M)$, using \cref{homotopy product co-orientation convention} for the co-orientations.
	The analogous fact holds for $V, W \in PC^i_\Gamma(M)$.
\end{corollary}

\begin{proof}
	We know $V, W \in PC_i^\Gamma(M)$ represent the same element of $C_i^\Gamma(M)$ if and only if $V \sqcup -W \in Q_*(M)$, so the corollary follows from \cref{L: dessicated homotopy}.
	Similarly for cochains.
\end{proof}

\begin{corollary}\label{C: homotopy}
	Suppose $F \colon M \times I \to N$ is a smooth homotopy between maps $f,g \colon M \to N$.

	If $\uW \in C_*^\Gamma(M)$ is a cycle then $\underline{f(W)}$ and $\underline{g(W)}$ are homologous cycles in $N$.

	If $\uW \in C^*_\Gamma(M)$ is a cocycle and $F$ is proper and a co-oriented homotopy from $f$ to $g$ (see \cref{D: co-oriented homotopy}), then $\underline{f(W)}$ and $\underline{g(W)}$ are cohomologous cocycles in $N$.
\end{corollary}

\begin{proof}
	By the preceding corollary we may work with any representative $r_W \colon W \to M$ of $\uW$.
	First suppose $W \in PC_*^\Gamma(M)$.
	Identifying $W \times I$ with $W \times_{pt} I$, by \cref{P: oriented fiber boundary} we have $\bd (W \times I) = \bd W \times I \bigsqcup (-1)^{w} W \times \bd I$.
	As $W$ is a cycle, $\bd W \in Q_*(M)$, and hence so is $F \circ ((\bd r_W) \times \id_I) \colon \bd W \times I \to N$ by the preceding lemma.
	Thus in $C^\Gamma_*(N)$ the boundary of $F \circ (r_W \times \id_I) \colon W \times I \to N$ is represented up to signs by the restriction of
	$F \circ (r_W \times \id_I)$ to $W \times \bd I = W \times (\{1\} \sqcup \{-0\})$, applying our conventions for the boundary of $I$ with its standard orientation.
	As $F \circ (r_W \times \id_I)|_{W \times 1} = gr_W$ and $F \circ (r_W \times \id_I)|_{W \times 0} = fr_W$, we see that up to the overall sign $(-1)^{w}$, the boundary of $F \circ (r_W \times \id_I)$ is represented in $C_*^\Gamma(N)$ by the disjoint union of $gr_W \colon W \to N$ and the negative of $fr_W \colon W \to N$.
	Thus $fr_W \colon W \to N$ and $gr_W \colon W \to N$ represent homologous cycles.

	Next consider the case of $W$ a cocycle.
	Employing \cref{homotopy product co-orientation convention}, by \cref{C: universal homotopy boundary co-orientation,D: co-oriented homotopy} we have $$\bd (F \circ (r_W \times \id_I)) = gr_W \amalg -fr_W \amalg H,$$ where $H$ has underlying map $F \circ (r_W \times \id_I) \circ i_{\bd W \times \id_I}$.
	So, up to co-orientation, $H = F \circ (r_{\bd W} \times \id_I) \colon \bd W \times I \to N$.
	But $\bd W \in Q^*(M)$, and hence so is $H \colon (\bd W) \times I \to N$ by \cref{L: dessicated homotopy}.
	Thus we have
	that $fr_W \colon W \to N$ and $gr_W \colon W \to N$ represent cohomologous cocycles.
\end{proof}

We now come to the culmination of this section.
The first part of the following theorem is Lipyanskiy's Theorem 5 in \cite{Lipy14}.

\begin{theorem}\label{T: homology homotopy functor}\index{geometric homology!functoriality}\index{geometric cohomology!functoriality}\index{geometric homology!functoriality!homotopy}\index{geometric cohomology!functoriality!homotopy}
	Given a \textnormal{continuous} map $f \colon M \to N$ of manifolds without boundary\footnote{Possibly we can still let $M$ and $N$ be manifolds with corners and obtain a true statement, but we simplify our assumption here to avoid treating the question of smooth approximations in that setting.}, it induces a map $f_* \colon H_*^\Gamma(M) \to H_*^\Gamma(N)$ that depends only on the homotopy class of $f$.
	If $f$ is proper and co-oriented, it also induces $f_* \colon H^*_\Gamma(M) \to H^{*+n-m}_\Gamma(N)$, which depends only on the proper homotopy class and co-orientation of $f$.

	The constructions are both functorial.
	In other words:
	\begin{enumerate}
		\item If $\id \colon M \to M$ is the identity map, then $\id_*$ is the identity map on $H_*^\Gamma(M)$ and $H^*_\Gamma(M)$, using the tautological co-orientation of $\id$ for the latter.
		\item If $f_1 \colon M \to N$ and $f_2 \colon N \to X$ are continuous maps of manifolds without boundary, then $(f_2f_1)_* = f_{2*}f_{1*} \colon H_*^\Gamma(M) \to H_*^\Gamma(X)$; if $f_1$ and $f_2$ are also co-oriented and proper, then $(f_2f_1)_* = f_{2*}f_{1*} \colon H^*_\Gamma(M) \to H^*_\Gamma(X)$.
	\end{enumerate}
\end{theorem}

\begin{proof}
	First consider the case of homology.
	Let $g$ be any smooth approximation to $f$.
	Then by \cref{C: homology chain map}, $g$ induces a chain map $C_*^\Gamma(M) \to C_*^\Gamma(N)$ and hence a map $H_*^\Gamma(M) \to H_*^\Gamma(N)$.
	We show that this map is independent of the choice of $g$.
	Let $h$ be any other smooth map homotopic to $f$ (and so also homotopic to $g$).
	The continuous homotopy from $g$ to $h$ can be smoothly approximated by a smooth homotopy $H \colon M \times I \to N$ from $g$ to $h$ \cite[Theorem III.2.5]{Kos93}.
	The cycles represented by $g(W)$ and $h(W)$ are homologous by \cref{C: homotopy}.

	The cohomological case is the same by taking proper smooth approximations, which we show can be found in the proof of \cref{T: basic trans}, and co-orienting the homotopies using \cref{L: co-orientable homotopies} so that the approximations $g$ and $h$ are co-oriented homotopic to $f$.

	The functoriality is clear after replacing both $f_1$ and $f_2$ by (proper) smooth approximations.
\end{proof}

If $f \colon M \to N$ as above, we will often abuse notation by writing simply $f$ to denote the map on geometric chains or geometric homology, rather than $f_*$.

With functoriality of geometric homology established, we can introduce reduced homology; cf.\ \cite[page 200]{Span81}.

\begin{definition}\label{D: reduced}\index{geometric homology!reduced|textbf}
	The \textbf{reduced geometric homology} $\td H_*^\Gamma(M)$ of the manifold without corners $M$ is defined to be the kernel of the homomorphism $H_*^\Gamma(M) \to H_*^\Gamma(pt)$ induced by the unique map of $M$ to a point.
\end{definition}

\begin{remark}\label{R: reduced h}
Given the observations of \cref{E: first examples} and the splitting of surjective maps to free modules, we have the standard facts that $\td H_i^\Gamma(M) = H_i^\Gamma(M)$ when $i \neq 0$, while $$H_0^\Gamma(M) \cong \td H_0^\Gamma(M) \oplus H_0^\Gamma(pt) \cong \td H_0^\Gamma(M) \oplus \Z.$$

As for reduced singular homology, generators of $\td H_0^\Gamma(M)$ can be represented by prechains consisting of pairs of points with opposite orientations in separate components of $M$, as if they were in the same component then we could connect them with a path and they would represent $0$ in $\td H_0^\Gamma(M)$.

We recall that $M$ must have at most a countable number of components by the footnote on page \pageref{F: countable}, so let us write $M = \amalg M_i$ for connected components $M_i$ of $M$.
Then if $g_{ij}$, $i \neq j$, is a homology class represented by a prechain with a negative point in component $M_i$ and a positive point in component $M_j$, then we also have the relations $g_{ij}+g_{jk}=g_{ik}$.
So $\td H_0^\Gamma(M)$ can be identified with either the free abelian group on generators $g_{1i}$, $i>1$, or on generators $g_{i,i+1}$ for $i \geq 1$ so long as $i+1$ exists as an index.
\end{remark}

We can also now establish the standard additivity property.

\begin{proposition}\label{P: additive h}\index{geometric homology!additivity}
Let $M = \amalg{M_a}$, and let $i_a \colon M_a \into M$ be the inclusion.
Define $i \colon \oplus H_*^\Gamma(M_a) \to H_*^\Gamma(M)$ by $i(\oplus \underline{V_a}) = \sum i_a(\underline{V_a})$.
Then $i$ is an isomorphism.
\end{proposition}
\begin{proof}
Note that we have written an element of $\oplus H_*^\Gamma(M_a)$ as $\oplus \underline{V_a}$, but to be an element of $\oplus H_*^\Gamma(M_a)$, we tacitly assume that all but a finite number of classes $\underline{V_a}$ are trivial.
The map $i$ is then well defined.

Suppose $\uW \in  H_*^\Gamma(M)$ is represented by $r_W \colon W \to M$.
Let $W_a$ be the union of those components of $W$ that map to $M_a$, and let $r_{W_a} \colon W_a \to M_a$ be defined by restricting the codomain of $r_W|_{W_a}$.
As $W$ is compact, there must be a finite number of indices $a$ such that $W_a$ is non-empty.
Each $W_a$ represents a geometric cycle in $M_a$, else $W$ could not represent a cycle in $W$, and $i(\oplus \underline{W_a}) = \uW$, so $i$ is surjective.

Next, suppose $i(\oplus \underline{V_a})=0$ with $V_a \in PC_*^\Gamma(M_a)$ representing $\underline{V_a}$.
By \cref{R: cycles and boundaries}, there is a $W \in PC_*^\Gamma(M)$ such that $\bd W \sqcup (\sqcup_a -V_a) \in Q_*(M)$.
Let $W_a$ be as above.
Then we must have $\bd W_a \sqcup -V_a \in Q_*(M_a)$, which implies that $\underline{V_a}$ represents $0$ in $H_*^\Gamma(M_a)$.
So $i$ is injective.
\end{proof}

There is another very useful application of \cref{C: homotopy} that will be needed frequently below.
To explain, we first observe that if we have a map $r_W \colon W \to M$ representing an element of $H_*^\Gamma(M)$, then, in general, this homology class is not preserved under homotopies of $r_W$.
For example, a map representing a cycle might be homotopic to maps that do not represent cycles by pulling apart cancelling boundary components during the homotopy, or maps of small rank may be homotopic to maps that do not have small rank.
So, in general, geometric homology and cohomology do not behave well with respect to homotopies of the elements of $PC(M)$; this is the case also with ordinary singular chains modeled on simplices.
This will cause difficulties below, for example when we want to alter a map representing a cycle by a homotopy to make it transverse to some other cycle.
The solution will be to define and use what we call universal homotopies.

\begin{definition}\label{D: universal homotopy}\index{universal homotopy|textbf}
	Given two maps $f,g \colon W\to M$, we say that there is a \textbf{universal homotopy} from $f$ to $g$ if there is a smooth homotopy $H \colon M\times I\to M$ with $H(-,0)=\id_M$ and such that $g=H(-,1)\circ f$; in other words, if $f$ and $g$ are homotopic by a composition of the form $W \times I \xr{f\times \id} M \times I \xr{H} M$ with $H(-,0)$ the identity. In this case, we call the composition a \textbf{universal homotopy} from $f$ to $g$. We say that the universal homotopy is proper if $f$, $g$, and $H$ are proper maps.
\end{definition}

If one thinks of an ordinary homotopy from $f$ to $g$ as a way of obtaining $g$ by deforming $f$, then we think of a universal homotopy as deforming $f$ by first performing $f$ and then deforming $M$.

\begin{proposition}\label{P: universal homotopy}
	Suppose there is a universal homotopy $H$ from $f \colon W\to M$ to $g \colon W\to M$. If $f \colon W\to M$ represents an element of $H_*^\Gamma(M)$, then $g = H(-,1)f \colon W\to M$ represents the same element. Similarly, if $f \colon W\to M$ represents an element of $H^*_\Gamma(M)$ and the universal homotopy is proper, then $g \colon W\to M$ represents the same element, assuming we choose co-orientations for $g$ and $H$ such that $H$ is a homotopy from $f$ to $g$ in the sense of \cref{D: co-oriented homotopy}.
\end{proposition}

\begin{proof}
	The proposition is immediate from \cref{C: homotopy}, taking $M=N$ there and letting $F$ be the map $H \colon M \times I \to M$ of \cref {D: universal homotopy} that realizes the universal homotopy. In the case of cohomology, we note that we can always find such co-orientations; see \cref{L: co-orientable homotopies,D: homotopy co-orientation}.
\end{proof}

\subsubsection{Contravariant functoriality of geometric homology and cohomology}\label{S: cohomology pullback}\index{geometric homology!functoriality|(}\index{geometric cohomology!functoriality|(}

In this section, we assign to a continuous map $f \colon M \to N$ of manifolds without boundary a map $H^*_\Gamma(N) \to H^*_\Gamma(M)$.
If $f$ is proper and $M$ and $N$ are oriented, we also have a map $f^* \colon H_*^\Gamma(N) \to H_{*+m-n}^\Gamma(M)$.
As in the preceding section, we first consider a smooth map $f$ and then generalize to the continuous case using smooth approximations.
For simplicity of notation in what follows, we will not always explicitly write the degree shift for the homology map.

First, suppose $W \in PC^*_\Gamma(N)$ is represented by a proper co-oriented map $r_W \colon W \to N$, and let $g \colon M \to N$ be a \textit{smooth} map such that $g$ is transverse to $r_W$.
We emphasize that there is no need for assumptions that $g$ be co-oriented or proper.
We define $g^*(W)$\index{precochain!pullback of} to be the pullback $r_W \times_N g \colon W \times_N M \to M$, co-oriented by our standard convention from \cref{D: pullback coorient}.
The pullback is proper by \cref{L: co-orientable pullback}.
Also, using our standard convention for notating dimensions, the pullback has dimension $w+m-n$, and we have $m - (w+m-n) = n - w$. In other words, the index of $W$ as a precochain of $N$ is the same as the index of $W \times_N M$ as a precochain of $M$.
So, when defined, $g^*$ takes elements of $PC_\Gamma^i(N)$ to elements of $PC_\Gamma^i(M)$.

Similarly, if $W \in PC_*^\Gamma(N)$ and $g \colon M \to N$ is smooth and transverse to $r_W$ and if $M$ and $N$ are oriented, the pullback $W \times_N M$ has the pullback orientation of \cref{S: orientation of fiber products} and maps to $M$ by projection.
Furthermore, if $g$ is proper, then $W \times_N M$ is compact, as we show in the following lemma, so that we can pull back elements of $PC_*^\Gamma(N)$ to $PC_*^\Gamma(M)$ when they are transverse to $g$.\index{prechain!pullback of}

\begin{lemma}\label{L: compact pullback}
	Suppose $g \colon M \to N$ is a proper map from a manifold with corners to a manifold without boundary.
	Suppose $W$ is a compact manifold with corners and that $r_W \colon W \to N$ is transverse to $g$.
	Then $W \times_N M$ is compact.
\end{lemma}

\begin{proof}
	As $W$ is compact, so is $r_W(W)$, and as $g$ is proper, $g^{-1}(r_W(W))$ is compact.
	Now we observe that we must have $W \times_N M \subset W \times g^{-1}(r_W(W)) \subset W \times M$.
	So $W \times_N M$ is compact.
\end{proof}

The following lemma is similar to \cref{L: pullback with Q}, as is its proof, though there the focus was on fiber products, not pullbacks.
It will be useful later to allow $M$ to be a manifold with corners here.

\begin{lemma}\label{L: pullback map Q}\index{pullback!of trivial is trivial}\index{pullback!of small rank is small rank}
	Suppose $M$ is a manifold with corners and $N$ is a manifold without boundary.
	Suppose $r_S \colon S \to N$ is trivial (oriented or co-oriented) or has small rank, and let $g \colon M \to N$ be a smooth map transverse to $r_S$.
	Then the pullback $g^*(S) = S \times_N M \to M$ is also trivial or has small rank, respectively.
	Consequently, if $S \in Q(N)$ then $S \times_N M \to M$ is in $Q(M)$.
\end{lemma}

\begin{proof}
	If $\rho$ is a (co\nobreakdash-)orientation reversing diffeomorphism of $S$ over $N$, then the restriction of $\rho \times \id_M$ to $S \times_N M \subset S \times M$ is a (co\nobreakdash-)orientation reversing diffeomorphism of $S \times_N M$ over $M$ by \cref{L: product preserves iso}.

	Next, suppose $r_S \colon S \to N$ has small rank.
	Suppose at the point $x \in S$ we have $v \in \ker(D_x r_S)$.
	By \cref{L: tangent of pullbacks}, if $(x,y) \in S \times M$, the tangent space $T_{(x,y)}(S \times_N M)$ is the pullback $T_xS\times_{T_{(r_S(x),g(y))}N} T_y M$, and we are interested in the map from this space under the derivative of the projection $\pi \colon S \times M \to M$.
	As $v \in \ker (D_xr_S)$, we have $(v,0) \in T_{(x,y)}(S \times_N M)$, and this projects to $0$ in $T_yM$.
	So $(v,0) \in \ker(D_{(x,y)}\pi)$.
	This shows that the pullback has small rank.

	The final statement is then a consequence of the definitions as in the proof of \cref{L: Q preservation}.
\end{proof}

The requirement that $g \colon M \to N$ must be transverse to $W$ means that even with the preceding lemmas we cannot define chain maps $g^* \colon C^*_\Gamma(N) \to C^*_\Gamma(M)$ or $g^* \colon C_*^\Gamma(N) \to C_*^\Gamma(M)$ because for any fixed $g$ there may be geometric cochains none of whose representatives are transverse to $g$.
Nonetheless, given any continuous $f \colon M \to N$ of manifolds without boundary we can define a map in geometric cohomology $f^* \colon H^*_\Gamma(N) \to H^*_\Gamma(M)$, and if $M$ and $N$ are oriented and $f$ is proper, we can further define a map in geometric homology $f^* \colon H_*^\Gamma(N) \to H_{*+m-n}^\Gamma(M)$.
The constructions are as follows.

\begin{definition}\label{D: cohomology pullback and homology transfer}\index{geometric cohomology!pullback of}\index{geometric homology!pullback of}
	Suppose $f \colon M \to N$ is a continuous (not necessarily smooth, co-oriented, or proper) map of manifolds without boundary and $\uW \in H^*_\Gamma(N)$ is represented by $r_W \colon W \to N$.
	Let $g \colon M \to N$ be any smooth map homotopic to $f$ that is transverse to $r_W$ (which we know exists by \cref{T: basic trans}).
	Then we define $f^*(\uW) \in H^*_\Gamma(M)$ as the cohomology class represented by the pullback of $W$ by $g$, i.e.\ $g^*(W) = W \times_N M \to M$.

	If also $M$ and $N$ are oriented and $f$ is proper, and if $\uW \in H_*^\Gamma(N)$ is represented by $r_W \colon W \to N$, then let $g$ be a smooth proper map properly homotopic to $f$ and transverse to $r_W$ (also by \cref{T: basic trans}).
	Then we define $f^*(\uW) \in H_*^\Gamma(M)$ as the homology class represented by the pullback of $W$ by $g$, i.e.\ $g^*(W) = W \times_N M \to M$, with the pullback oriented as in \cref{S: orientation of fiber products}.
\end{definition}

\begin{proposition}\label{P: cohomology pullback}\index{geometric homology!functoriality}\index{geometric cohomology!functoriality}\index{geometric homology!functoriality!homotopy}\index{geometric cohomology!functoriality!homotopy}
	Suppose $f \colon M \to N$ is a continuous (not necessarily smooth, co-oriented, or proper) map of manifolds without boundary.
	The map $f^* \colon H^k_\Gamma(N) \to H^k_\Gamma(M)$ is well defined and depends only on the homotopy class of $f$.
	Similarly, if $M$ is oriented and $f$ is proper, then the map $f_* \colon H_k^\Gamma(N) \to H_{k+m-n}^\Gamma(M)$ is well defined and depends only on the proper homotopy class of $f$.
\end{proposition}

\begin{proof}
	We treat only the cohomology case, as the homology arguments are analogous taking the appropriate maps and homotopies to be proper, which is consistent with our applications of \cref{T: basic trans,T: homotopy trans} below.

	First, let $r_W \colon W \to N$ represent a cocycle, and suppose $g_0 \colon M \to N$ is smooth, homotopic to $f$, and transverse to $r_W$.
	By \cref{D: cohomology pullback and homology transfer}, $f^*(\uW)$ is represented by the pullback $g^*(W) = W \times_N M \to M$.
	We first observe this represents a cocycle.
	Because $M$ is without boundary, we have by \cref{leibniz} that $\bd(W\times_N M) = (\bd W) \times_N M$.
	Then as $\bd W \in Q^*(N)$ by assumption, we have $(\bd W) \times_N M \in Q^*(M)$ by \cref{L: pullback map Q}.
	Hence $g^*(W) = W \times_N M \to M$ represents a cocycle.

	Next suppose we have two maps $g_0, g_1 \colon M \to N$ that are smooth, homotopic to $f$, and transverse to $r_W$.
	Let $G'$ be a smooth homotopy from $g_0$ to $g_1$, which exists by \cite[Proposition 9.2.33]{MaDo92} (or, if $f$ is proper and $g_0$ and $g_1$ are smooth proper maps properly homotopic to $f$, there is a smooth proper such $G'$ by \cite[Proposition 9.2.35]{MaDo92}).
	Then by \cref{T: homotopy trans} there exists a smooth homotopy $G \colon M \times I \to N$ from $g_0$ to $g_1$ such that $G$ is transverse to $r_W$.
	In particular, we apply \cref{T: homotopy trans} with the map $f$ of that theorem being our smooth homotopy $G'$, whose boundary consists of $g_0$ and $g_1$; then the resulting map $h(-,1)$ of that theorem will be our desired $G$.

	Now, consider the composition $W \times_N (M \times I) \xr{r_W \times_N G} M \times I \to M$ of the pullback of $r_W$ by $G$ and the projection $M \times I \to M$.
	As $r_W$ is proper and co-oriented, so is the pullback, and the projection $M \times I \to M$ has its standard co-orientation $(\beta_M \wedge \beta_I,\beta_M)$ and is proper.
	So this composition represents a geometric cochain in $M$.
	Furthermore, the boundary of this composition is the disjoint union of $\bd W \times_N (M \times I) \to M \times I \to M$ and, up to sign, $$W \times_N \bd(M \times I) = \left( W \times_N (M \times \{1\})\right) \sqcup \left(- W \times_N (M \times \{0\})\right) \to M \times I \to M$$ as $M$ is a manifold without boundary.
	Since $W$ is assumed to be a cocycle, $\bd W \in Q^*(N)$.
	So $\bd W \times_N (M \times I) \to M \times I$ is in\footnote{Note that even if $M$ is without boundary, $M \times I$ will have boundary (unless $M$ is empty).
	This instance is the main reason we have allowed such constructions as $PC^*(-)$ and $Q^*(-)$ to take manifolds with non-empty boundary as inputs up to this point.
	In this particular case, since $N$ does not have boundary, all of our ``smooth'' maps to $N$ are automatically smooth in the strong sense of \cite{Joy12}, as are then $W \times_N (M \times I) \to M \times I$ and $M \times I = M \times_{pt} I \to M$ as these are both pullbacks in the smooth category.}
	$Q^*(M \times I)$ by \cref{L: pullback map Q}, and then by the proof of \cref{L: Q preservation} the projection to $M$ preserves triviality and small rank.
	Therefore, the composite $\bd W \times_N (M \times I) \to M \times I \to M$ is in $Q^*(M)$.
	The other terms correspond to the pullbacks of $W$ via $g_0$ and $g_1$, so $g_0^*(W)$ and $g_1^*(W)$ represent cohomologous cocycles.

	Finally, as we have noted that $f^*$ is not a fully-defined chain map, we must show directly that $f^*$ as a cohomology map does not depend on the choice of precocycle representing a cohomology class.
	So suppose $r_{W_0} \colon W_0 \to N$ and $r_{W_1} \colon W_1 \to N$ represent the same cohomology class.
	In this case there will be a $Z \in PC^*_\Gamma(N)$ such that $\bd \underline{Z} = \underline{W_0} - \underline{W_1}.$
	As $\bd \underline{Z} = \underline{\bd Z}$ and $\underline{W_0} - \underline{W_1} = \underline{W_0 \sqcup - W_1}$ by \cref{D: chains and cochains,L: co/chains well defined}, the above equality means by \cref{L: co/chains well defined} that $\bd Z \sqcup -(W_0 \sqcup -W_1) \in Q^*(N)$.
	By \cref{T: basic trans} and its proof, we may choose a smooth $g$ homotopic to $f$ that is transverse to all of $Z$ (and so also $\bd Z$), $W_0$, and $W_1$.
	Then by \cref{L: pullback map Q}, the pullback by $g$ given by
	$$\left[\bd Z \sqcup -(W_0 \sqcup -W_1)\right] \times_N M = ((\bd Z) \times_N M) \sqcup (-W_0 \times_N M) \sqcup (W_1 \times_N M)$$
	is in $Q^*(M)$.
	As $M$ has no boundary, $(\bd Z) \times_N M = \bd (Z \times_N M)$.
	It follows that $$\bd \underline{(Z \times_N M)} = \underline{(W_0 \times_N M)} \sqcup \underline{(-W_1 \times_N M)}.$$
	So $W_0 \times_N M$ and $W_1 \times_N M$ represent cohomologous cocycles, and $f^*$ does not depend on the choice of representing precocycle.

	Thus $f^*$ does not depend on our choices and depends on $f$ only through its homotopy class.
\end{proof}

\begin{corollary}\label{C: contra funct}\index{geometric homology!functoriality}\index{geometric cohomology!functoriality}\index{geometric homology!functoriality!homotopy}\index{geometric cohomology!functoriality!homotopy}
	The assignment taking $f$ to $f^*$ is functorial.
	In other words:
		\begin{enumerate}
			\item If $\id \colon M \to M$ is the identity map, then $\id^*$ is the identity map on $H^*_\Gamma(M)$; if $M$ is also oriented, then $\id^*$ is the identity map on $H_*^\Gamma(M)$.
			\item If $f_1 \colon M \to N$ and $f_2 \colon N \to X$ are continuous maps of manifolds without boundary, then $(f_2f_1)^* = f_1^*f_2^* \colon H^*_\Gamma(X) \to H^*_\Gamma(M)$; if $f_1$ and $f_2$ are also proper maps and $M$, $N$, and $X$ are oriented, then $(f_2f_1)^* = f_1^*f_2^* \colon H_*^\Gamma(X) \to H_*^\Gamma(M)$.
		\end{enumerate}
\end{corollary}
\begin{proof}
	We first consider cohomology.
	In this case, both statements follow from \cref{P: pullback functoriality}.
	For the second statement, given a cocycle represented by $r_W \colon W \to X$, we first apply \cref{T: basic trans} to obtain a smooth approximation $g_2$ to $f_2$ that is transverse to $r_W$ and then again to obtain a smooth approximation $g_1$ to $f_1$ that is transverse to the pullback $W \times_X N \to N$.

	For the homology case, we create the pullbacks the same way using that \cref{T: basic trans} allows us to approximate proper maps by proper maps.
	The properties then follow from \cref{P: oriented fiber product basic properties} and \cref{P: oriented fiber mixed associativity}, in the latter case taking the map labeled $g$ in the statement of that proposition to be the identity.
\end{proof}

\begin{example}\label{E: contractible}
By \cref{E: dimension}, $H^i_\Gamma(pt)=0$ for $i\neq 0$, while it follows from \cref{E: first examples} and \cref{T: PD} that  $H^0_\Gamma(pt) \cong \Z$, generated by the class of the identity $pt \to pt$ with either co-orientation.
It follows now from the homotopy functoriality of pullbacks that if $D$ is any contractible manifold then similarly $H^i_\Gamma(D)=0$ for $i\neq 0$ and $H^0_\Gamma(D) \cong \Z$.
The representatives of generators of $H^0_\Gamma(D)$ are the pullbacks of the generators of $H^0_\Gamma(pt)$ by the homotopy equivalence $D \to pt$.
These are the identity maps $\id \colon D \to D$ with either co-orientation.

Next consider the projection $M \times D \to M$, which is a homotopy equivalence for any manifold $M$ and any contractible manifold $D$.
Then $H^i_\Gamma(M \times D) \cong H^i_\Gamma(M)$, and if $r_W \colon W \to M$ (with some co-orientation) represents a cohomology class in $H^i_\Gamma(M)$, the corresponding class in $H^i_\Gamma(M \times D)$ is represented by the map $W \times D \to M \times D$ given by $(w,d) \mapsto (r_W(w),d)$, with its pullback co-orientation.
\end{example}
\index{geometric homology!functoriality|)}\index{geometric cohomology!functoriality|)}

The functoriality of geometric cohomology allows us to establish the standard additivity property for cohomology.
The argument is analogous to that for \cref{P: additive h}.

\begin{proposition}\label{P: additive c}\index{geometric cohomology!additivity}
Let $M = \amalg{M_a}$, and let $i_a \colon M_a \into M$ be the inclusion.
Define $i^* \colon H^*_\Gamma(M) \to \prod_a H^*_\Gamma(M_a)$ by $i^*(\uV) = \prod i^*_a(\uV)$.
Then $i^*$ is an isomorphism.
\end{proposition}
\begin{proof}
We observe that if $r_V \colon V \to M$ represents a cochain, then $i^*_a(\uV)$ is simply the restriction of $r_V$ to those connected components of $V$ that map to $M_a$, which we write as $r_{V_a} \colon V_a \to M_a$.

So first consider $\prod_a \underline{W_a} \in \prod_a H^*_\Gamma(M_a)$ with each $\underline{W_a}$ represented by $r_{W_a} \colon W_a \to M_a$.
As noted in the footnote on page \pageref{F: countable}, $M$ and each representative precochain $W_a$ must have a countable number of components, and so $W = \sqcup W_a$ can be embedded in $\R^\infty$.
We define $r_W \colon W \to M$ by letting $r_W = r_{W_a}$ on $W_a$, and $r_W$ is proper since it maps properly to each $M_a$.
Then $i^*(\uW) = \prod i^*_a(\uW) = \prod \underline{W_a}$.
So $i^*$ is surjective.

Next, suppose $i^*(\uV) = 0$.
Let $V \in PC^*_\Gamma(M)$ represent $\uV$, and let $V_a$ be the union of those connected components of $V$ that map to $M_a$.
By \cref{R: cycles and boundaries}, there are $W_a \in PC^*_\Gamma(M_a)$ such that $\bd W_a \sqcup -V_a \in Q^*(M_a)$.
We let $W = \sqcup W_a$ as above, and it follows that $(\sqcup_a \bd W_a) \sqcup (\sqcup -V_a) = W \sqcup -V \in Q^*(M)$.
So $V$ represents $0$ in $H^*_\Gamma(M)$, and $i^*$ is injective.
\end{proof}

We can now discuss reduced cohomology; cf.\ \cite[page 240]{Span81}.

\begin{definition}\label{D: reduced c}\index{geometric cohomology!reduced}
The \textbf{reduced geometric cohomology} $\td H^*_\Gamma(M)$ of the manifold without corners $M$ is defined to be the cokernel of the homomorphism $H^*_\Gamma(pt)\to H^*_\Gamma(M)$ induced by the unique map of $M$ to a point.

As $H^i_\Gamma(pt) = 0$ for $i \neq 0$, we have the standard fact that  $\td H^i_\Gamma(M) = H^i_\Gamma(M)$ for $i \neq 0$.
When $i=0$, the generator of $H^0_\Gamma(pt) \cong \Z$, which is represented by the identity map of the point, maps to the class of $H^0_\Gamma(M)$ represented by the pullback of identity map of the point, which is the identity map of $M$ with its tautological co-orientation.
As noted in \cref{E: first examples}, we will show below that when $M$ is connected this map generates $H^0_\Gamma(M) \cong \Z$.
So by \cref{P: additive c}, we have $H^0_\Gamma(M) \cong \prod \Z$, where the product is indexed over the components of $M$ with each factor generated by the identity map of that component with its tautological co-orientation.
Then $\td H^0_\Gamma(M)$ is the quotient of this $\prod \Z$ by the subgroup generated by $(1,1,\ldots)$.
\end{definition}

\subsection{Mayer--Vietoris sequences}\label{S: MV}

In this section, we consider Mayer--Vietoris sequences for homology and cohomology.
In particular, we present covariant homology and cohomology sequences and a contravariant cohomology sequence.

Throughout our treatment, we will assume that $U$ and $V$ are open subsets of a manifold without boundary $M$.
Note that any of $U$, $V$, or $U \cap V$ may be empty.

For the covariant homology sequence, the maps $H_*^\Gamma(U) \to H_*^\Gamma(M)$, etc., induced by inclusion will be those defined in \cref{S: covariant functoriality}.
For a covariant cohomology sequence, however, the maps of the preceding section will not work, as the inclusion of an open set into a manifold is not generally proper.
Rather, to obtain a covariant cohomology sequence we will need to use in at least some of the terms a variant of geometric cohomology supported on open subsets, denoted $H^*_\Gamma(M|_U)$.
We will not provide an in-depth study of this variant cohomology, but it will not take us too far afield to introduce it here and treat its Mayer--Vietoris sequence in parallel with the homology sequence.
We justify this added effort by observing that we will have $H^*_\Gamma(M|_M) = H^*_\Gamma(M)$, so this sequence may be useful in studying our usual geometric cohomology.
However, the reader more interested in our standard geometric homology and cohomology can bypass this next section and go right to the following section.

\subsubsection{Cohomology supported on open sets}

In this brief section we provide the definition and some immediate properties for geometric cohomology supported on open sets.

\begin{definition}\index{geometric cohomology!supported on an open set|textbf}
	Let $M$ be a manifold without boundary and $U \subset M$ an open subset.
	Define $C_\Gamma^*(M|_U) \subset C_\Gamma^*(M)$\index{$C_\Gamma^*(M \pipe_U)$} to be the subcomplex consisting of elements of $C_\Gamma^*(M)$ that can be represented by $r_W \colon W \to M$ in $PC^*_\Gamma(M)$ with the image of $W$ contained in $U$.
	Let $H_\Gamma^*(M|_U) = H^*(C_\Gamma^*(M|_U))$.\index{$H_\Gamma^*(M \pipe_U)$}
\end{definition}

It is easy to observe that $C_\Gamma^*(M|_U)$ is a chain complex as it is closed under addition and taking boundaries using that the sum $\underline{W_1}+\underline{W_2}$ of two elements of $C^*_\Gamma(M|_U)$ represented by elements of $PC^*_\Gamma(M)$ with image in $U$ can be represented by the disjoint union $W_1 \sqcup W_2$ with image in $U$, and similarly the boundary of a precochain with image in $U$ has image in $U$.

In fact, we can reformulate this definition as follows.

\begin{lemma}
	Let $PC_\Gamma^*(M|_U) \subset PC_\Gamma^*(M)$ be the subset consisting of $r_W \colon W \to M$ with image in $U$, and let $Q^*(M|_U)$ be the elements of $Q^*(M)$ with image in $U$.
	Then $C_\Gamma^*(M|_U)$ is isomorphic to the set of equivalence classes of $PC_\Gamma^*(M|_U)$ under the relation $V\sim_U W$ if $V \sqcup -W \in Q^*(M|_U)$.
\end{lemma}

\begin{proof}
	We first observe that $\sim_U$ is an equivalence relation on $PC_\Gamma^*(M|_U)$ by the same proof as \cref{L: cancel Q}, assuming all maps have image in $U$.
	Furthermore, if $V,W \in PC_\Gamma^*(M|_U)$ and $V\sim_UW$ then $V\sim W$ in the sense of \cref{L: cancel Q}.
	So letting $\hat C_\Gamma^*(M|_U)$ temporarily denote the
	equivalence classes under $\sim_U$, we have a well-defined map $f \colon \hat C_\Gamma^*(M|_U) \to C_\Gamma^*(M|_U)$.
	By definition, elements of $C_\Gamma^*(M|_U)$ can be represented by elements of $PC_\Gamma^*(M|_U)$, so $f$ is surjective.
	Now suppose $W_1,W_2 \in PC^*_\Gamma(M|_U)$ represent the same element of $C^*_\Gamma(M|_U)$.
	Then $W_1\sim W_2$ in the sense of \cref{L: cancel Q}, i.e.\ $W_1 \sqcup -W_2 \in Q^*(M)$, but $W_1$ and $W_2$ each have support in $U$, so $W_1 \sqcup -W_2 \in Q^*(M|_U)$ and $W_1\sim_U W_2$.
	So $f$ is injective.
\end{proof}

The cohomology groups $H_\Gamma^*(M|_U)$ are functorial in the sense that if $U \subset W \subset M$ are open sets then we have $C_\Gamma^*(M|_U) \subset C_\Gamma^*(M|_W) \subset C_\Gamma^*(M|_M) = C_\Gamma^*(M)$, and these induce maps $H_\Gamma^*(M|_U) \to H_\Gamma^*(M|_W) \to H_\Gamma^*(M)$.

As Mayer--Vietoris sequences are often the engines of induction arguments and as many inductions start with Euclidean balls, we provide the following computation, which is akin to the dimension axiom or the Poincar\'e lemma; see \cref{E: dimension}.

\begin{proposition}\label{P: P lemma}\index{geometric cohomology!supported on an open set!for an open ball}
	Let $M$ be an $m$-dimensional smooth manifold, and $U \subset M$ an open set that is diffeomorphic to an open ball and whose closure in $M$ is compact.
	Then $H^n_\Gamma(M|_U)$ is $\Z$ if $n = m$ and is zero otherwise.
\end{proposition}

\begin{proof}
	If $r_W \colon W \to M$ is proper and has image in $U$, then $W$, as the preimage of the closure of $U$, must be compact.
	Also, $U$ is orientable, and so by choosing an orientation of $U$ every co-oriented $r_W \colon W \to U$ determines an orientation of $W$ as explained in \cref{S: co-orientations}.
	Thus $C^*_\Gamma(M|_U) \cong C_{m-*}^\Gamma(U)$, and the result follows from \cref{E: dimension} (the dimension axiom), homotopy functoriality of geometric homology, and that $U$ is homotopy equivalent to a point.
\end{proof}

\subsubsection{Covariant Mayer--Vietoris sequences}

We now turn to our covariant Mayer--Vietoris sequences.

A key tool in proving the existence of these sequences will be creasing, which we introduced in \cref{S: creasing}.
For geometric homology and cohomology, creasing replaces the role that subdivision plays in the classical singular theories.
In some sense creasing is simpler, as creasing only needs to be applied once while subdivision often needs to be iterated.
In order to perform creasing, the following definitions will be useful.

\begin{definition}\label{D: separating function}
	Let $U,V$ be open subsets of a manifold $M$.
	We will call a smooth function $\phi \colon U \cup V \to [-1/2,1/2]$ a \textbf{separating function for $\mathbf U$ and $\mathbf V$}\index{separating function|textbf} if $U - V \subset \phi^{-1}(-1/2)$ and $V - U \subset \phi^{-1}(1/2)$.
	As $U  - V$ and $V  - U$ are disjoint closed subsets of $U \cup V$ and all of our spaces are normal by assumption, such a function can always be found by an application of Urysohn's Lemma to find a continuous function with this property and then applying the Smooth Approximation Theorem \cite[Theorem III.2.5]{Kos93}.
	We use $\pm 1/2$ as our endpoints to be consistent with the maps $\phi \colon M \to (-1,1)$ used in the definition of creasing.
	Note that we only assume that $U  - V \subset \phi^{-1}(-1/2)$ and not that $U  - V = \phi^{-1}(-1/2)$, and similarly for $V  - U$.
	However, it will be useful to keep in mind that $\phi^{-1}([-1/2,1/2)) \subset U$ and $\phi^{-1}((-1/2,1/2]) \subset V$.

	Suppose now $\phi$ is a separating function for $U$ and $V$ and that we have a map $r_W \colon W \to U \cup V$.
	We will say that $\phi$ \textbf{separates $\mathbf W$ over $\mathbf U$ and $\mathbf V$} if $0$ is a regular value for $\phi$ and for the composite $\phi r_W$, meaning that $\phi$ and $\phi r_W$ are transverse to the inclusion of $0$ into $[-1/2,1/2]$.
	In particular, $0$ is a regular value in the usual sense for the restriction of $\phi r_W$ to each stratum of $W$.
	Consequently, if $\phi$ is separating for $W$ over $U$ and $V$, we may perform creasing of $W$ along $(\phi r_W)^{-1}(0)$ as defined in \cref{S: creasing}.
	More generally, we can form $W^+$, $W^-$, and $W^0$ as discussed in \cref{E: manifold decomposition,S: codim 0 and 1 co-or}.
	Note that some of these spaces may be empty depending on the configuration of $U$ and $V$.
\end{definition}

\begin{lemma}\label{L: existence of separating}\index{separating function!existence of}
	Given $U,V \subset M$ as above and $r_W \colon W \to U \cup V$.
	There exists a function separating $W$ over $U$ and $V$.
	More generally, for any countable collection $r_{W_i} \colon W_i \to U \cup V$, there is a $\phi$ that is separating for all $W_i$ over $U$ and $V$.
\end{lemma}

\begin{proof}
	Let $\phi$ be a separating function for $U$ and $V$; we have already observed these exist.
	By Sard's Theorem, when $M$ is any smooth manifold without boundary the set of regular values of $\phi \colon M \to [-1/2,1/2]$ is a set whose complement has measure zero.
	So for any countable set of maps $W_j \to U \cup V$ from smooth manifolds without boundary we can find a value $p \in (-1/2, 1/2)$ that is a regular value for $\phi$ and for all $\phi r_{W_j}$, as a countable union of sets of measure zero still has measure $0$.
	In particular, we can take the $W_j$ to be the strata of a single manifold with corners or even all the strata of a countable collection of manifolds with corners.
 Once we have found such a common regular value $p$, we can then replace $\phi$ by its composition with an orientation-preserving diffeomorphism of $[-1/2,1/2]$ that takes $p$ to $0$, for example the linear fractional transformation $x \mapsto \frac{x-p}{-4px+1}$.
	This composition is our desired function.
\end{proof}

We are now ready to demonstrate the existence of covariant Mayer--Vietoris sequences.
In the following statements, we let $i_U \colon U \cap V \to U$ and $i_V \colon U \cap V \to V$ be the inclusion maps, and we write $(i_U, -i_V)$ for the map $C_*^{\Gamma}(U \cap V) \to C_*^{\Gamma}(U) \oplus C_*^{\Gamma}(V)$ that takes $\uW$ to $(i_U(\uW), -i_V(\uW))$.
More generally, we use the notation $(a,b)$ for elements of a direct sum.

\begin{theorem}\label{T: relative MV}\index{Mayer-Vietoris sequence!covariant}
	Let $M$ be a manifold without boundary.
	For any pair of open sets $U$ and $V$ in $M$ there are Mayer--Vietoris exact sequences
	\begin{equation*}
		\begin{tikzcd}[
			column sep = small,
			arrow style = math font,
			cells = {nodes = {text height = 2ex,text depth = 0.75ex}}
			]
			\cdots \arrow[r] & H^{\Gamma}_{k}(U \cap V) \arrow[r, "(i_U{,} -i_V)"] & [25pt] H^{\Gamma}_k(U) \oplus H^{\Gamma}_{k}(V)
			\arrow[r] & H_{k}^{\Gamma}(U \cup V) \arrow[r] & H_{k-1}^{\Gamma}(U \cap V) \arrow[r] & \cdots
		\end{tikzcd}
	\end{equation*}
	and
	\begin{equation*}
		\begin{tikzcd}[
			column sep = small,
			arrow style = math font,
			cells = {nodes = {text height = 2ex,text depth = 0.75ex}}
			]
			\cdots \arrow[r] & H_{\Gamma}^{k}(M|_{U \cap V}) \arrow[r, "(i_U{,} -i_V)"] & [25pt] H_{\Gamma}^k(M|_U) \oplus H_{\Gamma}^{k}(M|_V)
			\arrow[r] & H^{k}_{\Gamma}(M|_{U \cup V}) \arrow[r] & H^{k+1}_{\Gamma}(M|_{U \cap V}) \arrow[r] & \cdots.
		\end{tikzcd}
	\end{equation*}
\end{theorem}

\begin{proof}
	The proof parallels standard proofs of the existence of Mayer--Vietoris sequences for singular homology.
	The proofs for homology and cohomology are analogous, so we give the cohomological case, which is slightly more exotic.

	Let $S^*$ denote the quotient of $C_{\Gamma}^*(M|_U) \oplus C_{\Gamma}^{*}(M|_V)$ by the image of $C_{\Gamma}^*(M|_{U \cap V})$ under the map $(i_U, -i_V)$, with $i_U$ and $i_V$ being the inclusions.
	Then we have a short exact sequence
	\begin{equation}\label{E: homology MV SES}
		\begin{tikzcd}[
			column sep = small,
			arrow style = math font,
			cells = {nodes = {text height = 2ex,text depth = 0.75ex}}
			]
			0 \arrow[r] &
			C_{\Gamma}^{*}(M|_{U \cap V}) \arrow[r, "(i_U{,}-i_V)"] &[25pt]
			C_{\Gamma}^*(M|_U) \oplus C_{\Gamma}^{*}(M|_V) \arrow[r] &
			S^* \arrow[r] &
			0.
		\end{tikzcd}
	\end{equation}
	This short exact sequence generates a long exact cohomology sequence, and the theorem will follow from showing there is a quasi-isomorphism $\psi \colon S^* \to C_{\Gamma}^*(M|_{U \cup V})$.
	Our quasi-isomorphism will be induced by the map $C_{\Gamma}^*(M|_U) \oplus C_{\Gamma}^{*}(M|_V) \to C_{\Gamma}^*(M|_{U \cup V})$ that takes $(\underline{W_1}, \underline{W_2})$ to $\underline{W_1}+\underline{W_2}$, letting context determine whether we think of $\underline{W_1}$ as an element of $C_{\Gamma}^*(M|_U)$ or $C_{\Gamma}^*(M|_{U \cup V})$ and similarly for $\underline{W_2}$.
	This induces a well-defined map $\psi \colon S^* \to C_{\Gamma}^*(M|_{U \cup V})$ as it takes elements in the image of $(i_U,-i_V)$ to $0$.
	To establish the quasi-isomorphism, we use creasing.

	First suppose a cocycle $\uW \in C_{\Gamma}^*(M|_{U \cup V})$ represented by $r_W \colon W \to U \cup V$.
	Let $\phi \colon U \cup V \to [-1/2,1/2]$ separate $W$ over $U$ and $V$.
	By \cref{T: cohomology creasing}, we have $\uW = \underline{W^-}+\underline{W^+} \in H^*_\Gamma(M)$.
	But by the creasing construction, if $W$ has image in $U \cup V$ then so does $\Cre(W)$, so the equality $\uW = \underline{W^-}+\underline{W^+}$ also holds in $H_*^\Gamma(M|_{U \cup V})$.
	Similarly, $\underline{W^-} \in C_{\Gamma}^*(M|_{U})$ and $\underline{W^+} \in C_{\Gamma}^*(M|_{V})$.
	We also have $\bd W \in Q^*(M|_{U \cup V}) \subset Q^*(M)$ by assumption that $W$ is a cocycle.
	By the computation in \cref{E: codim 1 pullbacks}, we have $\bd(W^-) = -(W^0) \sqcup (\bd W)^-$, and $\bd (W^+) = W^0 \sqcup (\bd W)^+$.
	As $\bd W \in Q^*(M)$ with image in $U \cup V$, we have $(\bd W)^\pm \in Q^*(M)$ with respective images in $U$ and $V$ by \cref{C: creasing Q}, so $(\bd W)^-$ represents in $0$ in $C^*_\Gamma(M|_U)$ and similarly $(\bd W)^+$ represents $0$ in $C^*_\Gamma(M|_V)$.
	Thus $(\underline{W^-}, \underline{W^+})$ is an element of $C_{\Gamma}^*(M|_U) \oplus C_{\Gamma}^{*}(M|_V)$ whose boundary is $(-\underline{W^0},\underline{W^0})$, which is in the image of $(i_U,-i_V)$.
	Therefore, $(\underline{W^-}, \underline{W^+})$ represents an element of
	$H^*(S^*)$ that maps to $\uW \in H^*_\Gamma(M_{U \cup V})$.
	Thus $\psi$ is surjective.

	Next, suppose $\uW$ is a cocycle in $S^*$ that maps to zero in $H_{\Gamma}^*(M|_{U \cup V})$, and that $\uW$ is represented by $(\underline{W_1},\underline{W_2}) \in C^*_\Gamma(M|_U) \oplus C^*_\Gamma(M|_V)$.
	As $\uW$ maps to zero in $H_{\Gamma}^*(M|_{U \cup V})$, there is some $Z \in PC^*_\Gamma(M|_{U \cup V})$ such that $\bd Z \sqcup - W_1 \sqcup -W_2 = T \in Q^*(M|_{U \cup V})$.
	Using that $\underline T = 0$ in $C^*_\Gamma(M|_{U \cup V})$, it will be useful to write this as $\bd \underline{Z} = \underline{\bd Z} = \underline{W_1} + \underline {W_2}$ in $C^*_\Gamma(M|_{U \cup V})$.
	By \cref{L: existence of separating} there is a function $\phi \colon U \cup V \to [-1/2,1/2]$ that is separating for $Z$, $W_1$, and $W_2$ over $U$ and $V$.
 	We claim that $(\underline{Z^-}+\underline{\Cre(W_1)},\underline{Z^+}+\underline{\Cre(W_2)}) \in C^*_\Gamma(M|_U) \oplus C^*_\Gamma(M|_V)$ represents an element of $S^*$ whose boundary is $\uW$.
	We observe again that $Z^-$ and $Z^+$ represent respective elements of $PC^*_\Gamma(M|_U)$ and $PC^*_\Gamma(M|_V)$, as they are closed subsets of $Z$, so their inclusions to $Z$ composed with $Z \to M$ are proper, and we utilize the tautological co-orientations of $Z^\pm \into Z$ as in \cref{E: codim 0 and 1 co-or as fiber products}.

	As $\bd Z \sqcup - W_1 \sqcup -W_2 \in Q^*(M|_{U \cup V})$, \cref{C: creasing Q} implies that $(\bd Z)^- \sqcup (- W_1)^- \sqcup (-W_2)^- \in Q^*(M|_U)$ and $(\bd Z)^+ \sqcup (- W_1)^+ \sqcup (-W_2)^+ \in Q^*(M|_V)$, so $\underline{(\bd Z)^- }= \underline{W_1^-} + \underline{W_2^-}$ and $\underline{(\bd Z)^+} = \underline{W_1^+} + \underline{W_2^+}$ in $C^*_\Gamma(M|_U)$ and $C^*_\Gamma(M|_V)$, respectively.
	So we can compute using \cref{E: codim 1 pullbacks,E: bd crease}:
	\begin{align*}
		\bd (\underline{Z^-}&+\underline{\Cre(W_1)},\underline{Z^+}+\underline{\Cre(W_2)})\\
		&= (\bd (\underline{Z^-})+\bd \underline{\Cre(W_1)},\bd (\underline{Z^+})+\bd\underline{\Cre(W_2)})\\
		& = (-\underline{Z^0} + \underline{(\bd Z)^-}+\underline{W_1} -\underline{W_1^-} - \underline{W_1^+} -\underline{\Cre(\bd W_1)},\underline{Z^0} + \underline{(\bd Z)^+} +\underline{W_2} -\underline{W_2^-} -\underline{W_2^+} -\underline{\Cre(\bd W_2)})\\
		& = (-\underline{Z^0} + \underline{W_1^-} + \underline{W_2^-} +\underline{W_1} -\underline{W_1^-} -\underline{W_1^+} -\underline{\Cre(\bd W_1)},\underline{Z^0}+ \underline{W_1^+}+ \underline{W_2^+} +\underline{W_2} -\underline{W_2^-} -\underline{W_2^+} -\underline{\Cre(\bd W_2)})\\
		& = (-\underline{Z^0} + \underline{W_2^-} + \underline{W_1} -\underline{W_1^+} -\underline{\Cre(\bd W_1)},\underline{Z^0}+ \underline{W_1^+} + \underline{W_2} -\underline{W_2^-} -\underline{\Cre(\bd W_2))})\\
		& = (\underline{W_1},\underline{W_2})+ (-\underline{Z^0},\underline{Z^0}) + (\underline{W_2^-}-\underline{W_1^+}, \underline{W_1^+} -\underline{W_2^-}) -(\underline{\Cre(\bd W_1)},\underline{\Cre(\bd W_2)}).
	\end{align*}
	To see that this is equal to $\uW = (\underline{W_1},\underline{W_2})$ in $S^*$, we first note the middle two terms are each in the image of $(i_U, -i_V)$, so they represent $0$ in $S^*$.
	It remains to show that $(\underline{\Cre(\bd W_1)},\underline{\Cre(\bd W_2)})=0$ in $S^*$.

 This term is obtained by applying the creasing construction to the entries of $(\bd W_1,\bd W_2)$, which by assumption represents $0$ in $S^*$.
	In general, suppose $A \in PC^*_{\Gamma}(M|_U)$ and $B \in PC^*_\Gamma(M|_V)$ such that $(\underline{A}, \underline{B}) \in C^*_\Gamma(M|_U) \oplus C^*_\Gamma(M|_V)$ represents zero in $S^*$.
	This means that $(A,B)$ is equivalent in $C^*_\Gamma(M|_U) \oplus C^*_\Gamma(M|_V)$ to an element of the form $(i_U,-i_V)(\underline{C})$ for some $C \in PC^*_{\Gamma}(M|_{U \cap V})$, i.e.\ that $A \sqcup -i_U(C) \in Q^*(M|_U)$ and $B \sqcup i_V(C) \in Q^*(V)$.
	Creasing preserves both support and membership in $Q^*$ by definition and by \cref{C: creasing Q}, so we have
	$$\Cre(A \sqcup -i_U(C)) = \Cre(A) \sqcup -\Cre(i_U(C)) = \Cre(A) \sqcup -i_U(\Cre(C)),$$
	and this is an element of $Q^*(M|_U)$.
	The equivalent computations holds for the other term in $V$.
	So in $C^*_\Gamma(M|_U)\oplus C^*_\Gamma(M|_U)$, we have
 $$0= (\underline{\Cre(A)} -\underline{i_U(\Cre(C))}, \underline{\Cre(B)} + \underline{i_V(\Cre(C))}) = (\underline{\Cre(A)},\underline{\Cre(B)}) +(i_U(\underline{\Cre(C)}), -i_V(\underline{\Cre(C)})).$$
	Consequently, $(\underline{\Cre(A)},\underline{\Cre(B)})$ is in the image of $(i_U, -i_V)$ and so represents $0$ in $S^*$.

	Applying this argument to $(\underline{\Cre(\bd W_1)},\underline{\Cre(\bd W_2))})$, we conclude that
	$\bd (\underline{Z^-} + \underline{\Cre(W_1)},\underline{Z^+} + \underline{\Cre(W_2)}) = (\underline{W_1},\underline{W_2}) = \uW$ in $S^*$, and so $\psi$ is injective.

	The one difference in the homological case is that in our above computation of $$\bd (\underline{Z^-} + \underline{\Cre(W_1)},\underline{Z^+} + \underline{\Cre(W_2)})$$ the term $(-\underline{Z^0},\underline{Z^0})$ will instead be $(\underline{Z^0},-\underline{Z^0})$ due to the difference in signs in the boundary formulas in \cref{L: W0 cochain,L: W0 chain}.
	However, this does not alter the overall argument.
\end{proof}

\begin{remark}\label{R: MV boundary}
	By the computations above and the standard zig-zag construction of long exact sequences from short exact sequences, the connecting homomorphism in our Mayer--Vietoris sequence takes the cohomology class represented by the cocycle $\uW \in C^*_\Gamma(M|_{U \cup V})$ to the class represented by $-\underline{W^0} \in C^*_\Gamma(M|_{U \cap V})$, where $W^0$ is determined by a function that separates $W$ over $U$ and $V$.

	The description for the homology sequence is similar, except that, using \cref{L: W0 chain}, the connecting homomorphism takes the homology class represented by the cycle $\uW \in C_*^\Gamma(U \cup V)$ to the class represented by $\underline{W^0} \in C_*^\Gamma(U \cap V)$.
\end{remark}

\begin{example}[Homology of spheres]\label{E: sphere homology}\index{geometric homology!of spheres}
Given the dimension axiom (\cref{E: dimension}), homotopy invariance, and Mayer-Vietoris sequences, the geometric homology groups of spheres can now be computed formally to agree with the singular homology of spheres, i.e.\ for $n>0$ we have $H_0^{\Gamma}(S^n)\cong H_n^\Gamma(S^n)\cong \Z$ and all other homology groups are $0$.
By \cref{E: first examples}, a generator of $H_0^{\Gamma}(S^n)$ is represented by the inclusion of any point.
By the same example, a generator of $H_n^{\Gamma}(S^n)$ is represented by the oriented identity map $\id \colon S^n \to S^n$, although the argument there for general closed, oriented manifolds relies on results we have not yet reached.
For spheres, however, we can see this by induction.

By applying creasing at the equator $S^{n-1}\subset S^n$, we can write the homology class of $\id: S^n \to S^n$ as the union of the inclusions of the closed northern and southern hemispheres. Then writing $S^n$, $n>1$, as the union of the two overlapping disks $U$ and $V$ that extend the hemispheres in the usual way, we see by a diagram chase that the Mayer-Vietoris boundary map $H_n^\Gamma(S^n)\xr{\cong} H_{n-1}^\Gamma(U \cap V)$ is an isomorphism with the image of the class of $\id: S^n \to S^n$ represented, up to sign, by the inclusion of the equator.
But $U \cap V$ has a deformation retraction to the equator $S^{n-1}$.
So if we assume inductively that $\id \colon S^{n-1} \to S^{n-1}$ represents a generator of $H_{n-1}^\Gamma(S^{n-1})\cong \Z$, the same will be true for $H_n^\Gamma(S^n)$.
The base case for $n=1$ is a bit more involved as in this case $H_1^\Gamma(S^1) \to H_{0}^\Gamma(U \cap V) \cong H_0^\Gamma(S^0) \cong \Z^2$ is only injective, not an isomorphism, but the necessary exact sequence computations in degree $0$ to show that the boundary map of the sequence takes the class represented by $\id \colon S^1 \to S^1$ onto a generator of a $\Z$ summand is nonetheless a standard exercise that we leave for the reader.
\end{example}

\subsubsection{Contravariant Mayer--Vietoris sequence}

In this section, we show that geometric cohomology possesses a contravariant Mayer--Vietoris sequence on manifolds.
This does not seem possible by deriving a long exact sequence from a short exact sequence.
For example, in general the restriction map $C^*_\Gamma(M) \to C^*_\Gamma(U)$ will not be surjective.
Consequently, this sequence will take more work than the covariant ones.

Instead, we proceed with a direct analysis at each term of the exact sequence as in Kreck's proof of the Mayer--Vietoris sequence for his cohomology theory using stratifolds in \cite{Krec10}.
In fact, our argument for exactness parallels Kreck's fairly closely.
However, our proof that the connecting map is well defined is more complicated than the analogous proof in Kreck (ignoring the extra complications in \cite{Krec10} arising from considerations of collars that we do not need).
This is because Kreck is able to define his version of separating functions directly on his representing objects as $\varphi \colon W \to [0,1]$, which then makes it possible to compare two different $\varphi$s using a cylinder $W \times I$ with different versions of $\varphi$ on each end.
Here, however, we must always use splitting maps that factor through $M$ or risk that even if $W$ is trivial $W^{\pm}$ may not be.

Throughout, there will be no loss of generality in taking $U \cup V = M$, and it will simplify notation, though we often continue to write $U \cup V$ for emphasis.

\medskip
\noindent\textbf{The connecting morphism.} We first define the connecting morphism and demonstrate its naturality.
This will require quite a bit of technical work.

\begin{definition}\label{D: connecting}\index{Mayer-Vietoris sequence!for geometric cohomology!connecting morphism}\index{Mayer-Vietoris sequence!for geometric cohomology!connecting morphism|(}\index{connecting morphism for geometric cohomology|(}
	Suppose that $U$ and $V$ are open subsets of a manifold without boundary $M = U \cup V$ and that $\uW \in H^k_\Gamma(U \cap V)$ is represented by $r_W \colon W \to U \cap V$.
	Of course we can also consider $r_W$ to have image in $U \cup V$.
	Let $\phi \colon U \cup V \to [-1/2,1/2]$ separate $W$ over $U$ and $V$.
	Define $\delta(\uW) \in H^{k+1}_\Gamma(U \cup V) = H^{k+1}_\Gamma(M)$ to be represented by $-W^0 = -(\phi r_W)^{-1}(0)$, which we recall can be identified with $-\phi^{-1}(0)\times_M W$ by \cref{L: pm0 as fiber products}.
\end{definition}

The choice to use $-W^0$ rather than $W^0$ in the definition is explained by \cref{R: MV boundary}.
An alternative convention that uses $W^0$ rather than $-W^0$ would be to use $(-i_U,i_V)$ rather than $(i_U,-i_V)$ in our covariant Mayer--Vietoris sequences and the map $\uW \to (-\uW|_U, \uW|_V)$ rather than $\uW \to (\uW|_U, -\uW|_V)$ in our contravariant Mayer--Vietoris sequence below.

Our next goal is to show that $\delta$ is well defined, i.e.\ that $W^0$ is a precocycle in $M = U \cup V$ and that its cohomology class does not depend on the choice of representing precochain $W$ or the choice of $\phi$.
We begin with a series of lemmas.
Although we use $0$ as our regular value in our construction of $\delta$, it will be useful in our arguments to allow our regular value to be anywhere in $(-1/2,1/2)$, and we immediately work in this generality.
So for $p \in (-1/2,1/2)$, we let $M^p = \phi^{-1}(p)$ and $W^p = (\phi r_W)^{-1}(p) \cong M^p \times_{M} W$.
Note that as $\phi$ separates $U$ and $V$, we have $M^p \subset U \cap V \subset U \cup V = M$ for all $p \in (-1/2,1/2)$, so the space $M^p$ is unambiguous whether we consider the domain of $\phi$ to be $U \cup V$ or $U \cap V$.

We first show our construction does give an element of $H^*_\Gamma(M)$.

\begin{lemma}
	Suppose that $U$ and $V$ are open subsets of a manifold without boundary $M = U \cup V$ and that $\uW \in H^k_\Gamma(U \cap V)$ is represented by $r_W \colon W \to U \cap V$.
	Consider $r_W$ to have image in $U \cup V = M$.
	Let $\phi \colon U \cup V \to [-1/2,1/2]$ separate $U$ and $V$, and let $p \in (-1/2, 1/2)$ be a regular value for $\phi$ and $\phi r_W$.
	Then $W^p = M^p \times_{M} W$ is an element of $PC^*_\Gamma(M)$ with boundary in $Q^*(M)$, so it represents an element of $H^*_\Gamma(M)$.
\end{lemma}

\begin{proof}
	For simplicity and consistency with the notation in \cref{S: codim 0 and 1 co-or,S: splitting}, we take $p = 0$, but the same argument works for any $p \in (-1/2, 1/2)$ satisfying the hypotheses.
	By \cref{L: W0 cochain}, $M^0 \times_{U\cap V} W$ is an element of $PC^*_\Gamma(U \cap V)$, but we must show that the map including this fiber product into $M$ is proper.

	We know that $W \to U \cap V$ is proper by hypothesis and that the inclusion $M^0 \into M$ has image in $U \cap V$, so by \cref{L: co-orientable pullback} the pullback $W \times_{U \cap V} M^0 \to M^0$ is proper.
	And as $M^0$ is a closed subset of $M$, the inclusion $M^0 \into M$ is proper.
	So the composition map $W \times_{U \cap V} M^0 \to M^0 \into M$ is proper.
	But now clearly $W \times_{U \cap V} M^0 \cong W \times_M M^0$ as extending the codomain does not alter the fiber product, and $W \times_M M^0$ and $M^0 \times_M W$ are diffeomorphic over $M$.
	We now co-orient $M^0 \times_M W$ as a fiber product as in \cref{L: W0 cochain}, giving $W^0 = M^0 \times_M W \in PC^*_\Gamma(M)$.

	Finally, as $M^0$ has no boundary, we compute using \cref{leibniz} that
	$$\bd (W^0) = - M^0 \times_{M} \bd W.$$
	But $W$ represents a cocycle in $U \cap V$, so $\bd W$ is a union of trivial and degenerate precochains in $U \cap V$.
	The definitions of trivial and degenerate do not rely on the properness of the map, so $\bd W$ will also be a union of trivial and degenerate manifolds over $M$.
	Consequently, the same is true of $M^0 \times_{M} \bd W$ by \cref{L: pullback with Q}, noting the proof of that result does not rely on any map being proper.
\end{proof}

The next lemma shows that our construction of $\delta$ does not depend on the use of $0$ as our regular value.
The proof uses $\phi^{-1}([p,q])$, so the reader should compare with \cref{L: Wpq cochain}.

\begin{lemma}\label{L: different point}
	Suppose that $U$ and $V$ are open subsets of a manifold without boundary $M = U \cup V$.
	Suppose $W \in PC^*_\Gamma(U \cap V)$ represents an element of $H^*_\Gamma(U \cap V)$ and that $\phi \colon M = U \cup V \to [-1/2,1/2]$ separates $U$ and $V$.
	Let $p, q\in (-1/2,1/2)$ be any regular values for $\phi$ and $\phi r_W$.
	Then $W^p$ and $W^q$, as defined just above, represent the same element of $H^*_\Gamma(M)$.
\end{lemma}

\begin{proof}
	Suppose without loss of generality that $p<q$.
	Then the submanifold $\phi^{-1}([p,q]) \subset U \cap V \subset M$ will be transverse to $r_W \colon W \to M$.
	As $\phi^{-1}([p,q])$ is closed in $M$, its inclusion is proper, and the inclusion can be tautologically co-oriented as a codimension $0$ embedding (see \cref{D: tautological co-orientation}).
	Thus, as in \cref{L: Wpq cochain}, $$\bd(\varphi^{-1}([p,q])\times_{M} W ) = (M^p \times_{M} W)\sqcup -(M^q \times_{M} W) \sqcup (\phi^{-1}([p,q])\times_{M} \bd W).$$
	As $W$ is a cocycle, the last term is in $Q^*(M)$ (with support in $U \cap V$) by \cref{L: pullback with Q}, so $M^p \times_{M} W$ and $M^q \times_{M} W$ represent cohomologous cocycles.
\end{proof}

Next, we show independence of the representative $W$, assuming a fixed $\phi$.

\begin{lemma}\label{L: different W}
	Suppose that $U$ and $V$ are open subsets of a manifold without boundary $M = U \cup V$.
	Suppose $W_1, W_2 \in PC^*_\Gamma(U \cap V)$ represent the same element of $H^*_\Gamma(U \cap V)$.
	Let $\phi \colon U \cup V \to [-1/2,1/2]$ separate $U$ and $V$ and suppose $q \in (-1/2,1/2)$ is a regular value for $\phi$, $\phi r_{W_1}$, and $\phi r_{W_2}$.
	Then $W_1^q$ and $W_2^q$ represent the same element of $H^*_\Gamma(M)$.
\end{lemma}

\begin{proof}
	By assumption, there is a $Z \in PC^*_{\Gamma}(U \cap V)$ with $\bd Z \sqcup -W_1 \sqcup W_2 \in Q^*(U \cap V)$.
	For our fixed $\phi$, there is by Sard's Theorem a $p \in (-1/2,1/2)$ such that $p \into (-1/2,1/2)$ is transverse to $\phi$, $\phi r_Z$, $\phi r_{W_1}$, and $\phi r_{W_2}$.
	By the preceding lemma, it suffices to show that $W_1^p$ and $W_2^p$ represent the same element of $H^*_\Gamma(M)$.

	As in the proofs of the preceding lemmas, the fiber product with $M^p$ over $M$ takes an element of $PC^*_\Gamma(U \cap V)$ to $PC^*_\Gamma(M)$ and an element of $Q^*(U \cap V)$ to $Q^*(M)$.
	So
	$$M^p \times_{M}(\bd Z \sqcup -W_1 \sqcup W_2) = (\bd Z)^p \sqcup -W_1^p \sqcup W_2^p$$
	is an element of $Q^*(M)$.
	Furthermore, by an obvious modification to \cref{C: co-orient W0}, we have $\bd (Z^p) = -(\bd Z)^p$.
	So $$-\bd (Z^p) \sqcup -W_1^p \sqcup W_2^p \in Q^*(M),$$
	which implies that $W_1^p$ and $W_2^p$ represent the same cohomology class in $M$.
\end{proof}

\begin{proposition}\label{P: connecting}\index{connecting morphism for geometric cohomology!is well defined}
	Suppose that $U$ and $V$ are open subsets of a manifold without boundary $M = U \cup V$.
	The map $H^k_\Gamma(U \cap V) \xr{\delta} H^{k+1}_\Gamma(U \cup V) = H^{k+1}_\Gamma(M)$ is well defined.
	In particular, it does not depend on the choice of separating function.
\end{proposition}

\begin{proof}
	Suppose $W_1$ and $W_2$ represent the same element of $H^*_\Gamma(U \cap V)$ and that we choose functions $\phi_1,\phi_2$ that are respectively separating for $W_1$ and $W_2$ over $U$ and $V$.
	We must show that $\phi_1^{-1}(0)\times_{M} W_1$ and $\phi_2^{-1}(0)\times_{M} W_2$ represent the same element of $H^*_\Gamma(M)$.

	By Sard's Theorem, we can find points $p_1<p_2$ in $(-1/2,1/2)$ such that $p_1$ is a regular value for $\phi$ and $\phi_1 r_{W_1}$ and $p_2$ is a regular value for all of $\phi$, $\phi_1 r_{W_1}$, $\phi_2 r_{W_1}$, and $\phi_2 r_{W_2}$.
	By \cref{L: different point}, $\phi_1^{-1}(0) \times_{M} W_1$ and $\phi_2^{-1}(0) \times_{M} W_2$ represent the same elements of $H^*_\Gamma(M)$ respectively as $\phi_1^{-1}(p_1) \times_{M} W_1$ and $\phi_2^{-1}(p_2) \times_{M} W_2$.
	Further, by \cref{L: different W}, $\phi_2^{-1}(p_2) \times_{M} W_2$ represents the same cohomology class as $\phi_2^{-1}(p_2) \times_{M} W_1$.
	So it suffices to show that $\phi_1^{-1}(p_1) \times_{M} W_1$ represents the same cohomology class as $\phi_2^{-1}(p_2) \times_{M} W_1$.

 	So for simplicity of notation, let us suppose a single precochain $W$ and two separating functions $\phi_1,\phi_2$ such that $p_1 < p_2$ are respective regular values for $\phi$, $\phi_1 r_W$, and $\phi_2 r_W$.

	First, suppose there exists $q$ with $-1/2<q<p_2$ such that $\phi_1^{-1}(p_1) \subset \phi_2^{-1}([-1/2,q])$.
	By postcomposing $\phi_1$ with a diffeomorphism of $[-1/2,1/2]$, we may suppose that $p_1<q$.
	Let us choose $u_1,u_2$ such that $-1/2<u_1<p_1<u_1<q$.
	Using the Urysohn lemma we can find a continuous $\hat \Phi$ on $M$ such that:

	\begin{enumerate}
		\item $\hat \Phi$ is equal to $\phi_1$ on $\phi_1^{-1}[u_1,u_2]$,
		\item $\hat\Phi$ is equal to $\phi_2$ on $\phi_2^{-1}([q,1/2])$,

		\item $\hat\Phi$ takes $U - V$ to $-1/2$,

		\item $\hat\Phi^{-1}(p_1) = \phi_1^{-1}(p_1)$,

		\item $\hat\Phi^{-1}(p_2) = \phi_2^{-1}(p_2)$.
	\end{enumerate}
	Furthermore, by a sufficiently small homotopy $\hat\Phi$ can be approximated by a smooth function $\Phi$ that is still separating, preserves the last two properties, and, for some $\epsilon >0$, agrees with $\phi_1$ on $\phi_1^{-1}((p_1-\epsilon,p_2+\epsilon))$ and with $\phi_2$ on $\phi_2^{-1}((p_2-\epsilon,p_2+\epsilon))$.
	Then $p_1$ and $p_2$ are regular values for $\Phi$ and $\Phi r_W$, and
	we have $\Phi^{-1}(p_1) \times_{M} W \cong \phi_1^{-1}(p_1) \times_{M} W$ and similar for $p_2$ and $\phi_2$.
	So now $\Phi^{-1}(p_1) \times_{M} W$ and $\Phi^{-1}(p_2) \times_{M} W$ represent the same element in $H^*_\Gamma(M)$ by \cref{L: different point}.

	While this argument is written for $\phi_1^{-1}(p_1)$ ``below some $q<p_2$,'' clearly an analogous argument holds for $\phi_1^{-1}(p_1)$ ``above some $q>p_2$,'' or with the roles of $\phi_1$ and $\phi_2$ reversed.

	Lastly, we have to consider the case where there is not a $q$ as in the preceding paragraph that allows us to separate $\phi_1^{-1}(p_1)$ from $\phi_2^{-1}(p_2)$.
	In this case we will modify $\phi_2$ in an appropriate way.
	Let use choose $r$ with $p_2 <r<1/2$.
	Define $\hat \Phi$ as follows:

	\begin{enumerate}
		\item On $\phi_2^{-1}([-1/2,r])$, take $\hat \Phi = \phi_2$,

		\item On $\phi_2^{-1}([r,1/2])$ use the Urysohn Lemma to construct a continuous function $\phi_2^{-1}([r,1/2]) \to [r,1/2]$ that takes $\phi_2^{-1}(r) \cup [\phi_2^{-1}([r,1/2]) \cap \phi_1^{-1}(p_1)]$ to $r$ and $V - U$ to $1/2$; let $\hat \Phi$ be this constructed function on $\phi_2^{-1}([r,1/2]$.
	\end{enumerate}

	Now, we can modify $\hat \Phi$ by an $\epsilon$-homotopy fixing $\phi_2^{-1}([-1/2,r])$ and $\Phi^{-1}(1/2)$ to a smooth $\Phi$ with $\Phi^{-1}(p_2) = \phi_2^{-1}(p_2)$ and such that there is a $q$, $r<q<1$, with $\phi_1^{-1}(p_1) \subset \Phi^{-1}([-1/2,q])$.
	We can now choose a $p_2'$ with $q<p_2'<1$ that separates $W$ via $\Phi$.
	From our previous cases, we know $\Phi^{-1}(p_2')\times_{M}W$ is cohomologous to $\Phi^{-1}(p_2)\times_{M}W = \phi_2^{-1}(p_2)\times_{M}W$, but $\Phi^{-1}(p_2')$ and $\phi_1^{-1}(p_1)$ are now related as in the preceding case.

	Altogether we have shown that $\delta$ is independent of choices.
\end{proof}

\begin{theorem}\label{T: natural connection}\index{connecting morphism for geometric cohomology!naturality}
	The connecting map $\delta$ is natural, i.e.\ for a continuous map $f \colon U' \cup V' \to U \cup V$ such that $f(U') \subset U$ and $f(V') \subset V$, the following diagram commutes for all $k$:
	\[
	\begin{tikzcd}
		H^k_\Gamma(U \cap V) \arrow[r, "\delta"] \arrow[d, "f^*"] & H^{k+1}_\Gamma(U \cup V) \arrow[d, "f^*"] \\
		H^k_\Gamma(U' \cap V') \arrow[r, "\delta"] & H^{k+1}_\Gamma(U' \cup V').
	\end{tikzcd}
	\]
\end{theorem}

We start with a lemma that will be needed in the proof of \cref{T: natural connection} and is of somewhat independent interest.

\begin{lemma}\label{L: cut pullback}
	Suppose $M$ and $N$ are manifolds without corners and that $\phi \colon N \to \R$ is smooth with $0$ a regular value.
	Suppose $W \in PC^*_\Gamma(N)$ and that $r_W \colon W \to N$ is transverse to $N^0$.
	Further, suppose $g \colon M \to N$ is smooth and transverse to both $W \to N$ and $W^0 \to N$.
	Then the pullback $W^0 \times_N M \to M$ is co-oriented isomorphic to $(W \times_N M)^0 \to M$, where this splitting of $W \times_N M$ is with respect to the composite $\phi g$.
\end{lemma}

\begin{proof}
	We first note that $W^0 \times_N M$ is well defined as we assume $g$ is transverse to $W^0 \to N$.

	Now, let us write the pullback of $W$ by $g$, which exists by the transversality of $g$ and $r_W$, as $r_W^* \colon W \times_N M \to M$.
	To see that $(W \times_N M)^0$ is well defined, it suffices to show that the transversality of $g$ with $r_{W^0}$ implies that $0$ is a regular value for the composition $\phi g r_W^*$.
	So we must show that the restriction of $\phi gr_W^*$ to each stratum of $W \times_N M$ is transverse to $0$.
	Recall by \cref{pullback} that $S^k(W \times_N M) = S^k(W) \times_N M$, as $M$ is without boundary.
	Similarly, we note for a bit later that $$S^k(W^0) = S^k(N^0 \times_N W) = N^0 \times_N S^k(W) = (S^k(W))^0.$$
	So let $(x,y) \in S^k(W) \times_N M$ with $x \in S^k(W)$ and $y \in M$ such that $\phi g(y)=0$, and let $z = r_W(x) = g(y) \in N$.
	As $D_x(r_W|_{S^k(W)}) = (D_xr_W)|_{T_xS^k(W)}$, we will simply write $D_xr_W$ rather than $D_x(r_W|_{S^k(W)})$ in what follows.
	We recall that the tangent bundle of the pullback is the pullback of the tangent bundles by \cref{L: tangent of pullbacks}, so the tangent space of $S^k(W) \times_N M$ is $T_x S^k(W) \times_{T_{z}N} T_yM$.
	As $r_W$ is transverse to $N^0$, we know $0$ is a regular value for $\phi r_W$; see \cref{S: splitting}.
	So, in particular, zero is a regular value for $\phi \circ r_W|_{S^k(W)}$, so there must be a vector $\xi \in T_xS^k(W)$ with $D_x(\phi r_W)(\xi)\neq 0$.
	Thus $D_xr_W(\xi) \neq 0$.
	As $g$ is assumed transverse to $W^0 \to N$, there must be $\alpha \in T_xS^k(W^0)$ and $\beta \in T_yM$ such that $D_xr_W(\xi) = D_xr_W(\alpha) + D_yg (\beta)$.
	Rewriting, $D_yg (\beta) = D_xr_W(\xi-\alpha)$.
	As $\alpha \in T_xW^0 \in \ker(D_x (\phi r_W))$, we have $$D_{z}\phi(D_yg (\beta)) = D_{z}\phi \circ D_xr_W(\xi-\alpha) = D_x(\phi r_W)(\xi)-D_x(\phi r_W)(\alpha) = D_x(\phi r_W)(\xi)\neq 0.$$
	So, recalling that $Dr_W^*$ is simply the projection to $T_yM$, the pair $( \xi-\alpha, \beta)$ is a non-zero vector in $$T_{(x,y)}(S^k(W) \times_N M) = T_x S^k(W) \times_{T_{z}N} T_yM$$ that maps by $D(\phi gr_W^*)$ to $D_{z}\phi(D_yg (\beta)) \neq 0 \in T_{0}[-1/2, 1/2]$.
	As $(x,y)$ was an arbitrary point of $(\phi gr^*_{W})^{-1}(0) \subset W \times_N M$, this shows that $0$ is a regular value for $\phi gr^*_W$.
	So $(W \times_N M)^0$ is well defined.

	Now, from the definitions and \cref{S: splitting}, $W^0 \times_N M = (N^0 \times_N W) \times_N M$, treating $N^0 \times_N W$ as a fiber product and $W^0 \times_N M$ as a pullback over $M$.
	On the other hand, $(W \times_N M)^0 = M^0 \times_M (W \times_N M)$, where $W \times_N M$ is the pullback to $M$ and then we take the fiber product with $M^0$.
	In both cases, these translate to the pairs of points $(x,y) \in W \times M$ with $r_W(x) = g(y)$ and $\phi(r_W(x)) = \phi(g(y)) = 0$.
	In other words, as spaces these are both precisely the limit of the following diagram, together with its map to $M$:
	\[
	\begin{tikzcd}
		& M \arrow[d,"g"] & \\
		W \arrow[r, "r_W"] & N \arrow[r, "\phi"] & \R & \arrow[l,hook'] 0.
	\end{tikzcd}
	\]
	Thus $(W \times_N M)^0$ and $W^0 \times_N M$ are diffeomorphic over $M$, and it remains to check the co-orientations.

	We return to the definitions of the pullback and fiber product co-orientations.
	It suffices to compare $(W \times_N M)^0$ and $W^0 \times_N M$ at an arbitrary point of the top stratum.
	We first consider the co-orientation of $(W \times_N M)^0$.
	Let $(x,y) \in (W \times_N M)^0$ with $r_W(x) = g(y) = z$.
	Choose $e \colon W \into N \times \R^K$, fix $\beta_N$ at $z$, and choose $\beta_W$ at $x$ so that $\omega_{r_W} = (\beta_W,\beta_N)$.
	Let $\nu$ be the Quillen-oriented normal bundle of $W$ in $N \times \R^K$ so that $\beta_W \wedge \beta_\nu = \beta_N \wedge \beta_E$, where $\beta_E$ is the standard orientation of $\R^K$.
	Then, by definition, the pullback map from $P = W \times_N M \subset M \times \R^K$ to $M$ is co-oriented by $(\beta_P,\beta_M)$ if we choose $\beta_P$ and $\beta_M$ such that $\beta_P \wedge \beta_\nu = \beta_M \times \beta_E$, recalling that we let $\nu$ also denote the pulled back normal bundle of $P$ in $M \times \R^K$.
	We suppose we have chosen such $\beta_P$ and $\beta_M$.
	Next, let $M^0 \subset M$ have normal co-orientation $\beta_\phi$ determined by pulling back the standard orientation from $\R$, and similarly let $\beta_\phi$ denote the pullback normal co-orientation to $(W \times_N M)^0$.
	Then by \cref{P: codim 1 co-orient}, the co-orientation of the inclusion $(W \times_N M)^0 = M^0 \times_M (W \times_N M) \hookrightarrow W \times_N M$ is $(\beta_Q,\beta_Q \wedge \beta_\phi)$ for any $\beta_Q$.
	The fiber product $(W \times_N M)^0 \into M$ is then co-oriented by the composition $(\beta_Q,\beta_Q \wedge \beta_\phi)*(\beta_P,\beta_M)$.
	So if we choose $\beta_Q$ so that $\beta_Q \wedge \beta_\phi = \beta_P$ (or equivalently $\beta_Q$ and $\beta_M$ so that $\beta_Q \wedge \beta_\phi \wedge \beta_\nu = \beta_M \wedge \beta_E$), the co-orientation is $(\beta_Q,\beta_M)$.

	On the other hand, consider $(N^0 \times_N W) \times_N M = W^0 \times_N M \to M$.
	Once again, at the same points, we fix $\beta_N$ and $\beta_W$ so that $\omega_{r_W} = (\beta_W,\beta_N)$.
	Again by \cref{P: codim 1 co-orient}, the co-orientation of the pullback $N^0 \times_N W \to W$ is $(\beta_{W^0},\beta_{W^0} \wedge \beta_\phi)$ for any $\beta_{W^0}$, continuing to let $\beta_\phi$ denote any normal co-orientation pulled back via $\phi$.
	If we choose $\beta_{W^0}$ so that $\beta_{W^0} \wedge \beta_\phi = \beta_W$ then we have $r_{W^0} \colon W^0 \to N$ co-oriented by $(\beta_{W^0},\beta_N)$.
	As $W^0 \subset W$, we can embed $W^0$ in $N \times \R^K$ via the composition $W^0 \into W \xhookrightarrow{e}N \times \R^K$, using the same $e$ and $K$ as above.
	We also assume the same oriented normal bundle $\nu$.
	As $\beta_W \wedge \beta_\nu = \beta_N \wedge \beta_E$ and $\beta_W = \beta_{W^0} \wedge \beta_\phi$, we have $\beta_{W^0} \wedge \beta_\phi \wedge \beta_\nu = \beta_N \wedge \beta_E$ so that $\beta_\phi \wedge \beta_\nu$ is the Quillen orientation for the normal bundle of $W^0$ in $N \times \R^K$.
	Using this to pull back $W^0 \to N$ to $W^0 \times_N M \to M$, by definition the pullback co-orientation is $(\beta_Q,\beta_M)$ when $\beta_Q$ and $\beta_M$ are chosen so that $\beta_Q \wedge \beta_\phi \wedge \beta_\nu = \beta_M \wedge \beta_E$.
	But this is exactly the same co-orientation we arrived at in the preceding paragraph.
	Thus the co-orientations of the two constructions agree.
\end{proof}

\begin{proof}[Proof of \cref{T: natural connection}]
	We first establish some notation.
	Let $M = U' \cup V'$ and $N = U \cup V$.
	Let $A = U' \cap V'$ and $B = U \cap V$.
	Let $r_W \colon W \to B \subset N$ represent an element $\uW \in H^*_\Gamma(B)$.
	Let $\psi \colon M \to [-1/2, 1/2]$ be separating for $U'$ and $V'$, and let $\phi \colon N \to [-1/2, 1/2]$ be separating for $U$ and $V$.
	We further suppose $0$ is a regular value for $\psi$ and that $\phi$ separates $W$ over $U$ and $V$, treating $W$ as having image in $N$.
	Let $N^0$ and $W^0$ be determined by $\phi$ as in \cref{S: splitting}, and similarly, let $M^0$ be determined by $\psi$.
	We let $r_{W^0} \colon W^0 \to N$ be the reference map for $W^0$.
	By postcomposing $\psi$ with an orientation preserving-diffeomorphism of $[-1/2, 1/2]$, we may assume that $\pm 1/4$ are regular values for $\psi$, and we let $K = \psi^{-1}([-1/4, 1/4]) \subset M$.

	Now, noting that $f(A) \subset B$, apply\footnote{To obtain transversality to both $W$ and $W^0$, it suffices to choose $g_1$ transverse to the disjoint union of $r_{W^+} \colon W^+ \to M$ and $r_{W^-} \colon W^- \to M$. For both $g_1$ and $g_1|_{M^0}$ to have the desired transversality, we note that the restriction of the function $H$ constructed in the proof of \cref{T: basic trans} to $M^0 \times D$ has the desired transversality at almost every $s \in D$ by the same argument as for $H$ on all of $M \times D$. As the intersection of two almost everywhere sets also has measure zero complement, we can choose an $s \in D$ that gives transversality for both $M$ and $M^0$. Compare the proof of the Transversality Homotopy Theorem in \cite{GuPo74}, which allows that our analogue of $M$ have boundary.} \cref{T: basic trans} to obtain $g_1 \colon A \to B$, a smooth approximation to $f|_A$ that is transverse to $r_W$ and $r_{W^0}$ and so that the restriction to $M^0$ is also transverse to $W$.
	Then $f^*(\uW) \in H^*_\Gamma(A)$ is represented by the pullback of $W$ by $g_1$ by \cref{D: cohomology pullback and homology transfer}.
	For clarity, we write this pullback as $g_1^*(W) = W \times_B A$.
	By \cref{L: transverse to pullback}, as $g_1|_{M^0}$ is transverse to $W$, the inclusion of $M^0$ into $M$ is transverse to $g_1^*(W)$.
	The fiber product $ - M^0 \times_M g_1^*(W)$ represents $\delta f^* (W)$ by \cref{D: connecting}.

	Next, to consider $f^* \delta (\uW)$, we wish to construct a homotopy from $f$ to a map $g_2 \colon M \to N$ such that
	\begin{enumerate}
		\item $g_2 = g_1$ on $K$,
		\item $g_2$ is smooth,
		\item $g_2$ is transverse to $W$ and $W^0$
		\item $g_2$ takes $\psi^{-1}([-1/2,1/4])$ to $U$ and $\psi^{-1}([-1/4,1/2])$ to $V$.
	\end{enumerate}
	We will show that we can construct such a $g_2$ in a lemma following the remainder of the proof of \cref{T: natural connection}.

	Assuming $g_2$ exists, by definition $f^* \delta (\uW)$ is represented by the pullback of $-W^0 \to N$ by $g_2$, which we will write $g_2^*(-W^0) = - g_2^*(W^0) \to M$.
	Now we apply \cref{L: cut pullback} to see that $g_2^*(W^0)$ is isomorphic to $(g_2^*(W))^0$, where the latter expression is the splitting with respect to the composite $\phi g_2$.
	In other words, $g_2^*(W^0) = M^0_{\phi g_2} \times_M g_2^*(W)$, writing $M^0_{\phi g_2}$ for the splitting of $M$ determined by the function $\phi g_2$.

	So to complete the proof we must show that this $M^0_{\phi g_2} \times_M g_2^*(W)$ represents the same cohomology class in $M$ as $M^0_\psi \times_M g_1^*(W)$, where $M^0_\psi$ is the splitting of $M$ using $\psi$.

	We next observe that as $M^0_\psi$ is in the interior of $K$, where we know $g_1 = g_2$, we must have $M^0_\psi \times_M g_1^*(W) = M^0_\psi \times_M g_2^*(W)$.
	So we compare $M^0_{\phi g_2} \times_M g_2^*(W)$ with $M^0_\psi \times_M g_2^*(W)$.
	The rough idea of the remainder of the proof is to show that $\psi$ and $\phi g_2$ are both separating functions for $g_2^*(W)$ over $g_2^{-1}(U)$ and $g_2^{-1}(V)$, which will allow us to invoke \cref{P: connecting}, though there remain some technicalities, including replacing $\psi$ with a slight modification.

	So, consider $g_2^{-1}(U)$ and $g_2^{-1}(V)$. As $N = U \cup V$, we must have $M = g_2^{-1}(U) \cup g_2^{-1}(V)$.
	We have $g_2^{-1}(U) - g_2^{-1}(V) = g_2^{-1}(U  - V)$, so $\phi g_2(g_2^{-1}(U) - g_2^{-1}(V)) = \phi g_2(g_2^{-1}(U  - V)) \subset \phi(U  - V) = -1/2$, and similarly $\phi g_2(g_2^{-1}(V)  - g_2^{-1}(U)) = 1/2$.
	So $\phi g_2$ is a separating function on $M$ for $g_2^{-1}(U)$ and $g_2^{-1}(V)$.
	We also know $M^0_{\phi g_2}$ is transverse to $g_2^*(W)$, so $\phi g_2$ separates $g_2^*(W)$ over $g_2^{-1}(U)$ and $g_2^{-1}(V)$.

	As we chose $g_2$ so that $g_2 = g_1$ on $K$ and so that it takes $\psi^{-1}([-1/2,1/4])$ to $U$ and $\psi^{-1}([-1/4,1/2])$ to $V$, we have $g_2^{-1}(U  - V) \subset \psi^{-1}([-1/2,-1/4])$ and $g_2^{-1}(V  - U) \subset \psi^{-1}([1/4,1/2])$.
	Let $\psi_1$ be the composition of $\psi$ with a smooth non-decreasing map $a \colon [-1/2,1/2] \to [-1/2,1/2]$ such that:
	\begin{enumerate}
		\item $a([-1/2,-1/4]) = -1/2$,
		\item $a(x)=x$ on a neighborhood of $0$,
		\item $a([1/4,1/2]) = 1/2$.
	\end{enumerate}
	In this case, $M^0_\psi = M^0_{\psi_1}$, so we have $M^0_\psi \times_M g_2^*(W) = M^0_{\psi_1} \times_M g_2^*(W)$.
	Furthermore, as $g_2^{-1}(U) - g_2^{-1}(V) = g_2^{-1}(U - V)$, and similarly reversing the roles of $U$ and $V$, we see $\psi_1$ separates $g_2^*(W)$ over $g_2^{-1}(U)$ and $g_2^{-1}(V)$.

	Lastly, we observe that $g_2^{-1}(U) \cup g_2^{-1}(V) = M$ and $g_2^{-1}(U) \cap g_2^{-1}(V) = g_2^{-1} (U \cap V) = g_2^{-1}(B)$, so $g_2$ maps this intersection into $B = U \cap V$.
	Therefore, the pullback $(g_2|_{g_2^{-1} (B)})^* (W)$ is an element of $PC^*_\Gamma(g_2^{-1} (B))$.
	But, $M^0_{\phi g_2} \times_M	(g_2|_{g_2^{-1} (B)})^* (W) = M^0_{\phi g_2} \times_M g_2^*(W)$ and $M^0_{\psi_1} \times_M (g_2|_{g_2^{-1} (B)})^* (W) = M^0_{\psi_1} \times_M g_2^*(W)$, as $M^0_{\phi g_2}$ and $M^0_{\psi_1}$ are both in $g_2^{-1} (B)$.
	Now, as desired, $M^0_{\phi g_2} \times_M	(g_2|_{g_2^{-1} (B)})^* (W)$ and $M^0_{\psi_1} \times_M (g_2|_{g_2^{-1} (B)})^* (W)$ are two splittings of the same element of $PC^*_\Gamma(g_2^{-1} (B)) = PC^*_\Gamma(g_2^{-1}(U) \cap g_2^{-1}(V)) $ with respect to different separating functions over $g_2^{-1}(U)$ and $g_2^{-1}(V)$, and so they represent the same element of $H^*_\Gamma(M)$ by \cref{P: connecting}.
\end{proof}

\begin{lemma}
	There exists a $g_2$ as claimed in the proof of \cref{T: natural connection}.
\end{lemma}

\begin{proof}
	As $\psi^{-1}(\pm 1/4)$ are bicollared submanifolds of $M$, the submanifold $K \subset A$ is collared, so the inclusion map $K \into A$ is a cofibration, e.g. via\footnote{We define a map $M \to [0,1]$ as in the cited theorem by taking $K$ to $0$, projecting the parts of the collars isomorphic to $\psi^{-1}(\pm 1/4) \times [0,1]$ to the second factor, and taking everything else to $1$.} \cite[Theorem VII.1.5]{Bred97}.
	Therefore, the restriction to $K \times I$ of the homotopy from $f|_A$ to $g_1$ used at the start of the proof of \cref{T: natural connection} can be extended to a continuous homotopy $M \times I \to N$ from $f$ to a continuous map we will call $f' \colon M \to N$.
	So by construction, $f' = g_1$ on $K$.
	Furthermore, choosing the collars to be contained in $A$, which is an open neighborhood of $K$ in $M$, then by the usual construction of a homotopy extension using a neighborhood deformation retraction (see \cite[Section VII.1]{Bred97}), we can construct the homotopy such that we have a neighborhood $L$ of $K$ such that $f'$ takes $L$ to $B$ (because our homotopy of $f|_A$ maps to $B$) and so that $f' = f$ on $M - L$.
	In particular, $f'$ takes $\psi^{-1}([-1/2,1/4])$ to $U$ and $\psi^{-1}([-1/4,1/2])$ to $V$.

	From here, we will modify the proof of \cref{T: basic trans}.
	As in that argument, we suppose $N$ properly embedded as a closed manifold in some Euclidean space; then $N$ has a distance function inherited from Euclidean space, which we will denote $d$.

	First, on $\psi^{-1}([-1/2,-1/4])$, we consider the function $x \mapsto d(f'(x), V - U)$.
	As distance to a closed set is a continuous function and $f'$ takes $\psi^{-1}([-1/2,1/4])$ to $U$, this is a continuous function $\psi^{-1}([-1/2,-1/4]) \to (0,\infty)$.
	Similarly, $x \mapsto d(f'(x), U - V)$ is a positive function on $\psi^{-1}([1/4, 1/2])$.
	As $(0,\infty) \cong \R$, by the Tietze Extension Theorem we can extend these functions to a single function $\delta \colon M \to (0, \infty)$.
	We now apply smooth approximation\footnotemark to replace $f'$ with a homotopic function $g'$ that is smooth, such that $g' = f' = g_1$ on $K$, and such that $d(f'(x),g'(x))<\delta(x)/2$ for all $x \in M$.
	\footnotetext{A good proof that $f'$ is homotopic rel $K$ to a smooth map can be found in \cite[Theorem 6.26]{Lee13}. To obtain the distance bound, we modify the argument there to require that the function called $\td \delta$ in that proof satisfies $\td \delta(x) < \delta(x)/4$ in addition to the other requirements on $\td \delta$ in the proof. Noting that the tubular neighborhood in \cite{Lee13} is the same as an $\epsilon$-neighborhood in \cite{GuPo74}, we can now use a triangle inequality argument as in our proof of \cref{T: basic trans} to see that the smooth approximation satisfies $d(f'(x),g'(x))<\delta(x)/2$.}

	Next, we modify the transversality portion of the proof of \cref{T: basic trans} to obtain $g_2$ from $g'$.
	For the modification, we replace the smooth function $\eta$ in the definition of $H(x,s)$ with one that satisfies all the properties stated there except
	\begin{enumerate}
		\item we take $\eta(x)=0$ on $K$, and
		\item on $M - K$, we take $\eta(x) > 0$ and, in addition to its other size constraints, we take $\eta$ sufficiently small that $d(H(x,s),g'(x)) < \delta(x)/2$.
		This can be done using the same argument as in the proof of \cref{T: basic trans} for the case of proper maps, using $\delta/2$ in the role of $\varepsilon$ (that part of the argument did not rely on any map being proper).
	\end{enumerate}
	Then, continuing the transversality argument of \cref{T: basic trans}, $H(-,s)|_{M-K}$ is transverse to $W$ and $W^0$ for almost every $s \in D$, so we again choose such an $s_0$ and let $$h(x,t) = H(x, ts_0) = \pi(g'(x) + \eta(x)ts_0),$$ for all $x \in M$ and $t \in [0,1]$.
	Let $g_2 = h(-,1)$.

	Altogether now, $g_2$ is smooth, and it is equal to $g_1$ on $K$ by construction.
	It is transverse to $W$ and $W^0$ on $K$ because $g_1$ is, and it is transverse to $W$ and $W^0$ on $M - K$ by the argument from \cref{T: basic trans} with our choice of $s_0$.
	Finally, we see $g_2$ takes $\psi^{-1}([-1/2,1/4])$ to $U$ and $\psi^{-1}([-1/4,1/2])$ to $V$: On $K$, we have $g_2(K) = g_1(K) \subset B = U \cap V$ from the construction of $g_1$.
	For $x\in \psi^{-1}([-1/2,-1/4])$, our definition of the bounding function $\delta$ above and the choice of $\delta/2$-small homotopies in the steps from $f'$ to $g'$ and $g'$ to $g_2$ ensure that $d(f'(x), g_2(x)) < \delta(x)$, so $g_2(x)$ is a positive distance from $V - U$.
	Similarly for $\psi^{-1}([1/4,1/2])$ and $U - V$.
\end{proof}
\index{Mayer-Vietoris sequence!for geometric cohomology!connecting morphism|)}\index{connecting morphism for geometric cohomology|)}

\medskip

\noindent\textbf{The cohomology Mayer-Vietoris sequence.} We now come to the Mayer-Vietoris sequence itself.
The following observations and notation will be useful when working with the restriction maps.

Suppose $\uW \in H^*_\Gamma(M)$ represented by $r_W \colon W \to M$ and that $f \colon U \into M$ is the inclusion of an open subset.
As $f$ is smooth and necessarily transverse to $r_W$, we have $f^*(\uW)$ represented by the pullback $W \times_M U \to U$ by \cref{D: cohomology pullback and homology transfer}.	By \cref{P: codim 0 pullback}, this pullback is simply the restriction of $r_W$ to $r_W^{-1}(U)$.
In what follows, we will write the pullback $W \times_M U$ as $W|_U$, and the corresponding cohomology class, represented by $W|_U$, as $\uW|_U$.

\begin{theorem}\label{T: absolute MV}\index{Mayer-Vietoris sequence!geometric cohomology}\index{geometric cohomology!Mayer-Vietoris sequence}
	Let $U,V \subset M$ be open subsets.
	There is a long exact Mayer--Vietoris sequence
	\[
	\cdots \to H^k_\Gamma(U \cup V) \xr{i} H^k_\Gamma(U) \oplus H^k_\Gamma(V) \xr{j} H^k_\Gamma(U \cap V) \xr{\delta} H^{k+1}_\Gamma(U \cup V)\to\cdots
	\]
	with $i(\uW) = (\uW|_U, -\uW|_V)$, $j(\uW_1,\uW_2) = \uW_1|_{U \cap V}+\uW_2|_{U \cap V}$, and $\delta$ as given in \cref{D: connecting}.
\end{theorem}

\begin{proof}
	We first make some general observations that will apply throughout the proof.
	Notationally, we always assume $\uW$ is represented by $W$, etc.
	In the following, we typically assume an appropriate separating function $\phi \colon U \cup V \to [-1/2,1/2]$.
	We further assume $\phi$ has been chosen to separate any precochain under discussion;
	see \cref{D: separating function}.

	We observe that $(U \cap V)^+ = \phi^{-1}([0,1/2])\cap(U \cap V)$ is a closed subspace both of $U \cap V$ and of $U$, as $\phi^{-1}([0,1/2]) \cap (U \cap V) = \phi|_U^{-1}([0,1/2])$ given that $\phi(U - V) = -1/2$.
	So the inclusions $(U \cap V)^+ \into U \cap V$ and $(U \cap V)^+ \into U$ are both proper maps.
	Thus if $W \in PC^*_\Gamma(U \cap V)$, in which case in particular $r_W \colon W \to U \cap V$ is proper, and assuming $\phi$ separates $W$ over $U$ and $V$, both the fiber products $(U \cap V)^+\times_{U \cap V} W \to U \cap V$ and $(U \cap V)^+\times_{U} W \to U$ will be proper maps.
	In fact, they both have the same domain, which is just $W^+$, and we will write $W^+$ for the corresponding elements of both $PC^*_\Gamma(U \cap V)$ or $PC^*_\Gamma(U)$, determining which is meant by context.
	Similarly, $W^-$ can represent an element of either $PC^*_\Gamma(U \cap V)$ or $PC^*_\Gamma(V)$.
	See \cref{F: MV1}.
	On the other hand, $(U \cap V)^- = \phi^{-1}([-1/2,0])\cap(U \cap V)$ is not generally closed in $U$, and so we do not obtain an element $W^-$ in $PC^*_\Gamma(U)$.
	Similarly, there is generally no $W^+$ in $PC^*_\Gamma(V)$.

	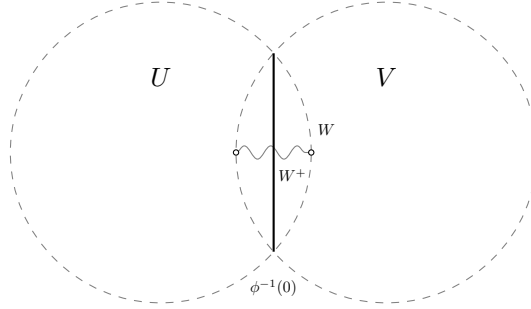
\begin{figure}[h]
			\begin{tikzpicture}
	\draw[gray, dashed] (0,0) circle (2);
	\node at (0,1){$U$};
	\draw[gray, dashed] (3,0) circle (2);
	\node at (3,1){$V$};
	\draw[gray, snake it] (1,0) -- (2,0);
	\draw[black, fill=white] (1,0) circle (1pt);
	\draw[black, fill=white] (2,0) circle (1pt);
	\draw[thick] (1.5, 1.32) -- (1.5, -1.32);
	\node[scale=.6] at (1.5, -1.8){$\phi^{-1}(0)$};
	\node[scale=.6] at (1.75, -.3){$W^+$};
	\node[scale=.6] at (2.2, .3){$W$};
	\end{tikzpicture}
		\caption{For a proper map $r_W \colon W \to U \cap V$, the restriction to $W^+$ is proper into $U$, while the restriction to $W^-$ is proper into $V$.}
		\label{F: MV1}
	\end{figure}

	Analogously, we have $\phi^{-1}([-1/2,0]) = \phi|_U^{-1}([-1/2,0])$, and so $\phi^{-1}([-1/2,0])$ is a closed subset of $U$ and of $U \cup V$.
	Therefore, the inclusions $\phi^{-1}([-1/2,0]) \into U$ and $\phi^{-1}([-1/2,0]) \into U \cup V$ are proper so that if $W \in PC^*_\Gamma(U)$, then $W^-$ can represent an element of either $PC^*_\Gamma(U)$ and $PC^*_\Gamma(U \cup V)$.
	See \cref{F: MV2}.
	Likewise, if $W \in PC^*_\Gamma(V)$ then we obtain $W^+$ in both $PC^*_\Gamma(V)$ and $PC^*_\Gamma(U \cup V)$.

	\begin{figure}[h]
			\begin{tikzpicture}
	\draw[gray, dashed] (0,0) circle (2);
	\node at (0,1){$U$};
	\draw[gray, dashed] (3,0) circle (2);
	\node at (3,1){$V$};
	\draw[gray, snake it] (-2,0) -- (2,0);
	\draw[black, fill=white] (-2,0) circle (1pt);
	\draw[black, fill=white] (2,0) circle (1pt);
	\draw[thick] (1.5, 1.32) -- (1.5, -1.32);
	\node[scale=.6] at (1.5, -1.8){$\phi^{-1}(0)$};
	\node[scale=.6] at (.7, -.3){$W^-$};
	\node[scale=1] at (-1, .3){$W$};
	\end{tikzpicture}
		\caption{For a proper map $r_W \colon W \to U$, the restriction to $W^-$ is proper into $U \cup V$.}
		\label{F: MV2}
	\end{figure}
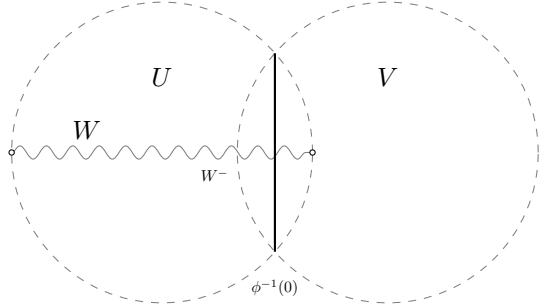

	We also note that $\phi^{-1}(0)$ is closed in $U \cap V$, $U$, $V$, and $U \cup V$, so if $W$ is an element of any of $PC_\Gamma^*(U \cap V)$, $PC_\Gamma^*(U)$, $PC_\Gamma^*(V)$, or $PC_\Gamma^*(U \cup V)$, then $W^0$ represents an element of all of these sets.

	Finally, let $A \subset B$ be any subsets of $M$, and suppose $W \xr{r_W} A$ and $r_W \xr{r_W} A \into B$ both represent precochains.
	Suppose further that $W \in Q^*(A)$; then also $W \in Q^*(B)$.
	This follows directly from the definitions, as being trivial or degenerate with image in $A$ implies being trivial or degenerate with image in $B$.

	These observations will all be used freely in the remainder of the proof, in which we proceed to consider exactness at each term of the sequence in turn.

	\textbf{Exactness at $\mathbf{H^k_\Gamma(U \cup V)}$.}
	Let $\uW \in H^{k-1}(U \cap V)$.
	Then $i\delta(\uW)$ is represented by $(-W^0|_U, W^0|_V)$.
	As $W^0$ already maps to $U \cap V$, we will write simply $(-W^0, W^0)$.
	By the above discussion, we have $W^+ \in PC^*_\Gamma(U)$, and by \cref{L: W0 cochain} as well as applying the discussion above to $(\bd W)^+$, we have $\bd (W^+) = W^0 \sqcup (\bd W)^+$ in $PC^*_\Gamma(U)$.
	Now $\bd W \in Q^*(U \cap V)$ as $W$ is a cocycle, so $(\bd W)^+ \in Q^*(U \cap V)$ by \cref{C: creasing Q} and so also $(\bd W)^+ \in Q^*(U)$ by the preceding discussion.
	Passing to cochains, we thus have $\bd (\underline{W^+}) = \underline{W^0}$ in $C^*_\Gamma(U)$, so $W^0$ represents $0 \in H^k_\Gamma(U)$.
	Similarly, $W^0$ represents $0 \in H^k_\Gamma(V)$ using $W^-$.
	So $i\delta = 0$.

	Next suppose $\uW \in H^k_\Gamma(U \cup V)$ with $i(\uW) = 0$.
	So $\uW|_U$ and $\uW|_V$ cobound in $U$ and $V$, respectively.
	By \cref{R: cycles and boundaries} there exists $A \in PC^*_\Gamma(U)$ with $\bd A \sqcup - W|_U \in Q^*(U)$ and $B \in PC^*_\Gamma(V)$ with $\bd B \sqcup - W|_V \in Q^*(V)$.
	We assume the separating function $\phi$ chosen so that it separates $W$, $A$, and $B$, considering each as mapping into $U \cup V$.
	Then $(\bd A \sqcup - W|_U)^- = (\bd A)^- \sqcup - (W|_U)^-$ is in $Q^*(U)$ by \cref{C: creasing Q}.
	By the above discussion, we can consider this equality holding in $PC^*_\Gamma(U \cup V)$, where also $(W|_U)^- = W^-$.
	So $(\bd A)^- \sqcup - W^- \in Q^*(U \cup V)$.
	Similarly $(\bd B)^+ \sqcup - W^+ \in Q^*(U \cup V)$.

	Now by \cref{L: W0 cochain}, we have $\bd(A^-) = -A^0 \sqcup (\bd A)^-$ in $PC^*_\Gamma(U)$ and $\bd (B^+) = B^0 \sqcup (\bd B)^+$ in $PC^*_\Gamma(U)$.
	Again, we may consider both of these equalities holding in $PC^*_\Gamma(U \cup V)$.
	So taking the disjoint union of $\bd(A^- \sqcup B^+)$ with $- W^- \sqcup - W^+$ and passing to $C^*_\Gamma(U \cup V)$, we have
	\begin{align*}
		\underline{	\bd(A^- \sqcup B^+) \sqcup - W^- \sqcup - W^+} &= \underline{-A^0 \sqcup (\bd A)^- \sqcup B^0 \sqcup (\bd B)^+ \sqcup - W^- \sqcup - W^+ }\\
		&= \underline{-A^0 \sqcup B^0 \sqcup (\bd A)^- \sqcup - W^- \sqcup (\bd B)^+ \sqcup - W^+ }\\
		&= \underline{-A^0 \sqcup B^0},
	\end{align*}
	as elements of $Q^*(U \cup V)$ represent $0$ in $C^*_\Gamma(U \cup V)$.
	So
		\[	\bd (\underline{A^- \sqcup B^+}) = \underline{W^-} + \underline{W^+} -(\underline{A^0}-\underline{B^0}).\]
	But $\underline{W^-} + \underline{W^+}$ is cohomologous to $\uW$ in $U \cup V$ by creasing, \cref{T: cohomology creasing}, so $W$ is cohomologous to $\underline{A^0}-\underline{B^0}$, which represents $\delta(\underline{B}|_{U \cap V}-\underline{A}|_{U \cap V})$ in $H^*_\Gamma(U \cup V)$.
	So $\uW \in H^*_\Gamma(U \cup V)$ is in the image of $\delta$.

	\textbf{Exactness at $\mathbf{H^k_\Gamma(U) \oplus H^k_\Gamma(V)}$.}
	It is immediate that the composition $ji$ is the $0$ map.

	Suppose $(\uW_1,\uW_2) \in \ker j \subset H^k_\Gamma(U) \oplus H^k_\Gamma(V)$.
	Using representatives $W_1,W_2$, this means that there is a $Z \in C^k_\Gamma(U \cap V)$ with $\bd Z \sqcup - W_1|_{U \cap V} \sqcup - W_2|_{U \cap V} \in Q^*(U \cap V)$.
	We choose $\phi$ to separate $Z$, $W_1$, and $W_2$ over $U$ and $V$.
	We claim that $\gamma = W_1^- \sqcup -Z^0 \sqcup -W_2^+$ represents an element of $H^k_\Gamma(U \cup V)$ whose image under $i$ is $(\uW_1,\uW_2)$.
	We compute using \cref{L: W0 cochain}:
	\begin{align*}
	\bd \gamma& = \bd (W_1^- \sqcup - Z^0 \sqcup -W_2^+)\\
		& = -W_1^0 \sqcup (\bd W_1)^- \sqcup (\bd Z)^0 \sqcup - W_2^0 \sqcup - (\bd W_2)^+\\
		& = (\bd W_1)^- \sqcup - (\bd W_2)^+ \sqcup (\bd Z \sqcup -W_1 \sqcup - W_2)^0.
	\end{align*}
	As $W_1$ and $W_2$ are cycles in $U$ and $V$ respectively, $(\bd W_1)^-$ and $(\bd W_2)^+$ are in $Q^*(U)$ and $Q^*(V)$ by \cref{R: cycles and boundaries,C: creasing Q}, and so they are also in $Q^*(U \cup V)$ by our discussion above.
	We also know $\bd Z \sqcup - W_1|_{U \cap V} \sqcup - W_2|_{U \cap V} \in Q^*(U \cap V)$, so its splitting at $0$ is also in $Q^*(U \cap V)$ by \cref{C: creasing Q}; hence $(\bd Z \sqcup -W_1 \sqcup - W_2)^0 \in Q^*(U \cup V)$.
	So $\bd \gamma$ represents $0$ in $C^*_\Gamma(U \cup V)$, and $\gamma$ represents an element of $H^k_\Gamma(U \cup V)$.

	Next we show that $\gamma|_U$ is cohomologous to $W_1$ in $U$.
	In fact, in $U$ we have
	\[\bd (Z^+) \sqcup \gamma|_U = Z^0 \sqcup (\bd Z)^+ \sqcup W_1^- \sqcup -Z^0 \sqcup -(W_2^+)|_U,\]
	using that $W_1^-$ and $Z^0$ are already precochains in both $U \cap V$ and $U$.
	So using \cref{L: W0 cochain}, passing to cochains, and undertaking some further manipulations, we have in $C^*_\Gamma(U)$ that
	\begin{align*}
	\bd (\underline{Z^+}) + \underline{\gamma}|_U &
		= \underline Z^0 + \underline{(\bd Z)^+} + \underline{W_1^-} - \underline Z^0 -\underline{W_2^+}|_U + (\underline{W_1^+} - \underline{W_1^+})\\
		&= \underline{(\bd Z)^+} + \underline{W_1^-} -\underline{W_2^+}|_U + \underline{W_1^+} - \underline{W_1^+}\\
		&= \underline{(\bd Z \sqcup - W_1|_{U \cap V} \sqcup - W_2|_{U \cap V})^+} + \underline{W_1^-} + \underline{W_1^+}\\
		&=\underline{W_1^-} + \underline{W_1^+}\\
		&=\underline{W_1}.
	\end{align*}
	Here we have used that $\bd Z \sqcup - W_1|_{U \cap V} \sqcup - W_2|_{U \cap V} \in Q^*(U \cap V)$, which implies that it is also in \cref{C: creasing Q},
	and that $\underline{W_1^-} + \underline{W_1^+}= \underline{W_1}$ in $H^*_\Gamma(U)$ by creasing, \cref{T: cohomology creasing}.
	So $\underline{W_1} = \underline{\gamma}|_U \in H^*_\Gamma(U)$, as claimed.

	Similarly, in $V$ we have
	\begin{align*}
	\bd (\underline{Z^-}) - \underline{\gamma}|_V &
		= - \underline Z^0 + \underline{(\bd Z)^-} - \underline{W_1^-}|_V + \underline Z^0 +\underline{W_2^+} + (\underline{W_2^-} - \underline{W_2^-})\\
		&= \underline{(\bd Z)^-} - \underline{W_1^-}|_V + \underline{W_2^+} + \underline{W_2^-} - \underline{W_2^-}\\
		&= \underline{(\bd Z \sqcup - W_1|_{U \cap V} \sqcup - W_2|_{U \cap V})^-} + \underline{W_2^-} + \underline{W_2^+}\\
		&=\underline{W_2^-} + \underline{W_2^+}\\
		&=\underline{W_2}.
	\end{align*}
	So $\underline{W_2} = -\underline{\gamma}|_V$ in $H^*_\Gamma(V)$.

	Altogether, we have $(\uW_1,\uW_2) = \left(\underline{\gamma}|_U, -\underline{\gamma}|_V\right) = i\left(\underline{\gamma}\right)$.

	\textbf{Exactness at $\mathbf{H^k_\Gamma(U \cap V)}$.}
	Consider a cocycle $W$ in $C^k_\Gamma(U)$.
	Then $\delta j(\uW)$ is represented by $-W^0$.
	By the discussion above, $W^- \in PC^*_\Gamma(U \cup V)$, and $$\bd (W^-) = -W^0 \sqcup (\bd W)^-.$$ As $W$ is a cocycle, $\bd W \in Q^*(U)$, and so $(\bd W)^-$ is in $Q^*(U)$ and $Q^*(U \cup V)$ via \cref{C: creasing Q} and the preceding discussion.
	So $\bd (\underline{W^-}) = - \underline{W^0}$ in $C^*_\Gamma(U \cup V)$, and consequently
 $\delta j(\uW) = 0 \in H^*_\Gamma(U \cup V)$.
	A similar argument holds for elements of $H^k_\Gamma(V)$.

	Now suppose $\uW \in H^k_\Gamma(U \cap V)$ and $\delta(\uW) = 0$.
	Representing $\uW$ by $W$, this means there is a $Z$ in $U \cup V$ such that $\bd Z \sqcup W^0 \in Q^*(U \cup V)$.
	Let $A = Z|_U \sqcup W^+ \in PC^*_\Gamma(U)$ and $B = -Z|_V \sqcup W^- \in PC^*_\Gamma(V)$.
	Then by \cref{L: W0 cochain},
	\begin{align*}
		\bd A& = \phantom{-}\bd Z|_U \sqcup \bd( W^+) = \phantom{-}\bd Z|_U \sqcup \phantom{-} W^0 \sqcup (\bd W)^{+} \in PC^*_\Gamma(U)\\
		\bd B& = -\bd Z|_V \sqcup \bd (W^-) = -\bd Z|_V \sqcup -W^0 \sqcup (\bd W)^{-} \in PC^*_\Gamma(V).
	\end{align*}
	As $\bd Z \cup W^0 \in Q^*(U \cup V)$, its restrictions to $U$ and $V$ are in $Q^*(U)$ and $Q^*(V)$, respectively, by \cref{L: pullback map Q}.
	Furthermore, as $W$ represents a cocycle, $\bd W \in Q^*(U \cap V)$, so $(\bd W)^\pm \in Q^*(U\cap V)$ by \cref{C: creasing Q} and consequently $(\bd W)^+$ and $(\bd W)^-$ are in $Q^*(U)$ and $Q^*(V)$, respectively.
	So by \cref{R: cycles and boundaries}, the pair $(A,B)$ represents an element of $H^k_\Gamma(U) \oplus H^k_\Gamma(V)$.

	We then have in $C^*_\Gamma(U \cap V)$ that
	\begin{align*}
		j(\underline{A}, \underline{B})& = \underline{A}|_{U \cap V} + \underline{B}|_{U \cap V}\\
		& = \underline{Z}|_{U \cap V}+ \underline{W^+} - \underline{Z}|_{U \cap V} + \underline{W^-}\\
		& = \underline{W^+} + \underline{W^-}.
	\end{align*}
	The cochain $\underline{W^+} + \underline{W^-}$ represents $\uW \in H^*_\Gamma(U \cap V)$ via creasing, \cref{T: cohomology creasing}, so $\uW \in \im(j)$.
\end{proof}

\begin{example}[Cohomology of spheres]\label{E: cohomology of spheres}\index{geometric cohomology!of spheres}
Analogous to our observation about homology of spheres in \cref{E: sphere homology}, once one has a dimension axiom (\cref{E: dimension}), homotopy invariance, and Mayer-Vietoris sequences, the geometric cohomology groups of spheres follow from formal computations, and we have $H^0_{\Gamma}(S^n)\cong H^n_\Gamma(S^n)\cong \Z$ with all other cohomology groups zero.
For $n>0$, by \cref{E: sphere homology} and our strong form of Poincar\'e duality in \cref{T: PD}, $H^0_{\Gamma}(S^n)$ is generated by the class of $\id \colon S^n \to S^n$ and $H^n_{\Gamma}(S^n)$ is generated by the class of the inclusion of a point, each with either co-orientation.

We can also see this by induction, again analogously to the arguments in \cref{E: sphere homology}, but with some interesting geometric differences. Let us start with $n>1$ and write $S^n$ again as the union of open disks $U$ and $V$ that extend the hemispheres. Then the Mayer-Vietoris sequence starts out as
\begin{equation*}
0\to H^0_\Gamma(S^n) \to H^0_\Gamma(U) \oplus H^0_\Gamma(V) \to H^0_\Gamma(U\cap V) \to H^1_\Gamma(S^n) \to H^1_\Gamma(U) \oplus H^1_\Gamma(V) \to.
\end{equation*}
The extended hemispheres $U$ and $V$ are contractible and so by \cref{E: contractible} the only nontrivial cohomology of each is $H^0_\Gamma(U) \cong H^0_\Gamma(V) \cong \Z$, in each case generated by the class of the tautologically co-oriented identity map.
The tautologically co-oriented identity map $\id_{S^n} \colon S^n \to S^n$ thus represents an element of $H^0_\Gamma(S^n)$ that restricts to generators of $H^0_\Gamma(U)$ and $H^0_\Gamma(V)$.
Furthermore, the identity maps of $U$ and $V$ each restrict to the identity map $\id_{U \cap V} \colon U \cap V \to U \cap V$.
But $U \cap V \cong S^{n-1} \times (-1,1)$, so if we assume inductively that $H^0_\Gamma(S^{n-1}) \cong \Z$ with the generator represented by the identity $\id_{S^{n-1}}$, then by \cref{E: contractible} the map $\id_{U \cap V}$ generates $H^0_\Gamma(U \cap V) \cong \Z$.
As $H^1_\Gamma(U) = H^1_\Gamma(V) = 0$, it follows that $H^1_\Gamma(S^1) = 0$ and $H^0_\Gamma(S^n) \cong \Z$, generated by $\id_{S^n}$.

Continuing to assume $n>1$ and that we know $H^*_\Gamma(S^{n-1})$, all $H^i_\Gamma(S^n)$ for $i>0$ must be zero from the long exact sequence until we reach degree $n$, where the Mayer-Vietoris sequence collapses to the isomorphism $H^{n-1}_\Gamma(U \cap V) \xr{\delta} H^n_\Gamma(S^n)$.
If we assume inductively that $H^{n-1}_\Gamma(S^{n-1}) \cong \Z$ is generated by the inclusion of a point $f \colon pt \into S^{n-1}$, then, by \cref{E: contractible}, if we identify $U \cap V$ again with $S^{n-1} \times (-1,1)$ we have $H^{n-1}_\Gamma(U \cap V) \cong H^{n-1}_\Gamma(S^{n-1} \times (-1,1))\cong \Z$ generated by $f \times \id_{(-1,1)}$.
We can choose a separating function $\phi$ so that $\phi^{-1}(0) = S^{n-1} \times 0 \subset S^n$.
Then by the definition of $\delta$ in \cref{D: connecting}, the image of this generator in $H^n_\Gamma(S^n)$ is represented up to sign by the composition $f \times \id_{0} \colon pt \times 0 \to S^{n-1} \times (-1,1) \into S^n$.

The base case $n=1$ is again a bit more complicated and left as a good exercise for the reader.
\end{example}

\subsection{Geometric homology and cohomology are singular homology and cohomology}\label{S: homology is homology}

In this section, we apply a theorem of Kreck and Singhof to show that geometric homology and cohomology are isomorphic to singular homology and cohomology on smooth manifolds.

\begin{theorem}\label{T: geometric is singular}\index{geometric homology!is singular homology}\index{geometric cohomology!is singular cohomology}
	On the category of smooth manifolds (without boundary) and continuous maps, geometric homology and cohomology are respectively isomorphic to singular homology and cohomology with integer coefficients, i.e.\ $H_*^\Gamma(-) \cong H_*(- ;\Z)$ and $H^*_\Gamma(-) \cong H^*(-;\Z)$ as functors.
\end{theorem}

\begin{proof}
	This is a consequence of \cite[Proposition 10]{Krec10b} once we verify that $H_*^\Gamma$ and $H^*_\Gamma$ are respectively an ordinary homology theory and an ordinary cohomology theory on the category of smooth manifolds and continuous maps as defined in \cite{Krec10b}.
	Putting together the definitions and axioms of Sections 4 and 5 of \cite{Krec10b}, we must show the following:

	\begin{enumerate}
		\item\label{I: homotopy functor} $H_*^\Gamma$ is a covariant homotopy functor on the category of smooth manifolds (without boundary) and continuous maps between them, and $H^*_\Gamma$ is a contravariant homotopy functor on the same category.

		\item For each triple $(M;U,V)$ with $M$ a smooth manifold and $U,V$ open subsets there are exact (homological or cohomological) Mayer--Vietoris sequences with natural connecting maps $\delta$.

		\item\label{I: neg dim} For all $M$, $H_k^\Gamma(M) = H^k_\Gamma(M) = 0$ for $k<0$.

		\item The Dimension Axiom: $H_k^\Gamma(pt) = H^k_\Gamma(pt) = 0$ for $k\neq 0$.

		\item $H_*^\Gamma$ and $H^*_\Gamma$ are additive: for a manifold $M$ of dimension $0$, each $H_k^\Gamma(M)$ is canonically isomorphic to $\oplus_{x \in M} H_k^\Gamma(x)$ and each $H^k_\Gamma(M)$ is canonically isomorphic to $\prod_{x \in M} H^k_\Gamma(x)$.
	\end{enumerate}

	Axiom \ref{I: homotopy functor} holds from the definitions, \cref{T: homology homotopy functor,P: cohomology pullback,C: contra funct}.
	We have Mayer--Vietoris sequences by \cref{T: relative MV,T: absolute MV}.
	The connecting map for the cohomology sequence is natural by \cref{T: natural connection}.
	The connecting map for the homology sequence is natural just as in the standard argument for singular homology: given a map of triples $(M;U',V') \to (N;U,V)$ there is a map of Mayer--Vietoris sequences induced by a map of short exactly sequences of chain complexes of the form of diagram \eqref{E: homology MV SES} (replacing supported cochain complexes with chain complexes), itself induced by functoriality from the maps of chain complexes $C_*^{\Gamma}(U') \to C_*^{\Gamma}(U)$ and similarly for $V$ and $U \cap V$.
	The resulting map of Mayer--Vietoris sequences shows that the connecting map is natural.

	Axiom \ref{I: neg dim} holds trivially for homology as there are no non-empty prechains of degree $<0$.
	It also holds for cohomology because for $k<0$ any representing precocycle must have small rank and its boundary must be in $Q^*(M)$ by \cref{R: cycles and boundaries}.
	Thus such a precocycle is degenerate and must represent $0 \in C^k_\Gamma(M)$.

	The Dimension Axiom has been proven in \cref{E: dimension}.

	The Additivity Axiom is apparent.

	As $H_*^\Gamma$ and $H^*_\Gamma$ satisfy all the axioms, \cite[Proposition 10]{Krec10b} says that for any isomorphism $H_0(pt) \xr{\cong} H^\Gamma_0(pt)$ there is a unique isomorphism of homology theories $H_*(- ;\Z) \xr{\cong} H_*^\Gamma(-)$, and similarly for cohomology.
	As we know that $H_0^\Gamma(pt) \cong H_\Gamma^0(pt) \cong \Z$ by \cref{E: dimension}, our theorem follows.
\end{proof}

\begin{example}\index{geometric cohomology!of the projective plane}
	Let us consider $H^*_\Gamma(\R P^2)$.
	We know from the standard computations of $H^*(\R P^2)$ and the above theorem that
	\begin{equation*}
		H_\Gamma^i(\R P^2) =
		\begin{cases}
			\Z_2,&i = 2,\\
			\Z,&i = 0,\\
			0,&\text{otherwise.}
		\end{cases}
	\end{equation*}
	From the viewpoint of geometric cohomology, $H^0(\R P^2)$ is generated by the identity map $\R P^2 \to \R P^2$ with its tautological co-orientation.

	In degree $1$, cochains are represented by co-oriented maps from (unions of) closed intervals or circles.
	For a map from the circle to be co-oriented it cannot represent the non-trivial element $\alpha \in \pi_1(\R P^2)$, and so it must be contractible.
	Any map from a disk must be co-orientable, as the disk is contractible (this can be considered an extension of \cref{L: co-orientable homotopies}).
	Thus by smooth approximation, any smooth co-oriented map from the circle to $\R P^2$ is the boundary of a smooth co-oriented map from the disk, and so represents $0$ in cohomology.
	The same is true for any circle decomposed into intervals that does not represent $\alpha \in \pi_2(\R P^2)$; via homotopy and creasing, such a ``circle'' is the boundary of a polygon.
	On the other hand, consider a map $g \colon [0,1] \to \R P^2$ that does represent $\alpha$.
	As $[0,1]$ is contractible, such a map can always be co-oriented.
	We can let $e$ represent the standard unit vector of $[0,1]$ and suppose that $g$ is co-oriented so that at $0 \in [0,1]$ the co-orientation is represented by $(e,\beta)$ for some local orientation $\beta$ of $\R P^2$ at $g(0)$.
	Traversing the path, the representation for the co-orientation at $1$ is then $(e,-\beta)$.
	Recalling our boundary conventions, the boundary of $g$ therefore consists of the maps $g|_0:0 \to \R P^2$ co-oriented by $(1,e)*(e,\beta) = (1,\beta)$ and $g|_1:1 \to \R P^2$ co-oriented by $(1,-e)*(e,-\beta) = (1,\beta)$.
	Thus with these co-orientations, $g|_0:0 \to \R P^2$ and $g|_1:1 \to \R P^2$ represent isomorphic manifolds over $\R P^2$.
	So $g$ does not represent a cocycle.
	On the other hand, as a co-oriented map from any single point to $\R P^2$ is a $2$-cocycle, this shows that by choosing the image to be the basepoint for $\pi_1(\R P^2)$ that twice any such map is a boundary.
	We leave it to the reader to verify that any two such points generate the same cohomology class and that a map from a single point cannot represent a coboundary, and so we have verified that $H^1_\Gamma(\R P^2) = 0$ and $H_\Gamma^2(\R P^2) \cong \Z_2$.
\end{example}

\subsubsection{A more direct comparison of singular and geometric homology}\label{S: homology direct}

While \cref{T: geometric is singular} says that geometric homology and cohomology are isomorphic to singular homology and cohomology, it does not provide actual isomorphisms, which we would like to have.
In this section, we provide a second proof that singular and geometric homology are isomorphic, in this case providing the isomorphism.
We will consider the case of cohomology in \cref{S: transversality}.

The basic idea for homology is that, as simplices are manifolds with corners and as the standard model simplex comes equipped with an orientation, smooth singular simplices are geometric chains.
In other words, if we let $S^{sm}_*(M)$\index{$S^{sm}_*(M)$}\index{smooth singular chain|(} denote the chain complex of smooth singular chains on the manifold without boundary $M$, then the generators, which are smooth singular simplices $\sigma \colon \Delta^i \to M$, are prechains, and so they represent elements of $C_*^\Gamma(M)$.
As $S^{sm}_*(M)$ is freely generated by the smooth singular simplices, this determines a map of graded abelian groups $S^{sm}_*(M) \to C_*^\Gamma(M)$.
We will see that there is consistency of the boundary orientations so that this is a chain map $S^{sm}_*(M) \to C_*^\Gamma(M)$.
We will show that it is, in fact, a quasi-isomorphism, i.e.\ it induces isomorphisms of homology groups.
Combined with the well-known fact that the inclusion $S^{sm}_*(M) \into S_*(M)$ is a chain homotopy equivalence if $S_*(M)$ is the full singular chain complex on $M$ \cite[Theorem 18.7]{Lee13}, this will provide a concrete chain of isomorphisms $H_*(M) \cong H_*^\Gamma(M)$.

It will be useful in what follows to have the analogous result for cubical singular homology in addition to the more common simplicial singular homology.
Details of cubical singular homology can be found, for example, in \cite{Mas91} or \cite[Section 8.3]{HW60}, though we provide a brief overview.

Rather than maps of simplices $\Delta^k \to M$, the cubical singular chain complex,\index{cubical singular chain} which we denote $SK_*(M)$,\index{$SK_*(M)$} is generated by maps $\interval^k \to M$\index{$I$@$\interval^k$} with $\interval$ being the standard interval $\interval = [0,1]$ (it will be convenient henceforth to use this notation for the interval in this context).
The boundary formula is defined so that if $\sigma: \interval^k \to M$ is a singular cube, then
\begin{equation}\label{E: cube bd}
	\bd \sigma = \sum_{i = 1}^k (-1)^i(\sigma \delta_i^0-\sigma \delta^1_i),
\end{equation}
where for $\epsilon\in\{0,1\}$, the map $\delta_i^\epsilon \colon \interval^{k-1} \to \interval^k$ is defined by
$$\delta_i^\epsilon(x_1,\ldots,x_{k-1}) = (x_1,\ldots,\epsilon,\ldots, x_k)$$
with $\epsilon$ in the $i$th slot.
The homology of the chain complex $SK_*(M)$ is not isomorphic to singular homology, as it does not satisfy the dimension axiom, so one instead forms the normalized complex $NK_*(M)$\index{$NK_*(M)$}\index{cubical singular chain!normalized} by quotienting out the subcomplex of \textit{degenerate singular cubes} that is generated by singular cubes $\sigma \colon \interval^k \to M$ such that $\sigma$ does not depend on at least one of the variables.
In other words, the degenerate singular cubes are those maps $\sigma \colon \interval^k \to M$ that factor through one of the coordinate projections $\interval^k \to \interval^{k-1}$.
It then holds that $NK_*(M)$ is chain homotopy equivalent to $S_*(M)$, the complex of singular simplicial chains.
In fact, this holds for $M$ any space and not just a manifold \cite[Theorem 8.4.7]{HW60}.
Of course we will need the smooth version $NK^{sm}_*(M)$ generated by smooth singular cubes and modulo degenerate smooth singular cubes.
We defer to an appendix to this section the proof that $NK^{sm}_*(M) \into NK_*(M)$ is a chain homotopy equivalence.

As the cubes $\interval^k$ are compact manifolds with corners equipped with their standard orientations,
we have, analogously to the simplicial case, maps $SK^{sm}_i(M) \to C_i^\Gamma(M)$ for each $i$.
We will show that this is a chain map and that it induces a chain map $NK^{sm}_*(M) \to C_*^\Gamma(M)$.

\begin{lemma}\label{L: singular and cubical chain maps}
	The maps $S^{sm}_*(M) \to C_*^\Gamma(M)$ and $SK^{sm}_*(M) \to C_*^\Gamma(M)$ are well-defined chain maps.
	The latter induces a chain map $NK^{sm}_*(M) \to C_*^\Gamma(M)$.
\end{lemma}

\begin{proof}
	We first consider the simplicial case, for which we need only verify that our boundary orientation convention for manifolds with corners is consistent with the simplicial boundary formula.

	We will think of the model $k$-simplex $\Delta^k$ as being the convex hull in $\R^k$ of the origin and the points $(0, \ldots , 1, \ldots, 0)$ as the $1$ varies through the $k$ coordinates.
	We will also use the abstract simplicial notation $[0,\ldots,k]$ when convenient, identifying the vertex $0$ with the origin and the vertex $i$ with the point $(0, \ldots , 1, \ldots, 0)$ with $1$ in the $i$th place.
	Then, letting $e_i$ denote the $i$th standard unit basis vector, the standard orientation of $\Delta^k$ is $e_1 \wedge \cdots \wedge e_k$.
	The simplicial boundary formula, written in terms of manifolds with corners, is then
	\[\bd [0,\ldots,k] = \amalg_i (-1)^i [0,\ldots, \hat i,\ldots, k].\]
	Let $F_i$ denote the face $[0,\ldots, \hat i,\ldots, k]$.
	We first consider the $F_i$ for $i > 0$.
	In this case, the vector $-e_i$ can be taken as our outward pointing normal to $F_i$ and the remaining $e_j$ are in the tangent space to $F_i$.
	So, by definition, the boundary orientation $\beta_i$ of $F_i$ as a manifold with corners boundary is the one such that $-e_i \wedge \beta_i$ agrees with $e_1 \wedge \cdots \wedge e_k$.
	We observe that
	\[ -e_i \wedge e_1 \wedge \cdots \wedge \hat e_i \wedge \cdots \wedge e_k = (-1)^ie_1 \wedge \cdots \wedge e_k,\]
	so $\beta_i = (-1)^i e_1 \wedge \cdots \wedge \hat e_i \wedge \cdots \wedge e_k$.
	Thus the corresponding simplicial boundary term is $(-1)^i [0,\ldots, \hat i,\ldots, k]$, as desired, observing that the standard simplicial orientation of $[0,\ldots, \hat i,\ldots, k]$ is $e_1 \wedge \cdots \wedge \hat e_i \wedge \cdots \wedge e_k$.
	Now, we consider $F_0$.
	For $2\leq j\leq k$, let $w_j$ be the vector from the vertex $(1,0,\ldots,0)$ to the $j$th vertex of $\Delta^k$.
 	These $w_j$ span the tangent space of $F_0$, and the simplicial orientation of $F_0$ is $w_2 \wedge \cdots \wedge w_k$.
	Taking $e_1$ to be an outward pointing normal at $F_0$, the orientation $\beta_0$ of $F_0$ as a manifold with corners boundary is the one such that $e_1 \wedge \beta_0 = e_1 \wedge \cdots \wedge e_k$.
	But now we observe that projection to the $x_1 = 0$ plane takes $w_j$ to $e_j$, which is sufficient to conclude that $e_1 \wedge w_2 \wedge \cdots \wedge w_k = e_1 \wedge \cdots \wedge e_k$, so $\beta_0 = w_2 \wedge \cdots \wedge w_k$, and this corresponds simplicially to $[1,\ldots, k]$, as desired in the simplicial boundary formula.

	We turn now to cubes.

	Any degenerate singular cube $\sigma \colon \interval^k \to M$ is also degenerate as a prechain in the sense of \cref{D: equiv triv and small}.
	In fact, it will have small rank as it filters through a projection.
	Furthermore, if that projection collapses the $i$th coordinate, then each face $\sigma \delta_j^\epsilon$ for $j\neq i$ will also be a degenerate small cube and so have small rank, while the term $\pm (\sigma \delta_i^0-\sigma \delta_i^1)$ will be trivial with the trivializing map $\rho$ being the interchange of the two faces.
	Thus, the degenerate smooth singular cubes are elements of $Q_*(M)$, and our map is well defined in each degree.

	For compatibility of boundary orientations, consider the $k-1$ dimensional face $F$ of $\interval^k$ given by $x_i = j$ with $j\in\{0,1\}$.
	Then we have an outward pointing vector given by $(-1)^{j+1}e_i$, where $e_i$ is the standard unit vector in the $i$th direction.
	The boundary orientation $\beta_F$ for $F$ as a manifold with corners is the one such that
	$(-1)^{j+1}e_i \wedge \beta_F$ is the orientation $e_1 \wedge\cdots\wedge e_k$ of $\interval^k$.
	Thus $\beta_F = (-1)^{j+1+i-1}e_1 \wedge \cdots \wedge \hat{e}_i \wedge \cdots e_j$.
	On the other hand, $e_1 \wedge \cdots \wedge \hat{e}_i \wedge \cdots e_k$ is precisely the standard orientation of $F$ when considering $\interval^k$ as a cubical complex.
	So the manifold with corners boundary orientation of $F$ is $(-1)^{i+j}$ times its standard orientation as a cube.
	But this corresponds precisely to the formula for the boundary in $K^{sm}_*(M)$ given by \eqref{E: cube bd}.
\end{proof}

\begin{theorem}\label{T: hom iso map}\index{geometric homology!via smooth singular simplices}\index{geometric homology!via smooth singular cubes}\index{geometric homology!is singular homology}
	Let $M$ be a manifold without boundary.
	The maps $$H_*(S^{sm}_*(M)) \to H_*^\Gamma(M)$$ and $$H_*(NK^{sm}_*(M)) \to H_*^\Gamma(M)$$ obtained by treating smooth singular simplices and smooth singular cubes as elements of $C_*^\Gamma(M)$ are isomorphisms.
\end{theorem}

\begin{proof}
	Both $H_*(S^{sm}_*(M))$ and $H_*(NK^{sm}_*(M))$ are isomorphic to the standard singular homology groups, so in the following we simply write $H_*$ for either of these theories and provide a uniform proof.
	We let $\Phi$ denote either of the chain maps of \cref{L: singular and cubical chain maps}, using the normalized chains in the cubical case.

	We apply \cite[Theorem 5.1.1]{Frie20}, which is based on standard Mayer--Vietoris techniques.
	We first observe that on the category consisting of the open subsets of $M$, the maps $\Phi: H_*(-) \to H_*^\Gamma(-)$ determined by the chain maps of \cref{L: singular and cubical chain maps} provide natural transformations of functors, which is clear from the definitions.
	We need to verify the following three properties:

	1.
	On $\emptyset$ or $U \subset M$ with $U$ homeomorphic to $\R^m$, the map $\Phi: H_*(U) \to H_*^\Gamma(U)$ is an isomorphism.
	As both $H_*$ and $H_*^\Gamma(-)$ are homotopy functors, we know from the respective Dimension Axioms (see \cref{E: dimension}) that in this case $H_k(U) = H_k^\Gamma(U) = 0$ for $k\neq 0$, while for $k = 0$ we have the commutative diagram
	\[
	\begin{tikzcd}
		H_0(pt) \arrow[r] \arrow[d, "\Phi"] & H_0(U) \arrow[d, "\Phi"] \\
		H_0^\Gamma(pt) \arrow[r] & H^\Gamma_0(U).
	\end{tikzcd}
	\]
	The horizontal maps are isomorphisms because these are homotopy functors, and the left hand vertical map is an isomorphism because $\Phi$ takes a generator of $H_0(pt) \cong \Z$ to a generator of $H_0^\Gamma(pt) \cong \Z$; see again \cref{E: dimension}.
	So the right hand map is also an isomorphism.

	The case of the empty set is clear as in this case both groups are trivial.

	2.
	$\Phi$ induces a commutative diagram of long exact Mayer--Vietoris sequences.
	This follows from basic homological algebra given the commutativity of the following diagram and its analogue for normalized smooth singular cubical chains
	\[
	\begin{tikzcd}[column sep=large]
		S^{sm}_*(U \cap V) \arrow[r, "{(i_U, -i_V)}", hook] \arrow[d, "\Phi"] & S^{sm}_*(U) \oplus S^{sm}_*(V) \arrow[d, "\Phi \oplus \Phi"] \\
		C_*^\Gamma(U \cap V) \arrow[r, "{(i_U, -i_V)}", hook] & C_*^\Gamma(U) \oplus C_*^\Gamma(V).
	\end{tikzcd}
	\]

	3.
	If $\{U_\alpha\}$ is an increasing collection of open submanifolds of $M$ such that $\Phi \colon H_*(U_\alpha) \to H_*^\Gamma(U_\alpha)$ is an isomorphism for all $\alpha$, then $\Phi \colon H_*(\cup_\alpha U_\alpha) \to H_*^\Gamma(\cup_\alpha U_\alpha)$ is an isomorphism.
	This argument is standard given that both singular (simplicial or cubical) chains and geometric chains are represented by compact spaces: If $W$ represents a cycle in $C_*^\Gamma(\cup_\alpha U_\alpha)$, then $W \to \cup_\alpha U_\alpha$ factors through some particular $U_\beta$, so, as $H_*(U_\beta) \xr{\Phi}H_*^\Gamma(U_\beta)$ is an isomorphism, $\uW$ is in the image of $H_*(U_\beta) \xr{\Phi}H_*^\Gamma(U_\beta) \to H_*^\Gamma(\cup_\alpha U_\alpha)$.
	But then $\uW$ is in the image of the composition $H_*(U_\beta) \to H_*(\cup_\alpha U_\alpha) \to H_*^\Gamma(\cup_\alpha U_\alpha)$, so $\Phi \colon H_*(\cup_\alpha U_\alpha) \to H_*^\Gamma(\cup_\alpha U_\alpha)$ is surjective.
	Similarly, if $\Phi \colon H_*(\cup_\alpha U_\alpha) \to H_*^\Gamma(\cup_\alpha U_\alpha)$ maps a class represented by a smooth singular cycle $\xi$ to $0$, then $\xi$ bounds as a geometric cycle, so by \cref{R: cycles and boundaries} there is some $Z \in PC_*^\Gamma(\cup_\alpha U_\alpha)$ so that $\bd Z \amalg -\xi \in Q_*(\cup_\alpha U_\alpha)$.
	But by compactness, there is some $\beta$ so that $Z$ and $\xi$ both have image in $U_\beta$.
	So $\xi$ represents a class in $H_*(U_\beta)$ that maps to $0$ in $H_*^\Gamma(U_\beta)$.
	As $\Phi$ is assumed an isomorphism on $U_\beta$, it must be that $\xi$ represents $0$ in $H_*(U_\beta)$, and so it also represents $0$ in $H_*(\cup_\alpha U_\alpha)$.

	It now follows from Theorem \cite[5.1.1]{Frie20} that $\Phi: H_*(M) \to H_*^\Gamma(M)$ is an isomorphism.
\end{proof}

\cref{T: hom iso map} is stated without proof in \cite[Section 10]{Lipy14}, where
Lipyanskiy writes, ``The fact that the natural maps induce isomorphisms follow from the standard Mayer--Vietoris arguments.''

Unfortunately, providing a direct comparison for cohomology theories is not so straight forward as there is no obvious map between $C^*_\Gamma(M)$ and $S^*(M) = \Hom(S_*(M),\Z)$.
It will take some work in the following sections to develop a geometric connection between these cohomology theories.
\index{smooth singular chain|)}

\subsection{Appendix: smooth singular cubes}\label{S: smooth cubes}\index{smooth singular cube|(}

In this appendix, we show that normalized singular cubical chains and normalized smooth singular cubical chains provide the same homology theory, analogous to the case for the perhaps more familiar simplicial chains.

\begin{proposition}\label{P: singular smooth cubes}
	The inclusion $\psi: NK^{sm}_*(M) \into NK_*(M)$ is a chain homotopy equivalence.
\end{proposition}

\begin{proof}
	The proof is analogous to the simplicial case as given in detail in \cite[Theorem 18.7]{Lee13}, though we need to take care with degenerate cubes, which is not an issue in the simplicial case.
	To account for this, we sketch the proof in \cite{Lee13} but provide some detailed modifications.

	We first observe that the map $SK^{sm}_*(M) \to NK_*(M)$ takes a smooth singular cubical chain to $0$ only if all of its cubes (with non-zero coefficient) are degenerate, and so we do have an injection $NK^{sm}_*(M) \into NK_*(M)$.
	We will define cube-wise a chain homotopy inverse $s \colon NK_*(M) \to NK^{sm}_*(M)$ by starting with a map $\td s \colon SK_*(M) \to SK^{sm}_*(M)$ and passing to quotients.

	Recall from \cref{S: homology is homology} that we write $\delta_i^\epsilon$ for the face inclusions of the standard cubes.
	If $\sigma \colon \interval^k \to M$ is a singular cube, we define homotopies $H_\sigma \colon \interval^k \times \interval = \interval^{k+1} \to M$ so that the following properties hold:
	\begin{enumerate}
		\item\label{I: smooth} $H_\sigma$ is a homotopy from $\sigma$ to a smooth map $\td \sigma \colon \interval^k \to M$.

		\item\label{I: faces} $H_{\sigma \delta_i^\epsilon} = H_\sigma \circ (\delta_i^\epsilon \times \id_\interval)$ so that the construction is compatible along faces.
		More explicitly, $H_{\sigma \delta_i^\epsilon}(x,t) = H_\sigma(\delta_i^\epsilon(x),t)$.

		\item\label{I: already smooth} If $\sigma$ is already smooth then $H_{\sigma}(x,t) = \sigma(x)$, i.e.\ the homotopy is constant.

		\item\label{I: degen} If $\sigma$ is independent of the coordinate $x_i$, then so is $H_\sigma$.
	\end{enumerate}

	The last condition, which we have added for cubes, ensures that if $\sigma$ is degenerate so will be $H_\sigma$ and $\td \sigma$.

	The construction is by induction on dimension.
	If $\sigma$ is a $0$-cube, then we define $H_\sigma(x,t) = \sigma(x)$, the constant homotopy.
	This satisfies the conditions.
	We then assume $H_\sigma$ defined with these properties for all cubes of dimension $<k$ and extend the definition to $k$-cubes.
	If $\sigma$ is already smooth, then the constant homotopy $H_\sigma(x,t) = \sigma(x)$ satisfies the conditions, noting that if $\sigma$ is smooth then so is each $\sigma \circ \delta_i^\epsilon$.
	If $\sigma$ is not smooth, we consider separately the two cases when $\sigma$ is degenerate or nondegenerate.

	First suppose $\sigma$ is not degenerate.
	By the induction hypothesis and Condition \eqref{I: faces}, $H_\sigma$ is determined on the subspace $(\interval^k \times 0) \cup ((\bd \interval^k) \times \interval)$.
	One can check as in the proof of \cite[Lemma 18.8]{Lee13} that Condition \eqref{I: faces} guarantees that the faces glue to form a continuous map.
	As $\bd \interval^k \into \interval^k$ is a the inclusion of a CW subcomplex, it is a cofibration, so there is a retraction $\interval^k \times \interval \to (\interval^k \times 0) \cup ((\bd \interval^k) \times \interval)$ \cite[Theorem VII.1.3]{Bred97}.
	The composition of the retraction with $H_\sigma$ as defined on the subspace determines a homotopy $F \colon \interval^k \times \interval \to M$ such that $F(-,1)$ is smooth on each $k-1$ face of $\bd \interval^k$.
	In fact, this implies that $F(-,1)$ is smooth on all of $\bd \interval^k$ by a minor modification of \cite[Lemma 18.9]{Lee13}.
	So by the Whitney Approximation Theorem \cite[Theorem 6.26]{Lee13}, there is a homotopy rel $\bd \interval^k$ from $F(-,1)$ to a smooth map $\td \sigma \colon \interval^k \to M$; we denote this homotopy $G$.
	Finally, let $u \colon \interval^k \to (0,1]$ be a continuous function that takes $\bd\interval^k$ to $1$ and the interior of the cube to $(0,1)$.
	Then we can define
	\begin{equation*}
		H_\sigma(x,t) =
		\begin{cases}
			F\left(x,\frac{t}{u(x)}\right),&x \in \interval^k, 0 \leq t \leq u(x),\\
			G\left(x,\frac{t-u(x)}{1-u(x)}\right),&x \in \text{Int}(\interval^k), u(x) \leq t \leq 1.\\
		\end{cases}
	\end{equation*}
	One can check as in the proof of \cite[Lemma 18.8]{Lee13} that this is a continuous homotopy that satisfies the first two conditions above, as required.

	Next, suppose $\sigma$ is degenerate, i.e.\ there is some coordinate $x_i$ so that $\sigma$ does not depend on $x_i$.
	Let $\pi_i \colon \interval^k\to\interval^{k-1}$ be given by $\pi_i(x_1,\ldots, x_k) = (x_1,\ldots, \hat x_i, \ldots, x_k)$ with the $x_i$ term omitted.
	In this case, we define $H_\sigma(x,t)$ by
	\[H_\sigma(x,t) = H_{\sigma \delta_i^0}(\pi_i(x),t) = H_{\sigma \delta_i^1}(\pi_i(x),t).\]
	The two expressions on the right are equal due to $\sigma$ not depending on $x_i$.
	We claim that if there are multiple coordinates of which $\sigma$ is independent then this definition is independent of the choice of such coordinate.
	This is clear for $1$-cubes for which there is only one possible coordinate.
	Suppose then the claim proven in dimensions $<k$ and that $\sigma \colon \interval^k \to M$ is independent of $x_j$ and $x_i$ with $j<i$.
	Since $\sigma$ is independent of $x_j$, so is $\sigma \circ \delta_i^0$, so inductively $H_{\sigma \delta_i^0}(\pi_i(x),t) = H_{\sigma \delta_i^0\delta_j^0}(\pi_j\pi_i(x),t)$.
	Similarly, using that the $i$th coordinate of the cube is the $i-1$-st coordinate of the $j$th faces, we have $H_{\sigma \delta_j^0}(\pi_j(x),t) = H_{\sigma \delta_j^0\delta_{i-1}^0}(\pi_{i-1}\pi_j(x),t)$.
	But $\delta_i^0\delta_j^0$ and $\delta_j^0\delta_{i-1}^0$ determine the same $k-2$ face of $\interval^k$, and $\pi_j\pi_i(x) = \pi_{i-1}\pi_j(x)$.
	So both constructions give the same $H_\sigma$.

	With this definition of $H_\sigma$, Conditions \eqref{I: smooth}, \eqref{I: already smooth}, and \eqref{I: degen} hold by construction and by induction.
	We must verify Condition \eqref{I: faces}.
	If $\sigma$ is independent of $x_i$, the condition is clear by construction for the faces $\sigma\delta_i^0$ and $\sigma\delta_i^1$.
	For $j\neq i$, first suppose $i<j$.
	As $\sigma$ is independent of $x_i$, so is $\sigma\delta_j^\epsilon$, so
	\begin{align*}
		H_{\sigma\delta_j^\epsilon}(x,t)& = H_{\sigma\delta_j^\epsilon \delta_i^0}(\pi_i(x),t)\\
		& = H_{\sigma \delta_i^0\delta_{j-1}^\epsilon}(\pi_i(x)),t)\\
		& = H_{\sigma \delta_i^0}(\delta_{j-1}^\epsilon\pi_i(x)),t)\\
		& = H_{\sigma \delta_i^0}(\pi_i(\delta_j^\epsilon(x)),t)\\
		& = H_\sigma(\delta_j^\epsilon(x),t).
	\end{align*}
	Here the first equality uses our definition of $H_{\sigma\delta_j^\epsilon}$ as $\sigma\delta_j^\epsilon$ is independent of $x_i$.
	The second equality is an identity for cubical face inclusions; see \cite[Section 4]{GrMa03} or below in \cref{S: cubes}.
	The third is Condition \eqref{I: faces} for $H_{\sigma \delta_i^0}$, which holds by induction hypothesis.
	The fourth equality is another cubical identity, and the last is the definition of $H_\sigma$.

	Similarly, if $j<i$, then $\sigma\delta_j^\epsilon$ is independent of its $i-1$-st coordinate, and we compute analogously:

	\begin{align*}
		H_{\sigma\delta_j^\epsilon}(x,t)& = H_{\sigma\delta_j^\epsilon \delta_{i-1}^0}(\pi_{i-1}(x),t)\\
		& = H_{\sigma \delta_{i}^0\delta_{j}^\epsilon}(\pi_{i-1}(x)),t)\\
		& = H_{\sigma \delta_{i}^0}(\delta_{j}^\epsilon\pi_{i-1}(x)),t)\\
		& = H_{\sigma \delta_{i}^0}(\pi_{i}(\delta_j^\epsilon(x)),t)\\
		& = H_\sigma(\delta_j^\epsilon(x),t).
	\end{align*}

	This completes our construction of the homotopies $H_\sigma$ satisfying the desired conditions.
	We can now define $\td s \colon SK_*(M) \to SK^{sm}_*(M)$ by $\td s(\sigma) = H_\sigma(-,1)$.
	This is a chain map thanks to Condition \eqref{I: faces}:
	\begin{align*}
	\bd \td s(\sigma) &= \bd H_\sigma(-,1)\\
	&=\sum_{i=1}^k (-1)^i (H_\sigma(-,1) \delta_i^0 - H_\sigma(-,1) \delta_i^1)\\
	&=\sum_{i=1}^k (-1)^i (H_{\sigma \delta_i^0}(-,1) - H_{\sigma \delta_i^1}(-,1))\\
	&=\sum_{i=1}^k (-1)^i (\td s(\sigma \delta_i^0) - \td s(\sigma \delta_i^1))\\
	&= \td s(\bd \sigma).
	\end{align*}

	If $\td \psi: SK_*^{sm}(M) \to SK_*(M)$ is the inclusion, we have by definition that $\td s \td \psi = \id$.
	We show that $\td\psi\td s$ is chain homotopic to the identity.\footnote{Here, finally, is a step that is easier in the cubical setting as we do not need to subdivide prisms into simplices.}
	Indeed, if $\sigma$ is a singular $k$-cube then treating $H_\sigma$ as a singular $k+1$ cube we have

	\begin{align*}
		\bd H_\sigma& = \sum_{i = 1}^{k+1} (-1)^i\left(H_\sigma \delta_i^0-H_\sigma \delta^1_i\right)\\
		& = \left(\sum_{i = 1}^{k} (-1)^i\left(H_\sigma (\delta_i^0 \times \id_\interval)-H_\sigma (\delta^1_i \times \id_\interval)\right)\right) +(-1)^{k+1}(H_\sigma(-,0)-H_\sigma(-,1))\\
		& = \left(\sum_{i = 1}^{k} (-1)^i\left(H_{\sigma \delta_i^0}-H_{\sigma\delta^1_i}\right)\right) +(-1)^{k+1}(\sigma -\td \psi\td s(\sigma)).
	\end{align*}
	So if we define $\td J(\sigma) = (-1)^{k+1}H_\sigma$ for a singular $k$-cube $\sigma$, we obtain
	\begin{align*}
		(-1)^{k+1}\bd \td J(\sigma)& = \left(\sum_{i = 1}^{k} (-1)^i\left( (-1)^k\td J(\sigma \delta_i^0)-(-1)^k\td J(\sigma\delta^1_i)\right)\right) +(-1)^{k+1}\left(\sigma(-)-\td \psi\td s(\sigma)\right)\\
		& = (-1)^k\td J(\bd \sigma)+(-1)^{k+1}(\sigma - \td \psi\td s(\sigma)),
	\end{align*}
	so $$\bd \td J(\sigma)+\td J(\bd \sigma) = \sigma - \td \psi\td s(\sigma),$$
	which shows that $\td \psi\td s$ is chain homotopic to the identity.

	Finally, we note that $\td \psi$, $\td s$, and $\td J$ all take degenerate simplices to degenerate simplices by definition and construction, so $\td \psi$ and $\td s$ descend to chain maps $\psi \colon NK_*^{sm}(M) \to NK_*(M)$ and $s \colon NK_*(M) \to NK_*^{sm}(M)$ with $s\psi = \id$ and $\td J$ descends to a chain homotopy $J \colon NK_*(M) \to NK_{*+1}(M)$.
\end{proof}
\index{smooth singular cube|)}

\subsection{Geometric homology with coefficients}

In this brief section, we consider geometric homology with coefficients in an arbitrary abelian group $G$.

\begin{definition}\index{geometric chain!with coefficients|textbf}\index{geometric homology!with coefficients}
	Let $M$ be a manifold with corners and $G$ an abelian group.
	We define the chain complex $C_*^\Gamma(M;G)$ of \textbf{geometric chains with coefficients in $G$} to be $C_*^\Gamma(M) \otimes_\Z G$.
	The corresponding \textbf{geometric homology groups with coefficients}\index{geometric homology!with coefficients} will be denoted $H_*^\Gamma(M;G)$.
\end{definition}

In the remainder of the section, we will simply write $\otimes$ rather than $\otimes_Z$.

As $C_*^\Gamma(M)$ and its subgroups are torsion-free, and hence flat as $\Z$-modules, the Universal Coefficient Theorem applies, and so for each $i$, $$H_i^\Gamma(M;G) \cong \left(H_i^\Gamma(M) \otimes G \right) \oplus \left(H_{i-1}^\Gamma(M) * G \right),$$
letting $*$ denote the torsion product as in \cite{Munk66}.

Generalizing \cref{T: hom iso map,P: singular smooth cubes}, we have the following, which says that geometric homology with coefficients is ordinary homology with coefficients.

\begin{theorem}\label{T: hom iso map G}\index{geometric homology!with coefficients!is singular homology with coefficients}
	Let $M$ be a manifold without boundary.
	The maps $H_*(S_*(M;G)) \leftarrow H_*(S^{sm}_*(M;G)) \to H_*^\Gamma(M;G)$ and $H_*(NK_*(M;G)) \leftarrow H_*(NK^{sm}_*(M;G)) \to H_*^\Gamma(M;G)$ obtained by tensoring the corresponding integer-coefficient chain complex inclusions with $G$ are isomorphisms.
\end{theorem}
\begin{proof}
	We note again that $C_*^\Gamma(M)$ and its subgroups are torsion-free, and hence flat as $\Z$-modules.
	The other chain complexes are all free; in the case of the normalized cubical complexes, they are generated by the nondegenerate singular cubes, which in each case constitute a free summand of the (smooth) group of singular cubes.
	So all of our chain maps induce maps of universal coefficient short exact sequences, and the claim follows from the Five Lemma.
\end{proof}


\section{Interaction with cubical structures}\label{S: transversality}

In this section, we bring in some auxiliary structures that will help us further develop geometric cohomology and its connections to singular cohomology.
In particular, we equip our manifolds with smooth cubulations.
Many of our results would apply just as well with the more familiar smooth triangulations, but we find cubulations to be more convenient.
In particular, in \cite{FMS-flows}, we have considered geometric cochains in the presence of cubulations, demonstrating how to obtain a fully-defined cochain-level cup product via intersection using certain flows developed in terms of the cubulation.
Cup products in geometric cohomology will be discussed in the following section.

Smooth cubulations are analogous to smooth triangulations in that they involve a homeomorphism $M \cong |X|$ between a manifold $M$ and the geometric realization $|X|$ of a cubical complex $X$, such that the restriction to each cubical face is a smooth embedding.
What is slightly different, aside from substituting cubes (i.e.\ copies of $\interval^k$) for simplices, is that cubical complexes require a bit more combinatorial structure; simplicial complexes can always be constructed just by gluing together simplices along faces, while gluing of cubical faces requires a certain combinatorial compatibility among the faces being glued.

The issue is that any simplicial complex can be given a total ordering of its vertices, and this ordering provides a canonical identification between any simplicial face and the standard model simplex of the same dimension.
By contrast, the natural combinatorial structure on the vertices of the standard cube is not a total ordering but rather a partial ordering.
In particular, if we take the standard cube to be $\interval^k = [0,1]^k \subset \R^k$, then we have $v \leq w$ for two vertices if each coordinate of $v$ is less than or equal to the corresponding coordinate of $w$.
There turn out to be spaces obtained from naively gluing cubes that do not support compatible partial orderings of this type among their cubes; for example, see \cref{F: cubical structure}, below.

So when we speak of cubical complexes we will restrict ourselves to complexes that do admit such combinatorial data.
Consequently, each cubical $k$-face comes equipped with an identification with the standard $k$-cube, and hence also a standard orientation.
As we will note below, smooth cubulations of this form exist for any smooth manifold.
In the remainder of this work, ``cubulation'' will always mean a smooth cubulation.

Analogously to simplicial complexes, cubical complexes possess algebraic cubical chain and cochain complexes and so cubical homology and cohomology that coincide with singular homology and cohomology\footnote{We will show below that cubical homology coincides with singular cubical homology, which coincide with simplicial singular homology by \cite{EM53}.
As all of the involved chain complexes are free, the corresponding cohomologies are also isomorphic by basic homological algebra \cite[Theorem 45.5]{Mun84}.\label{FN: cubical and singular}}.
Our primary goal in this section is to see that there are direct geometrically-defined isomorphisms between cubical (co)homology and geometric (co)homology.

For this, we first provide some background on cubical complexes and cubical homology and cohomology in \cref{S: cubes,S: cubical cochains}.
Then in \cref{S: cubical and geometric homology} we show that the obvious map that takes a face of a cube complex to its corresponding geometric chain induces an isomorphism from cubical homology to geometric homology.
Next, in \cref{S: transverse cochains}, we consider those geometric cochains that are transverse to a given cubulation and show that their cohomology agrees with the geometric cohomology obtained without that constraint.

The motivation for our interest in cochains that are transverse to the cubulation is that they allow us to define an \textit{intersection map} $\mc I$ from these transverse geometric cochains to the cubical cochains.
If $F$ is a face of the cubulation, $F^*$ its dual cochain, and $W$ represents a geometric cochain of complementary dimension to $F$, then the coefficient of $F^*$ in $\mc I(W)$ is simply the geometric intersection number of $W$ with $F$.
This intersection map is defined in \cref{S: intersection map}, which also contains our proof that the intersection map induces a cohomology isomorphism when $H^*(M)$ is finitely generated in each degree.
To implement this proof, we include in \cref{S: dual cubes} a discussion of what we call ``central subdivisions'' of cubical complexes, which are analogous to barycentric subdivisions of simplicial complexes and allow us to construct the cubical dual cells to faces of the cubulation.

\subsection{Cubical complexes and cubulations}\label{S: cubes}

We begin by recalling some notation from \cite{FMS-flows}.
In the context of cubical complexes we write the unit interval as $\interval = [0,1]$ and define the \textbf{standard $n$-cube}\index{cube} to be
\begin{equation*}
	\interval^n = \big\{ (x_1, \dots, x_n) \in \R^n\ |\ 0 \leq x_i \leq 1 \big\}.
\end{equation*}
Denote the set $\{1, \dots, n\}$ by $\overline{n}$.

Given a partition $F = (F_0, F_{01}, F_1)$\index{cube!partition of coordinates} of $\overline n$, it determines a \textbf{face}\index{cube!face|textbf} of $\interval^n$ given by
\begin{equation*}
	\{(x_1, \dots, x_n) \in \interval^n\ |\ \forall \varepsilon \in \{0, 1\},\ i \in F_\varepsilon \Rightarrow x_i = \varepsilon\}.
\end{equation*}
We abuse notation and write $F$ for both the partition and its associated face.
We refer to coordinates $x_i$ with $i \in F_{01}$ as \textbf{free}\index{cube!free coordinates|textbf} and to the others as \textbf{bound}\index{cube!bound coordinates|textbf}.
The \textbf{dimension} of $F$ is its number of free coordinates, and as usual the faces of dimension $0$ and $1$ are called vertices and edges, respectively.
The set of vertices of $\interval^n$ is denoted by $\vertices(\interval^n)$.
Given any face $F$, its \textbf{initial vertex|textbf}\index{cube!initial vertex} is obtained by setting all free coordinates to $0$, and its \textbf{terminal vertex|textbf}\index{cube!terminal vertex} is obtained by setting all free coordinates to $1$.
We note that $F$ is determined completely by its initial and terminal vertices: $F_0$ is the set of coordinates that are $0$ for both vertices, $F_1$ is the set of coordinates that are $1$ for both vertices, and $F_{01}$ is the set of coordinates that disagree.

\begin{example}
	The partition $(\{1,4,6\},\{3,5\}, \{2\})$ corresponds to the $2$-dimensional face of $\interval^6$ determined by $x_1=x_4=x_6=0$, $x_2=1$, and with $x_3$ and $x_5$ a free to vary in $\interval = [0,1]$.
	Its initial vertex is $(0,1,0,0,0,0)$ and its terminal vertex is $(0,1,1,0,1,0)$.
\end{example}

For $\varepsilon \in \{0, 1\}$ and $i \in \overline{n}$, we define maps $\delta_i^\varepsilon \colon \interval^{n-1} \to \interval^{n}$ by
\begin{align*}
	\delta_i^\varepsilon(x_1, \dots, x_{n-1}) & = (x_1, \dots, x_{i-1}, \varepsilon, x_i, \dots, x_{n-1}).
\end{align*}
Any composition of these is referred to as a \textbf{face inclusion map}.\index{cube!face inclusion}

We also have projection maps $\pi_i \colon \interval^n \to \interval^{n-1}$ such that
\[\pi_i(x_1, \ldots, x_n) = (x_1, \ldots, \hat x_i, \ldots, x_n),\]
with $\hat x_i$ as usual denoting the omission of the $x_i$ term.
Analogous to the face and degeneracy identities for simplicial sets, these operators satisfy the following relations \cite[Section 4]{GrMa03} (it is also easy and illuminating to work these out on one's own):\index{cube!face/degeneracy identities}

\[
\begin{array}{rlc}
	\delta_j^\eta \delta_i^\varepsilon &= \delta_{i+1}^\varepsilon \delta_j^\eta, &j\leq i,\\
	\pi_i \pi_j &= \pi_j \pi_{i+1}, & j \leq i,\\
	\pi_j \delta^\varepsilon_i &=
	\begin{cases}
		\delta_{i-1}^\varepsilon \pi_j, \\
		\id,\\
		\delta_i^\varepsilon \pi_{j-1},
	\end{cases}
	&\begin{array}{lll}j<i, \\ j=i, \\ j>i. \end{array}
\end{array}
\]

For $v \in \vertices(\interval^n)$ all coordinates are bound -- that is, $v_{01} = \emptyset$.
Thus $v$ is determined by the partition of $\overline n$ into $v_0$ and $v_1$, so
we have a bijection from the set of vertices of $\interval^n$ to the power set $\mathcal P(\overline n)$ of $\overline n$, sending $v$ to $v_1$.
The inclusion relation in the power set induces a poset structure on $\vertices(\interval^n)$ given explicitly by
\begin{equation*}
	v = (\epsilon_1, \dots, \epsilon_n) \leq w = (\eta_1, \dots, \eta_n) \iff \forall i,\ \epsilon_i \leq \eta_i.
\end{equation*}
We will freely use the identification of these two posets, thinking of $\mathcal P(\overline n)$ as a combinatorial model for $\interval^n$ in the same way that one identifies the totally ordered set $[0,\dots, n]$ with the $n$-simplex $\Delta^n$ in the simplicial setting.
Note that face embedding maps induce order-preserving maps at the level of vertices.

An \textbf{interval subposet}\index{cube!interval poset|textbf} of $\mathcal P(\overline n)$ is one of the form $[v, w] = \{u \in \mathcal P(\overline n)\ |\ v \leq u \leq w\}$ for a pair of vertices $v \leq w$.
This is precisely the set of vertices of the unique face $F$ of $\interval^n$ with $v$ as its initial vertex and $w$ as its terminal vertex.
This association determines a canonical bijection between faces of $\interval^n$ and such subposets; explicitly, to $[v, w]$ we associate the face $F$ defined by $F_\varepsilon = \{i \in \overline{n}\ |\ v_i = w_i = \varepsilon\}$ for $\varepsilon \in \{0, 1\}$.

The posets $\{\mathcal P(\overline n)\}_{n \geq 1}$ play the role for cubical complexes that finite totally ordered sets play for simplicial complexes.
Recall for comparison that an abstract ordered simplicial complex can be defined as a pair $(\vertices(X), X)$, where $X$ is a collection of subsets of the set $\vertices(X)$, called simplices, such that every singleton subset (vertex) is in $X$, every subset of a simplex (face) is again in $X$, and each simplex is equipped with a total order compatible with its faces.
We have the following cubical analogue.

\begin{definition}\label{D:cubical}\index{cubical complex|textbf}
	A \textbf{cubical complex} $X$ is a collection $\{ \sigma \}$ of finite non-empty subsets of a set $\vertices(X)$, together with, for each $\sigma \in X$, a bijection $\iota_\sigma \colon \sigma \to \mathcal P(\overline n)$ for some $n$, such that:
	\begin{enumerate}
		\item For all $v \in \vertices(X)$, $\{v\} \in X$,
		\item For all $\sigma \in X$ and all $[u,w] \subset \mathcal P(\overline n)$ the set $\rho = \iota_\sigma^{-1}([u,w])$ is in $X$ and the following diagram commutes, with the two diagonal maps being order preserving:
		\begin{equation*}
			\begin{tikzcd} [row sep = tiny, column sep = small]
				\sigma \arrow[rr, "\iota_\sigma"] && \mathcal P(\overline n) \\
				& [-5pt] {[}u,w{]} \arrow[ur, hook] & \\
				\rho \arrow[uu, hook] \arrow[rr, "\iota_\rho"'] && \mathcal P(\overline m).
				\arrow[ul, "\cong"'] \arrow[uu, dashed]
			\end{tikzcd}
		\end{equation*}
	\end{enumerate}
	We refer to an element $\sigma \in X$ as a \textbf{cube} or \textbf{face}\index{cubical complex!face} of $X$, refer to $\iota_\sigma \colon \sigma \to \mathcal P(\overline{n})$ as its \textbf{characteristic map},\index{cubical complex!characteristic map}
	and refer to $n$ as its \textbf{dimension}.
	If $\rho \subseteq \sigma \in X$, we say that $\rho$ is a \textbf{face} of $\sigma$ in $X$, and this determines a poset structure on $X$ itself.
	We identify elements in $\vertices(X)$ with the singleton subsets in $X$, referring to them as vertices.
\end{definition}

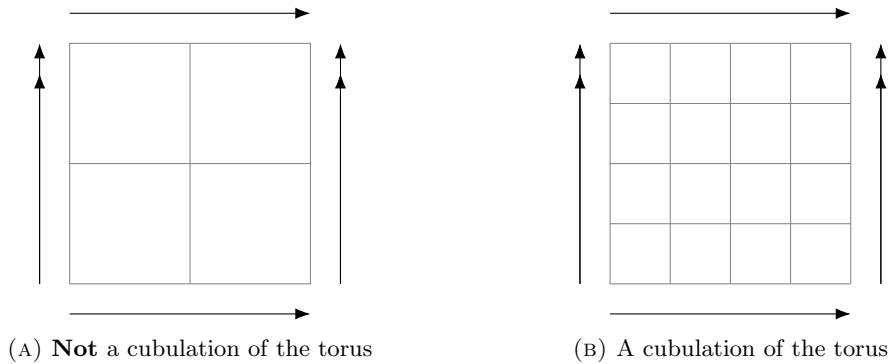
\begin{figure}
	\newcommand*{\xMin}{0}%
\newcommand*{\xMax}{4}%
\newcommand*{\yMin}{0}%
\newcommand*{\yMax}{4}%

\begin{subfigure}{.4\textwidth}
	\centering
	\begin{tikzpicture}[scale=.8]
		\draw[-{Latex[length=2mm]}] (-.5,\yMin)--(-.5,\yMax);
		\draw[-{Latex[length=2mm]}] (-.5,\yMin)--(-.5,\yMax-.5);
		\draw[-{Latex[length=2mm]}] (4.5,\yMin)--(4.5,\yMax);
		\draw[-{Latex[length=2mm]}] (4.5,\yMin)--(4.5,\yMax-.5);

		\draw[-{Latex[length=2mm]}] (\xMin, -.5)--(\xMax, -.5);
		\draw[-{Latex[length=2mm]}] (\xMin, 4.5)--(\xMax, 4.5);

		\draw [very thin,gray] (\xMin, \yMin) -- (\xMin, \yMax) -- (\xMax, \yMax) -- (\xMax, \yMin) -- (\xMin, \yMin);

		\draw [very thin,gray] (0.5*\xMax, \yMin) -- (0.5*\xMax, \yMax);

		\draw [very thin,gray] (\xMin, 0.5*\yMax) -- (\xMax, 0.5*\yMax);
	\end{tikzpicture}
	\caption{\textbf{Not} a cubulation of the torus}
\end{subfigure}\qquad
\begin{subfigure}{.4\textwidth}
	\centering
	\begin{tikzpicture}[scale=.8]
		\draw[-{Latex[length=2mm]}] (-.5,\yMin)--(-.5,\yMax);
		\draw[-{Latex[length=2mm]}] (-.5,\yMin)--(-.5,\yMax-.5);
		\draw[-{Latex[length=2mm]}] (4.5,\yMin)--(4.5,\yMax);
		\draw[-{Latex[length=2mm]}] (4.5,\yMin)--(4.5,\yMax-.5);

		\draw[-{Latex[length=2mm]}] (\xMin, -.5)--(\xMax, -.5);
		\draw[-{Latex[length=2mm]}] (\xMin, 4.5)--(\xMax, 4.5);

		\foreach \i in {0,...,4} {
			\foreach \j in {0,...,3} {
				\draw [-{Latex[length=1mm]}, very thin, gray] (\i,\j) -- (\i,\j+1);
				\draw [-{Latex[length=1mm]}, very thin, gray] (\j,\i) -- (\j+1,\i);
			}
		}
	\end{tikzpicture}
%
%
	\caption{A cubulation of the torus}
\end{subfigure}
	\caption{The first cellular decomposition of a torus pictured above does not represent the geometric realization of a cubical complex, as each square has the same set of vertices.
		On the right, each square has been coherently identified with the standard square and has a unique set of vertices.
		Therefore, (B) depicts a cubical structure on the torus.}
	\label{F: cubical structure}
\end{figure}

In the definition, the role of the characteristic maps $\iota_\sigma \colon \sigma \to \mathcal P(\overline{n})$ is to endow each cube $\sigma$ with a poset structure on its vertices that mirrors the poset structure on the topological cube $\interval^n$ that we previously established by identifying the vertices of $\interval^n$ with $\mathcal P(\overline n)$.
The first condition then says that every vertex is a cube of the cubical complex, while the second says that every face of a cube is also a cube in the complex with its compatible subposet structure.
We also note that the definition guarantees that no two cubes of a cubical complex can share the same set of vertices, i.e.\ a cube is completely determined by its set of vertices; cf.\ \cref{F: cubical structure}.

In analogy with the usual terminology in the simplicial setting, one could call these ``ordered cubical complexes," but we will not need the unordered version.
Our definition sits between cubical sets \cite{jardine2002cubical} and cellular subsets of the cubical lattice of $\R^\infty$ \cite{kaczynski2006computational}, analogously to the way that abstract ordered simplicial complexes sit between simplicial sets and simplicial complexes.

\begin{remark}
	In \cite{FMS-flows} we assumed $\vertices(X)$ to be a poset, but this extra condition was not used and is not necessary for what follows.
	It does, however, sometimes occur naturally, such as when a cubical complex is a subcomplex of a lattice, as will be the case in our construction below of cubulations of smooth manifolds.
\end{remark}

\noindent Given two cubical complexes $X$ and $Y$, their \textbf{product}\index{cubical complex!product of} is the cubical complex
\[
X \times Y = \{ \sigma_X \times \sigma_Y \mid \sigma_X \in X,\ \sigma_Y \in Y \},
\]
with vertex set $\vertices(X \times Y)$ and characteristic maps defined by
\[
\iota_{\sigma_X \times \sigma_Y} = \iota_{\sigma_X} \times \iota_{\sigma_Y} \colon \sigma_X \times \sigma_Y \to \mathcal{P}(\overline{n}) \times \mathcal{P}(\overline{m}) \cong \mathcal{P}(\overline{n + m}),
\]
where
\[
\mathcal{P}(\overline{n}) \times \mathcal{P}(\overline{m}) \cong \mathcal{P}(\overline{n + m}),
\]
is the canonical order-preserving bijection, which is compatible with the identification
\[
\vertices(\interval^n) \times \vertices(\interval^m) \cong \vertices(\interval^{n + m}).
\]

\subsection{Cubulations}

Let ${\tt Top}$ denote the category of topological spaces and continuous maps, and let ${\tt Cube}$\index{$Cu$@$\mathtt{Cube}$} be the subcategory whose objects are the $n$-cubes, identified with $\interval^n$, and whose morphisms are face inclusions.
The poset of faces of a cubical complex $X$ also determines a category, and the characteristic maps of $X$ determine a functor from this poset category to $\mathtt{Cube}$.
We define the \textbf{geometric realization}\index{cubical complex!geometric realization|textbf} $|X|$ of $X$ as the colimit of this functor; in other words, we glue topological cubes together according to the combinatorial data of the cubical complex $X$ in the evident way.
A \textbf{cubical structure} or \textbf{cubulation}\index{cubulation|textbf} on a space $S$ is a homeomorphism $h \colon |X| \to S$ for some cubical complex $X$.

A smooth cubulation is one for which the restriction of $h \colon |X| \to S$ to each cube is a smooth map of manifolds with corners.
For a given face of $F$ of $X$, we will also refer to the corresponding map $\interval^n \to S$ as the characteristic map of the face.
Smooth cubulations exist for any smooth manifold, as in the following construction of \cite{ShSh92}.\index{cubulation!smooth cubulations exist}

Start with a smooth triangulation (see for example \cite[Theorem 10.6]{Munk66} for the existence of such).
Consider the cell complex that is dual to its barycentric subdivision.
Intersecting those dual cells with each simplex in the triangulation provides a subdivision of the simplex into cells that are linearly isomorphic to cubes.
Moreover, starting with an ordered triangulation -- obtained for example by taking a barycentric subdivision -- such a cubical decomposition embeds cellularly into the cubical lattice of $\R^\infty$, and thus it is the geometric realization of a cubical complex.
See \cite{ShSh92} for details.

Using the product of cubical complexes, it follows immediately that the product of cubulated manifolds is again a cubulated manifold.

\subsection{Cubical chains and cochains}\label{S: cubical cochains}

We can also define an ``algebraic realization" for a cubical complex in analogy to its geometric realization.
Let $K_*(\interval^1)$ be the usual cellular chain complex of the interval with integer coefficients.
Explicitly, $K_0(\interval^1)$ is generated by the vertices, which we write in this context as $[\underline{0}]$ and $[\underline{1}]$, and $K_1(\interval^1)$ is generated by the unique 1-dimensional face, denoted $[\underline{0},\underline{1}]$ in the interval subposet notation.
The boundary map is $\bd [\underline{0},\underline{1}] = [\underline{1}]-[\underline{0}]$.

Let $K_*(\interval^n) = K_*(\interval^1)^{ \otimes n}$, with differential defined by the graded Leibniz rule.
Given a face inclusion $\delta_i^{\varepsilon} \colon \interval^n \to \interval^{n+1}$ the natural chain map $K_*(\delta_i^{\varepsilon}) \colon K_*(\interval^1)^{ \otimes n} \to K_*(\interval^1)^{ \otimes n+1}$ is defined on basis elements by
\begin{equation*}
	x_1 \otimes \cdots \otimes x_n \mapsto
	x_1 \otimes \cdots \otimes [\underline{\varepsilon}] \otimes \cdots \otimes x_n.
\end{equation*}
Regarding a cubical complex $X$ as a functor to $\mathtt{Cube}$, we can compose it with the chain functor above to obtain a functor to chain complexes.
The complex of \textbf{cubical chains}\index{cubical chain|textbf} of $X$, denoted $K_*(X)$,\index{$K_*(X)$} is defined to be the colimit of this composition.
As one would expect, in each degree it is a free abelian group generated by the cubes of that dimension, and its boundary homomorphism sends the
generator associated to a cube to a sum of generators associated to its codimension-one faces with appropriate signs.
Furtheremore, the functor of cubical chains is strongly monoidal in the sense that $K_*(X \times Y)$ is canonically isomorphic to $K_*(X) \otimes K_*(Y)$ for any pair of cubical complexes.
We will write $H_*(X)$\index{$H_*(X)$} for $H_*(K_*(X))$, the cubical homology of $X$.\index{cubical homology}

By abuse, we will use the same notation and terminology for an element in $X$, its geometric realization in $|X|$,
and the corresponding basis element in $K_*(X)$.
Most commonly we will write $F$ and refer simply to a ``face of $X$.''

We note that for each $\interval^n$ we have the ordered set $\{\e_1, \dots, \e_n\}$ where $\e_i = \frac{\bd\ }{\bd x_i}$.
For any face $F$ of $\interval^n$, the ordered subset $\beta_F = \{\e_i\ |\ i \in F_{01}\}$ defines the \textbf{canonical orientation} of $F$.\index{cubical complex!face!canonical orientation of}
In forming the cubical complex $X$, these orientations are preserved, and so each face of $X$ carries an orientation.
These orientations are compatible with our standard generators of $K_*(X)$ in the sense that if we identify $[0,1]^{ \otimes k}$ with $\interval^k$ with its standard orientation then, in the boundary formula, $k-1$ faces appear with sign $1$ or $-1$ according to whether or not their standard orientations agree with the boundary orientation of $\interval^k$ as a manifold with corners.

The \textbf{cubical cochain complex}\index{cubical cochain|textbf} of $X$ (with $\Z$ coefficients) is the chain complex
\[
K^*(X) = \Hom_\Z(K_*(X), \Z).
\]\index{$K^*(X)$}
If $F$ is a face of $X$, and correspondingly a generator of $K_*(X)$, then we will write $F^*$ for the dual, i.e.\ the element of $K^*(X)$ such that $F^*(F) = 1$ and $F^*(E) = 0$ for all other faces $E\neq F$ of $X$.
We will use the convention as in \cite[Section 42]{Mun84} that
$$(dF^*)(\xi) = F^*(\bd \xi).$$
As discussed in \cite[Section 42]{Mun84} in the simplicial context, these formulas are sufficient to determine the coboundary operator on all of $K^*(X)$, thinking of elements of $K^*(X)$ as possibly infinite formal sums $\sum n_i F_i^*$.
We will write $H^*(X)$\index{$H^*(X)$} for $H^*(K^*(X))$, the cubical cohomology of $X$.\index{cubical cohomology}

\subsection{Comparing cubical and geometric homology}\label{S: cubical and geometric homology}

Suppose $h \colon |X| \to M$ is a smooth cubulation.
As the cubes of $X$ are compact oriented manifolds with corners, the composition of the inclusion of a cube into $X$ with the map $h$ gives an element of $PC_*^\Gamma(M)$ and hence an element of $C_*^\Gamma(M)$.
Furthermore, the boundary formula for cubes in $K_*(X)$ agrees with the geometric boundary formula.
So we obtain a chain map, in fact a subcomplex inclusion, $\mc J \colon K_*(X) \into C^\Gamma_*(M)$.
We show in \cref{T: cubical homology iso} that the homology of $K_*(X)$ is both the standard homology of $M$ and the geometric homology $M$.

Recall from \cref{S: homology is homology} the normalized singular cubical chain complex of $M$, denoted $NK_*(M)$, as well as the subcomplex $NK^{sm}_*(M)$ generated by smooth cubes.
The homology groups of each are isomorphic to the standard simplicial singular homology of $M$, for example by our work in \cref{S: homology is homology,S: smooth cubes}.
Furthermore, by \cref{T: hom iso map,P: singular smooth cubes}, there are quasi-isomorphisms $NK_*(M) \xleftarrow{\psi} NK^{sm}_*(M) \xr{\phi} C^\Gamma_*(M)$, the latter induced by observing that smooth singular cubes are elements of $PC^\Gamma_*(M)$ and that degenerate cubes are elements of $Q_*(M)$.

There is also a map $\eta: K_*(X) \to NK^{sm}_*(M)$ that takes each cube into its embedding into $M$.
Then $\mc J = \phi\eta \colon K_*(X) \to C^\Gamma_*(M)$.

\begin{theorem}\label{T: cubical homology iso}\index{cubical homology!is singular cubical homology}
	The map $\eta: K_*(X) \to NK^{sm}_*(M)$ is a chain homotopy equivalence, and the map $\mc J \colon K_*(X) \to C^\Gamma_*(M)$ induces an isomorphism of homology groups $H_*(K_*(X)) \to H_*^\Gamma(M)$.
\end{theorem}

\begin{proof}

	As $\phi \colon NK^{sm}_*(M) \to C^\Gamma_*(M)$ is a quasi-isomorphism by \cref{T: hom iso map}, it suffices to show that $\eta$ is a chain homotopy equivalence.
	In fact, as $NK^{sm}_*(M)$ is freely generated by the nondegenerate smooth singular cubes, $K_*(X)$ and $NK^{sm}_*(M)$ are both free chain complexes bounded below, so it suffices by \cite[Theorem 46.2]{Mun84} to show that $\eta$ induces homology isomorphisms.

	The proof is analogous to the proof of \cite[Proposition V.8.3]{Dol72}, which provides an isomorphism between simplicial and singular homology.

	Consider have the diagram
	\[
	\begin{tikzcd}
		H_*(K_*(X)) \arrow[r, "\eta_*"] \arrow[d, "\cong"] & H_*(NK^{sm}_*(M)) \arrow[r, "\psi_*", "\cong"'] & H_*(NK_*(M)) \\
		H_*(CW_*(M)) \arrow[rru, "\Theta", "\cong"', out=0, in=200] & &
	\end{tikzcd}
	\]
	in which $CW_*(M)$ is the CW chain complex of $M$ corresponding to the CW complex structure given by the cubulation and $\Theta$ is the standard isomorphism between CW homology and singular homology as in Dold \cite[Proposition V.1.9]{Dol72}.
	The isomorphism in Dold is developed using simplicial singular homology, but as simplicial singular and cubical singular homology are isomorphic, the argument there goes through identically using singular cubes.
	The map on the left is an isomorphism at the chain level as there is an evident isomorphism in this case between the cubical chain complex and the CW chain complex that takes an embedding of a $k$-cube to the corresponding generator of $CW_k(M) = H_k(X^k, X^{k-1})$ (where we assume the expression on the right is singular cubical homology).
	As in Dold, the map $\Theta$ takes a class in $H_k(CW_*(M))$ represented by a $k$-cycle $z$ in $CW_k(M)$ to the class in $H_k(NK_*(M))$ represented by a singular (cubical) cycle in $NK_k(M)$ that represents the same class as $z$ in $H_k(X^k,X^{k-1})$.
	But all cycles in the image of the vertical map of the diagram are already represented by singular cycles in $NK^{sm}_k(M)$, so the diagram commutes, and it follows that $\eta_*$ is an isomorphism.
\end{proof}

As a corollary of the proof, we have the following useful result concerning the cohomology groups of the complexes, which we obtain using the following definitions:
\begin{align*}
	NK^*(M)& \defeq \Hom(NK_*(M),\Z)\\
	NK_{sm}^*(M)& \defeq \Hom(NK^{sm}_*(M),\Z).
\end{align*}

\begin{corollary}\index{cubical cohomology!is cubical singular cohomology}
	The following maps on cohomology induced by restrictions are isomorphisms: $$H^*(NK^*(M)) \to H^*(NK_{sm}^*(M)) \to H^*(X).$$
\end{corollary}

\begin{proof}
	This follows from basic homological algebra \cite[Theorem 45.5]{Mun84} as $NK_*(M)$, $NK^{sm}_*(M)$, and $K_*(X)$ are all free chain complexes, observing that even though $NK_*(M)$ and $NK^{sm}_*(M)$ are defined by taking quotients of the groups of singular cubical chains $SK_*(M)$ and $SK^{sm}_*(M)$ by the subgroups of degenerate cubes, the degenerate cubes correspond to generators of $SK_*(M)$ and $SK^{sm}_*(M)$, and so each $NK_i(M)$ and $NK^{sm}_i(M)$ is freely generated by the nondegenerate, respectively nondegenerate and smooth, singular $i$-cubes.
\end{proof}

\subsection{Cubically transverse geometric cohomology}\label{S: transverse cochains}

In this section we consider the cochains on $M$ represented by maps $W \to M$ that are transverse to a given cubulation of $M$.

\begin{definition}\label{D: trans cube}
	Let $M$ be equipped with a smooth cubulation $|X| \to M$.
	We say that $r_W \colon W \to M$ is \textbf{transverse}\index{cubulation!transversality to|textbf} to $X$ if $r_W \colon W \to M$ is transverse to each characteristic map of the cubulation.
	In particular, this implies by \cref{L: simple trans} that each induced $\bd^kW \to M$ is naively transverse to each face of the cubulation.
	If $r_W \colon W \to M$ is transverse to $X$ then the same is true for any $r_V \colon V \to M$ isomorphic to $r_W \colon W \to M$, and so we can define $PC^*_{\Gamma \pf X}(M)$\index{$PC^*_{\Gamma \pf X}(M)$} to be the subset of $PC^*_{\Gamma}(M)$ consisting of those precochains with reference maps transverse to $X$.

	We let $Q^*_{\Gamma \pf X}(M) = PC_{\Gamma \pf X}^*(M) \cap Q^*(M)$ and note that the equivalence relation of \cref{L: cancel Q} descends to an equivalence relation on $PC_{\Gamma \pf X}^*(M)$ such that $V\sim W$ if and only if $V \sqcup -W \in Q^*_{\Gamma \pf X}(M)$.
	The \textbf{geometric cochains of $M$ transverse to $X$},\index{geometric cochain!transverse to cubulation|textbf} denoted $C_{\Gamma \pf X}^*(M)$,\index{$C_{\Gamma \pf X}^*(M)$} are the equivalence classes in $PC_{\Gamma \pf X}^*(M)$.
	The set $C_{\Gamma \pf X}^*(M)$ is a chain complex under the operation $\sqcup$ and with boundary map $\bd$.
	If $V, W \in PC_{\Gamma \pf X}^*(M)$, then $V \sqcup -W \in Q^*_{\Gamma \pf X}(M)$ if and only if as elements of $PC_\Gamma^*(M)$ we have $V \sqcup -W \in Q^*(M)$, so $V$ and $W$ represent the same element of $C_{\Gamma \pf X}^*(M)$ if and only if they represent the same element of $C^*_\Gamma(M)$.
	Thus $C_{\Gamma \pf X}^*(M)$ is a subcomplex of $C_\Gamma^*(M)$.
	The \textbf{geometric cohomology transverse to $X$}\index{geometric cohomology!transverse to cubulation|textbf} is $H_{\Gamma \pf X}^*(M) \defeq H^*(C_{\Gamma \pf X}^*(M))$.\index{$H_{\Gamma \pf X}^*(M)$}

	When the specific cubulation $X$ is understood, we sometimes simplify the notation to $PC_{\Gamma\pf}^*(M)$, $C_{\Gamma\pf}^*(M)$, and $H_{\Gamma\pf}^*(M)$.
\end{definition}

The proof of \cref{L: co/chains well defined} continues to hold for transverse cochains, and so $r_W \colon W \to M$ represents $0$ in $C^*_{\Gamma \pf X}(M)$ if and only if it is in $Q^*_{\Gamma \pf X}(M)$.
Therefore, the evident map $C^*_{\Gamma \pf X}(M) \to C^*_\Gamma(M)$, which takes the element of $C^*_{\Gamma \pf X}(M)$ represented by $r_W \colon W \to M$ to the element of $C^*_\Gamma(M)$ represented by the same map, is a monomorphism of chain complexes, for such an $r_W$ is transverse to $X$ by definition and if it is also in $Q^*(M)$ then it is in $Q^*_{\Gamma \pf X}(M)$.
Thus we will think of $C^*_{\Gamma \pf X}(M)$ as a subcomplex of $C^*_\Gamma(M)$.
A key technical result, which will take the remainder of this section to prove,
is that this inclusion induces a cohomology isomorphism.
In other words, the cochains that are transverse to $X$ are sufficient to compute the cohomology of $M$.

\begin{theorem}\label{T: transverse complex}\index{geometric cohomology!transverse to cubulation!is geometric cohomology}
	The inclusion $C^*_{\Gamma \pf X}(M) \into C^*_\Gamma(M)$ is a quasi-isomorphism.
\end{theorem}

To show that the inclusion $C^*_{\Gamma \pf X}(M) \into C^*_\Gamma(M)$ is a quasi-isomorphism, it will be necessary to consider the following scenario.
Suppose we have a map $r_V \colon V \to M$ with $V$ a manifold with corners and $M$ a manifold with a cubulation.
Let $\bd V = W$, and suppose $W$ is already transverse to the cubulation.
We will construct a homotopy $h \colon V \times I \to M$ such that $g(-,0) = r_V$, $g(-,1)$ is transverse to the cubulation, and the restriction of $h$ to $W \times I$ is transverse to the cubulation. However, as noted at the end of \cref{S: covariant functoriality}, such a homotopy might not preserve cohomology classes, as we will need below. So we must instead use the universal homotopies of \cref{D: universal homotopy,P: universal homotopy}.

The technique for constructing such homotopies will be modeled on a variety of results in \cite{GuPo74}.
We use the Transversality Theorem and Transversality Homotopy Theorem of \cite[Section 2.3]{GuPo74} as stated.
However, for the Stability Theorem of \cite[Section 1.6]{GuPo74} we will provide details of the proof because the proof is only sketched in \cite{GuPo74} and we will need the result to be generalized in several ways.
Also, the Stability Theorem was not stated correct in early printings of \cite{GuPo74}, where the requirement that the submanifold of the target be a closed set was omitted\footnote{As stated in early printings of \cite{GuPo74}, the claim was that if $f \colon X \to Y$ is transverse to any submanifold $Z$ of $Y$ then this property is stable under small homotopies of $f$; more specifically that if $f_t:X \times I \to Y$ is a homotopy with $f_0$ transverse to $Z$ then there is an $\epsilon>0$ such that $f_t$ is transverse to $Z$ for all $t\in[0,\epsilon)$.
	Here is a counterexample:

	In the plane $\R^2$, let $Z = \{(x,y)|y = x^2, x\neq 0\}$ and consider maps $g_t: \R \to \R^2$ with
	$g_t(x) = (x,t^2+2t(x-t))$.
	For each fixed $t$, the image is the line given by $y-t^2 = 2t(x-t)$, which has slope $2t$ and passes through the point $(t,t^2)$.
	So the map $g_0$ embeds $\R$ as the x-axis, and as the image does not intersect $Z$, the map $g_0$ is transverse to $Z$.
	But for all $t\neq 0$, $g_t$ takes $\R$ to a line that is tangent to $Z$, and so $g_t$ is not transverse to $Z$ for $t\neq 0$, violating the Stability Theorem as stated on page 35 of \cite{GuPo74}.

	The error in the original proof comes from considering only what happens in neighborhoods of points $x$ such that $f(x) \in Z$ but not points $x$ with $f(x)\notin Z$.
	As we can see, the claim breaks down when $f(x)\notin Z$ but every neighborhood of $(x,0)$ in $X \times I$ has a point with image in $Z$.
	However, this can be avoided if $Z$ is a closed set in $Y$, as is the case for the statement of the theorem in later printings of \cite{GuPo74}.}.
The needed versions of these results is established in the following proposition:

\begin{proposition}\label{P: ball stability}\index{transversality!Transversality Theorem!universal with respect to cubulation}
	Suppose $r_V \colon V \to M$ is a proper map from a manifold with corners to a cubulated manifold without boundary.
	Then there is a proper homotopy $H \colon M \times I \to M$ such that $H(-,0) = \id_M$ and $H(-,1) r_V \colon V \to M$ is transverse to the cubulation.
	In other words, there is a proper universal homotopy from $r_V$ to a map that is transverse to the cubulation.

	Furthermore, given another proper map $r_W \colon W \to M$ from a manifold with corners that is transverse to the cubulation, we can choose the homotopy $H$ above so that also the resulting universal proper homotopy of $W$ given by $W \times I \xr{r_W \times \id_I} M \times I \xr{H} M$ is transverse to the cubulation.
\end{proposition}

Before proving the proposition, which is somewhat technical, we use it to prove \cref{T: transverse complex}, which states that $H^*(C^*_{\Gamma \pf X}(M)) \to H^*(C_\Gamma^*(M))$ is an isomorphism.

\begin{proof}[Proof of \cref{T: transverse complex}]
	The idea of the argument that $H^*(C^*_{\Gamma \pf X}(M)) \to H^*(C_\Gamma^*(M))$ is a surjection is contained already in the proof of \cite[Lemma 15]{Lipy14}, which involves constructing a homotopy to move a cycle into transverse position.
	We elaborate upon that argument.

	Suppose $\uV \in C_\Gamma^*(M)$ is a cocycle represented by $r_V \colon V \to M$.
	By \cref{P: ball stability}, there is a proper universal homotopy $h \colon V \times I \to M$ from $r_V$ to $H(-,1)r_{V} \colon V \to M$ and this latter map is transverse to the cubulation.
	Let us call the transverse map $r_{V'} \colon V' \to M$ with $V' \cong V$.
	By \cref{C: homotopy}, $r_V$ and $r_{V'}$ represent the same cohomology class in $H^*_{\Gamma}(M)$, but the class represented by $r_{V'}$ is in the image of $H^*(C^*_{\Gamma \pf X}(M))$.

	For injectivity, suppose $W \in PC^*_{\Gamma \pf X}(M)$ is transverse to the cubulation and represents zero in $H^*(C_\Gamma^*(M))$.
	Then by definition there is a $V \in PC^*_\Gamma(M)$ with $\bd V \sqcup -W \in Q^*(M)$.
	By \cref{P: ball stability}, there is a proper homotopy $H \colon M \times I \to M$ such that $H(-,0) = \id_M$ and both $H(-,1)r_V$ and $H\circ (r_W \times \id_I)$ are transverse to the cubulation; co-orient $H$ and $H(-,1)$ as in \cref{D: homotopy co-orientation} based on the tautological co-orientation of $H(-,0)$.
	Let $V'$ be the precochain $V \xr{H(-,1)r_V} M$, let $W'$ be the precochain $W \xr{H(-,1)r_W} M$, and let $Y$ be the precochain $W \times I \xr{H\circ (r_W \times \id_I)} M$, co-orienting $r_W \times \id_I$ using the naive pre-homotopy co-orientation convention, \cref{homotopy product co-orientation convention}.
	Let $Z$ be the precochain $V' \sqcup -Y$, which is transverse to the cubulation.
	We note that $\bd V' \sqcup -W'$ is the image of $\bd V \sqcup -W$ after composing with $H(-,1)$, and so it is in $Q^*(M)$ by \cref{L: Q preservation}.
	Also by \cref{C: universal homotopy boundary co-orientation,L: co-oriented homotopy}, we have $\bd Y = W' \sqcup -W \sqcup B$, where $B$ is the map $\bd W \times I \xr{H\circ (r_{\bd W} \times \id_I)} M$ co-oriented as a homotopy from $-\bd W$ to $-\bd W'$.
	Note that $B$ is also transverse to the cubulation as part of the boundary of $Y$.
	We now compute
	\begin{align*}
		\bd Z \sqcup -W &= \bd V' \sqcup -\bd Y \sqcup -W\\
		&= \bd V' \sqcup -(W' \sqcup -W \sqcup B) \sqcup -W\\
		&= \bd V' \sqcup -W' \sqcup W \sqcup -W \sqcup -B.
	\end{align*}
	We have already noted $\bd V' \sqcup -W' \in Q^*(M)$ and $W \sqcup -W$ is trivial.
	Since $W$ represents a cycle, $\bd W \in Q^*(M)$ and hence $B \in Q^*(M)$ by \cref{L: dessicated homotopy}.
	So $\bd Z \sqcup -W \in Q^*(M)$ and $Z$ is transverse to the cubulation, which tells us that $W$ represents $0$ in $H^*_{\Gamma \pf X}(M)$.
\end{proof}

It remains to prove \cref{P: ball stability}, which will require the following technical lemma that is also useful below in the proof of \cref{T: intersection qi}.

\begin{lemma}\label{L: minimizer}\index{function subordinate to a cover with constants}
	Let $M$ be a manifold without boundary, and let $\mc U = \{U_j\}$ be a locally finite open cover such that each $\bar U_j$ is compact.
	Suppose given $\varepsilon_j>0$ for each $j$.
	Then there exists a smooth function $\phi \colon M \to \R$ such that $0<\phi(x)<\varepsilon_j$ if $x \in \bar U_j$.
\end{lemma}

\begin{proof}
	Let $\eta_j = \min\{\varepsilon_k \mid \bar U_j \cap \bar U_k\neq \emptyset\}$.
	By the local finiteness and compactness conditions, $\{k \mid \bar U_j \cap \bar U_k\neq \emptyset\}$ is a finite set and so the $\eta_j$ are well defined.
	Let $\{\psi_j\}$ be a partition of unity subordinate to $\mc U$ and let $\phi_1 = \sum \eta_j\psi_j$.
	For any $x \in M$, this sum is positive.
	If $x \in \bar U_j$ then $\phi_1(x) = \sum_{\{k \mid \bar U_j \cap \bar U_k\neq \emptyset\}} \eta_k\psi_k$.
	But for any such $k$, we have $\eta_k \leq \varepsilon_j$.
	Thus $\phi_1(x) \leq \varepsilon_j$.
	Now take $\phi = \frac{1}{2}\phi_1$.
\end{proof}

We can now prove \cref{P: ball stability}.
In the following $D^N$ is the open unit ball in $\R^N$ and, more generally, $D^N_r$ is the open ball of radius $r$.

\begin{proof}[Proof of \cref{P: ball stability}]
	We begin with the case that $W$ is compact, and then we will show how to use the arguments of the compact case to obtain the general case.
	We first construct a map $F \colon M \times D^N \to M$, for some $N$, such that

	\begin{enumerate}
		\item $F(-,0) = \id \colon M \to M$,
		\item for almost all $s \in D^N$ the composition $V \xr{r_V} M \xr{F(-,s)}M$ is transverse to the cubulation,
		\item there is a ball neighborhood $D_r^N$ of $0$ in $D^N$ such that for all $s \in D_r^N$ the composition $W \xr{r_W} M \xr{F(-,s)}M$ is transverse to the cubulation.
	\end{enumerate}

	This will suffice to provide the transversality required in the $W$ compact case as then we can let $s_0$ be any point in $D_r^N$ such that the composition $V \xr{r_V} M \xr{F(-,s_0)}M$ is transverse to the cubulation and define $H(-,t) = F(-,ts_0)$.
	Then $H(-,0) = \id_M$ since $F(-,0) = \id$.
	We will have $H(-,1) r_V \colon V \to M$ transverse to the cubulation by our choice of $s_0$.
	Finally, as $ts_0 \in D_r^N$ for all $t \in I$, each $F(-,ts_0)r_W$ is transverse to the cubulation, which then implies that $H \circ (r_{W} \times \id)$ is transverse to it as well.
	This does not provide the properness of $H$, but we consider that below.

	The construction of $F$ is a small variation of the construction in the Transversality Homotopy Theorem of \cite[Section 2.3]{GuPo74}:
	Let $M_\epsilon$ be an $\epsilon$-neighborhood of $M$ in some $\R^N$ in the sense of the $\epsilon$-Neighborhood Theorem of \cite[Section 2.3]{GuPo74}; in particular,
	$M_\epsilon$ is an $\epsilon$-neighborhood of a proper embedding of $M$ into $\R^N$ that possesses a submersion $\pi \colon M_\epsilon \to M$.
	We may also assume $\epsilon$ is a bounded function.
	If $M$ is not compact, then $\epsilon$ is a smooth bounded positive function of $M$ and $M_\epsilon = \{z \in \R^N \mid |z-y|<\epsilon(y) \text{ for some }y \in M\}$.
	Let $f \colon M \times D^N \to M_\epsilon$ be given by $f(y, s) = y + \eta(y) s$, where $\eta \colon M \to \R$ is a smooth function such that $0 < \eta(y) < \epsilon(y)$ for all $y \in M$.
	As $\eta > 0$, this is clearly a submersion (onto its image) at all points.
	We let $F \colon M \times D^N \to M$ be the composition $M \times D^N \xr{f}M_\epsilon \xr{\pi}M$.
	Furthermore, the map $\ms F \colon V \times D^N \to M$ given by the composition $$V \times D^N \xr{r_V \times \id} M \times D^N \xr{F} M$$ as well as all the restrictions $\ms F|_{S^k(V)}$
	are submersions.
	In particular, each $\ms F|_{S^k(V)}$ is transverse to any submanifold of $M$, so it follows by the Transversality Theorem of \cite[Section 2.3]{GuPo74} that for any fixed submanifold $Z$ of $M$, each $\ms F|_{S^k(V)}(-,s)$ is transverse to $Z$ for almost all $s \in D^N$.
	In particular, we may take $Z$ to be the interior of any cube $E$ (of any dimension) of the cubulation.
	There are countably many cubes in the cubulation of $M$ and finitely many manifolds $S^k(V)$.
	As the countable union of measure zero sets has measure zero, for almost all $s \in D^N$ we have for all $k$ that $\ms F|_{S^k(V)}(-,s) = F(-,s)r_V|_{S^k(V)}$ is transverse to all cubical faces.

	We now show how to ensure $H$ is proper.
	We consider $M$ endowed with the the distance function $|\cdot|$ obtained by restricting the standard distance function on $\R^N$ to $M$.
	The projection $p \colon M \times I \to M$ is proper, and by
	\cref{L: nearby proper homotopy}, there is a continuous function $\varepsilon \colon M \times I \to (0,\infty)$ such that if $H \colon M \times I \to M$ satisfies $|p(y,t)-H(y,t)| < \varepsilon(y,t)$ for all $(y,t) \in M \times I$, then $H$ is proper.
	By essentially the same computation performed in the proof of \cref{T: basic trans},
	$$|p(y,t)-H(y,t)| = |y-H(y,t)| = |y - \pi(y + \eta(y) ts_0)| < 2 \eta(y).$$
	So for $H$ to be proper, it suffices to have $\eta(y) < \varepsilon(y,t)/2$ for all $(y,t) \in M \times I$.
	But if this is not already the case, we can simply replace $\eta$ in the above definition of $f$ with an $\eta \colon M \to (0,\infty)$ such that $\eta(y) < \min\{\epsilon(y), \min\{\varepsilon(y,t)/2 \mid t \in I\}\}$ for all $y$, analogously to the construction in the proof of \cref{T: basic trans}.
	Such an $\eta$ can be found by considering a locally finite open cover $\mc U = \{U_j\}$ of $M$ such that each $\bar U_j$ is compact.
	Then $\min\{\epsilon(y) \mid y \in \bar U_j\}$ and $\min\{\varepsilon(y,t)/2 \mid (y,t) \in \bar U_j \times I\}$ are both defined and positive as $\bar U_j$ and $\bar U_j \times I$ are compact and $\epsilon$ and $\varepsilon$ are both positive continuous functions.
	Now choose $b_j > 0$ such that $b_j$ is less than the minimum of these two minima for each $j$, and then apply \cref{L: minimizer}.

	It remains in the compact $W$ setting to show that if we are given compact $W$ with $r_W \colon W \to M$ transverse to the cubulation then $F(-,s)r_W$ is transverse to the cubulation for all $s$ in some neighborhood $D_r^N$ of $0$ in $D^N$.
	It is here that we need to generalize the Stability Theorem of \cite[Section 1.6]{GuPo74}.
	As the Stability Theorem is not necessarily true when the manifolds involved are not closed submanifolds, compact, or controlled in some other way, it is more convenient here to work with the closed cubical faces of the cubulation and with $\bd^kW$ rather than $S^k(W)$.
	We recall that by \cref{L: simple trans}, to prove that two maps of manifolds with corners are transverse it is sufficient to show that their compositions with all pairs of boundary inclusions are naively transverse (see \cref{D: naive transversality}).

	Let $\Upsilon_k$ denote the composition $$\Upsilon_k \colon \bd^kW \times D^N \xr{i_{\bd^ kW} \times \id_{D^N}} W \times D^N \xr{r_W \times \id_{D^N} } M \times D^N \xr{F} M.$$
	We provide the details primarily for $\Upsilon_0$, with $\bd^0 W$ being $W$ itself, the other cases being analogous.

	Let $E$ be a (closed) face of the cubulation.
	Recall that $F(-,0) = \id_M$, so $\Upsilon_0(-,0) = r_W$ is transverse to the cubulation.
	Thus for any $x \in W$, either $r_W(x)\notin E$ or $r_W$ is (naively) transverse to $E$ at $r_W(x)$.
	In the former case, as $E$ is closed, there is an open neighborhood $A_x$ of $(x,0) \in W \times D^N$ such that $\Upsilon_0(A_x) \cap E = \emptyset$.
	Now suppose that $r_W(x) \in E$ and is transverse there.
	By appealing to charts, we can suppose without loss of generality (at least locally in neighborhoods of $x$ and $r_W(x)$) that $M = \R_j^m$ for some $j$ with $m = \dim(M)$ and $x = 0$ and that $E = \{(y_1,\ldots,y_m) \mid y_i\geq 0\text{ for } i \leq \dim E\text{ and } y_i = 0 \text{ for } i>\dim(E)\}$.
	The transversality assumption means that the composition of $D_xr_W \colon T_xW \to T_{r_W(x)}M$ with the projection to the last $m-\dim(E)$ coordinates is a linear surjection.
	As this is an open condition on the Jacobian matrix of $r_W$ at $x$, it follows again that there is an open neighborhood $A_x$ of $(x,0) \in W \times D^N$ such that for each $(x',s)$ in the neighborhood $\Upsilon_0(-,s)$ is transverse to $E$ at $x'$ (it is possible that $\Upsilon_0(x',s)$ no longer intersects $E$, but this is fine).
	Taking the union of the $A_x$ over all $x \in W$ gives a neighborhood $B_E$ of $W \times 0$ in $W \times D^N$, and by the Tube Lemma, as $W$ is compact there is a neighborhood of $W \times 0$ of the form $W \times U_E \subset B_E$.
	For each $s \in U_E$, we have $\Upsilon_0(-,s)$ transverse to $E$.
	Now let $D^N_{1/2}$ be the open ball of radius $1/2$ and $\bar D^N_{1/2}$ its closure.
	As $W \times \bar D^N_{1/2}$ is compact, its image under $\Upsilon_0$ can intersect only a finite number of faces of the cubulation of $M$; call this collection $\mc E$.
	Then let $U_0$ be the finite intersection $U_0 = D^N_{1/2} \cap \bigcap_{E\in\mc E} U_E$.
	Then $W \times U_0 \subset W \times D^N$ is a neighborhood of $W \times 0$ on which $\Upsilon_0(-,s)$ is transverse to every cubical face that its image intersects.
	Let $U_k$ be defined similarly for each $k\geq 1$ using $\bd^kW$ in place of $W$ and $\Upsilon_k$ in place of $\Upsilon_0$.
	As $W$ has finite depth, $U = \cap U_k$ is a neighborhood of $0$ in $D^N$.
	Let $r>0$ be such that $D^N_r \subset U$.
	Then for every $s \in D^N_r$ we have $\Upsilon_k(-,s) \colon \bd^kW \to M$ transverse to all $E$ for all $k$ as required.

	This completes the proof of the proposition for $W$ compact.

	Next suppose that $W$ is no longer necessarily compact.
	We show how to apply and extend the preceding arguments.
	We will define a new homotopy $\hat H \colon M \times I \to M$ with the desired properties of the proposition.

	To begin we construct $F \colon M \times I \to M$ exactly as above, as its definition did not depend on the compactness of $W$.
	Let $K \subset W$ be compact and $E$ be a closed cube.
	Taking the union of the $A_x$ as defined above over all $x \in K$ and intersecting with $K \times D^N$ gives an open neighborhood $B_E$ of $K \times 0$ in $K \times D^N$, such that
	$\Upsilon_0(-,s) \colon W \to M$ is transverse to $E$ at all $x \in K$.
	Furthermore, as $K \times \bar D^N_{1/2}$ is compact, its image under $\Upsilon_0$ can intersect only a finite number of faces of the cubulation of $M$, so again by applying the tube lemma and then intersecting tubular neighborhoods, we find an open ball $D_{r,K}^N \subset D^N$ centered at $0$ such that for all $x \in K$ and $s \in D_{r,K}^N$ we have $\Upsilon_0(-,s)$ transverse to the cubulation at $x$.
	As the maps $\bd^kW \into W$ for $k\geq 1$ are all proper, we can similarly find $D_{r,K}^N$ so that for all $s \in D_{r,K}^N$ and all $k\geq 0$, we have $\Upsilon_k(-,s) \colon \bd^kW \to M$ transverse to the cubulation at any $x \in \bd^{k}W$ whose image in $W$ is in $K$.

	Let $\{\mc U_j\}$ be a locally finite covering of $M$ such that each $\bar{\mc U_j}$ is compact.
	As $r_W$ is proper, each $r^{-1}_W(\bar {\mc U_j})$ is compact in $W$.
	Proceeding as just above with $r_W^{-1}(\bar U_j)$ in place of $K$, we can find for each $j$ an $\varepsilon_{j,0}$ with $0 < \varepsilon_{j,0} \leq 1$ so that for every $s \in D^N_{\varepsilon_{j,0}}$ we have $\Upsilon_0(-,s)$ transverse to all cubical faces of $M$ at every $x \in r^{-1}_W(\bar {\mc U_j})$.
	Analogously, we have $\varepsilon_{j,k}$ for all $k\geq 1$ using $(r_{\bd^k W})^{-1}(\bar {\mc U_j})$.
	Let $\varepsilon_j = \min\{\varepsilon_{j,k} \mid k\geq 0\}$.
	These minima exist as $W$ has finite depth.

	Now, using \cref{L: minimizer}, we choose a smooth function $\phi \colon M \to \R$ such that for all $y \in M$ we have $0<\phi(y)<\epsilon_j$ if $y \in \bar{\mc U_j}$.
	Let $M\times_\phi D^N = \{(y,s) \in M \times D^N \mid |s|<\phi(y)\}$.
	By our construction, $\Upsilon_k(-,s) \colon \bd^{k}W \to M$ is transverse to the cubulation at each $x$ such that $(x,s)\in(r_{\bd^{k}W} \times \id)^{-1}(M\times_\phi D^N) = \{(x,s) \in \bd^{k}W \times D^N \mid |s|<\phi(r_{\bd^kW}(x))\}$.
	Unfortunately, however, while $M\times_\phi D^N$ is a neighborhood of $M \times 0$ in $M \times D^N$, there is not necessarily a $U \in D^M$ so that $M \times U \subset M\times_\phi D^M$.
	Thus,
	we cannot construct $H$ from $F$ as above using a fixed $s_0$ as there may be no single $s_0\neq 0$ so that $W \times s_0 \subset M\times_\phi D^N$.

	To account for this, we modify our functions above as follows: Let $\hat f \colon M \times D^N \to M_\epsilon$ be given by $\hat f(y, s) = y + \phi(y) \eta(y) s$; as $\phi(y)\eta(y)>0$, this is again a submersion onto its image at all points.
	Let $\hat F \colon M \times D^N \to M$ be the composition $M \times D^N \xr{\hat f}M_\epsilon \xr{\pi}M$, and let $\hat \Upsilon_k$ be the composition $\bd^kW \times D^N \xr{r_{\bd^kW} \times \id}M \times D^N \xr{\hat F} M$ for $k\geq 0$.
	Once again by the Transversality Theorem of \cite[Section 2.3]{GuPo74}, for almost all $s \in D^N$ we have $V \xr{r_V} M \xr{\hat F(-,s)} M$ transverse to all cubical faces of $M$.
	Letting $s_0 \neq 0$ be any such point\footnote{If $s_0 = 0$ satisfies the conditions, we can finish immediately by letting $H$ be the constant homotopy.} we define $\hat H \colon M \times I \to M$ by $\hat H(x,t) = \hat F(x,ts_0)$, and we claim that this $\hat H$ satisfies the conditions required for the proposition.

	The map $\hat H$ is proper again by \cref{L: nearby proper homotopy} because $\phi(y) \eta(y) \leq \eta(y)$.
	The conditions of the proposition for $V$ follow immediately from the construction.
	It remains to verify that each $\hat h_k \defeq \hat H \circ (r_{\bd^k W} \times \id_I) \colon \bd^k W \times I \to M$ is transverse to the cubulation.
	Again, we focus primarily on $\hat h_0$.

	We begin by observing that for $(x,t) \in W \times I$ we can write $\hat h_0(x,t) \in M$ explicitly as
	$$\hat h_0(x,t) = \pi(r_W(x)+\phi(r_W(x))\eta(r_W(x))ts_0).$$
	So, alternatively, we can observe that $\hat h_0$ is the composition
	\begin{equation}\label{E: alt hat h}
		W \times I \xr{\Phi} W \times I \xhookrightarrow{\Psi} W \times D^N \xr{r_W \times \id} M \times D^N \xr{F} M,
	\end{equation}
	with $\Phi(x,t) = (x,\phi(r_W(x))t)$, $\Psi(x,t) = (x,ts_0)$, and noting that on the right we do mean our original $F$ and not $\hat F$.

	The first map $\Phi$ is a diffeomorphism onto its image, which is a neighborhood of $W \times 0$ in $W \times I$, and the map $\Psi$ embeds this into $W \times D^N$ by a product map that is constant in the $W$ direction and nontrivial linear in the second factor.
	The composition of the last two maps is just our earlier map $\Upsilon_0$.
	By construction, the map $r_W \times \id$ now takes the image of $\Psi\Phi$ into $M\times_\phi D^N$, and so at each point $(z,s)$ in the image of $\Psi\Phi$ if we fix $s$ and consider $\Upsilon_0(-,s)$ we get by construction a map on $W$ that is transverse at $z$ to the cubulation of $M$.
	Let $(x,t) \in W \times I$, let $\Psi \Phi(x,t) = (z,s)$, and let $\R s_0$ denote the line in $\R^N = T_sD^N$ spanned by the position vector of $s_0$.
	As $\Phi$ is a diffeomorphism onto its image and $\Psi$ is an embedding that is the identity with respect to $W$ and nontrivial linear on $I$, we see
	that the derivative of $\Psi\Phi$ maps the tangent space $T_{(x,t)}(W \times [0,1])$ onto $ T_zW \times \R s_0 \subset T_{(z,s)}(W \times D^N)$.
	In particular, this image contains $T_zW \times 0$, and by construction $D\Upsilon_0$ takes $T_zW \times 0$ to a tangent subspace in $M$ at $\hat h_0(x,t)$ that is transverse to the tangent space there of any closed face of the cubulation containing $\hat h_0(x,t)$.
	The same holds for $k>0$ replacing $W$ with $\bd^{k}W$ in \eqref{E: alt hat h} and $r_W$ with $r_{\bd^{k}W}$.
	So we see that $\hat H$ satisfies all the requirements of the proposition.
\end{proof}

\subsection{The intersection map and the isomorphism between cubical and geometric cohomology}\label{S: intersection map}

To define the intersection map, we introduce an augmentation map as in singular homology theory.
For this we first need a quick lemma.

\begin{lemma}\label{L: Q0}
	If $W \in Q_0(M)$, then $W$ has the same number of positively and negatively oriented points.
\end{lemma}

\begin{proof}
	As elements of $PC_0^\Gamma(M)$ cannot be degenerate, if $W \in Q_0(M)$ then $W$ must be trivial, and so there is an orientation-reversing diffeomorphism $\rho$ of $W$ such that $r_W\rho = r_W$.
	But a compact $0$-manifold has an orientation-reversing diffeomorphism if and only there are the same number of points with each orientation.
\end{proof}

\begin{definition}\label{D: aug}
	We define the \textbf{augmentation map}\index{augmentation map}\index{$A$@$\aug$} $\aug \colon PC^\Gamma_0(M) \to \Z$ as follows: If $W \in PC^\Gamma_0(M)$, then $W$ is the disjoint union of a finite number of points, each with orientation denoted $1$ or $-1$.
	We let $\aug(W)$ be the sum of the orientations of the points in $W$, interpreting $1$ and $-1$ as integers.
	By \cref{L: Q0}, an element of $PC^\Gamma_0(M)$ can be in $Q_0(M)$ only if this sum is $0$, so the augmentation descends to a homomorphism $\aug \colon C^\Gamma_0(M) \to \Z$.
	Furthermore, if $W \in PC_1^\Gamma(M)$ then $\aug(\bd W) = 0$, so $\aug$ further descends to a homomorphism $\aug \colon H_0^\Gamma(M) \to \Z$.
\end{definition}

Later, we will construct in general a partially-defined intersection map $C^*_\Gamma(M) \otimes C_*^\Gamma(M) \to C_*^\Gamma(M)$.
In general, this is delicate, as geometric chains and cochains do not have fixed representatives.
However, at the moment we do not need this full generality to define the intersection map we will need to compare geometric cohomology and cubical cohomology.
This is reflected in the following more limited definition:

\begin{definition}\label{D: intersection number}\index{intersection number|(}
	Suppose $M$ is an $m$-manifold without boundary and $W \in PC_\Gamma^i(M)$ and $N \in PC_{i}^\Gamma(M)$ are transverse.
	We define the \textbf{intersection number} $I_M(W,N)$ (or simply $I(W,N)$ if $M$ is clear from context) by $$I_M(W,N) = \aug(W \times_M N),$$ with $W \times_M N$ as defined in \cref{D: PC products}.
\end{definition}

We observe that this definition makes sense as $W$ and $N$ are transverse with complementary dimensions and $W \times_M N$ is an element of $PC_0^\Gamma(M)$.
In fact, in this case in order for transversality to hold the maps $r_W \colon W \to M$ and $r_N \colon N \to M$ must have full rank at each $x \in W$ and $y \in N$ such that $r_W(x) = r_N(y)$.
As having full rank is an open condition, the maps will also have full rank on neighborhoods of these points.
In particular, by the Implicit Function Theorem, they must be immersions on neighborhoods of these points.
So, locally, the orientation of $N$ determines an orientation of $T_yN$, which we can consider to be a subspace of $T_{y}M$, slightly abusing notation to identify $y$ and $r_N(y)$ via the local immersion.
Furthermore, in a neighborhood of $x$ the co-orientation of $r_W$ determines an orientation of the normal bundle of the local immersion of $W$, and we can take the fiber of the normal bundle at $r_W(x)$ to be $T_yN$.

\begin{lemma}\label{L: intersection number}
	The intersection number $I_M(W,N)$ is equal to signed count of intersection points of $W$ and $N$, counting an intersection point with $+1$ if the normal co-orientation of $W$ agrees with the orientation of $N$ and $-1$ otherwise.
\end{lemma}

\begin{proof}
	This follows directly from \cref{C: complementary cap}.
\end{proof}

\begin{lemma}\label{L: Q-trivial intersection}
	Suppose $W \in PC_\Gamma^i(M)$ and $N \in PC_{i}^\Gamma(M)$ are transverse and that $W \in Q^i(M)$.
	Then $I(W,N) = 0$.
\end{lemma}

\begin{proof}
	By \cref{L: pullback with Q}, we know $W \times_M N \in Q_0(M)$, so $I(W,N) = \aug(W \times_M N) = 0$ by definition and \cref{L: Q0}.
\end{proof}

\index{intersection number|)}

\begin{definition}\label{D: intersection homomorphism}\index{intersection map|textbf}
	Given the manifold without boundary $M$ cubulated by $X$, we define the \textbf{intersection map} $\mc I \colon C^*_{\Gamma \pf X}(M) \to K^*(X)$ by the possibly infinite formal sum $$\mc I(\uW) = \sum_F I_M(W,F)F^*,$$ where the sum is taken over faces $F$ of the cubulation $X$ such that $\dim(F)+\dim(W) = \dim(M)$ and the $W$ on the right hand side is any element of $PC^*_{\Gamma \pf X}(M)$ representing $\uW$.

	In particular, for a face $F$ of dimension $\dim(M)-\dim(W)$, we have $$\mc I(\uW)(F) = I_M(W,F) = \aug(W \times_M F).$$
\end{definition}

\begin{proposition}\label{P: I is well defined}\index{intersection map!is well defined}
	The intersection map $\mc I$ is a well-defined chain map.
\end{proposition}

\begin{proof}
	If $W, W' \in PC^*_{\Gamma \pf X}(M)$ are two representatives of $\uW$ then $W \sqcup -W' \in Q^*(M)$ and it is transverse to $X$.
	So for any face $F$ we have $\aug(W \times_M F)-\aug(W' \times_M F) = \aug((W \sqcup -W') \times_M F) = 0$ by \cref{L: Q-trivial intersection}.
	So $\mc I(\uW)$ does not depend on the choice of $W$.

	To see that $\mc I$ is a chain map, let $W$ be transverse to the cubulation and representing $\uW$.
	Then we compute for any face $f$ of $X$ that
	\begin{align*}
		\mc I(\bd\uW)(f)& = \aug((\bd W) \times_M f)\\
		& = \aug(W \times_M \bd f)\\
		& = \mc I(\uW)(\bd f).
	\end{align*}
	For the second equality, we use \cref{P: Leibniz cap} together with the facts that the augmentations are both trivial unless $\dim(W \times_M f) = 1$ and that $\aug(\bd (W \times_M f)) = 0$.
\end{proof}

\begin{example}\label{E: coho 0 generator}
	Let $M$ be any connected manifold without boundary and given a cubulation $X$, let $\uW \in C^*_{\Gamma \pf X}(M)$ be represented by the tautologically co-oriented identity map $M \to M$, and let $v$ be a (positively-oriented) $0$-dimensional vertex of $X$.
	By \cref{P: cap with 1}, we have $M \times_M v = v$, and so
	\[\mc I(\uW)(v) = \aug(M \times_M v) = \aug(v) = 1.\]
	In other words, $\mc I(\uW)$ takes the value $1$ on each vertex of $X$, so it is the standard generator of $H^0(X) \cong \Z$.
	It will therefore follow from \cref{T: intersection qi}, below, that the tautologically co-oriented identity map $M \to M$ generates $H^0_\Gamma(M) \cong \Z$, as promised in \cref{E: first examples}.
\end{example}

\begin{example}\label{E: uncountable}\index{geometric cochain!complex need not be free abelian}
	Let $Z$ denote the integers, considered as a space with the discrete topology.
	This space has a unique cubulation, which we will also write simply as $Z$.
	In this example, we demonstrate that the map $\mc I \colon C^0_\Gamma(Z) \to K^0(Z)$ is an isomorphism.
	As $K^0(Z) = \Hom(K_0(Z),\Z) \cong \Hom\left(\displaystyle\oplus_{i \in \Z} \Z, \Z\right) \cong \prod_{i\in \Z} \Z$, this shows that the geometric cochain groups $C^*_\Gamma(M)$ need not be free abelian \cite{Sch08}.

	To demonstrate the isomorphism, first consider an element $\alpha \in K^0(Z)$.
	The cochain $\alpha$ is determined by the value it assigns to each element of $Z$.
	For $x \in Z$, let $V(x)$ denote the discrete space with $|\alpha(x)|$ points, which we think of as mapping to $x \in Z$ with the co-orientation chosen according to whether $\alpha(x)$ is positive or negative (if it is zero we let $V(x)$ be empty).
	Then $V = \displaystyle\sqcup_{i \in \Z} V(x)$ is a countable collection of points with a proper co-oriented map to $Z$, and $\mc I(\uV) = \alpha$.
	So $\mc I$ is surjective.

	Next, suppose that $\uV \in C^0_\Gamma(Z)$, represented by $V \in PC^0_\Gamma(Z)$, and suppose $\mc I(Z) = 0$.
	As $V \to Z$ is proper, the preimage of each point $x \in Z$ is a finite set of points, and as $\mc I(\uV) = 0$, the restriction to the preimage $r_V^{-1}(x)$ must have the same number of points mapping to $x$ with each co-orientation.
	By matching the points in each $r_V^{-1}(x)$ into pairs mapping with opposite co-orientations, we obtain a co-orientation reversing automorphism of $V$ by swapping the elements of each pair.
	So $V$ is trivial and $\uV = 0$.
	This $\mc I$ is also injective.
\end{example}

Our goal now is to show that the intersection map $\mc I$ induces an isomorphism $H^i_{\Gamma \pf X}(M) \to H^i(X)$ whenever $H^i_\Gamma(M)$ and $H^i(X)$ are finitely generated.
Recall that we already know these groups are abstractly isomorphic by \cref{T: transverse complex,T: geometric is singular} and the footnote on page \pageref{FN: cubical and singular}.
Here is the formal statement of the theorem:

\begin{theorem}\label{T: intersection qi}\index{intersection map!is a cohomology surjection}\index{intersection map!is a cohomology isomorphism when finitely generated}
	If $M$ is a manifold without boundary cubulated by $X$, the intersection map $\mc I \colon H^i_{\Gamma \pf X}(M) \to H^i(X)$ is a surjection.
	If, in addition, $H^i(M)$ is finitely generated, then $\mc I$ is an isomorphism.
\end{theorem}

The proof will take us the remainder of \cref{S: transversality}.
We begin in the next subsection by using the cubical structure to start building an inverse map to $\mc I$.

We conjecture that $\mc I$ is always an isomorphism, but we have not been able to demonstrate this.

\begin{remark}\label{R: intersection map extension}
	Putting together the map $\mc I \colon H^*_{\Gamma \pf X}(M) \to H^*(X)$ with the inverse of the isomorphism $H^*_{\Gamma \pf X}(M) \to H^*_\Gamma(M)$ of \cref{T: transverse complex}, it is sometimes useful to abuse notation and speak of the ``intersection map'' $\mc I \colon H^*_\Gamma(M) \to H^*(X)$.
	This map is given by taking a cohomology class representative that is transverse to the cubulation and applying the intersection map $\mc I$ to find a cubical cocycle representing the target cohomology class.
\end{remark}

Before starting toward the proof, we provide an application to cohomology with coefficients.

\begin{corollary}\index{geometric cohomology!with coefficients}\label{C: cohomology with coefficients}
	If $H_*(M)$ is finitely-generated in each degree and $G$ is an abelian group, then $H^*(C^*_\Gamma(M) \otimes G) \cong H^i(M;G)$.
\end{corollary}
\begin{proof}
	First choose a cubulation $X$ of $M$.
	Next, as each $H_*(M)$ is finitely generated, by \cite[Lemma 56.3 and Theorem 45.1]{Mun84} we can choose a chain complex $F_*$ and a map $\phi \colon F_* \to K_*(X)$ such $F_*$ vanishes in negative degrees, is finitely-generated in each degree, and $\phi$ is a quasi-isomorphism.

	Now consider the maps
	$$C^*_\Gamma(M) \otimes G \hookleftarrow C^*_{\Gamma \pf X}(M) \otimes G \xr{\mc I \otimes \id} K^*(X) \otimes G = \Hom(K_*(X),\Z)\otimes G \xr{\phi^* \otimes \id} \Hom(F_*, \Z) \otimes G.$$
	Omitting the tensor products with $G$, these maps are all quasi-isomorphisms by \cref{T: transverse complex}, \cref{T: intersection qi}, and \cite[Theorem 45.5]{Mun84}.
	Furthermore, continuing to omit $\otimes G$, these complexes are all torsion free and so flat, and any submodule is thus also flat.
	So it follows from the Universal Coefficient Theorem for Homology \cite[Theorem 7.55 ]{ROTMAN} and the Five Lemma that this remains a chain of quasi-isomorphisms after tensoring with $G$.
	Now by Step 2 of the proof of \cite[Theorem 56.1]{Mun84}, $\Hom(F_*, \Z) \otimes G \cong \Hom(F_*, G)$, and this is quasi-isomorphic to $\Hom(K_*(X), G) = K^*(X; G)$, again by \cite[Theorem 45.5]{Mun84}.
\end{proof}

In this argument, the requirement that $H_*(M)$ be finitely generated is used not just to ensure $\mc I$ is a quasi-isomorphism.
It is needed to construct an $F_*$ and have $\Hom(F_*, \Z) \otimes G \cong \Hom(F_*, G)$.

\subsubsection{Dualization in cubes}\label{S: dual cubes}

Analogous to barycentric subdivisions of simplices, we will need to consider standard subdivisions of cubes.
For this we let $\jinterval$ denote the interval $\interval = [0,1]$ thought of as the (non-disjoint) union $[0,1/2] \cup [1/2,1]$.
We can then write \index{$J$@$\jinterval^n$} $\jinterval^n = \left([0,1/2] \cup [1/2,1]\right)^n$, with the idea being that we consider $\interval^n$ as the union of $2^n$ \textit{subcubes} of side length $1/2$, each of which is the product within $\interval^n$ of $n$ factors, each factor equal to either $[0,1/2]$ or $[1/2,1]$.
We refer to $\jinterval^n$ with this structure as the \textbf{central subdivision}\index{cube!central subdivision|textbf} of $\interval^n$.
It is a cubical complex with the evident partial orders obtained from $0 < 1/2 < 1$.
The resulting subdivision of each cubical face of $\interval^n$ is precisely the central subdivision obtained by treating that face as some $\interval^k$ in its own right.

Analogously to the case with $\interval^n$, each cube $S$ of $\jinterval^n$ possesses faces (some of which it shares with other $n$-cubes) consisting of subsets of $S$ in which some variables have been bound to the values $0$, $1$, or $1/2$.
In general we refer to such faces as \textbf{faces of $\jinterval^n$}.
If no variable of such a face is bound to $0$ or $1$, we say that we have an \textbf{internal face of $\jinterval^n$},\index{$J$@$\jinterval^n$!internal face|textbf} otherwise we call it an \textbf{external} face.\index{$J$@$\jinterval^n$!external face|textbf}
External faces are all subsets of $\bd \interval^n$; internal faces are not subsets of $\bd \interval^n$.

To each face $F$ of $\interval^n$, we refer to the point $\hat F$ \index{$F$@$\hat F$} at which all its free variables are equal to $1/2$ as the \textbf{center}\index{cube!center|textbf} of $F$.
Each vertex is its own center.
To each face $F$ of $\interval^n$ we define its \textbf{dual face in $\interval^n$},\index{cube!dual face|textbf} or simply its \textbf{dual}, to be the face $F^\vee$ \index{$F^\vee$} of $\jinterval^n$ determined as follows:
\begin{itemize}
	\item If $i \in F_{01}$ (i.e.\ the coordinate $x_i$ is free in $F$), then $x_i = 1/2$ in $F^\vee$.

	\item If $i \in F_0$ (i.e.\ the coordinate $x_i$ is bound to $0$ in $F$), then $x_i$ is free in $[0,1/2]$ in $F^\vee$.

	\item If $i \in F_1$ (i.e.\ the coordinate $x_i$ is bound to $1$ in $F$), then $x_i$ is free in $[1/2,1]$ in $F^\vee$.
\end{itemize}

It is clear that $F$ and $F^\vee$ have complementary dimensions and that they intersect naively transversely in the center of $F$.

\begin{lemma}
	The set function $F \to F^\vee$ is a bijection between the faces of $\interval^n$ and the internal faces of $\jinterval^n$.
\end{lemma}

\begin{proof}
	Injectivity is clear as two different faces of $\interval^n$ will have different partitions of $\bar n$.

	Next consider an internal face of $\jinterval^n$.
	By definition this is a set in which some set of variables $A$ has been bound to $1/2$ while two other sets of variables, $B$ and $C$, are free on $[0,1/2]$ or $[1/2,1]$, respectively.
	But this is $F^\vee$ for the face $F$ with partition $(B,A,C)$.
	So our function is surjective.
\end{proof}

As each face $F$ of $\interval^n$ carries a natural orientation $\beta_F$ determined by the order of its free variables (or the orientation $1$ for vertices), this provides $F^\vee$ with the corresponding normal co-orientation.
In other words, $F^\vee$ is co-oriented at all points by the co-orientation $(\beta_{F^\vee}, \beta_{F^\vee} \wedge \beta_F)$, where $\beta_{F^\vee}$ is an arbitrary orientation of $F^\vee$.
Assigning $F^\vee$ this co-orientation, we interpret its embedding in $\interval^n$ as representing a cochain in $\interval^n$ of index $\dim(F)$.
We define a map $\Psi \colon K^*(\interval^n) \to C_\Gamma^*(\interval^n)$ given on generators by $\Psi(F^*) = F^\vee$, recalling that geometric chains and cochains can be defined on manifolds with corners --- see \cref{D: PC,D: chains and cochains}.
To avoid overburdening the notation too much, throughout the remainder of this section we will simply write $F^\vee$ for the output geometric cochain, rather than $\underline{F^\vee}$; it should be clear from context when we refer to the cube face and when we refer to its corresponding geometric cochain.

We introduce one more piece of notation for the following lemma.
If $f$ is a face in $\jinterval^n$, we let $\bd_{\text{int}}f$ \index{$bd$@$\bd_{\text{int}}f$} denote the union of the internal $\dim(f)-1$ faces of $\bd f$ and $\bd_{\text{ext}}f$ \index{$bd$@$\bd_{\text{ext}}f$} denote the union of the $\dim(f)-1$ external faces of $f$.
If $f$ is co-oriented, then we interpret the terms of $\bd_{\text{int}}f$ and $\bd_{\text{ext}}f$ as co-oriented with the boundary co-orientations.
When $f$ is interpreted as a geometric cochain, we interpret $\bd_{\text{int}}f$ and $\bd_{\text{ext}}f$ as sums of geometric cochains, in which case $\bd f = \bd_{\text{int}}f + \bd_{\text{ext}}f$.
We can then extend $\bd_{\text{int}}$ and $\bd_{\text{ext}}$ to linear operators in the obvious way.

\begin{lemma}\label{L: dualizing bijection}
	For any face $F$ of $\interval^n$, we have $\Psi(d F^*) = \bd_{\text{int}}\Psi(F^*).$
\end{lemma}

\begin{proof}
	We first show that there is a bijection between the internal faces of $\Psi(F^*)$ and the set of $E^\vee$ such that $E^*$ has non-zero coefficients in $dF^*$.
	Then we will return to carefully consider the co-orientations.

	So let $F$ be a face of $\interval^n$.
	Recall that $F$ is determined by the partition $(F_0,F_{01}, F_1)$ of $\overline{n}$ corresponding respectively to variables set to $0$, free variables of $F$, and variables set to $1$.
	By definition, $(dF^*)(E) = F^*(\bd E)$ and so the only faces participating in $dF^*$ are those $\dim(F)+1$ faces that have $F$ as a boundary, in other words those faces whose free variables are those in $F_{01}$ plus one more variable from $F_0$ or $F_1$.

	Now let us consider again $F^\vee$, which we recall is obtained by setting all variables in $F_{01}$ to $1/2$ and letting the variables in $F_0$ and $F_1$ become free variables on $[0,1/2]$ and $[1/2,1]$, respectively.
	Each boundary face of $F^\vee$ is then obtained by either setting one of the variables in $F_0$ to $0$ or $1/2$ or one of the variables in $F_1$ to $1/2$ or $1$.
	Furthermore, the internal faces are those where the variable in $F_0$ or $F_1$ has been set to $1/2$.
	So, in summary, an internal face of $F^\vee$ has all of the variables in $F_{01}$ as well as exactly one other variable set to $1/2$ and the rest remain free over the appropriate domains, namely $[0,1/2]$ for those in $F_0$ and $[1/2,1]$ for those in $F_1$.
	But, from the definition of dualization $E \to E^\vee$, this exactly describes the duals $E^\vee$ of the faces participating in $dF^*$, which is sufficient for our desired bijection due to the bijection of \cref{L: dualizing bijection}.

	It remains now to consider the signs.
	We recall that the boundary formula for cubes has the form
	$$\bd E = \sum_{i = 1}^k (-1)^i(E \delta_i^0-E \delta^1_i),$$ where the $\delta$s denote the embeddings of the faces.
	In $K_*(\interval^n)$, we can shorten this notation to $$\bd E = \sum_{i = 1}^k (-1)^i(E_i^0-E^1_i),$$ letting $E_i^j$, $j \in \{0,1\}$, denote the $i$th ``front or back face'' according to $j = 1$ or $j = 0$.
	In particular, we note that there are two factors, both $i$ and $j$, affecting sign.

	So now let us again fix a face $F$ of $\interval^n$ and let $E$ be a face of dimension $\dim(F)+1$ that includes $F$ in its boundary as $F = E_i^j$.
	Representing $E$ as $E = (E_0,E_{01},E_1)$, we obtain $F$ by setting the $i$th variable, which is in $E_{01}$, to $j$.
	From the coboundary formula, we have
	\begin{equation*}
		(dF^*)(E) = F^*(\bd E) = (E_i^j)^*\left(\sum_{i = 1}^k (-1)^i(E_i^0-E^1_i)\right)
		= (E_i^j)^*((-1)^{i+j}E_i^j) = (-1)^{i+j}.
	\end{equation*}
	So 	$E^*$ occurs in $dF^*$ with coefficient $(-1)^{i+j}$.
	Thus we must show that $\Psi(E^*) = E^\vee$ disagrees by a sign of $(-1)^{i+j}$ from the co-orientation assigned to $E^\vee$ as a member of $\bd_{\text{int}} \Psi(F^*) = \bd_{\text{int}} (F^\vee)$.

	We continue to assume $F = E_i^j$ and now consider $F^\vee$.
	If $A=(a,b,\ldots)$ is an ordered subset of $\overline{n}$, we let $\beta_A = e_a \wedge e_b \wedge \cdots$, where $e_a$ is, as usual, the standard unit vector in the $a$th coordinate.
	If $A, B$ are disjoint ordered subsets sets, we write $A \cup B$ for their concatenation.
	With this notation, we can write the normal co-orientation of $F^\vee$ as $(\beta_{F_0 \cup F_1},\beta_{F_0 \cup F_1} \wedge \beta_{F_{01}})$.
	Let $k \in F_0 \cup F_1$ be the unique index that is free in $E$ but bound in $F$ and so also free in $F^\vee$.
	Then $E^\vee$ is the boundary component of $F^\vee$ with $x_k$ set to $1/2$.
	In this case the boundary co-orientation of $E^\vee$ in $F^\vee$ is $(\beta_{F_0 \cup F_1 - \{k\}},\beta_{F_0 \cup F_1 - \{k\}} \wedge (-1)^{j+1} e_k)$.
	The sign $(-1)^{j+1}$ is because if $j = 1$ then in $F^\vee$ the variable $x_k$ is free on $[1/2,1]$ so the inward normal vector at $1/2$ points in the direction of $e_k$, and if $j = 0$ then in $F^\vee$ the variable $x_k$ is free on $[0,1/2]$ so the inward normal vector at $1/2$ points in the direction of $-e_k$.
	So the boundary co-orientation for $E^\vee$ in $\interval^n$ as a piece of $\bd F^\vee$ is the composition $$(\beta_{F_0 \cup F_1 - \{k\}},\beta_{F_0 \cup F_1 - \{k\}} \wedge (-1)^{j+1} \beta_{e_k})*(\beta_{F_0 \cup F_1},\beta_{F_0 \cup F_1} \wedge \beta_{F_{01}}).$$

	Meanwhile, the natural co-orientation for $E^\vee$ is $(\beta_{E_0 \cup E_1},\beta_{E_0 \cup E_1} \wedge \beta_{E_{01}}).$
	But $E_0 \cup E_1 = F_0 \cup F_1 - \{k\}$, so
	\begin{align*}
		(\beta_{F_0 \cup F_1 - \{k\}},\beta_{F_0 \cup F_1 - \{k\}}& \wedge (-1)^{j+1} \beta_{e_k})*(\beta_{F_0 \cup F_1},\beta_{F_0 \cup F_1} \wedge \beta_{F_{01}})\\
		& = (-1)^{j+1}(\beta_{E_0 \cup E_1},\beta_{E_0 \cup E_1} \wedge \beta_{e_k})*(\beta_{E_0 \cup E_1} \wedge \beta_{e_k},\beta_{E_0 \cup E_1} \wedge \beta_{e_k} \wedge \beta_{F_{01}})\\
		& = (-1)^{j+1}(\beta_{E_0 \cup E_1},\beta_{E_0 \cup E_1} \wedge \beta_{e_k} \wedge \beta_{F_{01}})\\
		& = (-1)^{j+i}(\beta_{E_0 \cup E_1},\beta_{E_0 \cup E_1} \wedge \beta_{E_{01}}).
	\end{align*}
	In the second line, we use that expressions of the form $(\beta_A, \beta_A \wedge \beta_B)$ are independent of the choice of $\beta_A$, so we replace $\beta_{F_0 \cup F_1}$ with $\beta_{E_0 \cup E_1} \wedge \beta_{e_k}$.
	But as $F = E_i^j$, we know that $k$ is the $i$th variable of $E_{01}$, so
	$\beta_{e_k} \wedge \beta_{F_{01}} = (-1)^{i-1}\beta_{E_{01}}$, which we use in the last equality.
	Thus the two co-orientations disagree by $(-1)^{i+j}$ as required.
\end{proof}

We need some additional notation for the next lemma.
Suppose $F$ is a face of $E$, which we will write $F \subseteq E$, recalling that we identify faces of a cubical complex with their vertex sets.
Let $F_E^\vee$ denote the dual of $F$ in $E$.
In other words, we identify $E$ with $\interval^k$ for an appropriate $k$ and form the dual of the face corresponding to $F$ in $E$.

\begin{lemma}\label{L: ext faces}
	Let $F$ be a face of $\interval^n$.
	Then the external faces $\bd_{\text{ext}}F^\vee$ of $F^\vee$ correspond exactly to the $F_E^\vee$ as $E$ ranges over the $n-1$ faces of $\interval^n$ containing $F$.
	Furthermore, if $F^\vee$ is given the normal co-orientation induced by $F$ as above, then the co-orientation of $F^\vee_E \to F^\vee \to \interval^n$ considered as a piece of the boundary of $F^\vee$ is given by $(\beta_{F_E^\vee},\beta_{F_E^\vee} \wedge \beta_v \wedge \beta_F)$, where $\beta_{F_E^\vee}$ is any arbitrary orientation for $F_E^\vee$, $v$ is an inward pointing normal to $\interval^n$ at $E$, and $\beta_F$ is the orientation of $F$ as a face of $\interval^n$.
\end{lemma}

\begin{proof}
	We first observe that each $n-1$ cube $E$ is determined by its single bound variable, and any cube $F$ in $E$ must share that bound variable.
	Conversely, given $F$, to obtain an $n-1$ cube containing it, we free all but one of its bound variables.
	So the $n-1$ cube $E$ contains $F$ if and only if there is an $i \in F_0$ such that $i \in E_0$ or an $i \in F_1$ such that $i \in E_1$.

	Now from the definitions, we recall that $F^\vee$ is obtained by setting the free variables of $F$ to $1/2$ and the bound variables of $F$ to be free, with the variables of $F_0$ running over $[0,1/2]$ and the variables of $F_1$ running over $[1/2,1]$.
	So if $f$ is an external face of $F^\vee$, it is obtained from the conditions defining $F^\vee$ by setting $x_i=0$ for an $i \in F_0$ or $x_i$ to $1$ for an $i \in F_1$.
	Thus $x_i$ is the only variable of $f$ bound to $0$ or $1$, and so the corresponding $E$ determined by $x_i = 0$ or $x_i = 1$ is the unique $n-1$ cube of $\interval^n$ containing $f$, and it also contains $F$ (since $i$ is in both $F_0$ and $E_0$ or both $F_1$ and $E_1$).
	Conversely, if $F \subseteq E$ and $i \in F_0$ with $x_i = 0$ defining $E$, then there is a unique external face of $F^\vee$ for which $x_i = 0$, namely the one with $x_i=0$, $x_k = 1/2$ for all $k \in F_{01}$, and all other variables free.
	Similarly, if $F \subseteq E$ and $i \in F_1$ with $x_i = 1$ defining $E$.
	So we have a bijection between external faces of $\bd F^\vee$ and the $n-1$ faces $E$ containing $F$.
	Furthermore, after we have set $x_i$ to $0$ or $1$ as appropriate, the behavior in the remaining variables shows that the external face of $\bd F^\vee$ contained in $E$ is precisely $F^\vee_E$.

	It remains to consider the co-orientations.
	Suppose $f$ is an external face of $F^\vee$ in the $n-1$ face $E$, i.e.\ $f = F^\vee_E$.
	The co-orientation of $F^\vee$ in $\interval^n$ and the co-orientation of $f$ in $E$ are both the normal co-orientations determined by the orientation of $F$.
	As $f$ is part of $\bd F^\vee$, its boundary co-orientation in $\interval^n$ is therefore the composite $$(\beta_f,\beta_f \wedge \beta_v)*(\beta_f \wedge \beta_v,\beta_f \wedge \beta_v \wedge \beta_F) = (\beta_f,\beta_f \wedge \beta_v \wedge \beta_F),$$ where $\beta_f$ is any arbitrary orientation for $f$, $v$ is the inward pointing normal of $f$ in $F^\vee$, which in this case is also an inward pointing normal of $E$ in $\interval^n$, and we use $\beta_f \wedge \beta_v$ as a convenient orientation for $F^\vee$.
\end{proof}

\subsubsection{Dualization in cubical complexes}\label{S: dualization of complexes}

Now suppose that $X$ is any cubical complex.
We obtain the \textbf{centrally subdivided cubical complex $sd(X)$}\index{cubical complex!centrally subdivided} by replacing each cube $\interval^n$ with $\jinterval^n$.
This is well defined as we have noted that the central subdivision of each cube is consistent with the central subdivision of each cube of which it is a face.

Now suppose that $M$ is an $n$-manifold without boundary cubulated by the cubical complex $X$ via $\phi \colon |X| \to M$.
For a given face $F$ of an $n$-cube $B$ of $X$, we can extend the above definitions, and abuse notation, by letting $F$ refer also to the image of $F$ in $M$ and letting $F_B^\vee$ also denote the composition $F_B^\vee \into \interval^n \into |X| \xr{\phi} M$, where $F_B^\vee$ is the dual of $F$ in $B$ as above.
Similarly, we consider this version of $F_B^\vee$ in $M$ to be co-oriented by the normal co-orientation obtained from the orientation of $F$, and we can then further abuse notation and consider $F_B^\vee$ to also represent a geometric cochain.
Moreover, in the cubulated $M$ we extend our earlier definition and write $F^\vee = \sqcup_{F \subseteq B} F^\vee_B$, where the union is taken over all $n$-cubes $B$ having $F$ as a face; when considering these as geometric cochains, we can replace union with sum.
We similarly define $\Psi_B(F^*) = F_B^\vee$, which is consistent with the dualization within a single cube defined above, and then set $\Psi(F^*) = \sum_{F \subseteq B} \Psi_B(F^*) = \sum_{F \subseteq B} F_{B}^\vee = F^\vee$.
This is sufficient to determine a homomorphism $\Psi \colon K^*(X) \to C_\Gamma^*(M)$, thinking of elements of $K^*(X)$ as possibly infinite formal sums $\sum n_i F^*$ with the sum taken over all faces of a given dimension; compare \cite[Section 42]{Mun84}. The image of such an element of $K^*(X)$ under $\Psi$ corresponds to a map from a manifold with corners to $M$ whose domain is a countable collection of cubes, and this map is proper as every point of $M$ has a neighborhood that intersects only finitely many image cubes.

\begin{lemma}\label{L: dual chain map}
	$\Psi$ gives a chain map $K^*(X) \to C_\Gamma^*(M)$.
\end{lemma}

\begin{proof}
	It suffices to show that for any face $F$ of the cubulation $X$ of $M$ we have $\Psi(d F^*) = \bd\Psi(F^*)$.

	By \cref{L: dualizing bijection} and the definition of the extended $\Psi$, we know that
	$$\Psi(dF^*) = \sum_{F \subseteq B} \Psi_B(dF^*) = \sum_{F \subseteq B}\bd_{\text{int}} \Psi_B(F^*) = \sum_{F \subseteq B} \bd_{\text{int}}F^\vee_B,$$
	where $\bd_{\text{int}}$ refers to the internal boundary in the appropriate $B$, so we must show that $\sum_{F \subseteq B} \bd_{\text{int}}(F^\vee_B) = \sum_{F \subseteq B} \bd (F^\vee_B)$, in other words that all of the external faces of the terms of $\bd\Psi(F^*) = \sum_{F \subseteq B} \bd (F^\vee_B)$ cancel in the sum.

	For this, let us fix an $n$-cube $B$ containing $F$, and let $E$ be an $n-1$ face of $B$ containing $F$.
	As $M$ is an $n$-manifold, there is exactly one other $n$-cube, say $B_1$, in the cubulation that shares the face $E$ with $B$.
	Furthermore, by \cref{L: ext faces}, $E$ contains external faces of $\bd(F^\vee_B)$ and $\bd(F^\vee_{B_1})$ and they both correspond geometrically to $F^\vee_E$.
	Also by the lemma, the co-orientation as a boundary piece of $F^\vee_B$ is $(\beta_{F_E^\vee},\beta_{F_E^\vee} \wedge \beta_{v_B} \wedge \beta_F)$, where $v_B$ is an inward pointing vector of $B$ at $E$, and the co-orientation as a boundary piece of $F^\vee_{B_1}$ is $(\beta_{F_E^\vee},\beta_{F_E^\vee} \wedge \beta_{v_{B_1}} \wedge \beta_F)$, where $v_{B_1}$ is an inward pointing vector of $B_1$ at $E$.
	Since a normal vector to $E$ that is inward pointing with respect to $B_1$ is outward pointing with respect to $B$, and vice versa, these are opposite co-orientations, so, as precochains, these two boundary pieces constitute a trivial pair and so cancel in $C^*_\Gamma(M)$.

	Applying this argument to all $n-1$ cubes that contain $F$ completes the argument by \cref{L: ext faces}.
\end{proof}

Morally speaking, $\Psi$ will be our inverse to the intersection map $\mc I$, but unfortunately the elements of $PC^*_\Gamma(M)$ that represent the image of $\Psi$ are not transverse to the cubulation.
So we need another step.

\subsubsection{Pushing the dual cubulation}\index{cubulation!pushing the dual|(}

To remedy the problem with the image of $\Psi$ consisting of elements of $PC^*_\Gamma(M)$ that are not transverse to the cubulation, we construct a homotopy of $M$ to itself that pushes these cochains into transverse position with respect to the cubulation, which we regard as fixed.
Moreover, we want to do so in a way such that the intersection number of the pushed $F^\vee$ with $F$ will be $1$, as expected.
Constructing such a homotopy is the purpose of the following technical lemma.
As we will see in the proof, the actual construction is a bit fiddly in order to ensure all the needed properties.

\begin{lemma}\label{L: push dual}
	Suppose $M$ is a manifold without boundary with a cubulation $X$.
	Then there is a proper smooth homotopy from the identity to a map $h \colon M \to M$ such that for each face $F$ of $X$ the following hold:
	\begin{enumerate}
		\item $h|_{F^\vee} \colon F^\vee \to M$ is transverse to $X$
		\item the only $\dim(F)$-face of $X$ that intersects $h(F^\vee)$ is $F$
		\item $I_M(h(F^\vee),F) = 1$.
	\end{enumerate}
\end{lemma}

\begin{proof}
	We will construct our homotopy in multiple steps.
	The basic idea is first to construct small isotopies near the centers $\hat F$ of the faces $F$ of $X$ that push the corresponding $F^\vee$ into transverse position with $F$ and so that the intersection number of the shifted $F^\vee$ with $F$ is $1$, giving the third condition.
	These will also be designed to ensure enough transversality in neighborhoods of the $\hat F$.
	Then we perform a more global shift to push the rest of $sd(X)$ into transverse position with $X$ to achieve the first condition, while leaving fixed the isotopies already constructed near the $\hat F$.
	We also ensure the global shift is small enough to provide the second condition while also not disrupting the third.

	Throughout we may assume that $M$ is embedded properly in some $\R^N$ with an $\epsilon$-neighborhood $M^\epsilon$ (with $\epsilon$ a function of $x \in M$) and a submersion $\pi \colon M^\epsilon \to M$; see \cite[Section 2.3]{GuPo74}.
	We may further assume that $\epsilon(x) < 1$ for all $x \in M$.
	As a subspace of $\R^N$, we can consider $M$ as a metric space with the subspace metric.

	To create a template and explain the basic idea, we first consider the standard $m = \dim(M)$-cube $\interval^m = [0,1]^m \subset \R^m$.
	Let $e_i$ denote the unit tangent vector of $\R^m$ pointing in the positive $i$th direction.
	Let $f = (f_0,f_{01},f_1)$ be a face of $\interval^m$, let $\hat f$ be its center, and define the vector $v_f$ to be the unit vector in the direction $\sum_{i \in f_1} e_i-\sum_{i \in f_0} e_i$ (if $f_0 = f_1 = \emptyset$, let $v = 0$; we leave the reader to adjust for this case in what follows, noting that in this case $f^\vee = \hat f$ is already transverse to $f = \interval^m$).
	At points of $f$, the vector $v_f$ points outward from the cube $\interval^m$.
	The vector $v_f$ is also tangent to $f^\vee$ at any point of $f^\vee$, and for small $\delta>0$ the points $\hat f-\delta v_f$ lie in the interior of the face $f^\vee$.
	So, if we were to translate $f^\vee$ by $\delta v_f$ for a sufficiently small positive $\delta$, we would obtain a translate of $f^\vee$ within the $m-\dim(f)$ plane containing it.
	As this plane is orthogonal to $f$ and intersects it at the single point $\hat f$, the only intersection of $f$ with the translated $f^\vee+\delta v_f$ is at $\hat f$ with a point in the image of the interior of $f^\vee$.
	As the translation preserves the tangent plane of $f^\vee$, this is a transverse intersection, and the translation preserves the normal co-orientation so that the intersection number of $f^\vee+\delta v_f$ with $f$ is $1$.
	Furthermore, under a sufficiently small perturbation of $\delta v_f$ to a nearby $z_f$ it remains the case that the translate $f^\vee + z_f$ intersects $f$ transversely in a single point, and the only intersection of $f$ with $f^\vee + z_f$ is with the translate of a point in the interior of $f^\vee$ (though this intersection need not be at $\hat f$) and the intersection number is 1.
	Additionally, by the Transversality Theorem of \cite{GuPo74}, for almost all vectors $u$ in $\R^m$, translation by $u$ takes all open faces of $\jinterval^m$ into transverse position with respect to all open faces of $\interval^m$.
	So almost all vectors $z_f$ sufficiently close to $\delta v_f$ have the property that translation by $z_f$ both makes all open faces of $\jinterval^m$ transverse to all open faces of $F$ and satisfies that the intersection number of $f^\vee$ and $f$ is $1$.

	Next, suppose we are given a small closed ball $D^f_r$ centered at $\hat f$ in $\R^m$ with $r$ sufficiently small that $D^f_r$ does not intersect any faces of $\interval^m$ that do not contain $f$ (and in particular $D^f_r$ does not intersect any other faces of the same dimension as $f$ besides $f$ itself).
	We choose numbers $0<r_4 < r_3 < r_2<r_1<r$ and let $D^f_{r_j}$ be closed balls of radius $r_j$ centered at $\hat f$.
	In the first part of the argument we will be concerned primarily with $r$, $r_1$, and $r_2$, but we will need $r_3$ and $r_4$ later.
	We choose a smooth function $\eta \colon \R^m \to [0,1]$ such that $\eta(x) = 1$ on $D^f_{r_2}$ and $\eta(x) = 0$ outside of $D^f_{r_1}$.
	Let $\delta, \delta_1 > 0$ satisfy $0 < \delta \ll \min\left\{\frac{r_4}{2}, \frac{r-r_1}{2}\right\}$, $0 < \delta_1 < \delta$, and $\delta_1$ is smaller than any of the non-zero coordinates of $\delta v_f$.
	We will impose an additional condition on $\delta_1$ below.
	Next, let $z_f$ be a vector in $\R^m$ such that $|z_f-\delta v_f|\ll \delta_1$; by the condition on $\delta_1$, the signs of the coordinates of $z_f$ in $f_0$ and $f_1$ agree with those of $\delta v_f$, so $z_f$ continues to point outward from $\interval^m$ at all points of $f$.
	Finally, let $\Phi \colon \R^m \to \R^m$ be defined by $\Phi(x) = x+ \eta(x) z_f$.

	The map $\Phi$ has the following properties:
	\begin{enumerate}
		\item Outside of $D^f_{r_1}$, the map $\Phi$ is the identity.

		\item As $|z_f|\leq |v_f| + |z_f-\delta v_f| < \delta+\delta_1 < 2 \delta < r-r_1$ and no point moves further than $|z_f|$ under $\Phi$, we have $\Phi(D^f_{r_1})\subset D^f_r$.

		\item On $D_{r_2}^f$, the map $\Phi$ is translation by $z_f$ as in the preceding paragraph.
		As $\delta < r_4 < r_2$, if we translate $f^\vee \cap D_{r_2}^f$ by $\delta v_f$, then there is exactly one intersection point between $(f^\vee \cap D_{r_2}^f) + \delta v_f$ and $f$, and the resulting intersection number is $1$.
		As in the preceding paragraph, we choose $\delta_1$ so small that this property is maintained replacing $\delta v_f$ with $z_f$.
		Then there is exactly one intersection point between $\Phi(f^\vee \cap D_{r_2}^f)$ and $f$, and the resulting intersection number is $1$.

		\item\label{I: pushout} As $\Phi$ only moves points in $D_{r_1}^f$ and acts on all of $\R^m$ by translation by a non-negative scalar multiple of $z_f$ (which is a positive scalar multiple on $D_{r_2}^f$) and as $z_f$ points outward from $\interval^m$ at $f$ in the $f_0$ and $f_1$ directions, it holds for any $x \in \R^m - \interval^m$ or any $x \in f \cap D_{r_2}^f$ that $\Phi(x) \in \R^m - \interval^m$.

		\item As $|z_f| < 2\delta < r_4$, together with the above facts, the distance between the set $\Phi(f^\vee-(f^\vee \cap D_{r_4}^f))$ and the union of $\dim(f)$-faces of $\interval^m$ is positive.

		\item\label{I: Phi homotopy} $\Phi$ is homotopic to the identity via the homotopy $(x,t) \mapsto x+ t \eta(x) z_f$.
	\end{enumerate}

	We now translate this procedure to the cubulation $X$ of $M$.
	Suppose $F$ is a face of the cubulation, and choose an $m$-cube $B$ of the cubulation containing $F$ as a face.
	By assumption there is a smooth diffeomorphism $\phi$ from $\interval^m$ to $B$ that is compatible with the cubical structure.
	Let $f$ be the face of $\interval^m$ that maps to $F$.
	Then $\phi$ takes $\hat f$ to the center $\hat F$ and $f^\vee_{\interval^m}$, the dual of $f$ in $\interval^m$, to the component of $F_B^\vee$ of $F^\vee$.
	By the definition of smooth maps, there is a neighborhood $U$ of $\hat f$ in $\R^m$ on which there is defined a smooth map $\phi_U \colon U \to M$ and such that $\phi_U$ agrees with $\phi$ on $U \cap \interval^m$.
	As $\phi$ is an embedding at $\hat f$, by the Inverse Function Theorem we may choose $U$ such that $\phi_U$ is a diffeomorphism from $U$ onto its image, and by making $U$ smaller if necessary, we can assume its image in $M$ intersects only faces of $X$ that contain $F$, and in particular no other $\dim(F)$ face of $X$ besides $F$.
	We then choose a ball $D_r^f \subset U$ centered at $\hat f$.
	With this ball $D_r^f$, we can now choose quantities and define a $\Phi$ as in the preceding paragraph that is the identity outside of $D_r^f$ and that pushes $f^\vee_{\interval^m}$ into transverse position with $f$ with intersection number $1$.
	Identifying $D_r^f$ in $\R^m$ with its image in $M$ via $\phi_U$ and using that $\Phi$ is fixed outside $D_{r_1}^f$, we obtain a map of $M$ to itself that corresponds to $\Phi$ under the identification this.
	We will henceforth abuse notation and use the same symbols for subsets of $\R^m$ and the corresponding subsets of $M$; under this identification, we identify $F$ with $f$ in a neighborhood of their centers.
	By property \eqref{I: pushout} of $\Phi$ listed above, the corresponding map $\Phi$ on $M$ pushes the other components $F^\vee_C$ of $F^\vee$ for $C\neq B$ away from $B$, and so it continues to hold that $I(\Phi(F^\vee),F) = 1$ and also the distance between the set $\Phi(F^\vee-(F^\vee \cap D_{r^4}^F))$ and the union of $\dim(F)$-faces of $X$ is positive.
	To see this last fact, we can use that $\Phi$ fixes the complement of $D_{r_1}^f$, and consider the intersections $D_{r_1}^f \cap F^\vee_C$ as subsets of $\R^m$ via the identification above.
	Once again, as $|z_f| < 2 \delta < r_4$ and as $\Phi(D_{r_1}^f)\subset D_r^f$, we can find a fixed positive lower bound for the distances between $\Phi(F^\vee-(F^\vee \cap D_{r^4}^F))$ and $F$ or any other $\dim(F)$ face of $X$.
	As in the proof of the Transversality Theorem of \cite{GuPo74}, by rechoosing our $z_f$ if necessary (still with $|z_f - \delta v_f| < \delta_1$), we can maintain the above properties and also ensure that the restriction of $\Phi$ to the intersection of any face of $sd(X)$ with $int(D_{r_2}^F)$ is transverse to the cubulation\footnotemark.

	\footnotetext{Let $e$ and $c$ be respectively the intersections of faces of $sd(X)$ and $X$ with $\phi_U$; via $\phi_U$ we identify these as submanifolds of $\R^m$.
		On $D_{r_2}^F$ in $\R^m$, the map $\Phi$ acts by translation.
		But the map $\tau \colon e \times \R^m \to \R^m$ given by translation by the second coordinate is a submersion, so by the Transversality Theorem the map $\tau(-,y)$ is transverse to $c$ for almost all $y$.
		As a finite number of cubes of $X$ and $sd(X)$ intersect $U$, almost all vectors in any neighborhood of $v_f$ will thus satisfy the claim.}

	Next, we can apply this procedure at all faces $F$ of $X$ of all dimensions simultaneously by choosing a sufficiently small $D^F_r$ neighborhood of each $\hat F$ so that they are all pairwise disjoint; we let each $r$, $r_1$, $r_2$, $\delta$ and $\delta_1$ depend on $F$ though we do not include this in the notation.
	We then generalize $\Phi$ by allowing a corresponding map on all balls $D_r^F$ simultaneously, which is well defined letting $\Phi$ be the identity outside of the disjoint $D_{r_1}^F$.
	This provides a smooth map $\Phi \colon M \to M$ that satisfies the second two conditions of the theorem.
	We also observe that $\Phi$ is homotopic to the identity by extending the homotopies observed in property \eqref{I: Phi homotopy} of $\Phi$ to be the identity outside of the $D_{r_1}^F$.
	To complete the proof, we must further modify $\Phi$ to ensure transversality beyond the interiors of the $D_{r_2}^F$ and without destroying the other required properties.

	Let $W_j = \cup_F D^F_{r_j}$, where the union is taken over all faces of $X$ of all dimensions.
	By construction, the restriction of $\Phi$ to the intersection of the interior of $W_2$ with any face of $sd(X)$ is transverse to all faces of $X$.

	As in the proof of \cref{T: transverse complex}, we next follow the construction in \cite[Section 2.3]{GuPo74}.
	We have $M$ embedded in some Euclidean space $\R^N$ with an $\epsilon$-neighborhood $M^\epsilon$ and a submersion $\pi \colon M^\epsilon \to M$; we also assume $\epsilon(x)<1$ for all $x \in M$.
	As we are happy with the map $\Phi$ as constructed so far on the interior of $W_2$, we let $\rho \colon M \to [0,1)$ be a smooth function that is $0$ on $W_3$ and positive on $M-W_3$.
	We will fine tune $\rho$ a bit more soon.
	Let $S$ be the open unit ball in $\R^N$.
	We now consider the map $H \colon M \times S \to M$ defined by $H(x,s) = \pi(\Phi(x)+\rho(x)\epsilon(\Phi(x)) s)$.
	At all points $(x,s)$ such that $\rho(x)>0$, i.e.\ on $M-W_3$, this is a submersion (and so transverse to all faces of $X$), and for all $(x,s)$ such that $\rho(x) = 0$, i.e.\ on $W_3$, we have $H(x,s) = \pi(\Phi(x)) = \Phi(x)$.

	Now let $G$ be the interior of any face of the cubical subdivision $sd(X)$.
	At any point $x \in G \cap W_3$, we already have that for any fixed $s_0 \in S$ the map $H|_G(-,s_0) = \Phi|_G(-)$ is transverse at $x$ to any face of $X$.
	Furthermore, by the Transversality Theorem of \cite{GuPo74}, for almost every $s_0 \in S$ and for any face $F$ of $X$, the map $H|_G(-,s_0)$ is transverse to the interior of $F$ at all points on the submanifold $G-G \cap W_3$ of $G$.
	But there are a countable number of faces of $X$, so for almost every $s_0 \in S$ and for every face $F$ of $X$, the map $H|_G(-,s_0)$ is transverse to the interior of $F$.
	But there are also only a countable number of faces of $sd(X)$, and so for almost every $s_0 \in S$, $H|_G(-,s_0)$ is transverse to the interior of $F$ at every point in the interior of $G$, for every $G$ and every $F$.
	Lastly, $H(-,s_0)$ is homotopic to $\Phi$ via $H(-, ts_0)$.

	To complete the proof we must do one last thing: we must fine tune $\rho$ to ensure that in forming $H(-,s_0)$ to obtain the required transversality we have not pushed any $F^\vee$ so far as to create new intersections with $\dim(F)$ faces of $X$ beyond the single desired intersection of $H(F^\vee,s_0)$ with $F$ that was obtained earlier using the map $\Phi$.
	For this we do the following.

	Recall that $H(-,s_0) = \Phi$ on $W_3$, and that $\Phi$ already creates the desired intersection between each $F$ and a point in the image of $W_4$.
	As noted above, if $F$ is any face of $X$ then the distance between $\Phi(F^\vee-(F^\vee \cap W_4))$ and the union of $\dim(F)$-faces of $X$ is positive.
	Now suppose $x \in M-W_4$ and consider a compact neighborhood $\bar U_x$ of $x$ in $M-W_4$.
	By the above, if $\bar U_x \cap F^\vee \neq \emptyset$, then there is a positive distance between $\Phi(\bar U_x \cap F^\vee)$ and the $\dim(F)$-skeleton of the cubulation.
	So using the Tube Lemma as in the proof of \cref{P: ball stability}, there is an $\varepsilon_{x,F} > 0$ such that $H(z,s)$ does not intersect the $\dim(F)$-skeleton of the cubulation for any $z \in \bar U_x \cap F^\vee$ and any $s$ with $|s| \leq \varepsilon_{x,F}$.
	As $\bar U_x$ is compact, it intersects only a finite number of $F^\vee$ as $F$ ranges over all faces of $X$, and we let $$\varepsilon_x = \min\{\varepsilon_{x,F} \mid \bar U_x \cap F^\vee\neq \emptyset\}.$$
	So, by construction, any map $g \colon M \to M$ such that $d(z,g(z))<\varepsilon_x$ for all $z \in \bar U_x$ satisfies the property that if $z \in F^\vee \cap \bar U_x$ then $g(z)$ is not contained in a $\dim(F)$ face of $X$.

	Now, suppose we have constructed such an $\varepsilon_x$ for all $x$ in $M - W_4$.
	Then the interiors $U_x$ of the $\bar U_x$ cover $M - W_4$, and since $M - W_4$ is a subspace of the Euclidean metric space $\R^N$, we can take a locally finite refinement $\mc U$.
	By \cref{L: minimizer}, we can find a smooth function $\rho_1 \colon M - W_4 \to \R$ so that $0<\rho_1(z)<\varepsilon_x$ if $z \in \bar U_x$ for $U_x \in \mc U$.
	Let $\rho_2 = \rho\rho_1 \colon M \to \R$, where $\rho$ is as chosen above.
	This is smooth and well defined on all of $M$ as $\rho(x) = 0$ for all $x \in W_3$ and $W_4 \subset int(W_3)$.
	We also see that $0<\rho_2(x)<\rho_1(x)$ for all $x \in M - W_4$.
	So now if we replace $H$ with $H_2 \colon M \times S \to M$ defined by $H_2(x,s) = \pi(\Phi(x)+\rho_2(x)\epsilon(\Phi(x)) s)$, then for any $s_0 \in S$ and for any $x \in (M - W_4) \cap F^\vee$, we have that $H_2(x,s_0)$ does not intersect any face of dimension $\dim(F)$ as desired.
	We can now set $h(-) = H_2(-,s_0)$ for almost all $s_0$ to achieve all three required conditions, noting that the conclusions of the preceding paragraphs did not depend on the choice of $\rho$.

	Lastly, we observe that the homotopies from the identity to $\Phi$ and from $\Phi$ to $h$ (given by $H_2(-,ts_0)$) can be taken to be proper by an application of \cref{L: nearby proper homotopy} as in the proof of \cref{P: ball stability}, assuming we make all of our homotopies above yet smaller if necessary to meet the requirements of \cref{L: nearby proper homotopy}.
\end{proof}
\index{cubulation!pushing the dual|)}

\subsubsection{The intersection map is an isomorphism for finitely-generated cohomology groups}

We are nearly ready to prove \cref{T: intersection qi}, which we restate here:

\begin{theorem*}[\cref{T: intersection qi}]\index{intersection map!is a cohomology surjection}\index{intersection map!is a cohomology isomorphism when finitely generated}
	If $M$ is a manifold without boundary cubulated by $X$, the intersection map $\mc I \colon H^i_{\Gamma \pf X}(M) \to H^i(X)$ is a surjection.
	If, in addition, $H^i(M)$ is finitely generated, then $\mc I$ is an isomorphism.
\end{theorem*}

The last piece is a basic fact about finitely-generated abelian groups:

\begin{lemma}\label{L: surjection of isomorphic groups}
	Every surjection of isomorphic finitely-generated abelian groups is an isomorphism.
\end{lemma}
\begin{proof}
	First, let $A$ be a finitely-generated abelian group, and let $f \colon A \to A$ be a surjection.
	Let $f^p$ denote the $p$-th iteration of $f$.
	We have $\ker(f)\subseteq \ker (f^2) \subseteq \cdots$, but since $A$ is a finitely-generated module over the Noetherian ring $\Z$, this sequence stabilizes; say $\ker(f^n)=\ker(f^{n+1})$ for some $n$.
	Suppose $x \in A$ with $f(x) = 0$. As $f$ is surjective, so is each $f^p$ and thus $x = f^n(y)$ for some $y \in A$.
	Then $0 = f(x) = f(f^n(y)) = f^{n+1}(y)$.
	But $\ker(f^n) = \ker(f^{n+1})$, so $x= f^n(y) = 0$.
	Thus $f$ is also injective and so an isomorphism.

	If instead $f: A \to B$ is a surjection of isomorphic abelian groups with $\phi: B \to A$ the isomorphism, then the above applies to $\phi f \colon A \to A$, and so $\phi^{-1}(\phi f) = f$ is a composition of isomorphisms and so an isomorphism.
\end{proof}

At last we prove \cref{T: intersection qi}.

\begin{proof}[Proof of \cref{T: intersection qi}]
	In \cref{S: dualization of complexes} we constructed a chain map $\Psi \colon K^*(X) \to C_\Gamma^*(M)$ by taking each $F^* \in K^*(X)$ to the geometric cochain represented by the inclusion of $F^\vee$ into $M$, but the image did not lie in $C^*_{\Gamma \pf X}(M)$.
	To remedy this, we define $\psi \colon K^*(X) \to C_{\Gamma \pf X}^*(M)$ as the composition of $\Psi$ with the chain map on $C^*_{\Gamma}(M)$ determined by the map $h$ constructed in \cref{L: push dual}.
	This last is a chain map by \cref{C: proper cofunctoriality}, using \cref{D: homotopy co-orientation} to co-orient $h$ with the co-orientation induced by $h$ being homotopic to the identity map of $M$, which we give its tautological co-orientation.
	Furthermore, for each face $F^*$ of $K^*(X)$, we have by construction that $\psi(F^*)$ can be represented by a precochain that is transverse to the cubulation.
	So $\psi$ has image in $C_{\Gamma \pf X}^*(M)$ as desired.

	By \cref{L: push dual} and the definitions, we have for each pair of faces $F$ and $f$ of $X$ of the same dimension that
	\begin{equation*}
		(\mc I \psi(F^*))(f) = (\mc I h \Psi(F^*))(f)\\
		= (\mc I h (F^\vee))(f)\\
		= I( h (F^\vee) , f) \\
		=
		\begin{cases}
			1, & F=f,\\
			0, & \text{otherwise.}
		\end{cases}
	\end{equation*}
	So the composition $\mc I \psi$ is the identity on $K^*(X)$, and it follows that $\mc I$ induces a cohomology surjection $H^*_{\Gamma \pf X}(M)\onto H^*(X)$.

	We know from \cref{T: geometric is singular} and the isomorphisms among cubical, singular cubical, and singular simplicial cohomologies that $H^i_{\Gamma \pf X}(M)$ and $H^i(X)$ are both isomorphic to the singular simplicial cohomology group $H^i(M)$.
	So the last claim holds by \cref{L: surjection of isomorphic groups}.
\end{proof}

We conjecture that $\mc I$ is an isomorphism in the general case but have not been able to prove this.

\section{Products of geometric chains and cochains}\label{S: products}

In this section we consider products of geometric chains and cochains, first simply as chains and cochains and then as pairings on homology and cohomology.
These pairings are all built from the fiber products and exterior products of maps as defined in \cref{S: orientations and co-orientations}.
However, while the exterior products were fully defined, the fiber products required transversality of $f \colon V \to M$ and $g \colon W \to M$ in order for $V \times_M W$ to be a well-defined manifold with corners possessing an oriented or co-oriented map to $M$.
Consequently, the fiber products do not induce fully-defined chain- and cochain-level products such as a pairing $C^*_\Gamma(M) \otimes C^*_\Gamma(M) \to C^*_\Gamma(M)$.
At best we can hope for a partially-defined (co)chain-level pairing, though even this is not clear once we take into account that a geometric chain or cochain is not represented simply by a single isomorphism class of a map $V \to M$ but is rather an equivalence class of such mappings up to triviality and degeneracy.

In \cref{S: chain products}, we address this issue and show that there is a natural notion of transversality among chains and cochains, despite the ambiguity in the representative prechains and precochains.
We use this to provide partially-defined cup, cap, and intersection pairings among geometric chains and cochains.
We consider it important to have such pairings, even when only partially defined, as cochain algebras contain much information that is lost on passage to cohomology.
For example, the singular cochains of a space are what carry the $E_\infty$-algebra structure, while passage to cohomology often contains just shadows of this structure, such as the Steenrod squares.
We also provide fully-defined chain and cochain exterior products, though these give us less trouble.
In \cref{S: (co)chain properties}, we collect the various properties of these partially-defined products, mostly based on properties we have established for fiber products of maps in earlier sections.

In \cref{S: homology products}, we then turn to the resulting products in geometric homology and cohomology, which we show are fully defined, providing cup, cap, intersection, and exterior products.
In \cref{S: kroneker}, we utilize this machinery along with our cubulations to construct a Universal Coefficient Theorem when the cohomology is finitely generated.
In \cref{S: usual cup}, we show that the cup product in geometric cohomology is abstractly isomorphic to the singular cohomology cup product, which we use to prove a K\"unneth Theorem in \cref{S: kunneth}.
Finally, in \cref{S: cubical cup and cap}, we relate geometric (co)homology cup and cap products to their singular (co)homology counterparts more concretely, using cubulations and the intersection maps introduced in \cref{S: intersection map}.
This includes some applications to Poincar\'e Duality and umkehr maps.

\subsection{Chain- and cochain-level products and transversality}\label{S: chain products}

In this section, we develop chain- and cochain-level products, as well as study some other aspects of chain- and cochain-level transversality.
We begin in \cref{S: simple products} with the case of ``simple products,'' in which two chains or cochains can be represented by transverse prechains or precochains, in which case we can take the fiber product.
Building on this case, we then define a more general notion of transversality for chains and cochains that allows for some amount of bilinear behavior.
After introducing some more transversality results in \cref{S: homology products} in the context of considering (co)homology products, we will return to this theme in \cref{S: product pullbacks}, where we consider pullbacks of cochains, and in \cref{S: Kronecker}, where we consider a Kronecker pairing for geometric chains and cochains.
In \cref{S: exterior chain products}, we define exterior products of chains and cochains.

\subsubsection{Transversality and products}\label{S: simple products}
In this section we define a notion of transversality of geometric chains and cochains, which allows us to define (co)chain-level cup, cap, and intersection products.

We begin in \cref{D: cochain trans} with the naive case in which our (co)chains possess transverse representing pre(co)chains.
We call such (co)chains \textbf{simply transverse} (see \cref{D: cochain trans}), and this is sufficient to define cup, cap, and intersection products via fiber products.
However, as we will discuss below, this definition is insufficient to obtain a product that behaves bilinearly.
To obtain this property, we introduce the more general concept of compound transversality in \cref{D: multicup}.
This allows for bilinear behavior of products of simply transverse (co)chains, although compound transversality itself does not seem to provide a satisfactory bilinear product.
The key problem has to do with demonstrating the existence of products that are independent of choices of representing pre(co)chains.
We will discuss the difficulties further after establishing some definitions and results.

We begin with simple transversality.

\begin{definition}\label{D: cochain trans}\index{transversality!simple|textbf}\index{geometric cochain!simply transverse|textbf}
	Let $M$ be a manifold without boundary.
	We say that $\uV, \uW \in C^*_{\Gamma}(M)$ are \textbf{simply transverse} as geometric cochains if there exist representatives $V,W \in PC^*_\Gamma(M)$ such that $V$ and $W$ are transverse as manifolds with corners mapping to $M$.
	We call the data of such a pair $(V,W)$ a \textbf{simple transverse representation}\index{simple transverse representation} for the pair $(\uV,\uW)$.

	We define simple transversality similarly if $\uV \in C^*_{\Gamma}(M)$ and $\uW \in C_*^{\Gamma}(M)$ or if $M$ is oriented and $\uV, \uW \in C_*^{\Gamma}(M)$.
\end{definition}

\begin{definition}\label{D: cochain products}
	Let $M$ be a manifold without boundary.
	For $\uV, \uW \in C^*_{\Gamma}(M)$ simply transverse, we define the \textbf{cup product}\index{cup product|textbf}\index{geometric cochain!cup product} $\uV \uplus \uW \in C^*_\Gamma(M)$ to be the geometric cochain represented by the fiber product $V \times_M W$ for some simple transverse representation $(V,W)$ of $(\uV,\uW)$.

	Analogously, if $\uV \in C^*_{\Gamma}(M)$ and $\uW \in C_*^{\Gamma}(M)$ are simply transverse, we define the \textbf{cap product}\index{cap product|textbf}\index{geometric cochain!cap product}\index{geometric chain!cap product} $\uV \nplus \uW \in C_*^\Gamma(M)$ by $V \times_M W$ for some simple transverse representation $(V,W)$.

	If $M$ is oriented and $\uV,\uW \in C_*^\Gamma(M)$ are simply transverse, we define the \textbf{intersection product}\index{intersection product|textbf}\index{geometric chain!intersection product} $\uV \bullet \uW \in C_*^\Gamma(M)$ by $V \times_M W$ for some simple transverse representation $(V,W)$.

	In each context, the given product $V \times_M W$ is as defined in \cref{D: PC products}, as $V$ and $W$ are transverse by assumption.
\end{definition}

The first main result of this section is that these products are well defined as operations on simply transverse geometric chains or cochains, independent of the prechain or precochain representatives chosen.
This is not immediately clear, as a geometric (co)chain $\uV$ has in general an infinite number of representatives in $PC(M)$ that may or may not be transverse to any other given element of $PC(M)$; see \cref{E: bad transversality}.

\begin{theorem}\label{T: cochain product}\index{cup product|textbf}\index{geometric cochain!cup product}\index{cap product|textbf}\index{geometric cochain!cap product}\index{geometric chain!cap product}\index{intersection product|textbf}\index{geometric chain!intersection product}\index{simple transverse representation}
	Given simply transverse $\uV$ and $\uW$, the cup, cap, or intersection products of \cref{D: cochain products} are well defined, independent of choice of simple transverse representation.
\end{theorem}

Rather than prove this theorem here, we first provide some further discussion and development.
\cref{T: cochain product} will then follow directly as a special case of \cref{T: multicup}, below.

To motivate our next step, suppose now a geometric cochain $\uV$ is simply transverse to two other geometric cochains of the same degree, $\underline{W_1}$ and $\underline{W_2}$.
This means we can form $\uV\uplus \underline{W_1}+\uV\uplus \underline{W_2}$, and we would like for this to equal $\uV\uplus (\underline{W_1} + \underline{W_2})$.
The problem is that it is not apparent from the definitions whether or not $\uV$ is simply transverse to $\uW_1+\uW_2$, as the simple transversality of the pairs $(\uV,\underline{W_1})$ and $(\uV,\underline{W_2})$ might be realized by representatives $V_1,V_2, W_1,W_2 \in PC^*_\Gamma(M)$ with $\underline{V_1} = \underline{V_2}$ so that $V_1$ and $W_1$ are transverse as spaces mapping to $M$ and $V_2$ and $W_2$ are transverse as spaces mapping to $M$, but neither $V_1$ nor $V_2$ is transverse to $W_1 \sqcup W_2$.
It is also not apparent how to find a $V_3$ representing $\uV$ that is transverse to both $W_1$ and $W_2$.
The simplest solution would then seem to be to just define
$\uV\uplus (\underline{W_1}+\underline{W_2})$ to be represented by $(V_1 \times_M W_1) \sqcup (V_2 \times_M W_2)$, so long as this is well defined.
The next definition builds on this idea.

\begin{definition}\label{D: multicup}
	Let $M$ be a manifold without boundary.
	Let $\uV, \uW \in C_\Gamma^*(M)$, and suppose $\uV$ and $\uW$ can be written as finite sums $\uV = \sum_i \underline{V_i} \in C_\Gamma^*(M)$ and $\uW = \sum_j \underline{W_j} \in C_\Gamma^*(M)$ such that each pair $(\underline{V_i},\underline{W_j})$ is simply transverse.
	Then we say that $\uV$ and $\uW$ are \textbf{compound transverse}\index{transversality!compound} and define the cup product $\uV\uplus\uW$ as $$\uV\uplus\uW = \sum_{i,j} \underline{V_i}\uplus \underline{W_j},$$
	where the cup products on the right are those of \cref{D: cochain products}.

	We extend the definition of the cap and intersection products analogously.
	In particular, for each simply transverse pair $(\underline{V_i},\underline{W_j})$ as above, there is a simple transverse representation $(V_{ij}, W_{ji})$, and the product $\uV\uplus \uW$, $\uV\nplus \uW$, or $\uV\bullet\uW$ is represented by $\sum_{ij}V_{ij} \times_M W_{ji}$.
\end{definition}

We now demonstrate these products are well defined, noting that \cref{T: cochain product} occurs as the special case in which the sums for $\uV$ and $\uW$ in \cref{D: multicup} each have only one term.
The proof utilizes the essential decompositions of prechains and precochains developed in \cref{S: essential decomp}.

\begin{theorem}\label{T: multicup}
	The products of compound transverse chains and cochains of \cref{D: multicup} are well defined.
	In particular, they do not depend on the decompositions of $\uV$ and $\uW$ into sums of geometric chains or cochains.
\end{theorem}

\begin{proof}
	We provide the argument for the cup product, the other proofs being analogous.
	Suppose $\uV = \sum_i \underline{V_i} = \sum_k \underline{V'_k}$ and $\uW = \sum_j \underline{W_j} = \sum_\ell \underline{W'_\ell}$.
	Suppose the pairs $(\underline{V_i},\underline{W_j})$ and $(\underline{V'_k},\underline{W'_\ell})$ are simply transverse.
	We must show that $\sum_{i,j} \underline{V_i}\uplus \underline{W_j} = \sum_{k,\ell} \underline{V'_k}\uplus \underline{W'_\ell}$.
	The assumptions mean that for each pair $(\underline{V_i},\underline{W_j})$, there are simple transverse representatives we can choose and call $(V_{ij}, W_{ji})$, and similarly for the primed versions.
	Then we must show that $$\left(\bigsqcup_{i,j} V_{ij} \times_M W_{ji}\right) \sqcup \left(-\bigsqcup_{k,\ell} V'_{k\ell} \times_M W'_{\ell k}\right) \in Q^*(M).$$

	For each $V_{ij}$, we have its essential decomposition\index{essential decomposition} $$V_{ij} = V_{ij,E} \sqcup V_{ij,TI} \sqcup V_{ij,NI},$$
	and by \cref{T: minimal rep}, the cochain $\underline{V_{ij}}$ is also represented by a precochain of the form $Z_{ij} \sqcup V_{ij,NI}$, where $Z_{ij}$ is the minimal essential precochain\index{minimial essential prechain or precochain|textbf} of $\underline{V_{ij}}$.
	As $Z_{ij} \sqcup V_{ij,NI}$ is obtained from $V_{ij}$ by removing some components, we may assume that all $V_{ij}$ in fact have this form without disturbing our transversality assumptions, and similarly for the $W_{ji}$, $V'_{k\ell}$, and $W'_{\ell k}$.
	Furthermore, as $V_{ia}$ and $V_{ib}$ represent the same $\underline{V_i}$, they will have the same minimal essential cochain, which we can therefore write simply as $V_{i,E}$ rather than $Z_{ij}$, and again similarly for the other precochains.

	We have $\left(\bigsqcup_j W_{j}\right) \sqcup -\left(\bigsqcup_\ell W'_{\ell}\right) \in Q^*(M)$ by assumption, so by \cref{L: Q essential}, $\left(\bigsqcup_j W_{j,E}\right) \sqcup -\left(\bigsqcup_\ell W'_{\ell,E}\right)$ must be trivial.
	Therefore, by \cref{L: trivial structure}, each connected component, say $\mc W$, appearing in one of the $W_{j,E}$ or $W'_{\ell, E}$ either has a co-orientation-reversing automorphism or appears zero times in all of $\left(\bigsqcup_j W_{j,E}\right) \sqcup -\left(\bigsqcup_\ell W'_{\ell,E}\right)$ when counting with sign.
	If $\mc W$ has a co-orientation-reversing automorphism, then $\mc W$ is trivial and so cannot appear in $W_{j,E}$ or $W'_{\ell,E}$.
	So, for each occurrence of $\mc W$ in $\left(\bigsqcup_j W_{j,E}\right) \sqcup -\left(\bigsqcup_\ell W'_{\ell,E}\right)$, there is an occurrence of $-\mc W$.
	Suppose $\mc W \subset W_{a,E}$ and $-\mc W \subset W_{b,E}$.
	Then, in particular, $\mc W \subset W_{ai}$ and $-\mc W \subset W_{bi}$ for all $i$, and so $\mc W$ is simply transverse to each $V_{ia}$ and $V_{ib}$.
	So $\bigsqcup_i V_{ia}$ and $\bigsqcup_i V_{ib}$ both represent $\uV \in C^*_\Gamma(M)$ and are transverse to $\mc W$.
	Thus
	$$\left(\bigsqcup_i V_{ia} \times_M \mc W\right) \sqcup \left(\bigsqcup_i V_{ib} \times_M -\mc W\right) = \left(\bigsqcup_i V_{ia} \times_M \mc W\right) \sqcup \left(-\bigsqcup_i V_{ib} \times_M \mc W\right) \in Q^*(M)$$
	by \cref{L: pullback with Q}, and similarly if one or both occurrences of $\pm \mc W$ are components of one of the $W'_{\ell, E}$, in which case the corresponding $V_{ia}$ or $V_{ib}$ is replaced with $V'_{k\ell}$ or $V'_{k\ell}$ .
	Continuing in this way, all of $$\left(\bigsqcup_{i,j} V_{ij} \times W_{j,E}\right) \sqcup \left(-\bigsqcup_{k,\ell} V'_{k\ell} \times_M W'_{\ell,E}\right)$$
	is in $Q^*(M)$.

	So it remains to show that
	\begin{equation}\label{E: multicup NI}
		\left(\bigsqcup_{i,j} V_{ij} \times W_{ji,NI} \right) \sqcup \left(- \bigsqcup_{k,\ell} V'_{k\ell} \times_M W'_{\ell k,NI}\right)
	\end{equation}
	is in $Q^*(M)$.
	As each $W_{ji,NI}$ and $W'_{\ell k,NI}$ has small rank, each component of \eqref{E: multicup NI} is of small rank by \cref{L: pullback with Q}.
	So it suffices to show that the boundary of \eqref{E: multicup NI} is a union of trivial and small rank precochains.
	The boundary terms of the form $(\bd V_{ij}) \times W_{ji,NI}$ and $(\bd V'_{k\ell}) \times_M W'_{\ell k,NI}$ all have small rank, again by \cref{L: pullback with Q}.
	So we consider
	$$\left(\bigsqcup_{i,j} V_{ij} \times \bd W_{ji,NI} \right) \sqcup \left(-\bigsqcup_{k,\ell} V'_{k\ell} \times_M \bd W'_{\ell k,NI}\right)$$
	(we can ignore the sign from the boundary formula, as all terms are multiplied by the same sign $(-1)^{m-v}$ in taking the boundary).

	We now consider the essential decompositions of the $\bd W_{ji,NI}$ and $\bd W'_{\ell k,NI}$.
	By \cref{L: pullback with Q}, any fiber product involving a TI component will be trivial and any fiber product involving an NI component will have small rank.
	So we must consider the terms $V_{ij} \times \left(\bd W_{ji,NI}\right)_E$ and $V'_{ij} \times \left(\bd W'_{ji,NI}\right)_E$.
	By \cref{L: same NI}, since $W_{ji}$ and $W_{ja}$ represent the same cochain for any $i,a$, we have that $\bd W_{ji,NI}$ and $\bd W_{ja,NI}$ represent the same cochain (and similar for the $W'$).
	So by \cref{C: Q essential}, these have the same minimal essential part and $\left(\bd W_{j,NI}\right)_E$ is the minimal essential part together with something trivial.
	As any fiber product with something trivial is trivial, we can concentrate on the minimal essential part; we now abuse notation and let $\left(\bd W_{j,NI}\right)_E$ stand just for the minimal essential part.
	So if we can show that
	$$\left(\bigsqcup_j \bd W_{j,NI}\right) \sqcup -\left(\bigsqcup_\ell \bd W'_{\ell,NI}\right)$$
	is in $Q^*(M)$ then we can proceed by the same argument we used above for $\left(\bigsqcup_j W_{j,E}\right) \sqcup -\left(\bigsqcup_\ell W'_{\ell,E}\right)$ to show that $\left(\bigsqcup_{i,j} V_{ij} \times W_{j,E}\right) \sqcup \left(-\bigsqcup_{k,\ell} V'_{k\ell} \times_M W'_{\ell,E}\right) \in Q^*(M)$.

	But, again, we know that
	$\left(\bigsqcup_j W_{j}\right) \sqcup \left(-\bigsqcup_\ell W'_{\ell}\right) \in Q^*(M)$, so by \cref{L: Lipy12}
	$$\left(\bigsqcup_j \left(W_{j,E} \sqcup W_{j,NI}\right)\right) \sqcup \left(-\bigsqcup_\ell \left(W'_{\ell,E} \sqcup W'_{\ell,NI}\right)\right) \in Q^*(M),$$
	as the $W_{j,TI}$ and $W'_{\ell, TI}$ are in $Q^*(M)$.
	Thus by \cref{L: Q essential}, $$\left(\bigsqcup_j W_{j,NI}\right) \sqcup \left(-\bigsqcup_\ell W'_{\ell,NI}\right) \in Q^*(M).$$
	So by \cref{L: bd defined}, the boundary
	$$\left(\bigsqcup_j \bd W_{j,NI}\right) \sqcup \left(-\bigsqcup_\ell \bd W'_{\ell,NI}\right)$$ is in $Q^*(M)$, as required.
\end{proof}

\cref{T: multicup} can be read to imply bilinear behavior of \textit{simply} transverse chains and cochains, although this is somewhat definitional: \cref{D: multicup} \textit{defines} $\uV \uplus (\underline{W_1} + \underline{W_2})$ as  $(\uV \uplus \underline{W_1}) + (\uV \uplus \underline{W_2})$, where the cup products on the right are those of \cref{D: cochain products}, and then \cref{T: multicup} tells us this is well defined.

The conundrum is that compound transverse chains and cochains again do not necessarily behave bilinearly, for essentially the same reasons that simple transversality does not provide bilinear products without extending to a broader definition of transversality:
The definition of compound transversality assumes fixed decompositions $\uV = \sum_i \underline{V_i}$ and $\uW = \sum_j \underline{W_j}$.
But now suppose we have $\uV \in C_\Gamma^*(M)$ that is \textit{compound} transverse to both $\underline{W_1}, \underline{W_2} \in C_\Gamma^*(M)$.
Again we would like $\uV \uplus (\underline{W_1} + \underline{W_2})  = (\uV \uplus \underline{W_1}) + (\uV + \underline{W_2})$, but $\uV \uplus \underline{W_1}$ and $\uV \uplus \underline{W_2}$ might be computed by the formula above using two different ways of writing $\uV$ as a sum, and \cref{T: multicup} does not cover this situation.

Conjecturally, we could mirror the preceding program by introducing an appropriate even broader notion of transversality, showing it is well defined, and then declaring it to provide a notion of bilinearity for products of compound transverse chains and cochains.
Even more conjecturally, this process can be repeated to all levels of ``$n$-transversality'' (with $n=1$ being simple transversality, $n=2$ being compound transversality, etc.) to provide well defined ``level $n$ products'' that definitionally provide bilinearity at the preceding the level.
But we will not pursue this project here.

\subsubsection{Exterior products}\label{S: exterior chain products}

As a complement to \cref{S: simple products}, we observe in this section that the exterior products defined in \cref{S: exterior products} give rise to well-defined products for geometric chains and cochains.
Unlike the cup, cap, and intersection products, these are fully defined, as the exterior products do not require any transversality assumptions.

\begin{definition}\label{D: exterior chain}\index{exterior product}\index{geometric cochain!exterior product}\index{geometric chain!exterior product}
	Suppose $\uV \in C_*^{\Gamma}(M)$ and $\uW \in C_*^{\Gamma}(N)$ are represented by $V \in PC_*^{\Gamma}(M)$ and $W \in PC_*^{\Gamma}(N)$.
	Then we define the \textbf{exterior chain product (or chain cross product)} $$\times \colon C_*^{\Gamma}(M) \otimes C_*^{\Gamma}(N) \to C_*^{\Gamma}(M \times N)$$ by $\uV \times \uW = \underline{V \times W}$, giving the product of oriented manifolds the standard product orientation, as in \cref{S: oriented product}.

	Similarly, suppose $\uV \in C^*_{\Gamma}(M)$ and $\uW \in C^*_{\Gamma}(N)$ are represented by $V \in PC^*_{\Gamma}(M)$ and $W \in PC^*_{\Gamma}(N)$.
	Then we define the \textbf{exterior cochain product (or cochain cross product)}
	$$\times \colon C^*_{\Gamma}(M) \otimes C^*_{\Gamma}(N) \to C^*_{\Gamma}(M \times N)$$ by $\uV \times \uW = \underline{V \times W}$, using the co-orientation of a product of co-oriented maps as defined in \cref{D: co-oriented exterior}.

	These are chain maps by \cref{P: oriented fiber boundary,P: boundary of exterior product}.
\end{definition}

As is standard for singular homology and cohomology, we use the symbol $\times$ for both products, allowing context to determine which product is meant.

\begin{proposition}
	The exterior chain and cochain products are well defined.
\end{proposition}

\begin{proof}
	We note that product of proper maps is proper by \cite[Proposition I.10.1.4]{Bou98}.

	It remains to show that if $V'$ and $W'$ are alternative representatives of $V$ and $W$ then $(V \times W) \sqcup -(V' \times W') \in Q(M \times N)$.
	We will show that $(V \times W) \sqcup -(V' \times W) \in Q(M \times N)$, then the general case follows from an equivalent argument with $W$.
	But we need only observe that $(V \times W) \sqcup -(V' \times W) = (V \sqcup -V') \times W$ and then apply
	\cref{L: exterior Q}.
\end{proof}

\subsection{Properties of the chain and cochain products}\label{S: (co)chain properties}\index{cup product!properties|(}\index{cap product!properties|(}\index{intersection product!properties|(}\index{exterior product!properties|(}

Now that we have defined cup, cap, intersection, and exterior products of geometric chains and cochains and shown that these products are well defined, at least when the necessary transversality and orientation conditions hold, they immediately inherit many of the properties demonstrated in \cref{S: orientations and co-orientations}.
We provide below some tables listing these properties and the locations of the previous results that support them.
The references are typically to results that involve only transversality of a pair of prechains or precochains, but in the chain and cochain setting they generalize to the more general products of \cref{D: multicup} by applying them to each summand.

For example, suppose $\uV, \uW \in C^*_\Gamma(M)$ are compound transverse.
This means we can write $\uV = \sum_i \underline{V_i} \in C_\Gamma^*(M)$, $\uW = \sum_j \underline{W_j} \in C_\Gamma^*(M)$ with each pair $(\underline{V_i},\underline{W_j})$ simply transverse.
And this means that there are representatives $V_{ij},W_{ji} \in PC_\Gamma^*(M)$ such that for all $i$ and $j$, we have $\underline{V_{ij}} = \underline{V_i}$, $\underline{W_{ji}} = \underline{W_j}$, and $V_{ij}$ transverse to $W_{ji}$.
We then have $\uV\uplus \uW$ represented by
$\sum_{i,j} V_{ij} \times_M W_{ji}$.
By \cref{P: graded comm}, we have
$$\sum_{i,j} V_{ij} \times_M W_{ji} = \sum_{i,j} (-1)^{(m-v)(m-w)}W_{ji} \times_M V_{ij} = (-1)^{(m-v)(m-w)}\sum_{i,j} W_{ji} \times_M V_{ij},$$
and this last expression represents $(-1)^{(m-v)(m-w)} \uW\uplus \uV$.
So we obtain the cup product commutativity formula $$\uV\uplus\uW = (-1)^{(m-v)(m-w)}\uW\uplus\uV$$ for compound transverse cochains.

The more complicated exceptions to this inheritance of properties from the pre(co)chain properties concern associativity and naturality, which we will address below in a separate section.

In the tables that follow, we assume to hold all transversality required for each expression to be defined.
For intersection products, we assume that the underlying manifold is oriented.
Unless stated otherwise, our default notations for cup, cap, and intersection products will have manifolds with corners $V$ and $W$ mapping to a manifold without boundary $M$.
Our default notations for chain and cochain cross products will assume $V \to M$ and $W \to N$.
We explain the further assumptions and notations prior to each table of formulas

\subsubsection{Boundary formulas}\index{cup product!boundary formula}\index{cap product!boundary formula}\index{intersection product!boundary formula}\index{exterior product!boundary formulas}

For our first table, with formulas involving boundaries, we also invoke the well-definedness of boundaries of geometric chains and cochains, see \cref{L: co/chains well defined}.
The cup, cap, and intersection products require (simple or compound) transversality of $\uV$ and $\uW$; the exterior products have no transversality requirements.
In the first line, we use that $\uV \times \uW = \uV \times_M \uW$, when $M$ is a point.

\begin{center}
	\begin{tabular}{|l|c|l|}
		\hline
		Chain cross product &$\bd(\uV \times \uW) = (\bd \uV) \times \uW+ (-1)^{v}\uV \times \bd \uW$&\cref{P: oriented fiber boundary}\\
		\hline
		Cochain cross product&$\bd(\uV \times \uW) = (\bd \uV) \times \uW+ (-1)^{m-v}\uV \times \bd \uW$&\cref{P: boundary of exterior product}\\
		\hline
		Cup product&$\bd (\uV \uplus \uW) = (\bd \uV) \uplus \uW+ (-1)^{m-v} \uV \uplus \bd \uW$&\cref{leibniz}\\
		\hline
		Cap product&$\bd(\uV\nplus \uW) = (-1)^{v+w-m} (\bd \uV)\nplus \uW + \uV\nplus\bd \uW$&\cref{P: Leibniz cap}\\
		\hline
		Intersection product &$\bd (\uV \bullet \uW) = (\bd \uV) \bullet \uW + (-1)^{m-v}\uV \bullet \bd \uW$&\cref{P: oriented fiber boundary}\\
		\hline
	\end{tabular}
\end{center}

\subsubsection{Commutativity formulas}\index{cup product!commutativity}\index{cap product!commutativity}\index{intersection product!commutativity}\index{exterior product!commutativity}

For the commutativity properties listed below, $\tau$ is the transposition map $\tau \colon N \times M \to M \times N$.
The cup and intersection products require transversality of $\uV$ and $\uW$; the exterior products have no transversality requirements.
In the first line, we again use that $\uV \times \uW = \uV \times_M \uW$, when $M$ is a point.

\begin{center}
	\begin{tabular}{|l|c|l|}
		\hline
		Chain cross product&$\tau(\uV \times \uW) = (-1)^{vw}\uW \times \uV$&\cref{P: commute oriented fiber}\\
		\hline
		Cochain cross product&$\tau^*(\uV \times \uW) = (-1)^{(m-v)(n-w)}\uW \times \uV$&\cref{P: exterior commutativity}\\
		\hline
		Cup product&$\uV\uplus \uW = (-1)^{(m-v)(m-w)} \uW\uplus \uV$&\cref{P: graded comm}\\
		\hline
		Intersection product&$\uV\bullet \uW = (-1)^{(m-v)(m-w)}\uW\bullet \uV$&\cref{P: commute oriented fiber}\\
		\hline
	\end{tabular}
\end{center}

\subsubsection{Unital properties}\label{S: unital properties}\index{cup product!unitality}\index{cap product!unitality}\index{intersection product!unitality}\index{exterior product!unitality}

For the following unital properties, we write $pt$ to refer to the point with its positive orientation.
We will write $\underline{pt}$ for the geometric chain given by $\id_{pt}:pt \to pt$ or for the geometric cochain given by the canonically co-oriented identity map $\id_{pt}:pt \to pt$.
Similarly, $\uM$ represents the geometric chain or cochain determined by $\id_M \colon M \to M$, canonically co-oriented in the cochain case.
Technically, $M$ must be compact for $\id_M$ to represent a chain, but the corresponding formulas hold more broadly at the referenced locations and so these identities could be taken as statements involving a broader class of geometric chains.
Note that, as a cochain, $\uM \in C^0_\Gamma(M)$, and \cref{P: projection pullbacks} shows that these behave like the singular cochain $1$.
We also let $\pi_1 \colon M \times N \to M$ and $\pi_2 \colon N \times M \to M$ denote the projections.
In the first formula for the cap product with $\uM$, $M$ is assumed oriented, and the first $\uV$ in the formula is represented by $V \to M$ as a cochain while the second $\uV$ in the formula is $V \to M$ as a chain with the induced orientation on $V$; see \cref{P: cap with identity M}.
In the second cap product formula, both instances of $\uV$ are as chains.
As $\id_M$ is transverse to all other maps, the following hold for all $\uV$.

\begin{center}
	\begin{tabular}{|l|c|l|}
		\hline
		Identity for chain cross product&$\uV \times \underline{pt} = \underline{pt} \times \uV = \uV$& Straightforward\\
		\hline
		Identity for cochain cross product&$\uV \times \underline{pt} = \underline{pt} \times \uV = \uV$& \cref{P: co-oriented exterior unit}\\
		\hline
		Cochain cross product with $1$&\begin{tabular}{c}$\pi_1^*\uV = \uV \times \underline{N}$\\$\pi_2^*\uV = \underline{N} \times \uV$ \end{tabular} &\cref{P: projection pullbacks}\\
		\hline
		Cup product with $1$&$\uV\uplus\uM = \uM\uplus \uV = \uV$&\cref{C: cup with identity}\\
		\hline
		Cap product with $\underline M$&$\uV\nplus \uM = \uV$ &\cref{P: cap with identity M}\\
		\hline
		Cap with product with 1&$\uM\nplus \uV = \uV$&\cref{P: cap with 1}\\
		\hline
		Intersection product with $\uM$ &$\uM\bullet \uV = \uV\bullet \uM = \uV$&\cref{P: oriented fiber product basic properties}\\
		\hline
	\end{tabular}
\end{center}

\subsubsection{Mixed properties}\label{S: mixed formulas}\index{cup product!mixed product properties}\index{cap product!mixed product properties}\index{intersection product!mixed product properties}\index{exterior product!mixed product properties}

The next grouping concerns properties that involve multiple kinds of products.
We recall that $\diag \colon M \to M \times M$ is the diagonal map.
For these properties we assume maps $V,W \to M$ and $X,Y \to N$.
We also have projections $\pi_M \colon M \times N \to M$ and $\pi_N \colon M \times N \to N$.
As in \cref{S: unital properties}, $\uM$ represents the geometric chain or cochain determined by $\id_M \colon M \to M$.
Then the last formula follows from \cref{P: compare cup and intersection orientations} by recalling from \cref{P: cap with identity M} that, when $M$ is oriented, the fiber product with $\id_M \colon M \to M$ takes a precochain to the corresponding chain with the induced orientation.
Again, $M$ must be compact for $\id_M$ to represent a chain, but the corresponding formulas hold more broadly at the referenced locations and so these identities could be taken as statements involving a broader class of geometric chains.

The first and last properties require that $\uV$ and $\uW$ be transverse.
The second holds for all $\uV, \underline{X}$.
The next three require that $\uV$ be transverse to $\uW$ and that $\underline{X}$ be transverse to $\underline{Y}$.

\begin{center}\index{cup product!is diagonal pullback of exterior product}\index{exterior product!is cup product of pullback}
	\begin{tabular}{|l|c|l|}
		\hline
		Cup from cross& $\uV\uplus \uW = \diag^*(\uV \times \uW)$&\ref{P: cross to cup}	\\
		\hline
		Cross from cup&$\uV \times \underline{X} = \pi_M^*(\uV)\uplus\pi_N^*(\underline{X})$& \ref{C: cross is cup}	\\
		\hline
		Cup of crosses&$(\uV \times \underline{X})\uplus (\uW \times \underline{Y}) = (-1)^{(m-w)(n-x)} (\uV\uplus \uW) \times (\underline{X}\uplus \underline{Y})$ &	\ref{C: criss cross}\\
		\hline
		Cap of crosses &$(\uV \times \underline{X})\nplus (\uW \times \underline{Y}) = (-1)^{(x+y-n)(m-v)} (\uV \nplus \uW) \times (\underline{X}\nplus \underline{Y})$ & \ref{P: cap cross}\\
		\hline
		Intersect. of crosses &$(\uV \times \underline{X})\bullet (\uW \times \underline{Y}) = (-1)^{(m-w)(n-x)}(\uV\bullet \uW) \times (\underline{X}\bullet \underline{Y})$&\ref{P: oriented interchange}\\
		\hline
		Cup and intersect. &$(\uV\uplus \uW)\nplus \uM = (-1)^{(m-v)(m-w)}(\uV\nplus \uM)\bullet(\uW\nplus \uM)$&\ref{P: compare cup and intersection orientations}\\
		\hline
	\end{tabular}
\end{center}

\subsubsection{Immersion formulas}\index{cup product!of immersions}\index{cap product!of immersions}\index{intersection product!of immersions}\index{exterior product!of immersions}

While geometric chains and cochains do not have unique representatives by prechains and precochains, we recall that if the two terms can be represented by transverse embeddings, then the cup, cap, and intersection products are represented by their intersection and we have simple formulas for the orientations or co-orientations; see \cref{P: normal pullback,P: cap of immersions,P: orient intersection}, respectively.
For cap and intersection products, the special cases where the (co)chains have complementary dimensions are further specified in \cref{C: complementary cap,C: orient complementary intersection}.

\subsubsection{Naturality and associativity formulas}

Formulas for associativity and naturality of geometric chain and cochain products are more delicate than our preceding formulas because they require sufficient transversality of more than two objects.
This would require some careful assumptions even just for maps of manifolds, but the ambiguity in representation of geometric chains and cochains makes the situation even more complicated.
To start, there is the question of which transversality we mean, as we have defined both simple and compound transversality for chains and cochains in \cref{S: simple products}.
If we limited ourselves to simple transversality, then we would be able to invoke results like \cref{C: fiber assoc} fairly directly, though even here there are a number of conditions that must be met.

If we want to work with compound transversality, the situation quickly becomes much more complicated.
For example, if $\uV$ and $\uW$ are compound transverse, then by definition we can write finite sums $\uV = \sum_i \underline{V_i}$ and $\uW = \sum_j \underline{W_j}$ and then find transverse representatives $V_{ij}$ and $W_{ji}$ of $\underline{V_i}$ and $\underline{W_j}$.
Then $\uV\uplus \uW$ is represented by $\sum_{ij}V_{ij} \times_M W_{ji}.$
Now suppose $Z$ is compound transverse to $\uV\uplus \uW$.
Then there must be similar decompositions of $Z$ and $\uV\uplus \uW$ into simply transverse pairs, but it is not clear that this condition can necessarily be written in terms of the $V_{ij}$ and $W_{ji}$ so that we can then take advantage of \cref{C: fiber assoc}.
So, rather than attempt to pursue the most general case, we will impose some extra restrictions in what follows so that we can utilize \cref{C: fiber assoc} and its analogues for the intersection product and for the cap product with a cup product.
These assumptions can be simplified if we only wish to consider simple transversality.

Similar concerns arise for our naturality formulas, as, for example, pulling back a cup product by a map $h$ requires that $h$ be transverse to the cup product, so we again have an interaction of three maps, leading to similar issues.

\textbf{Naturality.}\index{cup product!naturality}\index{cap product!naturality}\index{intersection product!naturality}\index{exterior product!naturality}
As noted above, naturality of cup and cap products requires some extra care to ensure not just that chains and cochains are appropriately transverse but that there are also the appropriate transversalities with respect to the maps we pull back by.
This requires a good number of further assumptions; see \cref{R: multiproducts,P: 3 out of 4 trans}.

So suppose first $h \colon N \to M$ is a map of manifolds without boundaries and that $\uV,\uW \in C^*_{\Gamma}(M)$.
For naturality of cup products we will assume not just that $\uV$ and $\uW$ are compound transverse, but we also require decompositions into finite sums $\uV = \sum_i \underline{V_i}$, $\uW = \sum_j \underline{W_j}$ such that each pair $(\underline{V_i},\underline{W_j})$ has representatives $V_{ij}$ and $W_{ji}$ such that:
\begin{itemize}
	\item $V_{ij}$ and $W_{ji}$ are transverse and
	\item $V_{ij}$, $W_{ji}$, and $V_{ij} \times_M W_{ji}$ are all transverse to $h$.
\end{itemize}
If we assume that $\uV$ and $\uW$ are simply transverse, then we need only representatives $V$ and $W$ such that $V$ and $W$ are transverse and, $V$, $W$, and $V \times_M W$ is transverse to $h$.

For naturality of cap products, we suppose $h \colon N \to M$, $\uV \in C^*_\Gamma(M)$, and $\uW \in C_*^\Gamma(N)$.
Then, via \cref{P: natural cap}, it is sufficient to assume decompositions $\uV = \sum_i \underline{V_i}$ and $\uW = \sum_j \underline{W_j}$ with representative pairs $(V_{ij}, W_{ji})$ such that each $V_{ij}$ is transverse to $h$ and each $W_{ji}$ is transverse to the pullback $V_{ij} \times_M N \to N$.
In the simple version we just need $V$ and $W$ with $V$ transverse $h$ and $W$ transverse to $V \times_M N \to N$.

The exterior products are simpler.
The naturality of the chain cross product requires no assumption, while the naturality of the cochain cross product requires only that $V$ and $W$ have representatives that are respectively transverse to $h$ and $k$.

With these assumptions, we have the following formulas.

\begin{center}
	\begin{tabular}{|l|c|l|}
		\hline
		Chain cross product&$(h \times k)(\uV \times \uW) = h(\uV) \times k(\uW)$ &Straightforward\\
		\hline
		Cochain cross product&$(h \times k)^*(\uV \times \uW) = h^*(\uV) \times k^*(\uW)$ &\cref{P: natural exterior}\\
		\hline
		Cup product &$h^*(\uV\uplus \uW) = h^*(\uV) \uplus h^*(\uW)$&\cref{C: fiber natural pullback}\\
		\hline
		Cap product &$\uV \nplus h(\uW) = h(h^*(\uV)\nplus \uW)$&\cref{P: natural cap}\\
		\hline
	\end{tabular}
\end{center}

\medskip\noindent\textbf{Associativity.}\index{cup product!associativity}\index{cap product!associativity}\index{intersection product!associativity}\index{exterior product!associativity}
For the associativity formulas we add a manifold with corners $X$ either mapping to $M$ for the cup, cap, and intersection products or to a third target manifold $Q$ for the cross products.
Once again, there are no special requirements for the exterior products in the first two formulas below.
For the other associativity formulas, while this might not encompass the most general possibility, in order to ensure associativity in a compound transversality setting, we assume decompositions into finite sums $\uV = \sum_i \underline{V_i}$, $\uW = \sum_j \underline{W_j}$, and $\underline{X} = \sum_k \underline{X_k}$ such that for each triple $(\underline{V_i},\underline{W_j}, \underline{X_k})$ there are representatives $V_{i,jk}$, $W_{j,ik}$, and $X_{k,ij}$ such that the following pairs are transverse: $(V_{i,jk}, W_{j,ik})$, $(W_{j,ik},X_{k,ij})$, $(V_{i,jk} \times_M W_{j,ik},X_{k,ij})$, and $(V_{i,jk}, W_{j,ik} \times_M X_{k,ij})$.
These assumptions will allow us to invoke \cref{C: fiber assoc,P: OC mixed associativity,P: oriented fiber mixed associativity}.
For simple transversality versions, we need only assume $V, W, X$ such that $(V, W)$, $(W,X)$, $(V \times_M W,X)$, and $(V, W \times_M X)$ are transverse pairs (and by \cref{P: 3 out of 4 trans}, the condition on the last pair is redundant).

We leave the reader to formulate associativity for products of larger collections of maps.

\begin{center}
	\begin{tabular}{|l|c|l|}
		\hline
		Chain cross product& $(\uV \times \uW) \times \uX = \uV \times (\uW \times \uX)$&Straightforward\\
		\hline
		Cochain cross product& $(\uV \times \uW) \times \uX = \uV \times (\uW \times \uX)$&\cref{P: exterior associativity}\\
		\hline
		Cup product &$(\uV\uplus \uW)\uplus\uX = \uV\uplus(\uW\uplus X)$&\cref{C: fiber assoc} \\
		\hline
		Cup/cap & $(\uV \uplus \uW)\nplus \uX = \uV\nplus(\uW\nplus\uX)$& \cref{P: OC mixed associativity}\\
		\hline
		Cap/intersection &
		$ (\uV \nplus \uW) \bullet \uX = (-1)^{(m-v)(m-x)}\uV \nplus (\uW \bullet \uX)$ &\cref{C: cap/intersect}\\
		\hline
		Intersection product &
		$(\uV\bullet\uW)\bullet\uX = \uV\bullet(\uW\bullet \uX)$&\cref{P: oriented fiber mixed associativity} and following\\
		\hline
	\end{tabular}
\end{center}

We note that with our assumptions, these triple products exhibit a form of linearity in each variable akin to that discussed in \cref{S: simple products}.

\index{cup product!properties|)}\index{cap product!properties|)}\index{intersection product!properties|)}\index{exterior product!properties|)}

\subsection{Homology and cohomology products}\label{S: homology products}

In this section we observe that the partially-defined cup, cap, and intersection products of geometric chains and cochains give rise to fully-defined products of geometric homology and cohomology classes.
Similarly, we obtain external homology and cohomology products, although this is more evident as external products are already fully defined for geometric chains and cochains.

\begin{theorem}\label{T: (co)homology products}\index{cup product|textbf}\index{geometric cohomology!cup product}\index{cap product|textbf}\index{geometric cohomology!cap product}\index{geometric cohomology!cap product}\index{intersection product|textbf}\index{geometric homology!intersection product}\index{simple transverse representation}
	Let $M$ and $N$ be manifolds without boundary.
	The chain cross product, cochain cross product, cup product, cap product, and, if $M$ is oriented, intersection product induce fully-defined bilinear maps
	\begin{align*}
		\times \colon & H^\Gamma_*(M) \otimes H^\Gamma_*(N) \to H^\Gamma_*(M \times N)\\
		\times \colon & H_\Gamma^*(M) \otimes H_\Gamma^*(N) \to H_\Gamma^*(M \times N)\\
		\uplus \colon & H_\Gamma^*(M) \otimes H_\Gamma^*(M) \to H_\Gamma^*(M)\\
		\nplus \colon & H_\Gamma^*(M) \otimes H^\Gamma_*(M) \to H^\Gamma_*(M)\\
		\bullet \colon & H^\Gamma_*(M) \otimes H^\Gamma_*(M) \to H^\Gamma_*(M).\\
	\end{align*}
\end{theorem}

We will prove \cref{T: (co)homology products} below.
The basic idea, in the case of the last three products, will be to show that we can represent pairs of homology or cohomology classes by simply transverse representatives, which is accomplished by \cref{T: transverse reps}, below, and then form the usual fiber products.
We therefore have the following immediate consequence.

\begin{theorem}\index{cup product!properties}\index{cap product!properties}\index{intersection product!properties}\index{exterior product!properties}
	The homology cross product, cohomology cross product, cup product, cap product, and, if $M$ is oriented, intersection product satisfy the properties enumerated in \cref{S: (co)chain properties}, except for the boundary formulas.
\end{theorem}

\begin{remark}
	For the naturality and associativity properties which require additional assumptions about transversality of representatives, those assumptions can all be met in the (co)homology setting.
	For example, for naturality of the cup product with respect to a map $h \colon N \to M$, once transverse representatives $V$ and $W$ of $\uV$ and $\uW$ have been found, we can use \cref{T: basic trans,R: countable trans} to replace $h$ with a homotopic map transverse to $V$, $W$, and $V \times_M W$, and by \cref{P: cohomology pullback}, we can use this transverse map to pull back our cohomology classes.
	For the cap product, we can again use \cref{T: basic trans} to assume $h$ transverse to $V$ and then \cref{T: transverse reps}, below, to choose a $W$ in the desired homology class transverse to the pullback $V \times_M N \to N$.
	Similarly, for associativity, we use \cref{T: transverse reps,R: countable trans2} to choose representatives $V$, $W$, and $X$ first so that $W$ is transverse to $V$ and then so that $X$ is transverse to $V \times_M W$ and $W$.
\end{remark}

To prove \cref{T: (co)homology products}, we will use the following important theorem.

\begin{theorem}\label{T: transverse reps}\index{geometric cohomology!transverse representative}\index{geometric homology!transverse representative}
Let $r_W \colon W \to M$ be a proper map from a manifold with corners to a manifold without boundary, and let $\uV \in H_*^\Gamma(M)$ or $\uV \in H^*_\Gamma(M)$.
Then there is an $r_V \colon V \to M$ representing $\uV$ such that $V$ is transverse to $W$.

Furthermore, if $V$ is transverse to $W$ and there is a pre(co)chain $Z$ such that $\bd Z \sqcup -V \in Q(M)$, i.e.\ $\uV = 0$ as a homology or cohomology class, then such a $Z$ can be chosen so that it is transverse to $W$.
\end{theorem}

We will prove \cref{T: transverse reps} below in \cref{S: transverse maps}.
For now we use it to prove \cref{T: (co)homology products}.

\begin{proof}[Proof of \cref{T: (co)homology products}]
	For the exterior products, by \cref{S: exterior chain products} we already have fully-defined maps
	\begin{align*}C^\Gamma_*(M) \times C^\Gamma_*(N)& \to C^\Gamma_*(M \times N)\\ C_\Gamma^*(M) \times C_\Gamma^*(N)& \to C_\Gamma^*(M \times N).
	\end{align*}
	These are easily seen to be bilinear and $\Z$-balanced (i.e.\ they satisfy $(a\uV) \times \uW = a(\uV \times \uW) = \uV \times a\uW$ for any $a \in \Z$).
	Moreover, these are chain maps: for the chain cross product this follows from the standard boundary formula for oriented products and our boundary conventions, \cref{Con: oriented boundary}, and for the cochain cross product this follows from \cref{P: boundary of exterior product}, recalling our indexing convention for cochains.
	The existence of the homology and cohomology cross products now follows from standard homological algebra.

	To define the other products, we utilize \cref{T: transverse reps}.
	Given (co)homology classes represented by $r_V \colon V \to M$ and $r_W \colon W \to M$, we can find by \cref{T: transverse reps} a map $r_V' \colon V' \to M$ that represents the same (co)homology class as $V$ and is transverse to $W$.
	We then represent the product by the (oriented or co-oriented) fiber product of $r_V'$ and $r_W$.

	To show that this gives a well-defined (co)homology class, we can suppose that $r_V'' \colon V'' \to M$ is another map transverse to $r_W$ representing the same (co)homology class as $r_V \colon V \to M$.
	Suppose $r_Z:Z \to M$ with $\bd Z \sqcup V' \sqcup -V'' \in Q(M)$, provides the (co)homology.
	We must show that $V' \times_M W$ and $V'' \times_M W$ are (co)homologous.

	Now by \cref{T: transverse reps}, we can assume that $Z$ is also transverse to $W$.
	So we can form $Z \times_M W$, and we have $$\bd (Z \times_M W) = \pm \left( (\bd Z) \times_M W \right) \sqcup \pm \left( Z \times_M \bd W \right)$$ via the appropriate boundary formulas, with the precise signs depending on which kind of product we are considering (see \cref{P: oriented fiber boundary,leibniz,P: Leibniz cap}).
	As $W$ represents a (co)cycle, $\bd W \in Q(M)$ by \cref{R: cycles and boundaries}, so $Z \times_M \bd W \in Q(M)$ by \cref{L: pullback with Q}.
	Meanwhile, we have
	$$((\bd Z) \times_M W) \sqcup (V' \times_M W) \sqcup -(V'' \times_M W) = \left((\bd Z) \sqcup V' \sqcup -V''\right) \times_M W,$$
	which is again in $Q(M)$ by \cref{L: pullback with Q}.
	Altogether then,
	\begin{multline*}
	\bd (Z \times_M W)  \sqcup (V' \times_M W) \sqcup -(V'' \times_M W)\\
	  = \pm \left( (\bd Z) \times_M W \right) \sqcup \pm \left( Z \times_M \bd W \right) \sqcup (V' \times_M W) \sqcup -(V'' \times_M W)\\
	 = \left[\pm Z \times_M \bd W\right] \sqcup \left[ \pm \left((\bd Z) \times_M W \right) \sqcup (V' \times_M W) \sqcup -(V'' \times_M W)\right].
	\end{multline*}
	We need to have $+(\bd Z) \times_M W$ in the bottom line of this computation, but we can obtain that by replacing $Z \times_M W$ with $-Z \times_M W$ if necessary.
	Then the bottom line will be an element of $Q(M)$.

	So $V' \times_M W$ and $V'' \times_M W$ are (co)homologous.
	\qedhere
\end{proof}

\subsubsection{Representing (co)homology classes by transverse maps}\label{S: transverse maps}

In this section, we prove \cref{T: transverse reps}.
The main tool will be \cref{P: perturb transverse to map}, stated below.
This proposition is analogous to \cref{P: ball stability} with the difference being that instead of making a map transverse to the faces of a cubulation we must make a map transverse to another map.
We will explain how to modify the proof of \cref{P: ball stability} to accomplish this.
We change notation slightly from that of \cref{T: transverse reps} to make it more consistent with \cref{P: ball stability}, which we hope will ease comparison of the two results for the reader.

\begin{proposition}\label{P: perturb transverse to map}\index{transversality!Transversality Theorem!universal}
	Suppose $r_V \colon V \to M$ and $r_X \colon X \to M$ are proper maps from manifolds with corners to a manifold without boundary.
	Then there is a proper homotopy $H \colon M \times I \to M$ such that $H(-,0) = \id$ and $H(-,1)r_V \colon V \to M$ is transverse to $r_X$.

	Furthermore, given another proper map from a manifold with corners $r_W \colon W \to M$ that is transverse to $r_X$, we can choose the homotopy $H$ above so that also the resulting universal proper homotopy of $W$ given by $W \times I \xr{r_W \times \id_I} M \times I \xr{H} M$ is transverse to $r_X$.
\end{proposition}

Proving \cref{P: perturb transverse to map} will involve \cref{L: all transversality is wrt embeddings}, which, in our current notation, says that two maps $r_V \colon V \to M$ and $r_X \colon X \to M$ are transverse if and only when we replace $r_X$ with an embedding $e \colon X \into M \times \R^n$ that projects to $r_X$, then $e$ is transverse to $r_V \times \id_{\R^n}$.
So this lemma allows us to replace transversality of arbitrary maps with transversality in which one map is an embedding.
The following lemma says, roughly speaking, that when we construct such an embedding $e$, then $e(X)$ does not run off to infinity in the $\R^n$ factors over compact subsets of $M$.

\begin{lemma}\label{L: compact preimage}
	Let $r_X \colon X \to M$ be a proper map from a manifold with corners to a manifold without boundary, let $\pi_M \colon M \times \R^n \to M$ be the projection, and let $e \colon X \to M \times \R^n$ be an embedding such that $\pi_Me = r_X$.
	Then if $L \subset M$ is compact, there exists a close ball $\bar B^n_L \subset \R^n$ such that $e(r_X^{-1}(L)) \subset L \times \bar B^n_L$.
\end{lemma}

\begin{proof}
	As $L$ is compact and $g$ is proper, $r_X^{-1}(L)$ is compact.
	So $e(r_X^{-1}(L))$ is compact, as is its image under the projection $\pi_{\R^n} \colon M \times \R^n \to \R^n$.
	Let $\bar B^n_L \subset \R^n$ be a closed ball containing this projection.
	Then $\pi_Me(r_X^{-1}(L)) = r_X(r_X^{-1}(L)) \subset L$ and $\pi_{\R^n}e(r_X^{-1}(L)) \subset \bar B^n_L$.
	So $e(r_X^{-1}(L)) \subset L \times \bar B^n_L$.
\end{proof}

We can now use \cref{L: all transversality is wrt embeddings,L: compact preimage} to augment the proof of \cref{P: ball stability} to a proof of \cref{P: perturb transverse to map}.

\begin{proof}[Proof of \cref{P: perturb transverse to map}]
	Suppose $e \colon X \to M \times \R^n$ is an embedding such that $\pi e = r_X$, with $\pi \colon M \times \R^n \to M$ the projection.
	Such an embedding always exists by \cref{C: embed V}.
	By \cref{L: all transversality is wrt embeddings}, it suffices to show that there is a proper  homotopy $H \colon M \times I \to M$ such that
	\begin{enumerate}
		\item $H(-,0) = \id$,
		\item $(H(-,1)r_V) \times \id_{\R^n} \colon V \times \R^n \to M \times \R^n$ is transverse to $e \colon X \to M \times \R^n$, and
		\item $(H (W \times \id_I)) \times \id_{\R^n}: W \times I \times \R^n \to M \times \R^n$ is transverse to $e$.
	\end{enumerate}
	To do so, we will run through the proof of \cref{P: ball stability} again, adapting it to this altered situation and referring back to that proof for some of the details.

	As in the proof of \cref{P: ball stability}, we begin with the case where $W$ is compact.
	We will construct $F \colon M \times D^N \to M$, with $D^N$ the unit ball centered at $0$ in $\R^N$ for some $N$, such that

	\begin{enumerate}
		\item $F(-,0) = \id \colon M \to M$,
		\item for almost all $s \in D^N$ the composition $V \times \R^n \xr{r_V \times \id_{\R^n}} M \times \R^n \xr{F(-,s) \times \id_{\R^n}} M \times \R^n$ is transverse to $e \colon X \to M \times \R^n$,

		\item there is a ball neighborhood $D_r^N$ of $0$ in $D^N$ such that for all $s \in D_r^N$ the composition $W \times \R^n \xr{r_W \times \id_{\R^n}} M \times \R^n \xr{F(-,s) \times \id_{\R^n}}M \times \R^n$ is transverse to $e$.
	\end{enumerate}

	Given such an $F$, we let $s_0$ be any point in $D_r^N$ such that $V \times \R^n \xr{r_V \times \id_{\R^n}} M \times \R^n \xr{F(-,s_0) \times \id_{\R^n}} M \times \R^n$ is transverse to $e \colon X \to M \times \R^n$.
	Then let $H(-,t) = F(-,ts_0)$, i.e.\ $H(y,t) = F(y,ts_0)$.
	Then $H(-,0) = F(-,0) = \id$, and $H(-,1)r_V$ will be transverse to $r_X$ by our choice of $s_0$ and \cref{L: all transversality is wrt embeddings}.
	Finally, $ts_0 \in D_r^N$ for all $t \in I$, each $F(-,ts_0)r_W \times \id_{\R^n}$ is transverse to $e$, so $F(-,ts_0)r_W$ is transverse to $r_X$.
	Each $F(-,ts_0)r_W$ is the restriction of $H \circ (r_W \times \id) \colon W \times I \to M$ to a fixed $W \times t$, so this implies $H \circ (r_W \times \id)$ is transverse to $r_X$.
	This does not provide the properness of $H$, which we will discuss below.

	We now claim that we can construct $F$ almost exactly as in \cref{P: ball stability}.
	Recall that we let $M_\epsilon$ be an $\epsilon$-neighborhood of a proper embedding of $M$ into some $\R^N$ in the sense of the $\epsilon$-Neighborhood Theorem of \cite[Section 2.3]{GuPo74}, with $\epsilon$ a smooth positive function of $M$ and $M_\epsilon = \{z \in \R^N \mid |z-y|<\epsilon(y) \text{ for some }y \in M\}$.
	We may assume $\epsilon$ is bounded by $1$ (see \cite[Exercise 2.3.1]{GuPo74}).
	Let $\pi_\epsilon: M_\epsilon \to M$ be the submersion.
	We define $f: M \times D^N \to M_\epsilon$ by $f(y, s) = y + \eta(y) s$, where $\eta \colon M \to \R$ is a smooth function such that $0 < \eta(y) < \epsilon(y)$ for all $y \in M$.
	As $\eta > 0$, this is clearly a submersion (onto its image) at all points.
	Then we let $F \colon M \times D^N \to M$ be the submersion $M \times D^N \xr{f}M_\epsilon \xr{\pi} M$ and let $\ms F \colon V \times D^N \to M$ be the composition $F \circ (r_V \times \id_{D^N})$.
	It is clear that $F(-,0) = \id_M$,
	Furthermore, the map $\ms F$, as well as each $\ms F|_{S^k(V)}$, is a submersion onto its image.
	Consequently the maps
	$$S^k(V) \times D^N \times \R^n \xr{\ms F \times \id_{\R^n}} M \times \R^n$$
	are also submersions onto their images.
	Thus they are all transverse to all the strata of $X$, embedded by $e$ into $M \times \R^n$.
	It follows by the Transversality Theorem of \cite[Section 2.3]{GuPo74} that for any stratum $S^j(X)$ of $X$, each $\ms F|_{S^k(V)}(-,s) \times \id_{\R^n}$ is transverse to $S^k(X)$ for almost all $s \in D^N$; note that $D^N$ remains our parameter space for invoking the Transversality Theorem, though we no longer write it as the last factor.
	As $X$ and $V$ each have finitely many strata and as the finite union of measure zero sets has measure zero, for almost all $s \in D^N$ we have all $\ms F|_{S^k(V)}(-,s) \times \id_{\R^n}$ transverse to all $S^j(X) \subset M \times \R^n$.
	So $\ms F(-,s) \times \id_{R^n}$ is transverse to $e$ for almost all $s \in D^N$.

	We also observe as in the proof of \cref{P: ball stability} that $H$ can be assumed proper, by replacing $\eta$ with a function with smaller values if necessary.

 	So it remains to show that if $W$ is compact (as we are currently assuming) and $r_W$ is transverse to $r_X$ (or, equivalently, $r_W \times \id_{\R^n}: W \times \id_{R^n} \to M \times \R^n$ transverse to $e$) then $(F(-,s) \circ r_W) \times \id_{\R^n} \colon W \times \R^n \to M \times \R^n$ is transverse to $e$ for all $s$ in some neighborhood $D_r^N$ of $0$ in $D^N$.
	Similarly to \cref{P: ball stability}, it is more convenient here to work with boundaries than with strata, recalling that, by \cref{L: simple trans}, to prove that two maps of manifolds with corners are transverse it is sufficient to show that their compositions with all pairs of boundary inclusions are naively transverse (see \cref{D: naive transversality}).
	Let $\Upsilon_k$ denote the composition $\Upsilon_k \colon \bd^k W \times D^N \xr{i_{\bd^k W} \times \id_{D^N}} W \times D^N \xr{r_W \times \id_{D^N}} M \times D^N \xr{F} M$.
	We must show that there is a $D_r^N$ such that for each $s\in D_r^N$ the maps $\Upsilon_k(-,s) \times \id_{\R^n}$ and $e i_{\bd^jX}$ are naively transverse for all $j,k$.
	We will sometimes provide the details only for $\Upsilon_0$, with $\bd^0 W$ being $W$ itself, when the other cases are analogous.

	We must start with two observations that were not needed in the proof of \cref{P: ball stability}.

	First, as we are currently assuming that $W$ is compact, $L_k = \Upsilon_k(\bd^k W \times \bar D^N_{1/2}) \subset M$ is compact, and so by \cref{L: compact preimage} we have $e(r_X^{-1}(L_k)) \subset M \times \bar B^n_k$ for some closed ball $\bar B^n_k \in \R^n$.
	In particular, this implies that for $|s| \leq 1/2$ only points in the compact set $\bd^k W \times \bar B^n_k$ can be taken by $\Upsilon_k(-,s) \times \id_{\R^n}$ to points of $e(X)$ in $M \times \R^n$.

	Second, we also have to be more careful here about the map $e \colon X \to M \times \R^n$, as its behavior can be more complicated than the embedding of a closed cube of a cubulation, in which each boundary component also embeds.
	Let $X_j = ei_{\bd^jX}(\bd^jX)$, the image of $\bd^jX$ in $M \times \R^n$.
	As the maps $i_{\bd^jX}$ are not necessarily embeddings, it will not generally be the case that $X_j \cong \bd^jX$.
	However, by \cite[Lemma 2.8]{Joy12}, the $i_{\bd^jX}$ are proper maps, so $X_j$ is a closed subset of $M \times \R^n$ (recall that proper maps are closed --- see \cite[Proposition I.10.1.1]{Bou98}).

	Now, suppose\footnote{In this argument we will use the symbol $w$ to refer to points in $W$ or $\bd^k W$, not the dimension of $W$.} $(w,z) \in W \times \bar B^n_0$, and let us fix $X_j$ as above.
	As $r_W$ is transverse to $r_X$ by assumption, $\Upsilon_0(-,0) \times \id_{\R^n} = r_W \times \id_{\R^n}$ is transverse to $e \colon X \to M \times \R^n$.
	So either $(r_W(w),z)\notin X_j$ or $r_W \times \id_{\R^n}$ is naively transverse to $ei_{\bd^jZ}$ at $(w,z)$ (recall \cref{D: naive transversality}).
	In the former case, as $X_j$ is closed, there is an open neighborhood $A_{(w,z)}$ of $(w,0,z) \in W \times D^N \times \R^n$ such that $(\Upsilon_0 \times \id_{\R^n})(A_{(w,z)}) \cap X_j = \emptyset$.
	On the other hand, suppose that $(r_W(w),z) \in X_j$ and $\Upsilon_0(-,0) \times \id_{\R^n}$ is naively transverse to $ei_{\bd^jX}$ at $(w,z)$.
	As $e$ is an embedding, the preimage of $(r_W(w),z)$ in $\bd^jX$ is the preimage of a point of $X$ under the boundary map $i_{\bd^jX}$, which is a finite set of points.
	Let $a \in \bd^jX$ be one point of this preimage.
	As the boundary maps are immersions, there is a neighborhood $C_a$ of $a$ in $\bd^jX$ on which $ei_{\bd^jX}$ restricts to an embedding from an $x-j$ dimensional manifold with corners into $M \times \R^n$.
	In fact, by choosing a chart that takes $0 \subset \R^{x-j}$ to $a$ and using the definition of a smooth map of manifolds with corners, $ei_{\bd^jX}$ (composed with the chart map) extends to a smooth immersion of a neighborhood of $0 \in \R^{x-j}$ into $M \times \R^n$.
	By further appealing to charts and local diffeomorphisms, we can identify a neighborhood of $(r_W(w),z)$ in $M \times \R^n$ with $\R^{m+n}$ and the image of the extension of $ei_{\bd^jX}$ with $\R^{x-j} \times 0 \subset \R^{m+n}$ (cf.\ the Local Immersion Theorem in \cite{GuPo74}).
	Making these identifications, the transversality assumption means that the composition of $D(\Upsilon_0(-,0) \times \id_{\R^n}) \colon T_wW \times T_z\R^n
	\to T_{(r_W(w),z)}(M \times \R^n)$ with the projection to the last $m+n-(x-j)$ coordinates is a linear surjection.
	As this is an open condition on the Jacobian matrix of $\Upsilon_0(-,0) \times \id_{\R^n}$ at $(w,z)$, it follows again that there is an open neighborhood $A_{(w,t),a}$ of $(w,0,z)$ in $W \times D_{1/2}^N \times \R^n$ such that for each $(w',s,z')$ in the neighborhood $\Upsilon_0(-,s) \times \id_{\R^n}$ is transverse to the restriction of $ei_{\bd^jX}$ to a neighborhood of $a$ (cf.\ the Stability Theorem in \cite{GuPo74}).
	As there are a finite number of possible choices for the point $a$ and the transversality assumptions must hold for all of them, by taking $A_{(w,z)} = \cap_a A_{(w,z),a}$ with the finite intersection running over all points of $\bd^jX$ that map to $(r_W(w),z)$, we obtain a neighborhood $A_{(w,z)}$ of $(w,0,z)$ in $W \times D_{1/2}^N \times \R^n$ such that for each $(w',s,z')$ in the neighborhood, $\Upsilon_0(-,s) \times \id_{\R^n} \colon W \times \R^n \to M \times \R^n$ is transverse to $ei_{\bd^jX}$ at $(w',z')$.

	Now, taking the union of the $A_{(w,z)}$ over all $(w,z) \in W \times \bar B^n_0$ gives a neighborhood $G_j$ of $W \times 0 \times \bar B^n_0$ in $W \times D_{1/2}^N \times \bar B^n_0$, and by the Tube Lemma, as $W \times \bar B^n_0$ is compact, there is a neighborhood of $W \times 0 \times \bar B^n_0$ of the form $W \times U_j \times \bar B^n_0$ in $G_j$.
	For each $s \in U_j$, we have $\Upsilon_0(-,s) \times \id_{\R^n} \colon W \times \bar B^n_0 \to M \times \R^n$ naively transverse to $ei_{\bd^jX}$.
	Furthermore, by the choice of $\bar B^n_0$, the map $\Upsilon_0(-,s) \times \id_{\R^n}$ takes no point of $W \times \R^n$ that is in the complement of $W \times \bar B^n_0$ to the image of $X$,
	so in fact $\Upsilon_0(-,s) \times \id_{\R^n} \colon W \times \R^n \to M \times \R^n$ is naively transverse to $ei_{\bd^jX}$.
	Repeating the argument for all of the finite $j$ such that $\bd^jX \neq \emptyset$ and taking $U^0 = \cap_j U_j$, we obtain a neighborhood of $0$ in $D_{1/2}^N$ on which $\Upsilon_0(-,s) \times \id_{\R^n} \colon W \times \R^n \to M \times \R^n$ is naively transverse to all $ei_{\bd^jX}$.
	By equivalent arguments, we can find open sets $U^k$ such that $\Upsilon_k(-,s) \times \id_{\R^n} \colon \bd^k W \times \R^n \to M \times \R^n$ is naively transverse to all $ei_{\bd^jX}$ for all $s \in U^k$ and all $k\geq 0$.
	Finally, taking $r$ sufficiently small so that $D_r^N \subset \cap_k U^k$, we obtain the desired $D_r^N$.

	This completes the proof of the proposition for $W$ compact.

	Next suppose that $W$ is no longer necessarily compact.

	We will continue to utilize $F \colon M \times I \to M$ as defined above, which did not rely on $W$ being compact.
	For $W$ not compact, the first two properties listed above for $F$ will continue to hold, but the third relied on compactness and so need not hold any longer in general.
	However, the only places where we required compactness above were in defining $L_k$ and hence $\bar B^n_k$ via \cref{L: compact preimage}, and then in applying the Tube Lemma.
	So now let $K \subset W$ be compact.
	Then each $i_{\bd^k W}^{-1}(K)$ is also compact as the $i_{\bd^k W}$ are proper.
	Proceeding exactly as in the argument above, using $L_k = \Upsilon_k(i_{\bd^k W}^{-1}(K) \times \bar D^N_{1/2})$, we can find an open neighborhood
	$D_{r_K}^N$ of $0$ in $D^N$ such that for all $s \in D_{r_K}^N$ and any $k\geq 0$ we have that $\Upsilon_k(-,s) \times \id_{\R^n} \colon W \times \R^n \to M \times \R^n$ is naively transverse to all $ei_{\bd^j X}$ at all points of $i_{\bd^k W}^{-1}(K) \times \R^n$.\footnote{We note that this condition still concerns transversality of $\Upsilon_k(-,s) \times \id_{\R^n}$ as a map with domain $\bd^k W \times \R^n$, not $i_{\bd^k W}^{-1}(K) \times \R^n$, which may not be a manifold with corners. But we only consider this condition at points of $i_{\bd^k W}^{-1}(K) \times \R^n$.}

	Let $\{\mc U_\ell\}$ be a locally finite covering of $M$ such that each $\bar{\mc U_\ell}$ is compact.
	As $r_W \colon W \to M$ is proper, each $r^{-1}_W(\bar {\mc U_\ell})$ is compact in $W$.
	Proceeding as just above with $r_W^{-1}(\bar U_\ell)$ in place of $K$, we can find for each $\ell$ an $\varepsilon_{\ell,k} \leq 1/2$ so that for every $s \in D^N_{\varepsilon_{\ell,k}}$ we have $\Upsilon_k(-,s) \times \id_{\R^n}$ naively transverse to all $ei_{\bd^j X}$ at every $(w,z) \in i_{\bd^k W}^{-1} r^{-1}_W(\bar {\mc U_\ell}) \times \R^n$.
	Let $\varepsilon_\ell = \min\{\varepsilon_{\ell,k} \mid k\geq 0\}$.
	These minima exist as $W$ has finite depth.

	Now, using \cref{L: minimizer}, we choose a smooth function $\phi \colon M \to \R$ such that for all $y \in M$ we have $0<\phi(y)<\epsilon_\ell$ if $y \in \bar{\mc U_\ell}$.
	Let $M \times_\phi D^N = \{(y,s) \in M \times D^N \mid |s|<\phi(y)\}$.
	By our construction, for all $k\geq 0$ we have $\Upsilon_k(-,s) \times \id_{\R^n} \colon \bd^k W \times \R^n \to M \times \R^n$ transverse to $ei_{\bd^j X}$ at each $(w,s,z) \in \bd^k W \times D^N \times \R^n$ such that $(w,s) \in (r_Wi_{\bd^{k}W} \times \id_{D^N})^{-1}(M\times_\phi D^N) = \{(w,s) \in \bd^k W \times D^N \mid |s|<\phi(r_Wi_{\bd^k W}(w))\}$.

	We now modify our above constructions as follows.
	Let $\hat f: M \times D^N \to M_\epsilon \subset \R^N$ be given by $\hat f(y, s) = y +\phi(y) \eta(y) s$; as $\phi(y)\eta(y)>0$, this is again a submersion onto its image at all points.
	Let $\hat F \colon M \times D^N \to M$ be the composition $M \times D^N \xr{\hat f} M_\epsilon \xr{\pi_\epsilon} M$, and let $\hat \Upsilon_k$ be the composition $\bd^k W \times D^N \xr{i_{\bd^k W} \times \id} W \times D^N \xr{r_W \times \id} M \times D^N \xr{\hat F} M$ for $k\geq 0$.
	Once again employing the Transversality Theorem of \cite[Section 2.3]{GuPo74} as above, for almost all $s \in D^N$ we have $V \times \R^n \xr{r_V \times \id_{\R^n}} M \times \R^n \xr{\hat F(-,s) \times \id_{\R^n}} M \times \R^n$ transverse to $e \colon X \to M \times \R^n$.
	Letting $s_0$ be any such point we define $\hat H \colon M \times I \to M$ to be $\hat H(y,t) = \hat F(y,ts_0)$, and we claim that this $\hat H$ satisfies the conditions required by the proposition.

	The map $\hat H$ is proper again by \cref{L: nearby proper homotopy} because $\phi(y) \eta(y) \leq \eta(y)$.
	The first two conditions of the proposition follow immediately from the construction and preceding observations.
	Let $$\hat h_k = \hat H \circ (r_W i_{\bd^k W} \times \id_I): \bd^k W \times I \to M.$$
	By \cref{L: all transversality is wrt embeddings}, it remains to verify that each $$\hat h_k \times \id_{\R^n} \colon \bd^k W \times I \times \R^n \to M \times \R^n$$ is naively transverse to each $ei_{\bd^j X}$.
	From here, the argument is essentially the same as the end of the proof of \cref{P: ball stability} with the $\R^n$ factor just along for the ride.
	Again we focus primarily on $k=0$ to simplify notation slightly.

	In detail, for $(w,t) \in W \times I$ we can write $\hat h_0 \colon W \times I \to M$ explicitly as
	$$\hat h_0(w,t) = \pi_\epsilon(r_W(w)+\phi(r_W(w))\eta(r_W(w))ts_0).$$
	So, alternatively, we can observe that $\hat h_0 (w,t)$ is the composition
	\begin{equation}\label{E: perturb transverse to map}
		W \times I \xr{\Phi} W \times I \xhookrightarrow{\Psi} W \times D^N \xr{r_W \times \id} M \times D^N \xr{F} M,
	\end{equation}
	with $\Phi(w,t) = (w,\phi(r_W(w))t)$, $\Psi(w,t) = (w,ts_0)$, and noting that on the right we do mean our original $F$ and not $\hat F$.

	As $0 < \phi(r_W(w)) < 1$ for all $w \in W$, the first map $\Phi$ is a diffeomorphism onto its image, which is a neighborhood of $W \times 0$ in $W \times I$, and the map $\Psi$ embeds this into $W \times D^N$ by a product map that is constant in the $W$ direction and nontrivial linear in the second factor.
	The composition of the last two maps is just our earlier map $\Upsilon_0$.
	By construction, the map $r_W \times \id_{D^N}$ now takes the image of $\Psi\Phi$ into $M\times_\phi D^N$ (as $|s_0|<1$), and so at each point $(w,s,z)$ in the image of $\Psi\Phi \times \id_{\R^n}$ if we fix $s$ and consider $\Upsilon_0(-,s) \times \id_{\R^n}$ we get by construction a map on $W \times \R^n$ that is naively transverse at $(w,s,z)$ to each $ei_{\bd^j X}$.
	Let $(w,t) \in W \times I$, let $\Psi \Phi(w,t) = (\xi,s)$, and let $\R s_0$ denote the line in $\R^N = T_sD^N$ spanned by the position vector of $s_0$.
	As $\Phi$ is a diffeomorphism onto its image and $\Psi$ is an embedding that is the identity with respect to $W$ and nontrivial linear on $I$ for each fixed $w$, we see
	that the derivative of $\Psi\Phi$ maps the tangent space $T_{(w,t)}(W \times I)$ onto $T_\xi W \times \R s_0 \subset T_{(\xi,s)}(W \times D^N)$.
	In particular, this image contains $T_\xi W \times 0$, and we have established that $D(\Upsilon_0 \times \id_{\R^n})$ takes $T_w W \times 0 \times T_z\R^n$ to a tangent subspace in $M \times \R^n$ at $(\hat h_0(w,t),z)$ that is transverse to the tangent space there of each $ei_{\bd^j X}$.

	The same argument holds for each $k>0$ replacing $W$ with $\bd^{k}W$ in \eqref{E: perturb transverse to map} and $r_W$ with $r_{\bd^k W}$.
	So we see that $\hat H$ satisfies all the requirements of the proposition.
	\qedhere
\end{proof}

Finally, we prove \cref{T: transverse reps}.

\begin{proof}[Proof of \cref{T: transverse reps}]
	Given \cref{L: compact preimage}, the proof is very analogous to that of \cref{T: transverse complex}.

	Let $r_W \colon W \to M$ be a proper map from a manifold with corners to a manifold without boundary, and let $\uV \in C_*^\Gamma(M)$ (or $C^*_\Gamma(M)$) be a cycle (or cocycle).
	Let $r_V \colon V \to M$ be any representative for $\uV$.
	By \cref{P: perturb transverse to map} there is a proper homotopy $H \colon M \to I$ such that $H(-,1)r_V$ is transverse to $r_W$.
	We co-orient $H$ and $H(-,1)$ as in \cref{D: homotopy co-orientation} based on the tautological co-orientation of $H(-,0)$.
	By \cref{D: universal homotopy,P: universal homotopy}, $H(-,1)r_V$ represents the same homology or cohomology class as $V$.

	Next, suppose $V$ is transverse to $W$ and there is a pre(co)chain $Z$ such that $\bd Z \sqcup -V \in Q(M)$.
	We consider first the case where $V$ is a precochain.
	By \cref{P: perturb transverse to map}, there is a proper homotopy $H$ such that $H(-,1)r_Z$ is transverse to $r_W$ and $V \times I \xr{r_V \times \id_I} M \times I \xr{H} M$ is transverse to $r_W$.
	Let $Z'$ denote the precochain $Z \xr{H(-,1)r_Z} M$, let $V'$ be the precochain $V \xr{H(-,1)r_V} M$, and let $Y$ be the precochain $V \times I \xr{H\circ (r_V \times \id_I)} M$, co-orienting $r_V \times \id_I$ using the naive pre-homotopy co-orientation convention, \cref{homotopy product co-orientation convention}.
	Let $A$ be the precochain $Z' \sqcup -Y$, which is transverse to $r_W$.
	We note that $\bd Z' \sqcup -V'$ is the image of $\bd Z \sqcup -V$ after composing with $H(-,1)$, and so it is in $Q(M)$ by \cref{L: Q preservation}.
	We also have $\bd Y = V' \sqcup -V \sqcup B$, where $B$ is (up to co-orientation) the precochain $\bd V \times I \xr{H\circ (r_{\bd V} \times \id_I)} M$ (cf. the proof of \cref{C: homotopy}).
	We now compute
	\begin{align*}
		\bd A \sqcup -V &= \bd Z' \sqcup -\bd Y \sqcup -V\\
		&= \bd Z' \sqcup -(V' \sqcup -V \sqcup B) \sqcup -V\\
		&= \bd Z' \sqcup -V' \sqcup V \sqcup -V \sqcup -B.
	\end{align*}
	We have already noted $\bd Z' \sqcup -V' \in Q(M)$ and $V \sqcup -V$ is trivial.
	Since $V$ represents a cycle, $\bd V \in Q^*(M)$ and hence $B \in Q^*(M)$ by \cref{L: dessicated homotopy}.
	So $\bd A \sqcup -V \in Q^*(M)$ and $A$ is transverse to $r_W$.
	So $A$ is our desired replacement for $Z$.

	In the oriented case, we instead have $\bd Y = (-1)^{v} V' \sqcup (-1)^{v+1} V \sqcup B$, identifying $V \times I$ with $V \times_{pt} I$ and applying \cref{P: oriented fiber boundary} (again cf.\ the proof of \cref{C: homotopy}).
	In this case, we let $A = Z' \sqcup (-1)^{v+1} Y$.
	Then
	\begin{align*}
		\bd A \sqcup -V &= \bd Z' \sqcup (-1)^{v+1}\bd Y \sqcup -V\\
		&= \bd Z' \sqcup (-1)^{v+1}((-1)^{v} V' \sqcup (-1)^{v+1} V \sqcup B) \sqcup -V\\
		&= \bd Z' \sqcup -V' \sqcup V \sqcup -V \sqcup (-1)^{v+1}B,
	\end{align*}
	and this is in $Q(M)$ for the same reasons as above.
\end{proof}

\begin{remark}\label{R: countable trans2}
As in \cref{R: countable trans}, and for the same reasons, the proof of \cref{T: transverse reps} allows us to extend the first part of theorem, concerning finding a representative $V$ of $\uV$ that is transverse to $W$, to finding a $V$ that is transverse to each of a countable collection $W_i$.
\end{remark}

While \cref{P: perturb transverse to map} as stated was what we needed to prove \cref{T: transverse reps}, we take the opportunity while in the midst of the machinery to prove the following generalization, which will be useful for other applications, particularly in \cite{GBF47}.

\begin{proposition}\label{P: perturb transverse to map; pairs}\index{transversality!Transversality Theorem!for collection of pairs}
	Suppose $r_V \colon V \to M$ and $r_X \colon X \to M$ are proper maps from manifolds with corners to a manifold without boundary.
	Then there is a proper homotopy $H \colon M \times I \to M$ such that $H(-,0) = \id$ and $H(-,1)r_V \colon V \to M$ is transverse to $r_X$.

	Furthermore, given a finite collection of pairs of proper maps of manifolds with corners $r_{W_i} \colon W_i \to M$ and $r_{Y_i} \colon Y_i \to M$ such that $r_{W_i}$ is transverse to $r_{Y_i}$ for all $i$, we can choose the homotopy $H$ above so that for each $i$ the resulting universal proper homotopy $W_i \times I \xr{r_{W_i} \times \id_I} M \times I \xr{H} M$ is transverse to $r_{Y_i}$ .
\end{proposition}
\begin{proof}
	First consider a single pair $(W,Y)$.
	We see in the part of the proof of \cref{P: perturb transverse to map} involving $W$ that we do not really need the $X$ there to be the same as the $X$ in the first part of the proof.
	If we replace $X$ in this part of the argument with $Y$ and choose the embedding $e$ to be an embedding $e \colon X \sqcup Y \into M \times \R^n$ for appropriate $n$, and if we assume that $W$ is transverse to $Y$ instead of $X$, then in the compact $W$ case the arguments are exactly the same to conclude there is a radius $r$ such that for all $s \in D_r^N$ the composition $W \times \R^n \xr{r_W \times \id_{\R^n}} M \times \R^n \xr{F(-,s) \times \id_{\R^n}}M \times \R^n$ is transverse to the restriction of $e$ to $Y$.
	Constructing $H$ using $s_0 \in D_r^N$, we conclude just as above that $W \times I \xr{r_W \times \id_I} M \times I \xr{H} M$ is transverse to $r_Y$.
	Similarly, if $W$ is not compact, we come to the same conclusions by modifying the argument to choose $\phi$ based on $Y$ instead of $X$.

	Taking this further, next suppose we have any finite collection of pairs $(W_i, Y_i)$ with each $W_i$ compact.
	By starting with an embedding $e \colon X \sqcup (\bigsqcup_i Y_i) \into M \times \R^n$ for appropriate $n$, the preceding paragraph says that we can find $D_{r_i}^N$ such that for all $s \in D_{r_i}^N$ the composition $W_i \times \R^n \xr{r_{W_i} \times \id_{\R^n}} M \times \R^n \xr{F(-,s) \times \id_{\R^n}}M \times \R^n$ is transverse to the restriction of $e$ to $Y_i$.
	Then choosing $r < \min\{r_i\}$, we end up with an $H$ such that that $W_i \times I \xr{r_{W_i} \times \id_I} M \times I \xr{H} M$ is transverse to $r_{Y_i}$ for each $i$, in addition to the original result that $H(-,1)r_V$ is transverse to $r_X$.

	For $W_i$ that are not necessarily compact, we further modify the argument as in the proof of the proposition by using open covers with compact closures to find functions $\phi_i \colon M \to \R$ as in that proof using the pairs $(W_i, Y_i)$.
	If we then let $\phi \colon M \to \R$ be such that $0 < \phi(y) < \min\{\phi_i(y)\}$ for all $i$ and all $y \in M$, the rest of the argument applies to show that again $W_i \times I \xr{r_{W_i} \times \id_I} M \times I \xr{H} M$ is transverse to $r_{Y_i}$ for all $i$, in addition to the original result that $H(-,1)r_V$ is transverse to $r_X$.
\end{proof}

\subsection{Transverse complexes}\label{S: product pullbacks}
\greg{Should this section go here?}

In \cref{S: cohomology pullback}, we showed that a continuous map $f \colon N \to M$ of manifolds without boundary yields a well-defined cohomology map $f^* \colon H^*_\Gamma(M) \to H^*_\Gamma(N)$.
In this section, we consider $f^*$ as a partially-defined map of cochain complexes $f^* \colon C^*_\Gamma(M) \to C^*_\Gamma(N)$.
In general, $f^*$ will not be fully defined as a chain map, as we cannot form pullbacks of cochains that are not transverse to $f$.
To solve this problem, we define a subcomplex of $C^*_\Gamma(M)$ that will be quasi-isomorphic to $C^*_\Gamma(M)$ if $f$ is proper.

\begin{definition}\label{D: transverse to map}\index{geometric cochain!transverse to a map}
	Let $f \colon N \to M$ be a smooth map from a manifold with corners to a manifold without boundary, and let $\uV \in C^*_\Gamma(M)$.
	We will say that $\uV$ is \textbf{transverse to $f$} if $\uV$ has a representative $r_V \colon V \to M$ such that $r_V$ is transverse to $f$.
	In this case we define the pullback $f^*(\uV)$ to be $\underline{f^*(V)} \in C^*_\Gamma(N)$.

	We will write the set of cochains transverse to $f$ as $C^*_{\Gamma \pf f}(M)$.\index{$C^*_{\Gamma \pf f}(M)$}
\end{definition}

We notice that the transversality situation here is simpler than the more general ones in \cref{S: chain products}, as $f$ is a fixed map.
We also remark that given a cubulation $X$ of $M$, it determines a proper map from the disjoint union of all the $\dim(M)$-cubes of $X$ to $M$, which we can also write $F \colon X \to M$ with an abuse of notation.
In this case, $C^*_{\Gamma \pf f}(M) = C^*_{\Gamma \pf X}(M)$.

\begin{proposition}\label{P: trans to f}\index{geometric cochain!transverse to a map}\index{geometric cochain!transverse to a map!pullback of}
	Given a smooth map from a manifold with corners to a manifold without boundary $f \colon N \to M$, the set $C^*_{\Gamma \pf f}(M)$ is a subcomplex of $C^*_{\Gamma}(M)$.
	If $\bd N = \emptyset$, then $f^*:C^*_{\Gamma \pf f}(M) \to C^*_{\Gamma}(N)$ is a well-defined chain map.
\end{proposition}

\begin{proof}
	If $\uV,\uW \in C^*_\Gamma(M)$ are represented by $r_V \colon V \to M$ and $r_W \colon W \to M$ that are transverse to $f$, then $\uV+\uW$ can be represented by $V \sqcup W$, which will also be transverse to $f$.
	So $C^*_{\Gamma \pf f}(M)$ is closed under addition.
	If $r_V \colon V \to M$ is transverse to $f$ then so is $-r_V$, i.e.\ $r_V$ with the opposite co-orientation, so $C^*_{\Gamma \pf f}(M)$ is closed under negation.
	The empty map $\emptyset \to M$ is always transverse to $f$ (since there are no points at which to check the tangent space condition), and so $0 \in C^*_{\Gamma \pf f}(M)$.
	Finally, if $\uV$ is represented by $r_V \colon V \to M$ transverse to $f$, then by definition $\bd V \to M$ is transverse to $f$, so $\bd \uV \in C^*_{\Gamma \pf f}(M)$.
	Therefore, $C^*_{\Gamma \pf f}(M)$ is a subcomplex of $C^*_{\Gamma}(M)$.

	Now suppose $\bd N = \emptyset$.
	To show that $f^*$ is well defined on $C^*_{\Gamma \pf f}(M)$ we must show that it does not depend on the choice of representative $V$.
	Suppose $V$ and $V'$ both represent $\uV$ and are transverse to $f$.
	Then $V \sqcup -V'$ is transverse to $f$ and an element of $Q^*(M)$.
	So by \cref{L: pullback map Q}, $f^*(V \sqcup -V')$, which is by definition $(V \sqcup -V') \times_M N = (V \times_M N) \sqcup (-V' \times_M N)$ mapping to $N$, is an element of $Q^*(N)$.
	So $f^*(V)$ and $f^*(V')$ represent the same element of $C^*_{\Gamma}(N)$.
	Thus $f^*$ is well defined.

	To see that $f^*$ is a homomorphism, let $V,W \to M$ represent elements of $C^*_{\Gamma}(M)$ that are transverse to $f$.
	Then
	$$f^*(\uV+\uW) = \underline{f^*(V \sqcup W)} = \underline{f^*(V) \sqcup f^*(W)} = \underline{f^*(V)}+\underline{f^*(W)},$$
	using the definitions and obvious properties of the pullback.
	Furthermore, as $\bd N = \emptyset$, $f^*$ is a chain map by \cref{leibniz}, .
\end{proof}

\begin{remark}
	While $C^*_{\Gamma \pf f}(N)$ is a subcomplex, it is not closed under taking cup products, even when they are well defined.
	As an example, let $f \colon M \to N$ be the inclusion of the $x$-axis into the plane $\R^2$.
	Let $V$ be represented by an embedding of $\R$ into $\R^2$ as the line $y = x$, and let $W$ similarly correspond to $y = -x$, with any co-orientations.
	Then $\uV \uplus \uW$ is represented by the embedding of the origin into $\R^2$, but this is not transverse to $f$, even though both $V$ and $W$ are transverse to $f$.
	For more general constructions that address this issue, see \cite{GBF47}.
\end{remark}

Generalizing \cref{T: transverse complex}, we have the following.

\begin{theorem}\label{T: transverse qi}\index{geometric cohomology!transverse to a map!is geometric cohomology}
	If $f \colon N \to M$ is a proper map from a manifold with corners to a manifold without boundary, then the inclusion $C^*_{\Gamma \pf f}(M) \into C^*_\Gamma(M)$ is a quasi-isomorphism.
\end{theorem}
\begin{proof}
	The proof is identical to that of \cref{T: transverse complex}, using \cref{P: perturb transverse to map} in place of \cref{P: ball stability}.
\end{proof}

The following corollary is immediate from \cref{{T: transverse qi}} and the construction in \cref{D: cohomology pullback and homology transfer}.
It says that, in the case at hand, the cohomology map of \cref{D: cohomology pullback and homology transfer} is induced by a chain map after restricting to a subcomplex of $C^*_\Gamma(M)$.

\begin{corollary}\label{C: transverse cohomology pullback}\index{geometric cohomology!transverse to a map!pullback of}
	Let $f \colon N \to M$ be a proper map between manifolds without boundary.
	The following diagram commutes, in which the top map is induced by the chain map $f^*$, the vertical map is induced by the inclusion $C^*_{\Gamma \pf f}(M) \into C^*_\Gamma(M)$, and the bottom map is that of \cref{D: cohomology pullback and homology transfer}.
	\begin{equation*}
	\begin{tikzcd}
		H^*_\Gamma(N) & H^*(C^*_{\Gamma \pf f}(M)) \arrow[l]  \arrow[d,"\cong"]\\
		& H^*_\Gamma(M). \arrow[lu]
	\end{tikzcd}
	\end{equation*}
	In particular, the image of the top map is independent of the choice of smooth $f$ within its (not necessarily proper) homotopy class.
\end{corollary}

\subsubsection{Kronecker pairing}\label{S: Kronecker}\index{Kronecker pairing|(}
\greg{Should this section go here?}
Using similar arguments to those in \cref{S: product pullbacks}, we consider a partially-defined Kronecker-type evaluation $C^*_\Gamma(M) \to \Hom(C_*^\Gamma(M),\Z)$.

The partially-defined cap product yields a partially-defined pairing
$$C^i_\Gamma(M) \times C_i^\Gamma(M) \xr{\nplus} C_0^\Gamma(M) \xr{\aug}\Z,$$
where $\aug \colon C_0^\Gamma(M) \to \Z$ is the augmentation map of \cref{D: aug}.
We consider here the extent to which this pairing corresponds to a function $C^i_\Gamma(M) \to \Hom(C_i^\Gamma(M),\Z)$.
This situation is closely related to the preceding discussion of pullbacks.

\begin{definition}\label{D: transverse to cohain}\index{geometric chain!transverse to a cochain}\index{$C_i^{\Gamma \pf \uV}(M)$}
	Let $\uV \in C^i_\Gamma(M)$ be a geometric cochain.
	We write $C_i^{\Gamma \pf \uV}(M)$ for the subgroup of $C_i^\Gamma(M)$ generated by geometric $i$-chains simply transverse to $\uV$.
\end{definition}

\begin{proposition}
	Given a geometric cochain $\uV \in C^i_\Gamma(M)$, the map $\aug(\uV\nplus -):C_i^{\Gamma \pf \uV}(M) \to \Z$ is a well-defined homomorphism.
\end{proposition}

\begin{proof}
	We first observe that $\uV\nplus -$ is defined on all elements of $C_i^{\Gamma \pf \uV}(M)$.
	If $\uW \in C_i^\Gamma(M)$ can be written as a sum $\uW = \sum \underline{W_i}$ with each $\underline{W_i}$ simply transverse to $\uV$, then $\uV\nplus \uW$ is well defined as $\sum \uV\nplus \underline{W_i}$ by \cref{T: multicup}.
	The element $0 \in C_i^\Gamma(M)$, as represented by the empty map, is simply transverse to $\uV$ with $\uV\nplus 0 = 0$, and if $\uW$ is transverse to $\uV$ then so is $-\uW$.
	\cref{D: multicup,T: multicup} imply that $\uV \nplus -$ is a homomorphism.
	We know that $\aug$ is a homomorphism, so the proposition follows.
\end{proof}

So given $\uV \in C^i_\Gamma(M)$, we obtain an element of $\Hom\left(C_i^{\Gamma \pf \uV}(M),\Z\right)$, but of course we will not in general obtain an element of $\Hom\left(C_i^{\Gamma}(M), \Z \right)$ due to transversality requirements.

\index{Kronecker pairing|)}

\subsection{The Kronecker pairing and the Universal Coefficient Theorem for geometric cohomology}\label{S: kroneker}

When applying the cap product to cohomology and homology classes of the same degree, we can compose with the augmentation map $\aug \colon H_0^\Gamma(M) \to \Z$ of \cref{D: aug} to obtain a bilinear Kronecker pairing
$$H^i_\Gamma(M) \otimes H_i^\Gamma(M) \xr{\nplus} H_0^\Gamma(M) \to \Z.$$
Taking the adjunct then provides a map
$$\alpha: H^i_\Gamma(M) \to \Hom(H_i^\Gamma(M),\Z).$$
Tracing through the definitions, this maps acts by counting the intersection number between a geometric chain and a geometric cochain in the sense of \cref{D: intersection number}.

When $H^i_\Gamma(M)$ is finitely generated, this map fits into a short exact sequence, just as for singular cohomology.

\begin{theorem}\label{T: UCT}\index{Universal Coefficient Theorem}
	If $H^i_\Gamma(M)$ is finitely generated, there is a short exact sequence
	\[
	0 \to \Ext\left(H_{i-1}^\Gamma(M),\Z\right) \to H^i_\Gamma(M) \xr{\alpha} \Hom\left(H_i^\Gamma(M),\Z\right) \to 0.
	\]
\end{theorem}

\begin{remark}
	The existence of a Universal Coefficient exact sequence holds even if $H^i_\Gamma(M)$ is not finitely generated, as we know by \cref{T: geometric is singular} that $H^i_\Gamma(M) \cong H^i(M)$, and then we have the usual singular cohomology Universal Coefficient Theorem.
	We can further identify $\Hom(H_i(M),\Z)$ and $\Ext(H_{i-1}(M),\Z)$ with $\Hom(H^\Gamma_i(M),\Z)$ and $\Ext(H^\Gamma_{i-1}(M),\Z)$, also using \cref{T: geometric is singular}.
	The reason we need to invoke finite generation in \cref{T: UCT} is that we use \cref{T: intersection qi}, which has that condition, in the proof. What we lose from \cref{T: UCT} by following instead the approach outlined in this remark is the identification of the map $H^i_\Gamma(M) \to \Hom(H_i^\Gamma(M),\Z)$ with the map $\alpha$ given by counting intersection numbers.
\end{remark}

\begin{proof}[Proof of \cref{T: UCT}]
	Let $M$ have a cubulation $X$, let $\mc I \colon C_{\Gamma \pf X}^*(M) \to K^*(X)$ be the intersection map of \cref{D: intersection homomorphism},
	and let $\mc J \colon K_*(X) \to C^\Gamma_*(M)$ be the map inducing the homology isomorphism of \cref{T: cubical homology iso}.
	We consider the diagram
	\[
	\begin{tikzcd}
		H^i_\Gamma(M) \arrow[r, "\alpha"] & \Hom(H_i^\Gamma(M),\Z) \arrow[dd, "\cong", "\mc J^*"'] \\
		H^i(C^*_{\Gamma \pf X}(M)) \arrow[u, "\cong"'] \arrow[d, "\cong", "\mc I"'] & \\
		H^i(K^*(X)) \arrow[r] & \Hom(H_i(K_*(X)), \Z).
	\end{tikzcd}
	\]
	The vertical maps on the left are isomorphisms by \cref{T: transverse complex,T: intersection qi}, while the right hand vertical map is an isomorphism by \cref{T: cubical homology iso}.
	We claim the diagram commutes.
	In fact, let $V \in PC_\Gamma^i(M)$ represent an element of $H^i(C^*_{\Gamma \pf X}(M))$.
	Then $V$ is transverse to the cubulation, and by definition the path clockwise around the diagram takes $\uV$ to a map that acts on an element $\xi$ of $H_*(K_*(X))$ represented by a $\Z$-linear combination of cubical faces $\sum_j c_j E_j$ by treating each $E_j$ as a geometric chain and forming
	$$\aug\left( V \times_M \sum_jc_jE_j\right) = \sum_jc_j\aug(V \times_M E_j).$$

	On the other hand, by \cref{D: intersection homomorphism}, the composition counterclockwise takes $\uV$ to a map that acts on $\xi$ by $\sum_j c_j I_M(V,E_j)$.
	But $I_M(V,E_j)$ is precisely $\aug(V \times_M E_j)$ by \cref{D: intersection number}.
	So the diagram commutes.

	We know that $K_i(X)$ is a free abelian group, so the bottom map of the diagram is a surjection by the algebraic Universal Coefficient Theorem, with kernel $\Ext(H_{i-1}(K_*(X)),\Z)$.
	The commutativity of the diagram thus implies that the top map of the diagram is a surjection with isomorphic kernel.
	To complete the proof, we again invoke \cref{T: cubical homology iso} to observe $\Ext(H_{i-1}(X),\Z) \cong \Ext(H^\Gamma_{i-1}(X),\Z)$.
\end{proof}

\begin{remark}
	Note that while we obtain the expected Universal Coeficient Theorem relating geometric cohomology and homology, we do not claim to have either an isomorphism or a quasi-isomorphism between $C^i_\Gamma(M)$ and $\Hom(C_i^\Gamma(M),\Z)$.
	In fact, as we do not know $C_*^\Gamma(M)$ to be a complex of free abelian groups (which we leave as an open question), it is not clear $\Hom(C_i^\Gamma(M),\Z)$ fits into a short exact Universal Coefficient-type sequence at all.
\end{remark}

\subsection{The geometric cup product is the usual cup product}\label{S: usual cup}

In this section we show that the geometric cup product agrees with the singular cup product in the sense that there is a natural ring isomorphism between singular cohomology $H^*(-)$ with the usual cup product and geometric cohomology $H^*_\Gamma(-)$ with the cup product $\uplus$.

\begin{theorem}\label{T: intersection is cup product}\index{geometric cohomology!cup product!is singular cup product}
	On the category of smooth manifolds without boundary and continuous maps, there are natural isomorphisms of functors $\Phi_p \colon H^p(-) \to H^p_\Gamma(-)$, $p \geq 0$, from singular cohomology to geometric cohomology that are also compatible with cup products.
	In other words, for each manifold without boundary $M$ there is a commutative diagram
	\[
	\begin{tikzcd}
		H^p(M) \otimes H^q(M) \arrow{r}{\smile} \arrow{d}{\Phi_p \otimes \Phi_q} &
		H^{p+q}(M) \arrow{d}{\Phi_{p+q}} \\
		H^p_\Gamma(M) \otimes H^q_\Gamma(M) \arrow{r}{\uplus} & H^{p+q}_\Gamma(M).
	\end{tikzcd}
	\]
\end{theorem}

Our proof is based on an axiomatic characterization of the cup product on manifolds due to Kreck and Singhof \cite[Proposition 12]{Krec10b}.\index{Kreck-Singhof Theorem}
As the proof of this proposition is only sketched in \cite{Krec10b}, we first fill in the details, restricting ourselves to $\Z$ coefficients and changing Kreck and Singhof's notation a bit to avoid conflicts with our earlier notation.
Before stating the result, we establish some further notation and conventions for this section.

In this section we assume the spheres $S^p$, $p>0$, to each have a fixed orientation.
We also want these orientations to be compatible in the sense that the composition $\nu \colon S^p \times S^q \to S^p \wedge S^q \cong S^{p+q}$ is orientation preserving away from the subspace that is collapsed to form the wedge product.
In particular, if $[S^p]$ and $[S^q]$ are the corresponding fundamental classes, the quotient should take $[S^p] \times [S^q]$ to $[S^{p+q}]$.
This can be arranged, for example, by modeling our spheres as the standardly-oriented cubes with their boundaries collapsed.
For each $p>0$, we let $s_p \in H^p(S^p) \cong \Hom(H_p(S^p), \Z)$ be the cohomology class that evaluates to $1$ on $[S^p]$.
Let $\pi_1 \colon S^p \times S^q \to S^p$ and $\pi_2 \colon S^p \times S^q \to S^q$ be the projections.

We let $K_p = K(\Z,p)$, $p>0$, be the Eilenberg-MacLane spaces, which we can assume have been constructed as CW complexes such that the $p+1$-skeleton of $K_p$ is $S^p$.
Let $\iota_p \in H^p(K_p)$ denote the fundamental class such that if $\phi_p \colon S^p \to K_p$ is the inclusion, then $\phi_p^*(\iota_p) = s_p$.
As the $p+1$ skeleton of $K_p$ is the image of $S^p$ under $\phi_p$, it is standard that $\phi_p^*$ is an isomorphism.
We also let $\mu \colon K_p \times K_q \to K_{p+q}$ be the unique-up-to-homotopy map that extends the collapse map $\nu$.

For $M$ connected, we always assume $H^0(M) \cong \Z$ generated by the class of the cochain $1 \in C^0(M)$,  i.e.\ the cochain that evaluates to $1$ on each singular $0$-simplex.
Then for $\lambda \in \Z$, we write $\lambda$ also for the class $\lambda 1$,.

\begin{proposition}[Kreck and Singhof, Proposition 12 of \cite{Krec10b}]\label{P: Kreck-Singhof pairing}\index{Kreck-Singhof Theorem}
	Consider singular cohomology $H^*(-)$ as a cohomology theory on smooth manifolds\footnote{As defined in \cite{Krec10b}; see the proof of \cref{T: geometric is singular} above.}.
	Suppose $\star$ is a natural multiplication on $H^*(-)$ such that if $M$ is connected and $\lambda \in H^0(M)$ then $\lambda\star \alpha = \alpha\star \lambda = \lambda\alpha$ for all $\alpha \in H^*(M)$ (and with the obvious extension when $M$ is not connected).
	Then if\footnote{Rather than $s_p \times s_q$, Kreck and Singhof require $\pi_1^*(s_p) \star \pi_2^*(s_q)$ to be the element of $H^{p+q}(S^p \times S^q)$ that evaluates to $1$ on the fundamental class of $S^p \times S^q$, but with our conventions that tensor products of cochains act by $(\alpha \otimes \beta)(x \otimes y) = \alpha(x)\beta(y)$, these are the same cohomology class (c.f.\ \cite[page 245]{Span81} and \cite[Section 3B]{Hatc02}.
	} $\pi_1^*(s_p) \star \pi_2^*(s_q) = s_p \times s_q \in H^{p+q}(S^p \times S^q)$ for all $p,q\geq 1$, the product $\star$ is the cup product.
\end{proposition}

\begin{proof}
	For a smooth manifold $M$, let $\alpha \in H^p(M)$ and $\beta \in H^q(M)$.
	The condition that $\lambda\star \alpha = \alpha\star \lambda = \lambda\alpha$ whem $M$ is connected already guarantees that $\star$ is the cup product when $p$ or $q$ is $0$, so we can suppose $p,q>0$.
	As $H^*$ is ordinary singular cohomology, we know that $\alpha$ and $\beta$ can be represented by maps $\bar \alpha \colon M \to K_p$ and $\bar\beta \colon M \to K_q$ with $\alpha = \bar \alpha^*(\iota_p)$ and $\beta = \bar\beta^*(\iota_q)$.
	Furthermore, $\alpha\smile \beta$ is the pullback of $\iota_{p+q}$ by the composition
	\begin{equation}\label{E: EM cross}
		M \xr{\diag} M \times M \xr{\bar\alpha \times \bar \beta} K_p \times K_q \xr{\mu} K_{p+q},
	\end{equation}
	while similarly $s_p \times s_q$ is the pullback of $\iota_{p+q}$ by
	\begin{equation}\label{E: sphere cross}
		S^p \times S^q \xr{\phi_p \times \phi_q} K_p \times K_q \xr{\mu} K_{p+q};
	\end{equation}
	see \cite[Section 4.3]{Hatc02}.

	As we will want to apply the naturality of $\star$ in the category of smooth manifolds, we will choose manifold replacements for $K_p$, $K_q$, and $K_{p+q}$.
	In particular, suppose we realize $K_p$ as a CW complex by the standard constructions and let $K_p^N$ be the $N$-skeleton of $K_p$ with $N$ much larger than $\dim(M)$.
	Then $K_p^N$ is homotopy equivalent to a finite simplicial complex \cite[Theorem 2C.5]{Hatc02}, and we can embed it simplicially into some Euclidean space and take an open regular neighborhood to get a smooth manifold $\mc K_p$ homotopy equivalent to the $N$-skeleton of $K_p$.
	We define $\mc K_q$ and $\mc K_{p+q}$ analogously, using a large enough skeleton $K^{N'}_{p+q}$ of $K_{p+q}$ for the restriction of $\mu$ to $K_p^N \times K_q^N \to K_{p+q}^{N'}$ to be defined.
	Abusing notation, we continue to write $\bar \alpha$, $\bar \beta$, $\phi_p$, $\mu$, etc.
	for the maps involving these manifold replacements of the $K_*$.
	These replacements will be sufficient for all cohomology and homotopy computations required in what follows.

	Next we make two more preliminary observations.
	The first is that it follows from $\star$ being natural with respect to pullbacks that, when $f^*$ is an isomorphism, the product $\star$ is also natural with respect to $(f^*)^{-1}$, as we see by applying the isomorphism $f^*$ to the claimed identity $(f^*)^{-1}(x)\star (f^*)^{-1}(y) = (f^*)^{-1}(x\star y)$.
	The second is that there is an evident commutative diagram
	\begin{equation}\label{D: projections}
		\begin{tikzcd}
			H^p(\mc K_p) \arrow{r}{\pi_1^*} \arrow[d, "\phi_p^*"] &
			H^p(\mc K_p \times \mc K_q) \arrow[d, "(\phi_p \times \phi_q)^*"] \\
			H^p(S^p) \arrow{r}{\pi_1^*} & H^p(S^p \times S^q)
		\end{tikzcd}
	\end{equation}
	and similarly for $\pi_2$, abusing notation to write $\pi_1$ and $\pi_2$ for the projections to the first and second factors for both pairs of spaces.
	Now we compute:
	\begin{align*}
		\alpha\smile \beta& = \diag^*(\bar \alpha \times \bar \beta)^*\mu^*(\iota_{p+q})&\text{see \eqref{E: EM cross}}\\
		& = \diag^*(\bar \alpha \times \bar \beta)^*((\phi_p \times \phi_q)^*)^{-1}(\phi_p \times \phi_q)^*\mu^*(\iota_{p+q})\\
		& = \diag^*(\bar \alpha \times \bar \beta)^*((\phi_p \times \phi_q)^*)^{-1}(s^p \times s^q)&\text{see \eqref{E: sphere cross}}\\
		& = \diag^*(\bar \alpha \times \bar \beta)^*((\phi_p \times \phi_q)^*)^{-1}(\pi_1^*(s_p) \star \pi_2^*(s_q))&\text{by assumption}\\
		& = \diag^*(\bar \alpha \times \bar \beta)^*(((\phi_p \times \phi_q)^*)^{-1}\pi_1^*(s_p) \star ((\phi_p \times \phi)^*)^{-1}\pi_2^*(s_q))&\text{by naturality}\\
		& = \diag^*(\bar \alpha \times \bar \beta)^*(\pi_1^*(\phi_p^*)^{-1}(s_p) \star \pi_2^*(\phi_q^*)^{-1}(s_q))&\text{by diagram \eqref{D: projections}}\\
		& = \diag^*(\bar \alpha \times \bar \beta)^*(\pi_1^*(\iota_p) \star \pi_2^*(\iota_q))\\
		& = (\diag^*(\bar \alpha \times \bar \beta)^*\pi_1^*(\iota_p)) \star (\diag^*(\bar \alpha \times \bar \beta)^*\pi_2^*(\iota_q))&\text{by naturality}\\
		& = \bar\alpha^*(\iota_p) \star \bar\beta^*(\iota_q)&\text{see below}\\
		& = \alpha\star\beta &\text{by definition}.
	\end{align*}
	For the penultimate equality, we have used that the composition of maps $$M \xr{\diag}M \times M \xr{\bar\alpha \times \bar \beta}K_p \times K_q \xr{\pi_1}K_p$$ is just $\bar \alpha$, and similarly for $\bar \beta$.
\end{proof}

	Now, recall that in the proof of \cref{T: geometric is singular}, which established an isomorphism between geometric and singular cohomology, we applied \cite[Theorem 10]{Krec10b}.
	That theorem of Kreck-Singhof shows that there is a natural isomorphism of these cohomology theories on the category of smooth manifolds and continuous maps.
	In fact, it shows there is such an isomorphism extending any given isomorphism of coefficients $\Phi_0 \colon H^0(pt) \to H^0_\Gamma(pt)$ to natural isomorphisms $\Phi_p \colon H^p(-) \to H^p_\Gamma(-)$ for all $p\geq 0$.
	We will here assume $\Phi_0$ chosen so that it takes $1 \in H^0(pt)$ to the element $\underline{pt} \in H^0_\Gamma(pt)$ represented by the identity $pt \to pt$ with tautological co-orientation (see \cref{E: first examples}).

	For the proof of \cref{T: intersection is cup product}, we need to arrange that for all $p\geq 1$ we have $\Phi_p(s_p) = s_p^\Gamma$, where
	$s_p \in H^p(S^p)$ is our preferred generator described above and $s_p^\Gamma \in H^p_\Gamma(S^p)$ is the generator represented by an embedded point with normal co-orientation agreeing with our chosen orientation of $S^p$.
	This will not necessarily be the case for the $\Phi_p$ output by \cref{T: geometric is singular}.
	However, part of the data for a cohomology theory in the Kreck-Singhof theory consists of the natural connecting maps $\delta$ of the Mayer--Vietoris sequence, and part of the output of the theorem is that the isomorphisms $\Phi_p$ commute with these connecting maps.
	Let us write the connecting map for a cohomology theory $h^*$ more explicitly as $h^p(U \cap V) \xr{\delta_p}h^{p+1}(U \cup V)$; we will generally write $\delta_p$ for the connecting map independent of which cohomology theory we are discussing.
	Of course if we replace a given $\delta_p$ by $-\delta_p$ for all spaces, then we still have a natural connecting map, and we will not have affected the exactness of the Mayer--Vietoris sequence.
	If we make such a change, we technically have a new cohomology theory with the same cohomology groups, but \cite[Theorem 10]{Krec10b} will output different isomorphisms $\Phi_p$.
	As \cref{T: intersection is cup product} does not particularly care about the signs of the connecting maps in the Mayer--Vietoris sequence, we will first tinker with the connecting maps in order to arrange that $\Phi_p(s_p) = s_p^\Gamma$ for all $p\geq 1$.

\begin{lemma}\label{L: connecting signs}
	Possibly by changing the signs of the connecting morphisms in the Mayer-Vietoris sequences for geometric cohomology, we can arrange for $\Phi_p(s_p) = s_p^\Gamma$ for all $p \geq 1$, where $\Phi_p \colon H^p(-) \to H^p_\Gamma(-)$ are the isomorphisms output by \cref{T: geometric is singular} given $\Phi_0(1) = \underline{pt}$, as above.
\end{lemma}
\begin{proof}
	Let
	\begin{align*}
		U_p& = \{(x_1,\ldots,x_{p+1}) \in S^p \mid x_{p+1}>-1/2\}\\
		V_p& = \{(x_1,\ldots,x_{p+1}) \in S^p \mid x_{p+1}<1/2\}.
	\end{align*}
	Then the equatorial inclusion $S^{p-1} \into U_p \cap V_p$ is a homotopy equivalence.
	We will abuse notation and let $s_p$ also denote its image under the isomorphism $H^{p}(S^p) \cong  H^{p}(U_{p+1} \cap V_{p+1})$ induced by the homotopy equivalence.
	For $p\geq 1$, we now choose the sign of $\delta_p \colon H^{p}(U_{p+1} \cap V_{p+1}) \to H^{p+1}(S^{p+1})$ so that $\delta_p(s_p) = s_{p+1}$.
	Similarly, for geometric cohomology we arrange for $\delta_p(s_p^\Gamma) = s_{p+1}^\Gamma$.
	For $p = 0$, to avoid confusion let us write $z_- = -1 \in \R$ and $z_+ = 1 \in \R$.
	We let $s_0$ be the element of $H^0(S^0) \cong \Z^2$ that restricts to $1 \in H^0(z_+)$ and $0 \in H^0(z_-)$.
	Similarly, let $s_0^\Gamma \in H^0_\Gamma(S^0)$ be represented by the identity map of $z_+$ with its tautological co-orientation, and then $s_0^\Gamma$ is the element of $H^0_\Gamma(S^0) \cong \Z^2$ that restricts to $\underline{pt} \in H^0_\Gamma(z_+)$ and $0 \in H^0_\Gamma(z_-)$.
	Then we choose the signs of $\delta_0$ so that $\delta_0(s_0) = s_{1}$ and $\delta_0(s_0^\Gamma) = s_{1}^\Gamma$.

	Taking $H^*(-)$ and $H^*_\Gamma(-)$ with these Mayer--Vietoris connecting maps and this $\Phi_0 \colon H^0(pt) \to H^0_\Gamma(pt)$, \cite[Theorem 10]{Krec10b} gives natural isomorphisms $\Phi_p \colon H^p(-) \to H^p_\Gamma(-)$ extending $\Phi_0$ on a point.
	The naturality implies that $\Phi_0(s_0) = s_0^\Gamma$.
	It now follows by induction, using the following diagram due to the commutativity of $\Phi_*$ with the connecting maps, that $\Phi_p(s_p) = s_p^\Gamma$ for all $p$:
	\[
	\begin{tikzcd}
		H^p(S^p) \cong H^p(U_{p+1} \cap V_{p+1}) \arrow{r}{\delta_p} \arrow[d, "\Phi_p"] &
		H^{p+1}(U_{p+1} \cup V_{p+1}) = H^{p+1}(S^{p+1}) \arrow[d, "\Phi_{p+1}"] \\
		H^p_\Gamma(S^p) \cong H^p_\Gamma(U_{p+1} \cap V_{p+1}) \arrow{r}{\delta_p} &
		H_\Gamma^{p+1}(U_{p+1} \cup V_{p+1}) = H^{p+1}_\Gamma(S^{p+1}).
	\end{tikzcd}
	\]
\end{proof}

\begin{corollary}\label{C: sphere product}
	Given $\Phi_p \colon H^p(-) \to H^p_\Gamma(-)$ as in \cref{L: connecting signs}, then $\Phi_{p+q}(s_p \times s_q) = s_p^\Gamma \times s_q^\Gamma$ for all $p,q \geq 1$.
\end{corollary}

\begin{proof}
	Recall the maps $\nu \colon S^p \times S^q \to S^p \wedge S^q \cong S^{p+q}$ and $\mu \colon K_p \times K_q \to K_{p+q}$, defined above.
	We consider the following diagram, which commutes by the naturality of $\Phi_{p+q}$:
	\[
	\begin{tikzcd}
		H^{p+q}(S^p \times S^q) \arrow{d}{\Phi_{p+q}} &
		H^{p+q}(S^{p+q}) \arrow[l, "\nu^*"'] \arrow{d}{\Phi_{p+q}} \\
		H^{p+q}_\Gamma(S^p \times S^q) &
		H_\Gamma^{p+q}(S^{p+q}) \arrow[l, "\nu^*"'].
	\end{tikzcd}
	\]
	Let $s_{p+q}^\Gamma$ be represented by the embedding of a point at $y \in S^{p+q}$, normally co-oriented consistently with the orientation of $S^{p+q}$.
	By possibly choosing a different $y$ if necessary, we can also choose a smooth map homotopic to $\nu$ that maps a Euclidean neighborhood of some point $x \in S^p \times S^q$ by an orientation-preserving diffeomorphism onto a neighborhood of $y$, taking $x$ to $y$ and the complement of the neighborhood of $x$ to the complement of the neighborhood of $y$.
	Then, from the definitions, the pullback of $s_{p+q}^\Gamma$ is the embedding of $x$ into $S^p \times S^q$ with normal co-orientation corresponding to the orientation of $S^p \times S^q$.
	By \cref{E: sphere product}, this is exactly $s_p^\Gamma \times s_q^\Gamma$, i.e.\ $\nu^*(s_{p+q}^\Gamma) = s_p^\Gamma \times s_q^\Gamma$.
	So, recalling that $\Phi_p(s_{p+q}) = s_{p+q}^\Gamma$, we have $\nu^*\Phi_{p+q}(s_{p+q}) = s_p^\Gamma \times s_q^\Gamma$.
	Thus, from the commutativity of the diagram and $\Phi_{p+q}$ being an isomorphism, it suffices to show that $\nu^*(s_{p+q}) = s_p \times s_q$.

	For this, consider the commutative diagram
	\[
	\begin{tikzcd}
		H^{p+q}(K_p \times K_q) \arrow[d, "(\phi_p \times \phi_q)^*"] &
		\arrow[l, "\mu^*"'] H^{p+q}(K_{p+q}) \arrow[d, "\phi_{p+q}^*"] \\
		H^{p+q}(S^p \times S^q) & \arrow[l, "\nu^*"'] H^{p+q}(S^{p+q}).
	\end{tikzcd}
	\]
	As the $p+q+1$ skeleton of $K_p \times K_q$ is $S^p \times S^q$, the vertical maps are isomorphisms.
	And we know from \eqref{E: sphere cross} and the definitions that
	$$(\phi_p \times \phi_q)^*\mu^*(\phi_{p+q}^*)^{-1}(s_{p+q}) = (\phi_p \times \phi_q)^*\mu^*(\iota_{p+q})\\
	= s_p \times s_q,$$
	so we must have $\nu^*(s_{p+q}) = s_p \times s_q$, as needed.
\end{proof}

We can now apply \cref{L: connecting signs,C: sphere product} to prove \cref{T: intersection is cup product}.

\begin{proof}[Proof of \cref{T: intersection is cup product}]
	Applying \cref{L: connecting signs}, we assume isomorphisms $\Phi_p \colon H^p(-) \to H^p_\Gamma(-)$, $p \geq 0$, such that $\Phi_p(s_p) = s_p^\Gamma$ for all $p \geq 0$.

	For each $M$, we can now define pairings $\star$ on $H^p(-) \otimes H^q(-) \to H^{p+q}(-)$ by the composition
	\[
	H^p(M) \otimes H^q(M) \xr{\Phi_p \otimes \Phi_q}
	H^p_\Gamma(M) \otimes H^q_\Gamma(M) \xr{\uplus}
	H^{p+q}_\Gamma(M) \xr{\Phi_{p+q}^{-1}}
	H^{p+q}(M).
	\]
	This pairing is natural, as all the maps are natural.
	We will apply \cref{P: Kreck-Singhof pairing} to show that this is really the cup product, which will prove the theorem.

	We first show that if $M$ is connected and $\lambda \in H^0(M) \cong \Z$ then $\lambda\star \alpha = \alpha\star \lambda = \lambda\alpha$ for all $\alpha \in H^*(M)$.
	Let $pt$ be an arbitrary point in $M$.
	By the naturality of $\Phi_0$, the following diagram, in which the horizontal maps are induced by the inclusion $pt \into M$, commutes:
	\[
	\begin{tikzcd}
		H^0(M) \arrow{r} \arrow[d, "\Phi_0"] & H^0(pt) \arrow[d, "\Phi_0"] \\
		H_\Gamma^0(M) \arrow{r} & H_\Gamma^0(pt)
	\end{tikzcd}
	\]
	The vertical maps are isomorphisms by our application of \cite[Theorem 10]{Krec10b}, and it is standard that the top map is an isomorphism.
	In fact, we can consider $H^0(M)$ as generated by the cochain $1_M$, and this pulls back to the generator $1_{pt} \in H^0(pt)$ (each represented by the map that takes a positively oriented point considered as a singular $0$-chain to $1$).
	It follows that the bottom map is an isomorphism.
	Consider the generator $\uM \in H_\Gamma^0(M)$ given by the identity map of $M$ with its tautological co-orientation $(\beta_M,\beta_M)$; see \cref{E: first examples}.
	This has normal orientation given by the positively-oriented $0$-dimensional normal bundle, so by the pullback construction of \cref{D: pullback coorient}, the pullback to $H^0_\Gamma(pt)$ is similarly represented by $\underline{pt}$, the identity map of $pt$ with its canonical co-orientation.
	As $\Phi_0(1_{pt}) = \underline{pt}$ by assumption, it follows from the commutativity that $\Phi_0(1_M) = \uM$.

	So for $\alpha \in H^p(M)$, we have
	\begin{align*}
		\lambda \star \alpha
		& = \Phi_p^{-1}(\lambda\uM \uplus \Phi_p(\alpha)) \\
		& = \Phi_p^{-1}(\lambda\Phi_p(\alpha)) \\
		& = \lambda\alpha,
	\end{align*}
	using the unital property of $\uplus$ --- see \cref{S: (co)chain properties}.
	The same argument holds for $\alpha\star \lambda$.
	If $M$ has multiple components, then these properties clearly hold component-wise, as needed.

	To apply \cref{P: Kreck-Singhof pairing}, it remains to show that $\pi_1^*(s_p) \star \pi_2^*(s_q) = s_p \times s_q$ for all $p,q\geq 1$.
	For this, we have
	\begin{align*}
		\pi_1^*(s_p) \star \pi_2^*(s_q)& = \Phi_{p+q}^{-1}(\Phi_p(\pi_1^*(s_p))\uplus\Phi_q(\pi_2^*(s_q)))&\text{by definition of $\star$}\\
		& = \Phi_{p+q}^{-1}(\pi_1^*\Phi_p(s_p)\uplus\pi_2^*\Phi_q(s_q))&\text{by naturality of the $\Phi$}\\
		& = \Phi_{p+q}^{-1}(\pi_1^*(s_p^\Gamma)\uplus\pi_2^*(s_q^\Gamma))&\text{by \cref{L: connecting signs}}\\
		& = \Phi_{p+q}^{-1}((s_p^\Gamma \times \underline{S^q})\uplus(\underline{S^p} \times s_q^\Gamma))&\text{by Prop.
			\ref{P: projection pullbacks}}\\
		& = \Phi_{p+q}^{-1}((s_p^\Gamma\uplus \underline{S^p})\times( \underline{S^q}\uplus s_q^\Gamma))&\text{by Cor.
			\ref{C: criss cross} }\\
		& = \Phi_{p+q}^{-1}(s_p^\Gamma \times s_q^\Gamma)&\text{by Cor.
			\ref{C: cup with identity}}\\
		& = s_p \times s_q & \text{by \cref{C: sphere product}}.
	\end{align*}
\end{proof}

\subsection{K\"unneth theorems}\label{S: kunneth}

Now that we know that the geometric cup product is naturally isomorphic to the singular chain cup product, we can use this to compare cohomology cross products and obtain the geometric cohomology K\"unneth Theorem.
We begin with the homology K\"unneth Theorem, which is simpler, and then address the cohomology one.

Recall from \cref{T: hom iso map,P: singular smooth cubes} that we have isomorphisms $$H_*(NK_*(M)) \xleftarrow{\cong} H_*(NK^{sm}_*(M)) \xr{\cong} H_*^\Gamma(M),$$ where $NK_*(M)$ is the complex of normalized singular cubical chains and $NK^{sm}_*(M) \subset NK_*(M)$ is the subcomplex generated by smooth singular cubes.
As elements of $NK^{sm}_*(M)$ are represented by linear combinations of smooth maps from cubes and the cross product is represented by taking geometric products, we have the following immediate compatibility of chain cross products. Note that the upper left vertical map is an inclusion as the $NK^{sm}_i(M)$ and $NK_i(M)$ are all free groups generated by the nondegenerate singular cubes.

\begin{lemma}\label{L: chain cross compare}\index{geometric homology!exterior product!is singular homology exterior product}
	Let $M$ and $N$ be manifolds without boundary.
	Then the following diagram commutes:
	\[
	\begin{tikzcd}
		NK_*(M) \otimes NK_*(N) \arrow{r}{\times}  & NK_*(M \times N)  \\
		NK^{sm}_*(M) \otimes NK^{sm}_*(N) \arrow{r}{\times} \arrow[d] \arrow[hookrightarrow]{u}& NK^{sm}_*(M \times N) \arrow[d]\arrow[hookrightarrow]{u} \\
		C_*^\Gamma(M) \otimes C_*^\Gamma(N) \arrow{r}{\times} & C_*^\Gamma(M \times N).
	\end{tikzcd}
	\]
	This induces the commutative diagram
	\[
	\begin{tikzcd}
		H_*(NK_*(M)) \otimes H_*(NK_*(N)) \arrow{r}{\times}  &
		H_*(NK_*(M \times N)) \\
		H_*(NK^{sm}_*(M)) \otimes H_*(NK^{sm}_*(N)) \arrow{r}{\times} \arrow[d, "\cong"] \arrow[u, "\cong"']&
		H_*(NK^{sm}_*(M \times N)) \arrow[d, "\cong"] \arrow[u, "\cong"'] \\
		H_*^\Gamma(M) \otimes H_*^\Gamma(N) \arrow{r}{\times} & H_*^\Gamma(M \times N).
	\end{tikzcd}
	\]
\end{lemma}

\begin{theorem}[K\"unneth Theorem]\label{T: homology kunneth}\index{K\"unneth Theorem!geometric homology}\index{geometric homology!K\"unneth Theorem}
	Let $M$ and $N$ be manifolds without boundary.
	There are natural short exact sequences
	\[
	0 \to \bigoplus_{p+q = a} H_p^\Gamma(M) \otimes H_q^\Gamma(N) \xr{\times} H_{p+q}^\Gamma(M \times N) \to \bigoplus_{p+q = a-1} H_p^\Gamma(M)* H_q^\Gamma(N) \to 0
	\]
	that split (non-naturally).
\end{theorem}

\begin{proof}
	As the $C_i^\Gamma(M)$ and $C_i^\Gamma(N)$ are flat by \cref{L: flat}, there is such a split short exact sequence with middle term $H_*(C_*^\Gamma(M) \otimes C_*^\Gamma(N))$ by the algebraic K\"unneth Theorem \cite[Theorem V.2.1]{HS}.
	We must show that $C_*^\Gamma(M) \otimes C_*^\Gamma(N) \xr{\times} C_*^\Gamma(M \times N)$ is a quasi-isomorphism.
	But now in the first diagram of \cref{L: chain cross compare}, the top map is a quasi-isomorphism by
	 \cite[Theorem XI.3.1]{Mas91}.
	Furthermore, the vertical maps on the right are quasi-isomorphisms by \cref{T: hom iso map,P: singular smooth cubes}, and so, as all modules are flat, the vertical maps on the left are also quasi-isomorphisms (e.g.\ apply the K\"unneth Theorem and then the Five Lemma to each vertical map).
	So the bottom horizontal map is a quasi-isomorphism.
\end{proof}

We now turn to cohomology.
We will provide two proofs of the cohomology K\"unneth Theorem, one using our maps $\Phi$ from the Kreck-Singhof argument of \cref{T: intersection is cup product} and one using our intersection maps $\mc I$ in the cubical setting of \cref{S: intersection map}.

For the first, recall that for $\uV \in H^*_\Gamma(M)$ and $\uW \in H^*_\Gamma(N)$ we have the relation
$\uV \times \uW = \pi_M^*(\uV)\uplus\pi_N^*(\uW)$, which follows from \cref{C: cross is cup}, while the same relation is well known to hold in singular cohomology \cite[Corollary 5.6.14]{Span81}.
So the following is immediate from the naturality of our comparison maps $\Phi$ and \cref{T: intersection is cup product}

\begin{proposition}\label{P: cross product is cross product}\index{geometric cohomology!exterior product!is singular cohomology exterior product}
	On the category of smooth manifolds without boundary and continuous maps, the isomorphisms $\Phi_p$ from singular cohomology to geometric cohomology are compatible with cross products.
	In other words, for manifolds without boundary $M$ and $N$ there are commutative diagrams
	\[
	\begin{tikzcd}
		H^p(M) \otimes H^q(N) \arrow{r}{\times} \arrow[d,"\Phi_p \otimes \Phi_q"', "\cong"] &
		H^{p+q}(M \times N) \arrow[d, "\Phi_{p+q}"', "\cong"] \\
		H^p_\Gamma(M) \otimes H^q_\Gamma(N) \arrow{r}{\times} &
		H^{p+q}_\Gamma(M \times N).
	\end{tikzcd}
	\]
\end{proposition}

\begin{theorem}[K\"unneth Theorem]\label{T: cohomology kunneth}\index{K\"unneth Theorem!geometric cohomology}\index{geometric cohomology!K\"unneth Theorem}
	If either $H^i_\Gamma(M)$ is finitely generated for all $i$ or $H^i_\Gamma(N)$ is finitely generated for all $i$, then there are natural short exact sequences
	\[
	0 \to \bigoplus_{p+q = a}H^p_\Gamma(M) \otimes H^q_\Gamma(N) \xr{\times} H^{p+q}_\Gamma(M \times N) \to \bigoplus_{p+q = a+1}H^p_\Gamma(M)* H^q_\Gamma(N) \to 0
	\]
	that split (non-naturally).
\end{theorem}

\begin{proof}
	Using \cref{P: cross product is cross product}, we can form the left part of the diagram
	\[
	\begin{tikzcd}
	\displaystyle\bigoplus_{p+q = a}H^p(M) \otimes H^q(N) \arrow[r,"\times",hook]\arrow[d,"\bigoplus \Phi_p \otimes \Phi_q"',"\cong"]& H^{p+q}(M \times N) \arrow[r,twoheadrightarrow]\arrow[d,"\Phi_{p+q}"',"\cong"]&\displaystyle\bigoplus_{p+q = a+1}H^p(M)* H^q(N) \arrow[d,dashed,"\bigoplus \Phi_p * \Phi_q"',"\cong"]\\
	\displaystyle\bigoplus_{p+q = a}H^p_\Gamma(M) \otimes H^q_\Gamma(N) \arrow[r,"\times"]& H^{p+q}_\Gamma(M \times N) \arrow[r,dashed,twoheadrightarrow]& \displaystyle\bigoplus_{p+q = a+1}H^p_\Gamma(M)* H^q_\Gamma(N).
	\end{tikzcd}
	\]
	Since we know there is a cohomology K\"unneth Theorem for singular cohomology \cite[Theorem 60.5]{Mun84}, and the vertical maps on the left and middle are isomorphisms, the bottom left horizontal map is injective.
	It follows that there is an isomorphism between the quotient terms of the two sequences, and the quotient term on the top is $\displaystyle\bigoplus_{p+q = a+1}H^p_\Gamma(M)* H^q_\Gamma(N)$, which is then isomorphic to $\displaystyle\bigoplus_{p+q = a+1}H^p_\Gamma(M)* H^q_\Gamma(N)$ via the maps $\Phi_p*\Phi_q$.

	As the exact sequences are isomorphic and the top one splits, they both split.

	This construction is natural since all of the morphisms involved (except the splitting morphisms) are natural.
\end{proof}

In the preceding argument we had to work directly with cohomology because our maps $\Phi_p$ are only defined on cohomology, not on cochains.
Alternatively, we can approach the cohomology K\"unneth theorem more explicitly and in closer analogy with the proof of \cref{T: homology kunneth} by utilizing cubulations and intersection maps.

One benefit of cubulations over triangulations is that the product of two cubical complexes is again a cubical complex: the product of two cubes $E$ and $F$ is simply the cube $E \times F$.
Furthermore, as noted in \cref{S: cubical cochains}, this product provides a canonical isomorphism $K_*(X) \otimes K_*(Y) \xr{\times} K_*(X \times Y)$ for any cubical complexes $X$ and $Y$.
Passing to cohomology, we thus have maps
\begin{multline*}
	K^*(X) \otimes K^*(Y) = \Hom(K_*(X), \Z) \otimes \Hom(K_*(Y), \Z) \xr{\theta} \Hom(K_*(X) \otimes K_*(Y),\Z) \\ \underset{\cong}{\xleftarrow{\times^*}} \Hom(K_*(X \times Y),\Z) = K^*(X \times Y),
\end{multline*}
where $\times^*$ is the Hom dual of $K_*(X) \otimes K_*(Y) \xr{\times} K_*(X \times Y)$ and $\theta$ is the canonical map so that $$\theta(\alpha \otimes \beta)(E \otimes F) = \alpha(E)\beta(F);$$ compare \cite[Sections 60,61]{Mun84}.
The composition left to right will be our \textbf{cubical cochain cross product},\index{exterior product!cubical cochain} which we denote simply $\times$, by the standard convention.

The cubical cochain cross product induces a cohomology cross product $H^*(X) \otimes H^*(Y) \xr{\times} H^*(X \times Y)$ that also takes the tensor product of cohomology classes represented by cocycles $\alpha$ and $\beta$ to the class of the cocycle $\alpha \times \beta$.

If either $H_*(X)$ or $H_*(Y)$ are finite in all dimensions, then $\theta \colon \Hom(K_*(X), \Z) \otimes \Hom(K_*(Y), \Z) \xr{\theta} \Hom(K_*(X) \otimes K_*(Y),\Z)$ is a quasi-isomorphism, as can be seen by replacing $K_*(X)$ or $K_*(Y)$ up to chain homotopy equivalence by a free complex of finite type; see \cite[Lemma 5.5.6 and 5.5.9]{Span81}.
Consequently, the cubical cochain cross product $K^*(X) \otimes K^*(Y) \xr{\times}  K^*(X \times Y)$ is a quasi-isomorphism.
Using this together with the algebraic K\"unneth theorem \cite[Theorem V.2.1]{HS} and assuming these finiteness hypotheses, there is a natural K\"unneth short exact sequence for cubical cohomology
	\[
	0 \to \bigoplus_{p+q = a}H^p(X) \otimes H^q(Y) \xr{\times} H^{p+q}(X \times Y) \to \bigoplus_{p+q = a+1} H^p(X)* H^q(Y) \to 0
	\]
that splits, but not naturally.

\begin{proposition}\label{P: cross product comparison}\index{geometric cohomology!exterior product!is singular cohomology exterior product}
	Let $M$ and $N$ be manifolds without boundary cubulated by cubical complexes $X$ and $Y$, and let $M \times N$ have the product cubulation.
	The following diagram commutes
	\[
	\begin{tikzcd}
		C^*_{\Gamma}(M) \otimes C^*_{\Gamma}(N) \arrow{r}{\times} & C^*_{\Gamma}(M \times N) \\
		C^*_{\Gamma\pf}(M) \otimes C^*_{\Gamma\pf}(N) \arrow{r}{\times} \arrow[u, hook] \arrow[d, "\mc I \otimes \mc I"] & C^*_{\Gamma\pf}(M \times N) \arrow[u, hook] \arrow[d, "\mc I"] \\
		K^*(X) \otimes K^*(Y) \arrow{r}{\times} & K^*(X \times Y).
	\end{tikzcd}
	\]
	Furthermore, if all $H^i(M)$ and $H^j(N)$ are finitely generated, then the cubical cohomology cross product is isomorphic to the geometric cohomology cross product.
	In particular, we have the following diagram with all vertical maps isomorphisms:
	\[
	\begin{tikzcd}
		H^*_{\Gamma}(M) \otimes H^*_{\Gamma}(N) \arrow{r}{\times} & H^*_{\Gamma}(M \times N) \\
		H^*_{\Gamma\pf}(M) \otimes H^*_{\Gamma\pf}(N) \arrow{r}{\times} \arrow[u, "\cong"'] \arrow[d, "\cong", "\mc I \otimes \mc I"'] &
		H^*_{\Gamma\pf}(M \times N) \arrow[u, "\cong"'] \arrow[d, "\cong", "\mc I"'] \\
		H^*(X) \otimes H^*(Y) \arrow{r}{\times} & H^*(X \times Y).
	\end{tikzcd}
	\]
\end{proposition}

\begin{proof}
	In the top diagram, the top square certainly commutes as the vertical maps are inclusions.
	Note that the product of two maps transverse to the cubulation will be transverse to the product cubulation, so the middle horizontal map is well defined.

	Let $V$ and $W$ represent elements of $C^*_{\Gamma\pf}(M)$ and $C^*_{\Gamma\pf}(N)$, and let $E$ and $F$ be cubes of the cubulations $X$ and $Y$, respectively.
	We check that the two ways around the bottom square evaluate the same on $E \times F$.
	Applying $\mc I(V \times W)$ to $E \times F$ gives $I_{M \times N}(V \times W, E \times F) = \aug((V \times W)\times_{M \times N}(E \times F))$ by \cref{D: intersection homomorphism,D: intersection number}, while going the other way around the diagram and applying the result to $E \times F$ yields $I_M(V,E)I_N(W,F) = \aug(V \times_M E)\aug(W \times_N F)$.

	We can now compute
	\begin{align*}
		\aug((V \times W)\times_{M \times N}(E \times F))& = (-1)^{(w+f-n)(m-v)}\aug((V \times_M E) \times (W \times_N F))\\
		& = (-1)^{(w+f-n)(m-v)}\aug(V \times_M E)\aug(W \times_N F).
	\end{align*}
	The first equality is due to \cref{P: cap cross}.
	Note that if either $V \times_M E$ or $W \times_N F$ is not $0$-dimensional then also $(V \times W)\times_{M \times N}(E \times F)$ is not $0$-dimensional, and all three expressions above are $0$.
	Otherwise, the second equality is apparent as the product of a $(-1)^a$-oriented point with a $(-1)^b$-oriented point is a $(-1)^{a+b}$-oriented point.
	In this case we also have $w+f = n$ so that $(-1)^{(w+f-n)(m-v)} = 1$.

	The bottom diagram now commutes as a consequence of the first diagram commuting, and the vertical maps are isomorphisms by \cref{T: transverse complex,T: intersection qi}.
\end{proof}

\cref{P: cross product comparison} now implies an alternative proof of \cref{T: cohomology kunneth} when all of the cohomology groups are finitely generated and the homology groups of either $X$ or $Y$ are all finitely generated.
In this case the vertical maps of the first diagram in \cref{P: cross product comparison} are all quasi-isomorphisms by \cref{T: transverse complex,T: intersection qi,L: flat}, and the
cubical cochain cross product $K^*(X) \otimes K^*(Y) \xr{\times} K^*(X \times Y)$ is a quasi-isomorphism by our discussion preceding the statement of the proposition.
Thus by the commutativity of the diagram, $C^*_{\Gamma}(M) \otimes C^*_{\Gamma}(N) \xr{\times}  C^*_{\Gamma}(M \times N)$ is a quasi-isomorphism, and the K\"unneth short exact sequence for $C^*_{\Gamma}(M) \otimes C^*_{\Gamma}(N)$ is again a purely algebraic matter \cite[Theorem V.2.1]{HS}.

While this argument for the cohomology K\"unneth Theorem requires more hypotheses than the proof of \cref{T: cohomology kunneth} given above, they will always be satisfied if $X$ and $Y$ are compact, and the benefit, as seen in the second diagram of \cref{P: cross product comparison}, is that we obtain isomorphisms between our geometric cohomology K\"unneth short exact sequence and the cubical cohomology K\"unneth short exact sequence given by the geometric intersection maps.

\subsection{Relating geometric and singular cup and cap products via intersection maps}\label{S: cubical cup and cap}

Our goal in this section is to relate the cup and cap products in geometric (co)homology to those in singular (co)homology using cubical (co)homology as an intermediary and via the intersection map $\mc I$ of \cref{D: intersection homomorphism}.
While we have already shown in Section \cref{S: usual cup} that the geometric cohomology cup product is abstractly isomorphic to the singular cohomology cup product, the approach through intersection provides a concrete comparison map, which, by \cref{T: intersection qi} also provides cohomology isomorphisms when $H^i(M)$ is finitely generated.
In the case of cap products this is especially necessary, as the work of Kreck and Singhof in \cite{Krec10b} does not directly apply to provide a more abstract comparison as we did for cup products.

In \cref{S: cubical products}, we first discuss formulas for the cubical cup and cap products, relying on known formulas for the singular cubical products. Then we apply the cubical formulas to the geometric world, first in \cref{S: cup via intersection} where we relate the geometric and cubical cup products.
In \cref{S: cap product via intersection} we show that the cubical cap product (and hence the singular cap product) determines the geometric cap product in general, while the geometric cap product determines the cubical cap product if all $H^i(M)$ are finitely generated.

We then turn to some applications in \cref{S: umkehr}, which concerns Poincar\'e Duality, and \cref{S: umkehr}, which concerns umkehr, or transfer, maps.

\subsubsection{Cubical cup and cap products}\label{S: cubical products}\index{cup product!cubical cohomology}\index{cap product!cubical cohomology}\index{cap product!cubical homology}

In this section we discuss cup and cap products for cubical and singular cubical homology and cohomology.
This will be needed below for comparing the geometric cap product with the classical cap products.

We first recall from Massey \cite[Chapter XI]{Mas91} some results about the normalized singular cubical chain complexes, which we have been denoting $NK_*(-)$, though we utilize some different notation from Massey.
Just as for the more familiar singular simplicial chains, there is an Eilenberg-Zilber theorem that provides a chain homotopy equivalence between $NK_*(X) \otimes NK_*(Y)$ and $NK_*(X \times Y)$ for any spaces $X$ and $Y$.
Explicit constructions of such homotopy inverse maps are given in \cite[Section XI.5]{Mas91}.
The map $\zeta: NK_*(X) \otimes NK_*(Y) \to NK_*(X \times Y)$ is simply the cross product that takes $S \otimes T$ for representative singular cubes $S: \interval^m \to X$ and $T \colon \interval^n \to Y$ to the product $S \times T \colon \interval^m \times \interval^n = \interval^{m+n} \to X \times Y$.
If $S$ or $T$ is degenerate, so is $S \times T$, so this product is well defined for the normalized complexes.
The homotopy inverse map\footnote{Massey sometimes writes this map as $\eta$.} $\xi: NK_*(X \times Y) \to NK_*(X) \otimes NK_*(Y)$ takes $S \colon \interval^n \to X \times Y$ to
$$\xi(S) = \sum \rho_{H,K}A_H(\pi_1S) \otimes B_K(\pi_2S),$$ where $\pi_i$ is the projection to the $i$th factor.
The precise definitions of $\rho_{H,K}$, $A_H$, and $B_K$ will not need to concern us except to note that $H$ and $K$ are complementary subsets of $\{1,\ldots, n\}$, the sum is over all such partitions, $\rho_{H,K}$ is either $1$ or $-1$ (in fact it is the sign of the permutation $HK$), and $A_H$ and $B_K$ are cubical faces of various dimensions of the singular cubes $\pi_1S$ and $\pi_2S$.
Again, this construction is sufficiently compatible with degeneracies to be well defined for the normalized singular cubical complexes.
We also observe that if $X$ is a smooth manifold and our input singular cubes are smooth, then all other cubes appearing in the constructions are smooth.

As usual, one then defines cup and cap products (up to one's favorite sign conventions) as follows (using our current sign conventions).
If $\alpha, \beta \in NK^*(X) = \Hom(NK_*(X),\Z)$, then $\alpha\smile \beta \in \Hom(NK_*(X),\Z)$ acts on a normalized singular cube $S$ by
\begin{equation}\label{E: cubical cup}
	(\alpha \smile \beta)(S) = (\alpha \otimes \beta)(\xi(\diag S)),
\end{equation}
with $\diag$ the diagonal map $X \to X \times X$,
while the cap product $\alpha\frown S$ is given by
\begin{equation}\label{E: cubical cap}
	\alpha\frown S = (\id \otimes \alpha) (\xi(\diag S)),
\end{equation}
identifying $NK_*(X) \otimes \Z$ with $NK_*(X)$.

So now suppose $M$ is a manifold without boundary with cubulation $X$.
As usual, let $K_*(X)$ and $K^*(X) = \Hom(K_*(X),\Z)$ be the cubical chain and cochain complexes.
If $E$ is a face of $X$, we also regard $E$ as an element of $NK^{sm}_*(M)$ given by the embedding of $E$ into $M$.
If $E, F \in K_*(X)$ are any cubical faces, then $E \times F$ is also a cubical face, and  we have $\pi_1\diag(E) = \pi_2\diag(E) = E$.
Then the following diagram commutes, with the upward arrows being inclusions (note that all complexes are free) and the lower horizontal maps being the restrictions of the top horizontal maps:
\[
\begin{tikzcd}
	NK_*(M) \otimes NK_*(M) \arrow{r}{\zeta} & NK_*(M \times M) \\
	NK^{sm}_*(M) \otimes NK^{sm}_*(M) \arrow{r}{\zeta} \arrow[u, hook] & NK^{sm}_*(M \times M) \arrow[u, hook] \\
	K_*(X) \otimes K_*(X) \arrow{r}{\zeta} \arrow[u, hook] & K_*(X \times X). \arrow[u, hook]
\end{tikzcd}
\]
\begin{equation}\label{E: cube products}
\begin{tikzcd}
	NK_*(M) \otimes NK_*(M) & \arrow[l, "\xi"'] NK_*(M \times M) \\
	NK^{sm}_*(M) \otimes NK^{sm}_*(M) \arrow[u, hook] & \arrow[l, "\xi"'] NK^{sm}_*(M \times M) \arrow[u, hook] \\
	K_*(X) \otimes K_*(X) \arrow[u, hook] & \arrow[l, "\xi"'] K_*(X \times X). \arrow[u, hook]
\end{tikzcd}
\end{equation}

The top map in each diagram is a homotopy equivalences, and the vertical maps on the right are all quasi-isomorphisms by \cref{P: singular smooth cubes,T: cubical homology iso}.
As these are all free modules, the vertical maps on the left are also quasi-isomorphisms.
It follows that the horizontal maps are all quasi-isomorphisms, and in fact chain homotopy equivalences \cite[Theorem 46.2]{Mun84}.

In fact, we can say something stronger:

\begin{lemma}\label{L: EZ for cubical}
	The maps $K_*(X) \otimes K_*(X) \xr{\zeta}  K_*(X \times X)$ and $K_*(X \times X) \xr{\xi} K_*(X) \otimes K_*(X)$ are homotopy inverses, and similarly for $NK^{sm}_*(M) \otimes NK^{sm}_*(M) \xr{\zeta}  NK^{sm}_*(M \times M)$ and $NK^{sm}_*(M \times M) \xr{\xi} NK^{sm}_*(M) \otimes NK^{sm}_*(M) $.
\end{lemma}
\begin{proof}
	Consider the commutative square
	\[
	\begin{tikzcd}
	A \arrow[r,"\td f{,} \td g"]\arrow[d,"\phi"]& B\arrow[d,"\psi"]\\
	C \arrow[r,"f{,} g"]& D,
	\end{tikzcd}
	\]
	meaning that this diagram commutes with the pair of maps $f, \td f$ or with the pair of maps $g,\td g$.
	Let $\eta \colon D \to B$ be a chain homotopy inverse to $\psi$ and suppose $f$ is homotopic to $g$.
	Then $f \phi = \psi \td f$, so $\eta f \phi = \eta \psi \td f \sim \td f$, writing $\sim$ to indicate a chain homotopy.
	Similarly,  $\eta g \phi \sim \td g$.
	So then $$\td f \sim \eta f \phi  \sim \eta g \phi \sim \td g.$$

	The lemma now follows by applying this argument to the horizontal compositions of the diagrams above, using that $\xi \zeta$ and $\zeta \xi$ are each homotopic to identity maps on the $NK_*$ lines.
\end{proof}

Putting this all together, for both the complexes of normalized smooth singular cubical chains and the cubical complexes $K_*$ coming from the smooth cubulations, we may define cup and cap products again by the formulas \eqref{E: cubical cup} and \eqref{E: cubical cap}.
Note that in the case of a geometric cube $E$ coming from a cubulation, $\diag E$ is not a cube in the cubical complex, but that does not matter as in the end formula for $\xi(\diag E)$ we work with $\diag E$ only through its projections $\pi_1(\diag E) = \pi_2(\diag E) = E$.
These products are then compatible with the constructions for normalized singular cubical chains and cochains, i.e.\ the restriction of the cup product is the cup product of the restriction, and the appropriate mixed functoriality version of that statement holds for cap products.

It follows directly from this definition that the singular cubical cup product satisfies the familiar formula $$\alpha \smile \beta = \diag^*(\alpha \times \beta),$$
where $\alpha \times \beta$ is the cochain cross product given by $$\alpha \times \beta = \xi^*\theta$$
and $\theta \colon NK^*(M) \otimes NK^*(M) \to \Hom(NK_*(M) \otimes NK_*(M),\Z)$ is the canonical map so that $$(\theta(\alpha \otimes \beta))(E \otimes F) = \alpha(E)\beta(F).$$
In other words, $$\alpha \smile \beta = \diag^* \xi^* \theta(\alpha \otimes \beta) = (\xi \diag)^*\theta (\alpha \otimes \beta).$$

In the setting of cubical complexes, we have the maps $\xi$ and $\theta$, but we do not have the map $\diag \colon K_*(X) \to K_*(X \times X) $, so to express the cup product similarly as a pullback of the cross product in that setting, we replace $\diag$ with a map we call $\Delta \colon K_*(X) \to K_*(X \times X)$.

\begin{definition}\label{D: Delta}\index{Serre diagonal}
	For a cubical complex $X$ with geometric realization $|X|$, let $\Delta \colon K_*(X) \to K_*(X \times X)$ be the chain map given by $\Delta = \zeta(\xi\diag)$.
	This is a chain map because it is the restriction of a chain map of the singular cubical complexes to subcomplexes, noting again that $\xi\diag(E)$ can be identified as an element of $K_*(X \times X)$ even though $\diag E$ is not.
\end{definition}

The composition $\xi \Delta = \xi \zeta (\xi \diag)$ is chain homotopic to $\xi \diag$ as a map $\diag \colon K_*(X) \to K_*(X) \otimes K_*(X)$ by \cref{L: EZ for cubical}.
The upshot is that $(\xi \Delta)^*\theta = \Delta^* \xi^* \theta \colon K^*(X) \otimes K^*(X) \to K^*(X)$ is chain homotopic to the cup product, with each map well defined in the cubical context.

\begin{remark}\label{R: smooth Massey}
This discussion remains true replacing each instance of $NK_*(M)$ with $NK_*^{sm}(M)$.
The proof of the Eilenberg-Zilber Theorem in \cite[Section XI.5]{Mas91} showing that $\xi$ and $\zeta$ are natural homotopy inverses works just as well when considering only smooth singular cubical complexes on the category of (smooth) manifolds with corners.
In particular, in that case we can take the chain homotopies to also be realized by smooth singular cubical chains.
\end{remark}

The following proposition demonstrates the compatibility between the cubical, cubical singular, and smooth cubical singular cup products.

\begin{proposition}\label{P: cubical and singular cups}\index{cup product!cubical cohomology}
	Let $X$ be a cubical complex and $|X|$ its geometric realization.
	The following diagram commutes, with the top line involving the cohomology of cubical cochain complexes $K^*$, and each horizontal composition being the cup product:
	\[
	\begin{tikzcd}
	H^*(X) \otimes H^*(X) \arrow[r,"\times"]& H^*(X \times X) \arrow[r,"\Delta^*"] & H^*(X)\\
	H^*(NK^*_{sm}(|X|)) \otimes H^*(NK^*_{sm}(|X|)) \arrow[r,"\times"]\arrow[u,"\cong"] & H^*(NK^*_{sm}(|X| \times |X|)) \arrow[r,"\diag^*"]\arrow[u,"\cong"] & H^*(NK^*_{sm}(|X|))\arrow[u,"\cong"]\\
	H^*(NK^*(|X|)) \otimes H^*(NK^*(|X|))\arrow[r,"\times"]\arrow[u,"\cong"] & H^*(NK^*(|X| \times |X|)) \arrow[u,"\cong"]\arrow[r,"\diag^*"] &H^*(NK^*(|X|))\arrow[u,"\cong"].
	\end{tikzcd}
	\]
\end{proposition}
\begin{proof}
The commutativity of the left side of the diagram is a consequence of diagram \eqref{E: cube products} on page \pageref{E: cube products} after applying $\Hom(-,\Z)$,  together with the naturality of our $\theta$ maps of the form $\Hom(C_*,\Z) \otimes \Hom(D_*,\Z) \to \Hom(C_* \otimes D_*,\Z)$.

Commutativity on the right side follows from our preceding discussion, noting that $\zeta \xi \colon NK_*^{sm}(|X| \times |X|) \to NK_*^{sm}(|X| \times |X|)$ is chain homotopic to the identity by \cref{L: EZ for cubical}.
\end{proof}

\subsubsection{Relating geometric and cubical cup products via the intersection map}\label{S: cup via intersection}

We now consider the relationship between the geometric cohomology cup product and the classical cup product, as mediated by our intersection map given a cubical structure on the manifold.
Our goal is to prove the following theorem; putting this together with \cref{P: cubical and singular cups}, we obtain the compatibility of the geometric cohomology and singular (cubical) cohomology cup products.
Recall the definition of the intersection map $\mc I$ in \cref{D: intersection homomorphism} and our convention from \cref{R: intersection map extension} by which we abuse notation and also write $\mc I$ for the composition $H^*_\Gamma(M) \xleftarrow{\cong} H^*_{\Gamma \pf X}(M) \xr{\mc I} H^*(X)$.

\begin{theorem}\label{T: cup compatibility}\index{geometric cohomology!cup product!is cubical cohomology cup product}
Let $M$ be a manifold without boundary cubulated by the cubical complex $X$.
The following diagram, in which the horizontal compositions are the geometric cohomology cup product and the cubical cohomology cup product, commutes.
If $H^*(M)$ is finitely generated, then the vertical maps are isomorphisms.
	\[
	\begin{tikzcd}
	H^*_\Gamma(M) \otimes H^*_\Gamma(M) \arrow[r,"\times"]\arrow[d,"\mc I \otimes \mc I"]& H^*_\Gamma(M \times M) \arrow[d,"\mc I"]\arrow[r,"\diag^*"] & H^*_\Gamma(M) \arrow[d,"\mc I"]\\
	H^*(X) \otimes H^*(X) \arrow[r,"\times"]& H^*(X \times X) \arrow[r,"\Delta^*"] & H^*(X).
	\end{tikzcd}
	\]
\end{theorem}

The proof of \cref{T: cup compatibility} will rely on some further work.

For the next proposition, we need to further limit the geometric cochains we allow.
So let $M$ be a manifold without boundary cubulated by $X$, and recall that there is a natural chain homotopy between $\zeta \xi \colon NK^{sm}_*(M \times M) \to  NK^{sm}_*(M \times M)$ and the identity  \cite[Section XI.5]{Mas91}; see also \cref{R: smooth Massey}.
The naturality here means that these chain homotopies can be chosen compatibly across our category of manifolds with corners, and we will fix these as the homotopies in what follows.
So $\zeta\xi \diag$, which we denoted $\Delta$ in \cref{D: Delta}, is chain homotopic to $\diag \colon K_*(X) \to NK^{sm}_*(M \times M)$, thinking of $K_*(X)$ as a subcomplex of $NK_*^{sm}(M)$.
Let $P \colon K_i(X) \to NK^{sm}_{i+1}(M \times M)$ be the chain homotopy such that $\diag - \Delta = \bd P + P \bd$, coming from our fixed chain homotopy between $\zeta \xi$ and the identity.

We now define the subcomplex $C^*_{\Gamma \pf P}(M \times M) \subset C^*_{\Gamma \pf (X \times X)}(M \times M) \subset C^*_\Gamma(M \times M)$.
The definition is analogous to our definition of $C^*_{\Gamma \pf X}(M)$ in \cref{D: trans cube} except instead of just being transverse to the cubulation, we require
$PC^*_{\Gamma \pf P} (M \times M)$ to consist of those elements of $V \in PC^*_{\Gamma}(M \times M)$ such that
\begin{enumerate}
\item $V$ is transverse to each characteristic map of the cubulation $X \times X$ of $M \times M$,
\item $V$ is transverse to $\diag(E)$ for every $E$ in $X$, which also implies that $V$ is transverse to $\diag \colon M \to M \times M$,
\item $V$ is transverse to $P(E)$ for every $E$ in $X$, thinking of the cubical chain $P(E)$ as an element of $PC^\Gamma_*(M \times M)$ in the usual way.
\end{enumerate}

Due to the last condition, the maps that each $V$ must be transverse to are not all embeddings.
Nevertheless, a directly analogous construction to the discussion in \cref{D: trans cube} allows us to define $C^*_{\Gamma \pf P} (M \times M)$.

\begin{lemma}\label{L: P-transverse iso}
	The inclusions $C^*_{\Gamma \pf P} (M \times M) \to C^*_{\Gamma \pf X} (M \times M) \to C^*_\Gamma (M \times M)$ are quasi-isomorphisms.
\end{lemma}
\begin{proof}
	The proof that $C^*_{\Gamma \pf P} (M \times M)  \to C^*_\Gamma (M \times M)$ is a quasi-isomorphism is identical to that of \cref{T: transverse complex}, though since we do not need to be transverse only to faces of the cubulation, we use \cref{P: perturb transverse to map} rather than \cref{P: ball stability} to provide the needed transversality arguments.
	For this, we just need to create an appropriate $r_{\mf X}  \colon \mf X \to M$ so that precochains transverse to $\mf X$ will be in $PC^*_{\Gamma \pf P} (M \times M)$.
	So we let $\mf X$ be the disjoint union of
	\begin{enumerate}
		\item the cubes in $X \times X$,
		\item the cubes of $X$,
		\item for each cube $E$ in $X$, the domain of $P(E)$, which will be a collection of cubes, thinking of the cubical chain $P(E)$ as an element of $PC^\Gamma_*(M)$ in the usual way.
	\end{enumerate}
	The map $r_{\mf X}$ then consists of
	\begin{enumerate}
		\item the embeddings of the cubes of $X \times X$ into $M \times M$ via the cubulation,
		\item the embeddings of the cubes of $X$ into $M \times M$ via the cubulation of $M$ and the diagonal map $\diag$,
		\item for each cube $E$ in $X$, the map realizing $P(E)$ as a cubical chain.
	\end{enumerate}
	We note that $\mf X$ is a manifold with corners as it is a countable disjoint union of cubes.
	To see that $r_{\mf X}$ is proper, we observe that each $P(E)$ is a finite cubical chain and, furthermore, by the naturality of the chain homotopies in the proof of the Eilenberg-Zilber theorem through acyclic models (see \cite[Section XI.5]{Mas91}), $P(E)$ will be supported in $E \times E$.
	Therefore, each compact set of $M \times M$ will intersect the image of at most a finite number of the cubes of $\mf X$.

	So, by construction, transversality to $\mf X$ implies membership in $PC^*_{\Gamma \pf P} (M \times M)$, and we conclude inclusion $C^*_{\Gamma \pf P} (M \times M) \to C^*_\Gamma (M \times M)$ is a quasi-isomorphism by the proof of \cref{T: transverse complex} using \cref{P: perturb transverse to map} rather than \cref{P: ball stability}.
	As we already know $C^*_{\Gamma \pf X} (M \times M) \to C^*_\Gamma (M \times M)$ is a quasi-isomorphism by \cref{T: transverse complex}, the lemma follows.
\end{proof}

Next, we notice that the chain map $\mc I \colon C^*_{\Gamma \pf (X\times X)}(M\times M) \to K^*(X \times X)$ of \cref{D: intersection homomorphism} restricts to a chain map on $C^*_{\Gamma \pf P}(M\times M)$.
Furthermore, if we have any map $f \colon N \to M \times M$ that is transverse to every element of $C^*_{\Gamma \pf P}(M\times M)$ and if $\bd N = \emptyset$ , then we have a pullback chain map $f^* \colon C^*_{\Gamma \pf P}(M\times M) \to C^*_\Gamma(N)$ given by $f^*(\uV) = \underline{V \times_{M \times M} N}$ as a pullback mapping to $N$, where $V$ is any representative of $\uV$; cf.\ \cref{S: cohomology pullback} and the statement and proof of \cref{P: trans to f}.
In particular, this is the case for $\diag \colon M \to M \times M$, in which case we have a well defined chain map $\diag^* \colon C^*_{\Gamma \pf P}(M\times M) \to C^*_\Gamma(M)$.
In fact, we can say a bit more: If $V \in  C^*_{\Gamma \pf P}(M\times M)$, then we know that $V$ is transverse to the composition $E \into M \xr{\diag} M \times M$ for any face $E$ of $X$.
So by \cref{L: transverse to pullback}, the pullback $V \times_{M \times M} M \to M$ is transverse to the inclusion $E \into M$.
So the image of $\diag^*$ is in the subcomplex $C^*_{\Gamma \pf X}(M)$.
Furthermore, $\diag^* \colon C^*_{\Gamma \pf P}(M \times M) \to C^*_{\Gamma \pf X}(M)$ is a chain map by \cref{leibniz}, as $M$ has no boundary.

This allows us to state our key proposition.

\begin{proposition}\label{P: diag/intersect}\index{Serre diagonal!versus geometric diagonal}
	Let $M$ be a manifold without boundary, cubulated by the cubical complex $X$, and let $C^*_{\Gamma \pf P}(M \times M)$ be as defined above.
	The following diagram commutes up to chain homotopy:
	\[
	\begin{tikzcd}
	C^*_{\Gamma \pf P}(M \times M) \arrow[r,"\diag^*"] \arrow[d,"\mc I"]& C^*_{\Gamma \pf X}(M) \arrow[d,"\mc I"]  \\
	 \arrow[r,"\Delta^*"] K^*(X\times X)&K^*(X).
	\end{tikzcd}
	\]
\end{proposition}
\begin{proof}
	Let $V \in PC^*_{\Gamma \pf P}(M \times M)$ be a representative of $\uV \in C^*_{\Gamma \pf P}(M \times M)$, and let $E$ be a cubical face of $X$.
	We first compute as follows:
	\begin{align*}
		(\mc I \diag^*(V))(E) & = \aug(\diag^*(V) \times_M  E) &\text{\cref{D: intersection homomorphism,D: intersection number}}\\
		&=\aug(V \times_{M \times M} \diag(E)) &\text{\cref{P: natural cap}},
	\end{align*}
	where $\diag(E)$ denotes the map $E \into M \xr{\diag} M \times M$ in the notation of \cref{S: covariant functoriality}.
	On the other hand,
		\begin{align*}
		(\Delta^*\mc I(V))(E) & = \mc I(V)(\Delta (E))\\
		&=\aug(V \times_{M \times M} \Delta (E)).
	\end{align*}
	Note that since $\Delta (E) = \diag(E) - \bd P(E) - P \bd(E)$ and $V \in PC^*_{\Gamma \pf P}(M \times M)$, we have $\Delta(E)$ transverse to $V$, and all other needed transversality for the preceding computations follows from the definitions and discussion above.

	Next, we define $L \colon C^j_{\Gamma \pf P}(M \times M) \to K^{j-1}(X)$ by $$L(\uV)(E) = \aug(V \times_{M \times M} P(E))$$
	for any representative $V$ of $\uV$.
	This is well defined by arguments analogous to those for \cref{P: I is well defined}.
	Then we have
	\begin{align*}
		((dL + L\bd)(V))(E) &=(dL(V))(E) + ((L\bd)(V))(E) \\
			& = L(V)(\bd E) +L(\bd V)(E)\\
			& = \aug(V \times_{M \times M} P(\bd E)) + \aug((\bd V) \times_{M \times M} P(E))\\
			& = \aug(V \times_{M \times M} P(\bd E)) + (-1)^{v+e+1}\aug(\bd(V \times_{M \times M} P(E)))+(-1)^{v+e}\aug(V \times_{M \times M} \bd P(E))\\
			&=\aug(V \times_{M \times M} P(\bd E)) +\aug(V \times_{M \times M} \bd P(E))\\
			&=\aug(V \times_{M \times M} (P\bd + \bd P)(E))\\
			&=\aug(V \times_{M \times M} (\diag - \Delta)(E))\\
			&=\aug(V \times_{M \times M} (\diag(E))) - \aug(V \times_{M \times M} \Delta(E))\\
			&=(\mc I \diag^*(V))(E)- 	(\Delta^*\mc I(V))(E)\\
			&=((\mc I \diag^*- \Delta^*\mc I)(V))(E).
	\end{align*}
	Here in the fourth line we use \cref{P: Leibniz cap}, and in the fifth line we use that $\aug \circ \bd =0$ and that $\aug$ will be trivial unless $\dim(V)+ \dim(E) = \dim(M \times M) \equiv 0 \mod 2$.
	As $E$ and $V$ were arbitrary, we see that we have a chain homotopy between $\mc I \diag^*$ and $\Delta^*\mc I$.
\end{proof}

We obtain the following corollary.

\begin{corollary}\label{C: diag commute}
	The following diagram commutes:
	\[
	\begin{tikzcd}
		&H^*_{\Gamma}(M \times M) \arrow[r,"\diag^*"]& H^*_{\Gamma}(M) \\
	H^*_{\Gamma \pf (X \times X)}(M \times M)  \arrow[rd,"\mc I"']\arrow[ru,"\cong"] & \arrow[l,"\cong"']H^*_{\Gamma \pf P}(M \times M) \arrow[r,"\diag^*"] \arrow[d,"\mc I"] \arrow[u,"\cong"']& H^*_{\Gamma \pf X}(M) \arrow[d,"\mc I"] \arrow[u,"\cong"'] \\
		&H^*(K^*(X\times X))	\arrow[r,"\Delta^*"] & H^*(K^*(X)) .
	\end{tikzcd}
	\]
\end{corollary}
\begin{proof}
	The upper left triangle commutes because the maps are all induced by subcomplex inclusions, and the maps are all isomorphisms by \cref{L: P-transverse iso}.
	The bottom left triangle similarly commutes because $C^*_{\Gamma \pf P}(M \times M) \subset C^*_{\Gamma \pf (X \times X)}(M \times M)$ and, at the chain level, the map $\mc I$ on the right is a restriction of the $\mc I$ map on the left.
	The bottom right square commutes by \cref{P: diag/intersect}.
	Finally, the upper right square commutes at the chain level thanks to our observation above that the chain map $\diag^*$ takes the subcomplex $C^*_{\Gamma \pf P}(M \times M) \subset C^*_\Gamma(M \times M)$ to the subcomplex $C^*_{\Gamma \pf X}(M) \subset C^*_\Gamma(M)$.
\end{proof}

We can now prove \cref{T: cup compatibility}.

\begin{proof}[Proof of \cref{T: cup compatibility}]
We obtain the proof by concatenating the diagrams of \cref{P: cross product comparison} (after taking cohomology) and \cref{C: diag commute} and considering the rectangle formed on the outside.
\end{proof}

\subsubsection{The cup product via flowing}
While we have \cref{T: cup compatibility} concerning the cup product on cohomology, the geometric cup product is not fully-defined as a map of cochains $C^*_\Gamma(M) \otimes C^*_\Gamma(M) \to C^*_\Gamma(M)$ due to the transversality requirements considered in detail in \cref{S: simple products}.
However, we do have the following result from \cite{FMS-flows}, which we rephrase for the current context and for which we refer to \cite{FMS-flows} for details.
This essentially says that, given two geometric cochains $V$ and $W$ on a compact cubulated manifold, there is a canonical flow, depending only on the cubulation, such that after enough time, $V$ and $W$ are transverse and then their co-oriented fiber product is compatible with the cubical cup product via the intersection map.
We note that this implies another proof of \cref{T: cup compatibility} when $M$ is compact, observing that in our language here the flow provides a universal homotopy of $M$ to itself.

\begin{theorem}[Theorem 1 of \cite{FMS-flows}]\label{T: FMMS}\index{cup product!via flows}
	Let $M$ be a closed manifold cubulated by the cubical complex $X$, let $\f_t$ be the time $t$ flow of the logistic vector field associated to the cubulation, and let $V,W \in PC^*_{\Gamma \pf X}(M)$.
	Then, for $t$ sufficiently large:
	\begin{enumerate}
		\item $\f_t(W)$ and $\f_{-t}(V)$ are transverse and
		\[
		\mc I \left(\f_t(W) \times_M \f_{-t}(V)\right) = \mc I \left(\f_t(W)\right) \sms \mc I \left(\f_{-t}(V)\right) \in K^*(X).
		\]
		\item $\f_{-t}(W)$ and $\f_t(V)$ are transverse and
		\[
		\mc I \left(\f_{-t}(W) \times_M \f_t(V)\right) = (-1)^{(m-v)(m-w)} \, \mc I \left(\f_t(V)\right) \sms \mc I \left(\f_{-t}(W)\right) \in K^*(X).
		\]
	\end{enumerate}
\end{theorem}

In general, the flow time $t$ required will depend on $V$ and $W$, but we can say the following.

\begin{corollary}[Corollary 46 of \cite{FMS-flows}]
	Let $M$ be a closed cubulated manifold, and let $F^*$ be a finitely-generated chain complex with chain map $g \colon F^*\to C^*_{\Gamma \pf}(M)$.
	Then, there is a $T \in \R$ such that for all $t > T$ the following diagram commutes:
	\[
	\begin{tikzcd} [row sep=tiny]
		& C^*_{\Gamma \pf}(M) \otimes C^*_{\Gamma \pf}(M) \arrow[r, "\mc I \otimes \mc I"] & K^*(X) \otimes K^*(X) \arrow[dd, "\sms"] \\
		F^*\otimes F^* \arrow[ur, in=180, out=45,"g\otimes g"] \arrow[dr, in=180, out=-45, "\f_t \circ g \; \uplus \; \f_{-t} \circ g"']& & \\
		& C^*_{\Gamma \pf}(M) \arrow[r, "\mc I"] & K^*(X).
	\end{tikzcd}
	\]
\end{corollary}

If $F^*$ is the subcomplex of $C^*_{\Gamma \pf}(M)$ generated by two cochains $\uV$ and $\uW$, this recapitulates \cref{T: FMMS} as a statement about cochains, not just precochains.
More generally, the corollary applies if $F^*$ is a finitely-generated subcomplex of $C^*_{\Gamma \pf}(M)$ whose inclusion induces a quasi-isomorphism, in which case we obtain a fully-defined geometric cochain model for a cup product on $C^*_\Gamma(M)$.

\subsubsection{Relating geometric and cubical cap products via the intersection map}\label{S: cap product via intersection}

In this section we show that the geometric and cubical cap products are compatible in the sense given below in \cref{T: equivalent cap,C: cap relation}.
Recall that in \cref{R: intersection map extension} we extended the definition of the intersection map $\mc I$ to give us a map $H^*_\Gamma(M) \to H^*(X)$ for the manifold $M$ cubulated by $X$.
We again let $\mc J \colon K_*(X) \to C^\Gamma_*(M)$ be the map that takes a cubical face of $X$ to its embedding into $M$; see \cref{T: cubical homology iso}.

\begin{theorem}\label{T: equivalent cap}\index{cap product!is the cubical cap product}
	Let $M$ be a manifold without boundary smoothly cubulated by the cubical complex $X$.
	Let $\uV \in H^*_\Gamma(M)$ and $W \in H_*(X )$.
	Then
	$$\uV \nplus \mc J(W) = \mc J(\mc I(\uV)\frown W).$$
\end{theorem}

Here the cap product on the left is our geometric cap product and the cap on the right is the cubical cap product.

Before proving the theorem, we note the following corollary.

\begin{corollary}\label{C: cap relation}\index{cap product!is the cubical cap product}
	The cubical cap product (and hence the singular cap product) determine the geometric cap product.
	If all $H^i(M)$ are finitely generated, then the geometric cap product determines the cubical cap product.
\end{corollary}

\begin{proof}
	Let $\uV \in H^*_\Gamma(M)$ and $\uW \in H_*^\Gamma(M)$.
	Choose a cubulation $X$ of $M$.
	As $\mc J \colon H_*(X) \to H_*^\Gamma(M)$ is an isomorphism by \cref{T: cubical homology iso}, we have by \cref{T: equivalent cap}
	$$\uV\nplus \uW = \uV\nplus \mc J(\mc J^{-1}(\uW)) = \mc J(\mc I(\uV)\frown \mc J^{-1}(\uW)).$$
	On the other hand, suppose $V \in H^*(X)$ and $W \in H_*(X)$.
	Then $\mc I \colon H^*_\Gamma(M) \to H^*(X)$ is an isomorphism by \cref{T: intersection qi,R: intersection map extension} when all $H^i(M)$ are finitely generated.
	So then
	$$V\frown M = \mc J^{-1}(\mc J(\mc I\mc I^{-1}(V)\frown W)) = \mc J^{-1}(\mc I^{-1}(V) \nplus \mc J(W)).$$
\end{proof}

We will approach the proof of \cref{T: equivalent cap} through a series of lemmas.
The first two concern transversality.
Then we have a series of lemmas that essentially consist of various reformulations of the cap products, eventually linking together the two terms of \cref{T: equivalent cap}.
Once we have all the lemma established, we explain how to tie them all together to prove the theorem.

\begin{lemma}\label{L: product transversal}
	Let $M$ be a cubulated manifold without boundary.
	Let $V \in PC^*_{\Gamma}(M)$, and suppose $W \in PC_*(M \times M)$ is represented by a collection of embeddings.
	Then there is a proper universal homotopy $h \colon V \times I \to M$ such that $h(-,1) \colon V \to M$ is transverse to the cubulation and $M \times V \xr{\id_M \times h(-,1)} M \times M$ is transverse to $W$ in $M \times M$.
\end{lemma}

\begin{proof}
	As in the proof of \cref{P: ball stability}, we use the transversality techniques of \cite[Section 2.3]{GuPo74}.
	Consider $M$ as embedded in some $\R^N$ with an $\epsilon$-neighborhood $M_\epsilon$ and proper submersion $\pi \colon M_\epsilon \to M$.
	We define $\ms F$ as in proof of \cref{P: ball stability} so that $\ms F \colon V \times D^N \to M$ is given by $\ms F(x,s) = \pi(r_V(x)+\eta(x)s)$.
	Then $\ms F$ is a submersion and so transverse to each face of the cubulation, and also $\id_M \times \ms F \colon M \times V \times D^N \to M \times M$ is a submersion and hence transverse to $W$.
	So now by the Transversality Theorem of \cite[Section 2.3]{GuPo74}, for any face $E$ of the cubulation, $\ms F(-,s)$ is transverse to $E$ for almost all $s \in D^N$ and, similarly, $\id_M \times \ms F(-,s)$ is transverse to $W$ for almost all $s \in D^N$ (and similarly for each stratum of $W$).
	As the cubulation must have a countable number of faces and $W$ has a finite number of strata, there is an $s_0 \in D^N$ such that $\ms F(-,s_0)$ is transverse to the cubulation and $\id_M \times \ms F(-,s_0)$ is transverse $W$.
	Now let $h(x,t) = \ms F(x,ts_0)$.
	This is a proper universal homotopy, and $h(-,1)$ has the required transversality properties.
	See the proof of \cref{P: ball stability} for additional details.
\end{proof}

\begin{lemma}\label{L: M times transverse diag}
	Suppose $r_V \colon V \to M$ and $r_W \colon W \to M$ are transverse maps from manifolds with corners to a manifold without boundary.
	Then $\id_M \times r_V \colon M \times V \to M \times M$ is transverse to $\diag r_W \colon W \to M \times M$ in $M \times M$, where $\diag \colon M \to M \times M$ is the diagonal map.
	In particular, $\id_M \times r_V \colon M \times V \to M \times M$ is transverse to $\diag \colon M \to M \times M$ for any $r_V$.
\end{lemma}

\begin{proof}
	Without loss of generality, we suppose $V$ and $W$ are manifolds without boundary; the general case then holds by applying the following to pairs of strata.

	Suppose $x \in V$ and $y \in W$ such that $r_V(x) = r_W(y)$.
	Then $Dr_V(T_xV)+Dr_W(T_yW) = T_{r_V(x)}M$ by assumption.
	Now suppose $(z,x) \in M \times V$ and $y \in W$ such that $(\id_M \times r_V) (z,x) = \diag r_W(y)$.
	This is equivalent to $r_V(x) = r_W(y) = z$.
	At any such $(z,x)$, the image of $D(\id_M \times r_V)$ acting on $T_{(z,x)}(M \times V) = T_zM \oplus T_xV$ is $$T_zM \oplus Dr_V(T_xV) \subset T_zM \oplus T_zM = T_{(z,z)}(M \times M).$$
	The image of $D(\diag r_W)$ acting on $T_yW$ is $\{(a,a) \in T_z \oplus T_z \mid a \in Dr_W(T_yW)\}$.
	As we know $Dr_V(T_xV)+Dr_W(T_yW) = T_{r_V(x)}M$, when $r_V(x) = r_W(y) = z$ these images together span $(T_{(z,z)}M \times M)$.

	The last statement follows by taking $r_W \colon W \to M$ to be $\id_M \colon M \to M$, which is certainly transverse to any $r_V$.
\end{proof}

For the next lemmas, recall \cref{D: Delta}.
We also make the following definition.

\begin{definition}
	Let $V \in PC^*_\Gamma(M)$.
	Below we write $M \times V$ for the element of $PC^*_\Gamma(M \times M)$ given by the co-oriented exterior product of $V$ with the identity $\id_M \colon M \times M$ given its tautological co-orientation.
\end{definition}

\begin{lemma}\label{L: image of cubical cap}
	Suppose $M$ is a manifold without boundary with cubulation $X$.
	Let $V \in PC^*_{\Gamma\pf}(M)$ represent a cocycle, and let $W \in K_*(X)$.
	Suppose $\id_M \times r_V \colon M \times V \to M \times M$ is transverse to $\mc J(\Delta(W))$.
	Then $$\mc J(\mc I(V)\frown W) = \pi_1( (M \times V)\times_{M \times M}\mc J(\Delta(W))) \in C_*^\Gamma(M),$$
	where $\pi_1 \colon M \times M \to M$ is the projection to the first factor.
\end{lemma}

\begin{proof}
	For a cubical face $E$ representing an element of $K_*(X)$, let us write $\xi(\diag (E)) = \sum_i E_{1i} \otimes E_{2i}$, analogously to Sweedler notation.
	By definition, at the chain/cochain level $\mc I(V)\frown E$ is given by
	$$(1 \otimes \mc I(V))(\xi(\diag (E))) = (1 \otimes \mc I(V))\left(\sum_i E_{1i} \otimes E_{2i}\right) = \sum_i E_{1i} \otimes \mc I(V)(E_{2i}) = \sum_i \mc I(V)(E_{2i})\cdot E_{1i},$$
	where $\mc I(V)(E_{2i})$ is the intersection number of $V$ with $E_{2i}$ by \cref{D: intersection homomorphism}.
	So $\mc J(\mc I(V)\frown W)$ is just the geometric cochain represented by $\sum_i I(V,E_{2i})E_{1i} = \sum_i \aug(V \times_M E_{2i})E_{1i} $, identifying the cubical face $E_{1i}$ with its embedding into $M$.
	Note that we have $I(V,E_{2i}) = 0$ if $V$ and $E_{2i}$ do not have complementary dimension in $M$, so we can take the sum $\sum_i I(V,E_{2i})E_{1i}$ to be over those $i$ such that $E_{2i}$ has complementary dimension to $V$.

	On the other hand, $\mc J(\Delta(E))$ is the geometric chain corresponding to $\sum_i E_{1i} \times E_{2i}$, and, applying our transversality assumption, we have
	\begin{align*}
		(M \times V)&\times_{M \times M}\mc J(\Delta(E))\\
		& = (M \times V)\times_{M \times M}\left(\sum_i E_{1i} \times E_{2i}\right)\\
		& = \sum_i (M \times_M E_{1i}) \times (V \times_M E_{2i})&\text{by \cref{P: cap cross}}\\
		& = \sum_i E_{1i} \times (V \times_M E_{2i})&\text{by \cref{P: cap with 1}}.
	\end{align*}
	Note that the signs in the formula for \cref{P: cap cross} all vanish in this setting.

	We now consider cases depending on the dimension of $V \times_M E_{2i}$.
	If $\dim(V)+\dim(E_{2i})<\dim (M)$, then $V \times_M E_{2i} = \emptyset$, and the corresponding terms in the above formula vanish.
	Similarly if $\dim(V)+\dim(E_{2i})\geq \dim (M)$ but $V$ and $E_{2i}$ do not intersect.
	For the remaining cases, suppose
	$V \times_M E_{2i}\neq \emptyset$.

	If $V$ and $E_{2i}$ have complementary dimension, then $V \times_M E_{2i}$ is $0$ dimensional, and $\pi_1(
	E_{1i} \times (V \times_M E_{2i}))$ is simply $\aug(V \times_M E_{2i})E_{1i}$.

	If $\dim(V \times_M E_{2i})\geq 2$, then when we take the projection,
	$\pi_1(E_{1i} \times (V \times_M E_{2i}))$ has small rank.
	In this case, $\dim(\bd (V \times_M E_{2i}))\geq 1$ (or is empty) and so
	$$\bd(\pi_1(E_{1i} \times (V \times_M E_{2i}))) = \pi_1(\bd E_{1i} \times (V \times_M E_{2i}))\pm \pi_1(E_{1i} \times \bd(V \times_M E_{2i}))$$
	also has small rank, and so these terms are degenerate and vanish in $C^\Gamma_*(M)$.

	Finally, suppose $\dim(V \times_M E_{2i}) = 1$.
	Then again $\pi_1(E_{1i} \times (V \times_M E_{2i}))$ has small rank, as does the boundary term $\pi_1(\bd E_{1i} \times (V \times_M E_{2i}))$.
	The second boundary summand $\pm\pi_1(E_{1i} \times \bd(V \times_M E_{2i}))$ may not have small rank.
	However, since $\dim(V \times_M E_{2i}) = 1$, it must consist of mappings of circles and compact intervals, and, therefore, its boundary consists of (maps to $M$ of) pairs of oppositely oriented points.
	So $E_{1i} \times \bd(V \times_M E_{2i})$ consists of pairs of oppositely oriented copies of $E_{1i}$ mapping to $M \times M$, and once we project via $\pi_1$, these pairs become trivial elements of $C^\Gamma_*(M)$.
	So $\pi_1(E_{1i} \times (V \times_M E_{2i}))$ is also degenerate in this case, and these terms are also $0$ in $C^\Gamma_*(M)$.

	We conclude that $\pi_1( (M \times V)\times_{M \times M}\mc J(\Delta(E)))$ as an element of $C_*^\Gamma(M)$ can be represented as the sum $\sum_i \aug(V \times_M E_{2i})E_{1i}$ over only those $i$ with $E_{2i}$ of complementary dimension to $V$.
	But this is the same formula we derived for $\mc J(\mc I(V)\frown E)$.
\end{proof}

\begin{lemma}\label{L: diagonal version of intersection}
	Let $M$ be a manifold without boundary.
	Let $V \in PC^*_\Gamma(M)$ and $W \in PC_*^\Gamma(M)$ be transverse.
	Let $\pi_1 \colon M \times M \to M$ be the projection to the first factor.
	Then $$V \times_M W = \pi_1((M \times V)\times_{M \times M} \diag(W)).$$
\end{lemma}

\begin{proof}
	By \cref{L: M times transverse diag}, $\id_M \times r_V \colon M \times V \to M \times M$ is transverse to $\diag r_W \colon W \to M \times M$ in $M \times M$, so both expressions are defined.
	We also have $V = M \times_M V = \diag^*(M \times V)$ by \cref{C: cup with identity} and \cref{P: cross to cup}.
	So we can compute
	\begin{align*}
		V \times_M W& = \pi_1\diag (V \times_M W)&\text{since $\pi_1\diag = \id_M$}\\
		& = \pi_1\diag (\diag^*(M \times V) \times_M W)&\text{by the above}\\
		& = \pi_1((M \times V)\times_{M \times M}\diag(W))&\text{by naturality of cap products.}
	\end{align*}
	For the last equality, see \cref{P: natural cap} and its interpretation in terms of naturality of the cap product in \cref{S: (co)chain properties}.
	\cref{P: natural cap} requires $\id_M \times r_V \colon M \times V \to M \times M$ to be transverse to $\diag \colon M \to M \times M$ and $W \to M$ to be transverse to the pullback of $M \times V$ by $\diag \colon M \to M \times M$ to
	$(M \times V)\times_{M \times M}M \to M$.
	The first requirement holds by \cref{L: M times transverse diag}.
	For the second transversality requirement, \cref{L: transverse to pullback} says that in the presence of the first transversality condition, this is equivalent to requiring $\id_M \times r_V \colon M \times V \to M \times M$ to be transverse to $\diag r_W \colon W \to M \times M$.
	But this also holds by \cref{L: M times transverse diag} as $V$ and $W$ are transverse.
\end{proof}

\begin{lemma}\label{L: diagonal equivalence}
	Let $M$ be a manifold without boundary with cubulation $X$.
	Let $W$ be a cycle in $K_*(X)$, and let $V \in PC^*_\Gamma(M)$ represent a cocycle such that $\id_M \times r_V \colon M \times V \to M \times M$ is transverse to $\diag(\mc J(W))$ and $\mc J(\Delta(W))$.
	Then $$\underline{\pi_1((M \times V) \times_{M \times M}\diag(\mc J(W)))} = \underline{\pi_1( (M \times V)\times_{M \times M}\mc J(\Delta(W)))} \in H_*^\Gamma(M).$$
\end{lemma}

\begin{proof}
	If we consider $W$ as an element of $NK_*^{sm}(M)$, then the geometric chain $\diag(\mc J(W))$ is represented by the singular cubical chain $\diag(W)$ and $\mc J(\Delta(M))$ is represented by the singular cubical chain $\zeta\xi\diag(W)$.

	As $\zeta\xi: NK_*^{sm}(M) \to NK_*^{sm}(M)$ is chain homotopic to the identity by \cite[Section XI.5]{Mas91}, $\diag(W)$ and $\zeta\xi\diag(W)$ must be homologous in $NK_*^{sm}(M)$, and so they are also homologous as geometric chains.
	In particular, $\diag(\mc J(W))$ and $\mc J(\Delta(M))$ represent the same element of $H_*^\Gamma(M \times M)$.

	As $V$ represents a cocycle, so does $M \times V$ by \cref{P: boundary of exterior product,L: exterior Q}.
	It now follows from \cref{T: (co)homology products} that
	$(M \times V)\times_{M \times M} \diag(\mc J(W))$ and $(M \times V)\times_{M \times M} \mc J(\Delta(W))$ represent the same geometric homology class, and so their images under $\pi_1$ represent the same geometric homology class.
\end{proof}

\begin{proof}[Proof of \cref{T: equivalent cap}]
	Let us first choose a cubical cycle $W \in K_*(X)$ representing our given cubical homology class.
	By \cref{L: product transversal}, we can choose a representative $V$ of our geometric cohomology class such that $V$ is transverse to the cubulation (and hence to $W$) and $M \times V$ is transverse to $\diag(\mc J(W)) \sqcup \mc J(\Delta(W))$, which is also represented by a union of embeddings.
	Then by \cref{L: image of cubical cap}, we have
	$$\underline{\mc J(\mc I(V)\frown W)} = \underline{\pi_1( (M \times V)\times_{M \times M}\mc J(\Delta(W)))} \in H_*^\Gamma(M),$$
	and by \cref{L: diagonal equivalence} this equals $\underline{\pi_1((M \times V)\times_{M \times M} \diag(\mc J(W)))}$.
	Then by \cref{L: diagonal version of intersection},
	$\pi_1((M \times V)\times_{M \times M} \diag(\mc J(W))) = V \times_M \mc J(W) \in PC_*^\Gamma(M)$.
	Finally, $V \times_M \mc J(W)$ represents $\uV\nplus \mc J(W)$ by definition.
\end{proof}

\subsubsection{Poincar\'e duality}\label{S: PD}

In \cref{T: PD}, we noticed that geometric homology and cohomology satisfy a very strong form of Poincar\'e duality, as for a closed oriented manifold $M$ we in fact have chain-level identities $C^{m-i}_\Gamma(M) = C_i^\Gamma(M)$ obtained by identifying co-oriented cochains with their corresponding oriented chains, using the orientations induced by the orientation of $M$.
\cref{T: equivalent cap} allows us to observe that this strong version of geometric Poincar\'e is compatible with the classical Poincar\'e duality:

\begin{corollary}[Poincar\'e duality]\label{C: PD}\index{Poincar\'e Duality!via cap product}
	Let $M$ be a closed oriented manifold with cubulation $X$.
	Let $\underline M \in C_m^\Gamma(M)$ be represented by the orientation-preserving identity map $\id_M \colon M \to M$, and let $[M] \in K_*(X)$ represent the cubical fundamental class.
	Then there is a commutative diagram of isomorphisms
	\[
	\begin{tikzcd}
		H^{n-i}_\Gamma(M) \arrow{r}{\nplus \uM} \arrow[d, "\mc I","\cong"'] & H_i^\Gamma(M) \\
		H^{n-i}(X) \arrow{r}{\frown [M]} & H_i(X). \arrow[u, "\mc J"',,"\cong"]
	\end{tikzcd}
	\]
\end{corollary}

The corollary follows immediately from \cref{T: equivalent cap} and the following lemma.

\begin{lemma}
	Let $M$ be closed, oriented, cubulated, and connected.
	Then $\mc J([M]) = \uM \in H_m^\Gamma(M)$.
\end{lemma}

\begin{proof}
	It suffices to prove the lemma for $M$ connected.
	Let $\uV \in H^m_\Gamma(M)$ be represented by a map $V = pt \into M$ taking the point to the center of an $m$-cube of the cubulation, co-oriented so that its normal co-orientation agrees with the orientation of $M$.
	By \cref{P: cap with identity M}, as $\uM$ and $\mc J([M])$ are both represented by embeddings with the same orientation in a neighborhood of the embedded point $V$, the cap products $\uV\nplus \uM$ and $\uV\nplus \mc J([M])$ in $H_0^\Gamma(M)$ are each represented by the same point with its induced orientation (which by \cref{P: cap of immersions} will be the positive orientation).
	This is a generator of $H_0^\Gamma(M) \cong \Z$, as we can see, for example, via our homology isomorphism $H_*(NK_*(M)) \to H_*^\Gamma(M)$.
	As $H^m_\Gamma(M) \cong H_0^\Gamma(M) \cong \Z$ by the isomorphisms between geometric and singular homology and cohomology, $\uV\nplus \colon H^m_\Gamma(M) \to H_0^\Gamma(M)$ must be injective as we have shown it is not the $0$ map.
	Since we have shown $\uV\nplus\uM = \uV\nplus\mc J([M])$, we have $\uM = \mc J([M])$.
\end{proof}

With $\nplus \uM$ as our Poincar\'e duality map, the relation of \cref{P: compare cup and intersection orientations}, which in \cref{S: mixed formulas} became the chain/cochain formula
$$(\uV\uplus \uW)\nplus \uM = (-1)^{(m-v)(m-w)}(\uV\nplus \uM)\bullet(\uW\nplus \uM) = (\uW\nplus \uM)\bullet(\uV\nplus \uM),$$
demonstrates the usual relationship between intersection products and cup products that is well known for homology classes represented by embedded manifolds, cf.\ \cite[Section VI.11]{Bred97}.
Here we see that this relationship extends not just for intersections of embedded manifolds but to all homology classes.
Of course this is always possible if one takes the above formula as a defining formula for the intersection product, but here we see that the intersection product can always be defined geometrically in terms of fiber products.

\subsubsection{Umkehr maps}\label{S: umkehr}\index{umkehr maps|(}

\Cref{C: PD} allows us to make some remarks about umkehr maps, also known as wrong-way or transfer maps, associated to maps of closed oriented manifolds $f \colon N \to M$.
These are maps
\begin{align*}
	f^! \colon H^{n-i}_\Gamma(N) \to H^{m-i}_\Gamma(M)\\
	f_! \colon H_{m-i}^\Gamma(M) \to H_{n-i}^\Gamma(N),
\end{align*}
typically defined by taking a homology or cohomology class, dualizing using Poincar\'e duality, applying $f$ or $f^*$, and then dualizing again; see \cite[Definition VI.11.2]{Bred97}.
We will show that when $M$ and $N$ are closed and oriented, these transfer maps correspond to the pullbacks and pushforwards already encountered in \cref{S: functoriality}, where we only required for cohomology pushforwards that $f$ be proper and co-oriented and for homology pullbacks that $f$ be proper and that $M$ and $N$ be oriented.

\begin{proposition}\index{Poincar\'e Duality!and umkehr maps}\index{umkehr maps!and Poincar\'e Duality}
	Let $f \colon N \to M$ be a map of closed oriented manifolds.
	We may consider $f$ co-oriented via the orientations of $M$ and $N$.
	Then the following diagrams commute:
	\[
	\begin{tikzcd}
		H^{n-i}_\Gamma(N) \arrow{r}{f} \arrow[d, "\nplus \uN"] & H^{m-i}_\Gamma(M) \arrow[d, "\nplus \uM"] & H^{i}_\Gamma(M) \arrow{r}{f^*} \arrow[d, "\nplus \uM"] &[.75in] H^{i}_\Gamma(N) \arrow[d, "\nplus \uN"] \\
		H_i^\Gamma(N) \arrow{r}{f} & H_i^\Gamma(M) & H_{m-i}^\Gamma(M) \arrow[r,"(-1)^{i(m-n)}f^*"] & H_{n-i}^\Gamma(N).
	\end{tikzcd}
	\]
\end{proposition}

\begin{proof}
	We start with the diagram on the left.
	Let $\uV \in H^{n-i}_\Gamma(N)$ be represented by a co-oriented map $r_v \colon V \to N$.
	Then $\uV\nplus\uN$ is represented by the same map to $N$ with its induced orientation; see \cref{S: co-orientations}.
	In particular, if $x \in V$ then $V$ is oriented at $x$ by the local orientation $\beta_V$ such that $(\beta_V,\beta_N)$ gives the co-orientation of $r_V$.
	The path down then right is then the composition $fr_V$, considering $V$ with its orientation $\beta_V$.
	On the other hand, by the definition in \cref{S: covariant functoriality}, the element $f(\uV)$ in $H^{m-i}_\Gamma(M)$ is represented by $fr_V$ co-oriented by composing the co-orientations of $r_V$ and $f$.
	So if the co-orientation of $r_V$ is again $(\beta_V,\beta_N)$, the co-orientation of $fr_V$ representing $f(\uV)$ is $(\beta_V,\beta_N)*(\beta_N,\beta_M) = (\beta_V,\beta_M)$.
	So $f(\uV)\nplus \uM$ is represented by $fr_V$ with $V$ oriented again by $\beta_V$.
	Thus the diagram commutes.

	For the second diagram, let $r_V \colon V \to M$ represent $\uV \in H^{i}_\Gamma(M)$.
	We can assume up to a homotopy that $f$ is smooth and transverse to $r_V$.
	Then $f^*(\uV)$ is represented by the co-oriented pullback $V \times_M N \to N$, and $f^*(\uV)\nplus\uN$ is represented by the same map, now with the orientation on $V \times_M N$ induced by the pullback co-orientation and the orientation of $N$.
	By the naturality of cap products, if we consider $f(f^*(\uV) \cap \uN)$, this is equal to $\uV \nplus f(\uN)$.
	As the pushforward $f \colon H_*^\Gamma(N) \to H_*^\Gamma(M)$ is determined by composition with $f$ (see \cref{S: covariant functoriality}), this means that the composition right then down in the diagram is represented by $V \times_M N$ with its cap orientation, thinking of $f \colon N \to M$ as representing an element of $C_*^\Gamma(M)$.
	If we consider $f \colon N \to M$ with the co-orientation $\hat N$ induced by the orientations of $M$ and $N$, then this prechain $f \colon N \to M$ is the cap orientation of $\hat N \times_M M$, again by \cref{P: cap with identity M}, so $V \times_M N$ with its cap orientation is $V \times_M (\hat N \times_M M)$, which by \cref{P: OC mixed associativity} is $(V \times_M \hat N) \times_M M$.
	In other words, going right then down gives us the chain represented by the co-oriented fiber product of $V$ and $N$, followed by taking its induced orientation, i.e.\ it is $(V \times_M \hat N)\,\check{\vrule height1.3ex width0pt}$.
	This is precisely the product that we denoted by $\check V \times^c_M N$ in \cref{S: c vs o}.
	On the other hand, the map down then right first gives $V$ its induced orientation $\check V$ by \cref{P: cap with identity M} and then forms $(-1)^{i(m-n)}$ times the oriented fiber product $\check V \times_M N$ by \cref{D: cohomology pullback and homology transfer}.
	This is precisely $(-1)^{i(m-n)}$ times the product that we denoted by $\check V \times^o_M N$ in \cref{S: c vs o}.
	By \cref{P: compare cup and intersection orientations},
	$$\check V \times_M^o N = (-1)^{(m-v)(m-n)} \check V \times_M^c N,$$
	hence the sign, noting that $i = m-v$.
\end{proof}
\index{umkehr maps|)} 

\section{Alexander duality and link theory}\label{S: Alexander}

In this section, as a sample application of geometric cohomology, we demonstrate some computations in link theory.
The first step is a proof of Alexander Duality in the language of geometric homology and cohomology.
This allows immediate applications to link theory, as Seifert surfaces of link components with their induced co-orientations are precochains in the link the complement.

\subsection{Alexander Duality}\index{Alexander Duality|(}

We saw in \cref{T: PD} that Poincar\'e duality takes a very strong chain-level form in the world of geometric chains and cochains.
It is not clear that we can provide quite the same chain-level version of Alexander duality, but nonetheless the Alexander duality isomorphisms can be expressed in geometric homology and cohomology by pleasingly geometric maps, at least when our subspaces of $S^n$ are taken to be embedded submanifolds.
This is the case, for example, when studying classical link theory.

So for this section let us assume an underlying sphere $S^n$ with its standard orientation and a subspace $L$ that is the union of finitely-many disjoint smoothly embedded closed manifolds $L_i$, not necessarily of the same dimension.
To avoid too many special cases, we will assume $n>1$, leaving the $S^1$ case to the reader.
In this context, we will prove the Alexander duality statement that $\td H_i^\Gamma(L) \cong \td H^{n-i-1}_\Gamma(S^n-L)$ for all $i \geq 0$.
The standard singular (co)homology proof proceeds by using relative homology, Lefschetz duality, excision, and some homotopy equivalences. We will instead define two mutually inverse maps.

Throughout the proof, we will use the orientation of $S^n$ to convert prechains to precochains by converting orientations to their induced co-orientations, and vice versa.
To facilitate this, recall \cref{N: hat check}:
If $V$ is a prechain, we let $\hat V$ denote the corresponding precochain so that $V = \hat V \times_{S^n} S^n$ by \cref{P: cap with identity M}; if $V$ is a precochain, we let $\check V$ be the corresponding prechain so that $\check V = V \times_{S^n} S^n$.
We think of the inverse operations $\hat{\phantom{a}}$ and $\check{\phantom{a}}$ as respectively raising and lowering the index.
If $V$ is replaced by a complicated expression, we sometimes write $V\,\hat{\vrule height1.3ex width0pt}$ or $V\,\check{\vrule height1.3ex width0pt}$.
By \cref{L: Q switch}, the operations $\hat{\phantom{a}}$ and $\check{\phantom{a}}$ preserve membership in $Q(M)$.
We also note, using \cref{P: Leibniz cap}, that if $V$ is a prechain then
$$\bd V = \bd(\hat V \times_{S^n} S^n) = (-1)^{v} (\bd \hat V) \times_{S^n} S^n,$$
and so $(\bd V)\hat{\vrule height1.3ex width0pt} = (-1)^{v} \bd \hat V$.
Similarly, if $V$ is a precochain, $$\bd(\check V) = \bd(V \times_{S^n} S^n) = (-1)^{v} (\bd V) \times_{S^n} S^n = (-1)^v (\bd V)\check{\vrule height1.3ex width0pt}.$$

We also commit the following abuse of notation.
Suppose $U$ is an open subset of $S^n$ and $W$ is a prechain supported in the interior of $U$.
Then we may of course consider $W$ as a prechain of either $U$ or $S^n$.
Furthermore, if $V$ is a precochain of either $S^n$ or $U$, then, we may write either $V \times_{S^n} W$ or $V \times_U W$ without ambiguity to result in the same prechain in $U$ and hence also in $S^n$.

\subsubsection{Homology to cohomology} We first define the map from homology to cohomology.

Let $V$ be a prechain representing an element of $\td H_i^\Gamma(L)$.
We first observe that $V$ represents the $0$ element of $H_i^\Gamma(S^n)$ when we consider $L$ as a subset of $S_n$.
This is trivial when $i \neq 0,n$ as then $H_i^\Gamma(S^n)=0$.
When $i = n$, we have $H_n^\Gamma(L) = 0$ by \cref{E: dimension range}, as all the manifold components of $L$ must have dimension $<n$, and so $V$ represents $0$ in $H_n^\Gamma(S^n)$.
For $i=0$, we recall from \cref{R: reduced h} that the generators of $\td H_0^\Gamma(L)$ can be represented by prechains $g_{ij}$, each consisting of a pair of oppositely-oriented points in separate components of $L$.
As we can connect such a pair with a path in $S^n$, the $g_{ij}$ each represent $0$ in $H_0^\Gamma(S^n)$, as well.

So each prechain $V$ that represents an element of $\td H_i^\Gamma(L)$ also represents the $0$ homology class in $S^n$. By \cref{R: cycles and boundaries}, there is thus some $Z \in PC^\Gamma_{i+1}(S^n)$ such that $\bd Z \sqcup -V \in Q_*(S^n)$.
We consider $\hat Z \in PC_\Gamma^{n-i-1}(S^n)$, which is $Z$ with its induced co-orientation, and then let $\mr Z \in PC_\Gamma^{n-i-1}(S^n - L)$ be the pullback of $\hat Z$ to $S^n-L$, namely $\hat Z \times_{S^n} S^n - L$.
As a space, $\mr Z = r_Z^{-1}(S^n-L)$, and its reference map, which we will write $\mr r_Z$, is the restriction of $r_Z$ to be a map from $\mr Z$ to $S^n-L$.

We next show that $\mr Z$ is a precocycle.
By \cref{leibniz}, we have $\bd{\mr Z} = \bd (\hat Z \times_{S^n} (S^n - L)) = (-1)^{i+1} (\bd \hat Z) \times_{S^n} (S^n - L)$.
On the other hand, we have $\bd Z \sqcup -V \in Q_*(S^n)$, and so $(\bd Z \sqcup -V)\hat{\vrule height1.3ex width0pt} = (\bd Z)\hat{\vrule height1.3ex width0pt} \sqcup - \hat V \in Q^*(S^n)$ by \cref{L: Q switch}.
When we pull back to $S^n - L$, we obtain an element of $Q^*(S^n - L)$ by \cref{L: pullback with Q}.
But $V$ is contained in $L$, so its pullback is empty, and the remaining piece of the pullback is again
$(\bd Z)\hat{\vrule height1.3ex width0pt} \times_{S^n} (S^n - L) = (-1)^{i+1} (\bd \hat Z) \times_{S^n} (S^n - L) = \bd{\mr Z}$.
So $\bd{\mr Z} \in Q^*(S^n - L)$, and $\mr Z$ represents an element of $\td H^{n-i-1}_\Gamma(S^n - L)$.

We will show that this construction gives a well-defined homomorphism $\mc Z \colon \td H_i^\Gamma(L) \to \td H^{n-i-1}_\Gamma(S^n-L)$.

\begin{proposition}
	The map $\mc Z \colon \td H_i^\Gamma(L) \to \td H^{n-i-1}_\Gamma(S^n-L)$ is a well-defined homomorphism for all $i \geq 0$.
	In particular, it does not depend on the choices involved in the definition.
\end{proposition}
\begin{proof}
	Once we show that $\mc Z$ is independent of the specific choices of $V$ (representing a given $\uV \in \td H_i^\Gamma(L)$) and of $Z$, it will follow immediately that $\mc Z$ is a homomorphism, as given $V,W\in PC_i^\Gamma(L)$ representing homology classes and $Z$ and $Y$ as constructed above from $V$ and $W$ respectively, then we may represent both $\mc Z(\uV + \uW)$ and $\mc Z(\uV) + \mc Z(\uW)$ by $\mr Z \sqcup \mr Y$.

	So suppose that $V,W \in PC_i^\Gamma(L)$ represent the same homology class $\uV$ in $\td H_i^\Gamma(L)$, possibly with $V = W$.
	Let $Z,Y \in PC_{i+1}^\Gamma(S^n)$ be such that $\bd Z \sqcup -V$ and $\bd Y \sqcup -W$ are in $Q_*(S^n)$.
	As $\uV=\underline{W} \in \td H_i^\Gamma(L)$, there is an $S \in PC_{i+1}^\Gamma(L)$ such that $\bd S \sqcup -V \sqcup W \in Q_*(L)$.
	We note that prechains in $Q_*(L)$ also represent elements of $Q_*(S^n)$.
	So now consider $S \sqcup -Z \sqcup Y$.
	By taking boundaries and rearranging, we have $$\bd (S \sqcup -Z \sqcup Y)\sqcup -V \sqcup W \sqcup V \sqcup -W = (\bd S \sqcup -V \sqcup W)\sqcup -(\bd Z \sqcup -V) \sqcup (\bd Y\sqcup -W)\in Q_*(S^n).$$
	But $-V \sqcup W \sqcup V \sqcup -W$ is trivial, so by \cref{L: Lipy12} we have $\bd(S \sqcup -Z \sqcup Y)\in Q_*(S^n)$, i.e.\ $S \sqcup -Z \sqcup Y$ represents an $i+1$ cycle by \cref{R: cycles and boundaries}.

	Next, suppose $i+1\neq n$, we have $H_{i+1}^\Gamma(S^n) = 0$, so there exists $R \in PC_{i+2}^\Gamma(S^n)$ such that $\bd R \sqcup -(S \sqcup -Z \sqcup Y) \in Q_*(S^n)$.
	When we convert this last element to a precochain, we obtain $(-1)^{i+2}\bd \hat R \sqcup -(\hat S \sqcup -\hat Z \sqcup \hat Y) \in Q^*(S^n)$.
	Pulling back to $S^n - L$ and using that $S$ is supported in $L$ and that $S^n - L$ has no boundary, we get
	$(-1)^{i+2}\bd \mr R \sqcup \mr Z \sqcup -\mr {Y} \in Q^*(S^n - L)$.
	Thus $\mr Z$ and $\mr {Y}$ represent the same cohomology class, again by \cref{R: cycles and boundaries}.

	Next, suppose $i+1 = n$; this is the case where our map is $\mc Z \colon \td H_{n-1}^\Gamma(L) \to \td H^0_\Gamma(S^n - L)$.
	We know $H_{i+1}^\Gamma(S^n) = H_{n}^\Gamma(S^n) \cong \Z$, and by \cref{E: sphere homology} we can take the identity map $\id \colon S^n \to S^n$ to be a generator (which, as usual, we represent by its domain $S^n$).
	So our cycle represented by $S \sqcup -Z \sqcup Y$ must be homologous to $k S^n$ for some $k \in \Z$; note that here $kS^n$ can be thought of as $|k|$ copies of $S^n$ or $-S^n$ mapping to $S^n$ by the identity map.
	By \cref{R: cycles and boundaries} there exists $R \in PC_{n+1}^\Gamma(S^n)$ such that $\bd R \sqcup -(S \sqcup -Z \sqcup Y) \sqcup kS^n \in Q_*(S^n)$.
	This time, the corresponding precochain is $(-1)^{n+1}\bd \hat R \sqcup -(\hat S \sqcup -\hat Z \sqcup \hat{Y}) \sqcup k \hat S^n \in Q^*(S^n)$, which pulls back to
	$(-1)^{n+1}\bd \mr R \sqcup \mr Z \sqcup -\mr Y \sqcup k\mr S^n \in Q^*(S^n - L)$.
	But $\mr S^n$ is the tautologically co-oriented identity map of $S^n - L$, which represents $0$ in reduced cohomology $\td H^0_\Gamma(S^n - L)$ by \cref{D: reduced c}.
	Thus $\mr Z$ and $\mr {Y}$ represent the same reduced cohomology class, as required.
\end{proof}

\subsubsection{Cohomology to homology}
Next we'll define a map $\mc A \colon \td H^{n-i-1}_\Gamma(S^n-L) \to \td H_i^\Gamma(L)$ that we will later show is an inverse to $\mc Z$.
We first define the nonreduced version $\mc A \colon H^{n-i-1}_\Gamma(S^n-L) \to H_i^\Gamma(L)$, demonstrate it is well defined, and then show that it naturally provides a map also in the cases where reduced (co)homology comes into play.

We suppose $S^n$ is the standard unit sphere in $\R^{n+1}$ with the induced metric; this assumption is convenient to start, though we will see below that it is not critical (see \cref{R: Alex tubular}).
Let $d \colon \R^{n+1} \to [0,\infty)$ be the distance function to $L$.
As a subset of $\R^n$, the subspace $L$ will have an $\epsilon$-neighborhood $Y_\epsilon$ in in the sense of \cite[Section 2.3]{GuPo74}, and so also a retraction $\pi \colon Y_\epsilon \to L$ that takes a point in $Y_\epsilon$ to the nearest point in $L$.
We let $U = Y_\epsilon \cap S^n$ and also write $\pi$ for the restriction to $U$.
We can further extend $\pi$ to a deformation retraction $h \colon U \times I \to L$ by setting
\begin{equation}\label{E: alex retraction}
	h(x,t) = \frac{(1-t)x + t\pi(x)}{|(1-t)x + t\pi(x)|},
\end{equation}
where the addition and scalar division are in $\R^{n+1}$ and $|-|$ is the norm in $\R^{n+1}$.
For small $\epsilon$, the denominator will never vanish, and for $x \in S^n -L$, we have $h(x,t) \in S^n -L$ for $t < 1$.

Next, consider the restriction of the distance-from-$L$ function $d$ to $S^n$, which we will write $\bar d \colon S^n \to [0,\infty)$.
This is smooth away from $0$, and by Sard's Theorem, almost every point of $[0,\infty)$ is a regular value for $\bar d$.
For such a regular value $x$, the subspace $S_x = \bar d^{-1}(x)$ will be a smooth $n-1$ dimensional submanifold of $S^n$ \cite[Section 2.1]{GuPo74}.
If we choose $x < \epsilon$, then $S_x$ will be contained in our $\epsilon$-neighborhood $U$.
Furthermore, $D_x^- = \bar d^{-1}([0,x])$ and $D_x^+ = \bar d^{-1}([x,\infty))$ will be smooth $n$-dimensional submanifolds with boundary of $S^n$.
We assume $D^{\pm}_x$ oriented consistently with $S^n$ and give $S_x$ the boundary orientation $S_x = \bd D_x^+ = -D_x^-$.

Define $\mc A \colon H^{n-i-1}_\Gamma(S^n-L) \to H_i^\Gamma(L)$ as follows.
We choose a regular value $x \in [0,\epsilon)$.
Let $\uV \in H^{n-i-1}_\Gamma(S^n-L)$ be represented by $r_V \colon V \to S^n - L$, which we can assume is transverse to $S_x$ by \cref{P: perturb transverse to map,P: universal homotopy}.
As $V$ is a precochain and $S_x$ is oriented, we obtain an orientation on $V \times_{S^n} S_x$ by \cref{D: PC products}.
Technically the fiber product happens in $S^n-L$ with the resulting prechain supported in $U-L$, but the notation $\times_{S^n}$ should be simpler and unambiguous.
As $V$ has dimension $i+1$, the prechain $V \times_{S^n} S_x$ has dimension $i$.
We define $\mc A(\uV)$ to be represented by $\pi(V \times_{S^n} S_x)$.

Once again, we need to show we have independence of the choices.
The following proposition establishes that $\mc A$ is independent of the choice of $V$ representing $\uV \in H^*_\Gamma(S^n-L)$ and of our choice of a regular value $x$.
It will follow later from \cref{R: Alex tubular} that $\mc A$ is also independent of the exact choice of distance function; in fact we could replace $S_x$ with the sphere bundle of a normal bundle identified with a tubular neighborhood of $L$.

\begin{proposition}
	The map $\mc A \colon H^{n-i-1}_\Gamma(S^n-L) \to H_i^\Gamma(L)$ is a well-defined homomorphism.
	In particular, it does not depend on the choices of $V$ or $x$ in the above definition.
\end{proposition}
\begin{proof}
	Once we show that $\mc A$ is independent of the specific choices of $V$ (representing a given $\uV$ in $H^*_\Gamma(S^n-L)$) and $x$, it will follow immediately that $\mc A$ is a homomorphism, as given $V,W\in PC^*_\Gamma(S^n-L)$ representing cohomology classes and transverse to $S_x$, then $\mc A(\uV + \uW)$ will be represented by $\pi((V \sqcup W) \times_{S^n} S_x)$, which is the same precochain as $\pi(V \times_{S^n} S_x) \sqcup \pi(W \times_{S^n} S_x)$.

	Next, suppose $V' \in PC^*_\Gamma(S^n-L)$ also represents $\uV \in H^*_\Gamma(S^n-L)$ and is transverse to $S_x$.
	As $V$ and $V'$ represent the same cohomology class, by \cref{R: cycles and boundaries} there exists $Z \in PC^*_\Gamma(S^n-L)$ so that $\bd Z\sqcup -V \sqcup V' \in Q^*(S^n-L)$; we may also assume $Z$ transverse to $S_x$ by \cref{T: transverse reps}.
	Then $(\bd Z \sqcup -V \sqcup V') \times_{S^n} S_x \in Q_*(S^n-L)$ by \cref{L: pullback with Q}.
	This prechain can also be considered to be in $Q_*(U) \subset PC_*^\Gamma(U)$, where $U$ is the $\epsilon$-neighborhood described above.
	So then $\pi((\bd Z \sqcup -V \sqcup V') \times_{S^n} S_x ) \in Q_*(L)$ by \cref{L: Q preservation}.
	But, as above, this is the same precochain as $\pi((\bd Z) \times_{S^n} S_x ) \sqcup -\pi(V \times_{S^n} S_x) \sqcup \pi(V' \times_{S^n} S_x)$.
	Furthermore, as $S_x$ is boundaryless, $\pi((\bd Z)) \times_{S^n} S_x) = \pm \bd \pi((Z \times_{S^n} S_x)$.
	So $\pi(V \times_{S^n} S_x)$ and $\pi(V' \times_{S^n} S_x)$ represent the same element of $H_*^\Gamma(L)$, again by \cref{R: cycles and boundaries}, and this shows $\mc A$ is independent of the choice of $V$.

	Next consider $0 < x, y < \epsilon$, for our $\epsilon$-neighborhood of $L$, such that $x$ and $y$ are both regular values of the distance function $\bar d$.
	Without loss of generality, suppose $x < y$, and let $D_{x,y} = \bar d^{-1}([x,y])$, which is a manifold with boundary in $S^n-L$.
	Given our convention for orienting $S_x$ and $S_y$, we have $\bd D_{x,y} = S_x \sqcup -S_y$.
	Let $V \in PC^*_\Gamma(S^n-L)$ represent $\uV \in H^*_\Gamma(S^n - L)$, and we suppose $V$ transverse to $S_x$ and $S_y$ by \cref{P: perturb transverse to map,P: universal homotopy}.
	As $V$ must represent a cocycle, we have $\bd V \in Q^*(S^n - L)$.
	We now compute using \cref{P: Leibniz cap} that
	\begin{align*}
		\bd \pi(V \times_{S^n} D_{x,y}) &= \pi(\bd (V \times_{S^n} D_{x,y}))\\
		&= \pi\left(\left[(-1)^{i+1} (\bd V) \times_{S^n}D_{x,y} \right] \bigsqcup (V \times_{S^n} \bd D_{x,y})\right)\\
		&= \pi\left(\left[(-1)^{i+1} (\bd V) \times_{S^n}D_{x,y} \right] \bigsqcup (V \times_{S^n} S_x)\bigsqcup (V \times_{S^n} -S_y)\right)\\
		&=(-1)^{i+1} \pi((\bd V) \times_{S^n}D_{x,y}) \sqcup \pi(V \times_{S^n} S_x) \sqcup -\pi(V \times_{S^n} S_y).
	\end{align*}
	Since $\bd V \in Q^*(S^n-L)$, we have the first term in $Q_*(L)$, and so pairing the other two terms with their opposite orientations,
	$$\bd \pi(V \times_{S^n} D_{x,y})\sqcup -\pi(V \times_{S^n} S_x) \sqcup \pi(V \times_{S^n} S_y) \in Q_*(L).$$
	Thus $\pi(V \times_{S^n} S_x)$ and $\pi(V \times_{S^n} S_y)$ are in the same homology class.
\end{proof}

Finally, let us verify that $\mc A$ determines a map of reduced (co)homology $\td H^{n-i-1}_\Gamma(S^n-L) \to \td H_i^\Gamma(L)$.
We need only check the cases $i = 0$ and $i = n-1$.

First suppose $i = 0$.
Then, as we assume $n > 1$, an element $\uV \in \td H^{n-i-1}_\Gamma(S^n-L) = H^{n-i-1}_\Gamma(S^n-L)$ is represented by a $1$-dimensional precochain $V$.
We compute
$$\bd(V \times_{S^n} D^+_x) = -(\bd V \times_{S^n} D^+_x) \sqcup (V \times_{S_n} S_x).$$
As $V$ represents a cocycle, $\bd V \in Q^*(S^n-L)$, and so $\bd V \times_{S^n} D^+_x \in Q_*(S^n-L) \subset Q_*(S^n)$ by \cref{L: pullback with Q}.
Thus $\bd(V \times_{S^n} D^+_x) \sqcup -(V \times_{S_n} S_x) \in Q_*(S^n)$.
But the only way for a $0$-dimensional prechain to be in $Q_*(S^n)$ is if it consists of pairs of ``cancelling'' points with opposite signs.
As $\bd(V \times_{S^n} D^+_x)$ is the boundary of a compact $1$-manifold, it must also consist of pairs of oppositely-signed points (though they do not necessarily cancel as they do not necessarily have the same image in $S^n$).
But then $V \times_{S_n} S_x$ must also consist of pairs of oppositely-signed points, and hence this is also true of $\pi(V \times_{S_n} S_x)$, which represents $\mc A(\uV)$.
This suffices to imply that $\mc A(\uV) \in \td H_0^\Gamma(L)$.

Next, suppose $i = n-1$ so that $\uV \in \td H^0_\Gamma(S^n-L)$.
Utilizing the preceding proposition and \cref{D: reduced c}, to show that $\mc A$ remains well defined with this domain we must show that $\mc A(S^n-L)$ represents $0 \in H_{n-1}^\Gamma(L)$, where here $S^n - L$ refers to the identity map $\id \colon S^n - L \to S^n - L$ as a precochain with the canonical orientation.
For this, we note that $S_x = - \bd D_x^-$, so by \cref{P: Leibniz cap},
$$\pi((S^n-L) \times_{S^n} S_x) = -\pi((S^n-L) \times_{S^n} \bd D_x^-) = -\bd \pi((S^n-L) \times_{S^n} D_x^-).$$
We conclude $\mc A(S^n - L)$ is a boundary of a prechain in $L$, which suffices.

Altogether, we see that $\mc A$ defines a homomorphism $\td H^{n-i-1}_\Gamma(S^n - L) \to \td H_i^\Gamma(L)$.

\subsubsection{The duality isomorphisms} We now show that $\mc A$ and $\mc Z$ are inverse isomorphisms.

\begin{theorem}\label{T: alex duality}
	For $i \geq 0$, the maps $\mc A \colon \td H^{n-i-1}_\Gamma(S^n-L) \to \td H_i^\Gamma(L)$ and $\mc Z \colon \td H_i^\Gamma(L) \to \td H^{n-i-1}_\Gamma(S^n-L)$ are inverse isomorphisms.
\end{theorem}

\begin{proof}
	First we consider the composition $\mc A \mc Z$.
	So suppose $\uV \in \td H_i^\Gamma(L)$ represented by $r_V \colon V \to L$.
	Then $\mc Z(\uV)$ is represented by $\mr r_Z \colon \mr Z \to S^n-L$, where $Z \in PC_{i+1}^\Gamma(S^n)$ is such that $\bd Z \sqcup -V \in Q_*(S^n)$.
	Recall that we choose an $x \in (0,\epsilon)$ that is a regular value for the distance function to $L$.
	By \cref{T: transverse reps}, we may suppose $Z$ is transverse to $S_x$, and then $\mc A \mc Z(\uV)$ is represented by $\pi(\mr Z \times_{S^n} S_x) = \pi(\hat Z \times_{S^n} S_x)$, which we denote $W$.

	Let $Y = \pi(\hat Z \times_{S^n} D_x^-)$.
	By \cref{P: Leibniz cap} and the above observation that $(\bd Z)\hat{\vrule height1.3ex width0pt} = (-1)^{i+1} (\bd \hat Z)$, we have
	\begin{equation*}
		\bd Y = \bd \pi(\hat Z \times_{S^n} D_x^-)
		= \pi\left(\left[(-1)^{i+1} (\bd \hat Z) \times_{S^n} D_x^-\right] \bigsqcup \hat Z \times_{S^n} (-S_x)\right)\\
		= \pi\left((\bd Z)\hat{\vrule height1.3ex width0pt} \times_{S^n} D_x^- \right) \bigsqcup -W.
	\end{equation*}
	We consider
	\begin{equation*}
		\bd Y \sqcup -V\sqcup W =\pi\left((\bd Z)\hat{\vrule height1.3ex width0pt} \times_{S^n} D_x^- \right) \bigsqcup -W \sqcup -V \sqcup W.
	\end{equation*}
	The terms involving $W$ form a trivial pair, so we consider the other terms $\pi\left((\bd Z)\hat{\vrule height1.3ex width0pt} \times_{S^n} D_x^- \right) \sqcup -V$.
	As $V$ is contained in $L$, we have $V = \pi(\hat V \times_{S^n} D_x^-)$, and so
	$$\pi\left((\bd Z)\hat{\vrule height1.3ex width0pt} \times_{S^n} D_x^- \right) \bigsqcup -V = \pi(((\bd Z)\hat{\vrule height1.3ex width0pt} \sqcup -\hat V) \times_{S^n} D_x^-).$$
	But now $\bd Z \sqcup -V \in Q_*(S^n)$, so $(\bd Z)\,\hat{\vrule height1.3ex width0pt} \sqcup -\hat V \in Q^*(S^n)$, and thus $(\bd Z)\hat{\vrule height1.3ex width0pt} \sqcup -\hat V) \times_{S^n} D_x^- \in Q_*(S^n)$ by \cref{L: pullback with Q}.
	Then $\pi$ takes this to $Q_*(L)$ by \cref{L: Q preservation}.
	So $\bd Y \sqcup -V \sqcup W \in Q_*(L)$, and therefore $\uV = \uW \in \td H_i^\Gamma(L)$.
	We conclude $\uV = \mc A \mc Z(\uV) \in \td H_i^\Gamma(L)$, as claimed.

	We note that this argument holds as written for those indices for which we consider reduced homology or cohomology.

	Next we consider $\uV \in \td H^{n-i-1}_\Gamma(S^n-L)$ represented by $r_V \colon V \to S^n-L$.
	By \cref{P: perturb transverse to map,P: universal homotopy}, we can assume $V$ is transverse to $S_x$ and then $\mc Z \mc A(\uV)$ is represented by $\mr Z$, where $Z \in PC_{i+1}^\Gamma(S^n)$ represents a chain that cobounds the homology class represented by $\pi(V \times_{S^n} S_x)$.
	As above, $\mr Z$ is the restriction of $\hat Z$ to $S^n-L$, where $\hat Z$ is $Z$ with the induced co-orientation.
	We will construct such a $Z$ and show that $\underline{\mr Z}$ is cohomologous to $\uV$.

	We first observe that by composing our distance function to $L$ with an appropriate order-preserving embedding $[0,\infty) \to (-1,1)$ that takes $x$ to $0$, we can split $V$ into $V^+ \sqcup V^-$ as in \cref{S: splitting}, where $V^-$ is the portion of $V$ close to $L$.
	By \cref{T: cohomology creasing}, $\underline{V^+} + \underline{V^-} = \uV$ as cohomology classes.
	Note that splitting $V$ and then taking the induced orientation is equivalent to taking the induced orientation and then splitting, so we can write $\check V^+$ unambiguously for both $(V^+)\check{\vrule height1.3ex width0pt}$ and $(\check V)^+$.
	We also note $\check V^+ = V \times_{S^n} D_x^+$.

	The precochain $V^0$ from this splitting will have the same underlying space (i.e. ignoring (co)orientations) as $V \times_{S^n} S_x$.
	To distinguish the notation, we will use $V^0$ to denote the precochain coming from the splitting of $V$ as in \cref{S: splitting}, and we will write $V_0$ for the prechain $V \times_{S^n} S_x$.
	Then $\pi(V_0)$ represents $\mc A(\uV)$.

	Later, we will need to compare $V_0$ with $V^0$, so we make the following computations.
	First, using \cref{P: Leibniz cap} and that $S_x = \bd D_x^+$, we have
	\begin{align*}
		\bd(V \times_{S^n} D^+_x) &= [(-1)^{i+1} (\bd V)\times_{S^n} D^+_x] \sqcup (V \times_{S^n} \bd D^+_x)\\
		&=[(-1)^{i+1} (\bd V)\times_{S^n} D^+_x] \sqcup V_0.
	\end{align*}
	On the other hand, by \cref{P: OC mixed associativity} and the formulas $D^+_x = \hat D^+_x \times_{S^n} S^n$ and $V^+ = V \times_{S^n} \hat D^+_x$, we see $$V \times_{S^n} D^+_x = V \times_{S^n}(\hat D^+_x \times_{S^n} S^n) = (V \times_{S^n} \hat D^+_x ) \times_{S^n} S^n = V^+ \times_{S^n} S^n.$$
	So we compute using \cref{P: Leibniz cap,L: W0 cochain} that
	\begin{align*}
		\bd(V \times_{S^n} D^+_x) &= \bd(V^+ \times_{S^n} S^n )\\
		&= (-1)^{i+1}\bd (V^+) \times_{S^n} S^n\\
		&= (-1)^{i+1} (V^0 \sqcup (\bd V)^+) \times_{S^n} S^n.
	\end{align*}
	Comparing the corresponding components, we see that $(V^0)\check{\vrule height1.3ex width0pt} = (-1)^{i+1}V_0$.

	We next define $Z$. The basic idea, up to orientation information, is to replace $V$ with $V^+ \sqcup T$, where we have better control over $T$ than $V^-$.
	More specifically, we create $T$ by considering a homotopy from $V_0$ to $\pi(V_0)$.
	For this, recall the deformation retraction $h$ of $U$ to $L$ of \eqref{E: alex retraction}.
	Let $r_{V_0} \colon V_0 \to S^n$ be the reference map for $V_0$, and let $T$ be prechain given by the space $V_0 \times I$ with the product orientation and reference map $r_T(y,t) = h(r_{V_0}(y),t)$.
	Then $h(-,0) = r_{V_0}$ and $h(-,1)$ is the reference map for $\pi(V_0) = \pi(V \times_{S^n} S_x)$, which represents $\mc A(\uV)$.
	By \cref{P: oriented fiber boundary},
	\begin{equation}\label{E: alex cylinder bdry}
		\bd T = \bd (V_0 \times I) = (\bd (V_0) \times I) \sqcup (-1)^{i} (V_0 \times \{1,-0\}) = (\bd (V_0) \times I) \sqcup (-1)^i\pi(V_0)\sqcup (-1)^{i+1} V_0,
	\end{equation}
	as $V_0 \times 1$ maps under $h$ by $\pi$.
	We define $Z = \check V^+ \sqcup (-1)^i T$.

	We now check that $\bd Z \sqcup -\pi(V_0) \in Q_*(S^n)$ so that $\bd \underline{Z} = \mc A(\uV)$ in $C_*^\Gamma(S^n)$, as desired.
	Using \cref{L: W0 chain,P: oriented fiber boundary} in the second line below, as well as our computations above of $\bd T$ and of $\bd \check V^+$, which is just $\bd (V \times_{S^n} D^+_x)$, we have
	\begin{align*}
		\bd Z \sqcup -\pi(V_0) & = \bd \check V^+ \sqcup \bd (-1)^i T \sqcup -\pi(V_0)\\
		&=\left[(-1)^{i+1} (\bd V) \times_{S^n} D_x^+\right] \sqcup V_0 \sqcup
		(-1)^i(\bd (V_0) \times I) \sqcup \pi(V_0)\sqcup -V_0
		\sqcup -\pi(V_0).
	\end{align*}
	We have two trivial pairs $V_0 \sqcup - V_0$ and $ -\pi(V_0) \sqcup \pi(V_0)$ in this expression.
	That leaves the terms involving $(\bd V) \times_{S^n} D_x^+$ and $(\bd V_0) \times I$.
	By assumption that $V$ represents a cocycle, $\bd V \in Q^*(S^n-L)$, so $(\bd V) \times_{S^n} D_x^+ \in Q_*(S^n-L)$ by \cref{L: pullback with Q}.
	Similarly, $$\bd (V_0) = \bd (V \times_{S^n} S_x) = (-1)^{i+1}(\bd V) \times_{S^n} S_x$$ by \cref{P: Leibniz cap}, and so this is in $Q_*(S^n-L)$ by \cref{L: pullback with Q}.
	Furthermore, the homotopy that defines $T$ behaves like a universal homotopy since the underlying homotopy $h$ is of $U \subset S^n$; so $(\bd V_0) \times I \in Q_*(S^n)$ by \cref{L: dessicated homotopy}.
	We conclude that $Z$ has the desired cobounding properties.

	It remains to show that $\mr Z = V^+ \sqcup (-1)^i\mr T$ represents the same cohomology class as $V$.
	We already know that $V$ represents the same class as $V^+ \sqcup V^-$, so instead we compare $\mr Z$ with $V^+ \sqcup V^-$, and we must show there is some precochain $Y$ such that
	$$\bd Y \sqcup -\mr Z\sqcup V^+\sqcup V^- = \bd Y \sqcup -V^+ \sqcup (-1)^{i+1}\mr T \sqcup V^+ \sqcup V^- \in Q^*(S^n-L).$$
	It is immediate that $-V^+ \sqcup V^+ \in Q^*(S^n-L)$, so it suffices to find a $Y$ such that
	$\bd Y \sqcup (-1)^{i+1} \mr T \sqcup V^- \in Q^*(S^n-L)$.
	For this we let the underlying space of $Y$ be $(\mr T \sqcup V^-) \times [0,1)$, and if we define $\rho \colon \mr T \sqcup V^- \to S^n -L$ to be the reference map for $(-1)^{i+1}\mr T \sqcup V^-$, we define a map $H \colon Y \to S^n-L$ by $H(u,t) = h(\rho(u),t)$.
	So $H$ is essentially a universal homotopy that pushes the image of $\mr T \sqcup V^-$ off of $S^n-L$ into $L$, but $H$ is well defined into $S^n-L$ because no point goes to $L$ until time $1$, which is omitted from the domain.
	Also, $H$ is proper because the reference map for $\mr T \sqcup V^-$ is proper and the homotopy moves points closer to $L$ (cf.\ \cite[Corollory I.10.3.7]{Bou98}), and $H$ behaves like a universal homotopy because it is built from a homotopy of $U$.
	We co-orient $H$ as a homotopy as in \cref{S: co-oriented homotopy} so that it is a homotopy from $(-1)^{i+1}\mr T \sqcup V^-$; in particular, its boundary will include the term $(-1)^{i}\mr T \sqcup - V^-$.
	So we can now compute
	\begin{align*}
		\bd Y \sqcup (-1)^{i+1} \mr T \sqcup V^-
		&= [(-1)^{i}\mr T \sqcup - V^-] \sqcup \mc H \sqcup [(-1)^{i+1} \mr T \sqcup V^-],
	\end{align*}
	where $\mc H$ is a one-ended homotopy of $\bd[ (-1)^i\mr T \sqcup - V^-]$ by \cref{D: co-oriented homotopy,L: co-oriented homotopy}.
	The two bracketed terms give us a trivial pair, so we consider
	$\mc H$.

	It suffices to show that $\bd[(-1)^i\mr T \sqcup - V^-] \in Q^*(S^n-L)$ as then the homotopy will be in $Q^*(S^n-L)$ by \cref{L: dessicated homotopy}.
	We recall from \cref{L: W0 cochain} that $\bd(V^-) = -(V^0) \sqcup (\bd V)^-$.
	Meanwhile $\mr T$ is the restriction to $S^n - L$ of $\hat T$, and using \eqref{E: alex cylinder bdry},
	\begin{align*}
		\bd (\hat T) &= (-1)^{i+1}(\bd T)\hat{\vrule height1.3ex width0pt}\\
		&=(-1)^{i+1}[(\bd (V_0) \times I) \sqcup (-1)^i\pi(V_0)\sqcup (-1)^{i+1} V_0]\,\hat{\vrule height1.3ex width0pt}\\
		&=[(-1)^{i+1}(\bd (V_0) \times I) \sqcup -\pi(V_0)\sqcup V_0]\,\hat{\vrule height1.3ex width0pt}.
	\end{align*}
	The terms of $\bd[ (-1)^i\mr T \sqcup - V^-]$ that come from these computations and involve $\bd V$ and $\bd V_0$ are all in $Q^*(S^n-L)$ by the assumption that $V$ represents a cocycle and \cref{L: pullback with Q,L: W0 cochain,L: dessicated homotopy}.
	The $\pi(V_0)$ term has image in $L$ and so vanishes when we restrict to $S^n - L$.
	Finally, the remaining terms are $(-1)^i(V_0)\hat{\vrule height1.3ex width0pt} = -V^0$ from the $T$ term and $V^0$ from the $V^-$ term, so these also form a trivial pair.

	This completes the argument that $\underline{\mr Z} = \uV \in \td H^{n-i-1}(S^n - L)$.
	Again we note that this argument holds as written for those indices for which we consider reduced homology or cohomology, so this completes the proof of the theorem.
\end{proof}

\begin{remark}\label{R: Alex tubular}
	So far for definiteness we have assumed $S^n$ to be the unit sphere in $\R^{n+1}$ and made use of an $\epsilon$ neighborhood of $L$.
	However, this is not always convenient to work with, and we now provide an alternative.
	Rather than the $\epsilon$ neighborhood $U$ of $L$, we can instead choose disjoint tubular neighborhoods of the components that collectively form a neighborhood $\mc U$ of $L$.
	We can identify the tubular neighborhoods with normal bundles and replace $S_x$ with an appropriate sphere bundle $\mc S$.
	As the sphere bundles, and their associated disk bundles, also permit deformation retractions to $L$, the constructions and (co)homological arguments above go through in this setting in complete analogy with the preceding proofs.
	As inverses to isomorphisms are unique, we establish that the map $\mc A$ constructed using tubular neighborhoods is identical to the map $\mc A$ as we originally defined it.
	Similarly, we could proceed with a different distance function on $S^n$ using an alternative retraction of a neighborhood of $L$ to $L$.
\end{remark}\index{Alexander Duality|)}

\subsection{Applications to spherical links}\index{Alexander Duality!application to links|(}

As a sample application of the geometric cohomology approach to Alexander duality, we consider linking of a union of embedded spheres $S^{n-2}$ in $S^n$, $n \geq 3$.
Of course, the classical case is embedded circles in $S^3$.

Here, the only nontrivial reduced homology groups of our link $L = \amalg_{a \in A} L_a$ are $\td H_0^\Gamma(L) \cong \Z^{|A|-1}$ and $\td H_{n-2}^\Gamma(L) \cong \Z^{|A|}$.

As observed in \cref{R: reduced h}, $\td H_0^\Gamma(L)$ is generated by homology classes represented by prechains $g_{ab}$, $a \neq b$, consisting of a negative point in component $L_a$ and a positive point in component $L_b$, subject to the relations $g_{ab}+g_{bc}=g_{ac}$.
Applying $\mc Z$, we see that $H^{n-1}_\Gamma(S^n - L)$ is generated by the restrictions to $S^n-L$ of maps $\gamma_{ab} \colon [0,1] \to S^n$, such that $\gamma_{ab}(0) \in L_a$ and $\gamma_{ab}(1) \in L_b$, given their induced co-orientations and subject to the same relations.

On the other hand, $\td H_{n-2}^\Gamma(L) = H_{n-2}^\Gamma(L)$ is generated by the embedding maps from $S^{n-2}$ to each component $L_a$ of $L$.
Recall that any knot $K \cong S^{n-2}$ in $S^n$ possesses a Seifert surface, i.e.\ a compact oriented embedded manifold whose boundary is $K$ \cite[Appendix B]{MC17}.
Seifert surfaces are usually assumed to be connected, but that will not be necessary for us.
If we let $S_a$ denote a Seifert surface for the component $L_a$ of $L$, then the corresponding generator of $H^1_\Gamma(S^n - L)$ is just the class of $\mr S^a$, the restriction to $S^n - L$ of $S_a$ with its induced co-orientation.

So, in summary, the generators of $H^1_\Gamma(S^n - L)$ can be represented by Seifert surfaces, while the elements of $H^{n-1}_\Gamma(S^n - L)$ can be realized as arcs between components.

Of course the real power of geometric cohomology is that we can also derive consequences for the {\it ring} $H^*_\Gamma(S^n-L)$ from geometry.
For example, consider a link in $S^3$.
If it is possible to choose Seifert surfaces $S_a$ simultaneously for all $a$ that are disjoint or intersect only in closed curves, then all cup product on $H^1_\Gamma(S^3-L)$ vanish.
On the other hand, if two Seifert surfaces $S_a$ and $S_b$ intersect in just an arc between $L_a$ and $L_b$, then this corresponds to the cup product of two generators of $H^1_\Gamma(S^3-L)$ being a generator of $H^2_\Gamma(S^3-L)$.
We develop this idea further in the next example.

\subsubsection{Linking numbers}

Let $L$ be a link of two components, $L_1$ and $L_2$, in $S^3$, and let $S_1$ and $S_2$ be corresponding Seifert surfaces that are transverse to each other.
As we have observed, $\mr S_1$ and $\mr S_2$, which we can interpret as $S_1$ and $S_2$ with their intersections with $L$ removed and co-oriented consistently with their orientations (which are determined by the orientations of $L_1$ and $L_2$), represent generators of $H^1_\Gamma(S^3-L) \cong \Z^2$.
As $\mr S_1$ and $\mr S_2$ are embedded in $S^3$, the fiber product $\mr S_1 \times_{S^3} \mr S_2$ is just the intersection $\mr S_1 \cap \mr S_2$ with its fiber product co-orientation, and this represents the cup product $\underline{\mr S_1} \uplus \underline{\mr S_2} \in H^2_\Gamma(S^3-L) \cong \Z$.
Since $\mr S_1 \times_{S^3} \mr S_2$ is a $1$-manifold with a proper co-oriented map to $S^3-L$, it must consist of closed loops in $S^3-L$ and the interiors of arcs with their endpoints on $L$.
Taking induced orientations, we can consider these to be oriented loops and arcs.
By our preceding computations, the loops represent trivial elements of $H^2_\Gamma(S^3-L)$, as does any arc with both endpoints on the same link component.
So, up to an overall choice of sign corresponding to a choice of generator of $H^2_\Gamma(S^3-L)$, the cup product corresponds to the count of arcs running from $L_1$ to $L_2$ minus the count of the arcs running from $L_2$ to $L_1$.

We claim that, using this choice of generator to identify $H^2_\Gamma(S^3-L)$ with $\Z$, the cup product $\underline{\mr S_1} \uplus \underline{\mr S_2}$ represents the linking number of $L_1$ and $L_2$, where we use the second (of eight!) characterization of linking numbers given by Rolfsen in \cite[Section 5D]{Ro90}.\index{linking number}
By this characterization, the linking number of $L_1$ and $L_2$ can be found by counting with sign the number of intersections of $L_1$ with $S_2$.
More precisely, we orient a bicollar of $S_2$ so that $L_2$ runs counterclockwise when seen from the positive side, and then count an intersection with $+1$ when $L_1$ runs from the negative to positive side of $S_2$ and with a $-1$ when $L_1$ runs from the positive to negative side of $S_2$.

To verify our claim, let us first develop carefully the orientations and co-orientations involved.
We assume that $S^3$ has its standard orientation and the components $L_i$ of $L$ have given orientations $\beta_{L_i}$.
Then, by \cref{Con: oriented boundary}, if $x \in L_i$ and $\mu_i$ is a normal vector to $L_i$ in $T_xS_i$ pointing outward from $S_i$, then the orientation of $S_i$ is such that $\beta_{S_i} = \beta_{\mu_i} \wedge \beta_{L_i}$.
The induced co-orientation on the inclusion of $S_i$ into $S^3$ is then $(\beta_{S_i}, \beta_{S^3}) = (\beta_{S_i}, \beta_{S_i} \wedge \beta_{\nu_i})$, where $\nu_i$ is a vector normal to $S_i$.
Comparing with Rolfsen, $\nu_i$ points in the positive direction of a bicollar on $S_i$.
Hence by \cref{normal co-or}, $\beta_{\nu_i}$ is the normal orientation of $S_i$.

Now, considering $S_1$ and $S_2$ as co-oriented in this way, then according to \cref{P: normal pullback} the intersection $S_1 \times_{S^3} S_2$ is normally co-oriented by $\beta_{\nu_1} \wedge \beta_{\nu_2}$.
So the induced orientation on each component curve $C$ of this intersection is such that $\beta_C \wedge \beta_{\nu_1} \wedge \beta_{\nu_2}$ is the orientation of $S^3$.

We now choose a point where $L_1$ intersects $S_2$ and suppose $L_1$ passes through $S_2$ in the direction of the oriented bicollar of $S_2$.
Then we can identify $\beta_{L_1}$ with $\beta_{\nu_2}$, and we can identify $\beta_{C}$ with $\pm\beta_{\mu_1}$.
But now we know that
\begin{align*}
	\beta_{S^3} &= \beta_{S_1} \wedge \beta_{\nu_1}\\
				&= \beta_{\mu_1} \wedge \beta_{L_1} \wedge \beta_{\nu_1}\\
				&= \beta_{\mu_1} \wedge \beta_{\nu_2} \wedge \beta_{\nu_1}\\
				&= -\beta_{\mu_1} \wedge \beta_{\nu_1} \wedge  \beta_{\nu_2}.
\end{align*}
So we must have $\beta_C = -\beta_{\mu_1}$, and $C$ is oriented to move away from this intersection point.
On the other hand, if $L_1$ points opposite the orientation of the bicollar of $S_2$, $C$ will be oriented toward the intersection point.

We conclude that the linking number of Rolfsen corresponds to the number of curves of $S_1 \times_{S^3} S_2$ that are oriented away from $L_1$ minus the number of curves that are oriented toward $L_1$.
Of course any curve that travels from $L_1$ to $L_1$ is counted as a plus and a minus, so these cancel.
The result is that the linking number counts curves that travel from $L_1$ to $L_2$ minus the number of curves that travel from $L_2$ to $L_1$.
But this was precisely our claim.

\subsubsection{A triple linking number}\index{triple linking number|(}\index{link!triple linking number|(}

If one is given a link and chooses Seifert surfaces, it is probably simpler to compute the linking number using the technique in \cite{Ro76} then the geometric cohomology computation above.
However, geometric cohomology allows us to easily compute a more sophisticated link invariant - the triple linking number provided by Massey products.

We first recall in general from \cite{UM57} that if we have three cohomology classes $\underline{\alpha} \in H^i(X)$, $\underline{\beta} \in H^j(X)$, and $\underline{\gamma} \in H^k(X)$ for some space $X$ such that $\underline{\alpha} \sms \underline{\beta} = \underline{\beta} \sms \underline{\gamma} =0$, then if we choose representative cocycles $\alpha, \beta, \gamma$ there must be cochains $R,T$ such that $dR = \alpha \sms \beta$ and $dT = \beta \sms \gamma$.
Then the Massey product $\langle \underline{\alpha}, \underline{\beta}, \underline \gamma \rangle$ is represented by $R \sms \gamma +(-1)^{\deg(\alpha) +1} \alpha \sms T$, which represents a well-defined class in $H^{i+j+k-1}(X)/(\underline{\alpha} \sms H^{j+k-1}(X) + H^{i+j-1}(X) \sms \underline{\gamma})$.

Now consider a link of three components $L = L_1 \sqcup L_2 \sqcup L_3$ in $S^3$, and let $S_1$, $S_2$, and $S_3$ denote Seifert surfaces, co-oriented as in the preceding section.
We further suppose that $\underline{\mr S_1} \uplus \underline{\mr S_2} = \underline{\mr S_2} \uplus \underline{\mr S_3} = \underline{\mr S_1} \uplus \underline{\mr S_3} = 0 \in H^2_\Gamma(S^3-L)$.
So all linking numbers vanish, and the cup product alone cannot detect whether or not $L$ is trivially linked in the sense that it is isotopic to a link with the components contained in disjoint balls in $S^3$.
But in this case we have a Massey triple product $\langle \underline{\mr S_1}, \underline{\mr S_2}, \underline{\mr S_3}\rangle$.
Furthermore, $H^1_\Gamma(S^3 - L)$ is generated by the $\underline{\mr S_i}$, and $\underline{\mr S_i} \uplus \underline{\mr S_i} = 0$ as we can find two disjoint Seifert surfaces for any knot (for example, push $S_i$ along its bicollar).
So the cup product between elements of $H^1_\Gamma(S^3-L)$ is trivial, and $\langle \underline{\mr S_1}, \underline{\mr S_2}, \underline{\mr S_3}\rangle$ is a well-defined element of $H^2_\Gamma(S^3 - L)$.

Now if $L$ is trivially linked, we can choose Seifert surfaces such that the product $S_i \times_{S^3} S_j$, $i \neq j$, is empty, and so we can take $R = T = \emptyset$, and $\langle \underline{\mr S_1}, \underline{\mr S_2}, \underline{\mr S_3}\rangle = 0$.
So if, for some link with trivial linking numbers, $\langle \underline{\mr S_1}, \underline{\mr S_2}, \underline{\mr S_3}\rangle \neq 0$, then $L$ is not trivially linked, even though its linking numbers all vanish.

As an example, we will demonstrate the nontrivial linking of the classic example of the Borromean rings.
Let us consider the following realization of the Borromean rings:
In $\R^3$, let $L_1$ be the set $\{x^2+y^2 = 4, z=0\}$.
Let $L_2$ be the set $\{x^2+z^2/9 = 1, y=0\}$.
Let $L_3$ be the set $\{y^2/9 + z^2 = 1 , x=0\}$.
Each of these is an ellipse in a standard coordinate plane.
We let $S_i$ be the interior of $L_i$ in its coordinate plane.
For convenience, we assume the $L_i$ oriented so that the $S_i$ are normally co-oriented by the positive coordinate directions normal to the planes containing them; e.g.\ $S_1$ is the disk of radius $2$ centered at the origin in the $x$-$y$ plane, normally co-oriented by the positive $z$ direction, which we denote $\beta_z$.

In this case $S_1 \times_{S^3} S_2$ will be the arc $-1 \leq x \leq 1$ along the $x$-axis, normally co-oriented by $\beta_z \wedge \beta y$, and $S_2 \times_{S^3} S_3$ will be the arc $-1 \leq z \leq 1$ along the $z$-axis, normally co-oriented by $\beta_y \wedge \beta_x$.
The first arc has both endpoints on $L_2$ and the second has both endpoints on $L_3$, so both represent $0$ in $H^2_\Gamma(S^3-L)$ after restricting to their interiors.
Let us construct surfaces $S_{12}$ and $S_{23}$ such that $\mathring S_{12}$ and $\mathring S_{23}$  cobound $\mr S_1 \times_{S^3} \mr S_2$ and $\mr S_2 \times_{S^3} \mr S_3$.

For $S_{12}$, consider the intersection of $S_2$ with the upper half space $z \geq 0$.
Then, geometrically, the only boundary component of $\mr S_{12}$ is $\mathring S_1 \times_{S^3} \mathring S_2$.
Let us consider the co-orientations.
Giving $S_{12}$ the normal co-orientation $\beta_y$, so that the co-orientation is $(\beta_{S_{12}}, \beta_{S_{12}} \wedge \beta y)$ for any choice of $\beta_{S_{12}}$, and writing $C = \bd \mr S_{12}$, the inward pointing normal from $C$ into $S_{12}$ is $\beta_z$ and so the co-orientation of $C$ is the composition
\begin{align*}
(\beta_C, \beta_C \wedge \beta_z) * (\beta_{S_{12}}, \beta_{S^3}) &=(\beta_C, \beta_C \wedge \beta_z) * (\beta_{S_{12}}, \beta_{S_{12}} \wedge \beta_y)\\
&=(\beta_C, \beta_C \wedge \beta_z) * (\beta_C \wedge \beta_z, \beta_C \wedge \beta_z\wedge \beta_y)\\
&=(\beta_C,  \beta_C \wedge \beta_z \wedge \beta_y).
\end{align*}
So the normal co-orientation is $\beta_z \wedge \beta_y$, which agrees with the co-orientation of $\mr S_1 \times_{S^3} \mr S_2$.
So $\mr S_{12}$ cobounds $\mr S_1 \times_{S^3} \mr S_2$.

Similarly, for $S_{23}$, consider the intersection of $S_3$ with the half space $y \geq 0$.
Then, geometrically, the only boundary component of $\mr S_{23}$ is $\mathring S_2 \times_{S^3} \mathring S_3$.
Let us consider the co-orientations.
Giving $S_{23}$ the normal co-orientation $\beta_x$, so that the co-orientation is $(\beta_{S_{23}}, \beta_{S_{23}} \wedge \beta_x)$ for any choice of $\beta_{S_{23}}$, and writing $D = \bd \mr S_{23}$, the inward pointing normal from $D$ into $S_{23}$ is $\beta_y$ and so the co-orientation of $D$ is the composition
\begin{align*}
(\beta_D, \beta_D \wedge \beta_y) * (\beta_{S_{23}}, \beta_{S^3}) &=(\beta_D, \beta_D \wedge \beta_y) * (\beta_{S_{23}}, \beta_{S_{23}} \wedge \beta_x)\\
&=(\beta_D, \beta_D \wedge \beta_y) * (\beta_D \wedge \beta_y, \beta_D \wedge \beta_y \wedge \beta_x)\\
&=(\beta_D,  \beta_D \wedge \beta_y \wedge \beta_x).
\end{align*}
So the normal co-orientation is $\beta_y \wedge \beta_x$, which agrees with the normal co-orientation of $\mr S_2 \times_{S^3} \mr S_3$.
So $\mr S_{23}$ cobounds $\mr S_2 \times_{S^3} \mr S_3$.

Thus our Massey product $\langle \mr S_1, \mr S_2, \mr S_3 \rangle$ can be represented by $\mr S_{12} \times_{S^3} \mr S_3 + \mr S_1 \times_{S^3} \mr S_{23}$.
The piece $\mr S_{12} \times_{S^3} \mr S_3$ is the arc $0 \leq z < 1$ on the $z$-axis, normally co-oriented by $\beta_y \wedge \beta_x$, so the induced orientation is in the $-\beta_z$ direction.
The piece $\mr S_1 \times_{S^3} \mr S_{23}$ is the arc $0 \leq y < 2$ on the $y$-axis, normally co-oriented by $\beta_z \wedge \beta_x$, so the induced orientation is in the $\beta_y$ direction.
Altogether then, the product is represented by a path from $L_2$ to $L_1$, which represents a nontrivial class in $H^2_\Gamma(S^3-L)$.
This allows us to conclude that the Borromean rings are nontrivially linked.

\index{Alexander Duality!application to links|)}
\index{triple linking number|)}\index{link!triple linking number|)}


	\bibliographystyle{alpha}
	\bibliography{auxy/foundations} 

	\index{naive transversality|see{transversality, naive}}
	\index{fiber product|seealso{pullback}}
	\index{pullback|seealso{fiber product}}
	\index{local orientation|see{orientation, local}}
	\index{tautological co-orientation|see{co-orientation, tautological}}
	\index{normal co-orientation|see{co-orientation, normal}}
	\index{pullback co-orientation|see{co-orientation, of pullback}}
	\index{product|seealso{external products}}
	\index{product orientation|see{orientation!of external product}}
	\index{product co-orientation|see{co-orientation!of external product}}
	\index{orientation!cap|seealso{cap orientation}}
	\index{Mayer-Vietoris sequence!for geometric cohomology!connecting morphism|seealso{connecting morphism for geometric cohomology}}
	\index{simple transversality|see{transversality, simple}}
	\index{compound transversality|see{transversality, compound}}
	\index{cross product|see{exterior product}}

	\printindex


\end{document}